%% file: main.tex
\documentclass[11pt,reqno,oneside]{amsbook}

\input{preamble}

\begin{document}

\frenchspacing
\raggedbottom

\frontmatter

\thispagestyle{empty}
\pdfbookmark[1]{Title}{titlepage}
\begin{center}
\vspace*{1.6cm}

{\LARGE Non-Abelian Hirota--Miwa Equations\\[0.45em]
for the KPZ Universality Class\par}

\vspace{1.5cm}

{\large\scshape C. Alexander Rodriguez\par}

\vspace{0.8cm}

{\large August 3, 2026\par}

\vspace{1.2cm}

{\scshape Abstract\par}
\end{center}

\vspace{0.1cm}

\begin{center}
\begin{minipage}{0.86\textwidth}\small
\setlength{\parindent}{1.2em}
This work introduces an algebraic framework yielding explicit, closed matrix differential-difference equations for eighteen models in the exactly solvable sector of the KPZ universality class across four scaling regimes. By organizing Fredholm determinant data into an overdetermined linear problem on a directed lattice graph, we derive a compatibility system termed the diamond equations. Elementary seed data extracted from the shift structure of the Fredholm kernel provides simple solutions to this system.

We then construct a Darboux transformation to compress the infinite-dimensional Fredholm data into a finite-dimensional matrix observable. We show this dressing procedure preserves the diamond equations; consequently the resulting matrix observable obeys the same nonlinear structure as the initial seed data. Verifying a closed nonlinear equation for any specific model thus reduces to checking a handful of linear conditions on its kernel data.

Under a scalar reduction, the framework produces variable-coefficient Hirota--Miwa equations for Fredholm determinants, recovering the one-point bilinear equations of the author's earlier work as specializations. To supply the necessary seed data, a product graph construction with admissible propagators builds multipoint data in the fully discrete regime, while Euclidean division in a polynomial quotient algebra handles vertex and polymer models. Finally, we demonstrate the diamond equations are a gauge-equivalent reparametrization of the non-abelian Hirota--Miwa system, a central system in classical integrability theory.
\end{minipage}
\end{center}

\vfill

\newpage

\setcounter{tocdepth}{0} 
\setcounter{secnumdepth}{2} 
\tableofcontents

\mainmatter

\include{chapters/1-introduction}

\part{The Diamond Framework}
\include{chapters/2-the-discrete-framework}

\include{chapters/3-the-semi-discrete-framework}

\include{chapters/4-the-parabolic-framework}

\include{chapters/5-the-continuum-framework}

\part{Model Verifications}
\settocdepth{1}
\include{chapters/6-discrete-particle-models}

\include{chapters/7-euclidean-division-vertex-models-and-positive-polymer-models}

\include{chapters/8-semi-discrete-model-verifications}

\include{chapters/9-parabolic-model-verifications}

\settocdepth{0}
\include{chapters/10-the-kpz-fixed-point}

\backmatter

\bibliographystyle{amsalpha}
\bibliography{references}

\end{document}

%% file: preamble.tex
\usepackage[utf8]{inputenc}
\usepackage[T1]{fontenc}
\usepackage{lmodern}
\usepackage{microtype}

\usepackage[letterpaper,margin=1in]{geometry}

\usepackage[dvipsnames]{xcolor}

\usepackage[english]{babel}
\usepackage{enumitem}
\usepackage{amsmath,amssymb}
\allowdisplaybreaks[4] 
\usepackage{thmtools}
\usepackage{graphicx}
\usepackage{xspace}
\usepackage{calc}

\usepackage{mathtools}
\usepackage{amsfonts}
\usepackage{bm}
\usepackage{latexsym}
\usepackage{mathrsfs}
\usepackage{upgreek}
\usepackage{dirtytalk}
\usepackage{multicol}
\usepackage{array}
\usepackage{tabularx}
\usepackage{booktabs}
\usepackage{subcaption}
\usepackage{tikz}
\usetikzlibrary{shapes,backgrounds}
\usetikzlibrary{arrows.meta,calc,positioning}
\usepackage{tikz-cd}
\usepackage{tikzsymbols}

\usepackage{hyperref}
\hypersetup{%
  colorlinks=true, linktocpage=true,%
  breaklinks=true, pageanchor=true,%
  plainpages=false, bookmarksnumbered, bookmarksopen=true, bookmarksopenlevel=1,%
  hypertexnames=true,%
  linkcolor=MidnightBlue, citecolor=MidnightBlue, urlcolor=MidnightBlue,%
  pdftitle={Non-abelian Hirota--Miwa Equations for the KPZ Universality Class},%
  pdfauthor={C. Alexander Rodriguez},%
  pdfkeywords={KPZ universality, integrable probability, Hirota--Miwa, Fredholm determinants, Darboux transformations}%
}

\makeatletter
\@ifpackageloaded{babel}%
  {%
    \addto\extrasenglish{%
    }%
    }{\relax}
\makeatother

\newtheorem{theorem}{Theorem}[chapter]
\newtheorem{corollary}[theorem]{Corollary}
\newtheorem{lemma}[theorem]{Lemma}
\newtheorem{proposition}[theorem]{Proposition}

\theoremstyle{definition}
\newtheorem{definition}[theorem]{Definition}

\newtheorem{example}[theorem]{Example}

\theoremstyle{remark}
\newtheorem{remark}[theorem]{Remark}

\newcommand{\N}{\mathbb{N}}
\newcommand{\PP}{\mathbb{P}}
\newcommand{\E}{\mathbb{E}}
\newcommand{\Z}{\mathbb{Z}}
\newcommand{\R}{\mathbb{R}}
\providecommand{\C}{}\renewcommand{\C}{\mathbb{C}}

\newcommand{\hh}{\mathfrak{h}}

\newcommand{\defeq}{\vcentcolon=}
\newcommand{\End}{\mathrm{End}}
\providecommand{\Hom}{}\renewcommand{\Hom}{\mathrm{Hom}}

\newcommand{\pa}{\partial}

\newcommand{\diff}{\mathrm{d}}

\newcommand{\vertiii}[1]{{\left\vert\kern-0.25ex\left\vert\kern-0.25ex\left\vert #1 \right\vert\kern-0.25ex\right\vert\kern-0.25ex\right\vert}}

\newcounter{bln}

\newenvironment{itemize*}
  {\begin{itemize}[topsep=-\parskip+\jot,itemsep=-\parskip-\jot]}
  {\end{itemize}}

\newenvironment{enumerate*}
  {\begin{enumerate}[label=(\alph*),topsep=-\parskip+\jot,itemsep=-\parskip-\jot]}
  {\end{enumerate}}

\newenvironment{enumerate**}
  {\begin{enumerate}[label=(\roman*),topsep=-\parskip+\jot,itemsep=-\parskip-\jot]}
  {\end{enumerate}}

\newenvironment{enumerate***}
  {\begin{enumerate}[label=(\alph*'),topsep=-\parskip+\jot,itemsep=-\parskip-\jot]}
  {\end{enumerate}}

\numberwithin{equation}{chapter}
\numberwithin{figure}{chapter}

\makeatletter
\newcommand{\settocdepth}[1]{%
  \immediate\write\@mainaux{\string\@writefile{toc}{\string\setcounter{tocdepth}{#1}}}%
}
\makeatother

%% file: chapters/1-introduction.tex
\chapter{Introduction}
\label{ch:introduction}

{
  \setlength{\parskip}{0pt}
}

\section{The KPZ universality class}\label{sec:KPZ-universality}

The KPZ universality class is a broad collection of
one-dimensional, asymmetric, randomly forced mathematical and
physical models linked through their shared universal scaling
behaviour.  This collection includes stochastic interface growth
on a one-dimensional substrate, directed polymer chains in a random
potential, driven lattice gas models, reaction-diffusion models in
two-dimensional random media, and randomly forced Hamilton--Jacobi
equations (see, e.g.~\cite{BarabasiStanley1995,
HalpinHealyZhang1995, Krug1997, Meakin1998, Spohn2012, Corwin2012,
QuastelSpohn2015, HalpinHealyTakeuchi2015}).  Despite disparate
microscopic descriptions, these models exhibit common non-Gaussian
asymptotic fluctuations with, for instance, interfaces moving at a
velocity proportional to time $t$, fluctuations of size $t^{1/3}$,
and decorrelation at a spatial scale of $t^{2/3}$.  The recognition
that these models share common scaling behaviour, and the proposal
of an equation to represent it, both date to the seminal 1986 paper
that gave the class its name~\cite{KPZ86}.

Kardar, Parisi, and Zhang introduced the stochastic partial
differential equation
\begin{equation}\label{eq:KPZ}
\partial_t h = \nu \partial_x^2 h
+ \lambda(\partial_x h)^2 + \sigma\xi,
\end{equation}
where $\xi$ is space-time white noise, as a paradigmatic model of
random interface growth in one dimension.  While each term captures
a distinct physical mechanism (diffusive smoothing, slope-dependent
lateral growth, and random forcing), it is the nonlinearity
$\lambda(\partial_x h)^2$, coupling growth to the local slope, that
drives the system away from Gaussian (Edwards--Wilkinson) behaviour.  Kardar, Parisi, and Zhang argued, building on the
renormalization-group analysis of Forster, Nelson, and
Stephen~\cite{FNS1977}, that~\eqref{eq:KPZ} exhibits characteristic
scaling exponents: fluctuations of height on the order $t^{1/3}$
and spatial correlations at scale $t^{2/3}$, and that these
exponents should be universal across a broad class of
one-dimensional growth processes.

As evidence for the scaling exponents predicted in~\cite{KPZ86}
accumulated over the decades that followed, a more ambitious
distributional picture gradually emerged: that the height
fluctuations of every model in the class, rescaled under the 1:2:3
scaling
\[
h^{\epsilon}(t, x) = \epsilon^{1/2}
\bigl(h(\epsilon^{-3/2}t, \epsilon^{-1}x) - c_\epsilon t\bigr),
\]
should converge as
$\epsilon \to 0$ to a single universal limit process, and that a
universal random metric should similarly govern metric and polymer
models.  These
conjectural limits are now known as the KPZ fixed
point~\cite{MQR21} and the directed
landscape~\cite{DOV22}, and they represent the central universal
content of the class.  Proving rigorously that they are the universal limits for the full
breadth of the class remains the field's defining open problem.

By the turn of the millennium, the KPZ scaling picture had been
established beyond serious doubt at the level of physical
evidence~\cite{HalpinHealyZhang1995, Krug1997, HalpinHealyTakeuchi2015}.  Yet the KPZ equation itself
had resisted rigorous mathematical treatment for nearly two decades.
Moreover, this intractability was not merely technical: the universal
scaling limits the equation was meant to describe could not be
extracted from an object that itself could not be rigorously
defined.  While a rigorous finite-time solution theory was eventually given
by Hairer~\cite{Hairer2013, Hairer2014} (see
also~\cite{BertiniGiacomin1997}), well-posedness is not
universality, and rigorous progress toward the universal scaling
limits had to come from a different direction.

Beginning in the late 1990s, several mathematical breakthroughs
were achieved through the study of a special subset of models within
the class, known as \emph{exactly solvable} models.  These models,
while representing specific points in the vast parameter space of
the universality class, possess a deep and rigid algebraic
structure, often rooted in representation theory or quantum
integrable systems, which allows for the exact computation of
physical observables.  It is through the asymptotic analysis of
these exact formulas that the universal scaling functions and limit
distributions of the KPZ class have been rigorously derived.  These
exactly solvable models, therefore, illuminate the universal
properties that are conjectured to hold for the entire class.

Yet for all the successes this program has produced, its resulting
formulas have often appeared fragmented and isolated by the specific
algebraic machinery used to obtain them: interacting particle
systems via the Bethe ansatz and biorthogonal ensembles,
zero-temperature last-passage percolation through the
Robinson--Schensted--Knuth correspondence, stochastic vertex models
through the Yang--Baxter equation and Macdonald processes, and
positive-temperature directed polymers through Whittaker functions
and geometric representation
theory \cite{TW2009, Corwin2012, BorodinCorwin2014,
CorwinOConnellSeppalainenZygouras2014, BP2016, MQR21, DOV22,
AggarwalCorwinHegde2024}.  It has been natural, then, to view the landscape
of exact solvability in the KPZ class as a collection of individual
achievements; a zoo of models, methods, and formulas, rather than
the output of a single organizing principle.  Whether a deeper
structure connects these results, a single algebraic framework of
the sort that organizes other domains of mathematical physics, is a
question at the heart of the present work.
\section{Overview of the main results}\label{sec:overview}
\definecolor{tileA}{RGB}{0,119,136}
\definecolor{tileB}{RGB}{204,68,0}
\definecolor{tileC}{RGB}{102,45,145}
\definecolor{tileD}{RGB}{70,90,120}

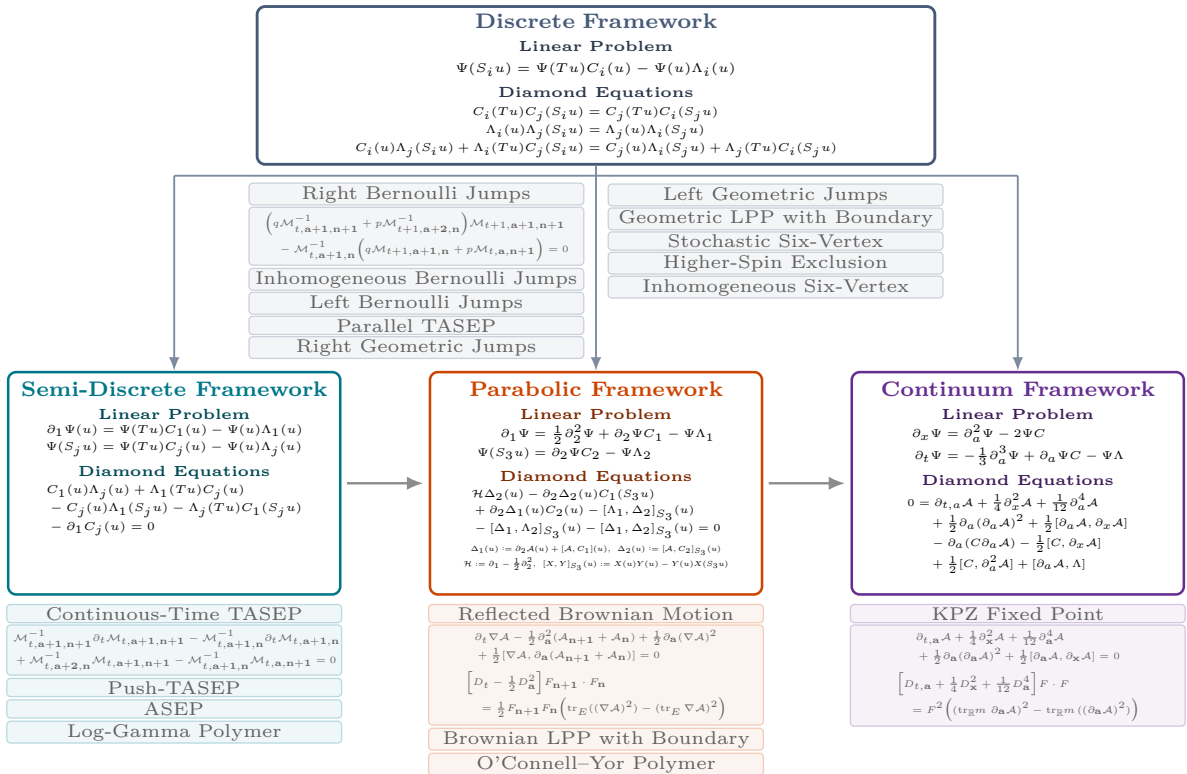
\begin{figure}[!htb]
\centering
\begingroup
\def\chipfont{\fontsize{5.8}{6.2}\selectfont}
\def\boxlab#1#2{{\fontsize{6.6}{7.4}\selectfont\bfseries\color{#1}#2}}
\def\subhead#1#2{{\fontsize{4.6}{5.0}\selectfont\bfseries\color{#1!62!black}#2}}
\tikzset{
  pbox/.style={draw=tileD, line width=0.9pt, rounded corners=2.5pt,
               inner xsep=4pt, inner ysep=2.2pt, align=center, text width=7.19cm},
  cbox/.style={draw=#1, line width=0.8pt, rounded corners=2.5pt,
               inner xsep=3pt, inner ysep=3pt},
  chipbase/.style={rounded corners=1.2pt, inner xsep=2pt, inner ysep=0.85pt,
                   align=center, font=\chipfont, line width=0.4pt, text width=3.55cm},
  pchip/.style={chipbase, draw=tileD!28, fill=tileD!6, text=black!62},
  chipA/.style={chipbase, draw=tileA!28, fill=tileA!6, text=black!62},
  chipB/.style={chipbase, draw=tileB!28, fill=tileB!6, text=black!62},
  chipC/.style={chipbase, draw=tileC!28, fill=tileC!6, text=black!62},
  repA/.style={chipA},
  repB/.style={chipB},
  repC/.style={chipC},
  eqbox/.style={rounded corners=1.2pt, inner xsep=2pt, inner ysep=1.4pt,
                align=center, line width=0.4pt, text width=3.55cm},
  eqboxA/.style={eqbox, draw=tileA!28, fill=tileA!6, text=black!62},
  eqboxB/.style={eqbox, draw=tileB!28, fill=tileB!6, text=black!62},
  eqboxC/.style={eqbox, draw=tileC!28, fill=tileC!6, text=black!62},
  eqboxD/.style={eqbox, draw=tileD!28, fill=tileD!6, text=black!62},
  bus/.style={draw=tileD!70, line width=0.6pt},
  down/.style={-{Latex[length=1.5mm,width=1.2mm]}, draw=tileD!70, line width=0.6pt},
  scal/.style={-{Latex[length=1.8mm,width=1.5mm]}, draw=black!45, line width=0.75pt},
}
\scalebox{1.20}{%
\begin{tikzpicture}[x=1cm,y=1cm,every node/.style={font=\fontsize{5}{5.8}\selectfont}]

\node[pbox, anchor=north] (P) at (0,0) {
  \boxlab{tileD}{Discrete Framework}\\[1.2pt]
  \subhead{tileD}{Linear Problem}\\[1.0pt]
  \scalebox{0.80}{$\Psi(S_i u) = \Psi(Tu)C_i(u) - \Psi(u)\Lambda_i(u)$}\\[2.0pt]
  \subhead{tileD}{Diamond Equations}\\[1.0pt]
  \scalebox{0.72}{$\setlength{\jot}{0.6pt}\begin{aligned}
     C_i(Tu) C_j(S_i u) &= C_j(Tu) C_i(S_j u)\\
     \Lambda_i(u) \Lambda_j(S_i u) &= \Lambda_j(u) \Lambda_i(S_j u)\\
     C_i(u) \Lambda_j(S_i u) + \Lambda_i(Tu) C_j(S_i u)
       &= C_j(u) \Lambda_i(S_j u) + \Lambda_j(Tu) C_i(S_j u)
   \end{aligned}$}
};

\coordinate (bus) at ($(P.south)+(0,-0.11)$);
\draw[bus] (P.south) -- (bus);
\draw[bus] ($(bus)+(-4.65,0)$) -- ($(bus)+(4.65,0)$);

\node[pchip, anchor=north] (pL1) at ($(bus)+(-1.98,-0.075)$) {Right Bernoulli Jumps};
\node[eqboxD, anchor=north] (pE) at ($(pL1.south)+(0,-0.012)$) {%
  \scalebox{0.60}{$\setlength{\jot}{0.6pt}\begin{aligned}
    & \bigl(q\mathcal{M}_{t,\mathbf{a}+\mathbf{1},\mathbf{n}+\mathbf{1}}^{-1}
      + p\mathcal{M}_{t+1,\mathbf{a}+\mathbf{2},\mathbf{n}}^{-1}\bigr)
      \mathcal{M}_{t+1,\mathbf{a}+\mathbf{1},\mathbf{n}+\mathbf{1}}\\
    & \quad - \mathcal{M}_{t,\mathbf{a}+\mathbf{1},\mathbf{n}}^{-1}
      \bigl(q\mathcal{M}_{t+1,\mathbf{a}+\mathbf{1},\mathbf{n}}
      + p\mathcal{M}_{t,\mathbf{a},\mathbf{n}+\mathbf{1}}\bigr) = 0
  \end{aligned}$}};
\node[pchip, anchor=north] (pL2) at ($(pE.south)+(0,-0.022)$) {Inhomogeneous Bernoulli Jumps};
\node[pchip, anchor=north] (pL3) at ($(pL2.south)+(0,-0.010)$) {Left Bernoulli Jumps};
\node[pchip, anchor=north] (pL4) at ($(pL3.south)+(0,-0.010)$) {Parallel TASEP};
\node[pchip, anchor=north] (pL5) at ($(pL4.south)+(0,-0.010)$) {Right Geometric Jumps};
\coordinate (pR0) at ($(bus)+(1.98,-0.075)$);
\foreach[count=\i, evaluate=\i as \j using int(\i-1)] \nm in {%
  {Left Geometric Jumps},
  {Geometric LPP with Boundary},
  {Stochastic Six-Vertex},
  {Higher-Spin Exclusion},
  {Inhomogeneous Six-Vertex}}{%
  \node[pchip, anchor=north] (pR\i) at ($(pR\j.south)+(0,-0.010)$) {\nm};
}

\coordinate (row) at ($(pL5.south)+(1.98,-0.13)$);

\node[cbox=tileA, anchor=north] (SD) at ($(row)+(-4.65,0)$) {
  \parbox[c][2.24cm][t]{3.45cm}{\centering
  \boxlab{tileA}{Semi-Discrete Framework}\\[1.2pt]
  \subhead{tileA}{Linear Problem}\\[0.8pt]
  \scalebox{0.72}{$\setlength{\jot}{0.6pt}\begin{aligned}
     \partial_1\Psi(u) &= \Psi(Tu)C_1(u) - \Psi(u)\Lambda_1(u)\\
     \Psi(S_j u) &= \Psi(Tu)C_j(u) - \Psi(u)\Lambda_j(u)
   \end{aligned}$}\\[2.0pt]
  \subhead{tileA}{Diamond Equations}\\[0.8pt]
  \scalebox{0.72}{$\setlength{\jot}{0.6pt}\begin{aligned}
     & C_1(u)\Lambda_j(u) + \Lambda_1(Tu)C_j(u)\\
     & - C_j(u)\Lambda_1(S_j u) - \Lambda_j(Tu)C_1(S_j u)\\
     & - \partial_1 C_j(u) = 0
   \end{aligned}$}}
};

\node[cbox=tileB, anchor=north] (PA) at ($(row)+(0,0)$) {
  \parbox[c][2.24cm][t]{3.45cm}{\centering
  \boxlab{tileB}{Parabolic Framework}\\[1.2pt]
  \subhead{tileB}{Linear Problem}\\[0.8pt]
  \scalebox{0.72}{$\setlength{\jot}{0.6pt}\begin{aligned}
     \partial_1 \Psi &= \tfrac{1}{2}\partial_2^2 \Psi + \partial_2\Psi C_1 - \Psi\Lambda_1\\
     \Psi(S_3 u) &= \partial_2\Psi C_2 - \Psi\Lambda_2
   \end{aligned}$}\\[2.0pt]
  \subhead{tileB}{Diamond Equations}\\[0.8pt]
  \scalebox{0.66}{$\setlength{\jot}{0.6pt}\begin{aligned}
     & \mathcal{H}\Delta_2(u) - \partial_2\Delta_2(u)C_1(S_3u)\\
     & + \partial_2\Delta_1(u)C_2(u) - [\Lambda_1,\Delta_2]_{S_3}(u)\\
     & - [\Delta_1,\Lambda_2]_{S_3}(u) - [\Delta_1,\Delta_2]_{S_3}(u) = 0
   \end{aligned}$}\\[1.4pt]
  \scalebox{0.46}{$\setlength{\jot}{0.6pt}\begin{gathered}
     \Delta_1(u) \defeq \partial_2\mathcal{A}(u) + [\mathcal{A}, C_1](u),\enspace
     \Delta_2(u) \defeq [\mathcal{A}, C_2]_{S_3}(u)\\
     \mathcal{H} \defeq \partial_1 - \tfrac{1}{2}\partial_2^2,\enspace
     [X,Y]_{S_3}(u) \defeq X(u)Y(u) - Y(u)X(S_3 u)
   \end{gathered}$}}
};

\node[cbox=tileC, anchor=north] (CO) at ($(row)+(4.65,0)$) {
  \parbox[c][2.24cm][t]{3.45cm}{\centering
  \boxlab{tileC}{Continuum Framework}\\[1.2pt]
  \subhead{tileC}{Linear Problem}\\[0.8pt]
  \scalebox{0.72}{$\setlength{\jot}{0.6pt}\begin{aligned}
     \partial_x\Psi &= \partial_a^2\Psi - 2\Psi C\\
     \partial_t\Psi &= -\tfrac{1}{3}\partial_a^3\Psi + \partial_a\Psi C - \Psi\Lambda
   \end{aligned}$}\\[2.0pt]
  \subhead{tileC}{Diamond Equations}\\[0.8pt]
  \scalebox{0.66}{$\setlength{\jot}{0.6pt}\begin{aligned}
     0 &= \partial_{t,a}\mathcal{A} + \tfrac{1}{4}\partial_x^2\mathcal{A}
          + \tfrac{1}{12}\partial_a^4\mathcal{A}\\
       &\quad + \tfrac{1}{2}\partial_a(\partial_a\mathcal{A})^2
          + \tfrac{1}{2}[\partial_a\mathcal{A}, \partial_x\mathcal{A}]\\
       &\quad - \partial_a(C\partial_a\mathcal{A})
          - \tfrac{1}{2}[C, \partial_x\mathcal{A}]\\
       &\quad + \tfrac{1}{2}[C, \partial_a^2\mathcal{A}]
          + [\partial_a\mathcal{A}, \Lambda]
   \end{aligned}$}}
};

\draw[down] ($(bus)+(-4.65,0)$) -- (SD.north);
\draw[down] (bus) -- (PA.north);
\draw[down] ($(bus)+(4.65,0)$) -- (CO.north);

\draw[scal] ($(SD.east)+(0.06,0)$) -- ($(PA.west)+(-0.06,0)$);
\draw[scal] ($(PA.east)+(0.06,0)$) -- ($(CO.west)+(-0.06,0)$);

\node[repA,  anchor=north] (a1) at ($(SD.south)+(0,-0.09)$) {Continuous-Time TASEP};
\node[eqboxA, anchor=north] (aE) at ($(a1.south)+(0,-0.012)$) {%
  \scalebox{0.60}{$\setlength{\jot}{0.6pt}\begin{aligned}
    & \mathcal{M}_{t,\mathbf{a}+\mathbf{1},\mathbf{n}+\mathbf{1}}^{-1}
      \partial_t \mathcal{M}_{t,\mathbf{a}+\mathbf{1},\mathbf{n}+\mathbf{1}}
      - \mathcal{M}_{t,\mathbf{a}+\mathbf{1},\mathbf{n}}^{-1}
      \partial_t \mathcal{M}_{t,\mathbf{a}+\mathbf{1},\mathbf{n}}\\
    & + \mathcal{M}_{t,\mathbf{a}+\mathbf{2},\mathbf{n}}^{-1}
      \mathcal{M}_{t,\mathbf{a}+\mathbf{1},\mathbf{n}+\mathbf{1}}
      - \mathcal{M}_{t,\mathbf{a}+\mathbf{1},\mathbf{n}}^{-1}
      \mathcal{M}_{t,\mathbf{a},\mathbf{n}+\mathbf{1}} = 0
  \end{aligned}$}};
\node[chipA, anchor=north] (a2) at ($(aE.south)+(0,-0.022)$) {Push-TASEP};
\node[chipA, anchor=north] (a3) at ($(a2.south)+(0,-0.022)$) {ASEP};
\node[chipA, anchor=north] (a4) at ($(a3.south)+(0,-0.022)$) {Log-Gamma Polymer};

\node[repB,  anchor=north] (b1) at ($(PA.south)+(0,-0.09)$) {Reflected Brownian Motion};
\node[eqboxB, anchor=north] (bE) at ($(b1.south)+(0,-0.012)$) {%
  \scalebox{0.60}{$\setlength{\jot}{0.6pt}\begin{aligned}
    & \partial_t \nabla\mathcal{A}
      - \tfrac{1}{2}\partial_{\mathbf{a}}^2
      (\mathcal{A}_{\mathbf{n}+\mathbf{1}} + \mathcal{A}_{\mathbf{n}})
      + \tfrac{1}{2}\partial_{\mathbf{a}}(\nabla\mathcal{A})^2\\
    & + \tfrac{1}{2}[\nabla\mathcal{A},
      \partial_{\mathbf{a}}(\mathcal{A}_{\mathbf{n}+\mathbf{1}}
      + \mathcal{A}_{\mathbf{n}})] = 0
  \end{aligned}$}\\[1.6pt]
  \scalebox{0.60}{$\setlength{\jot}{0.6pt}\begin{aligned}
    & \bigl[D_t - \tfrac{1}{2}D_{\mathbf{a}}^2\bigr]
      F_{\mathbf{n}+\mathbf{1}} \cdot F_{\mathbf{n}}\\
    & \quad = \tfrac{1}{2}F_{\mathbf{n}+\mathbf{1}} F_{\mathbf{n}}
      \bigl(\operatorname{tr}_E((\nabla\mathcal{A})^2)
      - (\operatorname{tr}_E \nabla\mathcal{A})^2\bigr)
  \end{aligned}$}};
\node[chipB, anchor=north] (b2) at ($(bE.south)+(0,-0.022)$) {Brownian LPP with Boundary};
\node[chipB, anchor=north] (b3) at ($(b2.south)+(0,-0.022)$) {O'Connell--Yor Polymer};

\node[repC, anchor=north] (c1) at ($(CO.south)+(0,-0.09)$) {KPZ Fixed Point};
\node[eqboxC, anchor=north] (cE) at ($(c1.south)+(0,-0.012)$) {%
  \scalebox{0.60}{$\setlength{\jot}{0.6pt}\begin{aligned}
    & \partial_{t,\mathbf{a}} \mathcal{A}
      + \tfrac{1}{4}\partial_{\mathbf{x}}^2 \mathcal{A}
      + \tfrac{1}{12}\partial_{\mathbf{a}}^4 \mathcal{A}\\
    & + \tfrac{1}{2}\partial_{\mathbf{a}}(\partial_{\mathbf{a}}\mathcal{A})^2
      + \tfrac{1}{2}[\partial_{\mathbf{a}}\mathcal{A},
      \partial_{\mathbf{x}}\mathcal{A}] = 0
  \end{aligned}$}\\[1.6pt]
  \scalebox{0.60}{$\setlength{\jot}{0.6pt}\begin{aligned}
    & \bigl[D_{t,\mathbf{a}} + \tfrac{1}{4}D_{\mathbf{x}}^2
      + \tfrac{1}{12}D_{\mathbf{a}}^4\bigr] F \cdot F\\
    & \quad = F^2\bigl((\operatorname{tr}_{\mathbb{R}^m}\partial_{\mathbf{a}}\mathcal{A})^2
      - \operatorname{tr}_{\mathbb{R}^m}((\partial_{\mathbf{a}}\mathcal{A})^2)\bigr)
  \end{aligned}$}};

\end{tikzpicture}}
\endgroup
\caption{The four frameworks and the models they govern.  Each
regime is formally a scaling limit of its predecessor, but the
frameworks are constructed independently; beneath each, a
representative model displays its multipoint equation.  The
catalogue of all eighteen equations appears in
\S\ref{sec:catalogue}.}
\label{fig:intro-regimes}
\end{figure}
The main results of this work are as follows.  We introduce
an algebraic framework, across four distinct scaling
regimes, that yields explicit, closed matrix
differential-difference equations for eighteen models in the
exactly solvable sector of the KPZ universality
class (\S\ref{sec:catalogue}).  We achieve this by
showing that the known Fredholm determinant data of these
models can be organized by an overdetermined linear problem on
a directed lattice graph, and as such embed the Fredholm data
into a general theory to provide a unified method for deriving
our matrix equations.

Our general theory is developed in four distinct regimes, each
formally a scaling limit of its predecessor, but nonetheless
each constructed independently as a standalone theory
(Figure~\ref{fig:intro-regimes}).  In
each regime, we construct an abstract overdetermined linear
problem whose compatibility conditions yield a system we call
the \emph{diamond equations}.  Our first observation is that
elementary data extracted from each model's Fredholm kernel,
such as the coefficients of its shift structure, already give
simple solutions to the diamond equations.  From these simple solutions, we construct what we call a
\emph{Darboux transformation}: a dressing procedure
that builds, using a \emph{dressing-compatible} kernel and its
operator-valued wave functions, a compressed
finite-dimensional matrix-valued observable.  As a hallmark of
integrability, we directly show that this dressing procedure
preserves the diamond equations, and so, from the simple
solutions of the diamond equations, it follows that our
compressed matrix observable itself again satisfies the same
diamond equations.  The consequence is that verifying a closed nonlinear equation
for any given model reduces entirely to checking a handful of
linear conditions on its Fredholm data: once these hold, the
Darboux theorem supplies the nonlinear equation.

The results recover the one-point equations
of~\cite{Rodriguez2025} as a special case.  When the matrix
observable is a scalar, the dressed diamond equations reduce,
under mild boundary conditions, to variable-coefficient
Hirota--Miwa equations for the Fredholm determinant.  The six
one-point bilinear equations derived in~\cite{Rodriguez2025}
are specializations of these, with coefficients determined by
each model's scalar seed data.

The diamond equations are not a new integrable system: for
invertible $C$-weights they are a gauge-equivalent
asymmetric reparametrization of the non-abelian
Hirota--Miwa system of Nimmo~\cite{Nimmo2006}.  The gauge connecting the two is
nonlocal: the asymmetric shift structure of each model's
Fredholm kernel falls directly into the diamond
parametrization, while the symmetric NAHM form is accessible
only by integrating this data along entire lattice paths.
The equivalence provides, nonetheless, a two-way conduit
between the exactly solvable sector of the KPZ universality
class and the classical theory of the Hirota--Miwa system.

We close this overview with several contributions outside the
core framework.  In the fully discrete regime, we develop a
mechanism, the product graph construction, which abstractly
builds dressing-compatible kernels, and hence solutions of the
diamond equations, from trivial linear data.  This is achieved
by algebraically abstracting a class of admissible propagators,
which allow one to lift typical one-point data into multipoint
data connected through the propagator.  For the vertex and polymer models of
\S\ref{sec:catalogue}, a different mechanism applies:
we show that on a polynomial quotient algebra, Euclidean
division naturally defines matrices satisfying the diamond
equations, providing the necessary seed data for the Darboux
transformation.  At the continuum level, we not only recover the matrix KP
equation for the KPZ fixed point, but show that the multipoint
distribution function itself satisfies the KP equation with a
forcing term built from our matrix observable.

The remainder of this introduction is organized as follows.
Section~\ref{sec:early-breakthroughs} surveys the
exact-solvability program in the KPZ class and the appearance
of Painlev\'e~II in its limiting distributions;
\S\ref{sec:classical-integrability} follows classical
integrability from the solitary wave to the Hirota--Miwa
equation; \S\ref{sec:equations-state-of-art} reviews the
closed equations known for KPZ distribution functions at the
outset of this work; and \S\ref{sec:one-point-equations}
recounts the one-point equations of~\cite{Rodriguez2025}, the
direct predecessor of the present framework.  Precise
statements of the main results, together with the catalogue of
all eighteen equations, occupy \S\ref{sec:main-results}.

\section{Early mathematical breakthroughs}\label{sec:early-breakthroughs}

The program of exact computation in the KPZ class began in 1999
with a classical problem in probabilistic combinatorics, going
back to Ulam in the
1960s.\footnote{Ulam~\cite{Ulam1961} raised the question of the
expected length of the longest increasing subsequence in the
1960s; Hammersley~\cite{Hammersley1972}, Logan and
Shepp~\cite{LoganShepp1977}, and Vershik and
Kerov~\cite{VershikKerov1977} progressively established that
$\E l_N \sim 2\sqrt{N}$, leaving the fluctuations as the central
open problem that Baik, Deift, and Johansson resolved.}  Baik,
Deift, and Johansson~\cite{BDJ1999} showed that the length $l_N$
of the longest increasing subsequence in a uniformly random
permutation of $\{1, \dotsc, N\}$ has fluctuations governed by an
explicit limiting distribution $F_2$,
\begin{equation}\label{eq:BDJ-limit}
\PP\left(\frac{l_N - 2\sqrt{N}}{N^{1/6}} \leq s\right)
\longrightarrow F_2(s) \qquad \text{as } N \to \infty.
\end{equation}
A direct mapping identifies $l_N$ with the zero-temperature free
energy of a directed polymer in a random environment, and through
last passage percolation the problem is connected to random growth
and interacting particle systems~\cite{Johansson2000, Romik2015}.
The proof combined a determinantal formula of
Gessel~\cite{Gessel1990} with the Deift--Zhou steepest descent
method for Riemann--Hilbert problems~\cite{DeiftZhou1993}.

The result was a landmark on several counts.  It demonstrated that
the large-$N$ asymptotics of the longest increasing subsequence
problem admitted a closed-form description, an achievement rarely
accessible for large-scale random systems far from equilibrium.
The function $F_2$, moreover, was not new.  Just five years
earlier, Tracy and Widom~\cite{TW1994} had identified it as the
limit law of the largest eigenvalue of a random matrix from the
Gaussian Unitary Ensemble, a completely different probabilistic
setting.  That the same distribution should govern both random
matrix eigenvalues and random growth was unexpected, and signalled
a connection whose depth would only become apparent over the
following decades.

Most significant for the perspective of this work is the form of
the Tracy--Widom distribution itself.  Tracy and
Widom~\cite{TW1994} had shown that
\begin{equation}\label{eq:TW-PII}
F_2(s) = \exp\biggl(
-\int_s^\infty (x - s) q(x)^2 dx\biggr),
\end{equation}
where $q$ is the Hastings--McLeod~\cite{HastingsMcLeod1980} solution
of the Painlev\'e~II equation,
\begin{equation}\label{eq:PainleveII}
q''(s) = sq(s) + 2q(s)^3,
\qquad q(s) \sim \operatorname{Ai}(s)
\text{ as } s \to +\infty.
\end{equation}
The Painlev\'e equations are among the central nonlinear ODEs of
mathematical physics.  Identified at the turn of the twentieth
century in the classification program of
Painlev\'e \cite{Painleve1900, Painleve1902} and
Gambier \cite{Gambier1910}, their connection to soliton theory
emerged in the 1970s, with Ablowitz and
Segur \cite{AblowitzSegur1977} deriving Painlev\'e~II as a
similarity reduction of the Korteweg--de Vries equation.  Jimbo,
Miwa, and Ueno \cite{JimboMiwaUeno1981} identified the
Painlev\'e equations as isomonodromic deformation equations for
linear systems with rational coefficients, placing them at the
center of the modern theory of integrable
systems \cite{ItsNovokshenov1986, FIKN2006, Clarkson2003}.  The
appearance
of a Painlev\'e transcendent in the fluctuation statistics of a
random growth model was a first, striking indication that the
distribution functions of KPZ models carry the imprint of
classical integrability.

Over the decade that followed, a rapid succession of results
delivered exact distribution formulas and Tracy--Widom limit laws
across the KPZ universality class:
Johansson~\cite{Johansson2000} established the Tracy--Widom limit
for TASEP through last passage percolation;
Pr\"ahofer and Spohn~\cite{PrahoferSpohn2002} showed that the
spatial fluctuation process of the polynuclear growth model
converges to the Airy$_2$ process, the first spatial limit process
in the class; Tracy and Widom~\cite{TW2008, TW2009} obtained exact
formulas for the full asymmetric simple exclusion process (ASEP)
via the Bethe ansatz; and Amir, Corwin, and
Quastel~\cite{ACQ2011} and Sasamoto and
Spohn~\cite{SasamotoSpohn2010} independently computed the
one-point distribution of the KPZ equation~\eqref{eq:KPZ} itself
with narrow wedge initial data, confirming Tracy--Widom
fluctuations for the continuum equation.  Running alongside these
advances, the Borodin--Corwin theory of Macdonald
processes~\cite{BorodinCorwin2014} provided a broad algebraic
framework encompassing many of the known solvable models.

The culmination of this program came with the construction of the
KPZ fixed point by Matetski, Quastel, and Remenik~\cite{MQR21}.
Starting from exact Fredholm determinant formulas for the
transition probabilities of TASEP with arbitrary initial data, they
carried out a scaling limit and obtained a well-defined Markov
process, the conjectural universal scaling limit of the entire KPZ
class.  The construction moreover established that the universal
object itself admits exact Fredholm determinant formulas for its
finite-dimensional distributions.  The Tracy--Widom distribution
$F_2$ that had appeared in Baik--Deift--Johansson is, from this
vantage point, the one-point marginal of the KPZ fixed point for
narrow wedge initial data, the most basic of the three canonical
initial conditions.  Shortly after, the directed landscape of
Dauvergne, Ortmann, and Vir\'ag~\cite{DOV22} provided a
complementary universal object, a random directed metric on the
plane conjectured to be the universal scaling limit for last
passage percolation and directed polymer models.

Yet the appearance of Painlev\'e~II in the Tracy--Widom
distribution raised a question that even this culminating
probabilistic program did not address.  The distribution of the
KPZ fixed point for the most basic initial condition is governed
by a closed nonlinear ODE with deep roots in classical
integrability.  Does this phenomenon extend?  Do the distributions
for other initial conditions satisfy analogous closed equations?
Do the distributions of prelimiting models, TASEP and the
stochastic six-vertex model and the log-gamma polymer at finite
time, satisfy analogous closed equations?  And are such equations,
when they exist, instances of the same classical integrable
structures that produced Painlev\'e~II, or are the appearances of
integrable equations in the KPZ setting unrelated coincidences?

Three recent results sharpened these questions considerably: the
Kadomtsev--Petviashvili (KP) equation was shown to govern the
distributions of the KPZ fixed point~\cite{QR2022}; coupled
matrix mKdV equations were established for the fixed point and
its periodic counterpart with narrow wedge initial
data~\cite{BPS2023}; and the non-abelian two-dimensional Toda
lattice was derived for polynuclear growth~\cite{MQR2024}.  These
results, and the questions they leave open, are surveyed in
\S\ref{sec:equations-state-of-art}.
\section{Classical integrability: from solitons to
Hirota--Miwa}\label{sec:classical-integrability}

The integrable equations that have surfaced in the study of KPZ
distributions, Painlev\'e~II, KP, the two-dimensional Toda
lattice, are not isolated objects.  They belong to a single
interconnected mathematical tradition, one whose central
achievement is the discovery that a vast collection of nonlinear
equations can be understood, and solved, through the spectral
theory of associated linear problems.  The story begins with a
solitary wave in a Scottish canal and reaches, for the purposes of
this work, a three-term discrete equation from which virtually
every known integrable hierarchy descends.

In August 1834, the Scottish engineer John Scott Russell was
observing boat experiments on the Union Canal near Hermiston,
Edinburgh, when a barge suddenly stopped and the water piled at
its prow rolled forward as a single rounded wave of elevation,
propagating along the channel without change of form.  Russell
followed the wave on horseback for one to two miles before losing
it in the windings of the canal, and reported his observation to
the British Association for the Advancement of Science, first in
1837 and then in his major 1844 \emph{Report on
Waves}~\cite{Russell1844}.  He called it the \emph{Wave of
Translation}, and argued that it was a genuine and stable
phenomenon in its own right.  The scientific establishment was
skeptical: George Airy and George Gabriel Stokes, foremost among
Russell's critics, regarded the wave as an artefact of linear
theory, and a theoretical explanation consistent with Russell's
observation was slow to emerge.  Partial progress came through the
work of Boussinesq~\cite{Boussinesq1871} in 1871 and
Rayleigh~\cite{Rayleigh1876} in 1876, but it was not until 1895,
more than sixty years after the original observation, that
Korteweg and de Vries~\cite{KdV1895} derived the equation
\begin{equation}\label{eq:KdV}
u_t + 6uu_x + u_{xxx} = 0
\end{equation}
as a long-wavelength asymptotic model for shallow water waves.
The Korteweg--de Vries equation admits an exact solitary-wave
solution, a sech$^2$ profile travelling at constant speed without
change of form, and this solution matched Russell's wave
quantitatively.  The long debate was resolved, after which the
equation, and the solitary wave along with it, lay essentially
dormant for another six decades.

A revival came in 1953, when Enrico Fermi, John Pasta, Stanis{\l}aw
Ulam, and Mary Tsingou, the latter as programmer of the new
MANIAC~I computer at Los Alamos, simulated a one-dimensional chain
of $32$ nonlinear oscillators, with the aim of observing the
gradual thermalization of energy initially concentrated in a
single Fourier mode~\cite{FPUT1955}.  The expectation, rooted in
the ergodic hypothesis of classical statistical mechanics, was
that nonlinear coupling would drive the system toward
equipartition.  Instead, after roughly $157$ oscillation periods,
the energy returned almost entirely to its initial mode.  The
recurrence was stunning, and flatly contradicted thermalization.

The FPUT recurrence was a puzzle for a decade.  Its resolution, in
1965, came from Norman Zabusky and Martin
Kruskal~\cite{ZabuskyKruskal1965}, who derived the KdV equation as
the continuum limit of the FPUT lattice and simulated it
numerically with periodic boundary conditions. They watched an initial cosine profile resolve into a train of solitary waves, each passing through the others without change of form: the system's energy was being carried by robust, particle-like solitary waves, and the FPUT recurrence was the consequence of their periodic reassembly.  To emphasize this
particle-like robustness, Zabusky and Kruskal coined the word
\emph{soliton}, combining the solitary-wave root with the suffix
\emph{-on} characteristic of elementary particles in physics.

The soliton's particle-like behaviour was striking, but even more striking was what came next. In 1967, Gardner, Greene, Kruskal, and Miura \cite{GGKM1967} announced something without precedent in the study of nonlinear partial differential equations: they solved the KdV equation exactly for a wide general class of initial data. The method, now called the \emph{inverse scattering transform}, interprets $u(x, t)$ as a potential in the time-independent Schr\"odinger equation and converts the nonlinear initial-value problem into a sequence of linear operations on the scattering data which then reconstructs a solution $u(x, t)$ via an inverse problem, the Gel'fand--Levitan--Marchenko equation~\cite{GelfandLevitan1955, Marchenko1955}.

Peter Lax~\cite{Lax1968} recognized that the inverse scattering
transform, remarkable as it was, was the consequence of a
structural property of the KdV equation that could be stated in
abstract terms and that turned out to generalize to many other nonlinear equations. Consider an overdetermined linear system
\begin{equation}\label{eq:Lax-system}
L\psi = \lambda\psi, \qquad \partial_t \psi = A\psi,
\end{equation}
in which $L = -\partial_x^2 + u$ is a linear differential operator depending on an
unknown function $u(x, t)$, and $A$ is a second linear
differential operator to be chosen.  The system is overdetermined
because it demands that $\psi$ satisfy two independent equations
simultaneously, and it is consistent only if these equations are
compatible.  Differentiating the spectral equation in $t$ and
using the flow equation to eliminate $\psi_t$ yields the
compatibility condition
\begin{equation}\label{eq:Lax-pair}
\frac{\partial L}{\partial t} = [A, L],
\end{equation}
the \emph{Lax equation}.  The Lax equation is a nonlinear equation
for $u$, and for the right choice of $L$ and $A$, it is the KdV equation. KdV thus arises as the compatibility condition of an
overdetermined linear problem.

The Lax pair elevated what had been a single remarkable
computation to a structural principle: a nonlinear evolution
equation is integrable when it arises as the compatibility
condition of an overdetermined linear system.  The existence of
such a system carries with it an infinite family of conserved
quantities, encoded in the spectral invariants of $L$, and the
possibility of exact solution through the inverse scattering
transform.  The recognition that KdV was not an isolated case
came rapidly.  Zakharov and
Shabat~\cite{ZakharovShabat1972} extended the method to the
nonlinear Schr\"odinger equation in 1972.  By 1974, Ablowitz,
Kaup, Newell, and Segur~\cite{AKNS1974} had
developed the AKNS scheme, a systematic framework encompassing the
nonlinear Schr\"odinger, sine-Gordon, and modified KdV equations
and establishing inverse scattering as, in their phrase, a
nonlinear Fourier transform.  A rapidly growing collection of
nonlinear equations turned out to be integrable, to possess Lax
pairs and soliton solutions and the full apparatus of inverse
scattering, and the integrable equations formed a rich,
interconnected hierarchy rather than a miscellaneous catalogue.

Running parallel to the inverse scattering program, and almost
entirely through the work of one person, a second approach to
integrability was emerging.  In 1971, Ryogo
Hirota~\cite{Hirota1971} published a short paper in
\emph{Physical Review Letters} that introduced the \emph{bilinear
method} for KdV.  Its central idea was a change of dependent
variable, $u = 2(\log \tau)_{xx}$,
that converts the KdV equation
into a quadratic equation in the new variable
$\tau$\footnote{In a series of papers in the early 1980s, the
Kyoto school, Date, Jimbo, Kashiwara, and
Miwa~\cite{DJKM1983}, showed that $\tau$ is the
$\tau$-\emph{function}, the central object of an algebraic
theory of integrable systems built on an
infinite-dimensional Grassmannian, and that the Hirota bilinear
identities are the Pl\"ucker relations defining the
Grassmannian embedding.  Segal and
Wilson~\cite{SegalWilson1985} arrived at the same picture
through an independent functional-analytic construction in
1985.} with all
derivatives appearing through a binary differential operator now
called the Hirota $D$-operator.  Unlike inverse scattering, the
bilinear method is purely algebraic and exceptionally efficient:
it produces $N$-soliton solutions as finite sums of exponentials,
without solving integral equations.  Hirota extended the method
rapidly over the next two years, deriving bilinear forms for
modified KdV, sine-Gordon, and the nonlinear Schr\"odinger
equation~\cite{Hirota1972a, Hirota1972b, Hirota1973}, and the
existence of a bilinear form became recognized as a signature of
integrability in its own right, alongside the Lax pair.

By the 1980s, Hirota had derived bilinear forms for more than a
dozen soliton equations, each through a separate application of
the $D$-operator, and had extended the method to discrete
analogues of KdV, the Toda lattice, and sine-Gordon.  The
accumulating catalogue strongly suggested a common algebraic
origin.  His 1981 paper~\cite{Hirota1981} identified it: a single
three-term bilinear relation on a multidimensional integer
lattice, now called the \emph{Hirota bilinear difference
equation}, or HBDE.\footnote{The equation is displayed in a form
due to Miwa~\cite{Miwa1982}, who proved that the lattice variables
of the HBDE can be identified with the higher times of the
Kadomtsev--Petviashvili hierarchy.  This identification shows that
the $\tau$-functions of the HBDE encode the entire KP hierarchy in
discrete coordinates, and conversely the $\tau$-functions of the
KP hierarchy satisfy the HBDE.  The equation is accordingly called
the \emph{Hirota--Miwa equation}.}
\begin{equation}\label{eq:HBDE}
\begin{aligned}
&\alpha\tau(l_1{+}1, l_2, l_3)\tau(l_1, l_2{+}1, l_3{+}1)
+ \beta\tau(l_1, l_2{+}1, l_3)\tau(l_1{+}1, l_2, l_3{+}1)\\
&\quad + \gamma\tau(l_1, l_2, l_3{+}1)\tau(l_1{+}1, l_2{+}1, l_3)
= 0,
\end{aligned}
\end{equation}
with $\alpha + \beta + \gamma = 0$.  Through a systematic
program of specializations and continuum limits, Hirota showed
that~\eqref{eq:HBDE} recovers the KdV equation, the KP equation,
modified KdV, the two-dimensional Toda lattice, sine-Gordon, and
their discrete analogues.  The entire catalogue of bilinear forms
that Hirota had built equation by equation through the 1970s was
contained in one discrete object.

The full significance of the Hirota--Miwa equation emerged over
the following two decades, as the same equation was independently
recognized in field after field.  In
algebraic geometry, the HBDE had already appeared, implicitly, in
Fay's 1973 \emph{Theta Functions on Riemann
Surfaces}~\cite{Fay1973}, as the trisecant identity for theta
functions on a Riemann surface: when the KP $\tau$-function is
expressed through Riemann theta functions, the Hirota bilinear
identities coincide with the trisecant relations of the Kummer
variety.  In discrete differential geometry, Doliwa and
Santini~\cite{DoliwaSantini1997} showed in 1997 that
multidimensional quadrilateral lattices, the discrete analogues
of classical conjugate nets, are governed by equations whose
$\tau$-function formulation is the HBDE, and the monograph of
Bobenko and Suris~\cite{BobenkoSuris2008} later established
multidimensional consistency as the defining principle of
discrete integrability, with the Hirota--Miwa equation as its
canonical instance.  In quantum integrable models, Krichever,
Lipan, Wiegmann, and Zabrodin~\cite{KLWZ1997} proved that the
fusion relations governing commuting quantum transfer matrices,
the functional equations satisfied by the eigenvalues of quantum
spin chains, are precisely the classical HBDE.  And in string
theory, Saito~\cite{Saito1987} showed as early as 1987 that
$N$-point bosonic string amplitudes satisfy the HBDE, with the
Miwa variables identified as Koba--Nielsen variables and external
momenta.

In 2006, Nimmo~\cite{Nimmo2006} introduced a \emph{non-abelian
Hirota--Miwa} system by dropping commutativity in the discrete Lax
pair: the dependent variables take values in a non-commutative
associative algebra, and the system is again the compatibility
condition of an overdetermined linear problem.  Nimmo constructed
explicit solutions using the quasideterminants of Gelfand and
Retakh~\cite{GelfandRetakh1991}, Gilson, Nimmo, and
Ohta~\cite{GilsonNimmoOhta2007} identified them as
quasi-Pl\"ucker coordinates extending the Grassmannian
interpretation, and Li, Nimmo, and
Tamizhmani~\cite{LiNimmoTamizhmani2009} studied continuum limits
to non-commutative KP.  Doliwa~\cite{Doliwa2010} gave a
complementary geometric interpretation through Desargues maps,
lattice maps into projective space whose consistency is equivalent
to Desargues' theorem.

Despite the long history of the Hirota--Miwa equation across
classical integrability, and despite the repeated appearance of
its reductions in KPZ distribution theory, no connection between
the Hirota--Miwa equation, in either its classical or its
non-abelian form, and the KPZ universality class had been made.
This work establishes that connection.
\section{Equations for KPZ distributions: state of the
art}\label{sec:equations-state-of-art}

The distinction between exact solvability and integrability in the sense of
\S\ref{sec:classical-integrability}\footnote{Integrability, in
the Lax sense, is a structural property of an equation: the
existence of a Lax pair, or equivalently an overdetermined linear
problem of which the equation is the compatibility condition.
Exact solvability is weaker: it asserts only that solutions can
be written in closed form, with no claim about the algebraic
machinery producing them.  The exactly solvable models of the KPZ
class described in
\S\ref{sec:KPZ-universality}--\S\ref{sec:early-breakthroughs}
were identified through diverse and apparently unrelated methods,
and whether they are integrable in the deeper structural sense is
among the questions this work addresses.} is sharp enough to
sort the KPZ literature into two unequal parts.
Every model in the exactly solvable cross-section of the class
possesses an exact formula for its distribution
functions, most commonly a Fredholm determinant.  Only a small
subset have been shown to carry a
closed nonlinear evolution equation for those distributions,
placing the model within reach of the classical integrable
hierarchies described in \S\ref{sec:classical-integrability}.  The
present section surveys what was known on this question at the
outset of this work.  The picture is one of scattered, deep,
and unconnected results: a short list of closed equations,
extracted by methods specific to individual models or specific
scaling regimes, with no mechanism explaining why the equations
should appear or predicting which integrable structure a new
model should carry.

The Painlev\'e~II equation~\eqref{eq:PainleveII} that
characterizes the Tracy--Widom distribution $F_2$
(\S\ref{sec:early-breakthroughs}) turned out to govern all
three canonical one-point limiting distributions of the class.
Tracy and Widom~\cite{TW1996} showed that the GOE and GSE
distributions $F_1$ and $F_4$ are expressible through the same
Hastings--McLeod solution $q$ and its integral, and Baik and
Rains~\cite{BaikRains2000, BaikRains2001a, BaikRains2001b},
in a series of papers on polynuclear growth with external
sources and on symmetrized random permutations, identified a
distribution $F_0$ governing KPZ models with stationary
initial data together with a one-parameter family $F_w$
interpolating between $F_2$, the superimposed GOE
distribution $F_1^2$, and $F_0$, all expressible
through the same Painlev\'e~II transcendent and thus all
governed by a single nonlinear ODE.

Each of these gap probabilities is the Fredholm
determinant of an operator of integrable type in the sense of
Its, Izergin, Korepin, and Slavnov~\cite{IIKS1990}: the
resolvent of such an operator satisfies a closed system of
total differential equations generalizing the
Jimbo--Miwa--M\^ori--Sato theory for correlation
functions~\cite{JMMS1980}, and the system reduces to the
Painlev\'e~II ODE for the Airy kernel on a single
semi-infinite interval.  The Painlev\'e representations
themselves were obtained by different routes, through the
operator-determinant method of Tracy and Widom~\cite{TW1994}
and the Riemann--Hilbert analysis of Baik and Rains.  The
structure is powerful, but it is specific to the asymptotic
limiting distributions $F_0$, $F_1$, $F_2$, not to the
finite-time distributions of any pre-limit model.

Closed nonlinear equations for exact, not merely asymptotic,
distribution functions first appeared in 2003, when
Borodin~\cite{Borodin2003} showed that the gap probability of
the Poissonized Plancherel measure, equivalently the Poissonized
distribution of the longest increasing subsequence, satisfies a
discrete Painlev\'e~II equation, and that the distribution of
the first row under the $z$-measures on partitions satisfies
discrete Painlev\'e~V.  The proof rests on a discrete analogue
of the Jimbo--Miwa--Ueno isomonodromy method: each gap
probability is expressed as a Fredholm determinant of a discrete
integrable operator, a discrete Riemann--Hilbert problem
produces a discrete Lax pair, and the discrete Painlev\'e
equation emerges as its compatibility
condition.\footnote{Borodin and
Boyarchenko~\cite{BorodinBoyarchenko2003} soon extended the
approach to classical discrete orthogonal polynomial ensembles,
treating in particular the Meixner ensemble, which corresponds
through the Robinson--Schensted--Knuth correspondence to
geometric last passage percolation, a central exactly solvable
KPZ model.}  Under a continuum limit, discrete Painlev\'e~II
degenerates to continuous Painlev\'e~II, and Borodin's pre-limit
result recovers the Tracy--Widom characterization of the
Baik--Deift--Johansson theorem from its discrete parent.

These were the first pre-limit KPZ distribution functions shown
to be integrable in the Lax sense of
\S\ref{sec:classical-integrability}.  Yet the discrete
isomonodromy method was tied closely to the rational structure
of the discrete log-derivative of each particular weight
function, and the equations it produced were model-specific.
Nothing in the approach suggested that discrete Painlev\'e
equations should appear more broadly across the exactly solvable
KPZ class, nor predicted which equation a given model ought to
satisfy.

The three results flagged at the end of
\S\ref{sec:early-breakthroughs}, all members of the
Hirota--Miwa hierarchy,\footnote{A separate line of
integro-differential equations for the KPZ equation itself,
beginning from the Amir--Corwin--Quastel
formula~\cite{ACQ2011} and developed by Cafasso, Claeys, and
Ruzza~\cite{CafassoClaeysRuzza2021} and
Krajenbrink~\cite{Krajenbrink2021}, has not been connected
to the Hirota--Miwa hierarchy and does not figure directly
in this work.} deserve closer attention.  Each is a
considerable individual advance, and each stands on its
own machinery.

In 2022, Quastel and Remenik~\cite{QR2022} proved that the
multipoint distributions of the KPZ fixed point, for general
initial data, are governed by the Kadomtsev--Petviashvili
equation.  Nothing in the prior literature had suggested that
a canonical nonlinear equation of mathematical physics would
govern the multipoint structure of the universal limit of
random interface growth; the Painlev\'e~II connections surveyed above had
reached only one-point limiting distributions for special
initial data.  For the one-point distribution
$F(t, x, r) = \PP(\mathfrak{h}(t, x) \leq r)$ of the KPZ
fixed point itself, however, the quantity
$\phi \defeq \partial_r^2 \log F$ satisfies the KP-II equation
\begin{equation}\label{eq:QR-KP}
\partial_t \phi + \tfrac{1}{2}\partial_r(\phi^2)
+ \tfrac{1}{12}\partial_r^3 \phi
+ \tfrac{1}{4}\partial_r^{-1}\partial_x^2 \phi = 0.
\end{equation}
When the distribution is independent of the spatial variable
$x$, as holds for flat initial data by translation invariance,
the equation reduces to KdV, and other self-similar reductions
of the KP-II equation recover the Painlev\'e~II equation
governing the Tracy--Widom distributions.  For $n$-point
equal-time multiposition distributions, an analogous matrix
KP-II equation governs $Q \defeq (I - K)^{-1}K(0,0)$, an
$n \times n$ matrix built from the resolvent;
for $n \geq 2$ the equation is intrinsically noncommutative.

The proof is purely algebraic.  It begins from the Fredholm
determinant formula of Matetski, Quastel, and
Remenik~\cite{MQR21} for the KPZ fixed point with general
initial data, and proceeds through an abstract
theorem: if a trace-class operator $K$ has a kernel satisfying
three differential identities relating derivatives in its
operator variables to derivatives in its physical variables,
then the resolvent of $K$ satisfies the matrix KP-II
equation.\footnote{The scalar case was established by
P\"oppe~\cite{Poppe1989}, who showed that Fredholm
determinants of kernels with appropriate linear evolution
properties yield solutions to the KP equation.}  For the KPZ fixed
point, the three kernel identities are linear
differential identities for the kernel, following from the
Airy differential equation.  From there, the passage
to matrix KP-II is substantial: it proceeds through an extended
sequence of resolvent manipulations and integration-by-parts
cancellations against the kernel, and the proof's depth rests as
much on that algebraic machinery as on the identities
themselves.

Powerful as the Quastel--Remenik theorem is, it leaves
significant territory uncovered.  It concerns the KPZ fixed
point itself, the universal scaling limit, and not the
distributions of any individual pre-limit model.  It addresses
equal-time, multiposition distributions, and does not reach
multitime distributions.  And the three kernel identities on
which the proof rests are specific to the analytic structure of
the KPZ fixed point kernel and to the Airy function ingredients
out of which it is built, with no obvious analogous identities
available for the kernels of pre-limit discrete models such as
TASEP or the stochastic six-vertex model.

In the periodic setting, Baik, Liu, and
Silva~\cite{BaikLiuSilva2022} showed that the one-point
distribution of the periodic KPZ fixed point, the scaling limit
of spatially periodic TASEP, satisfies coupled modified KdV
equations and coupled nonlinear heat equations, with the KP-II
equation emerging as the compatibility condition of the two
systems.  That the same families of integrable equations
should govern both the non-periodic and periodic cases, despite
the quite different analytical origins of the two Fredholm
determinant formulas, signalled that the integrable structure
did not depend on the particular analytic machinery used to
extract it.

Baik, Prokhorov, and Silva~\cite{BPS2023} extended this program
to the multi-point distributions of both the non-periodic and
periodic KPZ fixed points.  Their analysis is restricted to
narrow wedge initial data, in contrast with the general
initial data under which the Quastel--Remenik KP
equation holds; within that restriction, the multitime and
multiposition distributions of both fields, neither of which
the Quastel--Remenik framework reaches, satisfy coupled systems
of matrix modified KdV and matrix nonlinear Schr\"odinger
equations, with the matrix KP-II equation emerging as a
compatibility condition.  For equal-time distributions of
the non-periodic fixed point, the resulting matrix KP-II
coincides with that of Quastel and Remenik, though the two
derivations arrive at it through different objects.  The derivation proceeds
through the integrable operator formalism applied to a class of
kernels introduced under the name \emph{cubic integrable
operators}, into which both the non-periodic and periodic KPZ
fixed point kernels fall.

The third and most recent of the three results came from
Matetski, Quastel, and Remenik~\cite{MQR2024}, and is the first
in which an integrable lattice equation, rather than a continuum
PDE, governs the multipoint distribution functions of a KPZ
model at the pre-limit level: the multipoint distributions of
polynuclear growth are governed by the non-abelian
two-dimensional Toda lattice.  The equation was derived through
kernel identities and resolvent algebra closely analogous to
the Quastel--Remenik mechanism.  For $n$-point equal-time
distributions, the non-abelian 2D Toda equation is an
$n \times n$ matrix equation; for $n = 1$ it reduces to the
scalar 2D Toda lattice, and under specific initial conditions
the one-point distributions satisfy discrete Painlev\'e~II,
recovering earlier results of Baik and
Rains~\cite{BaikRains2001a, BaikRains2001c} and of
Borodin~\cite{Borodin2003}.  That the
distributions of a pre-limit KPZ model could be governed by a
canonical integrable lattice equation was a genuine surprise:
prior to this work, it had not been clear that such equations
should exist at the pre-limit level at all.  But the kernel identities driving the
non-abelian 2D Toda derivation are specific to the PNG
scattering kernel, and no obvious analogues are available for
other pre-limit models in the class.

Each of the results surveyed in this section is a major
advance in its own right, and each extracts a canonical member
of the Hirota--Miwa hierarchy from the distribution functions of
a specific KPZ model.  But no two of them share a common
mechanism.  Every closed nonlinear equation governing a KPZ
distribution has been derived by a method specific to a single
model, a single scaling regime, or a single analytic framework,
and nothing in the pattern suggests which framework a new model
should fit into.  For the central pre-limit discrete models in
the exactly solvable cross-section of the class, including
TASEP with general initial data, the stochastic six-vertex
model, and the log-gamma polymer at finite time, no closed
nonlinear equations for distribution functions are known at
all.  Whether the scattered appearances of Painlev\'e~II, KP, and the
non-abelian 2D Toda lattice in the KPZ setting are fragments of
a single algebraic framework, whether a common mechanism underlies their
existence and predicts which equation a new pre-limit model
should carry, and whether that mechanism brings the flagship
pre-limit discrete models within reach of the classical
integrability hierarchy traced in
\S\ref{sec:classical-integrability}, are the questions to which
the remainder of this work is addressed.
\section{Previous work: one-point equations and the
Hirota--Miwa connection}\label{sec:one-point-equations}
 
The open problem left at the end of
\S\ref{sec:equations-state-of-art}, the absence of closed
nonlinear equations for the finite-time distributions of
pre-limit discrete KPZ models, received its first answer in
2025, in the author's paper~\cite{Rodriguez2025}.\footnote{Among the models treated is the totally asymmetric
simple exclusion process (TASEP), a cornerstone of the
interacting particle systems literature since its introduction
by Spitzer~\cite{Spitzer1970} in 1970.}  For six
classical exactly solvable KPZ particle systems, the paper derived
closed nonlinear equations governing their finite-time
one-point distributions under arbitrary one-sided initial
data.  All six equations took Hirota bilinear form, a
classical signature of integrability in the sense of
\S\ref{sec:classical-integrability}.  For the three fully discrete-time models, explicit changes
of variables identified the equations as the Hirota--Miwa
equation itself.

The equations govern the one-point cumulative distributions
$F_{t, a, n} \defeq \PP(Y_n(t) > a)$, where $Y_n(t)$ is the
position of the $n$th particle of the model at time $t$, under
arbitrary one-sided initial data.  The
initial data enters only through the boundary condition
$F_{0, a, n} = \mathbf{1}_{y_n > a}$, and this generality was
itself new: no closed nonlinear equation had previously
governed the finite-time distributions of a fully discrete
pre-limit KPZ model under general initial data.  The paper also produced
Lax-pair and zero-curvature formulations for the bilinear
equations, in detail for reflected Brownian motion, TASEP, and
Parallel TASEP.  These are Lax pairs for the distribution
functions themselves, in the sense that the compatibility of
the associated linear system is equivalent to the bilinear
equation satisfied by $F_{t, a, n}$; in the Parallel TASEP
case, the Lax pair can alternatively be read off from the
classical zero-curvature formulation of the Hirota--Miwa
equation.

The technical engine behind all six equations was a general
solution theory operating on the known Fredholm determinant
formulas for the models.  Each of the models treated comes
equipped, from the prior exact-solvability literature, with a
representation
$F_{t, a, n} = \det(I - K_{t, a, n})_{L^2(X, \mu)}$
for a trace-class integral operator $K_{t, a, n}$ on a
separable Hilbert space, and the paper's solution theory takes
such a representation as input.  The theory identifies three
linear conditions on $K$ and on an associated pair of wave
functions $\psi_{t, a, n}, \phi_{t, a, n}$ that are jointly
sufficient for $F$ to satisfy a bilinear equation.  For
reflected Brownian motion, these are:
(i)~a rank-one decomposition
$\partial_a K_{t, a, n}(x, y)
= \psi_{t, a, n}(x)\phi_{t, a, n}(y)$;
(ii)~shift flows
$\nabla_n^+ \psi_{t, a, n} = \partial_a \psi_{t, a, n}$
and $\nabla_n^- \phi_{t, a, n}
= \partial_a \phi_{t, a, n}$; and
(iii)~dynamical flows
$\partial_t \psi_{t, a, n}
= \tfrac{1}{2}(\partial_a^2 + 2\partial_a + I)\psi_{t, a, n}$
and $\partial_t \phi_{t, a, n}
= -\tfrac{1}{2}(\partial_a^2 - 2\partial_a + I)\phi_{t, a, n}$.
Each condition is linear in the wave functions and in $K$
separately, and the bilinear equation for $F$ emerges from
their interaction through standard Fredholm determinant
identities.  The framework is in this sense a linearization:
the content is encoded in three linear conditions on an
integral kernel, and the nonlinear equation follows.
 
The framework provided a template for verification rather than
a mechanism for derivation.  The passage from the three
conditions to a closed bilinear equation had to be run, for
each model, through its own kernel decompositions and
resolvent manipulations, and no single calculation produced
more than one of the six equations.  The six proofs shared a
common outline: decompose $K$, extract $\psi$ and $\phi$,
identify the three conditions, run the resolvent algebra,
recover the bilinear equation.  What they did not share was a
common blueprint.

These results left the central questions of
\S\ref{sec:equations-state-of-art} only partly answered.  The
equations governed one-point distributions alone; the joint
law of several particles, where the richest content of the
integrable KPZ class lives, fell outside the reach of the
framework.  Moreover, given a new exactly solvable model with
a Fredholm determinant formula, the solution theory provided
no way to determine in advance if a bilinear equation, like
the HBDE, should govern the model's distribution function.
The framework developed in the remainder of this work
applies the central discovery of
\S\ref{sec:classical-integrability}, that nonlinear evolution
equations arise as compatibility conditions of overdetermined
linear systems, to produce every equation above, their
matrix-valued multipoint generalizations, and many more, from
a single master structure.
\section{Main results of this work}\label{sec:main-results}

We now state the main results precisely.  The diamond linear problem is built on a lattice
$\mathcal{V}$ generated by commuting invertible shifts
$T, S_1, S_2$, and on vector spaces $H$ and $E$ over a
field $\mathbb{F}$.  The shift $T$ plays a distinguished
role: in most of the KPZ models treated here it
corresponds to the spatial threshold parameter, the
variable $a$ in the cumulative distribution
$F_{t,a,n} = \PP(Y_n(t) > a)$, while $S_1$ and $S_2$
encode the particle index and the time
variable.\footnote{For Laplace-transform reductions,
such as the stochastic six-vertex model and the log-gamma
polymer, the role of $T$ is played by the Laplace
parameter, with analogous reinterpretations of
$S_1, S_2$.}  To each pair $(u, i)$ with
$u \in \mathcal{V}$ and $i \in \{1, 2\}$, assign
endomorphisms
$C_i(u), \Lambda_i(u) \in \End(E)$.  The \emph{diamond
linear problem} is the overdetermined system for
$\Psi \colon \mathcal{V} \to \Hom(E, H)$,
\begin{equation}\label{eq:intro-linear-problem}
\Psi(S_i u) = \Psi(Tu) C_i(u) - \Psi(u) \Lambda_i(u),
\qquad i \in \{1, 2\}.
\end{equation}

\definecolor{colT}{RGB}{0,119,136}
\definecolor{colS1}{RGB}{204,68,0}
\definecolor{colS2}{RGB}{102,45,145}

\begin{figure}[!htbp]
\centering

\begin{subfigure}[b]{0.35\textwidth}
\centering
\begin{tikzpicture}[>=Latex, scale=0.65, every node/.style={font=\scriptsize}]
  \node (u)  at (0,0)       {$(t,a,n)$};
  \node (T)  at (0,1.6)     {$(t,a{+}1,n)$};
  \node (S1) at (1.8,-1.2)  {$(t,a,n{+}1)$};
  \node (S2) at (-1.8,-1.2) {$(t{+}1,a,n)$};
  \draw[->, colT, very thick]
      (u) -- (T) node[midway, right] {$T$};
  \draw[->, colS1, thick]
      (u) -- (S1) node[midway, above right=-2pt] {$S_1$};
  \draw[->, colS2, thick]
      (u) -- (S2) node[midway, above left=-2pt] {$S_2$};
\end{tikzpicture}
\caption{Shift directions at a vertex.}
\label{fig:intro-triangle}
\end{subfigure}%
\hfill
\begin{subfigure}[b]{0.58\textwidth}
\centering
\begin{tikzpicture}[>=Latex, scale=0.8,
  every node/.style={font=\scriptsize},
  x={(1cm,0cm)},
  y={(0.5cm,0.38cm)},
  z={(0cm,1cm)}]

  \coordinate (u)    at (0,0,0);
  \coordinate (S1)   at (1.6,0,0);
  \coordinate (S2)   at (0,1.6,0);
  \coordinate (S12)  at (1.6,1.6,0);

  \coordinate (T)    at (0,0,1.4);
  \coordinate (TS1)  at (1.6,0,1.4);
  \coordinate (TS2)  at (0,1.6,1.4);
  \coordinate (TS12) at (1.6,1.6,1.4);

  \foreach \p in {u,S1,S2,S12,T,TS1,TS2,TS12}
    {\fill (\p) circle (0.04);}

  \node[below] at (u) {$u$};
  \node[below] at (S1) {$S_1 u$};
  \node[above left=-2pt] at (S2) {$S_2 u$};
  \node[below right=-2pt] at (S12) {$S_1 S_2 u$};
  \node[left] at (T) {$Tu$};
  \node[right] at (TS1) {$TS_1 u$};
  \node[above] at (TS2) {$TS_2 u$};
  \node[above] at (TS12) {$TS_1 S_2 u$};

  \draw[->, colS1, very thick]            (u)  -- (S1);
  \draw[->, colS2, very thick, dashed]    (u)  -- (S2);
  \draw[->, colS2, very thick]            (S1) -- (S12);
  \draw[->, colS1, very thick, dashed]    (S2) -- (S12);

  \draw[->, colS1, thick]            (T)   -- (TS1);
  \draw[->, colS2, thick, dashed]    (T)   -- (TS2);
  \draw[->, colS2, thick]            (TS1) -- (TS12);
  \draw[->, colS1, thick, dashed]    (TS2) -- (TS12);

  \draw[->, gray!70, thick] (u)    -- (T);
  \draw[->, gray!70, thick] (S1)   -- (TS1);
  \draw[->, gray!70, thick] (S2)   -- (TS2);
  \draw[->, gray!70, thick] (S12)  -- (TS12);
\end{tikzpicture}
\caption{Diamond face compatibility.}
\label{fig:intro-cube}
\end{subfigure}

\caption{The local lattice structure of the diamond linear
problem.  (A)~shows the three shift directions at a single
vertex; (B)~shows the diamond face, where compatibility
forces $\Psi(S_1 S_2 u)$ to agree along the two paths
from~$u$ (solid vs.\ dashed).}
\label{fig:intro-lattice}
\end{figure}
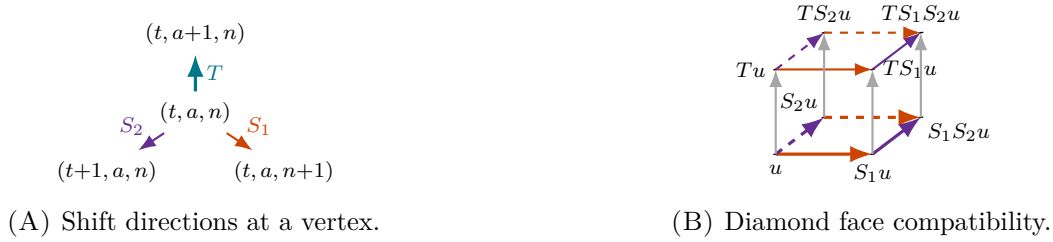

Compatibility of the linear
problem~\eqref{eq:intro-linear-problem} forces algebraic
constraints on the edge weights.

\begin{proposition}\label{prop:intro-CL-system}
The diamond linear problem~\eqref{eq:intro-linear-problem} is
compatible if and only if the edge weights $(C_i, \Lambda_i)$
satisfy, for all $u \in \mathcal{V}$ and all
$i, j \in \{1, 2\}$,
\begin{align}
C_i(Tu) C_j(S_i u) &= C_j(Tu) C_i(S_j u),
\label{eq:intro-C-diamond} \\
\Lambda_i(u) \Lambda_j(S_i u)
&= \Lambda_j(u) \Lambda_i(S_j u),
\label{eq:intro-L-diamond} \\
C_i(u) \Lambda_j(S_i u) + \Lambda_i(Tu) C_j(S_i u)
&= C_j(u) \Lambda_i(S_j u) + \Lambda_j(Tu) C_i(S_j u).
\label{eq:intro-mixed-diamond}
\end{align}
\end{proposition}

The proof, given in \S\ref{sec:linear-problem}, is a
direct computation: substitute the linear problem for
$\Psi(S_j S_i u)$ and $\Psi(S_i S_j u)$, subtract, and
equate coefficients.  We call
\eqref{eq:intro-C-diamond}--\eqref{eq:intro-mixed-diamond}
the \emph{diamond equations}.  The diamond equations
admit simple solutions from two sources.  For the
particle models, the bare edge weights are constant
scalars
$(c_1, c_2, \lambda_1, \lambda_2) \in \mathbb{F}^4$
determined by each model's jump probabilities and
auxiliary parameters.  Every product
in~\eqref{eq:intro-C-diamond}--\eqref{eq:intro-mixed-diamond}
commutes, so the diamond equations are satisfied
identically.\footnote{For inhomogeneous models such as
inhomogeneous discrete-time TASEP, the bare scalar weights
are $u$-dependent rather than constant, and the diamond
equations become genuine constraints on the functional
forms.  The framework handles both cases.}  For the vertex and polymer
models, a different source applies: Euclidean division
in a polynomial quotient algebra
$E = \mathbb{F}[w]/(\Pi)$ produces matrix-valued edge
weights satisfying the diamond equations identically.

\medskip

The Darboux theorem dresses these simple solutions into
the multipoint equations of
\S\ref{sec:catalogue}.  Two further ingredients
are needed: an adjoint linear problem and a
compatibility condition between $K$ and the wave
functions.

The adjoint linear problem is
\begin{equation}\label{eq:intro-adjoint-linear}
\Phi(Tu) = C_i(u) \Phi(S_i u)
- \Lambda_i(Tu) \Phi(TS_i u),
\qquad i \in \{1, 2\},
\end{equation}
for $\Phi \colon \mathcal{V} \to \Hom(H, E)$.  Its
compatibility conditions are again the diamond equations,
and $\Phi$ plays the role dual to $\Psi$ throughout the
construction.  The second ingredient is a kernel
$K \colon \mathcal{V} \to \End(H)$ satisfying
\emph{dressing compatibility}: the shifts of $K$ along
the lattice generators decompose through $\Psi$ and
$\Phi$,
\begin{equation}\label{eq:intro-dressing-compat}
K(Tu) - K(u) = \Psi(Tu) \Phi(Tu),
\quad
K(S_i u) - K(u) = \Psi(Tu) C_i(u) \Phi(S_i u),
\end{equation}
and the resolvent
$R(u) \defeq (I - zK(u))^{-1} \in \End(H)$ exists at
every $u \in \mathcal{V}$ for some fixed
$z \in \mathbb{F}$.  The dressed observable is the
endomorphism
\begin{equation}\label{eq:intro-M-def}
\mathcal{M}(u) \defeq I + z\Phi(u) R(u) \Psi(u)
\in \End(E),
\end{equation}
and the dressed edge weights are
\begin{equation}\label{eq:intro-dressed-weights}
\mathcal{M}_i(u) \defeq
\mathcal{M}(Tu)^{-1} C_i(u) \mathcal{M}(S_i u).
\end{equation}

\begin{theorem}[Darboux transformation]\label{thm:intro-Darboux}
Let $(C_i, \Lambda_i)$ be a solution of the diamond equations
\eqref{eq:intro-C-diamond}--\eqref{eq:intro-mixed-diamond},
let $\Psi, \Phi$ be solutions of the diamond linear problem
\eqref{eq:intro-linear-problem} and the adjoint linear
problem, and let $K$ be
dressing-compatible with $(\Psi, \Phi)$.  Then $\mathcal{M}(u)$
is invertible at every $u \in \mathcal{V}$, and the dressed
edge weights $\mathcal{M}_i$
of~\eqref{eq:intro-dressed-weights}, together with the
unchanged $\Lambda_i$, form a new solution
$(\mathcal{M}_i, \Lambda_i)$ of the diamond equations.
\end{theorem}

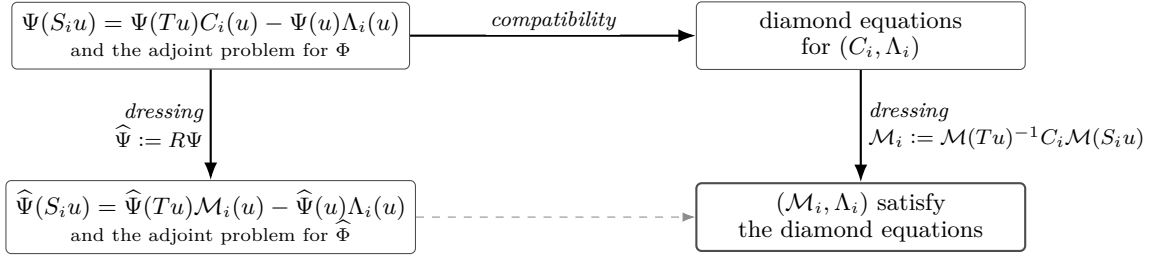
\begin{figure}[!ht]
\centering
\begin{tikzpicture}[>=Latex,
  box/.style={draw=black!70, rounded corners=2pt, align=center, inner sep=3.5pt, font=\footnotesize},
  albox/.style={draw=black!70, thick, rounded corners=2pt, align=center, inner sep=3.5pt, font=\footnotesize},
  arrlab/.style={font=\scriptsize\itshape, fill=white, inner sep=1.5pt}]
  \node[box] (LP) at (0,0)
    {$\Psi(S_i u) = \Psi(Tu)C_i(u) - \Psi(u)\Lambda_i(u)$\\[-1pt]
     {\tiny and the adjoint problem for $\Phi$}};
  \node[box, minimum width=4.35cm] (DE) at (8.6,0)
    {diamond equations\\[-1pt] for $(C_i, \Lambda_i)$};
  \draw[->, thick] (LP.east) -- (DE.west)
    node[arrlab, above, pos=0.5] {compatibility};
  \node[box] (DLP) at (0,-2.4)
    {$\widehat\Psi(S_i u) = \widehat\Psi(Tu)\mathcal{M}_i(u) - \widehat\Psi(u)\Lambda_i(u)$\\[-1pt]
     {\tiny and the adjoint problem for $\widehat\Phi$}};
  \draw[->, thick] (LP.south) -- (DLP.north)
    node[arrlab, left=1pt, pos=0.5, align=right] {dressing\\ $\widehat\Psi \defeq R\Psi$};
  \node[albox, minimum width=4.35cm] (DDE) at (8.6,-2.4)
    {$(\mathcal{M}_i, \Lambda_i)$ satisfy\\[-1pt] the diamond equations};
  \draw[->, dashed, black!45] (DLP.east) -- (DDE.west);
  \draw[->, thick] (DE.south) -- (DDE.north)
    node[arrlab, right=1pt, pos=0.5, align=left]
    {dressing \\ $\mathcal{M}_i \defeq \mathcal{M}(Tu)^{-1}C_i\mathcal{M}(S_i u)$};
\end{tikzpicture}
\caption{The Darboux transformation.  The left column dresses
the wave functions and the right column dresses the edge
weights; the dressed pair $(\mathcal{M}_i, \Lambda_i)$ again
satisfies the diamond equations.}
\label{fig:intro-darboux}
\end{figure}

The proof, given as Theorem~\ref{thm:Darboux} of
\S\ref{sec:Darboux-transformations}, uses only dressing
compatibility and the resolvent identity
$R(v) - R(w) = zR(v)(K(v) - K(w))R(w)$.  No Fredholm
determinant theory enters: the theorem assumes only
resolvent existence, not trace-class structure.

When $\dim E = m$, the dressed observable
$\mathcal{M}(u)$ is an $m \times m$ matrix whose
entries encode the joint distribution of $m$ observation
points.  The mixed dressed diamond
equation~\eqref{eq:intro-mixed-diamond} for
$(\mathcal{M}_i, \Lambda_i)$ is the closed multipoint
equation.  As a representative instance, the equation
governing the multipoint distribution of discrete-time
TASEP with right Bernoulli jumps reads
\begin{equation}\label{eq:intro-RBJ-multipoint}
\bigl(q\mathcal{M}_{t, \mathbf{a}+\mathbf{1},
\mathbf{n}+\mathbf{1}}^{-1}
+ p\mathcal{M}_{t+1, \mathbf{a}+\mathbf{2},
\mathbf{n}}^{-1}\bigr)
\mathcal{M}_{t+1, \mathbf{a}+\mathbf{1},
\mathbf{n}+\mathbf{1}}
- \mathcal{M}_{t, \mathbf{a}+\mathbf{1},
\mathbf{n}}^{-1}
\bigl(q\mathcal{M}_{t+1, \mathbf{a}+\mathbf{1},
\mathbf{n}}
+ p\mathcal{M}_{t, \mathbf{a},
\mathbf{n}+\mathbf{1}}\bigr) = 0,
\end{equation}
where $\mathcal{M}_{t, \mathbf{a}, \mathbf{n}} \in
\End(\R^m)$ encodes the joint distribution of
$m$ particles with labels $\mathbf{n}$ at thresholds
$\mathbf{a}$ at time $t$, and $p, q = 1 - p$ are the
right-jump and stay probabilities.  For $m \geq 2$ the
equation is intrinsically noncommutative.

\medskip

The scalar reduction of the Darboux theorem, the case
$\dim E = 1$, delivers the connection to the one-point
bilinear equations of~\S\ref{sec:one-point-equations}.
The bridge from the dressed observable to the Fredholm
determinant is Sylvester's determinant identity, which
gives
\begin{equation}\label{eq:intro-Woodbury}
F(u) = F(Tu) \det\nolimits_E\bigl(\mathcal{M}(Tu)\bigr),
\quad
F(u) \defeq \det\nolimits_H\bigl(I - zK(u)\bigr),
\end{equation}
whenever $K$ is dressing-compatible and $\dim E$ is
finite.  When $\dim E = 1$, the determinant collapses to
the scalar $\mathcal{M}(Tu)$, and the dressed weights
become explicit ratios of $F$ at neighbouring sites.
The mixed diamond reduces to a variable-coefficient
Hirota--Miwa equation.

\begin{corollary}[Scalar Hirota--Miwa reduction]
\label{cor:intro-scalar-HM}
Under the hypotheses of Theorem~\ref{thm:intro-Darboux}
with $\dim E = 1$, write
$\alpha_{ij}(u) \defeq c_i(Tu)\lambda_j(Tu)$ and
assume $F(T^\ell v) \to 1$ as $\ell \to -\infty$
together with the nondegeneracy
$\alpha_{12}(u) \neq \alpha_{21}(u)$.  Then the mixed
diamond for the dressed pair
$(\mathcal{M}_i, \lambda_i)$ is equivalent to
\begin{equation}\label{eq:intro-vcHM}
\alpha_{12}(u) F(S_1 u) F(T S_2 u)
- \alpha_{21}(u) F(S_2 u) F(T S_1 u)
- \bigl(\alpha_{12}(u) - \alpha_{21}(u)\bigr)
F(Tu) F(S_1 S_2 u)
= 0.
\end{equation}
\end{corollary}

The six one-point bilinear equations derived
in~\cite{Rodriguez2025} are scalar specializations
of~\eqref{eq:intro-vcHM}, recovered by substituting
each model's scalar weights for $\alpha_{12}$ and
$\alpha_{21}$.  Their shared Hirota form, observed
in~\cite{Rodriguez2025} through six parallel
computations and not explained by them, is forced: any
KPZ model accommodated by the framework produces a
Hirota-form bilinear equation by the same mechanism.

For models whose Fredholm determinant data is presently
available only in one-point form, the framework
nonetheless produces new equations governing their
observables.  For the stochastic six-vertex model with
parameters $0 < b_2 < b_1 < 1$ and $\tau = b_2/b_1$,
the observable is the $\tau$-Laplace transform of the
height function,
\[
\mathcal{L}_{t,x}(\zeta)
\defeq
\E_{\mathrm{step}}
\left[
\frac{1}{(\zeta\tau^{N_x(t)};\tau)_\infty}
\right],
\qquad
(a;\tau)_\infty
\defeq \prod_{k=0}^\infty (1 - a\tau^k),
\]
where the integer $m$, which discretizes the spectral
variable as $\zeta_m = \tau^m\zeta$, plays the role of
the distinguished T-shift.
Corollary~\ref{cor:intro-scalar-HM} specializes to the
bilinear equation
\begin{equation}\label{eq:intro-S6V-scalar-HM}
(1 - b_1)
\mathcal{L}_{t,x+1,m+1}
\mathcal{L}_{t-1,x,m}
- (1 - b_2)
\mathcal{L}_{t-1,x,m+1}
\mathcal{L}_{t,x+1,m}
+ (b_1 - b_2)
\mathcal{L}_{t,x,m+1}
\mathcal{L}_{t-1,x+1,m}
= 0,
\end{equation}
an instance of~\eqref{eq:intro-vcHM} with constant
Hirota coefficients built from the vertex weights
$b_1, b_2$.

\medskip

The Darboux theorem produces multipoint equations from
any dressing-compatible triple $(\Psi, \Phi, K)$, but
constructing such a triple for multipoint distributions
is substantial work.  The product graph construction
provides a systematic mechanism.

\begin{theorem}[Product graph construction]
\label{thm:intro-product-graph}
Let $(c_k, \lambda_k)$ be a scalar solution of the
diamond equations on a base graph $\mathcal{G}$, let
$(\psi, \phi)$ be solutions of the corresponding
linear and adjoint linear problems, and let $\{B_u\}$
be \emph{admissible propagators} on the product graph
$\mathcal{G}^m$: a family of strictly
lower-triangular endomorphism-valued functions
satisfying covariance, splitting, and compatibility
conditions with the diamond linear problem
\textup{(Definition~\ref{def:admissible-prop})}.  Then
the theorem constructs explicit wave functions
$\Psi, \Phi$ and a kernel $K$ on $\mathcal{G}^m$ that
are dressing-compatible, so that
Theorem~\ref{thm:intro-Darboux} applies and the
resulting dressed observable $\mathcal{M}$ satisfies
the diamond equations on $\mathcal{G}^m$.
\end{theorem}

The theorem is established as
Lemmas~\ref{lem:Lambda-diamond}, \ref{lem:Psi-Phi-linear},
\ref{lem:K-factorization} and
Theorem~\ref{thm:product-graph} of
\S\ref{sec:product-graph}.  It converts the diamond
framework from verification to construction: from a
solution of the diamond equations and a choice of
admissible propagator, the theorem manufactures the full
algebraic apparatus on the product graph.  The
construction iterates.  The concept of admissible
propagators is, to the author's knowledge, new.  In
every discrete KPZ model treated in this work, the
admissible propagators are constructed explicitly from
the transition operators of the underlying stochastic
process.

For each model whose Fredholm data admits a product
graph interpretation, verification of the multipoint
equation reduces to a uniform five-step procedure:
(1)~identify the lattice and shifts from the model's
parameter space; (2)~read off the bare scalar edge
weights and the seed functions $(\psi, \phi)$ from the
Fredholm determinant data; (3)~extract the admissible
propagator $B$ from the block-triangular structure of
the extended kernel; (4)~verify that $B$ satisfies the
three admissibility conditions; and (5)~apply the
product graph theorem and read off the multipoint
equation from the mixed dressed diamond.  The
model-specific content lies in steps (2)--(4); the
passage from admissible data to a closed multipoint
equation is the same theorem in every case.

The product graph construction applies to particle
models whose propagators have a natural
transition-operator interpretation.  For vertex and
polymer models, a different mechanism provides the seed
data: Euclidean division in a polynomial quotient
algebra.  Fix a polynomial
$\Pi \in \mathbb{F}[w]$ of degree $N$ and set
$E = \mathbb{F}[w]/(\Pi)$.  For any polynomial $p$
with $\deg p \leq N$, define endomorphisms
$C_p, \Lambda_p \colon E \to E$ by the division
identity
\begin{equation}\label{eq:intro-PQ-division}
p(w)f(w) = \Pi(w)(C_p f)(w) - (\Lambda_p f)(w),
\qquad f \in E.
\end{equation}
That is, $C_p f$ is the quotient and $-\Lambda_p f$
the remainder of Euclidean division of $pf$ by $\Pi$.
Commutativity of polynomial multiplication forces
these operators to satisfy the diamond
equations~\eqref{eq:intro-C-diamond}--\eqref{eq:intro-mixed-diamond}
for any pair $p, q$
(Lemma~\ref{thm:PQ-bare-diamond}): the product $pqf$
computed in two orders gives the same polynomial, and
uniqueness of the degree-graded decomposition forces
each diamond equation separately.  The Euclidean
division operators serve as the bare edge weights to
which the Darboux theorem is applied.

\medskip

The discrete framework covers ten of the eighteen
models.  The remaining eight require frameworks in which one or
more of the discrete shifts is replaced by a continuous
derivative.  This work develops three
such frameworks, each establishing diamond equations, a
Darboux theorem, and a scalar reduction in a
progressively more continuous setting.

In the \emph{semi-discrete framework}
(\S\ref{sec:sc-framework}), one of the discrete shifts
is replaced by a continuous derivative while the
remaining shifts stay discrete.  The linear problem
becomes
\begin{align}
\partial_1\Psi(u)
&= \Psi(Tu)C_1(u) - \Psi(u)\Lambda_1(u),
\label{eq:intro-sc-linear-1} \\
\Psi(S_j u)
&= \Psi(Tu)C_j(u) - \Psi(u)\Lambda_j(u),
\qquad j \geq 2.
\label{eq:intro-sc-linear-j}
\end{align}
Compatibility forces the \emph{semi-discrete diamond
equations} (Proposition~\ref{prop:sc-diamond}), which
augment the discrete
system~\eqref{eq:intro-C-diamond}--\eqref{eq:intro-mixed-diamond}
with derivative corrections.  The semi-discrete mixed
diamond, for instance, reads
\begin{equation}\label{eq:intro-sc-mixed}
C_1(u)\Lambda_j(u) + \Lambda_1(Tu)C_j(u)
- C_j(u)\Lambda_1(S_j u) - \Lambda_j(Tu)C_1(S_j u)
- \partial_1 C_j(u) = 0.
\end{equation}
The semi-discrete Darboux theorem
(Theorem~\ref{thm:sc-Darboux}) has the same character as
the discrete case: the $C$-weights are dressed by the
observable $\mathcal{M}$ while the $\Lambda$-weights
pass through unchanged.  The framework admits two
natural presentations, related by a change of lattice
basis (\S\ref{sec:sc-dual}): in the standard
presentation, the continuous derivative falls on
$\Psi(u)$; in the dual presentation, it falls on
$\Psi(Tu)$.  Among the models treated in this work,
Push-TASEP sits naturally in the standard
presentation, while continuous-time TASEP sits in the
dual.

The scalar reduction in the semi-discrete regime
produces a bilinear differential-difference equation
(Proposition~\ref{prop:sc-scalar-HM}):
\begin{equation}\label{eq:intro-sc-scalar}
F(Tu)\partial_1 F(S_j u)
- F(S_j u)\partial_1 F(Tu)
+ \alpha_j(u)
\bigl[F(Tu)F(S_j u) - F(u)F(TS_j u)\bigr] = 0,
\end{equation}
where $\alpha_j(u) = c_1(Tu)\lambda_j(Tu)/c_j(Tu)$.
As a representative instance, the multipoint equation
for continuous-time TASEP in the dual presentation reads
\begin{equation}\label{eq:intro-CT-TASEP}
\mathcal{M}_{t,\mathbf{a}+\mathbf{1},
\mathbf{n}+\mathbf{1}}^{-1}
\partial_t \mathcal{M}_{t,\mathbf{a}+\mathbf{1},
\mathbf{n}+\mathbf{1}}
- \mathcal{M}_{t,\mathbf{a}+\mathbf{1},
\mathbf{n}}^{-1}
\partial_t \mathcal{M}_{t,\mathbf{a}+\mathbf{1},
\mathbf{n}}
+ \mathcal{M}_{t,\mathbf{a}+\mathbf{2},
\mathbf{n}}^{-1}
\mathcal{M}_{t,\mathbf{a}+\mathbf{1},
\mathbf{n}+\mathbf{1}}
- \mathcal{M}_{t,\mathbf{a}+\mathbf{1},
\mathbf{n}}^{-1}
\mathcal{M}_{t,\mathbf{a},\mathbf{n}+\mathbf{1}} = 0,
\end{equation}
an $m \times m$ matrix differential-difference equation
coupling the dressed observable at four lattice points
with the time derivative $\partial_t$.

\medskip

In the \emph{parabolic framework}
(\S\ref{sec:dc-framework-clean}), a second continuous
limit is taken: two of the three lattice directions
become continuous derivatives, and the lattice reduces to
$I_1 \times I_2 \times \mathcal{V}'$ with two
continuous variables and one surviving discrete shift
$S_3$.  The linear problem is parabolic in the
continuous directions,
\begin{align}
\partial_1 \Psi(u)
&= \tfrac{1}{2}\partial_2^2 \Psi(u)
  + \partial_2\Psi(u)C_1(u) - \Psi(u)\Lambda_1(u),
\label{eq:intro-dc-linear-1} \\
\Psi(S_3 u)
&= \partial_2\Psi(u)C_2(u) - \Psi(u)\Lambda_2(u),
\label{eq:intro-dc-linear-S3}
\end{align}
coupling a heat operator to a first-order recurrence.
The compatibility conditions are the \emph{parabolic
diamond equations}
(Proposition~\ref{prop:dc-clean-diamond}), which involve
the heat operator
$\mathcal{H} \defeq \partial_1 - \frac{1}{2}\partial_2^2$
and the shifted bracket
$[X, Y]_{S_3}(u) \defeq X(u)Y(u) - Y(u)X(S_3 u)$.
The parabolic mixed diamond reads
\begin{equation}\label{eq:intro-dc-mixed}
\mathcal{H} C_2(u)
+ \partial_2 C_1(u)C_2(u)
- \partial_2 C_2(u)C_1(S_3 u)
+ \partial_2\Lambda_2(u)
- [C_1, \Lambda_2]_{S_3}(u)
- [\Lambda_1, C_2]_{S_3}(u) = 0.
\end{equation}

The parabolic Darboux theorem
(Theorem~\ref{thm:dc-clean-general-darboux}) reverses
the discrete pattern: the $C$-weights are fixed by the
dressing and the $\Lambda$-weights absorb the
corrections, expressed through the additive dressed
observable
$\mathcal{A}(u) \defeq z\Phi(u)R(u)\Psi(u)$:
\begin{align}
\widehat\Lambda_1(u)
&= \Lambda_1(u)
  + \partial_2\mathcal{A}(u)
  + [\mathcal{A}(u), C_1(u)],
\label{eq:intro-dc-L1-hat} \\
\widehat\Lambda_2(u)
&= \Lambda_2(u)
  + [\mathcal{A}, C_2]_{S_3}(u).
\label{eq:intro-dc-L2-hat}
\end{align}
By the Darboux theorem, the dressed pair
$(C_i, \widehat\Lambda_i)$ again satisfies the parabolic
diamond equations, and the mixed diamond for this pair is
a closed equation for $\mathcal{A}$.  Writing
$\nabla\mathcal{A} \defeq
\mathcal{A}_{\mathbf{n}+\mathbf{1}}
- \mathcal{A}_{\mathbf{n}}$, the multipoint equation
for reflected Brownian motion, with constant scalar
edge weights $C_1 = C_2 = I$ and
$\Lambda_1 = -\frac{1}{2}I$, $\Lambda_2 = -I$, reads
\begin{equation}\label{eq:intro-RBM}
\partial_t \nabla\mathcal{A}
- \tfrac{1}{2}\partial_{\mathbf{a}}^2
(\mathcal{A}_{\mathbf{n}+\mathbf{1}}
+ \mathcal{A}_{\mathbf{n}})
+ \tfrac{1}{2}\partial_{\mathbf{a}}
(\nabla\mathcal{A})^2
+ \tfrac{1}{2}[\nabla\mathcal{A},
\partial_{\mathbf{a}}(\mathcal{A}_{\mathbf{n}+\mathbf{1}}
+ \mathcal{A}_{\mathbf{n}})] = 0,
\end{equation}
an $m \times m$ matrix PDE, where
$\partial_{\mathbf{a}} \defeq \sum_{i=1}^m \partial_{a_i}$
denotes the collective threshold derivative.

The scalar reduction in the parabolic regime produces a
bilinear equation with a \emph{trace defect}: for the
multipoint Fredholm determinant
$F_{\mathbf{n}} = \det_H(I - zK_{\mathbf{n}})$,
\begin{equation}\label{eq:intro-RBM-Hirota}
\bigl[D_t - \tfrac{1}{2}D_{\mathbf{a}}^2\bigr]
F_{\mathbf{n}+\mathbf{1}} \cdot F_{\mathbf{n}}
= \tfrac{1}{2} F_{\mathbf{n}+\mathbf{1}}
F_{\mathbf{n}}
\bigl(\operatorname{tr}_E
((\nabla\mathcal{A})^2)
- (\operatorname{tr}_E
\nabla\mathcal{A})^2\bigr),
\end{equation}
where $D_t, D_{\mathbf{a}}$ denote Hirota derivatives.
At $m = 1$ the trace defect vanishes and
\eqref{eq:intro-RBM-Hirota} reduces to the scalar
Hirota equation
$[D_t - \frac{1}{2}D_a^2]
F_{n+1} \cdot F_n = 0$,
recovering the one-point bilinear equation
of~\cite{Rodriguez2025} for reflected Brownian motion.

\medskip

In the \emph{continuum framework}
(\S\ref{sec:kpz-reduced-framework}), the passage is
completed: every direction becomes continuous, leaving
three real variables $(a, x, t)$ and no discrete shifts.
The linear problem is
\begin{align}
\partial_x\Psi &= \partial_a^2\Psi - 2\Psi C,
\label{eq:intro-kpz-Lx} \\
\partial_t\Psi &= -\tfrac{1}{3}\partial_a^3\Psi
  + \partial_a\Psi C - \Psi\Lambda,
\label{eq:intro-kpz-Lt}
\end{align}
and its compatibility forces two PDE conditions on
$(C, \Lambda)$ (Proposition~\ref{prop:kpz-red-compat}):
\begin{gather}
2\partial_a\Lambda + \partial_x C + \partial_a^2 C = 0,
\label{eq:intro-kpz-constraint} \\
2\partial_t C - \partial_x\Lambda + \partial_a^2\Lambda
+ \tfrac{2}{3}\partial_a^3 C
- 2\partial_a C C + 2[C, \Lambda] = 0.
\label{eq:intro-kpz-evolution}
\end{gather}
The continuum Darboux theorem
(Theorem~\ref{thm:kpz-red-Darboux-compat}) dresses both
$C$ and $\Lambda$, in contrast to the parabolic regime
where only the $\Lambda$-weights are corrected:
\begin{equation}\label{eq:intro-kpz-dressed}
\widehat C = C - \partial_a\mathcal{A},
\qquad
\widehat\Lambda = \Lambda
+ \tfrac{1}{2}(\partial_x\mathcal{A}
+ \partial_a^2\mathcal{A}),
\end{equation}
where $\mathcal{A} = z\Phi R\Psi$ is the additive
dressed observable.  By the Darboux theorem, the dressed
pair $(\widehat C, \widehat\Lambda)$ again satisfies
the compatibility
system~\eqref{eq:intro-kpz-constraint}--\eqref{eq:intro-kpz-evolution}.
In particular, $\mathcal{A}$ satisfies
\begin{align}
0
&= \partial_{t,a}\mathcal{A}
+ \tfrac{1}{4}\partial_x^2\mathcal{A}
+ \tfrac{1}{12}\partial_a^4\mathcal{A}
+ \tfrac{1}{2}\partial_a(\partial_a\mathcal{A})^2
+ \tfrac{1}{2}[\partial_a\mathcal{A},
\partial_x\mathcal{A}]
\notag \\
&\quad
- \partial_a(C\partial_a\mathcal{A})
- \tfrac{1}{2}[C, \partial_x\mathcal{A}]
+ \tfrac{1}{2}[C, \partial_a^2\mathcal{A}]
+ [\partial_a\mathcal{A}, \Lambda].
\label{eq:intro-pKP}
\end{align}
The first line is the matrix potential KP equation; the
second line collects the background corrections, which
vanish when $C = \Lambda = 0$, the setting of the KPZ
fixed point (\S\ref{sec:FP}).  For $m$ observation
points, $\mathcal{A}$ is an $m \times m$ matrix; when
$m = 1$, the commutators vanish and the equation
reduces to scalar potential KP.

The Fredholm determinant
$F = \det_H(I - zK)$ satisfies KP in bilinear form
with a trace defect
(Corollary~\ref{cor:FP-trace-defect}):
\begin{equation}\label{eq:intro-trace-defect}
\Bigl[D_{t,\mathbf{a}}
+ \tfrac{1}{4}D_{\mathbf{x}}^2
+ \tfrac{1}{12}D_{\mathbf{a}}^4\Bigr] F \cdot F
= F^2\bigl((\operatorname{tr}_{\R^m}
\partial_{\mathbf{a}}\mathcal{A})^2
- \operatorname{tr}_{\R^m}
((\partial_{\mathbf{a}}\mathcal{A})^2)\bigr),
\end{equation}
where $D_{\mathbf{a}}, D_{\mathbf{x}}, D_t$ denote
Hirota derivatives.  The
left-hand side is KP-II in Hirota bilinear form.  At
$m = 1$ the trace defect vanishes, and
\eqref{eq:intro-trace-defect} reduces to
\begin{equation}\label{eq:intro-FP-scalar}
\Bigl[D_a D_t + \tfrac{1}{4}D_x^2
+ \tfrac{1}{12}D_a^4\Bigr] F \cdot F = 0.
\end{equation}

\medskip

Applied to eighteen exactly solvable KPZ models, the
framework produces the equations catalogued in
\S\ref{sec:catalogue}.

\begin{theorem}[Matrix equations for integrable KPZ
models]\label{thm:intro-catalog}
For each model listed in
\S\ref{sec:catalogue}, the Fredholm determinant
data of the model satisfies the hypotheses of the
Darboux theorem appropriate to its
regime\textup{:} the discrete Darboux
\textup{(Theorem~\ref{thm:intro-Darboux})} for the
models of Chapters~6 and~7\textup{,} the semi-discrete
Darboux \textup{(Theorem~\ref{thm:sc-Darboux})} for
the models of Chapter~8\textup{,} the parabolic Darboux
\textup{(Theorem~\ref{thm:dc-clean-general-darboux})}
for the models of Chapter~9\textup{,} and the continuum
Darboux
\textup{(Theorem~\ref{thm:kpz-red-Darboux-compat})}
for the KPZ fixed point.  In each case the dressed
observable $\mathcal{M}$ or its additive counterpart
$\mathcal{A}$ is governed by a closed matrix equation
given by the mixed dressed diamond of the corresponding
regime.  Every equation is new except for the KPZ fixed point, where the matrix KP equation was established by Quastel and Remenik~\cite{QR2022}; the framework provides an independent derivation.
\end{theorem}

Each verification occupies a dedicated chapter section
and consists of identifying the lattice shifts and edge
weights from the model's Fredholm determinant data,
verifying that the seed functions, propagators, and
kernel satisfy the dressing compatibility conditions,
and applying the Darboux theorem.

\medskip

\begin{figure}[!b]
\centering
\includegraphics[width=0.55\textwidth]{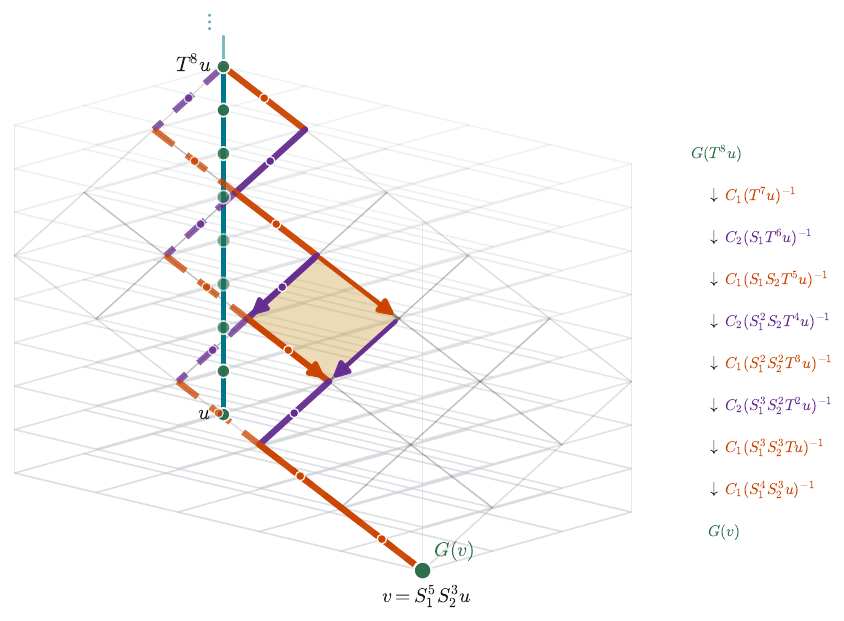}
\caption{The gauge function $G(v)$, built along two lattice paths by
the recursion $G(S_i w) = C_i(w)^{-1}G(Tw)$.  The paths agree by the
$C$-diamond equation.}
\label{fig:intro-gauge}
\end{figure}

When the $C$-weights are invertible, the diamond equations
are gauge equivalent to the non-abelian Hirota--Miwa
system of Nimmo~\cite{Nimmo2006}
(Proposition~\ref{prop:diamond-NAHM}): the
non-commutative master equation from which the KdV,
KP, Toda, and sine-Gordon hierarchies descend through
specialization.  The gauge connecting the two is
nonlocal (Figure~\ref{fig:intro-gauge}), built from bare
$C$-data along entire lattice paths, but the underlying
system is the same.  The
diamond parametrization is simply the form in which
integrable KPZ models present themselves.  The
invertibility hypothesis is genuine: the Euclidean-division
constructions of \S\ref{sec:PQ} produce singular
$C$-weights, and for those models the diamond equations and
all Darboux conclusions hold as stated, while the passage
to the Hirota--Miwa gauge is not asserted.

The equivalence carries content in both directions.
Inward: the separate appearances of Painlev\'e~II, KP,
and the two-dimensional Toda lattice in KPZ
distribution theory are not independent integrable
structures, but reductions of one system acting in
different regimes.  Outward: the Fredholm determinants
of stochastic kernels furnish a new class of solutions
to the non-abelian Hirota--Miwa equation, solutions
constructed from random growth rather than algebraic
geometry or representation theory.  The existence of
these solutions, and the constraints their stochastic
origin places on the spectral data, remain largely
unexplored.

Section~\ref{sec:classical-integrability} traced a line
from a solitary wave in a Scottish canal to a single
discrete equation governing classical soliton theory,
and asked whether that equation also governs the KPZ
universality class.  Across eighteen models and four
regimes, it does.

\subsection{Catalogue of equations}\label{sec:catalogue}

\begin{center}
\textsc{I. Discrete Particle Models}
\end{center}
\medskip

\noindent\textbf{1. Right Bernoulli Jumps.}
\hfill Theorem~\ref{thm:RBJ-H5}

\begin{equation*}
\bigl(q\mathcal{M}_{t, \mathbf{a}+\mathbf{1},
\mathbf{n}+\mathbf{1}}^{-1}
+ p\mathcal{M}_{t+1, \mathbf{a}+\mathbf{2},
\mathbf{n}}^{-1}\bigr)
\mathcal{M}_{t+1, \mathbf{a}+\mathbf{1},
\mathbf{n}+\mathbf{1}}
- \mathcal{M}_{t, \mathbf{a}+\mathbf{1},
\mathbf{n}}^{-1}
\bigl(q\mathcal{M}_{t+1, \mathbf{a}+\mathbf{1},
\mathbf{n}}
+ p\mathcal{M}_{t, \mathbf{a},
\mathbf{n}+\mathbf{1}}\bigr) = 0.
\end{equation*}

\medskip

\noindent\textbf{2. Inhomogeneous Bernoulli Jumps.}
\hfill Theorem~\ref{thm:IBJ-multipoint}

\begin{equation*}
\begin{aligned}
&\bigl(\Lambda_{\mathbf{a},
\mathbf{n}}
\mathcal{M}_{t, \mathbf{a}+\mathbf{1},
\mathbf{n}+\mathbf{1}}^{-1}
+ p_{t+1}
\mathcal{M}_{t+1, \mathbf{a}+\mathbf{2},
\mathbf{n}}^{-1}
q_{\mathbf{n}+\mathbf{1}}\bigr)
\mathcal{M}_{t+1, \mathbf{a}+\mathbf{1},
\mathbf{n}+\mathbf{1}}\\
&\quad -
\mathcal{M}_{t, \mathbf{a}+\mathbf{1},
\mathbf{n}}^{-1}
\bigl(p_{t+1}q_{\mathbf{n}+\mathbf{1}}
\mathcal{M}_{t, \mathbf{a},
\mathbf{n}+\mathbf{1}}
+ \mathcal{M}_{t+1, \mathbf{a}+\mathbf{1},
\mathbf{n}}
\Lambda_{\mathbf{a},
\mathbf{n}}\bigr) = 0.
\end{aligned}
\end{equation*}

\medskip

\noindent\textbf{3. Left Bernoulli Jumps.}
\hfill Theorem~\ref{thm:LBJ-H6}

\begin{equation*}
\bigl(q\mathcal{M}_{t+1,\mathbf{a}+\mathbf{1},
\mathbf{n}}^{-1}
+ p\mathcal{M}_{t,\mathbf{a}+\mathbf{1},
\mathbf{n}+\mathbf{1}}^{-1}\bigr)
\mathcal{M}_{t+1, \mathbf{a},
\mathbf{n}+\mathbf{1}}
- \mathcal{M}_{t,\mathbf{a}+\mathbf{1},
\mathbf{n}}^{-1}
\bigl(q\mathcal{M}_{t, \mathbf{a},
\mathbf{n}+\mathbf{1}}
+ p\mathcal{M}_{t+1, \mathbf{a},
\mathbf{n}}\bigr) = 0.
\end{equation*}

\medskip

\noindent\textbf{4. Parallel TASEP.}
\hfill Theorem~\ref{thm:PTASEP-H4}

\begin{equation*}
\bigl(q\mathcal{M}_{t+1, \mathbf{a}+\mathbf{1},
\mathbf{n}+\mathbf{1}}^{-1}
+ p\mathcal{M}_{t+1, \mathbf{a}+\mathbf{2},
\mathbf{n}}^{-1}\bigr)
\mathcal{M}_{t+2, \mathbf{a}+\mathbf{1},
\mathbf{n}+\mathbf{1}}
- \mathcal{M}_{t, \mathbf{a}+\mathbf{1},
\mathbf{n}}^{-1}
\bigl(q\mathcal{M}_{t+1, \mathbf{a}+\mathbf{1},
\mathbf{n}}
+ p\mathcal{M}_{t+1, \mathbf{a},
\mathbf{n}+\mathbf{1}}\bigr) = 0.
\end{equation*}

\medskip

\noindent\textbf{5. Right Geometric Jumps.}
\hfill Theorem~\ref{thm:RGP-H7}

\begin{equation*}
\bigl(\mathcal{M}_{t, \mathbf{a}+\mathbf{1},
\mathbf{n}+\mathbf{1}}^{-1}
- q\mathcal{M}_{t-1, \mathbf{a}+\mathbf{2},
\mathbf{n}}^{-1}\bigr)
\mathcal{M}_{t-1, \mathbf{a}+\mathbf{1},
\mathbf{n}+\mathbf{1}}
- \mathcal{M}_{t, \mathbf{a}+\mathbf{1},
\mathbf{n}}^{-1}
\bigl(\mathcal{M}_{t-1, \mathbf{a}+\mathbf{1},
\mathbf{n}}
- q\mathcal{M}_{t, \mathbf{a},
\mathbf{n}+\mathbf{1}}\bigr) = 0.
\end{equation*}

\medskip

\noindent\textbf{6. Left Geometric Jumps.}
\hfill Theorem~\ref{thm:LG-H8}

\begin{equation*}
\bigl(\mathcal{M}_{t, \mathbf{a}+\mathbf{1},
\mathbf{n}}^{-1}
- q\mathcal{M}_{t+1, \mathbf{a}+\mathbf{1},
\mathbf{n}+\mathbf{1}}^{-1}\bigr)
\mathcal{M}_{t, \mathbf{a},
\mathbf{n}+\mathbf{1}}
- \mathcal{M}_{t+1, \mathbf{a}+\mathbf{1},
\mathbf{n}}^{-1}
\bigl(\mathcal{M}_{t+1, \mathbf{a},
\mathbf{n}+\mathbf{1}}
- q\mathcal{M}_{t, \mathbf{a},
\mathbf{n}}\bigr) = 0.
\end{equation*}

\medskip

\noindent\textbf{7. Geometric LPP with Boundary.}
\hfill Theorem~\ref{thm:GLPP-H4p}

\begin{equation*}
\bigl(\mathcal{M}_{n-1, \mathbf{m},
\mathbf{a}-2\cdot\mathbf{1}}^{-1}
- q\mathcal{M}_{n, \mathbf{m}+\mathbf{1},
\mathbf{a}-\mathbf{1}}^{-1}\bigr)
\mathcal{M}_{n-1, \mathbf{m}+\mathbf{1},
\mathbf{a}-\mathbf{1}}
- \mathcal{M}_{n, \mathbf{m},
\mathbf{a}-\mathbf{1}}^{-1}
\bigl(\mathcal{M}_{n, \mathbf{m}+\mathbf{1},
\mathbf{a}}
- q\mathcal{M}_{n-1, \mathbf{m},
\mathbf{a}-\mathbf{1}}\bigr) = 0.
\end{equation*}

\medskip

\noindent\textbf{8. Stochastic Six-Vertex Model} (one-point).
\hfill Corollary~\ref{cor:S6V-scalar}

\begin{equation*}
(1 - b_1)
\mathcal{L}_{t,x+1,m+1}
\mathcal{L}_{t-1,x,m}
- (1 - b_2)
\mathcal{L}_{t-1,x,m+1}
\mathcal{L}_{t,x+1,m}
+ (b_1 - b_2)
\mathcal{L}_{t,x,m+1}
\mathcal{L}_{t-1,x+1,m}
= 0.
\end{equation*}

\bigskip
\begin{center}
\textsc{II. Euclidean Division, Vertex, and Polymer Models}
\end{center}
\medskip

\noindent\textbf{9. Higher-Spin Exclusion Process} (one-point).
\hfill Theorem~\ref{thm:HSEP-diamond}

\begin{equation*}
\begin{aligned}
&\mathcal{M}_{n+1,t+1,m+1}^{-1}
\bigl(C_{p_1}\mathcal{M}_{n+1,t,m}\Lambda_{p_2}
- C_{p_2}\mathcal{M}_{n,t+1,m}\Lambda_{p_1}\bigr)\\
&\quad +
\bigl(\Lambda_{p_1}\mathcal{M}_{n+2,t+1,m+1}^{-1}
C_{p_2}
- \Lambda_{p_2}\mathcal{M}_{n+1,t+2,m+1}^{-1}
C_{p_1}\bigr)\mathcal{M}_{n+1,t+1,m}
= 0.
\end{aligned}
\end{equation*}

\medskip

\noindent\textbf{10. Inhomogeneous Stochastic Six-Vertex Model} (one-point).
\hfill Theorem~\ref{thm:IS6V-diamond}

\begin{equation*}
\mathcal{M}_p(u) \Lambda_{p'}
+
\Lambda_p \mathcal{M}_{p'}(S_p u)
-
\mathcal{M}_{p'}(u) \Lambda_p
-
\Lambda_{p'} \mathcal{M}_p(S_{p'} u)
= 0.
\end{equation*}

\medskip

\noindent\textbf{11. Log-Gamma Polymer} (one-point).
\hfill Theorem~\ref{thm:LG-onepoint-matrix}

\begin{align*}
\mathcal{M}_1(u)\Lambda_j
+\Lambda_1\mathcal{M}_j(u)
-\mathcal{M}_j(u)\Lambda_1
-\Lambda_j\mathcal{M}_1(S_ju)
=s\partial_s\mathcal{M}_j(u).
\end{align*}

\medskip

\noindent\textbf{12. O'Connell--Yor Polymer} (one-point).
\hfill Theorem~\ref{thm:OY-finite-string-equation}

\begin{equation*}
\pa_\tau(\mathcal{A}_a - \mathcal{A}_{S_3a})
+\tfrac{1}{2}\pa_r^2(\mathcal{A}_a + \mathcal{A}_{S_3a})
-(\pa_r\mathcal{A}_a)(M_p + \mathcal{A}_a - \mathcal{A}_{S_3a})
+(M_p + \mathcal{A}_a - \mathcal{A}_{S_3a})\pa_r\mathcal{A}_{S_3a}
= 0.
\end{equation*}

\bigskip
\begin{center}
\textsc{III. Semi-Discrete Models}
\end{center}
\medskip

\noindent\textbf{13. Continuous-Time TASEP.}
\hfill Theorem~\ref{thm:CT-TASEP-multipoint}

\begin{equation*}
\mathcal{M}_{t,\mathbf{a}+\mathbf{1},
\mathbf{n}+\mathbf{1}}^{-1}
\pa_t \mathcal{M}_{t,\mathbf{a}+\mathbf{1},
\mathbf{n}+\mathbf{1}}
- \mathcal{M}_{t,\mathbf{a}+\mathbf{1},
\mathbf{n}}^{-1}
\pa_t \mathcal{M}_{t,\mathbf{a}+\mathbf{1},
\mathbf{n}}
+ \mathcal{M}_{t,\mathbf{a}+\mathbf{2},
\mathbf{n}}^{-1}
\mathcal{M}_{t,\mathbf{a}+\mathbf{1},
\mathbf{n}+\mathbf{1}}
- \mathcal{M}_{t,\mathbf{a}+\mathbf{1},
\mathbf{n}}^{-1}
\mathcal{M}_{t,\mathbf{a},\mathbf{n}+\mathbf{1}} = 0.
\end{equation*}

\medskip

\noindent\textbf{14. Push-TASEP.}
\hfill Theorem~\ref{thm:Push-TASEP-multipoint}

\begin{equation*}
 (\pa_t \mathcal{M}_{t,\mathbf{a},
\mathbf{n}+\mathbf{1}})
\mathcal{M}_{t,\mathbf{a},
\mathbf{n}+\mathbf{1}}^{-1}
- (\pa_t \mathcal{M}_{t,\mathbf{a}+\mathbf{1},
\mathbf{n}})
\mathcal{M}_{t,\mathbf{a}+\mathbf{1},
\mathbf{n}}^{-1}
 + \mathcal{M}_{t,\mathbf{a}+\mathbf{1},
\mathbf{n}}
\mathcal{M}_{t,\mathbf{a}+\mathbf{1},
\mathbf{n}+\mathbf{1}}^{-1}
- \mathcal{M}_{t,\mathbf{a},
\mathbf{n}}
\mathcal{M}_{t,\mathbf{a},
\mathbf{n}+\mathbf{1}}^{-1} = 0.
\end{equation*}

\medskip

\noindent\textbf{15. ASEP} (one-point).
\hfill Corollary~\ref{cor:asep-F-equation}

\begin{equation*}
\mathcal{L}_{t,x,m}\partial_t\mathcal{L}_{t,x,m+1}
-\mathcal{L}_{t,x,m+1}\partial_t\mathcal{L}_{t,x,m}
+\gamma\left[
\mathcal{L}_{t,x+1,m}\mathcal{L}_{t,x-1,m+1}
-\mathcal{L}_{t,x,m}\mathcal{L}_{t,x,m+1}
\right]=0.
\end{equation*}

\bigskip
\begin{center}
\textsc{IV. Parabolic Models}
\end{center}
\medskip

\noindent\textbf{16. Reflected Brownian Motion.}
\hfill Theorem~\ref{thm:RBM-multipoint}

\begin{equation*}
\pa_t (\mathcal{A}_{\mathbf{n}+\mathbf{1}}
- \mathcal{A}_{\mathbf{n}})
- \tfrac{1}{2}\pa_{\mathbf{a}}^2(\mathcal{A}_{\mathbf{n}+\mathbf{1}}
+ \mathcal{A}_{\mathbf{n}})
- \pa_{\mathbf{a}}(\mathcal{A}_{\mathbf{n}}
  \mathcal{A}_{\mathbf{n}+\mathbf{1}})
+ \mathcal{A}_{\mathbf{n}+\mathbf{1}}
  \pa_{\mathbf{a}} \mathcal{A}_{\mathbf{n}+\mathbf{1}}
+ \pa_{\mathbf{a}} \mathcal{A}_{\mathbf{n}}
  \mathcal{A}_{\mathbf{n}} = 0.
\end{equation*}

\medskip

\noindent\textbf{17. Brownian LPP with Boundary.}
\hfill Theorem~\ref{thm:BLPP-multipoint}

\begin{equation*}
\pa_{\mathbf{t}} (\mathcal{A}_{n+1} - \mathcal{A}_n)
- \tfrac{1}{2}\pa_{\mathbf{a}}^2(\mathcal{A}_{n+1}
+ \mathcal{A}_n)
+ \pa_{\mathbf{a}}(\mathcal{A}_n \mathcal{A}_{n+1})
- \mathcal{A}_{n+1} \pa_{\mathbf{a}} \mathcal{A}_{n+1}
- \pa_{\mathbf{a}} \mathcal{A}_n
  \mathcal{A}_n = 0.
\end{equation*}

\bigskip
\begin{center}
\textsc{V. The KPZ Fixed Point}
\end{center}
\medskip

\noindent\textbf{18. KPZ Fixed Point.}
\hfill Theorem~\ref{thm:FP-pKP}

\begin{equation*}
\pa_{t,\mathbf{a}} \mathcal{A}
+ \tfrac{1}{4}\pa_{\mathbf{x}}^2 \mathcal{A}
+ \tfrac{1}{12}\pa_{\mathbf{a}}^4 \mathcal{A}
+ \tfrac{1}{2}\pa_{\mathbf{a}}
  (\pa_{\mathbf{a}}\mathcal{A})^2
+ \tfrac{1}{2}[\pa_{\mathbf{a}}\mathcal{A},
\pa_{\mathbf{x}}\mathcal{A}] = 0.
\end{equation*}

\subsection*{Acknowledgments}
I thank my advisor, Jeremy Quastel, for generous support, patient
guidance, and many helpful discussions throughout the development of
this work.

%% file: chapters/2-the-discrete-framework.tex
\chapter{The Discrete Framework}
\label{ch:the-discrete-framework}

{
  \setlength{\parskip}{0pt}
}
\label{sec:discrete-framework}

This chapter introduces the diamond equations, which arise as compatibility conditions for an overdetermined linear problem on a directed lattice (Proposition~\ref{prop:CL-gen-system}), and the Darboux transformation (Theorem~\ref{thm:Darboux}) that dresses the edge weights while preserving the diamond structure. The scalar reduction extracts variable-coefficient Hirota--Miwa equations for Fredholm determinants (Proposition~\ref{prop:scalar-HM}), and the product graph construction (Theorem~\ref{thm:product-graph}) provides a systematic mechanism for building dressing-compatible data. For invertible edge weights $C_i$, the diamond equations are gauge equivalent to the non-abelian Hirota--Miwa system (Proposition~\ref{prop:diamond-NAHM}).

\section{The Linear Problem and Diamond Equations}\label{sec:linear-problem}
Let $E$ and $H \neq 0$ be vector spaces over a field $\mathbb{F}$, and let $\mathcal{G} = (\mathcal{V}, \mathcal{E})$ be a directed graph whose vertex set is a lattice  generated by invertible, mutually commuting shifts $S_i \colon \mathcal{V} \to \mathcal{V}$, $i \in \mathbb{Z}_{\geq 0}$, and whose edges take the form $e_i(u) \colon u \to S_i u$. To each edge $e_i(u)$, assign endomorphisms
\begin{equation*}
    C_i(u) \in \End(E), \qquad \Lambda_i(u) \in \End(E).
\end{equation*}
Fix a distinguished lattice direction with shift $T \defeq S_0 \colon \mathcal{V} \to \mathcal{V}$, so that the lattice is generated by $T, S_1, S_2, \dotsc$. Consider the overdetermined linear system for $\Psi:\mathcal V\to\Hom(E,H)$
\begin{equation}\label{eq:Linear-Psi-3}
\Psi(S_iu)=\Psi(Tu)C_i(u)-\Psi(u)\Lambda_i(u),\qquad  i \in \N.
\end{equation}
Compatibility of the linear system \eqref{eq:Linear-Psi-3} forces the following algebraic constraints on $(C_i, \Lambda_i)$.
\begin{proposition}\label{prop:CL-gen-system}
The system \eqref{eq:Linear-Psi-3} is compatible\footnote{That is, for any choice of initial data $\{\Psi(T^ku_0)\}_{k \in \Z}$ along a $T$-orbit, the value of $\Psi$ at any vertex reached by a finite composition of $S_i$-shifts is independent of the order of the composition.} if and only if for all $i, j\in \N$ and $u \in \mathcal V$,
\begin{align}
C_i(Tu)C_j(S_iu) &= C_j(Tu)C_i(S_ju), \label{eq:diamond-C-ij}\\
\Lambda_i(u)\Lambda_j(S_iu) &= \Lambda_j(u)\Lambda_i(S_ju), \label{eq:diamond-L-ij}\\
C_i(u)\Lambda_j(S_iu)+\Lambda_i(Tu)C_j(S_iu)&-C_j(u)\Lambda_i(S_ju)-\Lambda_j(Tu)C_i(S_ju)=0. \label{eq:diamond-mixed-ij}
\end{align}
\end{proposition}
\begin{proof}
Fix $i, j \in \N$ and $u \in \mathcal{V}$. We compute $\Psi(S_jS_iu)$ in two ways. Applying \eqref{eq:Linear-Psi-3} with index~$j$ at the vertex $S_iu$
gives
\[
\Psi(S_jS_iu) = \Psi(TS_iu)C_j(S_iu) - \Psi(S_iu)\Lambda_j(S_iu).
\]
Since $TS_i = S_iT$, the recurrence \eqref{eq:Linear-Psi-3} with
index~$i$ expresses $\Psi(S_iTu)$ and $\Psi(S_iu)$ in terms of
$\Psi(T^2u)$, $\Psi(Tu)$, and $\Psi(u)$.  Substituting and
collecting:
\begin{align*}
\Psi(S_jS_iu)
&=\bigl(\Psi(T^2u)C_i(Tu)-\Psi(Tu)\Lambda_i(Tu)\bigr)C_j(S_iu)
\\
&\quad -\bigl(\Psi(Tu)C_i(u)-\Psi(u)\Lambda_i(u)\bigr)\Lambda_j(S_iu)\\
&=\Psi(T^2u)C_i(Tu)C_j(S_iu)
-\Psi(Tu)\bigl[\Lambda_i(Tu)C_j(S_iu)
+C_i(u)\Lambda_j(S_iu)\bigr]
\\
&\quad +\Psi(u)\Lambda_i(u)\Lambda_j(S_iu).
\end{align*}
Performing the same calculation with $i$ and $j$ interchanged and subtracting the two results gives
\begin{align}
\Psi(S_jS_iu)-\Psi(S_iS_ju) 
&=\Psi(T^2u)\bigl[C_i(Tu)C_j(S_iu)-C_j(Tu)C_i(S_ju)\bigr]\notag\\
& -\Psi(Tu)\bigl[C_i(u)\Lambda_j(S_iu)+\Lambda_i(Tu)C_j(S_iu)
\notag\\
&\qquad\qquad\quad -C_j(u)\Lambda_i(S_ju)-\Lambda_j(Tu)C_i(S_ju)\bigr]\notag\\
& +\Psi(u)\bigl[\Lambda_i(u)\Lambda_j(S_iu)-\Lambda_j(u)\Lambda_i(S_ju)\bigr]. \label{eq:compat-identity}
\end{align}

Now suppose the system \eqref{eq:Linear-Psi-3} is compatible, so that
$\Psi(S_jS_iu) = \Psi(S_iS_ju)$ for every solution~$\Psi$.  Since the values $\Psi(u), \Psi(Tu), \Psi(T^2u)$ may be independently prescribed in $\Hom(E, H)$, each bracketed coefficient must vanish, giving \eqref{eq:diamond-C-ij}--\eqref{eq:diamond-mixed-ij}.

Conversely, assume \eqref{eq:diamond-C-ij}--\eqref{eq:diamond-mixed-ij}.
Then \eqref{eq:compat-identity} gives
$\Psi(S_jS_iu) = \Psi(S_iS_ju)$ for every $u \in \mathcal{V}$ and
every $i, j \in \N$.  Since any two permutations of a fixed finite
collection of shifts differ by a sequence of such pairwise
transpositions, the value of $\Psi(S_{i_1} \cdots S_{i_r}u)$ is
independent of the ordering.  Hence \eqref{eq:Linear-Psi-3} is
compatible.
\end{proof}
We call \eqref{eq:diamond-C-ij}--\eqref{eq:diamond-mixed-ij} the \emph{diamond equations}.
\begin{remark}
In KPZ applications, we will often take $E = \R^m$ and $H = L^2(\Z)$, so that $\Psi(u)$ is a row vector $(\Psi_1(u), \dotsc, \Psi_m(u))$ with entries in $H$ and the endomorphisms $C_i(u), \Lambda_i(u)$ are $m \times m$ matrices encoding the data of the physical model. The lattice $\mathcal V$ encodes the physical parameter space and the directed edges encode shifts in those parameters, e.g.\ $(t, a, n) \in \mathcal V$ with $T(t, a, n) = (t, a+1, n)$, $S_1(t, a, n) = (t+1, a, n)$, and $S_2(t, a, n) = (t, a, n+1)$. In each model of this work the lattice is generated by two shifts beyond $T$, so the diamond equations reduce to \eqref{eq:diamond-C-ij}--\eqref{eq:diamond-mixed-ij} with $i, j \in \{1, 2\}$.
\end{remark}

\section{Darboux Transformations}\label{sec:Darboux-transformations}
The Darboux transformation produces, from a solution $(C_i, \Lambda_i)$ of the diamond equations together with wave functions and a compatible kernel, new edge weights $(\mathcal{M}_i, \Lambda_i)$ satisfying the same diamond equations.  We now develop this construction.

Fix a solution $(C_i, \Lambda_i)$ of
\eqref{eq:diamond-C-ij}--\eqref{eq:diamond-mixed-ij}, and let
$\Psi \colon \mathcal{V} \to \Hom(E, H)$ be a solution of the
linear problem \eqref{eq:Linear-Psi-3}.  The formal adjoint
linear problem\footnote{Its compatibility conditions are again
equivalent to the diamond equations
\eqref{eq:diamond-C-ij}--\eqref{eq:diamond-mixed-ij}; the
verification is analogous to
Proposition~\ref{prop:CL-gen-system}.} for
$\Phi \colon \mathcal{V} \to \Hom(H, E)$ is
\begin{equation}\label{eq:Linear-Phi-Prob}
\Phi(Tu) = C_i(u)\Phi(S_i u) - \Lambda_i(Tu)\Phi(TS_i u),
\qquad i \in \mathbb{N}.
\end{equation}

\begin{definition}
We say that $K \colon \mathcal{V} \to \End(H)$ is \emph{dressing
compatible} with $(\Psi, \Phi)$ if:
\begin{enumerate}[label=\arabic*., leftmargin=*]
\item For every $u \in \mathcal{V}$ and $i \in \mathbb{N}$:
\begin{align}\label{eq:K-differences}
K(Tu) - K(u) &= \Psi(Tu)\Phi(Tu), \quad 
K(S_i u) - K(u) = \Psi(Tu)C_i(u)\Phi(S_i u). 
\end{align}
\item There exists $z \in \mathbb{F}$ such that the resolvent
exists at every vertex:
\[
R(u) \defeq (I - zK(u))^{-1} \in \End(H),
\qquad u \in \mathcal{V}.
\]
\end{enumerate}
\end{definition}
 
\begin{example}\label{ex:simple-kernel}
Under suitable analytic assumptions, the kernel
\[
K(u) = \sum_{p \leq 0} \Psi(T^p u)\Phi(T^p u)
\]
is dressing compatible with $(\Psi, \Phi)$; the verification
is a telescoping argument.
\end{example}
 
\begin{theorem}\label{thm:Darboux}
Let $K$ be dressing compatible with $(\Psi, \Phi)$, with resolvent $R(u) = (I - zK(u))^{-1}$.  Define the dressed observable and dressed
waves by
\[
\mathcal{M}(u) \defeq I + z\Phi(u)R(u)\Psi(u),
\]
\[
\widehat\Psi(u) \defeq R(u)\Psi(u), \qquad
\widehat\Phi(u) \defeq \Phi(u)R(T^{-1}u).
\]
Then $\mathcal{M}(u)$ is invertible for every
$u \in \mathcal{V}$. The dressed waves
$\widehat\Psi, \widehat\Phi$ satisfy the
linear system \eqref{eq:Linear-Psi-3} and adjoint linear
system \eqref{eq:Linear-Phi-Prob} respectively, with data
$(\mathcal{M}_i, \Lambda_i)$:
\begin{align*}
\widehat\Psi(S_i u)
&= \widehat\Psi(Tu)\mathcal{M}_i(u)
- \widehat\Psi(u)\Lambda_i(u), \\
\widehat\Phi(Tu)
&= \mathcal{M}_i(u)\widehat\Phi(S_i u)
- \Lambda_i(Tu)\widehat\Phi(TS_i u),
\end{align*}
for $i \in \mathbb{N}$, where the dressed edge weights are given by
\begin{equation}\label{eq:dressed-C}
\mathcal{M}_i(u) \defeq
\mathcal{M}(Tu)^{-1}C_i(u)\mathcal{M}(S_i u),
\end{equation}
and $\widehat K(u) \defeq K(u)R(u)$ is dressing compatible
with $(\widehat\Psi, \widehat\Phi)$.
Moreover, the pair $(\mathcal{M}_i, \Lambda_i)$ satisfies
the diamond equations
\eqref{eq:diamond-C-ij}--\eqref{eq:diamond-mixed-ij}.
\end{theorem}
 
\begin{remark}\label{rem:Darboux-via-wave}
One could alternatively derive the diamond equations for $(\mathcal{M}_i, \Lambda_i)$ from the compatibility of the dressed linear problem via the calculation of Proposition~\ref{prop:CL-gen-system}: the $\widehat\Psi(u)$-bracket vanishes because $\Lambda$ is unchanged, and the $\widehat\Psi(T^2u)$-bracket vanishes by telescoping, 
\[\mathcal{M}_i(Tu)\mathcal{M}_j(S_iu) = \mathcal{M}(T^2u)^{-1}C_i(Tu)C_j(S_iu)\mathcal{M}(S_iS_ju),\] 
which is symmetric in $i, j$ by the $C$-diamond equation. The remaining $\widehat\Psi(Tu)$-bracket, however, requires $\widehat\Psi(Tu)$ to be injective, which is not assumed here. The direct proof below circumvents the injectivity requirement entirely, using only resolvent existence.
\end{remark}
 
The proof of Theorem~\ref{thm:Darboux} rests on the following resolvent identities. For the remainder of this section, we suppress functional dependence on $z$.
 
\begin{lemma}\label{lem:resolvent-ids}
For every $u \in \mathcal{V}$ and $i \in \mathbb{N}$, the following resolvent identities hold:
\begin{align}
R(Tu) - R(u)
  &= zR(u)\Psi(Tu)\Phi(Tu)R(Tu),
  \label{eq:resolvent-T} \\
R(S_i u) - R(u)
  &= zR(u)\Psi(Tu)C_i(u)\Phi(S_i u)R(S_i u).
  \label{eq:resolvent-Si}
\end{align}
\end{lemma}
 
\begin{proof}
For any $v, w \in \mathcal{V}$, since $R(v)^{-1} = I - zK(v)$
we have $R(v)^{-1} - R(w)^{-1} = z(K(w) - K(v))$, and
therefore
\begin{equation}\label{eq:resolvent-general}
R(w) - R(v) = R(v)\bigl(R(v)^{-1} - R(w)^{-1}\bigr)R(w)
= zR(v)\bigl(K(w) - K(v)\bigr)R(w).
\end{equation}
Substituting the dressing compatibility conditions
\eqref{eq:K-differences} gives
\eqref{eq:resolvent-T}--\eqref{eq:resolvent-Si}.
\end{proof}
 
\begin{proof}[Proof of Theorem~\ref{thm:Darboux}]
We prove the theorem in five steps: invertibility of \(\mathcal M\),
the dressed linear problem, the dressed adjoint problem, the
dressed kernel, and the dressed diamond equations.

\noindent\textit{Invertibility of $\mathcal{M}$.}
We show that
\begin{equation}\label{eq:M-inverse}
\mathcal{M}(Tu)^{-1} = I - z\Phi(Tu)R(u)\Psi(Tu).
\end{equation}
Collecting the $z$-linear term as $z\Phi(Tu)(R(Tu) - R(u))\Psi(Tu)$ and using \eqref{eq:resolvent-T} to reduce the $z^2$-term,
\begin{align*}
    &(I - z\Phi(Tu)R(u)\Psi(Tu))(I + z\Phi(Tu)R(Tu)\Psi(Tu)) \\
    &\qquad = I + z\Phi(Tu)\bigl(R(Tu) - R(u)\bigr)\Psi(Tu) \\
    &\qquad\qquad - z^2\Phi(Tu)R(u)\Psi(Tu)\Phi(Tu)R(Tu)\Psi(Tu)   \\
    &\qquad = I + z\Phi(Tu)\bigl(R(Tu) - R(u)\bigr)\Psi(Tu) - z\Phi(Tu)\bigl(R(Tu) - R(u)\bigr)\Psi(Tu) \\
    &\qquad = I.
\end{align*}
The right inverse follows identically from the opposite factorization
\[
R(Tu) - R(u) = zR(Tu)\Psi(Tu)\Phi(Tu)R(u).
\]

\noindent\textit{The dressed linear problem.}
We begin by expanding $\widehat \Psi(Tu)$ using
the definition of $\mathcal{M}$:
\begin{align*}
    \widehat \Psi(Tu) &= R(Tu) \Psi(Tu) = R(u)\Psi(Tu) + (R(Tu) - R(u))\Psi(Tu)  \\
    &=  R(u) \Psi(Tu)+ zR(u)\Psi(Tu)\Phi(Tu)R(Tu)\Psi(Tu) \\
    &= R(u) \Psi(Tu) \mathcal{M}(Tu),
\end{align*}
where we used \eqref{eq:resolvent-T}. Therefore
\begin{equation}\label{eq:intertwine-wave}
 \widehat\Psi(Tu)\mathcal{M}(Tu)^{-1} = R(u)\Psi(Tu).
\end{equation}
 
Fix $i \in \mathbb{N}$ and expand the right-hand side of
the dressed linear problem.  Substituting the definition
\eqref{eq:dressed-C} and applying
\eqref{eq:intertwine-wave}:
\begin{align*}
\widehat\Psi(Tu)\mathcal{M}_i(u)
- \widehat\Psi(u)\Lambda_i(u)
&= R(u)\Psi(Tu)C_i(u)\mathcal{M}(S_i u)
- R(u)\Psi(u)\Lambda_i(u).
\end{align*}
Expanding
$\mathcal{M}(S_i u) = I +
z\Phi(S_i u)R(S_i u)\Psi(S_i u)$ and collecting
the $R(u)$-terms:
\begin{align*}
&= R(u)\bigl[\Psi(Tu)C_i(u) -
\Psi(u)\Lambda_i(u)\bigr]
+ zR(u)\Psi(Tu)C_i(u)\Phi(S_i u)R(S_i u)\Psi(S_i u).
\end{align*}
By the linear problem \eqref{eq:Linear-Psi-3}, the
bracket equals $\Psi(S_i u)$, and by
\eqref{eq:resolvent-Si} the second term equals
$\bigl(R(S_i u) - R(u)\bigr)\Psi(S_i u)$.  Therefore
\begin{align*}
&= R(u)\Psi(S_i u) + \bigl(R(S_i u) -
R(u)\bigr)\Psi(S_i u)
= R(S_i u)\Psi(S_i u) = \widehat\Psi(S_i u).
\end{align*}
 
\noindent\textit{The adjoint dressed linear problem.}
We begin by expanding
$\mathcal{M}(Tu)\widehat \Phi(Tu)$ using the definition of
$\mathcal{M}$:
\begin{align*}
\mathcal{M}(Tu)\widehat \Phi(Tu)
&= \Phi(Tu)R(u)
   + z\Phi(Tu)R(Tu)\Psi(Tu)\Phi(Tu)R(u) \\
&= \Phi(Tu)R(u)
   + \Phi(Tu)\bigl(R(Tu) - R(u)\bigr)
= \Phi(Tu)R(Tu),
\end{align*}
where we used \eqref{eq:resolvent-T} in the form
$R(Tu) - R(u) = zR(Tu)\Psi(Tu)\Phi(Tu)R(u)$.
Therefore
\begin{equation}\label{eq:adjoint-intertwine-wave}
\widehat \Phi(Tu) = \mathcal{M}(Tu)^{-1}\Phi(Tu)R(Tu).
\end{equation}

Fix $i \in \mathbb{N}$ and expand the right-hand side of
the adjoint dressed linear problem.  Substituting the
definition \eqref{eq:dressed-C} and applying
\eqref{eq:adjoint-intertwine-wave} at $T^{-1}S_i u$ to
collapse
$\mathcal{M}(S_i u)\Phi(S_i u)R(T^{-1}S_i u)
= \Phi(S_i u)R(S_i u)$:
\begin{align*}
\mathcal{M}_i(u)\widehat\Phi(S_i u)
- \Lambda_i(Tu)\widehat\Phi(TS_i u)
= \bigl[\mathcal{M}(Tu)^{-1}C_i(u)\Phi(S_i u)
- \Lambda_i(Tu)\Phi(TS_i u)\bigr]R(S_i u).
\end{align*}
Expanding
$\mathcal{M}(Tu)^{-1} = I -
z\Phi(Tu)R(u)\Psi(Tu)$ and collecting:
\begin{align*}
&= \bigl[C_i(u)\Phi(S_i u) -
\Lambda_i(Tu)\Phi(TS_i u)\bigr]R(S_i u)
- z\Phi(Tu)R(u)\Psi(Tu)C_i(u)\Phi(S_i u) R(S_i u).
\end{align*}
By the adjoint linear problem
\eqref{eq:Linear-Phi-Prob}, the bracket equals $\Phi(Tu)$:
\begin{align*}
&= \Phi(Tu)\bigl[I -
zR(u)\Psi(Tu)C_i(u)\Phi(S_i u)\bigr]R(S_i u).
\end{align*}
By \eqref{eq:resolvent-Si} rearranged as
$\bigl[I - zR(u)\Psi(Tu)C_i(u)\Phi(S_i u)\bigr]
R(S_i u) = R(u)$:
\begin{align*}
&= \Phi(Tu)R(u) = \widehat\Phi(Tu).
\end{align*}
 
\noindent\textit{The dressed kernel.}
Since $K(u)$ and
$R(u)$ commute at each vertex, $\widehat K(u) = R(u)K(u)$.
The identity
$(I - \eta\widehat K(u))(I - zK(u)) = I - (z+\eta)K(u)$
shows that $(I - \eta\widehat K(u))^{-1}$ exists whenever
$(I - (z+\eta)K(u))^{-1}$ does; in particular at
$\eta = -z$, where the condition is vacuous.  When $z = 0$,
$R = I$ and $\widehat K = K$, so the kernel identities
reduce to those of $K$ itself.  For $z \neq 0$, it remains
to verify the kernel identities.

\noindent\textit{The $T$-condition.}
Since $\widehat K(Tu) - \widehat K(u)
= z^{-1}(R(Tu) - R(u))$,
\eqref{eq:resolvent-general} and
\eqref{eq:K-differences} give
\[
\widehat K(Tu) - \widehat K(u)
= R(u)\bigl(K(Tu) - K(u)\bigr)R(Tu)
= R(u)\Psi(Tu)\Phi(Tu)R(Tu).
\]
The intertwining identity \eqref{eq:intertwine-wave}
gives $R(u)\Psi(Tu)
= \widehat\Psi(Tu)\mathcal{M}(Tu)^{-1}$,
and \eqref{eq:adjoint-intertwine-wave} gives
$\Phi(Tu)R(Tu)
= \mathcal{M}(Tu)\widehat\Phi(Tu)$.
Substituting:
\begin{align*}
R(u)\Psi(Tu)\Phi(Tu)R(Tu)
= \widehat\Psi(Tu)\mathcal{M}(Tu)^{-1}
   \mathcal{M}(Tu)\widehat\Phi(Tu)
= \widehat\Psi(Tu)\widehat\Phi(Tu).
\end{align*}

\noindent\textit{The $S_i$-condition.}
Fix $i \in \N$.  By the same argument as above,
with \eqref{eq:resolvent-Si} and
\eqref{eq:K-differences},
\[
\widehat K(S_iu) - \widehat K(u)
= R(u)\Psi(Tu)C_i(u)\Phi(S_iu)R(S_iu).
\]
Applying \eqref{eq:intertwine-wave} and
\eqref{eq:adjoint-intertwine-wave} at $T^{-1}S_iu$:
\begin{align*}
R(u)\Psi(Tu)C_i(u)\Phi(S_iu)R(S_iu)
&= \widehat\Psi(Tu)\mathcal{M}(Tu)^{-1}
   C_i(u)\mathcal{M}(S_iu)\widehat\Phi(S_iu) \\
&= \widehat\Psi(Tu)\mathcal{M}_i(u)\widehat\Phi(S_iu).
\end{align*}

\noindent\textit{The diamond equations.}
For \eqref{eq:diamond-C-ij}, fix $i, j \in \mathbb{N}$ and
expand the left-hand side:
\begin{align*}
&\mathcal{M}_i(Tu)\mathcal{M}_j(S_i u) \\
&\qquad = \mathcal{M}(T^2u)^{-1} C_i(Tu)
   \mathcal{M}(TS_i u)\mathcal{M}(TS_i u)^{-1}
   C_j(S_i u)\mathcal{M}(S_i S_j u) \\
&\qquad = \mathcal{M}(T^2u)^{-1} C_i(Tu)C_j(S_i u)
   \mathcal{M}(S_i S_j u).
\end{align*}
Similarly,
$\mathcal{M}_j(Tu)\mathcal{M}_i(S_j u) =
\mathcal{M}(T^2u)^{-1} C_j(Tu)C_i(S_j u)
\mathcal{M}(S_i S_j u)$, and \eqref{eq:diamond-C-ij} for
$(\mathcal{M}_i, \Lambda_i)$ follows from
\eqref{eq:diamond-C-ij} for $(C_i, \Lambda_i)$.
In addition, since $\Lambda$ is unchanged, \eqref{eq:diamond-L-ij} holds
immediately for $(\mathcal{M}_i, \Lambda_i)$.
 
It remains to verify the mixed equation. Fix
$i, j \in \mathbb{N}$; we must show
\begin{equation}\label{eq:diamond-mixed-dressed}
\mathcal{M}_i(u)\Lambda_j(S_i u)
+ \Lambda_i(Tu)\mathcal{M}_j(S_i u)
- \mathcal{M}_j(u)\Lambda_i(S_j u)
- \Lambda_j(Tu)\mathcal{M}_i(S_j u) = 0.
\end{equation}
 
\noindent\textit{Explicit formula for $\mathcal{M}_i(u)$.}
Fix $i \in \mathbb{N}$. Substituting the inverse formula
\eqref{eq:M-inverse} into the definition \eqref{eq:dressed-C}
and expanding gives
\begin{align*}
\mathcal{M}_i(u)
&= C_i(u)
   + zC_i(u)\Phi(S_i u)R(S_i u)\Psi(S_i u)
   - z\Phi(Tu)R(u)\Psi(Tu)C_i(u) \\
&\quad
   - z^2\Phi(Tu)R(u)\Psi(Tu)C_i(u)
     \Phi(S_i u)R(S_i u)\Psi(S_i u).
\end{align*}
By \eqref{eq:resolvent-Si}, the $z^2$-term equals
$-z\Phi(Tu)\bigl(R(S_i u) - R(u)\bigr)\Psi(S_i u)$.
Collecting terms based at $R(S_i u)$ and $R(u)$ separately
yields
\begin{align*}
\mathcal{M}_i(u)
&= C_i(u)
  + z\bigl(C_i(u)\Phi(S_i u) - \Phi(Tu)\bigr)
    R(S_i u)\Psi(S_i u)
\\
&\quad  + z\Phi(Tu)R(u)
    \bigl(\Psi(S_i u) - \Psi(Tu)C_i(u)\bigr).
\end{align*}
Substituting the linear problems
\eqref{eq:Linear-Psi-3} and \eqref{eq:Linear-Phi-Prob}:
\begin{equation}\label{eq:dressed-C-explicit}
\mathcal{M}_i(u)
= C_i(u)
  + z\Lambda_i(Tu)\Phi(TS_i u)R(S_i u)\Psi(S_i u)
  - z\Phi(Tu)R(u)\Psi(u)\Lambda_i(u).
\end{equation}
 
Substituting \eqref{eq:dressed-C-explicit} for each of the
four terms, the contributions independent of $z$ combine to
give the mixed diamond equation \eqref{eq:diamond-mixed-ij}
for $(C, \Lambda)$, which vanishes by assumption. Adjacent
pairs cancel among the $z$-linear contributions, and the
surviving terms group into $\Lambda\Lambda$-brackets that
vanish by \eqref{eq:diamond-L-ij}. Explicitly, the
$z$-linear contributions are:
\begin{align*}
&z\Lambda_i(Tu)\Phi(TS_i u)R(S_i u)
  \Psi(S_i u)\Lambda_j(S_i u) \tag*{(A)}
\\
&\quad - z\Phi(Tu)R(u)\Psi(u)
  \Lambda_i(u)\Lambda_j(S_i u),  \\
&z\Lambda_i(Tu)\Lambda_j(TS_i u)
  \Phi(TS_i S_j u)R(S_i S_j u)\Psi(S_i S_j u) \tag*{(B)}
\\
&\quad - z\Lambda_i(Tu)\Phi(TS_i u)R(S_i u)
  \Psi(S_i u)\Lambda_j(S_i u),  \\
&{-z}\Lambda_j(Tu)\Phi(TS_j u)R(S_j u)
  \Psi(S_j u)\Lambda_i(S_j u) \tag*{(C)}
\\
&\quad + z\Phi(Tu)R(u)\Psi(u)
  \Lambda_j(u)\Lambda_i(S_j u),  \\
&{-z}\Lambda_j(Tu)\Lambda_i(TS_j u)
  \Phi(TS_i S_j u)R(S_i S_j u)\Psi(S_i S_j u) \tag*{(D)}
\\
&\quad + z\Lambda_j(Tu)\Phi(TS_j u)R(S_j u)
  \Psi(S_j u)\Lambda_i(S_j u). 
\end{align*}
The first term of (A) cancels with the second term of (B),
and similarly the first term of (C) cancels with the second
term of (D). The surviving terms group as
\begin{align*}
&z\Phi(Tu)R(u)\Psi(u)
  \bigl[\Lambda_j(u)\Lambda_i(S_j u)
  - \Lambda_i(u)\Lambda_j(S_i u)\bigr] \\
&\quad + z\bigl[\Lambda_i(Tu)\Lambda_j(TS_i u)
  - \Lambda_j(Tu)\Lambda_i(TS_j u)\bigr]\Phi(TS_i S_j u)R(S_i S_j u)\Psi(S_i S_j u),
\end{align*}
and both brackets vanish by \eqref{eq:diamond-L-ij}. This
proves \eqref{eq:diamond-mixed-dressed}.
\end{proof}

\begin{remark}
One can also dress $\Lambda$ instead of $C$: defining
\[
\widetilde{\mathcal{M}}_i(u) \defeq
\mathcal{M}(u)\Lambda_i(u)\mathcal{M}(S_i u)^{-1},
\]
the pair $(C_i(u), \widetilde{\mathcal{M}}_i(u))$ again
satisfies the diamond equations.
\end{remark}
 
\begin{remark}\label{rem:resolvent-subset}
In applications, the resolvent may exist only on a
subset $\mathcal{V}' \subseteq \mathcal{V}$.  Each
diamond equation at a vertex~$u$ holds whenever the
resolvent exists at the finitely many shifts of~$u$
appearing in the proof.  For the mixed diamond with
indices~$i, j$, the explicit
formula~\eqref{eq:dressed-C-explicit}, which
uses~\eqref{eq:M-inverse} to avoid the resolvent at
$T$-shifted vertices, requires the resolvent
only at $u$, $S_i u$, $S_j u$, and~$S_i S_j u$.
In particular, $\mathcal{V}'$ need not be closed
under~$T$.
\end{remark}

\subsection{Fredholm determinants and the scalar reduction}\label{sec:scalar-HM}
In the probabilistic applications of subsequent chapters,
the Fredholm determinants
$F(u) = \det_H(I - zK(u))$ encode the
distributional data of each
model.\footnote{In these applications
$H$ is typically an $L^2$ space, so $F(u)$ is
obtained as a Fredholm determinant in the classical
sense.}
The corollary below establishes the fundamental
relation between $F$ and $\mathcal M$.
 
\begin{corollary}\label{cor:Woodbury}
Let $K$ be dressing compatible with $(\Psi, \Phi)$,
and suppose that $\dim E < \infty$ and that
$F(u) \defeq \det_H(I - zK(u))$
is well-defined for every $u \in \mathcal{V}$.  Then
\begin{equation}\label{eq:Woodbury-det}
F(u) = F(Tu)\det_E\bigl(\mathcal{M}(Tu)\bigr).
\end{equation}
In particular, for each $i \in \mathbb{N}$,
\begin{equation}\label{eq:det-dressed}
\det_E(\mathcal{M}_i(u))F(u)F(S_i u)
= \det_E(C_i(u))F(Tu)F(T^{-1}S_i u).
\end{equation}
\end{corollary}
 
\begin{proof}
Since $(I - zK(u))^{-1}$ exists by hypothesis,
$F(u) \neq 0$.  By dressing compatibility,
\begin{align*}
I - zK(u) = (I - zK(Tu)) + z\Psi(Tu)\Phi(Tu)
= (I-zK(Tu))\bigl(I + zR(Tu)\Psi(Tu)\Phi(Tu)\bigr).
\end{align*}
Taking $\det_H$ of both sides and using
multiplicativity of the determinant:
\[
F(u) = F(Tu) 
\det_H\bigl(I + zR(Tu)\Psi(Tu)\Phi(Tu)\bigr).
\]
Since $\dim E < \infty$, Sylvester's determinant
identity gives
\begin{align*}
  \det_H\bigl(I + zR(Tu)\Psi(Tu)\Phi(Tu)\bigr)
= \det_E\bigl(I + z\Phi(Tu)R(Tu)\Psi(Tu)\bigr)
= \det_E\bigl(\mathcal{M}(Tu)\bigr),
\end{align*}
establishing \eqref{eq:Woodbury-det}.
Identity~\eqref{eq:det-dressed} follows by
taking $\det_E$ of
$\mathcal{M}_i(u) = \mathcal{M}(Tu)^{-1}C_i(u)
\mathcal{M}(S_i u)$, applying
\eqref{eq:Woodbury-det}, and rearranging.
\end{proof}

When $E = \mathbb{F}$, the edge weights
$C_i, \Lambda_i, \mathcal{M}_i$ reduce to
scalar-valued functions $c_i, \lambda_i, \mathcal{M}_i$,
and Corollary~\ref{cor:Woodbury} gives
\begin{align}
    \mathcal{M}_i(u) = c_i(u)
    \frac{F(Tu)F(T^{-1}S_iu)}{F(u)F(S_iu)}.
    \label{eq:M-scalar}
\end{align}
We now show that the mixed diamond equation for $(\mathcal{M}_i, \lambda_i)$ reduces to a variable-coefficient Hirota--Miwa equation for $F$.  Write $\alpha_{ij}(u) \defeq c_i(Tu)\lambda_j(Tu)$.

\begin{proposition}\label{prop:scalar-HM}
Let $(c_i, \lambda_i)$ be a scalar solution of the
diamond equations
\eqref{eq:diamond-C-ij}--\eqref{eq:diamond-mixed-ij}, and
let $F(u) = \det_H(I - zK(u))$ be as in
Corollary~\ref{cor:Woodbury}.
Suppose that\footnote{This boundary condition is natural in the probabilistic applications of Chapter~\ref{ch:discrete-particle-models}.} $F(T^\ell v) \to 1$ as $\ell \to -\infty$ for all $v \in \mathcal V$, and that $\alpha_{ij}(u)$ is nonzero with $\alpha_{ij}(u) \neq \alpha_{ji}(u)$ for all $u \in \mathcal V$ and $i \neq j$, with the same conditions holding for the limits $\lim_{\ell \to -\infty} \alpha_{ij}(T^\ell v)$, which are assumed to exist.
Then for each pair
$i, j \in \mathbb{N}$ with $i \neq j$, the mixed
diamond equation \eqref{eq:diamond-mixed-ij} for
$(\mathcal{M}_i, \lambda_i)$ reduces to the scalar
Hirota--Miwa equation with variable coefficients:
    \begin{align}
        \alpha_{ij}(u) F(S_iu)F(TS_ju)
        - \alpha_{ji}(u) F(S_ju)F(TS_iu)
        - \bigl(\alpha_{ij}(u) - \alpha_{ji}(u)\bigr)
          F(Tu)F(S_iS_ju) = 0.
        \label{eq:HM-variable}
    \end{align}
\end{proposition}
 
\begin{proof}
We substitute \eqref{eq:M-scalar} into the mixed diamond
equation \eqref{eq:diamond-mixed-ij} for
$(\mathcal{M}_i, \lambda_i)$, use
\eqref{eq:diamond-C-ij}--\eqref{eq:diamond-L-ij} to
reduce the coefficients, and identify the resulting
expression as a quantity conserved along $T$-orbits whose
value is determined by the boundary conditions. Rearranging
will give \eqref{eq:HM-variable}.
 
Fix $i, j \in \mathbb{N}$ with $i \neq j$.
Since $(I - zK(v))^{-1}$ exists by hypothesis, $F(v) \neq 0$.
Substituting \eqref{eq:M-scalar} into
\eqref{eq:diamond-mixed-ij} gives
\begin{align}
    &c_i(u)\lambda_j(S_iu)
      \frac{F(Tu)F(T^{-1}S_iu)}{F(u)F(S_iu)}
    + \lambda_i(Tu) c_j(S_iu)
      \frac{F(TS_iu)F(T^{-1}S_iS_ju)}{F(S_iu)F(S_iS_ju)} \notag \\
    & - c_j(u)\lambda_i(S_ju)
      \frac{F(Tu)F(T^{-1}S_ju)}{F(u)F(S_ju)}
    - \lambda_j(Tu) c_i(S_ju)
      \frac{F(TS_ju)F(T^{-1}S_iS_ju)}{F(S_ju)F(S_iS_ju)} = 0.
      \label{eq:HM-expanded}
\end{align}
Multiplying  by
$F(S_iu)F(S_ju)/\bigl(F(Tu)F(T^{-1}S_iS_ju)\bigr)$
and grouping terms with common denominators:
\begin{align}
    &\frac{1}{F(u)F(T^{-1}S_iS_ju)}
      \bigl(c_i(u)\lambda_j(S_iu) F(T^{-1}S_iu)F(S_ju)
\notag \\
    &\qquad\qquad
      - c_j(u)\lambda_i(S_ju) F(T^{-1}S_ju)F(S_iu)\bigr) \notag \\
    &\quad = \frac{1}{F(Tu)F(S_iS_ju)}
      \bigl(\lambda_j(Tu) c_i(S_ju) F(TS_ju)F(S_iu)
\notag \\
    &\qquad\qquad
      - \lambda_i(Tu) c_j(S_iu) F(TS_iu)F(S_ju)\bigr).
      \label{eq:HM-cleared}
\end{align}
 
\noindent\textit{Coefficient reduction.}
Since $\alpha_{ij}(u) \neq 0$, equations
\eqref{eq:diamond-C-ij}--\eqref{eq:diamond-L-ij} give
$c_i(S_ju) = c_i(Tu)c_j(S_iu)/c_j(Tu)$ and
$\lambda_j(S_iu) = \lambda_j(u)\lambda_i(S_ju)/\lambda_i(u)$.
Substituting into \eqref{eq:HM-cleared}, replacing
$\lambda_j(S_iu)$ on the left and $c_i(S_ju)$ on the
right, factoring the resulting common scalar
prefactors from each side, and recalling
$\alpha_{ij}(u) = c_i(Tu)\lambda_j(Tu)$:
\begin{align}
    &\frac{\lambda_i(S_ju)}{\lambda_i(u)}
      \frac{\alpha_{ij}(T^{-1}u) F(T^{-1}S_iu)F(S_ju)
        - \alpha_{ji}(T^{-1}u) F(T^{-1}S_ju)F(S_iu)}
           {F(u)F(T^{-1}S_iS_ju)} \notag \\
    &\quad = \frac{c_j(S_iu)}{c_j(Tu)}
      \frac{\alpha_{ij}(u) F(TS_ju)F(S_iu)
        - \alpha_{ji}(u) F(TS_iu)F(S_ju)}
           {F(Tu)F(S_iS_ju)}.
           \label{eq:HM-reduced}
\end{align}
 
\noindent\textit{$T$-orbit invariance.}
Define
\begin{align}
\mathcal{J}(u) \defeq
\frac{\alpha_{ij}(u) F(S_iu)F(TS_ju)
      - \alpha_{ji}(u) F(S_ju)F(TS_iu)}
     {\bigl(\alpha_{ij}(u) - \alpha_{ji}(u)\bigr)
      F(Tu)F(S_iS_ju)},
      \label{eq:J-def}
\end{align}
which is well-defined by the nondegeneracy hypothesis.
Since $T$ and $S_i$ commute, evaluating
$\mathcal{J}$ at $T^{-1}u$ identifies the second fraction on
the left-hand side of \eqref{eq:HM-reduced} as
$(\alpha_{ij}(T^{-1}u) -
\alpha_{ji}(T^{-1}u))\mathcal{J}(T^{-1}u)$, while
the second fraction on the right-hand side equals
$(\alpha_{ij}(u) -
\alpha_{ji}(u))\mathcal{J}(u)$.
Rearranging, equation~\eqref{eq:HM-reduced} therefore becomes
\begin{align}
    \frac{c_j(Tu)\lambda_i(S_ju)}{c_j(S_iu)\lambda_i(u)}
    \bigl(\alpha_{ij}(T^{-1}u) -
    \alpha_{ji}(T^{-1}u)\bigr)
    \mathcal{J}(T^{-1}u)
    = 
    \bigl(\alpha_{ij}(u) -
    \alpha_{ji}(u)\bigr)
    \mathcal{J}(u).
    \label{eq:HM-J-form}
\end{align}
The scalar mixed diamond equation
\eqref{eq:diamond-mixed-ij} for the bare data
$(c, \lambda)$, after substituting and rearranging,
gives
$c_j(Tu)\lambda_i(S_ju)/(c_j(S_iu)\lambda_i(u))
= (\alpha_{ij}(u) - \alpha_{ji}(u))
/(\alpha_{ij}(T^{-1}u) - \alpha_{ji}(T^{-1}u))$,
so the scalar prefactors on both sides of
\eqref{eq:HM-J-form} are equal. Hence
$\mathcal{J}(T^{-1}u) = \mathcal{J}(u)$ so 
 $\mathcal{J}$ is constant along $T$-orbits and the
boundary conditions give
\[
\mathcal{J}(u)
= \lim_{m \to -\infty} \mathcal{J}(T^m u)
= 1,
\]
since every $F$-factor in \eqref{eq:J-def} tends to $1$
and $\alpha_{ij}(T^m u)$ converges to a nondegenerate
limit. Setting $\mathcal{J}(u) = 1$ in \eqref{eq:J-def}
and clearing the denominator gives
\eqref{eq:HM-variable}.
\end{proof}

\section{Gauge Equivalence and Non-Abelian Hirota--Miwa}\label{sec:gauge-NAHM}
 
The diamond equations are not a new integrable system.
They are gauge equivalent to the non-abelian Hirota--Miwa
system~\cite{Nimmo2006}, the non-commutative
generalization of the classical Hirota--Miwa equation.
 
To state the equivalence, we first formulate the non-abelian
Hirota--Miwa system.  Write $T_0 \defeq T$ and
$T_i \defeq S_i$ for $i \in \N$, so that the lattice
$\mathcal{V}$ is generated symmetrically by the commuting
shifts $\{T_i\}_{i \in \Z_{\geq 0}}$.  Following
Nimmo~\cite{Nimmo2006} (see also
Gilson--Nimmo--Ohta~\cite{GilsonNimmoOhta2007}), consider
the linear system for
$\widetilde\Psi \colon \mathcal{V} \to \Hom(E, H)$,
\begin{equation}\label{eq:NAHM-linear}
    \widetilde\Psi(T_i u) - \widetilde\Psi(T_j u)
    + \widetilde\Psi(u)U_{ij}(u) = 0,
   \qquad i,j \in \Z_{\geq 0},
\end{equation}
where $U_{ij}(u) \in \End(E)$.
 
The system \eqref{eq:NAHM-linear} is compatible if and only
if the following three conditions hold for all
$i, j, k \in \Z_{\geq 0}$: the \emph{cycle relation}
\begin{equation}\label{eq:NAHM-cycle}
    U_{ij}(u) + U_{jk}(u) + U_{ki}(u) = 0,
\end{equation}
the \emph{shifted cycle relation}
\begin{equation}\label{eq:NAHM-shifted-cycle}
    U_{jk}(T_i u) + U_{ki}(T_j u) + U_{ij}(T_k u) = 0,
\end{equation}
and the \emph{quadratic relation}
\begin{equation}\label{eq:NAHM-quadratic}
    U_{ij}(u)U_{ik}(T_j u)
    = U_{ik}(u)U_{ij}(T_k u).
\end{equation}
In particular, $U_{ii} = 0$ and $U_{ji} = -U_{ij}$.  As observed by
Nimmo~\cite{Nimmo2006}, the shifted cycle relation is
redundant when the $U_{ij}$ are invertible:
\eqref{eq:NAHM-quadratic} and its cyclic permutations,
summed and factored using \eqref{eq:NAHM-cycle}, imply
\eqref{eq:NAHM-shifted-cycle}.  The system
\eqref{eq:NAHM-cycle}--\eqref{eq:NAHM-quadratic} is the
\emph{non-abelian Hirota--Miwa equation}.
 
\begin{remark}[Derivation]\label{rem:NAHM-derivation}
    The derivation is analogous to
    Proposition~\ref{prop:CL-gen-system}.  For 
    $i, j, k \in \Z_{\geq 0}$, summing \eqref{eq:NAHM-linear} over the pairs
    $(i,j)$, $(j,k)$, and $(k,i)$ gives \eqref{eq:NAHM-cycle}.
    For the remaining conditions, apply
    \eqref{eq:NAHM-linear} for the pair $(i,j)$ at the vertex
    $T_k u$; eliminate $\widetilde\Psi(T_i T_k u)$ and
    $\widetilde\Psi(T_j T_k u)$ using the shifted equations
    for the pairs $(j,k)$ and $(k,i)$; then eliminate
    $\widetilde\Psi(T_j u)$ and $\widetilde\Psi(T_k u)$ via
    \eqref{eq:NAHM-linear} at $u$.  This yields
    \[
    \widetilde\Psi(T_i u)\bigl(U_{jk}(T_i u) + U_{ki}(T_j u)
      + U_{ij}(T_k u)\bigr)
    + \widetilde\Psi(u)\bigl(U_{ik}(u)U_{ij}(T_k u)
      - U_{ij}(u)U_{ik}(T_j u)\bigr) = 0.
    \]
    Since $\widetilde\Psi(u)$ and $\widetilde\Psi(T_i u)$ may
    be prescribed independently, each bracketed coefficient
    must vanish, giving \eqref{eq:NAHM-shifted-cycle} and
    \eqref{eq:NAHM-quadratic}.  Conversely, if
    \eqref{eq:NAHM-cycle}--\eqref{eq:NAHM-quadratic} hold
    (together with \eqref{eq:NAHM-shifted-cycle} when the
    $U_{ij}$ are not assumed invertible), the displayed
    identity vanishes for all $\widetilde\Psi$, so the system
    is compatible.
\end{remark}
 
\begin{remark}
    In the scalar case, the ansatz
    \[ U_{ij} = \alpha_{ij}\frac{\tau(T_i u)\tau(T_j u)}{\tau(u)\tau(T_i T_j u)}\]
    recovers
    the classical Hirota--Miwa equation
    of~\cite{Hirota1981, Miwa1982}.
\end{remark}
 
\begin{remark}
    Nimmo~\cite{Nimmo2006} writes the linear problem as
    $\phi_{,i} - \phi_{,j} + U_{ij}\phi = 0$ in an abstract
    associative algebra, so the coefficients act on the left.
    Our linear problem \eqref{eq:NAHM-linear} has $\End(E)$
    acting on $\Hom(E,H)$ by right multiplication, which
    corresponds to working in the opposite algebra
    $\End(E)^{\mathrm{op}}$.  The cycle relation is unchanged;
    the quadratic relation \eqref{eq:NAHM-quadratic} has
    the product order reversed relative to
    Nimmo's~\cite[equation~(11)]{Nimmo2006}.
\end{remark}
 
\begin{proposition}\label{prop:diamond-NAHM}
    Assume that $C_i(u)$ is invertible for every $i \in \N$
    and $u \in \mathcal{V}$.  Then the diamond equations
    \eqref{eq:diamond-C-ij}--\eqref{eq:diamond-mixed-ij}
    are gauge equivalent to the non-abelian Hirota--Miwa
    system
    \eqref{eq:NAHM-cycle}--\eqref{eq:NAHM-quadratic}.
\end{proposition}

The proof uses the following auxiliary construction.

\begin{lemma}\label{lem:gauge-existence}
    Let $(C_i, \Lambda_i)_{i \in \N}$ satisfy the diamond
    equations
    \eqref{eq:diamond-C-ij}--\eqref{eq:diamond-mixed-ij},
    and suppose that $C_i(u)$ is invertible for every
    $i \in \N$ and $u \in \mathcal{V}$.  Then there exists
    $G \colon \mathcal{V} \to \mathrm{GL}(E)$ satisfying
    \begin{equation}\label{eq:gauge-system}
        C_i(u)G(S_i u) = G(Tu), \qquad i \in \N.
    \end{equation}
    The compatibility condition for \eqref{eq:gauge-system}
    is precisely \eqref{eq:diamond-C-ij}.
\end{lemma}

\begin{proof}
    Since $C_i(u)$ is invertible, \eqref{eq:gauge-system} is
    equivalent to
    \[G(S_i u) = C_i(u)^{-1}G(Tu).\]
    Choose $G$ arbitrarily on one $T$-orbit and extend by
    this recursion.  Computing $G(T^2 u)$ via $S_i S_j u$ and
    $S_j S_i u$ gives
    \[
    C_i(Tu)C_j(S_i u)G(S_j S_i u)
    = C_j(Tu)C_i(S_j u)G(S_i S_j u),
    \]
    which is consistent for all initial data if and only if
    \eqref{eq:diamond-C-ij} holds.
\end{proof}

\begin{proof}[Proof of Proposition~\ref{prop:diamond-NAHM}]
Fix a gauge $G$ as in Lemma~\ref{lem:gauge-existence}, and define
\begin{equation}\label{eq:U0i-def}
    \widetilde\Psi(u) \defeq \Psi(u)G(u), \qquad
    U_{0i}(u) \defeq -G(u)^{-1}\Lambda_i(u)G(S_i u),
    \qquad i \in \N.
\end{equation}
Right-multiplying the diamond linear problem
\eqref{eq:Linear-Psi-3} by $G(S_i u)$ and using
\eqref{eq:gauge-system}:
\begin{align*}
    \widetilde\Psi(S_i u)
    &= \Psi(S_i u)G(S_i u)
    = \bigl[\Psi(Tu)C_i(u)
      - \Psi(u)\Lambda_i(u)\bigr]G(S_i u) \\
    &= \widetilde\Psi(Tu) + \widetilde\Psi(u)U_{0i}(u).
\end{align*}
Setting $U_{00} \defeq 0$ and
\begin{equation}\label{eq:Uij-from-U0i}
    U_{ij}(u) \defeq U_{0j}(u) - U_{0i}(u),
    \qquad i, j \in \Z_{\geq 0},
\end{equation}
we obtain $U_{ji} = -U_{ij}$ and
\begin{equation}\label{eq:gauged-NAHM-linear}
    \widetilde\Psi(T_i u) - \widetilde\Psi(T_j u)
    + \widetilde\Psi(u)U_{ij}(u) = 0,
    \qquad i, j \in \Z_{\geq 0},
\end{equation}
which is the non-abelian Hirota--Miwa linear problem
\eqref{eq:NAHM-linear}.  The cycle relation
\eqref{eq:NAHM-cycle} holds automatically, since
$U_{ij} + U_{jk} + U_{ki}
= (U_{0j} - U_{0i}) + (U_{0k} - U_{0j}) + (U_{0i} - U_{0k})
= 0$.

Since the gauge
$\Psi \mapsto \widetilde\Psi = \Psi G$ is
pointwise invertible, the two linear problems are
compatible for exactly the same initial data.  By
Proposition~\ref{prop:CL-gen-system}, the compatibility
conditions of \eqref{eq:Linear-Psi-3} are the diamond
equations; by the analogous derivation
(Remark~\ref{rem:NAHM-derivation}), those
of \eqref{eq:gauged-NAHM-linear} are
\eqref{eq:NAHM-cycle}--\eqref{eq:NAHM-quadratic}.
\end{proof}

The following corollary makes the substitution explicit.

\begin{corollary}\label{cor:diamond-gauge-equiv-NAHM}
    Under the substitution
    \begin{equation}\label{eq:CL-from-GU}
        C_i(u) = G(Tu)G(S_i u)^{-1}, \qquad
        \Lambda_i(u) = -G(u)U_{0i}(u)G(S_i u)^{-1},
    \end{equation}
    the $C$-diamond \eqref{eq:diamond-C-ij} reduces to the
    compatibility condition for \eqref{eq:gauge-system}, while the
    $\Lambda$-diamond \eqref{eq:diamond-L-ij} and mixed
    diamond \eqref{eq:diamond-mixed-ij} become, respectively,
    \begin{align}
        U_{0i}(u)U_{0j}(S_i u)
        &= U_{0j}(u)U_{0i}(S_j u),
        \label{eq:U0-quadratic} \\
        U_{0j}(S_i u) + U_{0i}(Tu)
        &= U_{0i}(S_j u) + U_{0j}(Tu).
        \label{eq:U0-shifted-cycle}
    \end{align}
    These are the specializations of
    \eqref{eq:NAHM-quadratic}
    and \eqref{eq:NAHM-shifted-cycle} to the triple
    $(0, i, j)$.  Moreover, the definition
    $U_{ij} = U_{0j} - U_{0i}$ extends these to the full
    quadratic relation \eqref{eq:NAHM-quadratic} for all
    distinct triples in $\Z_{\geq 0}$.
\end{corollary}
 
\begin{proof}
\noindent\textit{Explicit form.}
Inverting \eqref{eq:gauge-system} and \eqref{eq:U0i-def}
gives \eqref{eq:CL-from-GU}.  Substituting into
\eqref{eq:diamond-C-ij}--\eqref{eq:diamond-mixed-ij}, the
outer $G$-factors cancel: \eqref{eq:diamond-C-ij} becomes
the compatibility condition for \eqref{eq:gauge-system}, and the remaining
two equations reduce to
\eqref{eq:U0-quadratic}--\eqref{eq:U0-shifted-cycle}.
 
\noindent\textit{Full quadratic relation.}
Using $U_{ij} = U_{0j} - U_{0i}$, we must show
\begin{equation}\label{eq:quadratic-general}
    U_{ij}(u)U_{ik}(T_j u)
    = U_{ik}(u)U_{ij}(T_k u)
\end{equation}
for all distinct $i, j, k \in \N$.  Expanding both sides and
subtracting, we apply \eqref{eq:U0-quadratic} three times:
\begin{align*}
    U_{0j}(u)U_{0k}(T_j u)
    &= U_{0k}(u)U_{0j}(T_k u), \\
    U_{0i}(u)U_{0k}(T_i u)
    &= U_{0k}(u)U_{0i}(T_k u), \\
    U_{0j}(u)U_{0i}(T_j u)
    &= U_{0i}(u)U_{0j}(T_i u).
\end{align*}
The first identity shows that
$U_{0j}(u)U_{0k}(T_j u)$ and $-U_{0k}(u)U_{0j}(T_k u)$
cancel.  In the remaining six terms, applying the third
identity to $-U_{0j}(u)U_{0i}(T_j u)$ and the second to
$U_{0k}(u)U_{0i}(T_k u)$ produces the common left factor
$U_{0i}(u)$:
\begin{align*}
    &U_{ij}(u)U_{ik}(T_j u) - U_{ik}(u)U_{ij}(T_k u) \\
    &\qquad = U_{0i}(u)\Bigl[
      \bigl(U_{0k}(T_i u) - U_{0k}(T_j u)\bigr)
    + \bigl(U_{0j}(T_k u) - U_{0j}(T_i u)\bigr) \\
    &\qquad\qquad\quad
    + \bigl(U_{0i}(T_j u) - U_{0i}(T_k u)\bigr)
    \Bigr].
\end{align*}
Define $V_{ab}(u) \defeq U_{0a}(T_b u) - U_{0a}(Tu)$ for
$a, b \in \N$.  The additive relation
\eqref{eq:U0-shifted-cycle}, rearranged as
$U_{0a}(T_b u) - U_{0a}(Tu)
= U_{0b}(T_a u) - U_{0b}(Tu)$,
gives $V_{ab} = V_{ba}$.  The bracket becomes
\[
(V_{ki} - V_{kj}) + (V_{jk} - V_{ji})
+ (V_{ij} - V_{ik}),
\]
and each pair cancels by the symmetry of $V$:
$V_{ki} = V_{ik}$, $V_{kj} = V_{jk}$, $V_{ji} = V_{ij}$.
This establishes \eqref{eq:quadratic-general}.
\end{proof}
 
\begin{remark}
    \label{rem:gauge-homogeneous}
    With $H = E$, the function
    $\Psi_G(u) \defeq G(u)^{-1} \in \End(E) = \Hom(E, E)$
    satisfies
    \[
    \Psi_G(S_i u) = \Psi_G(Tu)C_i(u), \qquad i \in \N,
    \]
    which is \eqref{eq:Linear-Psi-3} with
    $\Lambda_i \equiv 0$.  Thus the gauge is a distinguished
    nonvanishing homogeneous solution of the diamond linear
    problem.  In the KPZ applications of subsequent chapters,
    the asymmetric $(C_i, \Lambda_i)$-form arises directly
    from the kernels of each model; constructing
    $G$ is neither necessary nor natural in that setting.
\end{remark}

\begin{remark}
    On each finite truncation to the directions
    $T, S_1, \dotsc, S_N$, the diamond formulation uses $2N$
    endomorphism-valued fields $C_i, \Lambda_i$, whereas the
    symmetric Hirota--Miwa formulation uses
    $\binom{N+1}{2}$ antisymmetric fields $U_{ij}$.
\end{remark}

\begin{remark}\label{rem:Doliwa}
    A closely related arbitrary-gauge formulation appears in
    Doliwa's treatment of Desargues maps~\cite{Doliwa2009}.
    There the linear problem is written in terms of
    coefficients $A_{ij}$ with no distinguished lattice
    direction, and a gauge $G$ solving
    $G(T_i u)A_{ij}(u) = -G(T_j u)A_{ji}(u)$ brings it
    to the Hirota--Miwa form.
\end{remark}

\section{Solution Construction on Product Graphs}\label{sec:product-graph}
The Darboux theorem of \S\ref{sec:Darboux-transformations} turns any
dressing-compatible triple \((\Psi,\Phi,K)\) into a new
solution of the diamond equations, but constructing such
a triple for multipoint distributions is nontrivial. We
now give a systematic construction on product graphs. Starting from seed linear data on a base graph
\(\mathcal G\), we form \(m\) block-diagonal copies on
the product graph \(\mathcal G^m\), one for each layer.
A strictly lower-triangular family of propagators
\(B_u\) couples the layers, and the admissibility
conditions imposed on \(B_u\) are precisely those needed
to produce dressing-compatible data on
\(\mathcal G^m\). Once this is done, the Darboux theorem
applies directly.\footnote{In the KPZ applications, the base graph carries the one-point kernel data of a model, while the product graph encodes the joint shifts of $m$ observation points.}

\medskip

Fix a solution $(c_k, \lambda_k)$ of the
diamond equations
\eqref{eq:diamond-C-ij}--\eqref{eq:diamond-mixed-ij}, and
let $\psi \colon \mathcal{V} \to \Hom(E, H)$,
$\phi \colon \mathcal{V} \to \Hom(H, E)$ be solutions of
the corresponding linear problems:
\begin{align}
    \psi(S_k v) &= \psi(Tv)c_k(v) - \psi(v)\lambda_k(v),
    \label{eq:seed-psi} \\
    \phi(Tv) &= c_k(v)\phi(S_k v) - \lambda_k(Tv)
    \phi(TS_k v). \label{eq:seed-phi}
\end{align}

For $m \geq 1$, let $\mathcal{G}^m$ denote the directed graph
on $\mathcal{V}^m$ generated by the diagonal shifts
\[
\mathcal{T}u \defeq (Tu_1, \dotsc, Tu_m), \qquad
\mathcal{S}_k u \defeq (S_k u_1, \dotsc, S_k u_m),
\]
where $u = (u_1, \dotsc, u_m) \in \mathcal{V}^m$.

Set $E^m \defeq E \times \cdots \times E$ ($m$ copies). For
each layer $j \in \{1, \dotsc, m\}$, let
$\psi^{(j)}, \phi^{(j)}$ be solutions of the seed linear
relations \eqref{eq:seed-psi}--\eqref{eq:seed-phi} on
$\mathcal{G}$ with common data $(c_k, \lambda_k)$. For
$u \in \mathcal{V}^m$, we assemble these into block-diagonal
data on the product graph: define
\[
\psi(u)
= \bigl(\psi^{(1)}(u_1), \dotsc, \psi^{(m)}(u_m)\bigr)
\in \Hom(E^m, H),
\]
\[
\phi(u)
= \begin{pmatrix} \phi^{(1)}(u_1) \\ \vdots \\
\phi^{(m)}(u_m) \end{pmatrix}
\in \Hom(H, E^m),
\]
and the diagonal edge weights
\[
C_k(u) \defeq \mathrm{Diag}\bigl(c_k(u_1), \dotsc,
c_k(u_m)\bigr) \in \End(E^m),
\]
\[
D_k(u) \defeq \mathrm{Diag}\bigl(\lambda_k(u_1), \dotsc,
\lambda_k(u_m)\bigr) \in \End(E^m).
\]
The seed linear relations then lift to $\mathcal{G}^m$ in
matrix form:
\begin{align}
    \psi(\mathcal{S}_k u)
    &= \psi(\mathcal{T}u)C_k(u)
    - \psi(u)D_k(u),
    \label{eq:seed-psi-product} \\
    \phi(\mathcal{T}u)
    &= C_k(u)\phi(\mathcal{S}_k u)
    - D_k(\mathcal{T}u)
    \phi(\mathcal{T}\mathcal{S}_k u).
    \label{eq:seed-phi-product}
\end{align}
Since $(c_k, \lambda_k)$ satisfies the diamond equations on
$\mathcal{G}$ and $C_k, D_k$ are diagonal, the pair
$(C_k, D_k)$ satisfies the diamond equations
\eqref{eq:diamond-C-ij}--\eqref{eq:diamond-mixed-ij} on
$\mathcal{G}^m$. We use lowercase $\psi, \phi$ for this block-diagonal lift
throughout.

\begin{definition}\label{def:admissible-prop}
 For $u \in \mathcal{V}^m$, let
\(
B_u \colon \mathcal{V}^m \times \mathcal{V}^m \to \End(E^m)
\)
be a family of endomorphism-valued functions, strictly
lower-triangular in the layer indices:
$[B_u(x, y)]_{ij} = 0$ for $i \leq j$, with each entry
depending only on the corresponding components,
$[B_u(x, y)]_{ij} = [B_u]_{ij}(x_i, y_j)$.
We call the collection $\{B_u\}_{u \in \mathcal{V}^m}$
\emph{admissible propagators} if the following conditions
hold for all $u \in \mathcal{V}^m$, all
$x, y$ on the $\mathcal{T}$-orbit of~$u$, and all
$k \in \N$.
\begin{enumerate}
    \item $\mathcal{T}$-Covariance:
    \begin{align*}
        B_{\mathcal{T}u}(\mathcal{T}x, \mathcal{T}y)
        = B_u(x, y).
    \end{align*}
    \item $\mathcal{T}$-Splitting:
    \begin{align*}
        B_{\mathcal{T}u}(x, y) - B_u(x, y)
        = B_u(x, \mathcal{T}u)
          B_{\mathcal{T}u}(\mathcal{T}u, y).
    \end{align*}
     \item $\mathcal{S}_k$-Compatibility:
    \begin{align*}
        &C_k(x)B_{\mathcal{S}_k u}(\mathcal{S}_k x, \mathcal{S}_k y)
        - D_k(\mathcal{T}x)
          B_{\mathcal{S}_k u}(\mathcal{S}_k \mathcal{T}x, \mathcal{S}_k y)
          \nonumber \\
        &\qquad
        - B_u(\mathcal{T}x, \mathcal{T}y)C_k(y)
        + B_u(\mathcal{T}x, y)D_k(y) \\
        &\qquad
        = B_u(\mathcal{T}x, \mathcal{T}u)
          C_k(u)B_{\mathcal{S}_k u}(\mathcal{S}_k u, \mathcal{S}_k y).
    \end{align*}
    \item Regularity:
All sums in \eqref{eq:Psi-def}--\eqref{eq:K-def} are
well-defined, and the formal operations used in the proofs below,
namely reindexing, termwise multiplication by \(C_k\), \(D_k\),
and \(B_u(u,u)\), telescoping along negative
\(\mathcal T\)-orbits with vanishing of tail terms, and
distributing the product of two convergent sums,
are valid.\footnote{These conditions are verified for each model
in the subsequent chapters.}
\end{enumerate}
\end{definition}

\begin{remark}
    The structure of the $\mathcal{S}_k$-compatibility
    condition reflects the linear problems it is designed
    to intertwine. On the left-hand side, the first pair
    of terms has the form of the adjoint linear problem
    \eqref{eq:seed-phi-product} applied to the left
    argument of $B$, and the second pair has the form of
    the linear problem \eqref{eq:seed-psi-product} applied
    to the right argument. The three-factor product on the
    right-hand side mirrors the dressing compatibility
    structure of \eqref{eq:K-differences}.
\end{remark}
The admissible propagator conditions ensure that the
diagonal edge weights $D_k$ can be corrected to new
weights $\Lambda_k$ that, together with $C_k$, satisfy
the diamond equations on $\mathcal{G}^m$. Proofs of the lemmas in this subsection are deferred to \S\ref{sec:technical-proofs}.

\begin{lemma}\label{lem:Lambda-diamond}
    For $v \in \mathcal{V}^m$, set $P(v) \defeq I + B_v(v, v) \in \End(E^m)$.  Then $P(v)$ is invertible with $P(v)^{-1} = I - B_v(\mathcal{T}v, \mathcal{T}v)$, and the edge weights $(C_k, \Lambda_k)$ with
    \begin{equation}
        \Lambda_k(u) \defeq P(u) D_k(u) P(\mathcal{S}_k u)^{-1}
        \label{eq:Lambda-def}
    \end{equation}
    satisfy the diamond equations \eqref{eq:diamond-C-ij}--\eqref{eq:diamond-mixed-ij} on $\mathcal{G}^m$.
\end{lemma}

The propagators also correct the seed solutions $\psi, \phi$
into dressed solutions $\Psi, \Phi$ that satisfy the linear
problems on $\mathcal{G}^m$ with edge weights
$(C_k, \Lambda_k)$. Define
$\Psi \colon \mathcal{V}^m \to \Hom(E^m, H)$ and
$\Phi \colon \mathcal{V}^m \to \Hom(H, E^m)$ by
\begin{align}
    \Psi(u) &= \psi(u) + \sum_{p < 0}
    \psi(\mathcal{T}^{p} u)
    B_{\mathcal{T}^{-1}u}(\mathcal{T}^{p} u, u),
    \label{eq:Psi-def} \\
    \Phi(u) &= \phi(u) + \sum_{q \leq 0}
    B_u(u, \mathcal{T}^{q} u)
    \phi(\mathcal{T}^{q} u).
    \label{eq:Phi-def}
\end{align}
Note that $\Psi(u)\Phi(u) =
\sum_{j=1}^{m} \Psi_j(u)\Phi_j(u) \in \End(H)$.

\begin{lemma}\label{lem:Psi-Phi-linear}
The functions $\Psi, \Phi$ defined by \eqref{eq:Psi-def}--\eqref{eq:Phi-def} satisfy the linear problems \eqref{eq:Linear-Psi-3} and \eqref{eq:Linear-Phi-Prob} with edge weights $(C_k, \Lambda_k)$.
\end{lemma}

Finally, we construct a dressing-compatible kernel from the
seed data and propagators. Define
$K \colon \mathcal{V}^m \to \End(H)$ by
\begin{align}
    K(u) = \sum_{p, q \leq 0}
    \psi(\mathcal{T}^{p} u)
    \bigl[\delta_{p,q}I_{E^m}
    + B_u(\mathcal{T}^{p} u,
    \mathcal{T}^{q} u)\bigr]
    \phi(\mathcal{T}^{q} u).
    \label{eq:K-def}
\end{align}
Explicitly, in components,
\begin{align*}
    K(u) = \sum_{1 \leq j \leq i \leq m}
    \sum_{p, q \leq 0}
    \psi^{(i)}(T^p u_i)
    \bigl[\delta_{p,q}\delta_{i,j}
    + [B_u]_{ij}(T^p u_i, T^q u_j)\bigr]
    \phi^{(j)}(T^q u_j).
\end{align*}
The Kronecker $\delta_{i,j}$ contributes the diagonal $(i=j)$ and the strictly lower-triangular $[B_u]_{ij}$ contributes the strict lower triangle $(i> j)$.
\begin{lemma}\label{lem:K-factorization}
The kernel $K$ defined by \eqref{eq:K-def} satisfies
\begin{align*}
    K(\mathcal{T}u) - K(u)
    &= \Psi(\mathcal{T}u)\Phi(\mathcal{T}u), \quad 
    K(\mathcal{S}_k u) - K(u)
    = \Psi(\mathcal{T}u)C_k(u)\Phi(\mathcal{S}_k u).
\end{align*}
\end{lemma}

Combining Lemmas~\ref{lem:Lambda-diamond},
\ref{lem:Psi-Phi-linear}, and \ref{lem:K-factorization}
with Theorem~\ref{thm:Darboux},
we obtain the main result of this section.
\begin{theorem}\label{thm:product-graph}
Let $\{B_u\}$ be admissible propagators, let $(C_k, \Lambda_k, \Psi, \Phi, K)$ be as constructed above, and assume that there exists $z \in \mathbb{F}$ such that the resolvent $R(u) \defeq (I - zK(u))^{-1} \in \End(H)$ exists for every $u \in \mathcal{V}^m$.  Then $K$ is dressing compatible with $(\Psi, \Phi)$, and the pair $(\mathcal{M}_k, \Lambda_k)$, with
\[
\mathcal{M}(u) \defeq I + z\Phi(u)R(u)\Psi(u), \quad
\mathcal{M}_k(u) \defeq \mathcal{M}(\mathcal{T}u)^{-1}C_k(u)\mathcal{M}(\mathcal{S}_k u),
\]
satisfies the diamond equations \eqref{eq:diamond-C-ij}--\eqref{eq:diamond-mixed-ij} on $\mathcal{G}^m$.
\end{theorem}
\begin{remark}
    The product graph construction and the notion of
    admissible propagators are, to the author's
    knowledge, new.
\end{remark}

\subsection{Technical Construction Proofs}\label{sec:technical-proofs}

We begin with preliminary identities that streamline the main proofs.

\begin{lemma}\label{lem:absorption}
For any $v, x, y \in \mathcal{V}^m$,
\begin{enumerate}[label=(\roman*)]
\item\label{it:inverse} $P(v)$ is invertible with
    $P(v)^{-1} = I - B_v(\mathcal{T}v, \mathcal{T}v)$,
\item\label{it:left-absorb}
    $P(v)^{-1}
    B_{\mathcal{T}v}(\mathcal{T}v, y)
    = B_v(\mathcal{T}v, y)$,
\item\label{it:right-absorb}
    $B_{\mathcal{T}^{-1}v}(x, v)P(v)
    = B_v(x, v)$.
\end{enumerate}
\end{lemma}

\begin{proof}
Proof of \ref{it:inverse}:
$\mathcal{T}$-splitting at the vertex
$\mathcal{T}^{-1}v$ with $x = y = v$ gives
\[
B_v(v, v) - B_{\mathcal{T}^{-1}v}(v, v)
= B_{\mathcal{T}^{-1}v}(v, v)
  B_v(v, v).
\]
By $\mathcal{T}$-covariance,
$B_{\mathcal{T}^{-1}v}(v, v)
= B_v(\mathcal{T}v, \mathcal{T}v)$.  Rearranging,
\[
\bigl(I - B_v(\mathcal{T}v, \mathcal{T}v)\bigr)
\bigl(I + B_v(v, v)\bigr) = I.
\]
Since $B_v(v, v)$ is strictly lower-triangular in the
layer indices, it is nilpotent, and so
$P(v) = I + B_v(v, v)$ is invertible.  A left inverse
of an invertible endomorphism is its inverse.

\medskip

Proof of \ref{it:left-absorb}:
$\mathcal{T}$-splitting at $v$ with $x = \mathcal{T}v$
gives
\[
B_{\mathcal{T}v}(\mathcal{T}v, y)
- B_v(\mathcal{T}v, y)
= B_v(\mathcal{T}v, \mathcal{T}v)
  B_{\mathcal{T}v}(\mathcal{T}v, y).
\]
Rearranging,
$\bigl(I - B_v(\mathcal{T}v, \mathcal{T}v)\bigr)
B_{\mathcal{T}v}(\mathcal{T}v, y)
= B_v(\mathcal{T}v, y)$, and the left-hand factor
is $P(v)^{-1}$ by~\ref{it:inverse}.

\medskip

Proof of \ref{it:right-absorb}:
$\mathcal{T}$-splitting at $\mathcal{T}^{-1}v$ with
$y = v$ gives
\[
B_v(x, v) - B_{\mathcal{T}^{-1}v}(x, v)
= B_{\mathcal{T}^{-1}v}(x, v)
  B_v(v, v).
\]
Rearranging,
$B_{\mathcal{T}^{-1}v}(x, v)
\bigl(I + B_v(v, v)\bigr)
= B_v(x, v)$.
\end{proof}

In particular, since $P(\mathcal{S}_k u)^{-1}
= I - B_{\mathcal{S}_k u}(\mathcal{T}\mathcal{S}_k u,
\mathcal{T}\mathcal{S}_k u)$,  the definition of
$\Lambda_k$ may equivalently be written as
\begin{equation}\label{eq:Lambda-expanded}
    \Lambda_k(u) = \bigl(I + B_u(u, u)\bigr)D_k(u)
    \bigl(I - B_{\mathcal{S}_k u}(\mathcal{T}\mathcal{S}_k u,
    \mathcal{T}\mathcal{S}_k u)\bigr).
\end{equation}

\begin{lemma}\label{lem:compat-00}
    For all $u \in \mathcal{V}^m$ and $k \in \N$,
    \begin{equation}\label{eq:compat-00}
        P(u)^{-1}C_k(u)P(\mathcal{S}_k u)
        = C_k(u) + D_k(\mathcal{T}u)
          B_{\mathcal{S}_k u}(\mathcal{T}\mathcal{S}_k u,
          \mathcal{S}_k u)
        - B_u(\mathcal{T}u, u)D_k(u).
    \end{equation}
\end{lemma}

\begin{proof}
    Expanding the left-hand side using Lemma~\ref{lem:absorption}\ref{it:inverse} and the
    definition of $P$ gives
    \begin{align*}
    P(u)^{-1}C_k(u)P(\mathcal{S}_k u)
    &= \bigl(I - B_u(\mathcal{T}u, \mathcal{T}u)\bigr)
       C_k(u)
       \bigl(I + B_{\mathcal{S}_k u}(\mathcal{S}_k u,
       \mathcal{S}_k u)\bigr) \\
    &= C_k(u)
       + C_k(u)B_{\mathcal{S}_k u}(\mathcal{S}_k u,
         \mathcal{S}_k u)
       - B_u(\mathcal{T}u, \mathcal{T}u)C_k(u) \\
    &\quad
       - B_u(\mathcal{T}u, \mathcal{T}u)C_k(u)
         B_{\mathcal{S}_k u}(\mathcal{S}_k u,
         \mathcal{S}_k u).
    \end{align*}
    $\mathcal{S}_k$-compatibility at $x = y = u$ gives
    \begin{align*}
    C_k(u)&B_{\mathcal{S}_k u}(\mathcal{S}_k u,
    \mathcal{S}_k u)
    - D_k(\mathcal{T}u)
      B_{\mathcal{S}_k u}(\mathcal{T}\mathcal{S}_k u,
      \mathcal{S}_k u) \\
    &- B_u(\mathcal{T}u, \mathcal{T}u)C_k(u)
    + B_u(\mathcal{T}u, u)D_k(u) \\
    &= B_u(\mathcal{T}u, \mathcal{T}u)C_k(u)
      B_{\mathcal{S}_k u}(\mathcal{S}_k u,
      \mathcal{S}_k u).
    \end{align*}
    Rearranging for the three non-$C_k(u)$ terms in the
    expansion:
    \begin{align*}
    C_k(u)B_{\mathcal{S}_k u}(\mathcal{S}_k u,
    \mathcal{S}_k u)
    &- B_u(\mathcal{T}u, \mathcal{T}u)C_k(u)
    - B_u(\mathcal{T}u, \mathcal{T}u)C_k(u)
      B_{\mathcal{S}_k u}(\mathcal{S}_k u,
      \mathcal{S}_k u) \\
    &= D_k(\mathcal{T}u)
      B_{\mathcal{S}_k u}(\mathcal{T}\mathcal{S}_k u,
      \mathcal{S}_k u)
    - B_u(\mathcal{T}u, u)D_k(u).
    \end{align*}
    Substituting into the expansion yields
    \eqref{eq:compat-00}.
\end{proof}

\begin{proof}[Proof of Lemma~\ref{lem:Lambda-diamond}]
We verify that $(C_k, \Lambda_k)$ satisfies the diamond
equations
\eqref{eq:diamond-C-ij}--\eqref{eq:diamond-mixed-ij} on
$\mathcal{G}^m$.  Since $C_k$ is unchanged,
equation~\eqref{eq:diamond-C-ij} holds automatically.

\noindent\textit{$\Lambda$-equation \eqref{eq:diamond-L-ij}.}
Since $(C_k, D_k)$ satisfies the diamond equations on
$\mathcal{G}^m$, we have
$D_i(u)D_j(\mathcal{S}_i u) = D_j(u)D_i(\mathcal{S}_j u)$
for all~$u$.  The corrected weight $\Lambda_k(u)$ is a
conjugation of $D_k(u)$: indeed,
$\Lambda_k(u) = P(u)D_k(u)P(\mathcal{S}_k u)^{-1}$.
Therefore
\begin{align*}
    \Lambda_i(u)\Lambda_j(\mathcal{S}_i u)
    &= P(u)D_i(u)P(\mathcal{S}_i u)^{-1}
       P(\mathcal{S}_i u)D_j(\mathcal{S}_i u)
       P(\mathcal{S}_j \mathcal{S}_i u)^{-1} \\
    &= P(u)D_i(u)D_j(\mathcal{S}_i u)
       P(\mathcal{S}_i \mathcal{S}_j u)^{-1},
\end{align*}
where the intermediate factors cancel by
Lemma~\ref{lem:absorption}\ref{it:inverse}.
The expression with $i$ and $j$ swapped gives
$P(u)D_j(u)D_i(\mathcal{S}_j u)
P(\mathcal{S}_j \mathcal{S}_i u)^{-1}$.  Since
$D_i(u)D_j(\mathcal{S}_i u) =
D_j(u)D_i(\mathcal{S}_j u)$ and
$\mathcal{S}_i \mathcal{S}_j =
\mathcal{S}_j \mathcal{S}_i$, the two expressions are
equal, establishing \eqref{eq:diamond-L-ij}.

\noindent\textit{Mixed equation \eqref{eq:diamond-mixed-ij}.}
Fix $i, j \in \N$ and $u \in \mathcal{V}^m$.  We
must show
\begin{equation}\label{eq:mixed-CL-target-app}
    C_i(u)\Lambda_j(\mathcal{S}_i u)
    + \Lambda_i(\mathcal{T}u)C_j(\mathcal{S}_i u)
    - C_j(u)\Lambda_i(\mathcal{S}_j u)
    - \Lambda_j(\mathcal{T}u)C_i(\mathcal{S}_j u) = 0.
\end{equation}

\noindent\textit{Step~1: Reduction by conjugation.}
We write out the four terms using
$\Lambda_k(v) = P(v)D_k(v)P(\mathcal{S}_k v)^{-1}$,
together with
$P(\mathcal{T}u) = P(u)$ (by $\mathcal{T}$-covariance) and
$P(\mathcal{S}_k \mathcal{T}u) = P(\mathcal{S}_k u)$
($\mathcal{T}\mathcal{S}_k = \mathcal{S}_k\mathcal{T}$
and $\mathcal{T}$-covariance):
\begin{align*}
    C_i(u)\Lambda_j(\mathcal{S}_i u)
    &= C_i(u)P(\mathcal{S}_i u)
       D_j(\mathcal{S}_i u)
       P(\mathcal{S}_i\mathcal{S}_j u)^{-1}, \\[4pt]
    \Lambda_i(\mathcal{T}u)C_j(\mathcal{S}_i u)
    &= P(u)D_i(\mathcal{T}u)
       P(\mathcal{S}_i u)^{-1}
       C_j(\mathcal{S}_i u), \\[4pt]
    -C_j(u)\Lambda_i(\mathcal{S}_j u)
    &= -C_j(u)P(\mathcal{S}_j u)
       D_i(\mathcal{S}_j u)
       P(\mathcal{S}_i\mathcal{S}_j u)^{-1}, \\[4pt]
    -\Lambda_j(\mathcal{T}u)C_i(\mathcal{S}_j u)
    &= -P(u)D_j(\mathcal{T}u)
       P(\mathcal{S}_j u)^{-1}
       C_i(\mathcal{S}_j u).
\end{align*}
Since $P(u)$ and $P(\mathcal{S}_i\mathcal{S}_j u)$ are
invertible, \eqref{eq:mixed-CL-target-app} vanishes
if and only if the expression obtained by left-multiplying
by $P(u)^{-1}$ and right-multiplying by
$P(\mathcal{S}_i\mathcal{S}_j u)$ vanishes:
\begin{align}\label{eq:conjugated-target}
    &\phantom{{}+{}}
    P(u)^{-1}C_i(u)P(\mathcal{S}_i u)
    D_j(\mathcal{S}_i u)
    \nonumber\\
    &+ D_i(\mathcal{T}u)P(\mathcal{S}_i u)^{-1}
       C_j(\mathcal{S}_i u)P(\mathcal{S}_i\mathcal{S}_j u)
    \nonumber\\
    &- P(u)^{-1}C_j(u)P(\mathcal{S}_j u)
       D_i(\mathcal{S}_j u)
    \nonumber\\
    &- D_j(\mathcal{T}u)P(\mathcal{S}_j u)^{-1}
       C_i(\mathcal{S}_j u)P(\mathcal{S}_i\mathcal{S}_j u)
    = 0.
\end{align}

\noindent\textit{Step~2: Expansion via
Lemma~\ref{lem:compat-00}.}
Each of the four lines of \eqref{eq:conjugated-target}
contains a factor of the form
$P(\cdot)^{-1}C_k(\cdot)P(\cdot)$, which we expand
using \eqref{eq:compat-00}.  For lines~1 and~3, this is
\eqref{eq:compat-00} at the vertex~$u$ (with $k = i$ and
$k = j$ respectively):
\[
    P(u)^{-1}C_k(u)P(\mathcal{S}_k u)
    = C_k(u) + D_k(\mathcal{T}u)
      B_{\mathcal{S}_k u}(\mathcal{T}\mathcal{S}_k u,
      \mathcal{S}_k u)
    - B_u(\mathcal{T}u, u)D_k(u).
\]
For line~2, this is \eqref{eq:compat-00} at the vertex
$\mathcal{S}_i u$ (with $k = j$):
\begin{align*}
    P(\mathcal{S}_i u)^{-1}&C_j(\mathcal{S}_i u)
    P(\mathcal{S}_j\mathcal{S}_i u) \\
    &= C_j(\mathcal{S}_i u)
    + D_j(\mathcal{T}\mathcal{S}_i u)
      B_{\mathcal{S}_i\mathcal{S}_j u}
      (\mathcal{T}\mathcal{S}_i\mathcal{S}_j u,
      \mathcal{S}_i\mathcal{S}_j u) \\
    &\quad - B_{\mathcal{S}_i u}(\mathcal{T}\mathcal{S}_i u,
      \mathcal{S}_i u)D_j(\mathcal{S}_i u),
\end{align*}
and line~4 is the analogous identity at $\mathcal{S}_j u$
(with $k = i$), obtained by swapping
$i \leftrightarrow j$.

\noindent\textit{Step~3: Cancellation.}
Substituting these four expansions into
\eqref{eq:conjugated-target} and distributing, each line
produces three terms, for a total of twelve.  We number
them sequentially (three per line of
\eqref{eq:conjugated-target}):
\begin{align}
    &
    C_i(u)D_j(\mathcal{S}_i u)
    + D_i(\mathcal{T}u)
      B_{\mathcal{S}_i u}(\mathcal{T}\mathcal{S}_i u,
      \mathcal{S}_i u)D_j(\mathcal{S}_i u)
    - B_u(\mathcal{T}u, u)D_i(u)D_j(\mathcal{S}_i u)
    \nonumber\\
    &+ D_i(\mathcal{T}u)C_j(\mathcal{S}_i u)
    + D_i(\mathcal{T}u)D_j(\mathcal{T}\mathcal{S}_i u)
      B_{\mathcal{S}_i\mathcal{S}_j u}
      (\mathcal{T}\mathcal{S}_i\mathcal{S}_j u,
      \mathcal{S}_i\mathcal{S}_j u)
    \nonumber\\
    &- D_i(\mathcal{T}u)
      B_{\mathcal{S}_i u}(\mathcal{T}\mathcal{S}_i u,
      \mathcal{S}_i u)D_j(\mathcal{S}_i u)- C_j(u)D_i(\mathcal{S}_j u)
    \nonumber\\
    &- D_j(\mathcal{T}u)
      B_{\mathcal{S}_j u}(\mathcal{T}\mathcal{S}_j u,
      \mathcal{S}_j u)D_i(\mathcal{S}_j u)+ B_u(\mathcal{T}u, u)D_j(u)D_i(\mathcal{S}_j u)
     \nonumber\\
    &- D_j(\mathcal{T}u)C_i(\mathcal{S}_j u)- D_j(\mathcal{T}u)D_i(\mathcal{T}\mathcal{S}_j u)
      B_{\mathcal{S}_i\mathcal{S}_j u}
      (\mathcal{T}\mathcal{S}_i\mathcal{S}_j u,
      \mathcal{S}_i\mathcal{S}_j u) \nonumber\\
      &+ D_j(\mathcal{T}u)
      B_{\mathcal{S}_j u}(\mathcal{T}\mathcal{S}_j u,
      \mathcal{S}_j u)D_i(\mathcal{S}_j u).
    \label{eq:12-term}
\end{align}
Terms~2 and~6 cancel, as do terms~8 and~12.  The remaining
eight terms group into three families:
\begin{enumerate}[label=(\roman*)]
    \item \emph{Terms independent of $B$ (positions 1, 4,
    7, 10):}
    \[
        C_i(u)D_j(\mathcal{S}_i u)
        + D_i(\mathcal{T}u)C_j(\mathcal{S}_i u)
        - C_j(u)D_i(\mathcal{S}_j u)
        - D_j(\mathcal{T}u)C_i(\mathcal{S}_j u) = 0,
    \]
    by the mixed diamond equation
    \eqref{eq:diamond-mixed-ij} for $(C_k, D_k)$.
    \item \emph{Terms with left factor
    $B_u(\mathcal{T}u, u)$ (positions 3, 9):}
    \[
        B_u(\mathcal{T}u, u)
        \bigl(D_j(u)D_i(\mathcal{S}_j u)
        - D_i(u)D_j(\mathcal{S}_i u)\bigr) = 0,
    \]
    by \eqref{eq:diamond-L-ij} for $(C_k, D_k)$ at~$u$.
    \item \emph{Terms with right factor
    $B_{\mathcal{S}_i\mathcal{S}_j u}
    (\mathcal{T}\mathcal{S}_i\mathcal{S}_j u,
    \mathcal{S}_i\mathcal{S}_j u)$
    (positions 5, 11):}
    \[
        \bigl(D_i(\mathcal{T}u)
        D_j(\mathcal{T}\mathcal{S}_i u)
        - D_j(\mathcal{T}u)
        D_i(\mathcal{T}\mathcal{S}_j u)\bigr)
        B_{\mathcal{S}_i\mathcal{S}_j u}
        (\mathcal{T}\mathcal{S}_i\mathcal{S}_j u,
        \mathcal{S}_i\mathcal{S}_j u) = 0,
    \]
    by \eqref{eq:diamond-L-ij} for $(C_k, D_k)$ at
    $\mathcal{T}u$.
\end{enumerate}
This completes the proof. 
\end{proof}
\begin{proof}[Proof of Lemma~\ref{lem:Psi-Phi-linear}]
First, we show that
\[\Psi(\mathcal{S}_k u) = \Psi(\mathcal{T}u)C_k(u)
- \Psi(u)\Lambda_k(u).\]
Since
$\Lambda_k(u) = P(u)D_k(u)P(\mathcal{S}_k u)^{-1}$,
right-multiplying by $P(\mathcal{S}_k u)$ gives the
equivalent identity
\begin{equation}\label{eq:Psi-target}
    \Psi(\mathcal{S}_k u)P(\mathcal{S}_k u)
    - \Psi(\mathcal{T}u)C_k(u)P(\mathcal{S}_k u)
    + \Psi(u)P(u)D_k(u) = 0.
\end{equation}
We compute each of the three terms separately.

\noindent\textit{Term~1: $\Psi(\mathcal{S}_k u)P(\mathcal{S}_k u)$.}
From definition~\eqref{eq:Psi-def} at $\mathcal{S}_k u$:
\[
    \Psi(\mathcal{S}_k u)
    = \psi(\mathcal{S}_k u)
    + \sum_{p < 0} \psi(\mathcal{T}^p \mathcal{S}_k u)
      B_{\mathcal{T}^{-1}\mathcal{S}_k u}
      (\mathcal{T}^p \mathcal{S}_k u, \mathcal{S}_k u).
\]
Right-multiplying by
$P(\mathcal{S}_k u)$ and applying
Lemma~\ref{lem:absorption}\ref{it:right-absorb} at
$v = \mathcal{S}_k u$ to each term of the sum converts the propagator factors 
into
$B_{\mathcal{S}_k u}(\mathcal{T}^p \mathcal{S}_k u,
\mathcal{S}_k u)$.  After this replacement, the leading
$\psi(\mathcal{S}_k u)B_{\mathcal{S}_k u}(\mathcal{S}_k u,
\mathcal{S}_k u)$ may be absorbed as the $p = 0$
contribution of the resulting sum, giving
\begin{equation}\label{eq:Psi-T1}
    \Psi(\mathcal{S}_k u)P(\mathcal{S}_k u)
    = \psi(\mathcal{S}_k u)
    + \sum_{p \leq 0}
      \psi(\mathcal{T}^p \mathcal{S}_k u)
      B_{\mathcal{S}_k u}(\mathcal{T}^p \mathcal{S}_k u,
      \mathcal{S}_k u).
\end{equation}

\noindent\textit{Term~2: $\Psi(\mathcal{T}u)C_k(u)P(\mathcal{S}_k u)$.}
From definition~\eqref{eq:Psi-def} at $\mathcal{T}u$
(reindexing $p \to p - 1$):
\begin{align*}
\Psi(\mathcal{T}u) &= \psi(\mathcal{T}u) +
\sum_{p \leq 0} \psi(\mathcal{T}^p u)
B_u(\mathcal{T}^p u, \mathcal{T}u).
\end{align*}
Right-multiplying by $C_k(u)P(\mathcal{S}_k u)$ gives
a leading term
$\psi(\mathcal{T}u)C_k(u)P(\mathcal{S}_k u)$ plus a
sum whose summands involve the factor
$B_u(\mathcal{T}^p u, \mathcal{T}u)C_k(u)
P(\mathcal{S}_k u)$.  To evaluate this factor, we
rearrange $\mathcal{S}_k$-compatibility at
$(x, y) = (\mathcal{T}^{p-1}u, u)$: moving the
$B_u C_k$-term to the right-hand side and recognizing
$P(\mathcal{S}_k u)$ gives
\begin{align}
    &B_u(\mathcal{T}^p u, \mathcal{T}u)C_k(u)
    P(\mathcal{S}_k u)
    \nonumber \\
    &= C_k(\mathcal{T}^{p-1} u)
      B_{\mathcal{S}_k u}(\mathcal{T}^{p-1}
      \mathcal{S}_k u, \mathcal{S}_k u)
\nonumber \\
    &\qquad - D_k(\mathcal{T}^p u)
      B_{\mathcal{S}_k u}(\mathcal{T}^p
      \mathcal{S}_k u, \mathcal{S}_k u)
    + B_u(\mathcal{T}^p u, u)D_k(u).
    \label{eq:compat-b0-app}
\end{align}
Substituting \eqref{eq:compat-b0-app} into the sum
produces three sums over $p \leq 0$.  We handle each in
turn.

The $C_kB_{\mathcal{S}_k u}$-sum is reindexed via
$p \mapsto p + 1$, shifting the range to $p \leq -1$:
\[
\sum_{p \leq -1}
\psi(\mathcal{T}^{p+1} u)C_k(\mathcal{T}^p u)
B_{\mathcal{S}_k u}(\mathcal{T}^p \mathcal{S}_k u,
\mathcal{S}_k u).
\]
Expanding the leading term,
\[
\psi(\mathcal{T}u)C_k(u)P(\mathcal{S}_k u)
= \psi(\mathcal{T}u)C_k(u)
+ \psi(\mathcal{T}u)C_k(u)
B_{\mathcal{S}_k u}(\mathcal{S}_k u,\mathcal{S}_k u),
\]
the latter piece is the $p = 0$ contribution of the
reindexed sum, extending the range to $p \leq 0$.
Combining this extended sum with the
$D_kB_{\mathcal{S}_k u}$-sum via the seed linear
relation
$\psi(\mathcal{T}^{p+1}u)C_k(\mathcal{T}^p u)
- \psi(\mathcal{T}^p u)D_k(\mathcal{T}^p u)
= \psi(\mathcal{T}^p \mathcal{S}_k u)$
yields
\begin{equation}\label{eq:Psi-T2}
\begin{aligned}
    \Psi(\mathcal{T}u)C_k(u)P(\mathcal{S}_k u)
    &= \psi(\mathcal{T}u)C_k(u) \\
    &\quad + \sum_{p \leq 0}
      \psi(\mathcal{T}^p \mathcal{S}_k u)
      B_{\mathcal{S}_k u}(\mathcal{T}^p \mathcal{S}_k u,
      \mathcal{S}_k u) \\
    &\quad + \sum_{p \leq 0}
      \psi(\mathcal{T}^p u)
      B_u(\mathcal{T}^p u, u)D_k(u).
\end{aligned}
\end{equation}

\noindent\textit{Term~3: $\Psi(u)P(u)D_k(u)$.}
From definition~\eqref{eq:Psi-def} at $u$:
\[
    \Psi(u)
    = \psi(u)
    + \sum_{p < 0} \psi(\mathcal{T}^p u)
      B_{\mathcal{T}^{-1}u}
      (\mathcal{T}^p u, u).
\]
Right-multiplying by $P(u) = I + B_u(u, u)$ and
applying Lemma~\ref{lem:absorption}\ref{it:right-absorb} at
$v = u$ to each term of the sum converts
$B_{\mathcal{T}^{-1}u}(\mathcal{T}^p u, u)P(u)$
into $B_u(\mathcal{T}^p u, u)$.  After this replacement,
the leading
$\psi(u)B_u(u, u)$ may be absorbed as the $p = 0$
contribution of the resulting sum.
Right-multiplying by $D_k(u)$ then gives
\begin{equation}\label{eq:Psi-T3}
    \Psi(u)P(u)D_k(u)
    = \psi(u)D_k(u)
    + \sum_{p \leq 0}
      \psi(\mathcal{T}^p u)
      B_u(\mathcal{T}^p u, u)D_k(u).
\end{equation}

\noindent\textit{Combining the three terms.}
From \eqref{eq:Psi-T1}, \eqref{eq:Psi-T2}, and
\eqref{eq:Psi-T3}:
\begin{align*}
    \Psi(\mathcal{S}_k u)P(\mathcal{S}_k u)
    - \Psi(\mathcal{T}u)C_k(u)P(\mathcal{S}_k u)
    + \Psi(u)P(u)D_k(u)
    = \psi(\mathcal{S}_k u)
    - \psi(\mathcal{T}u)C_k(u)
    + \psi(u)D_k(u),
\end{align*}
where the $B_{\mathcal{S}_k u}$-sums from
\eqref{eq:Psi-T1} and \eqref{eq:Psi-T2} cancel, and the
$B_u \cdot D_k$-sums from \eqref{eq:Psi-T2} and
\eqref{eq:Psi-T3} cancel.  The remaining identity is the
seed linear relation \eqref{eq:seed-psi-product} at~$u$,
which vanishes. 

\medskip

Next, we verify the adjoint linear problem
for $\Phi$.  The argument mirrors the $\Psi$-verification: we isolate three terms from the target identity, reduce each using left absorption and $\mathcal{S}_k$-compatibility, and combine them using the seed adjoint linear relation.  We must show that
\[\Phi(\mathcal{T}u) = C_k(u)\Phi(\mathcal{S}_k u)
- \Lambda_k(\mathcal{T}u)\Phi(\mathcal{T}\mathcal{S}_k u).\]
By $\mathcal{T}$-covariance,
$P(\mathcal{T}u) = P(u)$ and
$P(\mathcal{S}_k \mathcal{T}u) = P(\mathcal{S}_k u)$
(using $\mathcal{T}\mathcal{S}_k = \mathcal{S}_k\mathcal{T}$),
so
$\Lambda_k(\mathcal{T}u)
= P(u)D_k(\mathcal{T}u)P(\mathcal{S}_k u)^{-1}$.
Left-multiplying by $P(u)^{-1}$ gives the equivalent
identity
\begin{equation}\label{eq:Phi-target}
    P(u)^{-1}\bigl[C_k(u)\Phi(\mathcal{S}_k u)
    - \Phi(\mathcal{T}u)\bigr]
    - D_k(\mathcal{T}u)P(\mathcal{S}_k u)^{-1}
    \Phi(\mathcal{T}\mathcal{S}_k u) = 0.
\end{equation}
We compute each of the three terms separately.

\noindent\textit{Term~1:
$P(u)^{-1}C_k(u)\Phi(\mathcal{S}_k u)$.}
From definition~\eqref{eq:Phi-def} at $\mathcal{S}_k u$:
\[
\Phi(\mathcal{S}_k u) = \phi(\mathcal{S}_k u)
+ \sum_{q \leq 0}
B_{\mathcal{S}_k u}(\mathcal{S}_k u,
\mathcal{T}^q \mathcal{S}_k u)
\phi(\mathcal{T}^q \mathcal{S}_k u).
\]
Left-multiplying by $P(u)^{-1}C_k(u)$ gives a
term $P(u)^{-1}C_k(u)\phi(\mathcal{S}_k u)$ plus a
sum whose summands involve the factor
$P(u)^{-1}C_k(u)
B_{\mathcal{S}_k u}(\mathcal{S}_k u,
\mathcal{T}^q \mathcal{S}_k u)$.  To evaluate this
factor, we rearrange $\mathcal{S}_k$-compatibility at
$(x, y) = (u, \mathcal{T}^q u)$: moving the
$C_kB_{\mathcal{S}_k u}$-term to the right-hand side
and recognizing
$P(u)^{-1}$
gives
\begin{align}
    P(u)^{-1}C_k(u)
    &B_{\mathcal{S}_k u}(\mathcal{S}_k u,
    \mathcal{T}^q \mathcal{S}_k u) \nonumber \\
    &= D_k(\mathcal{T}u)
      B_{\mathcal{S}_k u}(\mathcal{T}\mathcal{S}_k u,
      \mathcal{T}^q \mathcal{S}_k u)
\nonumber \\
    &\quad + B_u(\mathcal{T}u, \mathcal{T}^{q+1}u)
      C_k(\mathcal{T}^q u)
    - B_u(\mathcal{T}u, \mathcal{T}^q u)
      D_k(\mathcal{T}^q u).
    \label{eq:compat-0b-app}
\end{align}
Substituting \eqref{eq:compat-0b-app} into the sum
produces three sums over $q \leq 0$.  We handle each in
turn.

The leading term expands as
$P(u)^{-1}C_k(u)\phi(\mathcal{S}_k u)
= C_k(u)\phi(\mathcal{S}_k u)
- B_u(\mathcal{T}u, \mathcal{T}u)C_k(u)
\phi(\mathcal{S}_k u)$.
The $q = 0$ contribution from the
$B_uC_k$-sum is
$B_u(\mathcal{T}u, \mathcal{T}u)C_k(u)
\phi(\mathcal{S}_k u)$, which cancels the correction
above.

In the remaining $B_uC_k$-terms (now summing over
$q \leq -1$), reindexing $q \mapsto q + 1$ shifts the
range to $q \leq 0$:
\[
\sum_{q \leq 0}
B_u(\mathcal{T}u, \mathcal{T}^q u)
C_k(\mathcal{T}^{q-1} u)
\phi(\mathcal{T}^{q-1} \mathcal{S}_k u).
\]
Combining with the $-B_uD_k$-sum via the seed adjoint
linear relation
$C_k(\mathcal{T}^{q-1}u)
\phi(\mathcal{T}^{q-1}\mathcal{S}_k u)
- D_k(\mathcal{T}^q u)
\phi(\mathcal{T}^q \mathcal{S}_k u)
= \phi(\mathcal{T}^q u)$
yields
\begin{equation}\label{eq:Phi-T1}
\begin{aligned}
    P(u)^{-1}C_k(u)\Phi(\mathcal{S}_k u)
    &= C_k(u)\phi(\mathcal{S}_k u) \\
    &\quad + \sum_{q \leq 0}
      D_k(\mathcal{T}u)
      B_{\mathcal{S}_k u}(\mathcal{T}\mathcal{S}_k u,
      \mathcal{T}^q \mathcal{S}_k u)
      \phi(\mathcal{T}^q \mathcal{S}_k u) \\
    &\quad + \sum_{q \leq 0}
      B_u(\mathcal{T}u, \mathcal{T}^q u)
      \phi(\mathcal{T}^q u).
\end{aligned}
\end{equation}

\noindent\textit{Term~2: $P(u)^{-1}\Phi(\mathcal{T}u)$.}
From definition~\eqref{eq:Phi-def} at $\mathcal{T}u$
(reindexing $q \to q - 1$):
\[
\Phi(\mathcal{T}u) = \phi(\mathcal{T}u)
+ \sum_{q \leq 1}
  B_{\mathcal{T}u}(\mathcal{T}u, \mathcal{T}^q u)
  \phi(\mathcal{T}^q u).
\]
Left-multiplying by
$P(u)^{-1} = I - B_u(\mathcal{T}u, \mathcal{T}u)$ and
applying Lemma~\ref{lem:absorption}\ref{it:left-absorb} at
$v = u$ to each term of the sum converts
$P(u)^{-1}B_{\mathcal{T}u}(\mathcal{T}u,
\mathcal{T}^q u)$ into
$B_u(\mathcal{T}u, \mathcal{T}^q u)$.  The $q = 1$
term produces
$B_u(\mathcal{T}u, \mathcal{T}u)\phi(\mathcal{T}u)$,
which combines with the leading
$P(u)^{-1}\phi(\mathcal{T}u)
= \phi(\mathcal{T}u)
- B_u(\mathcal{T}u, \mathcal{T}u)\phi(\mathcal{T}u)$
to give~$\phi(\mathcal{T}u)$:
\begin{equation}\label{eq:Phi-T2}
    P(u)^{-1}\Phi(\mathcal{T}u)
    = \phi(\mathcal{T}u)
    + \sum_{q \leq 0}
      B_u(\mathcal{T}u, \mathcal{T}^q u)
      \phi(\mathcal{T}^q u).
\end{equation}

\noindent\textit{Term~3:
$D_k(\mathcal{T}u)P(\mathcal{S}_k u)^{-1}
\Phi(\mathcal{T}\mathcal{S}_k u)$.}
From definition~\eqref{eq:Phi-def} at
$\mathcal{T}\mathcal{S}_k u$ (reindexing $q \to q - 1$):
\[
\Phi(\mathcal{T}\mathcal{S}_k u)
= \phi(\mathcal{T}\mathcal{S}_k u)
+ \sum_{q \leq 1}
  B_{\mathcal{T}\mathcal{S}_k u}
  (\mathcal{T}\mathcal{S}_k u,
  \mathcal{T}^q \mathcal{S}_k u)
  \phi(\mathcal{T}^q \mathcal{S}_k u).
\]
Left-multiplying by
$P(\mathcal{S}_k u)^{-1}
= I - B_{\mathcal{S}_k u}(\mathcal{T}\mathcal{S}_k u,
\mathcal{T}\mathcal{S}_k u)$ and applying
Lemma~\ref{lem:absorption}\ref{it:left-absorb} (left absorption) at
$v = \mathcal{S}_k u$ to each term of the sum converts
$P(\mathcal{S}_k u)^{-1}
B_{\mathcal{T}\mathcal{S}_k u}
(\mathcal{T}\mathcal{S}_k u,
\mathcal{T}^q \mathcal{S}_k u)$ into
$B_{\mathcal{S}_k u}(\mathcal{T}\mathcal{S}_k u,
\mathcal{T}^q \mathcal{S}_k u)$.
The $q = 1$ term produces
$B_{\mathcal{S}_k u}(\mathcal{T}\mathcal{S}_k u,
\mathcal{T}\mathcal{S}_k u)
\phi(\mathcal{T}\mathcal{S}_k u)$, which combines with
\[
P(\mathcal{S}_k u)^{-1}
\phi(\mathcal{T}\mathcal{S}_k u)
= \phi(\mathcal{T}\mathcal{S}_k u)
- B_{\mathcal{S}_k u}(\mathcal{T}\mathcal{S}_k u,
\mathcal{T}\mathcal{S}_k u)
\phi(\mathcal{T}\mathcal{S}_k u)
\]
to give $\phi(\mathcal{T}\mathcal{S}_k u)$.
Left-multiplying by $D_k(\mathcal{T}u)$:
\begin{equation}\label{eq:Phi-T3}
\begin{aligned}
    D_k(\mathcal{T}u)P(\mathcal{S}_k u)^{-1}
    &\Phi(\mathcal{T}\mathcal{S}_k u)
    = D_k(\mathcal{T}u)\phi(\mathcal{T}\mathcal{S}_k u) \\
    &\quad + \sum_{q \leq 0}
      D_k(\mathcal{T}u)
      B_{\mathcal{S}_k u}(\mathcal{T}\mathcal{S}_k u,
      \mathcal{T}^q \mathcal{S}_k u)
      \phi(\mathcal{T}^q \mathcal{S}_k u).
\end{aligned}
\end{equation}

\noindent\textit{Combining the three terms.}
From \eqref{eq:Phi-T1}, \eqref{eq:Phi-T2},
and \eqref{eq:Phi-T3}:
the $D_kB_{\mathcal{S}_k u}$-sums from
\eqref{eq:Phi-T1} and \eqref{eq:Phi-T3}
cancel, and the $B_u$-sums from
\eqref{eq:Phi-T1} and \eqref{eq:Phi-T2} cancel,
leaving
\[
    C_k(u)\phi(\mathcal{S}_k u)
    - \phi(\mathcal{T}u)
    - D_k(\mathcal{T}u)\phi(\mathcal{T}\mathcal{S}_k u)
    = 0,
\]
which is the seed adjoint linear relation
\eqref{eq:seed-phi-product} at~$u$. 
\end{proof}
\begin{proof}[Proof of Lemma~\ref{lem:K-factorization}]

We verify the two dressing compatibility conditions:
\begin{align}
    K(\mathcal{T}u) - K(u)
    &= \Psi(\mathcal{T}u)\Phi(\mathcal{T}u),
    \label{eq:K-T-app} \\
    K(\mathcal{S}_k u) - K(u)
    &= \Psi(\mathcal{T}u)C_k(u)\Phi(\mathcal{S}_k u).
    \label{eq:K-Sk-app}
\end{align}

\noindent\textit{Step~1: Setup.}
From definition~\eqref{eq:K-def}, shifting
$u \mapsto \mathcal{T}u$ and reindexing $p \to p - 1$,
$q \to q - 1$ gives
\[
K(\mathcal{T}u) = \sum_{p, q \leq 1}
\psi(\mathcal{T}^p u)
\bigl[\delta_{p,q}I + B_{\mathcal{T}u}(\mathcal{T}^p u,
\mathcal{T}^q u)\bigr]\phi(\mathcal{T}^q u),
\]
where $\delta_{p,q}$ denotes the Kronecker delta in the
summation indices.
Therefore, splitting into diagonal and off-diagonal parts: 
\begin{align}
K(\mathcal{T}u) - K(u)
&= \sum_{p \leq 1}
   \psi(\mathcal{T}^p u)\phi(\mathcal{T}^p u)
   - \sum_{p \leq 0}
   \psi(\mathcal{T}^p u)\phi(\mathcal{T}^p u)
\nonumber \\
&\quad + \sum_{p, q \leq 1}
   \psi(\mathcal{T}^p u)
   B_{\mathcal{T}u}(\mathcal{T}^p u, \mathcal{T}^q u)
   \phi(\mathcal{T}^q u)
\nonumber \\
&\quad - \sum_{p, q \leq 0}
   \psi(\mathcal{T}^p u)
   B_u(\mathcal{T}^p u, \mathcal{T}^q u)
   \phi(\mathcal{T}^q u).
\label{eq:KT-raw-app}
\end{align}

\noindent\textit{Step~2: Diagonal part.}
Since both diagonal sums have identical summands
$\psi(\mathcal{T}^p u)\phi(\mathcal{T}^p u)$, their
difference reduces to the single term at $p = 1$, namely
$\psi(\mathcal{T}u)\phi(\mathcal{T}u)$.

\noindent\textit{Step~3: Off-diagonal part.}
Applying $\mathcal{T}$-splitting to the first off-diagonal sum replaces
$B_{\mathcal{T}u}(\mathcal{T}^p u, \mathcal{T}^q u)$
by $B_u(\mathcal{T}^p u, \mathcal{T}^q u) +
B_u(\mathcal{T}^p u, \mathcal{T}u)
B_{\mathcal{T}u}(\mathcal{T}u, \mathcal{T}^q u)$.
Substituting and decomposing the index set
$\{p, q \leq 1\}$ as the disjoint union
$\{p, q \leq 0\} \sqcup \{p \leq 0, q = 1\}
\sqcup \{p = 1, q \leq 1\}$, the $B_u$-contributions
over $\{p, q \leq 0\}$ cancel with the second off-diagonal sum in
\eqref{eq:KT-raw-app}.  The remaining $B_u$-contributions
are:
\begin{align}
    (\mathrm{i})\colon &\quad \sum_{p \leq 0}
    \psi(\mathcal{T}^p u)
    B_u(\mathcal{T}^p u, \mathcal{T}u)
    \phi(\mathcal{T}u),
    & &\text{from $\{p \leq 0, q = 1\}$,}
    \label{eq:KT-bdry-q1} \\
    (\mathrm{ii})\colon &\quad \sum_{q \leq 1}
    \psi(\mathcal{T}u)
    B_u(\mathcal{T}u, \mathcal{T}^q u)
    \phi(\mathcal{T}^q u),
    & &\text{from $\{p = 1, q \leq 1\}$.}
    \label{eq:KT-bdry-p1}
\end{align}
Next, we similarly decompose the $B_u \cdot B_{\mathcal{T}u}$-contributions. The contribution from
$\{p = 1, q \leq 1\}$ is:
\begin{equation}\label{eq:KT-BuBT-p1}
\begin{aligned}
    \sum_{q \leq 1}
    \psi(\mathcal{T}u)
    B_u(\mathcal{T}u, \mathcal{T}u) B_{\mathcal{T}u}(\mathcal{T}u, \mathcal{T}^q u)
    \phi(\mathcal{T}^q u).
\end{aligned}
\end{equation}
Combining \eqref{eq:KT-bdry-p1} and
\eqref{eq:KT-BuBT-p1}, and using that $B_{\mathcal{T}u}(\mathcal{T}u, \mathcal{T}^q u)$ is equal to the factor
$B_u(\mathcal{T}u, \mathcal{T}^q u)
+ B_u(\mathcal{T}u, \mathcal{T}u)
B_{\mathcal{T}u}(\mathcal{T}u, \mathcal{T}^q u)$ by $\mathcal{T}$-splitting at~$u$ with
$x = \mathcal{T}u$.  This gives
\begin{equation}\label{eq:KT-combined-p1}
    \sum_{q \leq 1}
    \psi(\mathcal{T}u)
    B_{\mathcal{T}u}(\mathcal{T}u, \mathcal{T}^q u)
    \phi(\mathcal{T}^q u).
\end{equation}
The remaining $B_u \cdot B_{\mathcal{T}u}$-contributions
come from $\{p, q \leq 0\}$ and
$\{p \leq 0, q = 1\}$, which together form $\{p \leq 0, q \leq 1\}$ and so give:
\begin{equation}\label{eq:KT-cross}
    \sum_{\substack{p \leq 0 \\ q \leq 1}}
    \psi(\mathcal{T}^p u)
    B_u(\mathcal{T}^p u, \mathcal{T}u)
    B_{\mathcal{T}u}(\mathcal{T}u, \mathcal{T}^q u)
    \phi(\mathcal{T}^q u).
\end{equation}

\noindent\textit{Step~4: Factorization.}
From definition~\eqref{eq:Psi-def} at $\mathcal{T}u$
(reindexing $p \to p - 1$) and
definition~\eqref{eq:Phi-def} at $\mathcal{T}u$
(reindexing $q \to q - 1$):
\begin{align*}
\Psi(\mathcal{T}u)
&= \psi(\mathcal{T}u)
   + \sum_{p \leq 0} \psi(\mathcal{T}^p u)
     B_u(\mathcal{T}^p u, \mathcal{T}u), \\
\Phi(\mathcal{T}u)
&= \phi(\mathcal{T}u)
   + \sum_{q \leq 1}
     B_{\mathcal{T}u}(\mathcal{T}u, \mathcal{T}^q u)
     \phi(\mathcal{T}^q u).
\end{align*}
Expanding the product
$\Psi(\mathcal{T}u)\Phi(\mathcal{T}u)$ gives four
contributions:
\begin{align*}
\Psi(\mathcal{T}u)\Phi(\mathcal{T}u)
&= \psi(\mathcal{T}u)\phi(\mathcal{T}u) \\
&\quad + \sum_{p \leq 0}
  \psi(\mathcal{T}^p u)
  B_u(\mathcal{T}^p u, \mathcal{T}u)
  \phi(\mathcal{T}u) \\
&\quad + \sum_{q \leq 1}
  \psi(\mathcal{T}u)
  B_{\mathcal{T}u}(\mathcal{T}u, \mathcal{T}^q u)
  \phi(\mathcal{T}^q u) \\
&\quad + \sum_{\substack{p \leq 0 \\ q \leq 1}}
  \psi(\mathcal{T}^p u)
  B_u(\mathcal{T}^p u, \mathcal{T}u)
  B_{\mathcal{T}u}(\mathcal{T}u, \mathcal{T}^q u)
  \phi(\mathcal{T}^q u).
\end{align*}
These four terms coincide with the diagonal contribution
from Step~2, \eqref{eq:KT-bdry-q1},
\eqref{eq:KT-combined-p1}, and \eqref{eq:KT-cross}
respectively.

\noindent\textit{Verification of \eqref{eq:K-Sk-app}.}
The proof splits $K(\mathcal{S}_k u) - K(u)$ into diagonal and off-diagonal contributions, reduces each using the seed linear relations and $\mathcal{S}_k$-compatibility, and reassembles the result as $\Psi(\mathcal{T}u)C_k(u)\Phi(\mathcal{S}_k u)$.

\noindent\textit{Step~1: Setup.}
From definition~\eqref{eq:K-def}, where $\delta_{p,q}$
denotes the Kronecker delta in the summation indices:
\begin{align*}
K(\mathcal{S}_k u)
&= \sum_{p, q \leq 0}
   \psi(\mathcal{T}^p \mathcal{S}_k u)
   \bigl[\delta_{p,q}I
   + B_{\mathcal{S}_k u}(\mathcal{T}^p \mathcal{S}_k u,
     \mathcal{T}^q \mathcal{S}_k u)\bigr]
   \phi(\mathcal{T}^q \mathcal{S}_k u), \\
K(u)
&= \sum_{p, q \leq 0}
   \psi(\mathcal{T}^p u)
   \bigl[\delta_{p,q}I
   + B_u(\mathcal{T}^p u, \mathcal{T}^q u)\bigr]
   \phi(\mathcal{T}^q u).
\end{align*}
Note that $B_{\mathcal{S}_k u}$ has arguments on the
$\mathcal{T}$-orbit of $\mathcal{S}_k u$, and $B_u$ has
arguments on the $\mathcal{T}$-orbit of~$u$.  We split into
diagonal and off-diagonal parts:
\begin{align}
K(\mathcal{S}_k u) - K(u)
&= \sum_{p \leq 0}
   \bigl[\psi(\mathcal{T}^p \mathcal{S}_k u)
   \phi(\mathcal{T}^p \mathcal{S}_k u)
   - \psi(\mathcal{T}^p u)\phi(\mathcal{T}^p u)
   \bigr]
\nonumber \\
&\quad + \sum_{p, q \leq 0}
   \bigl[\psi(\mathcal{T}^p \mathcal{S}_k u)
   B_{\mathcal{S}_k u}(\mathcal{T}^p \mathcal{S}_k u,
   \mathcal{T}^q \mathcal{S}_k u)
   \phi(\mathcal{T}^q \mathcal{S}_k u)
\nonumber \\
&\qquad\qquad
   - \psi(\mathcal{T}^p u)
   B_u(\mathcal{T}^p u, \mathcal{T}^q u)
   \phi(\mathcal{T}^q u)\bigr].
\label{eq:KSk-split-app}
\end{align}

\noindent\textit{Step~2: Diagonal part.}
Using the seed linear relations
\eqref{eq:seed-psi-product}--\eqref{eq:seed-phi-product}:
\begin{align*}
\psi(\mathcal{T}^p \mathcal{S}_k u)
&= \psi(\mathcal{T}^{p+1} u)C_k(\mathcal{T}^p u)
   - \psi(\mathcal{T}^p u)D_k(\mathcal{T}^p u), \\
\phi(\mathcal{T}^p u)
&= C_k(\mathcal{T}^{p-1} u)
   \phi(\mathcal{T}^{p-1} \mathcal{S}_k u)
   - D_k(\mathcal{T}^p u)
   \phi(\mathcal{T}^p \mathcal{S}_k u).
\end{align*}
Substituting the first relation into
$\psi(\mathcal{T}^p \mathcal{S}_k u)
\phi(\mathcal{T}^p \mathcal{S}_k u)$
and the second into
$\psi(\mathcal{T}^p u)\phi(\mathcal{T}^p u)$,
the $D_k$-cross-terms cancel, giving the summand 
\[
\psi(\mathcal{T}^{p+1} u)C_k(\mathcal{T}^p u)
\phi(\mathcal{T}^p \mathcal{S}_k u)
- \psi(\mathcal{T}^p u)C_k(\mathcal{T}^{p-1} u)
\phi(\mathcal{T}^{p-1} \mathcal{S}_k u).
\]
Summing over $p \leq 0$ telescopes to:
\begin{equation}\label{eq:KSk-diag-app}
\mathrm{diagonal}
= \psi(\mathcal{T}u)C_k(u)\phi(\mathcal{S}_k u).
\end{equation}

\noindent\textit{Step~3: Off-diagonal part.}
Substituting the seed linear relations for
$\psi(\mathcal{T}^p \mathcal{S}_k u)$ and
$\phi(\mathcal{T}^q u)$ into the off-diagonal sum produces
four sums:
\begin{align}
\tag{A}
 &\sum_{p,q \leq 0}
\psi(\mathcal{T}^{p+1}u)C_k(\mathcal{T}^p u)
B_{\mathcal{S}_k u}(\mathcal{T}^p \mathcal{S}_k u,
\mathcal{T}^q \mathcal{S}_k u)
\phi(\mathcal{T}^q \mathcal{S}_k u),
 \\
\tag{B}
 & -\sum_{p,q \leq 0}
\psi(\mathcal{T}^p u)D_k(\mathcal{T}^p u)
B_{\mathcal{S}_k u}(\mathcal{T}^p \mathcal{S}_k u,
\mathcal{T}^q \mathcal{S}_k u)
\phi(\mathcal{T}^q \mathcal{S}_k u),
 \\
\tag{C}
 & -\sum_{p,q \leq 0}
\psi(\mathcal{T}^p u)
B_u(\mathcal{T}^p u, \mathcal{T}^q u)
C_k(\mathcal{T}^{q-1}u)
\phi(\mathcal{T}^{q-1}\mathcal{S}_k u),
 \\
\tag{D}
 & \sum_{p,q \leq 0}
\psi(\mathcal{T}^p u)
B_u(\mathcal{T}^p u, \mathcal{T}^q u)
D_k(\mathcal{T}^q u)
\phi(\mathcal{T}^q \mathcal{S}_k u).
\end{align}

\noindent\textit{Step~4: Extracting boundary terms.}
In~(A), separating $p = 0$ and reindexing $p \to p - 1$ in
the remainder:
\begin{align}
(\mathrm{A}) &= \sum_{q \leq 0}
\psi(\mathcal{T}u)C_k(u)
B_{\mathcal{S}_k u}(\mathcal{S}_k u,
\mathcal{T}^q \mathcal{S}_k u)
\phi(\mathcal{T}^q \mathcal{S}_k u)
\tag{A-bdry} \\
&\quad + \sum_{p,q \leq 0}
\psi(\mathcal{T}^p u)C_k(\mathcal{T}^{p-1}u)
B_{\mathcal{S}_k u}(\mathcal{T}^{p-1}\mathcal{S}_k u,
\mathcal{T}^q \mathcal{S}_k u)
\phi(\mathcal{T}^q \mathcal{S}_k u).
\tag{A-int}
\end{align}
To reindex~(C), set $q' = q - 1$ so that the
$C_k$ and $\phi$ arguments align with those in~(A)
and~(D). The original range $q \leq 0$ becomes
$q' \leq -1$; extending to $q' \leq 0$ and separating
the boundary contribution at $q' = 0$ gives
(dropping the primes):
\begin{align}
(\mathrm{C}) &= +\sum_{p \leq 0}
\psi(\mathcal{T}^p u)
B_u(\mathcal{T}^p u, \mathcal{T}u)
C_k(u)\phi(\mathcal{S}_k u)
\tag{C-bdry} \\
&\quad - \sum_{p,q \leq 0}
\psi(\mathcal{T}^p u)
B_u(\mathcal{T}^p u, \mathcal{T}^{q+1}u)
C_k(\mathcal{T}^q u)
\phi(\mathcal{T}^q \mathcal{S}_k u).
\tag{C-int}
\end{align}

\noindent\textit{Step~5: Interior cancellation.}
Fix $p, q \leq 0$.  The interior contributions
(A-int)~+ (B) + (C-int) + (D) assemble into the summands
\begin{align*}
\psi(\mathcal{T}^p u)\Big[
&C_k(\mathcal{T}^{p-1}u)
 B_{\mathcal{S}_k u}(\mathcal{T}^{p-1}\mathcal{S}_k u,
 \mathcal{T}^q \mathcal{S}_k u)
- D_k(\mathcal{T}^p u)
  B_{\mathcal{S}_k u}(\mathcal{T}^p \mathcal{S}_k u,
  \mathcal{T}^q \mathcal{S}_k u) \\
&- B_u(\mathcal{T}^p u, \mathcal{T}^{q+1}u)
   C_k(\mathcal{T}^q u)
+ B_u(\mathcal{T}^p u, \mathcal{T}^q u)
  D_k(\mathcal{T}^q u)
\Big]\phi(\mathcal{T}^q \mathcal{S}_k u).
\end{align*}
The expression in the square brackets is the left-hand side of $\mathcal{S}_k$-compatibility
at $(x, y) = (\mathcal{T}^{p-1}u, \mathcal{T}^q u)$.
Hence
\begin{equation}\label{eq:KSk-interior-app}
\mathrm{interior} = \sum_{p,q \leq 0}
\psi(\mathcal{T}^p u)
B_u(\mathcal{T}^p u, \mathcal{T}u)C_k(u)
B_{\mathcal{S}_k u}(\mathcal{S}_k u,
\mathcal{T}^q \mathcal{S}_k u)
\phi(\mathcal{T}^q \mathcal{S}_k u).
\end{equation}

\noindent\textit{Step~6: Factorization.}
Combining \eqref{eq:KSk-diag-app}, (A-bdry), (C-bdry), and
\eqref{eq:KSk-interior-app}:
\begin{align*}
&K(\mathcal{S}_k u) - K(u)
 \\
& = \psi(\mathcal{T}u)C_k(u)\phi(\mathcal{S}_k u)+ \sum_{q \leq 0}
  \psi(\mathcal{T}u)C_k(u)
  B_{\mathcal{S}_k u}(\mathcal{S}_k u,
  \mathcal{T}^q \mathcal{S}_k u)
  \phi(\mathcal{T}^q \mathcal{S}_k u) \\
&\quad + \sum_{p \leq 0}
  \psi(\mathcal{T}^p u)
  B_u(\mathcal{T}^p u, \mathcal{T}u)C_k(u)
  \phi(\mathcal{S}_k u) \\
&\quad + \sum_{p,q \leq 0}
  \psi(\mathcal{T}^p u)
  B_u(\mathcal{T}^p u, \mathcal{T}u)C_k(u)
  B_{\mathcal{S}_k u}(\mathcal{S}_k u,
  \mathcal{T}^q \mathcal{S}_k u)
  \phi(\mathcal{T}^q \mathcal{S}_k u).
\end{align*}
This is the expansion of the product
$\Psi(\mathcal{T}u)C_k(u)\Phi(\mathcal{S}_k u)$, since
\[
\Psi(\mathcal{T}u)
= \psi(\mathcal{T}u)
  + \sum_{p \leq 0} \psi(\mathcal{T}^p u)
    B_u(\mathcal{T}^p u, \mathcal{T}u)
\]
(from definition~\eqref{eq:Psi-def} at $\mathcal{T}u$,
reindexing $p \to p - 1$), and
\[
\Phi(\mathcal{S}_k u)
= \phi(\mathcal{S}_k u)
  + \sum_{q \leq 0}
    B_{\mathcal{S}_k u}(\mathcal{S}_k u,
    \mathcal{T}^q \mathcal{S}_k u)
    \phi(\mathcal{T}^q \mathcal{S}_k u)
\]
(from definition~\eqref{eq:Phi-def} at $\mathcal{S}_k u$).
This establishes \eqref{eq:K-Sk-app} and completes the
proof. 
\end{proof}

%% file: chapters/3-the-semi-discrete-framework.tex
\chapter{The Semi-Discrete Framework}
\label{ch:the-semi-discrete-framework}

{
  \setlength{\parskip}{0pt}
}

\label{sec:sc-framework}

This chapter develops the semi-discrete analogue of
Chapter~\ref{ch:the-discrete-framework}, replacing one discrete shift
by a continuous derivative. The semi-discrete diamond
equations, which arise as compatibility conditions for the
resulting linear problem
(Proposition~\ref{prop:sc-diamond}), admit a Darboux
transformation (Theorem~\ref{thm:sc-Darboux}) that dresses
the edge weights while preserving the diamond structure.
The scalar reduction extracts bilinear
differential-difference equations for Fredholm determinants
(Proposition~\ref{prop:sc-scalar-HM}), and a dual lattice
presentation (\S\ref{sec:sc-dual}) reformulates the linear
problem in the alternative basis used by several of the
continuous-time models. The semi-discrete theory arises from the scaling
$S_1 = e^{\epsilon\partial_1}$, $\epsilon \to 0$, but all
results are stated and proved independently.

\section{The Semi-Discrete Linear Problem}
\label{sec:sc-linear}

Let $E$ and $H \neq 0$ be vector spaces over $\mathbb{F} = \R$ or $\C$,
let $\mathcal{V}'$ be a lattice generated by invertible, mutually
commuting shifts $T$ and $S_j$, $j \geq 2$, and set
$\mathcal{V} = I \times \mathcal{V}'$ for an open interval
$I \subseteq \R$.  The shifts act on the
$\mathcal{V}'$-coordinate, leaving the $I$-coordinate
unchanged; write $\partial_1$ for the derivative in the
continuous variable.  Assign endomorphisms
\begin{equation*}
    C_k(u) \in \End(E), \qquad \Lambda_k(u) \in \End(E), \qquad k \geq 1,
\end{equation*}
all differentiable in the continuous variable.  Consider the
overdetermined linear system for
$\Psi \colon \mathcal{V} \to \Hom(E, H)$:
\begin{align}
\partial_1\Psi(u)
  &= \Psi(Tu)C_1(u) - \Psi(u)\Lambda_1(u),
  \label{eq:sc-linear-1} \\
\Psi(S_j u)
  &= \Psi(Tu)C_j(u) - \Psi(u)\Lambda_j(u),
  \qquad j \geq 2.
  \label{eq:sc-linear-j}
\end{align}
Compatibility of the linear system
\eqref{eq:sc-linear-1}--\eqref{eq:sc-linear-j} forces the
following conditions on $(C_k, \Lambda_k)$.
\begin{proposition}\label{prop:sc-diamond}
The system
\eqref{eq:sc-linear-1}--\eqref{eq:sc-linear-j} is
compatible\footnote{That is, for any choice of initial data
$\{\Psi(T^\ell u_0)\}_{\ell \in \Z}$ along a $T$-orbit,
the continuous equation and each $S_j$-recurrence give the
same value for $\partial_1\Psi(S_j u)$, and any two discrete
recurrences give the same value for $\Psi(S_i S_j u)$.}
if and only if two families of conditions hold.  First, for
each $j \geq 2$ and $u \in \mathcal{V}$, the following
semi-discrete diamond equations hold:
\begin{align}
C_1(Tu)C_j(u)
  &= C_j(Tu)C_1(S_j u),
  \label{eq:sc-diamond-C} \\
\Lambda_1(u)\Lambda_j(u)
  &= \Lambda_j(u)\Lambda_1(S_j u)
  + \partial_1\Lambda_j(u),
  \label{eq:sc-diamond-L} \\
C_1(u)\Lambda_j(u)
  + \Lambda_1(Tu)C_j(u)
  &- C_j(u)\Lambda_1(S_j u)
  - \Lambda_j(Tu)C_1(S_j u)
  - \partial_1 C_j(u) = 0.
  \label{eq:sc-diamond-mixed}
\end{align}
Second, for each pair $i, j \geq 2$ and
$u \in \mathcal{V}$, the fully discrete diamond equations
\eqref{eq:diamond-C-ij}--\eqref{eq:diamond-mixed-ij} of
Proposition~\ref{prop:CL-gen-system} hold.
\end{proposition}

\begin{proof}
Fix $j \geq 2$ and $u \in \mathcal{V}$.  We compute $\partial_1\Psi(S_j u)$ in
two ways.  
Differentiating \eqref{eq:sc-linear-j} gives
\[
\partial_1\Psi(S_j u) = \partial_1\Psi(Tu) C_j(u) + \Psi(Tu)\partial_1 C_j(u) - \partial_1\Psi(u)  \Lambda_j(u) - \Psi(u)\partial_1\Lambda_j(u).
\]
Substituting \eqref{eq:sc-linear-1} at $Tu$ and at $u$ to eliminate
$\partial_1\Psi(Tu)$ and $\partial_1\Psi(u)$, then
collecting by $\Psi(T^2 u)$, $\Psi(Tu)$, $\Psi(u)$:
\begin{align*}
\partial_1\Psi(S_j u)
&= \Psi(T^2 u)\bigl[C_1(Tu)C_j(u)\bigr] \\
&\quad - \Psi(Tu)\bigl[
  \Lambda_1(Tu)C_j(u) - \partial_1 C_j(u)
  + C_1(u)\Lambda_j(u)\bigr] \\
&\quad + \Psi(u)\bigl[
  \Lambda_1(u)\Lambda_j(u)
  - \partial_1\Lambda_j(u)\bigr].
\end{align*}
On the other hand, evaluating \eqref{eq:sc-linear-1} at
$S_j u$ gives
\[
\partial_1\Psi(S_j u)
= \Psi(TS_j u)C_1(S_j u)
  - \Psi(S_j u)\Lambda_1(S_j u).
\]
Since $TS_j = S_j T$, the recurrence \eqref{eq:sc-linear-j}
expresses $\Psi(S_j Tu)$ and $\Psi(S_j u)$ in terms of
$\Psi(T^2 u)$, $\Psi(Tu)$, and $\Psi(u)$.  Substituting and
collecting:
\begin{align*}
\partial_1\Psi(S_j u)
&= \Psi(T^2 u)\bigl[C_j(Tu)C_1(S_j u)\bigr] \\
&\quad - \Psi(Tu)\bigl[
  \Lambda_j(Tu)C_1(S_j u)
  + C_j(u)\Lambda_1(S_j u)\bigr] \\
&\quad + \Psi(u)\bigl[
  \Lambda_j(u)\Lambda_1(S_j u)\bigr].
\end{align*}
Subtracting the second from the first:
\begin{align}
0
&= \Psi(T^2 u)\bigl[
  C_1(Tu)C_j(u) - C_j(Tu)C_1(S_j u)\bigr]
  \notag \\
&\quad - \Psi(Tu)\bigl[
  C_1(u)\Lambda_j(u)
  + \Lambda_1(Tu)C_j(u)
  - C_j(u)\Lambda_1(S_j u) \notag \\
&\qquad\qquad\quad
  - \Lambda_j(Tu)C_1(S_j u)
  - \partial_1 C_j(u)\bigr]
  \notag \\
&\quad + \Psi(u)\bigl[
  \Lambda_1(u)\Lambda_j(u)
  - \Lambda_j(u)\Lambda_1(S_j u)
  - \partial_1\Lambda_j(u)\bigr].
  \label{eq:sc-compat-identity}
\end{align}
If the system is compatible, this identity holds for all
choices of initial data.  Since $\Psi(u)$, $\Psi(Tu)$,
$\Psi(T^2 u)$ may be independently prescribed in
$\Hom(E, H)$, each bracketed coefficient must vanish, giving
\eqref{eq:sc-diamond-C}--\eqref{eq:sc-diamond-mixed}.  The
compatibility conditions between two discrete directions
$S_i$ and $S_j$, $i, j \geq 2$, are obtained by the same
calculation as Proposition~\ref{prop:CL-gen-system}.

Conversely, assume
\eqref{eq:sc-diamond-C}--\eqref{eq:sc-diamond-mixed} and
the fully discrete diamond equations for all $i, j \geq 2$.
Then \eqref{eq:sc-compat-identity} vanishes for every
$j \geq 2$ and every $u$, giving $\partial_1\Psi(S_j u)$
the same value by either computation, and
$\Psi(S_j S_i u) = \Psi(S_i S_j u)$ for $i, j \geq 2$ by
Proposition~\ref{prop:CL-gen-system}.  Hence the system is
compatible.
\end{proof}

We call
\eqref{eq:sc-diamond-C}--\eqref{eq:sc-diamond-mixed}
\emph{semi-discrete diamond equations}.

\begin{remark}\label{rem:sc-KPZ}
In KPZ applications, the continuous direction $\partial_1$
and the discrete shifts $S_j$, $T$ encode the physical
parameters of the model.  For example, in continuous-time
Push-TASEP, $\partial_1 = \partial_t$
is the time derivative, $S_2$ shifts the particle index,
and $T$ shifts the spatial threshold parameter.
\end{remark}

\begin{remark}\label{rem:sc-scaling}
The semi-discrete framework arises from the fully discrete
theory of Chapter~\ref{ch:the-discrete-framework} under the scaling
$S_1 = e^{\epsilon\partial_1}$.  Expand the discrete edge
weights in the $S_1$-direction as
\[
C_1^{(\epsilon)}(u) = \epsilon C_1(u) + O(\epsilon^2),
\qquad
\Lambda_1^{(\epsilon)}(u) = -I + \epsilon\Lambda_1(u)
+ O(\epsilon^2).
\]
The leading-order $O(\epsilon)$ terms in the fully discrete
diamond equations for the pair $(1, j)$ then recover the
semi-discrete diamond equations
\eqref{eq:sc-diamond-C}--\eqref{eq:sc-diamond-mixed}; the
equations for $i, j \geq 2$ are unchanged.  The same scaling
yields the semi-discrete adjoint problem, dressing conditions,
and Darboux theorem of
\S\ref{sec:sc-Darboux}--\S\ref{sec:sc-scalar-HM}.
\end{remark}

\section{Semi-Discrete Darboux Transformations}
\label{sec:sc-Darboux}

The Darboux transformation produces, from a solution
$(C_k, \Lambda_k)$ of the semi-discrete diamond equations together with wave
functions and a compatible kernel, new edge weights
$(\mathcal{M}_k, \Lambda_k)$ satisfying the same semi-discrete diamond equations
of Proposition~\ref{prop:sc-diamond}.  We now develop this
construction.

Fix a solution $(C_k, \Lambda_k)$ of the semi-discrete diamond equations
of Proposition~\ref{prop:sc-diamond}, and
let $\Psi \colon \mathcal{V} \to \Hom(E, H)$ be a solution of
the linear problem
\eqref{eq:sc-linear-1}--\eqref{eq:sc-linear-j}.  The adjoint
linear problem\footnote{Its compatibility conditions are again
the diamond equations of
Proposition~\ref{prop:sc-diamond}; the
verification is analogous.} for
$\Phi \colon \mathcal{V} \to \Hom(H, E)$ is
\begin{align}
\partial_1\Phi(Tu)
  &= \Lambda_1(Tu)\Phi(Tu) - C_1(u)\Phi(u),
  \label{eq:sc-adjoint-1} \\
\Phi(Tu)
  &= C_j(u)\Phi(S_j u) - \Lambda_j(Tu)\Phi(TS_j u),
  \qquad j \geq 2.
  \label{eq:sc-adjoint-j}
\end{align}
\begin{definition}\label{def:sc-dressing}
We say that $K \colon \mathcal{V} \to \End(H)$ is
\emph{dressing compatible} with $(\Psi, \Phi)$ if:
\begin{enumerate}[label=\arabic*., leftmargin=*]
\item For every $u \in \mathcal{V}$ and $j \geq 2$,
\begin{align}
K(Tu) - K(u) &= \Psi(Tu)\Phi(Tu),
  \label{eq:sc-K-T} \\
\partial_1 K(u) &= \Psi(Tu)C_1(u)\Phi(u),
  \label{eq:sc-K-d1} \\
K(S_j u) - K(u)
  &= \Psi(Tu)C_j(u)\Phi(S_j u).
  \label{eq:sc-K-Sj}
\end{align}
\item There exists $z \in \mathbb{F}$ such that the resolvent
exists at every vertex:
\[
R(u) \defeq (I - zK(u))^{-1} \in \End(H),
\qquad u \in \mathcal{V}.
\]
\end{enumerate}
\end{definition}

\begin{example}\label{ex:sc-simple-kernel}
Under suitable analytic assumptions, the kernel
\[
K(u) = \sum_{p \leq 0} \Psi(T^p u)\Phi(T^p u)
\]
is dressing compatible with $(\Psi, \Phi)$; the $T$- and
$S_j$-conditions follow by telescoping.  The
$\partial_1$-condition follows by differentiating termwise,
cancelling the $\Lambda_1$-contributions via the linear and
adjoint problems, and telescoping the remainder.
\end{example}
\begin{theorem}\label{thm:sc-Darboux}
Let $K$ be dressing compatible with $(\Psi, \Phi)$, with
resolvent $R(u) = (I - zK(u))^{-1}$.  Define the dressed
observable and dressed waves by
\[
\mathcal{M}(u) \defeq I + z\Phi(u)R(u)\Psi(u),
\]
\[
\widehat\Psi(u) \defeq R(u)\Psi(u), \qquad
\widehat\Phi(u) \defeq \Phi(u)R(T^{-1}u).
\]
Then $\mathcal{M}(u)$ is invertible for every
$u \in \mathcal{V}$.  The dressed waves
$\widehat\Psi, \widehat\Phi$ satisfy the linear problem
\eqref{eq:sc-linear-1}--\eqref{eq:sc-linear-j} and adjoint
linear problem
\eqref{eq:sc-adjoint-1}--\eqref{eq:sc-adjoint-j}
respectively, with data $(\mathcal{M}_k, \Lambda_k)$: for each $j \geq 2$,
\begin{align*}
\partial_1\widehat\Psi(u)
  &= \widehat\Psi(Tu)\mathcal{M}_1(u)
     - \widehat\Psi(u)\Lambda_1(u), \\
\widehat\Psi(S_j u)
  &= \widehat\Psi(Tu)\mathcal{M}_j(u)
     - \widehat\Psi(u)\Lambda_j(u),
\shortintertext{and}
\partial_1\widehat\Phi(Tu)
  &= \Lambda_1(Tu)\widehat\Phi(Tu)
     - \mathcal{M}_1(u)\widehat\Phi(u), \\
\widehat\Phi(Tu)
  &= \mathcal{M}_j(u)\widehat\Phi(S_j u)
     - \Lambda_j(Tu)\widehat\Phi(TS_j u),
\end{align*}
where the dressed edge weights are given by
\begin{equation}\label{eq:sc-dressed-C}
\mathcal{M}_1(u) \defeq
  \mathcal{M}(Tu)^{-1}C_1(u)\mathcal{M}(u),
\quad
\mathcal{M}_j(u) \defeq
  \mathcal{M}(Tu)^{-1}C_j(u)\mathcal{M}(S_j u),
\end{equation}
and $\widehat K(u) \defeq K(u)R(u)$ is dressing compatible
with $(\widehat\Psi, \widehat\Phi)$.
Moreover, the pair $(\mathcal{M}_k, \Lambda_k)$ satisfies the
semi-discrete diamond equations of Proposition~\ref{prop:sc-diamond}.
\end{theorem}

The proof of Theorem~\ref{thm:sc-Darboux} rests on the
following resolvent identities.  For the remainder of this
subsection, we suppress functional dependence on $z$.

\begin{lemma}\label{lem:sc-resolvent-ids}
For every $u \in \mathcal{V}$ and $j \geq 2$, the
following resolvent identities hold:
\begin{align}
R(Tu) - R(u)
  &= zR(u)\Psi(Tu)\Phi(Tu)R(Tu),
  \label{eq:sc-resolvent-T} \\
\partial_1 R(u)
  &= zR(u)\Psi(Tu)C_1(u)\Phi(u)R(u),
  \label{eq:sc-resolvent-d1} \\
R(S_j u) - R(u)
  &= zR(u)\Psi(Tu)C_j(u)\Phi(S_j u)R(S_j u).
  \label{eq:sc-resolvent-Sj}
\end{align}
\end{lemma}

\begin{proof}
For any $v, w \in \mathcal{V}$, since
$R(v)^{-1} = I - zK(v)$ we have
$R(v)^{-1} - R(w)^{-1} = z(K(w) - K(v))$, and therefore
\begin{equation}\label{eq:sc-resolvent-general}
R(w) - R(v) = R(v)\bigl(R(v)^{-1} - R(w)^{-1}\bigr)R(w)
= zR(v)\bigl(K(w) - K(v)\bigr)R(w).
\end{equation}
Substituting the dressing compatibility conditions
\eqref{eq:sc-K-T} and \eqref{eq:sc-K-Sj} gives
\eqref{eq:sc-resolvent-T} and \eqref{eq:sc-resolvent-Sj}.
For \eqref{eq:sc-resolvent-d1}, differentiating
$R(u) R^{-1}(u) = I$ gives
$\partial_1 R(u) = -R(u)(\partial_1 R^{-1}(u))R(u)
= zR(u)(\partial_1 K(u))R(u)$, and substituting \eqref{eq:sc-K-d1}
completes the proof.
\end{proof}

\begin{remark}\label{rem:sc-Darboux-via-wave}
One could alternatively derive the diamond equations for
$(\mathcal{M}_k, \Lambda_k)$ from the compatibility of the
dressed linear problem via the calculation of
Proposition~\ref{prop:sc-diamond}: the $\widehat\Psi(u)$-bracket
vanishes because $\Lambda$ is unchanged, and the
$\widehat\Psi(T^2 u)$-bracket vanishes by telescoping.  The
remaining $\widehat\Psi(Tu)$-bracket, however, requires
$\widehat\Psi(Tu)$ to be injective, which is not assumed here.
The direct proof below circumvents the injectivity requirement
entirely, using only resolvent existence.
\end{remark}

\begin{proof}
We prove the theorem in six steps: invertibility of
$\mathcal{M}$, the dressed linear problem, the adjoint
dressed linear problem, the explicit formulas for
$\mathcal{M}_k$, the dressed kernel, and the dressed
diamond equations.

\noindent\textit{Invertibility of $\mathcal{M}$.}
We show that
\begin{equation}\label{eq:sc-M-inverse}
\mathcal{M}(Tu)^{-1}
= I - z\Phi(Tu)R(u)\Psi(Tu).
\end{equation}
Collecting the $z$-linear term as
$z\Phi(Tu)(R(Tu) - R(u))\Psi(Tu)$ and using
\eqref{eq:sc-resolvent-T} to reduce the $z^2$-term,
\begin{align*}
&\bigl(I - z\Phi(Tu)R(u)\Psi(Tu)\bigr)
  \bigl(I + z\Phi(Tu)R(Tu)\Psi(Tu)\bigr) \\
&\qquad = I + z\Phi(Tu)\bigl(R(Tu) - R(u)\bigr)\Psi(Tu) \\
&\qquad\quad - z^2\Phi(Tu)R(u)\Psi(Tu)\Phi(Tu)R(Tu)\Psi(Tu) \\
&\qquad = I + z\Phi(Tu)\bigl(R(Tu) - R(u)\bigr)\Psi(Tu)
  - z\Phi(Tu)\bigl(R(Tu) - R(u)\bigr)\Psi(Tu) \\
&\qquad = I.
\end{align*}
The right inverse follows from the opposite factorization
\[
R(Tu) - R(u) = zR(Tu)\Psi(Tu)\Phi(Tu)R(u).
\]

\noindent\textit{The dressed linear problem.}
We begin by expanding $\widehat\Psi(Tu)$ using the
definition of $\mathcal{M}$:
\begin{align*}
\widehat\Psi(Tu)
&= R(Tu)\Psi(Tu)
 = R(u)\Psi(Tu)
   + \bigl(R(Tu) - R(u)\bigr)\Psi(Tu) \\
&= R(u)\Psi(Tu)
   + zR(u)\Psi(Tu)\Phi(Tu)R(Tu)\Psi(Tu) \\
&= R(u)\Psi(Tu)\mathcal{M}(Tu),
\end{align*}
where we used \eqref{eq:sc-resolvent-T}.  Therefore
\begin{equation}\label{eq:sc-intertwine-wave}
\widehat\Psi(Tu)\mathcal{M}(Tu)^{-1} = R(u)\Psi(Tu).
\end{equation}
 
Fix $j \geq 2$.  The verification that
$\widehat\Psi(S_j u) =
\widehat\Psi(Tu)\mathcal{M}_j(u) -
\widehat\Psi(u)\Lambda_j(u)$
is identical to Theorem~\ref{thm:Darboux}, using
\eqref{eq:sc-intertwine-wave} and
\eqref{eq:sc-resolvent-Sj}. For the continuous direction, we compute $\partial_1\widehat\Psi$
directly.  By the Leibniz rule, the $\partial_1$-resolvent
identity \eqref{eq:sc-resolvent-d1}, and the linear problem
\eqref{eq:sc-linear-1}:
\begin{align*}
\partial_1\widehat\Psi(u)
&= (\partial_1 R(u))\Psi(u) + R(u)(\partial_1\Psi(u)) \\
&= zR(u)\Psi(Tu)C_1(u)\Phi(u)R(u)\Psi(u) \\
&\quad + R(u)\bigl[\Psi(Tu)C_1(u) - \Psi(u)\Lambda_1(u)\bigr] \\
&= R(u)\Psi(Tu)C_1(u)
   \bigl[I + z\Phi(u)R(u)\Psi(u)\bigr]
   - \widehat\Psi(u)\Lambda_1(u) \\
&= R(u)\Psi(Tu)C_1(u)\mathcal{M}(u) - \widehat\Psi(u)\Lambda_1(u).
\end{align*}
Applying \eqref{eq:sc-intertwine-wave} and the definition
$\mathcal{M}_1(u) =
\mathcal{M}(Tu)^{-1}C_1(u)\mathcal{M}(u)$:
\[
\partial_1\widehat\Psi(u)
= \widehat\Psi(Tu)\mathcal{M}_1(u)
  - \widehat\Psi(u)\Lambda_1(u).
\]
 
\noindent\textit{The adjoint dressed linear problem.}
We begin by expanding
$\mathcal{M}(Tu)\widehat\Phi(Tu)$ using the definition of
$\mathcal{M}$:
\begin{align*}
\mathcal{M}(Tu)\widehat\Phi(Tu)
&= \Phi(Tu)R(u)
   + z\Phi(Tu)R(Tu)\Psi(Tu)\Phi(Tu)R(u) \\
&= \Phi(Tu)R(u)
   + \Phi(Tu)\bigl(R(Tu) - R(u)\bigr) \\
&= \Phi(Tu)R(Tu),
\end{align*}
where we used \eqref{eq:sc-resolvent-T} in the form
$R(Tu) - R(u) = zR(Tu)\Psi(Tu)\Phi(Tu)R(u)$.
Therefore
\begin{equation}\label{eq:sc-adjoint-intertwine}
\widehat\Phi(Tu)
= \mathcal{M}(Tu)^{-1}\Phi(Tu)R(Tu).
\end{equation}
 
Fix $j \geq 2$.  The verification that
$\widehat\Phi(Tu) =
\mathcal{M}_j(u)\widehat\Phi(S_j u) -
\Lambda_j(Tu)\widehat\Phi(TS_j u)$
is identical to Theorem~\ref{thm:Darboux}, using
\eqref{eq:sc-adjoint-intertwine} and
\eqref{eq:sc-resolvent-Sj}. For the continuous direction, we compute
$\partial_1\widehat\Phi(Tu)$ directly.  By the Leibniz rule,
the adjoint linear problem \eqref{eq:sc-adjoint-1}, and the
$\partial_1$-resolvent identity \eqref{eq:sc-resolvent-d1}:
\begin{align*}
\partial_1\widehat\Phi(Tu)
&= (\partial_1\Phi(Tu))R(u)
   + \Phi(Tu)(\partial_1 R(u))\\
&= \bigl[\Lambda_1(Tu)\Phi(Tu) - C_1(u)\Phi(u)\bigr]R(u) \\
&\quad + z\Phi(Tu)R(u)\Psi(Tu)C_1(u)\Phi(u)R(u) \\
&= \Lambda_1(Tu)\Phi(Tu)R(u)
   - \bigl[I - z\Phi(Tu)R(u)\Psi(Tu)\bigr]
     C_1(u)\Phi(u)R(u) \\
&= \Lambda_1(Tu)\widehat\Phi(Tu)
   - \mathcal{M}(Tu)^{-1}C_1(u)\Phi(u)R(u).
\end{align*}
Applying the adjoint intertwining
$\Phi(u)R(u) = \mathcal{M}(u)\widehat\Phi(u)$
(from \eqref{eq:sc-adjoint-intertwine} at $T^{-1}u$)
and the definition
$\mathcal{M}_1(u) =
\mathcal{M}(Tu)^{-1}C_1(u)\mathcal{M}(u)$:
\[
\partial_1\widehat\Phi(Tu)
= \Lambda_1(Tu)\widehat\Phi(Tu)
  - \mathcal{M}_1(u)\widehat\Phi(u).
\]

\noindent\textit{Explicit formulas for $\mathcal{M}_k$.}
Fix $j \geq 2$.  The derivation of $\mathcal{M}_j$ is
identical to Theorem~\ref{thm:Darboux}: substitute
\eqref{eq:sc-M-inverse} into
$\mathcal{M}_j(u) =
\mathcal{M}(Tu)^{-1}C_j(u)\mathcal{M}(S_j u)$,
apply \eqref{eq:sc-resolvent-Sj} to the $z^2$-term, and
substitute the linear problem \eqref{eq:sc-linear-j} and the
adjoint \eqref{eq:sc-adjoint-j}.  The result is
\begin{equation}\label{eq:sc-Mj-explicit}
\mathcal{M}_j(u)
= C_j(u)
  + z\Lambda_j(Tu)\Phi(TS_j u)R(S_j u)\Psi(S_j u)
  - z\Phi(Tu)R(u)\Psi(u)\Lambda_j(u).
\end{equation}

Substituting \eqref{eq:sc-M-inverse} into
$\mathcal{M}_1(u) =
\mathcal{M}(Tu)^{-1}C_1(u)\mathcal{M}(u)$:
\begin{align*}
\mathcal{M}_1(u)
&= C_1(u)
  + zC_1(u)\Phi(u)R(u)\Psi(u)
  - z\Phi(Tu)R(u)\Psi(Tu)C_1(u) \\
&\qquad
  - z^2\Phi(Tu)R(u)\Psi(Tu)C_1(u)
    \Phi(u)R(u)\Psi(u).
\end{align*}
By \eqref{eq:sc-resolvent-d1}, the $z^2$-term equals
$-z\Phi(Tu)[\partial_1 R(u)]\Psi(u)$, giving
\begin{equation}\label{eq:sc-M1-intermediate}
\begin{aligned}
\mathcal{M}_1(u)
&= C_1(u)
  + zC_1(u)\Phi(u)R(u)\Psi(u) \\
&\quad - z\Phi(Tu)R(u)\Psi(Tu)C_1(u)
  - z\Phi(Tu)[\partial_1 R(u)]\Psi(u).
\end{aligned}
\end{equation}
Rearranging the adjoint \eqref{eq:sc-adjoint-1} as
$C_1(u)\Phi(u) =
\Lambda_1(Tu)\Phi(Tu) - \partial_1\Phi(Tu)$
and the linear problem \eqref{eq:sc-linear-1} as
$\Psi(Tu)C_1(u) =
\partial_1\Psi(u) + \Psi(u)\Lambda_1(u)$,
we substitute into \eqref{eq:sc-M1-intermediate}.  By the
Leibniz rule, the three derivative terms combine to
$-\partial_1[z\Phi(Tu)R(u)\Psi(u)]$, yielding
\begin{equation}\label{eq:sc-M1-explicit}
\begin{aligned}
\mathcal{M}_1(u)
&= C_1(u)
  + z\Lambda_1(Tu)\Phi(Tu)R(u)\Psi(u) \\
&\quad - z\Phi(Tu)R(u)\Psi(u)\Lambda_1(u)
  - \partial_1\bigl[z\Phi(Tu)R(u)\Psi(u)\bigr].
\end{aligned}
\end{equation}

\noindent\textit{The dressed kernel.}
Since $K(u)$ and
$R(u)$ commute at each vertex, $\widehat K(u) = R(u)K(u)$.
The identity
$(I - \eta\widehat K(u))(I - zK(u)) = I - (z+\eta)K(u)$
shows that $(I - \eta\widehat K(u))^{-1}$ exists whenever
$(I - (z+\eta)K(u))^{-1}$ does; in particular at
$\eta = -z$, where the condition is vacuous.  When $z = 0$,
$R = I$ and $\widehat K = K$, so the kernel identities
reduce to those of $K$ itself.  For $z \neq 0$, it remains
to verify the three kernel identities.
The $T$- and $S_j$-conditions are identical to the
corresponding steps in Theorem~\ref{thm:Darboux}, using
\eqref{eq:sc-resolvent-general},
\eqref{eq:sc-intertwine-wave}, and
\eqref{eq:sc-adjoint-intertwine} in place of their
discrete counterparts.

\noindent\textit{The $\partial_1$-condition.}
Differentiating $\widehat K(u) = z^{-1}(R(u) - I)$
and substituting \eqref{eq:sc-resolvent-d1}:
\[
\partial_1\widehat K(u) = z^{-1}\partial_1 R(u)
= R(u)\Psi(Tu)C_1(u)\Phi(u)R(u).
\]
Applying \eqref{eq:sc-intertwine-wave} to the left
factor and \eqref{eq:sc-adjoint-intertwine} at
$T^{-1}u$ to the right:
\begin{align*}
R(u)\Psi(Tu)C_1(u)\Phi(u)R(u)
= \widehat\Psi(Tu)\mathcal{M}(Tu)^{-1}
   C_1(u)\mathcal{M}(u)\widehat\Phi(u)
= \widehat\Psi(Tu)\mathcal{M}_1(u)\widehat\Phi(u).
\end{align*}

\noindent\textit{The dressed diamond equations.}
Since $\Lambda$ is unchanged,
\eqref{eq:sc-diamond-L} and the fully discrete
$\Lambda$-diamond of Proposition~\ref{prop:CL-gen-system} hold
for $(\mathcal{M}_k, \Lambda_k)$.  The $C$-diamonds
\eqref{eq:sc-diamond-C} and the fully discrete $C$-diamond
are preserved under dressing by the same telescoping argument
as in Theorem~\ref{thm:Darboux}, and the fully discrete mixed
diamonds for $i, j \geq 2$ follow by the same calculation,
since $\mathcal{M}_i, \mathcal{M}_j$ with $i, j \geq 2$ have
the same explicit formulas as in the fully discrete case.  It
remains to verify the semi-discrete mixed diamond equation for
each $j \geq 2$:
\begin{equation}\label{eq:sc-mixed-dressed}
\mathcal{M}_1(u)\Lambda_j(u)
+ \Lambda_1(Tu)\mathcal{M}_j(u)
- \mathcal{M}_j(u)\Lambda_1(S_j u)
- \Lambda_j(Tu)\mathcal{M}_1(S_j u)
-\partial_1\mathcal{M}_j(u) = 0.
\end{equation}

Fix $j \geq 2$.  To organize the cancellation, set
$Q(v) \defeq z\Phi(Tv)R(v)\Psi(v) \in \End(E)$,
so that the explicit formulas
\eqref{eq:sc-Mj-explicit}--\eqref{eq:sc-M1-explicit} read
\begin{align}
\mathcal{M}_j(u)
  &= C_j(u) + \Lambda_j(Tu)Q(S_j u)
     - Q(u)\Lambda_j(u),
  \label{eq:sc-Mj-A} \\
\mathcal{M}_1(u)
  &= C_1(u) + \Lambda_1(Tu)Q(u)
     - Q(u)\Lambda_1(u)
     - \partial_1Q(u).
  \label{eq:sc-M1-A}
\end{align}
Substituting \eqref{eq:sc-Mj-A}--\eqref{eq:sc-M1-A} (and the
formula for $\mathcal{M}_1(S_j u)$ obtained by evaluating
\eqref{eq:sc-M1-A} at $S_j u$) into each of the five terms of
\eqref{eq:sc-mixed-dressed}, the contributions independent of
$Q$ are
$C_1(u)\Lambda_j(u) + \Lambda_1(Tu)C_j(u)
- C_j(u)\Lambda_1(S_j u) - \Lambda_j(Tu)C_1(S_j u)
- \partial_1 C_j(u)$,
which vanishes by the bare mixed diamond
\eqref{eq:sc-diamond-mixed}.  The
remaining terms, all linear in $Q$, are:
\begin{gather*}
\Lambda_1(Tu)Q(u)\Lambda_j(u)
- Q(u)\Lambda_1(u)\Lambda_j(u)
- (\partial_1Q(u))\Lambda_j(u),
\tag*{(A)} \\
\Lambda_1(Tu)\Lambda_j(Tu)Q(S_j u)
- \Lambda_1(Tu)Q(u)\Lambda_j(u),
\tag*{(B)} \\
{-\Lambda_j(Tu)Q(S_j u)\Lambda_1(S_j u)}
+ Q(u)\Lambda_j(u)\Lambda_1(S_j u),
\tag*{(C)} \\
{-\Lambda_j(Tu)\Lambda_1(TS_j u)Q(S_j u)}
+ \Lambda_j(Tu)Q(S_j u)\Lambda_1(S_j u) + \Lambda_j(Tu)\partial_1Q(S_j u),
\tag*{(D)} \\
{-(\partial_1\Lambda_j(Tu))Q(S_j u)}
- \Lambda_j(Tu)\partial_1Q(S_j u) + (\partial_1Q(u))\Lambda_j(u)
+ Q(u)\partial_1\Lambda_j(u).
\tag*{(E)}
\end{gather*}
The first term of~(A) cancels with the second term of~(B),
the third term of~(A) with the third term of~(E),
the first term of~(C) with the second term of~(D),
and the third term of~(D) with the second term of~(E).
The surviving six terms group as
\begin{align*}
&Q(u)\bigl[
  -\Lambda_1(u)\Lambda_j(u)
  + \Lambda_j(u)\Lambda_1(S_j u)
  + \partial_1\Lambda_j(u)
\bigr] \\
&\quad + \bigl[
  \Lambda_1(Tu)\Lambda_j(Tu)
  - \Lambda_j(Tu)\Lambda_1(TS_j u)
  - \partial_1\Lambda_j(Tu)
\bigr]Q(S_j u).
\end{align*}
Both brackets vanish by \eqref{eq:sc-diamond-L} at $u$ and $Tu$, respectively. This proves \eqref{eq:sc-mixed-dressed}.
\end{proof}

\begin{remark}\label{rem:sc-domain-restriction}
In applications, the resolvent may exist only on a subset
$\widetilde{\mathcal{V}} \subseteq \mathcal{V}$.  Each
identity at a vertex~$u$ holds whenever the resolvent exists
at the finitely many shifts of~$u$ appearing in the proof.
For the dressed semi-discrete mixed diamond equation
\eqref{eq:sc-mixed-dressed} at a vertex~$u$, the explicit
formulas \eqref{eq:sc-Mj-A}--\eqref{eq:sc-M1-A}, which use
\eqref{eq:sc-M-inverse} to avoid the resolvent at
$T$-shifted vertices, require the resolvent only at $u$
and~$S_j u$.  In particular, $\widetilde{\mathcal{V}}$ need
not be closed under~$T$.
\end{remark}

\subsection{Fredholm determinants and the scalar reduction}
\label{sec:sc-scalar-HM}

In the probabilistic applications of subsequent chapters,
the Fredholm determinants
$F(u) = \det_H(I - zK(u))$ encode the distributional data of
each model.\footnote{In these applications $H$ is typically an
$L^2$ space, so $F(u)$ is obtained as a Fredholm determinant in
the classical sense.}  The corollary below establishes the
fundamental relation between $F$ and $\mathcal{M}$.

\begin{corollary}\label{cor:sc-Woodbury}
Let $K$ be dressing compatible with $(\Psi, \Phi)$, and
suppose that $\dim E < \infty$ and that
$F(u) \defeq \det_H(I - zK(u))$ is well-defined for every
$u \in \mathcal{V}$.  Then
\begin{equation}\label{eq:sc-Woodbury-det}
F(u) = F(Tu)\det_E\bigl(\mathcal{M}(Tu)\bigr).
\end{equation}
In particular, for each $j \geq 2$,
\begin{align}
\det_E(\mathcal{M}_j(u))F(u)F(S_j u)
  &= \det_E(C_j(u))F(Tu)F(T^{-1}S_j u),
  \label{eq:sc-det-dressed-j} \\
\det_E(\mathcal{M}_1(u))F(u)^2
  &= \det_E(C_1(u))F(Tu)F(T^{-1}u).
  \label{eq:sc-det-dressed-1}
\end{align}
\end{corollary}

\begin{proof}
The proof of \eqref{eq:sc-Woodbury-det} is identical to the
fully discrete case (Corollary~\ref{cor:Woodbury}): it depends
only on the $T$-dressing compatibility
$K(Tu) - K(u) = \Psi(Tu)\Phi(Tu)$, which is unchanged.
Identities \eqref{eq:sc-det-dressed-j}
and \eqref{eq:sc-det-dressed-1} follow by taking $\det_E$ of
$\mathcal{M}_j(u) =
\mathcal{M}(Tu)^{-1}C_j(u)\mathcal{M}(S_j u)$ and
$\mathcal{M}_1(u) =
\mathcal{M}(Tu)^{-1}C_1(u)\mathcal{M}(u)$ respectively,
applying \eqref{eq:sc-Woodbury-det}, and rearranging.
\end{proof}

When $E = \mathbb{F}$, the edge weights
$C_k, \Lambda_k, \mathcal{M}_k$ reduce to scalar-valued
functions $c_k, \lambda_k, \mathcal{M}_k$, and
Corollary~\ref{cor:sc-Woodbury} gives
\begin{align}
\mathcal{M}_j(u)
  &= c_j(u)
     \frac{F(Tu)F(T^{-1}S_j u)}{F(u)F(S_j u)},
     \qquad j \geq 2,
  \label{eq:sc-M-scalar-j} \\
\mathcal{M}_1(u)
  &= c_1(u)
     \frac{F(Tu)F(T^{-1}u)}{F(u)^2}.
  \label{eq:sc-M-scalar-1}
\end{align}
We now show that the mixed diamond equations for
$(\mathcal{M}_k, \lambda_k)$ reduce to bilinear equations
for $F$.  Define the coefficient ratios
\begin{equation}\label{eq:sc-alpha-def}
\alpha_{ij}(u) \defeq c_i(Tu)\lambda_j(Tu), \qquad
\alpha_j(u) \defeq \frac{c_1(Tu)\lambda_j(Tu)}{c_j(Tu)}.
\end{equation}

\begin{proposition}\label{prop:sc-scalar-HM}
Let $(c_k, \lambda_k)$ be a scalar solution of the
semi-discrete diamond equations
\eqref{eq:sc-diamond-C}--\eqref{eq:sc-diamond-mixed} and
the fully discrete diamond equations of
Proposition~\textup{\ref{prop:CL-gen-system}},
and let $F(u) = \det_H(I - zK(u))$ be as in
Corollary~\textup{\ref{cor:sc-Woodbury}}.
Assume the boundary
conditions\footnote{These boundary conditions are natural in
the probabilistic applications of Chapter~\ref{ch:semi-discrete-model-verifications}.}
$\lim_{\ell \to -\infty} F(T^\ell v) = 1$ and
$\lim_{\ell \to -\infty} \partial_1 F(T^\ell v) = 0$ for all
$v \in \mathcal{V}$.
\begin{enumerate}[label=\textup{(\alph*)}, leftmargin=*]
\item For each pair $i, j \geq 2$, assume
$\alpha_{ij}(u) \neq 0$ and
$\alpha_{ij}(u) \neq \alpha_{ji}(u)$ for all
$u \in \mathcal{V}$ and $i \neq j$, with the same
conditions holding for the limits
$\lim_{\ell \to -\infty} \alpha_{ij}(T^\ell v)$, which are
assumed to exist.  Then the fully discrete
mixed diamond equation reduces
to the scalar Hirota--Miwa equation with variable coefficients:
\begin{equation}\label{eq:sc-HM-dd}
\alpha_{ij}(u)F(S_i u)F(TS_j u)
- \alpha_{ji}(u)F(S_j u)F(TS_i u)
- \bigl(\alpha_{ij}(u) - \alpha_{ji}(u)\bigr)
  F(Tu)F(S_i S_j u) = 0.
\end{equation}

\item For each $j \geq 2$, assume that $c_j(u) \neq 0$,
$F$ is $C^1$ in the continuous variable, and that
$\lim_{\ell \to -\infty} \alpha_j(T^\ell u)$ exists and is
finite, for all $u \in \mathcal{V}$.  Then the
semi-discrete mixed
diamond equation \eqref{eq:sc-diamond-mixed} reduces to the
bilinear differential-difference equation:
\begin{equation}\label{eq:sc-HM-d1}
F(Tu)\partial_1 F(S_j u)
- F(S_j u)\partial_1 F(Tu)
+ \alpha_j(u)
\bigl[F(Tu)F(S_j u) - F(u)F(TS_j u)\bigr] = 0.
\end{equation}
\end{enumerate}
\end{proposition}
\begin{proof}
Part~(a) follows from Proposition~\ref{prop:scalar-HM}
applied to the unchanged discrete directions $i, j \geq 2$:
the dressed edge weights $\mathcal{M}_i, \mathcal{M}_j$ for
$i, j \geq 2$ are identical to the fully discrete case, and
the hypotheses of Proposition~\ref{prop:scalar-HM} hold by
assumption.

For part~(b), fix $j \geq 2$.  Since the resolvent
exists at every vertex, $F(u) \neq 0$ throughout
$\mathcal{V}$.  Set
\[
G_j(u) \defeq \frac{F(Tu)F(T^{-1}S_j u)}{F(u)F(S_j u)},
\qquad
G_1(u) \defeq \frac{F(Tu)F(T^{-1}u)}{F(u)^2},
\]
so that \eqref{eq:sc-M-scalar-j}--\eqref{eq:sc-M-scalar-1}
read $\mathcal{M}_j(u) = c_j(u)G_j(u)$ and
$\mathcal{M}_1(u) = c_1(u)G_1(u)$.  Substituting into the
dressed mixed diamond \eqref{eq:sc-diamond-mixed} for
$(\mathcal{M}_k, \lambda_k)$:
\begin{align}
c_1(u)G_1(u)\lambda_j(u)
&+ [\lambda_1(Tu) - \lambda_1(S_j u)]c_j(u)G_j(u)
- \lambda_j(Tu)c_1(S_j u)G_1(S_j u) \notag \\
&= [\partial_1 c_j(u)]G_j(u)
  + c_j(u)\partial_1 G_j(u).
  \label{eq:sc-var-expanded}
\end{align}
Substituting the bare mixed diamond
\eqref{eq:sc-diamond-mixed} for $(c_k, \lambda_k)$,
\[
\partial_1 c_j(u) = c_1(u)\lambda_j(u)
+ [\lambda_1(Tu) - \lambda_1(S_j u)]c_j(u)
- \lambda_j(Tu)c_1(S_j u),
\]
into \eqref{eq:sc-var-expanded} and cancelling
the $[\lambda_1(Tu) - \lambda_1(S_j u)]c_j(u)G_j(u)$
terms from each side eliminates all $\lambda_1$-dependence:
\begin{align}\label{eq:sc-var-reduced}
c_1(u)\lambda_j(u)\bigl[G_1(u) - G_j(u)\bigr]
- \lambda_j(Tu)c_1(S_j u)\bigl[G_1(S_j u) - G_j(u)\bigr]
= c_j(u)\partial_1 G_j(u).
\end{align}
 
\noindent\textit{Coefficient reduction.}
By \eqref{eq:sc-diamond-C}, $c_1(S_j u) = c_1(Tu)c_j(u)/c_j(Tu)$. Therefore 
$\lambda_j(Tu)c_1(S_j u) = \alpha_j(u)c_j(u)$ and
$c_1(u)\lambda_j(u) = \alpha_j(T^{-1}u)c_j(u)$.
Dividing \eqref{eq:sc-var-reduced} by $c_j(u)$:
\begin{equation}\label{eq:sc-var-divided}
\alpha_j(T^{-1}u)\bigl[G_1(u) - G_j(u)\bigr]
- \alpha_j(u)\bigl[G_1(S_j u) - G_j(u)\bigr]
= \partial_1 G_j(u).
\end{equation}
 
\noindent\textit{$T$-orbit invariance.}
To complete the reduction, define
\begin{equation}\label{eq:sc-J-def}
\mathcal{J}(u) \defeq
\alpha_j(u)\left[
\frac{F(u)F(TS_j u)}{F(Tu)F(S_j u)} - 1\right]
+ \partial_1\log F(Tu) - \partial_1\log F(S_j u).
\end{equation}
We claim $\mathcal{J}$ is constant along $T$-orbits.
Expanding the definitions of $G_1$ and $G_j$,
\begin{equation}\label{eq:sc-ratio-ids}
\frac{G_1(S_j u)}{G_j(u)}
= \frac{F(u)F(TS_j u)}{F(Tu)F(S_j u)},
\qquad
\frac{G_1(u)}{G_j(u)}
= \frac{F(T^{-1}u)F(S_j u)}{F(u)F(T^{-1}S_j u)}.
\end{equation}
The $F$-ratio in $\mathcal{J}(u)$ is
$G_1(S_j u)/G_j(u)$ by the first identity, and the $F$-ratio in $\mathcal{J}(T^{-1}u)$ 
is $G_1(u)/G_j(u)$ by the second identity.  Therefore
\begin{align*}
\mathcal{J}(u) - \mathcal{J}(T^{-1}u)
&= \frac{\alpha_j(u)[G_1(S_j u) - G_j(u)]
  - \alpha_j(T^{-1}u)[G_1(u) - G_j(u)]}{G_j(u)}\\
&\quad + \partial_1\log F(Tu)
  + \partial_1\log F(T^{-1}S_j u) \\
&\quad - \partial_1\log F(u)
  - \partial_1\log F(S_j u).
\end{align*}
The numerator of the first term equals the negative of the left-hand side
of \eqref{eq:sc-var-divided}, so the first line is
$-\partial_1 G_j(u)/G_j(u) = -\partial_1\log G_j(u)$, and 
the second and third lines combine to equal $+\partial_1\log G_j(u)$ by the
definition of $G_j$. Hence $\mathcal{J}$
is constant along $T$-orbits, and the  boundary conditions give
\[\mathcal{J}(u)
= \lim_{n \to \infty}\mathcal{J}(T^{-n}u) = 0,\]
since the $F$-ratio tends to $1$ and the
$\partial_1\log F$ terms vanish by assumption. 
Multiplying by $-F(Tu)F(S_j u)$ yields
\eqref{eq:sc-HM-d1}.
\end{proof}

\begin{remark}
\label{rem:sc-opposite-boundary}
The proof of Proposition~\ref{prop:sc-scalar-HM} uses only
that the $T$-orbit invariant $\mathcal{J}$ tends to zero at
one end of the orbit.  The same scalar equations
\eqref{eq:sc-HM-dd} and \eqref{eq:sc-HM-d1} therefore hold
under the opposite boundary normalization
\[
F(T^\ell v) \to 1, \qquad
\partial_1 F(T^\ell v) \to 0,
\qquad \ell \to +\infty,
\]
provided the coefficient ratios $\alpha_{ij}$ and $\alpha_j$
have finite limits as $\ell \to +\infty$, with the same
nondegeneracy conditions as in parts~\textup{(a)} and~\textup{(b)}.
\end{remark}

\begin{remark}
\label{rem:sc-product-graph}
The product graph construction of \S\ref{sec:product-graph}
extends to the semi-discrete setting: each statement below
is the $O(\epsilon)$ residue of its fully discrete
counterpart under the scaling
$\mathcal{S}_1 = e^{\epsilon\partial_1}$, illustrating the
role of the discrete framework as the master theory.  The
direct verifications are deferred to a future work.

The constructions of $\Psi, \Phi, K$ on the product graph
$\mathcal{G}^m$
(\eqref{eq:Psi-def}--\eqref{eq:K-def}) are unchanged in form:
they remain $\mathcal{T}$-orbit sums of seed data against a
strictly lower-triangular propagator $B_u$.  The admissible
propagator conditions carry over with one modification: the
discrete $\mathcal{S}_1$-compatibility condition is replaced by
a $\partial_1$-compatibility condition,
\begin{align*}
&\partial_1 B_u(\mathcal{T}x, y)
+ C_1(x)B_u(x, y)
- D_1(\mathcal{T}x)B_u(\mathcal{T}x, y) \\
& - B_u(\mathcal{T}x, \mathcal{T}y)C_1(y)
+ B_u(\mathcal{T}x, y)D_1(y)
= B_u(\mathcal{T}x, \mathcal{T}u)
  C_1(u)B_u(u, y).
\end{align*}
The $\mathcal{T}$-covariance, $\mathcal{T}$-splitting, and
$\mathcal{S}_j$-compatibility conditions for $j \geq 2$ are
inherited without modification.  The corrected
$\Lambda$-weights for discrete directions take the same form
as in \S\ref{sec:product-graph}:
$\Lambda_j(u) = P(u)D_j(u)P(\mathcal{S}_j u)^{-1}$, where
$P(v) = I + B_v(v, v)$.  For the continuous direction, the
corrected weight acquires an additional derivative term:
\[\Lambda_1(u) = P(u)D_1(u)P(u)^{-1}
+ (\partial_1 P(u))P(u)^{-1}.\]
\end{remark}

\section{Dual lattice presentation}
\label{sec:sc-dual}

Some KPZ models naturally present their semi-discrete linear
data with the derivative falling on $\Psi(Tu)$ rather than
$\Psi(u)$, corresponding to the scaling
$S_1 = Te^{\epsilon\partial_1}$ rather than
$S_1 = e^{\epsilon\partial_1}$.  This subsection develops
the resulting formulas in a form ready for direct use in the
model verifications.  No new algebraic structure is
introduced: the dual presentation is the standard
semi-discrete theory written in a different lattice basis.

Retain the setup of \S\ref{sec:sc-linear}, and consider the
overdetermined linear system for
$\Psi \colon \mathcal{V} \to \Hom(E, H)$:
\begin{align}
\partial_1\Psi(Tu)
  &= \Psi(Tu)C_1(u) - \Psi(u)\Lambda_1(u),
  \label{eq:sc-dual-linear-1} \\
\Psi(S_j u)
  &= \Psi(Tu)C_j(u) - \Psi(u)\Lambda_j(u),
  \qquad j \geq 2.
  \label{eq:sc-dual-linear-j}
\end{align}
\begin{lemma}
\label{lem:sc-dual}
Define $\widetilde T \defeq T^{-1}$,
$\widetilde S_j \defeq T^{-1}S_j$ for $j \geq 2$,
$\widetilde C_k(u) \defeq -\Lambda_k(T^{-1}u)$, and
$\widetilde\Lambda_k(u) \defeq -C_k(T^{-1}u)$.  Then the
dual system
\eqref{eq:sc-dual-linear-1}--\eqref{eq:sc-dual-linear-j}
is the standard semi-discrete linear problem
\eqref{eq:sc-linear-1}--\eqref{eq:sc-linear-j} with data
$(\widetilde C_k, \widetilde\Lambda_k)$ on the lattice
generated by $\widetilde T, \widetilde S_j$.
\end{lemma}

\begin{proof}
Evaluate \eqref{eq:sc-dual-linear-1} at $T^{-1}u$:
\[
\partial_1\Psi(u)
= \Psi(u)C_1(T^{-1}u)
  - \Psi(T^{-1}u)\Lambda_1(T^{-1}u).
\]
Since $\widetilde T = T^{-1}$, this is
\eqref{eq:sc-linear-1} with data
$(\widetilde C_1, \widetilde\Lambda_1)$ and shift
$\widetilde T$.  Similarly, evaluating
\eqref{eq:sc-dual-linear-j} at $T^{-1}u$ gives
\eqref{eq:sc-linear-j} with data
$(\widetilde C_j, \widetilde\Lambda_j)$ and shift
$\widetilde S_j = T^{-1}S_j$.
\end{proof}

\begin{corollary}
\label{cor:sc-dual-diamonds}
The dual linear system
\eqref{eq:sc-dual-linear-1}--\eqref{eq:sc-dual-linear-j}
is compatible if and only if two families of conditions hold.
First, for each $j \geq 2$ and $u \in \mathcal{V}$, the
following dual diamond equations hold:
\begin{align}
\partial_1 C_j(Tu) + C_1(Tu)C_j(Tu)
  &= C_j(Tu)C_1(S_j u),
  \label{eq:sc-dual-diamond-C} \\
\Lambda_1(u)\Lambda_j(Tu)
  &= \Lambda_j(u)\Lambda_1(S_j u),
  \label{eq:sc-dual-diamond-L} \\
\partial_1\Lambda_j(Tu)
  + C_1(u)\Lambda_j(Tu)
  &+ \Lambda_1(Tu)C_j(Tu)
  - C_j(u)\Lambda_1(S_j u) \notag \\
  &\quad - \Lambda_j(Tu)C_1(S_j u) = 0.
  \label{eq:sc-dual-diamond-mixed}
\end{align}
Second, for each pair $i, j \geq 2$ and
$u \in \mathcal{V}$, the fully discrete diamond equations
\eqref{eq:diamond-C-ij}--\eqref{eq:diamond-mixed-ij} of
Proposition~\ref{prop:CL-gen-system} hold.
\end{corollary}

\begin{proof}
By Lemma~\ref{lem:sc-dual}, the dual system
\eqref{eq:sc-dual-linear-1}--\eqref{eq:sc-dual-linear-j}
is the standard semi-discrete linear problem
\eqref{eq:sc-linear-1}--\eqref{eq:sc-linear-j} with data
$(\widetilde C_k, \widetilde\Lambda_k)$.  Its compatibility
conditions are therefore the diamond equations of
Proposition~\ref{prop:sc-diamond} for this data; expanding
the definitions and substituting $u \mapsto T^2 u$ yields
the stated equations.
\end{proof}

\begin{remark}
The system
\eqref{eq:sc-dual-linear-1}--\eqref{eq:sc-dual-linear-j}
arises from the fully discrete linear problem
\eqref{eq:Linear-Psi-3} under the scaling
$S_1 = Te^{\epsilon\partial_1}$, with
$C_1^{(\epsilon)} = I + \epsilon C_1
+ O(\epsilon^2)$ and
$\Lambda_1^{(\epsilon)} = \epsilon\Lambda_1
+ O(\epsilon^2)$.
\end{remark}

Fix a solution $(C_k, \Lambda_k)$ of the dual diamond
equations of Corollary~\ref{cor:sc-dual-diamonds}, and let
$\Psi \colon \mathcal{V} \to \Hom(E, H)$ solve the dual
linear problem
\eqref{eq:sc-dual-linear-1}--\eqref{eq:sc-dual-linear-j}.
The dual adjoint linear problem for
$\Phi \colon \mathcal{V} \to \Hom(H, E)$ is the system
\begin{align}
\partial_1\Phi(Tu)
  &= \Lambda_1(Tu)\Phi(T^2u) - C_1(u)\Phi(Tu),
  \label{eq:sc-dual-adjoint-1} \\
\Phi(Tu)
  &= C_j(u)\Phi(S_j u) - \Lambda_j(Tu)\Phi(TS_j u),
  \qquad j \geq 2.
  \label{eq:sc-dual-adjoint-j}
\end{align}
This is the standard adjoint problem
\eqref{eq:sc-adjoint-1}--\eqref{eq:sc-adjoint-j} written in
the dual basis of Lemma~\ref{lem:sc-dual}; its compatibility
conditions are again the dual diamond equations of
Corollary~\ref{cor:sc-dual-diamonds}.

\begin{definition}\label{def:sc-dual-dressing}
Let $\Psi$ and $\Phi$ solve the dual linear and adjoint
linear problems.  A map
$K \colon \mathcal{V} \to \End(H)$ is
\emph{dual dressing compatible} with $(\Psi, \Phi)$ if:
\begin{enumerate}[label=\arabic*., leftmargin=*]
\item For every $u \in \mathcal V$ and $j \geq 2$,
\begin{align}
K(Tu) - K(u) &= \Psi(Tu)\Phi(Tu),
  \label{eq:sc-dual-K-T} \\
\partial_1 K(u) &= \Psi(u)\Lambda_1(u)\Phi(Tu),
  \label{eq:sc-dual-K-d1} \\
K(S_j u) - K(u)
  &= \Psi(Tu)C_j(u)\Phi(S_j u),
  \qquad j \geq 2.
  \label{eq:sc-dual-K-Sj}
\end{align} 
\item There exists $z \in \mathbb{F}$ such that the resolvent
exists at every vertex:
\[
R(u) \defeq (I - zK(u))^{-1} \in \End(H),
\qquad u \in \mathcal{V}.
\]
\end{enumerate}
\end{definition}

\begin{example}\label{ex:sc-dual-simple-kernel}
Under suitable analytic assumptions, the kernel
\[
K(u) = \sum_{p \leq 0} \Psi(T^p u)\Phi(T^p u)
\]
is dual dressing compatible.  The $T$- and
$S_j$-conditions follow by telescoping.  The
$\partial_1$-condition follows by differentiating termwise:
the $C_1$-terms cancel by the dual linear and adjoint
problems, and the remaining $\Lambda_1$-terms telescope to
the endpoint $\Psi(u)\Lambda_1(u)\Phi(Tu)$.
\end{example}

\begin{theorem}
\label{thm:sc-dual-Darboux}
Let $K$ be dual dressing compatible with $(\Psi, \Phi)$, with
resolvent $R(u) = (I - zK(u))^{-1}$.  Define the dressed
observable and dressed waves by
\[
\mathcal{M}(u) \defeq I + z\Phi(u)R(u)\Psi(u),
\]
\[
\widehat\Psi(u) \defeq R(u)\Psi(u), \qquad
\widehat\Phi(u) \defeq \Phi(u)R(T^{-1}u).
\]
Then $\mathcal{M}(u)$ is invertible for every
$u \in \mathcal{V}$.  The dressed waves
$\widehat\Psi, \widehat\Phi$ satisfy the dual
linear problem
\eqref{eq:sc-dual-linear-1}--\eqref{eq:sc-dual-linear-j} and
dual adjoint problem
\eqref{eq:sc-dual-adjoint-1}--\eqref{eq:sc-dual-adjoint-j}
respectively, with data $(\mathcal{M}_k, \Lambda_k)$: for each $j \geq 2$,
\begin{align*}
\partial_1\widehat\Psi(Tu)
  &= \widehat\Psi(Tu)\mathcal{M}_1(u)
     - \widehat\Psi(u)\Lambda_1(u), \\
\widehat\Psi(S_j u)
  &= \widehat\Psi(Tu)\mathcal{M}_j(u)
     - \widehat\Psi(u)\Lambda_j(u),
\shortintertext{and}
\partial_1\widehat\Phi(Tu)
  &= \Lambda_1(Tu)\widehat\Phi(T^2u)
     - \mathcal{M}_1(u)\widehat\Phi(Tu), \\
\widehat\Phi(Tu)
  &= \mathcal{M}_j(u)\widehat\Phi(S_j u)
     - \Lambda_j(Tu)\widehat\Phi(TS_j u),
\end{align*}
where the dressed edge weights are given by
\begin{align}
\mathcal{M}_1(u)
  &\defeq \mathcal{M}(Tu)^{-1}C_1(u)\mathcal{M}(Tu)
  + \mathcal{M}(Tu)^{-1}\partial_1\mathcal{M}(Tu),
  \label{eq:sc-dual-dressed-C1} \\
\mathcal{M}_j(u)
  &\defeq \mathcal{M}(Tu)^{-1}C_j(u)\mathcal{M}(S_j u),
  \label{eq:sc-dual-dressed-Cj}
\end{align}
and $\widehat K(u) \defeq K(u)R(u)$ is dual dressing compatible
with $(\widehat\Psi, \widehat\Phi)$.
Moreover, $(\mathcal{M}_k, \Lambda_k)$ satisfies the dual
diamond equations of
Corollary~\ref{cor:sc-dual-diamonds}.
\end{theorem}

\begin{proof}
The strategy is to convert the dual dressing data to the
standard semi-discrete setting of \S\ref{sec:sc-Darboux}
via Lemma~\ref{lem:sc-dual}, apply
Theorem~\ref{thm:sc-Darboux}, and gauge the result back to
the dual setting.

\noindent\textit{Tilted dressing data.}
Lemma~\ref{lem:sc-dual} and
Corollary~\ref{cor:sc-dual-diamonds} ensure that
$(\widetilde C_k, \widetilde\Lambda_k)$ satisfies the
diamond equations of Proposition~\ref{prop:sc-diamond}.
Define $\widetilde K(u) \defeq -K(T^{-1}u)$ and
$\widetilde z \defeq -z$.  The dual dressing compatibility
conditions
\eqref{eq:sc-dual-K-T}--\eqref{eq:sc-dual-K-Sj} transport
to the dressing conditions of
Definition~\ref{def:sc-dressing} for $\widetilde K$ with
this data.  The $T$- and $\partial_1$-conditions are
immediate:
\[
\widetilde K(\widetilde Tu) - \widetilde K(u)
= \Psi(\widetilde Tu)\Phi(\widetilde Tu),
\qquad
\partial_1\widetilde K(u)
= \Psi(\widetilde Tu)\widetilde C_1(u)\Phi(u).
\]
For $j \geq 2$, write
\[
\widetilde K(\widetilde S_ju) - \widetilde K(u)
= \bigl[K(T^{-1}u) - K(T^{-2}u)\bigr]
  - \bigl[K(T^{-2}S_ju) - K(T^{-2}u)\bigr],
\]
where we substituted $\widetilde K(u) = -K(T^{-1}u)$ and
inserted $\pm K(T^{-2}u)$.
The two brackets are \eqref{eq:sc-dual-K-T} and
\eqref{eq:sc-dual-K-Sj} at $T^{-2}u$; substituting and
factoring:
$\Psi(T^{-1}u)\bigl[\Phi(T^{-1}u)
  - C_j(T^{-2}u)\Phi(T^{-2}S_ju)\bigr]$.
The dual adjoint equation \eqref{eq:sc-dual-adjoint-j} at
$T^{-2}u$ reduces the bracket to
$-\Lambda_j(T^{-1}u)\Phi(T^{-1}S_ju)
= \widetilde C_j(u)\Phi(\widetilde S_ju)$, giving the
$S_j$-condition of Definition~\ref{def:sc-dressing}.

Theorem~\ref{thm:sc-Darboux} therefore applies to
$(\Psi, \Phi, \widetilde K)$ in the tilted basis.
The tilted resolvent is
$\widetilde R(u) = R(T^{-1}u)$.  The tilted Darboux
observable is
$\widetilde{\mathcal M}(u)
= I - z\Phi(u)R(T^{-1}u)\Psi(u)
= \mathcal M(u)^{-1}$,
where the last equality is \eqref{eq:sc-M-inverse}
in the tilted basis.
Since $\widetilde{\mathcal M}(u)$ is well-defined by
Theorem~\ref{thm:sc-Darboux}, $\mathcal M(u)$ is invertible
at every vertex.

\noindent\textit{The dressed wave equations.}
Theorem~\ref{thm:sc-Darboux} produces dressed waves
$\Psi'(u) \defeq \widetilde R(u)\Psi(u)
= R(T^{-1}u)\Psi(u)$ and
$\Phi'(u) \defeq \Phi(u)\widetilde R(\widetilde T^{-1}u)
= \Phi(u)R(u)$
for the tilted system.  The intertwining identities
\eqref{eq:sc-intertwine-wave} and
\eqref{eq:sc-adjoint-intertwine} in the tilted basis give
\[
\widehat\Psi(u) = \Psi'(u)\mathcal M(u),
\qquad
\widehat\Phi(u) = \mathcal M(u)^{-1}\Phi'(u).
\]

Theorem~\ref{thm:sc-Darboux} gives the tilted dressed
discrete wave equation
\[
\Psi'(\widetilde S_j u)
= \Psi'(\widetilde Tu)\widetilde{\mathcal M}_j(u)
  - \Psi'(u)\widetilde\Lambda_j(u).
\]
Evaluating at $Tu$ and substituting the definitions of
Lemma~\ref{lem:sc-dual} together with
$\widetilde{\mathcal M} = \mathcal M^{-1}$:
\[
\Psi'(S_ju) = \Psi'(Tu)C_j(u)
  - \Psi'(u)\mathcal M(u)\Lambda_j(u)
    \mathcal M(S_ju)^{-1}.
\]
Substituting the gauge relation
$\widehat\Psi = \Psi'\mathcal M$ and multiplying on the
right by $\mathcal M(S_ju)$:
\[
\widehat\Psi(S_ju)
= \widehat\Psi(Tu)\mathcal M(Tu)^{-1}
  C_j(u)\mathcal M(S_ju)
  - \widehat\Psi(u)\Lambda_j(u).
\]
This is the discrete dressed wave equation with
dressed weight $\mathcal M_j(u)$ given by
\eqref{eq:sc-dual-dressed-Cj}.  For the continuous direction,
Theorem~\ref{thm:sc-Darboux} gives the tilted dressed
continuous wave equation
$\partial_1\Psi'(u)
= \Psi'(\widetilde Tu)\widetilde{\mathcal M}_1(u)
  - \Psi'(u)\widetilde\Lambda_1(u)$.
Evaluating at $Tu$ and substituting the definitions of
Lemma~\ref{lem:sc-dual} together with
$\widetilde{\mathcal M} = \mathcal M^{-1}$:
\[
\partial_1\Psi'(Tu) = \Psi'(Tu)C_1(u)
  - \Psi'(u)\mathcal M(u)\Lambda_1(u)
    \mathcal M(Tu)^{-1}.
\]
Differentiating the gauge relation
$\widehat\Psi(Tu) = \Psi'(Tu)\mathcal M(Tu)$ and
substituting the above for $\partial_1\Psi'(Tu)$:
\[
\partial_1\widehat\Psi(Tu)
= \Psi'(Tu)C_1(u)\mathcal M(Tu)
  + \Psi'(Tu)\partial_1\mathcal M(Tu)
  - \Psi'(u)\mathcal M(u)\Lambda_1(u).
\]
Substituting
$\Psi'(Tu) = \widehat\Psi(Tu)\mathcal M(Tu)^{-1}$ and
$\Psi'(u)\mathcal M(u) = \widehat\Psi(u)$:
\begin{align*}
\partial_1\widehat\Psi(Tu)
= \widehat\Psi(Tu)\bigl[
  \mathcal M(Tu)^{-1}C_1(u)\mathcal M(Tu)
  + \mathcal M(Tu)^{-1}\partial_1\mathcal M(Tu)\bigr]
- \widehat\Psi(u)\Lambda_1(u).
\end{align*}
This is the continuous dressed wave equation with
dressed weight $\mathcal M_1(u)$ given by
\eqref{eq:sc-dual-dressed-C1}.

\noindent\textit{The dressed adjoint equations.}
Theorem~\ref{thm:sc-Darboux} gives the
tilted dressed discrete adjoint equation
$\Phi'(\widetilde Tu)
= \widetilde{\mathcal M}_j(u)\Phi'(\widetilde S_ju)
  - \widetilde\Lambda_j(\widetilde Tu)
    \Phi'(\widetilde T\widetilde S_ju)$.
Evaluating at $T^2u$ and substituting the definitions of
Lemma~\ref{lem:sc-dual} together with
$\widetilde{\mathcal M} = \mathcal M^{-1}$:
\[
\Phi'(Tu) = C_j(u)\Phi'(S_ju)
  - \mathcal M(Tu)\Lambda_j(Tu)
    \mathcal M(TS_ju)^{-1}\Phi'(TS_ju).
\]
Substituting the gauge relation
$\Phi' = \mathcal M\widehat\Phi$ and multiplying on the
left by $\mathcal M(Tu)^{-1}$:
\[
\widehat\Phi(Tu)
= \mathcal M(Tu)^{-1}C_j(u)\mathcal M(S_ju)
  \widehat\Phi(S_ju)
  - \Lambda_j(Tu)\widehat\Phi(TS_ju).
\]
This is the discrete dressed adjoint equation with
dressed weight $\mathcal M_j(u)$ given by
\eqref{eq:sc-dual-dressed-Cj}.
For the continuous adjoint,
Theorem~\ref{thm:sc-Darboux} gives 
$\partial_1\Phi'(\widetilde Tu)
= \widetilde\Lambda_1(\widetilde Tu)\Phi'(\widetilde Tu)
  - \widetilde{\mathcal M}_1(u)\Phi'(u)$.
Evaluating at $T^2u$ and substituting the definitions of
Lemma~\ref{lem:sc-dual} together with
$\widetilde{\mathcal M} = \mathcal M^{-1}$:
\[
\partial_1\Phi'(Tu)
= \mathcal M(Tu)\Lambda_1(Tu)\mathcal M(T^2u)^{-1}\Phi'(T^2u)
  - C_1(u)\Phi'(Tu).
\]
Differentiating the gauge relation
$\widehat\Phi(Tu) = \mathcal M(Tu)^{-1}\Phi'(Tu)$ and
substituting the above for $\partial_1\Phi'(Tu)$:
\begin{align*}
\partial_1\widehat\Phi(Tu)
&= -\mathcal M(Tu)^{-1}\partial_1\mathcal M(Tu)
    \widehat\Phi(Tu) \\
&\quad + \Lambda_1(Tu)\mathcal M(T^2u)^{-1}\Phi'(T^2u)
  - \mathcal M(Tu)^{-1}C_1(u)\Phi'(Tu).
\end{align*}
Substituting
$\Phi'(T^2u) = \mathcal M(T^2u)\widehat\Phi(T^2u)$ and
$\Phi'(Tu) = \mathcal M(Tu)\widehat\Phi(Tu)$:
\begin{align*}
\partial_1\widehat\Phi(Tu)
= \Lambda_1(Tu)\widehat\Phi(T^2u)
- \bigl[\mathcal M(Tu)^{-1}C_1(u)\mathcal M(Tu)
    + \mathcal M(Tu)^{-1}\partial_1\mathcal M(Tu)
  \bigr]\widehat\Phi(Tu).
\end{align*}
This is the continuous dressed adjoint equation with
dressed weight $\mathcal M_1(u)$ given by
\eqref{eq:sc-dual-dressed-C1}.

\noindent\textit{The dressed kernel.}
Since $K(u)$ and
$R(u)$ commute at each vertex, $\widehat K(u) = R(u)K(u)$.
The identity
$(I - \eta\widehat K(u))(I - zK(u)) = I - (z+\eta)K(u)$
shows that $(I - \eta\widehat K(u))^{-1}$ exists whenever
$(I - (z+\eta)K(u))^{-1}$ does; in particular at
$\eta = -z$, where the condition is vacuous.  When $z = 0$,
$R = I$ and $\widehat K = K$, so the kernel identities
reduce to those of $K$ itself.  For $z \neq 0$, the $T$-
and $S_j$-conditions are identical to the
corresponding steps in Theorem~\ref{thm:Darboux}, using
the resolvent identity and the intertwining identities
established above in place of their discrete counterparts.
For the $\partial_1$-condition, differentiating
$\widehat K(u) = z^{-1}(R(u) - I)$ and substituting
$\partial_1 R(u) = zR(u)(\partial_1 K(u))R(u)$:
\[
\partial_1\widehat K(u)
= R(u)\Psi(u)\Lambda_1(u)\Phi(Tu)R(u)
= \widehat\Psi(u)\Lambda_1(u)\widehat\Phi(Tu).
\]

\noindent\textit{The dressed diamond equations.}
It remains to verify the dual diamond equations of
Corollary~\ref{cor:sc-dual-diamonds} for
$(\mathcal M_k, \Lambda_k)$.  Equation
\eqref{eq:sc-dual-diamond-L} holds because $\Lambda_k$ is
unchanged.

For \eqref{eq:sc-dual-diamond-C}, substitute
\eqref{eq:sc-dual-dressed-Cj}--\eqref{eq:sc-dual-dressed-C1}.
The Leibniz rule gives
\begin{align*}
\partial_1\mathcal M_j(Tu)
&= -\mathcal M(T^2u)^{-1}\partial_1\mathcal M(T^2u)
     \mathcal M(T^2u)^{-1}C_j(Tu)\mathcal M(TS_ju) \\
& + \mathcal M(T^2u)^{-1}\partial_1 C_j(Tu)
     \mathcal M(TS_ju)
   + \mathcal M(T^2u)^{-1}C_j(Tu)
     \partial_1\mathcal M(TS_ju),
\end{align*}
and the product $\mathcal M_1(Tu)\mathcal M_j(Tu)$ expands
as
\begin{align*}
\mathcal M(T^2u)^{-1}C_1(Tu)C_j(Tu)\mathcal M(TS_ju)
+ \mathcal M(T^2u)^{-1}\partial_1\mathcal M(T^2u)
  \mathcal M(T^2u)^{-1}C_j(Tu)\mathcal M(TS_ju).
\end{align*}
The $\partial_1\mathcal M(T^2u)$ terms cancel between the
two, leaving
\begin{align*}
&\partial_1\mathcal M_j(Tu)
  + \mathcal M_1(Tu)\mathcal M_j(Tu)
= \mathcal M(T^2u)^{-1}\bigl[
  \partial_1 C_j(Tu) + C_1(Tu)C_j(Tu)\bigr]
  \mathcal M(TS_ju) \\
&\quad + \mathcal M(T^2u)^{-1}C_j(Tu)
    \partial_1\mathcal M(TS_ju).
\end{align*}
The right-hand side
$\mathcal M_j(Tu)\mathcal M_1(S_ju)$ expands as
\[
\mathcal M(T^2u)^{-1}C_j(Tu)C_1(S_ju)\mathcal M(TS_ju)
+ \mathcal M(T^2u)^{-1}C_j(Tu)
  \partial_1\mathcal M(TS_ju).
\]
The $\partial_1\mathcal M(TS_ju)$ terms cancel, and
cancelling $\mathcal M(T^2u)^{-1}$ on the left and
$\mathcal M(TS_ju)$ on the right reduces
\eqref{eq:sc-dual-diamond-C} for
$(\mathcal M_k, \Lambda_k)$ to
\eqref{eq:sc-dual-diamond-C} for $(C_k, \Lambda_k)$.

For \eqref{eq:sc-dual-diamond-mixed}, substitute
\eqref{eq:sc-dual-dressed-Cj}--\eqref{eq:sc-dual-dressed-C1}
and multiply on the left by $\mathcal M(Tu)$ and on the
right by $\mathcal M(TS_ju)^{-1}$.  The five terms become:
\begin{gather*}
\mathcal M(Tu)\partial_1\Lambda_j(Tu)
  \mathcal M(TS_ju)^{-1}, \\
C_1(u)\mathcal M(Tu)\Lambda_j(Tu)
  \mathcal M(TS_ju)^{-1}
+ \partial_1\mathcal M(Tu)\Lambda_j(Tu)
  \mathcal M(TS_ju)^{-1}, \\
\mathcal M(Tu)\Lambda_1(Tu)\mathcal M(T^2u)^{-1}
  C_j(Tu), \\
-C_j(u)\mathcal M(S_ju)\Lambda_1(S_ju)
  \mathcal M(TS_ju)^{-1}, \\
-\mathcal M(Tu)\Lambda_j(Tu)
  \mathcal M(TS_ju)^{-1}C_1(S_ju) \\
\quad - \mathcal M(Tu)\Lambda_j(Tu)
  \mathcal M(TS_ju)^{-1}\partial_1\mathcal M(TS_ju)
  \mathcal M(TS_ju)^{-1}.
\end{gather*}
The terms involving $\partial_1\Lambda_j(Tu)$,
$\partial_1\mathcal M(Tu)$, and
$\partial_1\mathcal M(TS_ju)$ combine by the Leibniz rule
into
$\partial_1\bigl[\mathcal M(Tu)\Lambda_j(Tu)
  \mathcal M(TS_ju)^{-1}\bigr]$.
The standard mixed diamond equation for
$(\widetilde{\mathcal M}_k, \widetilde\Lambda_k)$,
written in the dual edge weights $(C_k, \Lambda_k)$ via
Lemma~\ref{lem:sc-dual} (evaluated at $T^2u$), reads
\begin{align*}
0
&= \partial_1\bigl[\mathcal M(Tu)\Lambda_j(Tu)
     \mathcal M(TS_ju)^{-1}\bigr] \\
& + C_1(u)\mathcal M(Tu)\Lambda_j(Tu)
     \mathcal M(TS_ju)^{-1}
   + \mathcal M(Tu)\Lambda_1(Tu)\mathcal M(T^2u)^{-1}
     C_j(Tu) \\
&- C_j(u)\mathcal M(S_ju)\Lambda_1(S_ju)
     \mathcal M(TS_ju)^{-1}
   - \mathcal M(Tu)\Lambda_j(Tu)
     \mathcal M(TS_ju)^{-1}C_1(S_ju),
\end{align*}
which is precisely the expression above.  This vanishes by
Theorem~\ref{thm:sc-Darboux}.  Since $\mathcal M(Tu)$ and
$\mathcal M(TS_ju)$ are invertible,
\eqref{eq:sc-dual-diamond-mixed} holds for
$(\mathcal M_k, \Lambda_k)$.

The proof of the fully discrete diamond equations for
$i, j \geq 2$ uses only the discrete-direction linear
problem, adjoint problem, and dressing conditions, all of
which are identical to the fully discrete framework.  The
proof of Theorem~\ref{thm:Darboux} therefore applies
directly.
\end{proof}

\begin{corollary}
\label{cor:sc-dual-Woodbury}
In the setting of Theorem~\ref{thm:sc-dual-Darboux}, suppose
that $\dim E < \infty$ and that
$F(u) \defeq \det_H(I-zK(u))$ is well-defined for every
$u \in \mathcal{V}$.  Then
\begin{equation}\label{eq:sc-dual-Woodbury-det}
F(u)=F(Tu)\det_E\bigl(\mathcal{M}(Tu)\bigr).
\end{equation}
Consequently, for $j \geq 2$,
\begin{equation}\label{eq:sc-dual-det-dressed-j}
\det_E(\mathcal{M}_j(u))F(u)F(S_j u)
=\det_E(C_j(u))F(Tu)F(T^{-1}S_j u).
\end{equation}
\end{corollary}

\begin{proof}
The identity \eqref{eq:sc-dual-Woodbury-det} depends only on
the $T$-dressing compatibility \eqref{eq:sc-dual-K-T},
so the proof is identical to
Corollary~\ref{cor:sc-Woodbury}.  Taking $\det_E$ of
\eqref{eq:sc-dual-dressed-Cj} and applying
\eqref{eq:sc-dual-Woodbury-det} at $u$ and $T^{-1}S_ju$ gives
\eqref{eq:sc-dual-det-dressed-j}.
\end{proof}

When $E = \mathbb{F}$, the edge weights
$C_k, \Lambda_k, \mathcal{M}_k$ reduce to scalar-valued
functions $c_k, \lambda_k, \mathcal{M}_k$, and
Corollary~\ref{cor:sc-dual-Woodbury} gives
\begin{align}
\mathcal{M}_j(u)
&= c_j(u)
   \frac{F(Tu)F(T^{-1}S_j u)}{F(u)F(S_j u)},
   \qquad j \geq 2,
   \label{eq:sc-dual-M-scalar-j} \\
\mathcal{M}_1(u)
&= c_1(u)
   + \partial_1\log F(u)
   - \partial_1\log F(Tu).
   \label{eq:sc-dual-M-scalar-1}
\end{align}
The first identity follows immediately from
\eqref{eq:sc-dual-det-dressed-j}.  For the second,
\eqref{eq:sc-dual-dressed-C1} gives
\[
\mathcal{M}_1(u) = c_1(u) + \partial_1\log\mathcal{M}(Tu),
\]
and \eqref{eq:sc-dual-Woodbury-det} gives
$\mathcal{M}(Tu) = F(u)/F(Tu)$.

Define the coefficient ratio
\begin{equation}\label{eq:sc-dual-beta-def}
\beta_j(u) \defeq \frac{c_j(u)\lambda_1(u)}{\lambda_j(u)}.
\end{equation}

\begin{corollary}
\label{cor:sc-dual-scalar}
Suppose $(c_k, \lambda_k)$ is a scalar solution of the
dual diamond equations of
Corollary~\ref{cor:sc-dual-diamonds},
and let $F(u) = \det_H(I - zK(u))$ be as in
Corollary~\textup{\ref{cor:sc-dual-Woodbury}}.  Assume
$\lambda_j(u) \neq 0$ and $F$ is $C^1$ in the continuous
variable for all $u \in \mathcal{V}$, and assume the
boundary normalization
$F(T^\ell v) \to 1$, $\partial_1 F(T^\ell v) \to 0$ as
$\ell \to -\infty$.  Assume also that
$\beta_j(T^\ell v)$
has a finite limit as $\ell \to -\infty$ for every
$v \in \mathcal{V}$.  Then the bilinear
differential-difference equation takes the form
\begin{equation}\label{eq:sc-dual-HM}
F(u)\partial_1 F(S_j u)
  - F(S_j u)\partial_1 F(u)
+ \beta_j(u)
\bigl[F(Tu)F(T^{-1}S_j u) - F(u)F(S_j u)\bigr] = 0.
\end{equation}
\end{corollary}

\begin{proof}
Proposition~\ref{prop:sc-scalar-HM}(b) applies in the
tilted basis of Lemma~\ref{lem:sc-dual}, with tilted
Fredholm determinant $\widetilde F(u) = F(T^{-1}u)$; the
boundary normalization $F(T^\ell v) \to 1$ as $\ell \to -\infty$ corresponds to
$\widetilde F(\widetilde T^\ell v) \to 1$ as $\ell \to +\infty$, which is
covered by Remark~\ref{rem:sc-opposite-boundary}.  Evaluating
the resulting equation at $T^2u$, the tilted coefficient is
\[
\widetilde\alpha_j(T^2u)
= \frac{(-\lambda_1(u))(-c_j(u))}{-\lambda_j(u)}
= -\beta_j(u),
\]
and substituting $\widetilde F(v) = F(T^{-1}v)$ into \eqref{eq:sc-HM-d1} yields
\eqref{eq:sc-dual-HM}.
\end{proof}

\begin{remark}
The same equation \eqref{eq:sc-dual-HM} holds under the
opposite boundary normalization $F(T^\ell v) \to 1$,
$\partial_1 F(T^\ell v) \to 0$ as $\ell \to +\infty$,
provided $\beta_j(T^\ell v)$
has a finite limit as $\ell \to +\infty$.
\end{remark}

%% file: chapters/4-the-parabolic-framework.tex
\chapter{The Parabolic Framework}
\label{ch:the-parabolic-framework}

{
  \setlength{\parskip}{0pt}
}
\label{sec:dc-framework-clean}

This chapter develops the parabolic analogue of the semi-discrete
framework (Chapter~\ref{ch:the-semi-discrete-framework}), replacing
the distinguished shift $T$ by a continuous derivative $\partial_2$
to produce a linear problem parabolic in $\partial_1$, with one
surviving discrete shift $S_3$.  The parabolic diamond equations,
which arise as compatibility conditions for the resulting linear
problem (Proposition~\ref{prop:dc-clean-diamond}), admit a Darboux
transformation (Theorem~\ref{thm:dc-clean-general-darboux}) with
inverted character: the $C$-weights remain fixed while the
$\Lambda$-weights absorb the dressing.  The normalized trace-defect
identity (Proposition~\ref{prop:dc-clean-normalized-trace-defect-direct})
extracts a bilinear identity for the Fredholm determinant that reduces
to a scalar Hirota identity when $\dim E = 1$
(Corollary~\ref{cor:dc-clean-scalar-trace-defect}).  The parabolic
framework arises from the semi-discrete theory under the scaling
$T^{(\epsilon)} = e^{\epsilon\partial_2}$,
$\partial_1^{(\epsilon)} = \epsilon^{-1}\partial_2 + \partial_1$,
$\epsilon \to 0$
(Remark~\ref{rem:dc-clean-scaling}), but all results are stated
and proved independently.

\section{The Parabolic Linear Problem}
\label{sec:dc-linear-clean}

Let $E$ and $H \neq 0$ be vector spaces over $\mathbb{F} = \R$ or $\C$,
let $\mathcal{V}'$ be a lattice equipped with an invertible shift $S_3$,
and set $\mathcal{V} \defeq I_1 \times I_2 \times \mathcal{V}'$ for
open intervals $I_1, I_2 \subseteq \R$.  Write $\partial_1$ and
$\partial_2$ for the derivatives in the two continuous coordinates;
the shift $S_3$ acts on $\mathcal{V}'$ and commutes with both.
Assign endomorphisms\footnote{Throughout, all maps of the continuous variables are assumed smooth.}
\begin{equation*}
C_k(u), \Lambda_k(u) \in \End(E), \qquad k = 1, 2.
\end{equation*}

Consider the overdetermined linear problem for
$\Psi \colon \mathcal V \to \Hom(E,H)$:
\begin{align}
\partial_1\Psi(u)
&= \frac{1}{2}\partial_2^2\Psi(u)
   + \partial_2\Psi(u)C_1(u)
   - \Psi(u)\Lambda_1(u),
\label{eq:dc-clean-linear-d1} \\
\Psi(S_3u)
&= \partial_2\Psi(u)C_2(u)
   - \Psi(u)\Lambda_2(u).
\label{eq:dc-clean-linear-S3}
\end{align}

Define the shifted bracket\footnote{The shifted bracket is not anti-symmetric: $[X,Y]_{S_3}(u)\neq -[Y,X]_{S_3}(u)$ in general, since the shift $S_3$ acts on different factors in the two expressions.} 
\begin{align*}
    [X, Y]_{S_3}(u) \defeq X(u)Y(u) - Y(u)X(S_3u)
\end{align*}
and the heat operator
\[\mathcal H\defeq \partial_1-\frac{1}{2}\partial_2^2.\]

Compatibility of the linear problem \eqref{eq:dc-clean-linear-d1}--\eqref{eq:dc-clean-linear-S3} forces the following conditions on $(C_k, \Lambda_k)$. 
\begin{proposition}
\label{prop:dc-clean-diamond}
The system \eqref{eq:dc-clean-linear-d1}--\eqref{eq:dc-clean-linear-S3}
is compatible\footnote{That is, the two computations of
\(\partial_1\Psi(S_3 u)\), obtained by applying the
\(\partial_1\)-flow and the \(S_3\)-shift in opposite
orders, agree for arbitrary local \(\partial_2\)-Cauchy
data. } if and only if
\begin{align}
\partial_2 C_2(u) -[C_1, C_{2}]_{S_3}(u)  
&= 0,
\label{eq:dc-clean-C-diamond} 
\end{align}
\begin{align}
    &\mathcal H \Lambda_2(u) 
+\partial_2\Lambda_1(u)C_2(u)-\partial_2\Lambda_2(u)C_1(S_3u) - [\Lambda_1, \Lambda_2]_{S_3}(u)  =0, \label{eq:dc-clean-L-diamond}
\end{align}
\begin{align}
\mathcal H C_2(u)
+ \partial_2 C_1(u)C_2(u)
-\partial_2C_2(u)C_1(S_3u) +\partial_2\Lambda_2(u)
- [C_1,\Lambda_2]_{S_3}(u)
- [\Lambda_1,C_2]_{S_3}(u)
=0.
\label{eq:dc-clean-mixed-diamond}
\end{align}
\end{proposition}

\begin{proof}
The compatibility equations are obtained by comparing
$\partial_1\Psi(S_3u)$ computed in two ways.  First, we differentiate
\eqref{eq:dc-clean-linear-S3} with respect to $\partial_1$ to obtain 
\begin{align}
\partial_1\Psi(S_3u)
= \partial_1\partial_2\Psi(u)C_2(u)
   +\partial_2\Psi(u)\partial_1C_2(u)
   -\partial_1\Psi(u)\Lambda_2(u)
   -\Psi(u)\partial_1\Lambda_2(u). \label{eq:dc-clean-firstway}
\end{align}
Differentiating \eqref{eq:dc-clean-linear-d1} with respect to $\partial_2$ yields an expression for $\partial_1\partial_2\Psi(u)$:
\begin{align*}
\partial_1\partial_2\Psi(u)
= \tfrac{1}{2}\partial_2^3\Psi(u)
   +\partial_2^2\Psi(u)C_1(u)
   +\partial_2\Psi(u)(\partial_2C_1(u)-\Lambda_1(u))
   -\Psi(u)\partial_2\Lambda_1(u).
\end{align*}
Substituting this into \eqref{eq:dc-clean-firstway} gives
\begin{align}
\partial_1\Psi(S_3u)
&= \frac{1}{2}\partial_2^3\Psi(u)[C_2(u)]
   \notag\\
   &\quad +\partial_2^2\Psi(u)[C_1(u)C_2(u)
      -\frac{1}{2}\Lambda_2(u)] \notag \\
&\quad
   +\partial_2\Psi(u)[(\partial_2C_1(u)-\Lambda_1(u))C_2(u)
      +\partial_1C_2(u)-C_1(u)\Lambda_2(u)] \notag\\
&\quad
   +\Psi(u)[-\partial_2\Lambda_1(u)C_2(u)
      +\Lambda_1(u)\Lambda_2(u)-\partial_1\Lambda_2(u)]. \label{eq:dc-compt-firstway}
\end{align}

For the second computation of $\partial_1 \Psi(S_3 u)$, we evaluate
\eqref{eq:dc-clean-linear-d1} at $S_3u$ 
\begin{align}
\partial_1\Psi(S_3u)
&= \frac{1}{2}\partial_2^2\Psi(S_3u)
   +\partial_2\Psi(S_3u)C_1(S_3u)
   -\Psi(S_3u)\Lambda_1(S_3u). \label{eq:dc-clean-linear-secondway} 
\end{align}
Next, differentiating \eqref{eq:dc-clean-linear-S3} by $\partial_2$ gives
\begin{align}
\partial_2\Psi(S_3u)
&= \partial_2^2\Psi(u)C_2(u)
   +\partial_2\Psi(u)(\partial_2C_2(u)-\Lambda_2(u))
   -\Psi(u)\partial_2\Lambda_2(u).
\label{eq:dc-clean-d2-S3}
\end{align}
Differentiating \eqref{eq:dc-clean-d2-S3} once more gives
\begin{align}
\partial_2^2\Psi(S_3u)
&= \partial_2^3\Psi(u)C_2(u)
   +\partial_2^2\Psi(u)(2\partial_2C_2(u)-\Lambda_2(u))
\notag \\
& \quad +\partial_2\Psi(u)(\partial_2^2C_2(u)-2\partial_2\Lambda_2(u))
   -\Psi(u)\partial_2^2\Lambda_2(u).
\label{eq:dc-clean-d22-S3}
\end{align}
Substituting \eqref{eq:dc-clean-linear-S3},
\eqref{eq:dc-clean-d2-S3}, and \eqref{eq:dc-clean-d22-S3}
into \eqref{eq:dc-clean-linear-secondway} gives
\begin{align}
\partial_1\Psi(S_3u)
&= \frac{1}{2}\partial_2^3\Psi(u)[C_2(u)] \notag\\
&\quad
   +\partial_2^2\Psi(u)[\partial_2C_2(u)
      -\frac{1}{2}\Lambda_2(u)+C_2(u)C_1(S_3u)] \notag\\
&\quad
   +\partial_2\Psi(u)[\tfrac{1}{2}\partial_2^2C_2(u)
      -\partial_2\Lambda_2(u) \notag \\
&\qquad\qquad
      +(\partial_2C_2(u)-\Lambda_2(u))C_1(S_3u)
      -C_2(u)\Lambda_1(S_3u)] \notag \\
&\quad
   +\Psi(u)[-\frac{1}{2}\partial_2^2\Lambda_2(u)
      -\partial_2\Lambda_2(u)C_1(S_3u)
      +\Lambda_2(u)\Lambda_1(S_3u)]. \label{eq:dc-compt-secondway}
\end{align}

If the system is compatible, then \eqref{eq:dc-compt-firstway} and \eqref{eq:dc-compt-secondway} must agree for all choices of arbitrary local $\partial_2$-Cauchy data.  Note that $\partial_2^3\Psi$ coefficients agree automatically. Equating the coefficients of $\partial_2^2\Psi$ gives
\eqref{eq:dc-clean-C-diamond}.  Equating the coefficients of
$\partial_2\Psi$ gives \eqref{eq:dc-clean-mixed-diamond}.
Equating the coefficients of $\Psi$ gives
\eqref{eq:dc-clean-L-diamond}.  

Conversely, if
\eqref{eq:dc-clean-C-diamond}--\eqref{eq:dc-clean-mixed-diamond}
hold, then \eqref{eq:dc-compt-firstway} and \eqref{eq:dc-compt-secondway} have identical coefficients and therefore agree for arbitrary local $\partial_2$-Cauchy data.
\end{proof}

We call
\eqref{eq:dc-clean-C-diamond}--\eqref{eq:dc-clean-mixed-diamond}
the \emph{parabolic diamond equations}.

\begin{remark}\label{rem:dc-clean-scaling}
The parabolic framework arises from the semi-discrete theory of
Chapter~\ref{ch:the-semi-discrete-framework} under the scaling
$T^{(\epsilon)} = e^{\epsilon\partial_2}$, with the semi-discrete
continuous derivative (written $\partial_1^{(\epsilon)}$ to
distinguish it from the parabolic $\partial_1$) decomposed as
$\partial_1^{(\epsilon)} = \epsilon^{-1}\partial_2 + \partial_1$.
Expand the semi-discrete edge weights as
\[
C_1^{(\epsilon)} = \epsilon^{-2}I + \epsilon^{-1}C_1, \qquad
\Lambda_1^{(\epsilon)} = \epsilon^{-2}I + \epsilon^{-1}C_1
  + \Lambda_1,
\]
\[
C_2^{(\epsilon)} = \epsilon^{-1}C_2, \qquad
\Lambda_2^{(\epsilon)} = \epsilon^{-1}C_2 + \Lambda_2.
\]
The leading-order terms in the semi-discrete diamond equations
then recover the parabolic diamond equations
\eqref{eq:dc-clean-C-diamond}--\eqref{eq:dc-clean-mixed-diamond}.
A direct scaling from the fully discrete theory
(Chapter~\ref{ch:the-discrete-framework}) also exists: scaling
$T^{(\epsilon)} = e^{\epsilon\partial_2}$ and a non-distinguished
shift $S_1^{(\epsilon)} = e^{p\epsilon\partial_2
+ p(1-p)\epsilon^2\partial_1}$ for a parameter
$p \notin \{0,1\}$, the parabolic linear problem emerges at
order $\epsilon^2$.
\end{remark}

\section{Darboux Transformations}
\label{sec:dc-clean-darboux}

The Darboux transformation produces, from a solution
$(C_k, \Lambda_k)$ of the parabolic diamond equations together with wave
functions and a compatible kernel, new edge weights
$(C_k, \widehat \Lambda_k)$ satisfying the same parabolic diamond equations
of Proposition~\ref{prop:dc-clean-diamond}.  We now develop this
construction.

Fix a solution $(C_1,\Lambda_1,C_2,\Lambda_2)$ of the
parabolic diamond equations.  Let
$\Psi \colon \mathcal V \to \Hom(E,H)$ solve
\eqref{eq:dc-clean-linear-d1}--\eqref{eq:dc-clean-linear-S3}.
The adjoint linear problem\footnote{Its compatibility conditions are
the diamond equations of
Proposition~\ref{prop:dc-clean-diamond}.} for
$\Phi \colon \mathcal V \to \Hom(H,E)$ is
\begin{align}
\partial_1\Phi(u)
&= -\frac{1}{2}\partial_2^2\Phi(u)
   +C_1(u)\partial_2\Phi(u)
   +(\Lambda_1(u)+\partial_2C_1(u))\Phi(u),
\label{eq:dc-clean-adjoint-d1} \\
\Phi(u)
&= -C_2(u)\partial_2\Phi(S_3u)
   -(\Lambda_2(u)+\partial_2C_2(u))\Phi(S_3u).
\label{eq:dc-clean-adjoint-S3}
\end{align}

\begin{definition}
\label{def:dc-clean-dressing}
We say that $K \colon \mathcal V \to \End(H)$ is \emph{dressing
compatible} with $(\Psi,\Phi)$ if: 
\begin{enumerate}[label=\arabic*., leftmargin=*]
\item For every $u \in \mathcal{V}$,
\begin{align}
\partial_2K(u)
&= \Psi(u)\Phi(u),
\label{eq:dc-clean-K-d2} \\
\partial_1K(u)
&= \frac{1}{2}\bigl(\partial_2\Psi(u)\Phi(u)
   -\Psi(u)\partial_2\Phi(u)\bigr)
   +\Psi(u)C_1(u)\Phi(u),
\label{eq:dc-clean-K-d1} \\
K(S_3u)-K(u)
&= \Psi(u)C_2(u)\Phi(S_3u),
\label{eq:dc-clean-K-S3}
\end{align}
\item There exists $z \in \mathbb{F}$ such that the resolvent
exists at every vertex:
\[
R(u) \defeq (I - zK(u))^{-1} \in \End(H),
\qquad u \in \mathcal{V}.
\]
\end{enumerate}
\end{definition}

The dressed waves inherit the linear problem with corrected
\(\Lambda\)-weights. 
\begin{proposition}
\label{prop:dc-clean-dressed-waves}
Let $K$ be dressing compatible with $(\Psi, \Phi)$, with resolvent $R(u) = (I-zK(u))^{-1}$. Define the dressed observable and dressed waves by 
\[
\mathcal A(u) \defeq z\Phi(u)R(u)\Psi(u),
\]
\[
\widehat\Psi(u) \defeq R(u)\Psi(u), \qquad
\widehat\Phi(u) \defeq \Phi(u)R(u).
\]
Then $\widehat K(u) \defeq K(u)R(u)$ is dressing compatible
with $(\widehat\Psi, \widehat\Phi)$ which satisfy the
linear problem
\eqref{eq:dc-clean-linear-d1}--\eqref{eq:dc-clean-linear-S3}
and the adjoint linear problem
\eqref{eq:dc-clean-adjoint-d1}--\eqref{eq:dc-clean-adjoint-S3}
respectively, with data $(C, \widehat\Lambda)$, where
\begin{align}
\widehat\Lambda_1(u)
&\defeq \Lambda_1(u)+\partial_2\mathcal A(u)+[\mathcal A(u),C_1(u)],
\label{eq:dc-clean-L1-dressed-general} \\
\widehat\Lambda_2(u)
&\defeq \Lambda_2(u)+[\mathcal A, C_2]_{S_3}(u)
\label{eq:dc-clean-L2-dressed-general}
\end{align}
\end{proposition}
The proof of Proposition~\ref{prop:dc-clean-dressed-waves} relies on the
following resolvent identities.  For the remainder of this
subsection, we suppress functional dependence on $z$.
\begin{lemma}
    For every $u \in \mathcal V$ the following resolvent identities hold: 
    \begin{align}
\partial_2R(u)
&= zR(u)\Psi(u)\Phi(u)R(u),
\label{eq:dc-clean-R-d2} \\
\partial_1R(u)
&= zR(u)\Bigl[\tfrac{1}{2}\bigl(\partial_2\Psi(u)\Phi(u)
   -\Psi(u)\partial_2\Phi(u)\bigr) \notag\\
&\qquad\qquad \quad +\Psi(u)C_1(u)\Phi(u)\Bigr]R(u),
\label{eq:dc-clean-R-d1} \\
R(S_3u)-R(u)
&= zR(u)\Psi(u)C_2(u)\Phi(S_3u)R(S_3u).
\label{eq:dc-clean-R-S3}
\end{align}
\end{lemma}
\begin{proof}
Differentiating (resp.\ evaluating) $(I-zK(u))R(u) = I$ gives the
general resolvent identities
$\partial_i R(u) = zR(u)(\partial_i K(u))R(u)$ and
$R(w) - R(v) = zR(v)(K(w) - K(v))R(w)$, as in
Lemma~\ref{lem:sc-resolvent-ids}.  Substituting
\eqref{eq:dc-clean-K-d2} and \eqref{eq:dc-clean-K-d1} proves
\eqref{eq:dc-clean-R-d2} and \eqref{eq:dc-clean-R-d1};
substituting \eqref{eq:dc-clean-K-S3} proves
\eqref{eq:dc-clean-R-S3}.
\end{proof}
\begin{proof}[Proof of Proposition \ref{prop:dc-clean-dressed-waves}]
Three preliminary identities are used throughout the proof.
Differentiating $\widehat\Psi=R\Psi$ and applying
\eqref{eq:dc-clean-R-d2} gives
\begin{align}
\partial_2\widehat\Psi(u)
&= R(u)\partial_2\Psi(u)+\widehat\Psi(u)\mathcal A(u).
\label{eq:dc-clean-d2-hatPsi}
\end{align}
The adjoint counterpart is
\begin{align}
\partial_2\widehat\Phi(u)
&= \partial_2\Phi(u)R(u)+\mathcal A(u)\widehat\Phi(u).
\label{eq:dc-clean-d2-hatPhi}
\end{align}
Differentiating $\mathcal A=z\Phi R\Psi$ and substituting
\eqref{eq:dc-clean-R-d2} gives
\begin{align}
\partial_2\mathcal A(u)
&= z\partial_2\Phi(u)R(u)\Psi(u)
   +\mathcal A(u)^2
   +z\Phi(u)R(u)\partial_2\Psi(u).
\label{eq:dc-clean-d2A-expanded}
\end{align}

For the $S_3$-equation \eqref{eq:dc-clean-linear-S3} of $\widehat\Psi$, write
$\widehat\Psi(S_3u)=R(S_3u)\Psi(S_3u)$ in the form
\begin{align*}
\widehat\Psi(S_3u)
&= R(u)\Psi(S_3u)+(R(S_3u)-R(u))\Psi(S_3u).
\end{align*}
Substituting the linear problem
\eqref{eq:dc-clean-linear-S3} for $\Psi(S_3u)$ in the first
term, and the resolvent identity \eqref{eq:dc-clean-R-S3} in
the second, gives
\begin{align*}
\widehat\Psi(S_3u)
= R(u)\bigl[\partial_2\Psi(u)C_2(u)
   -\Psi(u)\Lambda_2(u)\bigr]
+zR(u)\Psi(u)C_2(u)\Phi(S_3u)R(S_3u)\Psi(S_3u).
\end{align*}
The last term is $\widehat \Psi(u)C_2(u) \mathcal A(S_3u)$. Moreover, rearranging \eqref{eq:dc-clean-d2-hatPsi} in the form $R(u)\partial_2 \Psi(u) = \partial_2 \widehat \Psi(u) - \widehat \Psi(u) \mathcal A(u)$ and collecting, 
\begin{align*}
    \widehat \Psi(S_3u) &= \partial_2 \widehat \Psi(u)C_2(u) - \widehat \Psi(u) [\mathcal A(u)C_2(u) + \Lambda_2(u) - C_2(u)\mathcal A(S_3u)] \\
    &= \partial_2\widehat \Psi(u)C_2(u) - \widehat \Psi(u) \widehat \Lambda_2(u),
\end{align*}
where the last equality follows from inserting the definition of
$\widehat\Lambda_2$.

For the $\partial_1$-equation \eqref{eq:dc-clean-linear-d1} of $\widehat\Psi$, differentiating
$\widehat\Psi=R\Psi$ by $\partial_1$ and
substituting \eqref{eq:dc-clean-linear-d1} for $\partial_1\Psi$
and \eqref{eq:dc-clean-R-d1} for $\partial_1R$ gives
\begin{align*}
\partial_1\widehat\Psi
&= \frac{1}{2}R\partial_2^2\Psi
   +R\partial_2\Psi C_1
   -\widehat\Psi\Lambda_1
   +\frac{1}{2}R\partial_2\Psi\mathcal A
   -\frac{1}{2}z\widehat\Psi \partial_2\Phi R\Psi
   +\widehat\Psi C_1\mathcal A.
\end{align*}
The first three terms come from $R(\partial_1\Psi)$; the last
three from $(\partial_1R)\Psi$.
On the other hand, differentiating
\eqref{eq:dc-clean-d2-hatPsi} once more gives
\begin{align*}
\frac{1}{2}\partial_2^2\widehat\Psi
&= \frac{1}{2}R\partial_2^2\Psi
   +\frac{1}{2}R\partial_2\Psi\mathcal A
   +\frac{1}{2}z\widehat\Psi \partial_2\Phi R\Psi
   +z\widehat\Psi \Phi R\partial_2\Psi
   +\widehat\Psi\mathcal A^2,
\end{align*}
obtained by applying the Leibniz rule to both terms of
\eqref{eq:dc-clean-d2-hatPsi}, using \eqref{eq:dc-clean-R-d2}
to differentiate~$R$, and expanding $\partial_2\mathcal A$ via
\eqref{eq:dc-clean-d2A-expanded}.  
Evaluating $\partial_1 \widehat \Psi - \frac{1}{2}\partial_2^2 \widehat \Psi$ gives
\begin{align*}
\partial_1\widehat\Psi
-\frac{1}{2}\partial_2^2\widehat\Psi
&= R\partial_2\Psi C_1
   -\widehat\Psi\Lambda_1
   +\widehat\Psi C_1\mathcal A \\
&\quad
   -z\widehat\Psi \partial_2\Phi R\Psi
   -z\widehat\Psi \Phi R\partial_2\Psi
   -\widehat\Psi\mathcal A^2.
\end{align*}
The last line is $-\widehat \Psi \partial_2 \mathcal A$
 by \eqref{eq:dc-clean-d2A-expanded},
leaving
\begin{align*}
\partial_1\widehat\Psi
-\frac{1}{2}\partial_2^2\widehat\Psi
&= R\partial_2\Psi C_1
   +\widehat\Psi\bigl[C_1\mathcal A
      -\Lambda_1-\partial_2\mathcal A\bigr] \\
&= \bigl(R\partial_2\Psi+\widehat\Psi\mathcal A\bigr)C_1
   -\widehat\Psi\bigl[\Lambda_1+[\mathcal A,C_1]
      +\partial_2\mathcal A\bigr] \\
&= \partial_2\widehat\Psi C_1
   -\widehat\Psi\widehat\Lambda_1,
\end{align*}
where the second line uses $C_1\mathcal A = \mathcal A C_1 -[\mathcal A, C_1]$ and
rearranges the $C_1$ terms, and
the third uses \eqref{eq:dc-clean-d2-hatPsi} and the
definition of $\widehat\Lambda_1$.

For the $S_3$-equation \eqref{eq:dc-clean-adjoint-S3} of $\widehat\Phi$,
evaluate \eqref{eq:dc-clean-d2-hatPhi} at $S_3u$:
\[
\partial_2\widehat\Phi(S_3u)
=\partial_2\Phi(S_3u)R(S_3u)
+\mathcal A(S_3u)\widehat\Phi(S_3u).
\]
Then expand the right-hand side of the dressed version of
\eqref{eq:dc-clean-adjoint-S3} to give
\begin{align*}
&-C_2\partial_2\widehat\Phi(S_3u)
-(\widehat\Lambda_2+\partial_2C_2)\widehat\Phi(S_3u) \\
&\quad
=-C_2\partial_2\Phi(S_3u)R(S_3u)
  -C_2\mathcal A(S_3u)\widehat\Phi(S_3u)
  -(\Lambda_2+\partial_2C_2)\widehat\Phi(S_3u) \\
&\qquad\quad
  -\mathcal A(u)C_2\widehat\Phi(S_3u)
  +C_2\mathcal A(S_3u)\widehat\Phi(S_3u).
\end{align*}
The two $C_2\mathcal A(S_3u)\widehat\Phi(S_3u)$ terms cancel.
The bracketed term, after right-multiplication by $R(S_3u)$,
\[
\bigl[-C_2\partial_2\Phi(S_3u)
-(\Lambda_2+\partial_2C_2)\Phi(S_3u)\bigr]R(S_3u),
\]
equals $\Phi(u)R(S_3u)$ by the adjoint problem
\eqref{eq:dc-clean-adjoint-S3}, so the expression becomes
\begin{align*}
&\Phi(u)R(S_3u)
  -\mathcal A(u)C_2(u)\widehat \Phi(S_3u).
\end{align*}
Left-multiplying the resolvent identity
\eqref{eq:dc-clean-R-S3} by $\Phi(u)$ gives $\Phi(u)R(S_3u)-\Phi(u)R(u)
=\mathcal A(u)C_2\widehat\Phi(S_3u)$,
and so substituting into the expression above gives $\Phi(u)R(u)=\widehat\Phi(u)$.

For the $\partial_1$-equation \eqref{eq:dc-clean-adjoint-d1} of $\widehat\Phi$, differentiating
$\widehat\Phi=\Phi R$ by $\partial_1$ and
substituting \eqref{eq:dc-clean-adjoint-d1} for $\partial_1\Phi$
and \eqref{eq:dc-clean-R-d1} for $\partial_1R$ gives
\begin{align*}
\partial_1\widehat\Phi
&= -\frac{1}{2}\partial_2^2\Phi R
   +C_1\partial_2\Phi R
   +(\Lambda_1+\partial_2C_1)\widehat\Phi \\
&\quad
   +\frac{1}{2}z\widehat\Phi\partial_2\Psi \widehat \Phi
   -\frac{1}{2}z\widehat\Phi\Psi\partial_2\Phi R
   +\mathcal A C_1\widehat\Phi,
\end{align*}
where the first line is $(\partial_1\Phi)R$ and the second
is $\Phi(\partial_1R)$.
Separately, differentiating \eqref{eq:dc-clean-d2-hatPhi}
gives
\begin{align*}
\frac{1}{2}\partial_2^2\widehat\Phi
&= \frac{1}{2}(\partial_2^2\Phi)R
   +\frac{1}{2}z(\partial_2\Phi)R\Psi\widehat\Phi
   +\frac{1}{2}(\partial_2\mathcal A)\widehat\Phi
   +\frac{1}{2}\mathcal A(\partial_2\Phi)R
   +\frac{1}{2}\mathcal A^2\widehat\Phi,
\end{align*}
obtained by applying $\partial_2$ to both terms of
\eqref{eq:dc-clean-d2-hatPhi} and using
\eqref{eq:dc-clean-R-d2} in $\partial_2[(\partial_2\Phi)R]
=(\partial_2^2\Phi)R+z(\partial_2\Phi)R\Psi\widehat\Phi$.
Adding the two displays cancels the
$\frac{1}{2}(\partial_2^2\Phi)R$ terms.  The
$z\widehat\Phi\Psi(\partial_2\Phi)R$ contribution from
$\Phi(\partial_1R)$ cancels against the
$\mathcal A(\partial_2\Phi)R$ term from
$\frac{1}{2}\partial_2^2\widehat\Phi$.  The remaining three quadratic
terms
$\frac{1}{2}z\widehat\Phi(\partial_2\Psi)\widehat\Phi
+\frac{1}{2}z(\partial_2\Phi)R\Psi\widehat\Phi
+\frac{1}{2}\mathcal A^2\widehat\Phi$
are, by \eqref{eq:dc-clean-d2A-expanded}, equal to
$\frac{1}{2}(\partial_2\mathcal A)\widehat\Phi$.
Together with the explicit
$\frac{1}{2}(\partial_2\mathcal A)\widehat\Phi$ term, this gives
\begin{align*}
(\partial_1+\frac{1}{2}\partial_2^2)\widehat\Phi
&= C_1\partial_2\Phi R
   +(\Lambda_1+\partial_2C_1)\widehat\Phi
    +\partial_2\mathcal A\widehat\Phi +\mathcal A C_1\widehat\Phi
 \\
&= C_1[\partial_2\Phi R+\mathcal A\widehat\Phi ]
   +(
   \Lambda_1+\partial_2C_1+\partial_2\mathcal A+[\mathcal A,C_1])\widehat\Phi \\
&= C_1\partial_2\widehat\Phi
   +(\widehat\Lambda_1+\partial_2C_1)\widehat\Phi,
\end{align*}
where the second line uses
$\mathcal A C_1
=C_1\mathcal A
+[\mathcal A,C_1]$,
and the third applies \eqref{eq:dc-clean-d2-hatPhi} and the
definition of $\widehat\Lambda_1$.

\noindent\textit{The dressed kernel.}
Since $K(u)$ and
$R(u)$ commute at each vertex, $\widehat K(u) = R(u)K(u)$.
The identity
$(I - \eta\widehat K(u))(I - zK(u)) = I - (z+\eta)K(u)$
shows that $(I - \eta\widehat K(u))^{-1}$ exists whenever
$(I - (z+\eta)K(u))^{-1}$ does; in particular at
$\eta = -z$, where the condition is vacuous.  When $z = 0$,
$R = I$ and $\widehat K = K$, so the kernel identities
reduce to those of $K$ itself.  For $z \neq 0$, it remains
to verify the three kernel identities.

\noindent\textit{The $\partial_2$-condition.}
Differentiating $\widehat K(u) = z^{-1}(R(u) - I)$
and substituting \eqref{eq:dc-clean-R-d2}:
\[
\partial_2\widehat K(u)
= z^{-1}\partial_2 R(u)
= R(u)\Psi(u)\Phi(u)R(u)
= \widehat\Psi(u)\widehat\Phi(u).
\]

\noindent\textit{The $\partial_1$-condition.}
Differentiating $\widehat K(u) = z^{-1}(R(u) - I)$
and substituting \eqref{eq:dc-clean-R-d1}:
\begin{align*}
\partial_1\widehat K(u)
= \tfrac12\bigl(R(u)\partial_2\Psi(u)\Phi(u)R(u)
- R(u)\Psi(u)\partial_2\Phi(u)R(u)\bigr)
+ R(u)\Psi(u)C_1(u)\Phi(u)R(u).
\end{align*}
The last term equals
$\widehat\Psi(u)C_1(u)\widehat\Phi(u)$.
For the first term,
\eqref{eq:dc-clean-d2-hatPsi} gives
$R(u)\partial_2\Psi(u)
= \partial_2\widehat\Psi(u) - \widehat\Psi(u)\mathcal A(u)$,
and for the second term \eqref{eq:dc-clean-d2-hatPhi} gives
$\partial_2\Phi(u)R(u)
= \partial_2\widehat\Phi(u)
  - \mathcal A(u)\widehat\Phi(u)$.
Substituting:
\begin{align*}
&\tfrac12\bigl(R\partial_2\Psi\Phi R
  - R\Psi\partial_2\Phi R\bigr) \\
&\quad= \tfrac12\bigl(
  \partial_2\widehat\Psi\widehat\Phi
  - \widehat\Psi\mathcal A\widehat\Phi
  - \widehat\Psi\partial_2\widehat\Phi
  + \widehat\Psi\mathcal A\widehat\Phi\bigr) \\
&\quad= \tfrac12\bigl(
  \partial_2\widehat\Psi\widehat\Phi
  - \widehat\Psi\partial_2\widehat\Phi\bigr),
\end{align*}
where the $\pm\widehat\Psi\mathcal A\widehat\Phi$
terms cancel.  Collecting,
\[
\partial_1\widehat K(u)
= \tfrac12\bigl(\partial_2\widehat\Psi(u)\widehat\Phi(u)
  - \widehat\Psi(u)\partial_2\widehat\Phi(u)\bigr)
  + \widehat\Psi(u)C_1(u)\widehat\Phi(u).
\]

\noindent\textit{The $S_3$-condition.}
Since $\widehat K(S_3u) - \widehat K(u)
= z^{-1}(R(S_3u) - R(u))$,
\eqref{eq:dc-clean-R-S3} gives
\[
\widehat K(S_3u) - \widehat K(u)
= R(u)\Psi(u)C_2(u)\Phi(S_3u)R(S_3u)
= \widehat\Psi(u)C_2(u)\widehat\Phi(S_3u).
\]
\end{proof}

Proposition~\ref{prop:dc-clean-dressed-waves} shows that the dressed
waves satisfy the same linear problem with corrected
\(\Lambda\)-weights. Theorem~\ref{thm:dc-clean-general-darboux} next establishes that the corrected data \((C_k,\widehat\Lambda_k)\) satisfy the full diamond system. Moreover, the dressed observable $\mathcal A(u) = z \Phi(u) R(u)\Psi(u)$ itself satisfies an equation \eqref{eq:dc-clean-A-relative-Lambda-shifted-commutator}.
\begin{theorem}
\label{thm:dc-clean-general-darboux}
Under the hypotheses of
Proposition~\ref{prop:dc-clean-dressed-waves}, the corrected data
$(C_k,\widehat\Lambda_k)$
satisfies the parabolic diamond equations \eqref{eq:dc-clean-C-diamond}--\eqref{eq:dc-clean-mixed-diamond}. In particular, the dressed observable \(\mathcal A\) satisfies 
\begin{align}
&\mathcal H\Delta_2(u)
-\partial_2\Delta_2(u)C_1(S_3u)
+\partial_2\Delta_1(u)C_2(u)
\notag \\
&\qquad - 
[\Lambda_1,\Delta_2]_{S_3}(u)
-
[\Delta_1,\Lambda_2]_{S_3}(u)
-
[\Delta_1,\Delta_2]_{S_3}(u) =0, 
\label{eq:dc-clean-A-relative-Lambda-shifted-commutator}
\end{align}
where
\begin{align}
\Delta_1(u)
&\defeq
\partial_2\mathcal A(u)
+[\mathcal A, C_1](u), \qquad
\Delta_2(u)
\defeq
[\mathcal A, C_2]_{S_3}(u).
\label{eq:dc-clean-Delta-A}
\end{align}
\end{theorem}
The proof of Theorem~\ref{thm:dc-clean-general-darboux} reduces to three
identities.  Throughout, write
\[
X(u)\defeq z\partial_2\Phi(u)R(u)\Psi(u),
\qquad
Y(u)\defeq z\Phi(u)R(u)\partial_2\Psi(u),
\]
so that $\partial_2\mathcal A = X + \mathcal A^2 + Y$ by
\eqref{eq:dc-clean-d2A-expanded}.  We suppress the dependence on
$u\in\mathcal V$ throughout; all quantities are evaluated at the same
vertex unless otherwise indicated.  A subscript $S$ denotes evaluation
at $S_3u$.

The first identity describes the parabolic evolution of the dressed
observable $\mathcal A$ under the heat operator $\mathcal H$.
\begin{lemma}
\label{lem:dc-clean-general-A-heat}
The dressed observable satisfies
\begin{align}
\mathcal H\mathcal A
&=
\mathcal A C_1\mathcal A
-\mathcal A\partial_2\mathcal A 
+[\Lambda_1,\mathcal A]
+(\partial_2C_1)\mathcal A
-\partial_2X
+C_1X
+YC_1.
\label{eq:dc-clean-general-A-heat}
\end{align}
\end{lemma}

\begin{proof}
The non-commutative heat product rule,
\[
\mathcal H(fg) = (\mathcal Hf)g + f(\mathcal Hg)
  - \partial_2 f \partial_2 g,
\]
applied to the factorization $\mathcal A = z\Phi\widehat\Psi$, with
$\widehat\Psi = R\Psi$, gives
\begin{align}
\mathcal H\mathcal A
&= z(\mathcal H\Phi)\widehat\Psi
  + z\Phi(\mathcal H\widehat\Psi)
  - z(\partial_2\Phi)(\partial_2\widehat\Psi).
\label{eq:dc-clean-HA-product-rule}
\end{align}
The adjoint equation \eqref{eq:dc-clean-adjoint-d1} gives
$\mathcal H\Phi = -\partial_2^2\Phi + C_1\partial_2\Phi
+ (\Lambda_1+\partial_2C_1)\Phi$.  Substituting into the first term
of \eqref{eq:dc-clean-HA-product-rule} and combining with the third
gives
\begin{align*}
z(\mathcal H\Phi)\widehat\Psi
  - z(\partial_2\Phi)(\partial_2\widehat\Psi)
&= -\partial_2(z\partial_2\Phi\widehat\Psi)
  + C_1 X + (\Lambda_1+\partial_2C_1)\mathcal A \\
&= -\partial_2 X + C_1 X + (\Lambda_1+\partial_2C_1)\mathcal A.
\end{align*}

For the remaining contribution, the dressed form of
\eqref{eq:dc-clean-linear-d1} from
Proposition~\ref{prop:dc-clean-dressed-waves} gives
\[
\mathcal H\widehat\Psi
= \partial_2\widehat\Psi C_1
  - \widehat\Psi\bigl(\Lambda_1+\partial_2\mathcal A
    +[\mathcal A, C_1]\bigr).
\]
Left-multiplying by $z\Phi$ and using $z\Phi\widehat\Psi = \mathcal A$,
\[
z\Phi(\mathcal H\widehat\Psi)
= z\Phi\partial_2\widehat\Psi C_1
  - \mathcal A\bigl(\Lambda_1+\partial_2\mathcal A
    +[\mathcal A, C_1]\bigr).
\]
By \eqref{eq:dc-clean-d2-hatPsi},
$z\Phi\partial_2\widehat\Psi = Y+\mathcal A^2$.  The
$\mathcal A^2C_1$ term from $z\Phi\partial_2\widehat\Psi C_1$
cancels against $-\mathcal A^2C_1$ from the bracket,
and collecting by type gives \eqref{eq:dc-clean-general-A-heat}.
\end{proof}

The remaining two identities work together: the second shows that the
$\Lambda$-diamond residual factors through a single quantity
$\mathcal E$, and the third shows that $\mathcal E$ vanishes.

Substituting $\widehat\Lambda_k = \Lambda_k + \Delta_k$ into the
$\Lambda$-diamond \eqref{eq:dc-clean-L-diamond} and subtracting the
seed equation, the dressed $\Lambda$-diamond holds if and only if
\begin{align}
\mathcal H\Delta_2
+ \partial_2\Delta_1 C_2
- \partial_2\Delta_2 C_{1,S}
- [\Lambda_1,\Delta_2]_{S_3}
- [\Delta_1,\Lambda_2]_{S_3}
- [\Delta_1,\Delta_2]_{S_3}
= 0.
\label{eq:dc-clean-Lambda-residual}
\end{align}

\begin{lemma}
\label{lem:dc-clean-general-Lambda-residual}
The left-hand side of \eqref{eq:dc-clean-Lambda-residual} factors as
\[
\partial_2\mathcal E - C_1\mathcal E + \mathcal E C_{1,S},
\]
where
\begin{align}
\mathcal E
&\defeq
\mathcal AC_2\mathcal A_S
- [\mathcal A, \Lambda_2]_{S_3}
+ (\partial_2C_2)\mathcal A_S
+ YC_2 + C_2X_S.
\label{eq:dc-clean-E-definition}
\end{align}
\end{lemma}

\begin{lemma}
\label{lem:dc-clean-shifted-defect}
The quantity $\mathcal E$ from
Lemma~\ref{lem:dc-clean-general-Lambda-residual} vanishes identically.
\end{lemma}

\begin{proof}[Proof of Lemma~\ref{lem:dc-clean-shifted-defect}]
The linear $S_3$-equation \eqref{eq:dc-clean-linear-S3} gives
\[
\partial_2\Psi C_2
=
\Psi_S+\Psi\Lambda_2.
\]
Multiplying on the left by $z\Phi R$ gives
\[
YC_2
=
z\Phi R\Psi_S+\mathcal A\Lambda_2.
\]
The adjoint $S_3$-equation \eqref{eq:dc-clean-adjoint-S3} gives
\[
C_2\partial_2\Phi_S
=
-\Phi-(\Lambda_2+\partial_2C_2)\Phi_S.
\]
Multiplying on the right by $zR_S\Psi_S$ gives
\[
C_2X_S
=
-z\Phi R_S\Psi_S
-(\Lambda_2+\partial_2C_2)\mathcal A_S.
\]
Substituting both into \eqref{eq:dc-clean-E-definition}, the
$\mathcal A\Lambda_2$ term from $YC_2$ and the
$-\Lambda_2\mathcal A_S$ term from $C_2X_S$ cancel the shifted
commutator $-[\mathcal A,\Lambda_2]_{S_3}$, and the
$-(\partial_2C_2)\mathcal A_S$ term from $C_2X_S$ cancels the
remaining linear term, leaving
\[
z\Phi R\Psi_S
-z\Phi R_S\Psi_S
+\mathcal AC_2\mathcal A_S.
\]
The resolvent shift identity \eqref{eq:dc-clean-R-S3} gives
$R_S-R = zR\Psi C_2\Phi_S R_S$, so
\[
z\Phi(R-R_S)\Psi_S
=
-z\Phi(zR\Psi C_2\Phi_S R_S)\Psi_S
=
-\mathcal AC_2\mathcal A_S,
\]
and the three terms cancel.
Hence $\mathcal E=0$.
\end{proof}

\begin{proof}[Proof of Lemma~\ref{lem:dc-clean-general-Lambda-residual}]
Write $\mathcal P_\Lambda$ for the left-hand side of
\eqref{eq:dc-clean-Lambda-residual}, and define the linear operator
\[
L(Z)\defeq \partial_2 Z - C_1 Z + Z C_{1,S}.
\]
The proof expresses the difference
$\mathcal R \defeq \mathcal P_\Lambda - L(\mathcal E)$ as a linear
combination of quantities that vanish by the identities of the
preceding subsection.

The first pair encodes
the decomposition \eqref{eq:dc-clean-d2A-expanded}:
\begin{align*}
\mathcal D_{\partial\mathcal A}
&\defeq
\partial_2\mathcal A - X - \mathcal A^2 - Y,
\\
\mathcal D_{\partial\mathcal A,S}
&\defeq
\partial_2\mathcal A_S - X_S - \mathcal A_S^2 - Y_S.
\end{align*}
The third encodes the $C$-diamond \eqref{eq:dc-clean-C-diamond}:
\[
\mathcal D_C
\defeq
\partial_2 C_2 - C_1 C_2 + C_2 C_{1,S}.
\]
The fourth and fifth encode the $\mathcal A$-heat identity
(Lemma~\ref{lem:dc-clean-general-A-heat}) and its shifted evaluation:
\begin{align*}
\mathcal D_H
&\defeq
\mathcal H\mathcal A
+ \partial_2 X - C_1 X - YC_1
- (\Lambda_1 + \partial_2C_1)\mathcal A
+ \mathcal A\Lambda_1 \\
& \quad + \mathcal A \partial_2\mathcal A
- \mathcal AC_1\mathcal A,
\\
\mathcal D_{H,S}
&\defeq
\mathcal H\mathcal A_S
+ \partial_2 X_S - C_{1,S} X_S - Y_SC_{1,S}
- (\Lambda_{1,S} + \partial_2C_{1,S})\mathcal A_S
+ \mathcal A_S\Lambda_{1,S}
\\
&\quad
+ \mathcal A_S \partial_2\mathcal A_S
- \mathcal A_SC_{1,S}\mathcal A_S.
\end{align*}
The sixth encodes the mixed diamond \eqref{eq:dc-clean-mixed-diamond}:
\begin{align*}
   \mathcal D_M
\defeq
\mathcal HC_2
- \bigl[
  (\partial_2C_2)C_{1,S}
  - C_2\Lambda_{1,S}
  + (\Lambda_1 - \partial_2C_1)C_2
  + C_1\Lambda_2
  - \partial_2\Lambda_2
  - \Lambda_2C_{1,S}
\bigr].
\end{align*}
All six defects vanish identically: $\mathcal D_{\partial\mathcal A}$
and $\mathcal D_{\partial\mathcal A,S}$ by
\eqref{eq:dc-clean-d2A-expanded}, $\mathcal D_C$ by
\eqref{eq:dc-clean-C-diamond}, $\mathcal D_H$ and
$\mathcal D_{H,S}$ by Lemma~\ref{lem:dc-clean-general-A-heat},
and $\mathcal D_M$ by \eqref{eq:dc-clean-mixed-diamond}.  Since
the seed data are smooth, the $\partial_2$-derivatives of all six
defects vanish as well.

We claim the remainder $\mathcal R$ decomposes as
\begin{align}
\mathcal R
={}&
(\partial_2\mathcal D_{\partial\mathcal A})C_2
+ \mathcal D_HC_2
+ \mathcal A\mathcal D_M
- \mathcal D_M\mathcal A_S
- C_2\mathcal D_{H,S}
\notag\\
&- \mathcal A\mathcal D_C\mathcal A_S
- (\partial_2\mathcal D_C)\mathcal A_S
- Y\mathcal D_C
- \mathcal D_CX_S
\notag\\
&- C_1\mathcal D_{\partial\mathcal A}C_2
+ C_2\mathcal D_{\partial\mathcal A,S}C_{1,S}
- (\partial_2\mathcal A)\mathcal D_C.
\label{eq:lambda-gap-defect-factorization}
\end{align}
Since every term on the right-hand side is linear in a vanishing
defect, this identity gives $\mathcal R = 0$.  The remainder of the
proof verifies \eqref{eq:lambda-gap-defect-factorization} by
expanding both sides and matching term by term.

We first expand $\mathcal R$ without using any of the above defects.
The non-commutative heat product rule applied to
$\Delta_2 = \mathcal AC_2 - C_2\mathcal A_S$ gives
\begin{align*}
\mathcal H\Delta_2
=
(\mathcal H\mathcal A)C_2
+ \mathcal A(\mathcal HC_2)
- (\partial_2\mathcal A)(\partial_2C_2)
- (\mathcal HC_2)\mathcal A_S
- C_2(\mathcal H\mathcal A_S)
+ (\partial_2C_2)(\partial_2\mathcal A_S),
\end{align*}
and the ordinary Leibniz rule gives
\begin{align*}
\partial_2\Delta_2
&=
(\partial_2\mathcal A)C_2
+ \mathcal A(\partial_2C_2)
- (\partial_2C_2)\mathcal A_S
- C_2(\partial_2\mathcal A_S),
\\
\partial_2\Delta_1
&=
\partial_2^2\mathcal A
+ (\partial_2\mathcal A)C_1
+ \mathcal A(\partial_2C_1)
- (\partial_2C_1)\mathcal A
- C_1(\partial_2\mathcal A).
\end{align*}
Recall from \eqref{eq:dc-clean-E-definition} that
\[
\mathcal E
=
\mathcal AC_2\mathcal A_S
- [\mathcal A,\Lambda_2]_{S_3}
+ (\partial_2C_2)\mathcal A_S
+ YC_2
+ C_2X_S,
\]
so that
\begin{align*}
L(\mathcal E)
={}&
\partial_2(\mathcal AC_2\mathcal A_S
- \mathcal A\Lambda_2
+ \Lambda_2\mathcal A_S
+ (\partial_2C_2)\mathcal A_S
+ YC_2
+ C_2X_S)
\\
&- C_1(\mathcal AC_2\mathcal A_S
- \mathcal A\Lambda_2
+ \Lambda_2\mathcal A_S
+ (\partial_2C_2)\mathcal A_S
+ YC_2
+ C_2X_S)
\\
&+ (\mathcal AC_2\mathcal A_S
- \mathcal A\Lambda_2
+ \Lambda_2\mathcal A_S
+ (\partial_2C_2)\mathcal A_S
+ YC_2
+ C_2X_S)C_{1,S}.
\end{align*}
Substituting these three displays, together with
$\Delta_1 = \partial_2\mathcal A + \mathcal AC_1 - C_1\mathcal A$,
$\Delta_{1,S} = \partial_2\mathcal A_S
  + \mathcal A_SC_{1,S} - C_{1,S}\mathcal A_S$,
and $\Delta_2 = \mathcal AC_2 - C_2\mathcal A_S$, into
$\mathcal R = \mathcal P_\Lambda - L(\mathcal E)$ and preserving
factor order throughout gives the fully expanded form
\begin{align}
\mathcal R
={}&
(\mathcal H\mathcal A)C_2
+ \mathcal A(\mathcal HC_2)
- (\mathcal HC_2)\mathcal A_S
- C_2(\mathcal H\mathcal A_S)
+ (\partial_2^2\mathcal A)C_2
- (\partial_2^2C_2)\mathcal A_S
\notag\\
&- \Lambda_1\mathcal AC_2
+ \Lambda_1C_2\mathcal A_S
+ \mathcal AC_2\Lambda_{1,S}
- C_2\mathcal A_S\Lambda_{1,S}
\notag\\
&- \mathcal AC_1\mathcal AC_2
+ C_1\mathcal A^2C_2
+ \mathcal AC_1C_2\mathcal A_S
- \mathcal AC_2C_{1,S}\mathcal A_S
\notag\\
&+ \mathcal A(\partial_2C_1)C_2
- (\partial_2C_1)\mathcal AC_2
\notag\\
&- \mathcal A(\partial_2C_2)\mathcal A_S
- \mathcal A(\partial_2C_2)C_{1,S}
+ C_1(\partial_2C_2)\mathcal A_S
\notag\\
&- (\partial_2\mathcal A)\mathcal AC_2
+ (\partial_2\mathcal A)C_1C_2
- (\partial_2\mathcal A)C_2C_{1,S}
- (\partial_2\mathcal A)(\partial_2C_2)
- C_1(\partial_2\mathcal A)C_2
\notag\\
&- C_2\mathcal A_S^2C_{1,S}
+ C_2\mathcal A_SC_{1,S}\mathcal A_S
- C_2\mathcal A_S(\partial_2\mathcal A_S)
+ C_2(\partial_2\mathcal A_S)C_{1,S}
\notag\\
&- \mathcal AC_1\Lambda_2
+ \mathcal A(\partial_2\Lambda_2)
+ \mathcal A\Lambda_2C_{1,S}
+ C_1\Lambda_2\mathcal A_S
- (\partial_2\Lambda_2)\mathcal A_S
- \Lambda_2C_{1,S}\mathcal A_S
\notag\\
&+ C_1C_2X_S
- C_2X_SC_{1,S}
- C_2(\partial_2X_S)
+ C_1YC_2
- YC_2C_{1,S}
- Y(\partial_2C_2)
\notag\\
&- (\partial_2Y)C_2
- (\partial_2C_2)X_S.
\label{eq:lambda-gap-raw-R}
\end{align}
We now verify \eqref{eq:lambda-gap-defect-factorization}.  The
twelve terms on the right-hand side group naturally into four blocks,
each of which we expand separately and match against
\eqref{eq:lambda-gap-raw-R}.

The first block collects the defects that act on the left, with
$C_2$ as a right factor:
\begin{align*}
B_{\mathcal A}
&\defeq
(\partial_2\mathcal D_{\partial\mathcal A})C_2
+ \mathcal D_HC_2
- C_1\mathcal D_{\partial\mathcal A}C_2
\\
={}&
\bigl[
  \partial_2^2\mathcal A
  - \partial_2X
  - (\partial_2\mathcal A)\mathcal A
  - \mathcal A(\partial_2\mathcal A)
  - \partial_2Y
\bigr]C_2
\\
&+ \bigl[
  \mathcal H\mathcal A
  + \partial_2X
  - C_1X
  - YC_1
  - \Lambda_1\mathcal A
\\
&\qquad
  - (\partial_2C_1)\mathcal A
  + \mathcal A\Lambda_1
  + \mathcal A(\partial_2\mathcal A)
  - \mathcal AC_1\mathcal A
\bigr]C_2
\\
&- C_1\bigl[
  \partial_2\mathcal A - X - \mathcal A^2 - Y
\bigr]C_2
\\
={}&
(\partial_2^2\mathcal A)C_2
+ (\mathcal H\mathcal A)C_2
- \Lambda_1\mathcal AC_2
+ \mathcal A\Lambda_1C_2
\\
&- (\partial_2C_1)\mathcal AC_2
- \mathcal AC_1\mathcal AC_2
\\
&- (\partial_2\mathcal A)\mathcal AC_2
- (\partial_2Y)C_2
\\
&- C_1(\partial_2\mathcal A)C_2
+ C_1\mathcal A^2C_2 \\
&+ C_1YC_2
- YC_1C_2.
\end{align*}
The $(\partial_2 X)C_2$ and $\mathcal A(\partial_2\mathcal A)C_2$
terms from $(\partial_2\mathcal D_{\partial\mathcal A})C_2$ and
$\mathcal D_HC_2$ cancel pairwise, and the $C_1XC_2$ terms from
$\mathcal D_HC_2$ and $C_1\mathcal D_{\partial\mathcal A}C_2$
cancel, leaving twelve terms.

The second block collects the defects that act on the right, with
$C_2$ as a left factor:
\begin{align*}
B_S
&\defeq
- C_2\mathcal D_{H,S}
+ C_2\mathcal D_{\partial\mathcal A,S}C_{1,S}
\\
={}&
- C_2\bigl[
  \mathcal H\mathcal A_S
  + \partial_2X_S
  - C_{1,S}X_S
  - Y_SC_{1,S}
  - \Lambda_{1,S}\mathcal A_S
  - (\partial_2C_{1,S})\mathcal A_S
\\
&\hspace{3cm}
  + \mathcal A_S\Lambda_{1,S}
  + \mathcal A_S(\partial_2\mathcal A_S)
  - \mathcal A_SC_{1,S}\mathcal A_S
\bigr]
\\
&+ C_2\bigl[
  \partial_2\mathcal A_S - X_S - \mathcal A_S^2 - Y_S
\bigr]C_{1,S}
\\
={}&
- C_2(\mathcal H\mathcal A_S)
- C_2(\partial_2X_S)
+ C_2C_{1,S}X_S
- C_2X_SC_{1,S}
+ C_2\Lambda_{1,S}\mathcal A_S
- C_2\mathcal A_S\Lambda_{1,S}
\\
&+ C_2(\partial_2C_{1,S})\mathcal A_S
- C_2\mathcal A_S(\partial_2\mathcal A_S)
+ C_2(\partial_2\mathcal A_S)C_{1,S} \\
&- C_2\mathcal A_S^2C_{1,S}
+ C_2\mathcal A_SC_{1,S}\mathcal A_S.
\end{align*}
The $C_2Y_SC_{1,S}$ terms from $C_2\mathcal D_{H,S}$ and
$C_2\mathcal D_{\partial\mathcal A,S}C_{1,S}$ cancel, leaving
eleven terms.

The third block collects the mixed-diamond defect, sandwiched by
$\mathcal A$ and $\mathcal A_S$:
\begin{align*}
B_M
&\defeq
\mathcal A\mathcal D_M
- \mathcal D_M\mathcal A_S
\\
={}&
\mathcal A(\mathcal HC_2)
- (\mathcal HC_2)\mathcal A_S
- \mathcal A(\partial_2C_2)C_{1,S}
+ \mathcal AC_2\Lambda_{1,S}
- \mathcal A\Lambda_1C_2
+ \mathcal A(\partial_2C_1)C_2
\\
&- \mathcal AC_1\Lambda_2
+ \mathcal A(\partial_2\Lambda_2)
+ \mathcal A\Lambda_2C_{1,S}
+ C_1\Lambda_2\mathcal A_S
- C_2\Lambda_{1,S}\mathcal A_S
\\
&+ \Lambda_1C_2\mathcal A_S
- \Lambda_2C_{1,S}\mathcal A_S
- (\partial_2C_1)C_2\mathcal A_S
+ (\partial_2C_2)C_{1,S}\mathcal A_S
- (\partial_2\Lambda_2)\mathcal A_S.
\end{align*}

The fourth block collects the $C$-diamond defect with its five
coefficient operators:
\begin{align*}
B_C
&\defeq
- \mathcal A\mathcal D_C\mathcal A_S
- (\partial_2\mathcal D_C)\mathcal A_S
- Y\mathcal D_C
- \mathcal D_CX_S
- (\partial_2\mathcal A)\mathcal D_C
\\
={}&
- \mathcal A(\partial_2C_2)\mathcal A_S
+ \mathcal AC_1C_2\mathcal A_S
- \mathcal AC_2C_{1,S}\mathcal A_S
\\
&- (\partial_2^2C_2)\mathcal A_S
+ (\partial_2C_1)C_2\mathcal A_S
+ C_1(\partial_2C_2)\mathcal A_S \\
&- (\partial_2C_2)C_{1,S}\mathcal A_S
- C_2(\partial_2C_{1,S})\mathcal A_S- Y(\partial_2C_2)
+ YC_1C_2
- YC_2C_{1,S}
\\
&- (\partial_2C_2)X_S
+ C_1C_2X_S
- C_2C_{1,S}X_S
\\
&- (\partial_2\mathcal A)(\partial_2C_2)
+ (\partial_2\mathcal A)C_1C_2
- (\partial_2\mathcal A)C_2C_{1,S}.
\end{align*}

The sum $B_{\mathcal A} + B_S + B_M + B_C$ is exactly the
right-hand side of \eqref{eq:lambda-gap-defect-factorization}.
Seven cross-block pairs cancel before comparing with
\eqref{eq:lambda-gap-raw-R}: the $\mathcal A\Lambda_1C_2$ terms
between $B_{\mathcal A}$ and $B_M$, the $YC_1C_2$ terms between
$B_{\mathcal A}$ and $B_C$, the $C_2\Lambda_{1,S}\mathcal A_S$
terms between $B_S$ and $B_M$, the $C_2C_{1,S}X_S$ and
$C_2(\partial_2C_{1,S})\mathcal A_S$ terms between $B_S$ and $B_C$,
and the $(\partial_2C_1)C_2\mathcal A_S$ and
$(\partial_2C_2)C_{1,S}\mathcal A_S$ terms between $B_M$ and $B_C$.
After removing these pairs, the remaining terms match
\eqref{eq:lambda-gap-raw-R} term by term, so
\eqref{eq:lambda-gap-defect-factorization} holds.

Since all defects on the right-hand side of
\eqref{eq:lambda-gap-defect-factorization} vanish, together with
their $\partial_2$-derivatives, $\mathcal R = 0$.
\end{proof}

\begin{proof}[Proof of Theorem~\ref{thm:dc-clean-general-darboux}]
Equation \eqref{eq:dc-clean-C-diamond} is unchanged, since the
Darboux transformation fixes $C_1$ and $C_2$.

For \eqref{eq:dc-clean-mixed-diamond}, substituting
$\widehat\Lambda_k = \Lambda_k + \Delta_k$ into
\eqref{eq:dc-clean-mixed-diamond} and subtracting the seed equation
gives
\begin{align*}
&\partial_2\Delta_2
+\Delta_2C_{1,S}
+C_2\Delta_{1,S}
-\Delta_1C_2
-C_1\Delta_2 \\
&\qquad=
\mathcal A
\bigl[
\partial_2C_2+C_2C_{1,S}-C_1C_2
\bigr]
-
\bigl[
\partial_2C_2+C_2C_{1,S}-C_1C_2
\bigr]\mathcal A_S.
\end{align*}
The bracket vanishes by \eqref{eq:dc-clean-C-diamond}.

It remains to prove \eqref{eq:dc-clean-Lambda-residual}.  By
Lemma~\ref{lem:dc-clean-general-Lambda-residual}, the left-hand side
of \eqref{eq:dc-clean-Lambda-residual} equals
$\partial_2\mathcal E - C_1\mathcal E + \mathcal EC_{1,S}$.
Lemma~\ref{lem:dc-clean-shifted-defect} proves that
$\mathcal E$ vanishes identically, so this expression is zero.

Since the seed $\Lambda$-diamond holds and the residual
\eqref{eq:dc-clean-Lambda-residual} vanishes, the dressed
data satisfy
\eqref{eq:dc-clean-C-diamond}--\eqref{eq:dc-clean-mixed-diamond}
and \eqref{eq:dc-clean-L-diamond}.
\end{proof}
\section{Bilinear Equations}

In the probabilistic applications of subsequent chapters, the
Fredholm determinant $F(u) \defeq \det_H(I - zK(u))$ encodes the
distributional data of each model.  The proposition below extracts a
bilinear identity for $F$ from the matrix Darboux theory through the
relation $\partial_2 \log F = -\operatorname{tr}_E \mathcal A$.

\begin{proposition}
\label{prop:dc-clean-normalized-trace-defect-direct}
Under the hypotheses of
Proposition~\ref{prop:dc-clean-dressed-waves}, assume additionally
that $\dim E<\infty$, that $H$ is a real separable Hilbert space
with $K(u)\in\mathfrak S_1(H)$, that $K$ is $C^1$ in $\partial_1$
and $C^2$ in $\partial_2$ in trace norm, and that
$C_2(u)\in\operatorname{GL}(E)$ for every $u\in\mathcal V$.  Define
$F(u)\defeq \det_H(I-zK(u))$ and
$\mathcal B(u)\defeq \Delta_2(u)C_2(u)^{-1}$.
Then
\begin{align}
\frac{
\left(D_1-\frac{1}{2}D_2^2\right)
F(S_3u)\cdot F(u)
}{
F(S_3u)F(u)
}
&=
\frac{1}{2}
\left[
\operatorname{tr}_E\bigl(\mathcal B(u)^2\bigr)
-
\bigl(\operatorname{tr}_E\mathcal B(u)\bigr)^2
\right]
\notag \\
&\quad
+\operatorname{tr}_E
\left[
\mathcal B(u)
\left(
C_1(u)+\Lambda_2(u)C_2(u)^{-1}
\right)
\right].
\label{eq:dc-clean-normalized-trace-defect-direct}
\end{align}
\end{proposition}

\begin{corollary}
\label{cor:dc-clean-scalar-trace-defect}
If $C_1(u)+\Lambda_2(u)C_2(u)^{-1}=a(u)I$ for a scalar
function~$a$, then
\begin{align}
\frac{
\left(D_1-a(u)D_2-\frac{1}{2}D_2^2\right)
F(S_3u)\cdot F(u)
}{
F(S_3u)F(u)
}
&=
\frac{1}{2}
\left[
\operatorname{tr}_E\bigl(\mathcal B(u)^2\bigr)
-
\bigl(\operatorname{tr}_E\mathcal B(u)\bigr)^2
\right].
\label{eq:dc-clean-normalized-scalar-linear-trace-defect-direct}
\end{align}
When $\dim E = 1$, the trace-defect term vanishes and the
Fredholm determinant satisfies the scalar Hirota equation
\begin{equation}\label{eq:dc-clean-scalar-Hirota}
\left(D_1-a(u)D_2-\frac{1}{2}D_2^2\right)
F(S_3u)\cdot F(u) = 0.
\end{equation}
\end{corollary}

\begin{proof}[Proof of Proposition~\ref{prop:dc-clean-normalized-trace-defect-direct}]
The Hirota derivative identities
$D_1 F_S\cdot F/(F_SF) = \partial_1\log(F_S/F)$ and
$D_2^2 F_S\cdot F/(F_SF) = \partial_2^2\log(F_SF)
+ (\partial_2\log(F_S/F))^2$
reduce the proof to computing the log-derivatives of~$F$.

By the Jacobi formula and \eqref{eq:dc-clean-K-d2},
\[
\partial_2\log F
= -z\operatorname{Tr}_H(R\Psi\Phi)
= -\operatorname{tr}_E\mathcal A,
\qquad
\partial_2\log\frac{F_S}{F}
= \operatorname{tr}_E(\mathcal A-\mathcal A_S).
\]
By \eqref{eq:dc-clean-d2A-expanded},
$\partial_2^2\log F = -\operatorname{tr}_E(X+\mathcal A^2+Y)$
and
$\partial_2^2\log F_S = -\operatorname{tr}_E(X_S+\mathcal A_S^2+Y_S)$.

Substituting \eqref{eq:dc-clean-K-d1} into the Jacobi formula gives
$\partial_1\log F
= \frac{1}{2}\operatorname{tr}_E(X-Y)
- \operatorname{tr}_E(\mathcal AC_1)$,
and the difference of the shifted and unshifted evaluations is
\begin{align*}
\partial_1\log\frac{F_S}{F}
&=
\frac{1}{2}\operatorname{tr}_E(X_S-Y_S-X+Y)
-\operatorname{tr}_E(\mathcal A_SC_{1,S})
+\operatorname{tr}_E(\mathcal AC_1).
\end{align*}

Substituting all five expressions into the Hirota identities gives
\begin{align}
\frac{
\left(D_1-\frac{1}{2}D_2^2\right)F_S\cdot F
}{
F_SF
}
=
\operatorname{tr}_E
\bigl[
\mathcal A C_1-\mathcal A_SC_{1,S}
+Y+X_S
+\frac{1}{2}\mathcal A^2
+\frac{1}{2}\mathcal A_S^2
\bigr]
-\frac{1}{2}
\left[
\operatorname{tr}_E(\mathcal A-\mathcal A_S)
\right]^2.
\label{eq:dc-direct-before-E2}
\end{align}

It remains to eliminate $\operatorname{tr}_E(Y+X_S)$.
Lemma~\ref{lem:dc-clean-shifted-defect} gives
$YC_2+C_2X_S+\mathcal AC_2\mathcal A_S
+(\partial_2C_2)\mathcal A_S
-\mathcal A\Lambda_2+\Lambda_2\mathcal A_S = 0$.
Right-multiplying by $C_2^{-1}$ and taking
$\operatorname{tr}_E$ gives
\begin{align*}
\operatorname{tr}_E(Y+X_S)
=
-\operatorname{tr}_E(\mathcal AC_2\mathcal A_SC_2^{-1})
-\operatorname{tr}_E((\partial_2C_2)\mathcal A_SC_2^{-1})
+\operatorname{tr}_E(\mathcal A\Lambda_2C_2^{-1})
-\operatorname{tr}_E(\Lambda_2\mathcal A_SC_2^{-1}).
\end{align*}
Equation \eqref{eq:dc-clean-C-diamond} gives
$(\partial_2C_2)\mathcal A_SC_2^{-1}
= C_1C_2\mathcal A_SC_2^{-1} - C_2C_{1,S}\mathcal A_SC_2^{-1}$,
and under $\operatorname{tr}_E$ the second term becomes
$\operatorname{tr}_E(\mathcal A_SC_{1,S})$.
Substituting into \eqref{eq:dc-direct-before-E2}, the
$\operatorname{tr}_E(\mathcal A_SC_{1,S})$ terms cancel, leaving
\begin{align}
\frac{
\left(D_1-\frac{1}{2}D_2^2\right)F_S\cdot F
}{
F_SF
}
&=
\operatorname{tr}_E(\mathcal A C_1)
-\operatorname{tr}_E(C_1C_2\mathcal A_SC_2^{-1})
\notag \\
&\quad
+\operatorname{tr}_E(\mathcal A\Lambda_2C_2^{-1})
-\operatorname{tr}_E(\Lambda_2\mathcal A_SC_2^{-1})
\notag \\
&\quad
-\operatorname{tr}_E(\mathcal A C_2\mathcal A_SC_2^{-1})
+\frac{1}{2}\operatorname{tr}_E(\mathcal A^2)
+\frac{1}{2}\operatorname{tr}_E(\mathcal A_S^2)
\notag \\
&\quad
-\frac{1}{2}
\left[
\operatorname{tr}_E(\mathcal A-\mathcal A_S)
\right]^2.
\label{eq:dc-direct-expanded-final}
\end{align}

Expanding $\mathcal B = \mathcal A - C_2\mathcal A_SC_2^{-1}$,
trace cyclicity gives
$\operatorname{tr}_E\mathcal B = \operatorname{tr}_E(\mathcal A-\mathcal A_S)$.
The first four terms of \eqref{eq:dc-direct-expanded-final} are
$\operatorname{tr}_E[\mathcal B(C_1+\Lambda_2C_2^{-1})]$, and the
next three give $\frac{1}{2}\operatorname{tr}_E(\mathcal B^2)$, so
\eqref{eq:dc-direct-expanded-final} is exactly
\eqref{eq:dc-clean-normalized-trace-defect-direct}.
\end{proof}

\begin{proof}[Proof of Corollary~\ref{cor:dc-clean-scalar-trace-defect}]
If \(C_1+\Lambda_2C_2^{-1}=aI\), then
\[
\operatorname{tr}_E
\left[
\mathcal B
\left(C_1+\Lambda_2C_2^{-1}\right)
\right]
=
a\operatorname{tr}_E\mathcal B
=
a
\frac{D_2F_S\cdot F}{F_SF}.
\]
Moving this term to the left in
\eqref{eq:dc-clean-normalized-trace-defect-direct} gives
\eqref{eq:dc-clean-normalized-scalar-linear-trace-defect-direct}.
If $\dim E = 1$, then
\[
\operatorname{tr}_E(\mathcal B^2)
=
\bigl(\operatorname{tr}_E\mathcal B\bigr)^2,
\]
so the right-hand side vanishes.
\end{proof}

%% file: chapters/5-the-continuum-framework.tex
\chapter{The Continuum Framework}
\label{ch:the-continuum-framework}

{
  \setlength{\parskip}{0pt}
}

\label{sec:kpz-reduced-framework}

This chapter develops the continuum analogue of the parabolic
framework (Chapter~\ref{ch:the-parabolic-framework}), in which the
surviving discrete shift~$S_3$ degenerates to produce a fully
continuous linear problem in three real variables~$(t,x,a)$.
The compatibility conditions, which reduce from three diamond
families to a pair of coupled PDEs
(Proposition~\ref{prop:kpz-red-compat}), admit a Darboux
transformation (Theorem~\ref{thm:kpz-red-Darboux-compat}) that
dresses both~$C$ and~$\Lambda$ additively through the observable
$\mathcal{A} = z\Phi R\Psi$, which satisfies a matrix potential KP
equation with background corrections in~$(C,\Lambda)$ that vanish
when $C = \Lambda = 0$
(Corollary~\ref{cor:kpz-shift-vacuum-kp}).  The trace-defect Hirota
equation (Proposition~\ref{prop:kpz-shift-trace-defect}) extracts a
bilinear identity for the Fredholm determinant that reduces to a
scalar bilinear equation when $\dim E = 1$.  The continuum
framework arises directly from the fully discrete theory
(Chapter~\ref{ch:the-discrete-framework}) under the scaling
$T^{(\epsilon)} = e^{-\frac{1}{2}\epsilon^2\pa_x + \epsilon\pa_a}$, $S_1^{(\epsilon)} = e^{6\epsilon^3\pa_t - 2\epsilon^2\pa_x - 2\epsilon\pa_a}$,
$S_2^{(\epsilon)} = e^{-\frac{1}{2}\epsilon^2\pa_x - \epsilon\pa_a}$,
$\epsilon \to 0$
(Remark~\ref{rem:kpz-scaling}).  All results are stated
and proved independently.

\section{The Continuum Linear Problem}
\label{sec:kpz-red-compatibility}

Let $E$ and $H \neq 0$ be vector spaces over $\mathbb{F} = \R$ or $\C$, and
let $\mathcal V \subseteq \R^3$ be an open set with
coordinates $(t,x,a)$.  Write
$\pa_t$, $\pa_x$, and $\pa_a$ for the three
coordinate derivatives.  Assign endomorphisms\footnote{Throughout, all maps of $(t,x,a)$ are assumed smooth.}
\begin{equation*}
C(u), \Lambda(u) \in \End(E).
\end{equation*}
Consider the overdetermined linear problem for
$\Psi\colon\mathcal V\to\Hom(E,H)$:
\begin{align}
\Psi_x
  &= \Psi_{aa}-2\Psi C,
  \label{eq:kpz-red-Lx} \\
\Psi_t
  &= -\tfrac13\Psi_{aaa}+\Psi_a C-\Psi\Lambda.
  \label{eq:kpz-red-Lt}
\end{align}
Compatibility of the linear problem
\eqref{eq:kpz-red-Lx}--\eqref{eq:kpz-red-Lt} forces the
following conditions on $(C,\Lambda)$.

\begin{proposition}
\label{prop:kpz-red-compat}
The system
\eqref{eq:kpz-red-Lx}--\eqref{eq:kpz-red-Lt} is
compatible\footnote{That is, the two computations of $\Psi_{xt}$,
obtained by applying $\pa_t$ to \eqref{eq:kpz-red-Lx} and
$\pa_x$ to \eqref{eq:kpz-red-Lt}, agree for arbitrary Cauchy
data $\Psi$, $\Psi_a$, $\Psi_{aa}$.} if and only if
\begin{equation}
\label{eq:kpz-red-constraint}
2\Lambda_a+C_x+C_{aa} =0,
\end{equation}
\begin{equation}
\label{eq:kpz-red-evolution}
2C_t-\Lambda_x+\Lambda_{aa}
+\tfrac23 C_{aaa}-2C_aC+2[C,\Lambda] =0.
\end{equation}
\end{proposition}

\begin{proof}
The compatibility equations are obtained by computing
$\Psi_{xt}$ and $\Psi_{tx}$ separately.

First differentiate \eqref{eq:kpz-red-Lx} in $t$ and use
\eqref{eq:kpz-red-Lt} to eliminate $\Psi_t$ and its
$a$-derivatives.  This gives
\begin{align*}
\Psi_{xt}
&= -\tfrac13\Psi_{aaaaa}
   +\tfrac53\Psi_{aaa}C
   +\Psi_{aa}(2C_a-\Lambda) \\
&\quad
   +\Psi_a(C_{aa}-2\Lambda_a-2C^2)
   +\Psi(-\Lambda_{aa}+2\Lambda C-2C_t).
\end{align*}
Next differentiate \eqref{eq:kpz-red-Lt} in $x$ and use
\eqref{eq:kpz-red-Lx} to eliminate $\Psi_x$ and its
$a$-derivatives to give 
\begin{align*}
\Psi_{tx}
&= -\tfrac13\Psi_{aaaaa}
   +\tfrac53\Psi_{aaa}C
   +\Psi_{aa}(2C_a-\Lambda) \\
&\quad
   +\Psi_a(2C_{aa}-2C^2+C_x)
   +\Psi(\tfrac23 C_{aaa}-2C_aC+2C\Lambda-\Lambda_x).
\end{align*}
The coefficients of $\Psi_{aaaaa}$, $\Psi_{aaa}$, and
$\Psi_{aa}$ agree identically.  If the system is compatible,
the remaining coefficients must agree for arbitrary Cauchy
data $\Psi$, $\Psi_a$, $\Psi_{aa}$.  Equating the
coefficients of $\Psi_a$ gives
\eqref{eq:kpz-red-constraint}, and equating the coefficients
of $\Psi$ gives \eqref{eq:kpz-red-evolution}.  Conversely, if
\eqref{eq:kpz-red-constraint}--\eqref{eq:kpz-red-evolution}
hold, then $\Psi_{xt}$ and $\Psi_{tx}$ have identical
coefficients and therefore agree for arbitrary Cauchy data.
\end{proof}

\begin{remark}\label{rem:kpz-scaling}
The continuum framework arises from the fully discrete theory of
Chapter~\ref{ch:the-discrete-framework} under the scaling
$T^{(\epsilon)} = e^{-\frac{1}{2}\epsilon^2\pa_x + \epsilon\pa_a}$,
$S_1^{(\epsilon)} = e^{6\epsilon^3\pa_t
- 2\epsilon^2\pa_x - 2\epsilon\pa_a}$,
$S_2^{(\epsilon)} = e^{-\frac{1}{2}\epsilon^2\pa_x - \epsilon\pa_a}$.
Expand the discrete edge weights as
\[
C_1^{(\epsilon)} = -2I + O(\epsilon^4), \qquad
\Lambda_1^{(\epsilon)} = -3I - 6\epsilon^2 C
  + 6\epsilon^3(\Lambda + C_a) + O(\epsilon^4),
\]
\[
C_2^{(\epsilon)} = -I + O(\epsilon^4), \qquad
\Lambda_2^{(\epsilon)} = -2I - 2\epsilon^2 C + O(\epsilon^4).
\]
At order~$\epsilon^2$ the discrete equations recover the
$x$-equation~\eqref{eq:kpz-red-Lx}, and at order~$\epsilon^3$
the $t$-equation~\eqref{eq:kpz-red-Lt}.
The same scaling yields the adjoint problem, dressing
conditions, and Darboux transformation of
\S\ref{sec:kpz-red-Darboux}.  A sequential scaling through
the parabolic theory
(Chapter~\ref{ch:the-parabolic-framework}) also exists,
identifying $\pa_1 = \tfrac{1}{2}\pa_x + \epsilon\pa_t$,
$\pa_2 = \pa_a$,
$S_3 = e^{-\frac{1}{2}\epsilon^2\pa_x - \epsilon\pa_a}$.
\end{remark}

\section{Darboux Transformations}
\label{sec:kpz-red-Darboux}

The Darboux transformation produces, from a solution
$(C,\Lambda)$ of the compatibility equations together with
wave functions and a dressing-compatible kernel, new
coefficients $(\widehat C,\widehat\Lambda)$ satisfying the
same compatibility equations of
Proposition~\ref{prop:kpz-red-compat}.  We now develop this
construction.

Fix a solution $(C,\Lambda)$ of
\eqref{eq:kpz-red-constraint}--\eqref{eq:kpz-red-evolution}.
Let
$\Psi\colon\mathcal V\to\Hom(E,H)$ solve
\eqref{eq:kpz-red-Lx}--\eqref{eq:kpz-red-Lt}.
The adjoint linear problem\footnote{Its compatibility
conditions are those of
Proposition~\ref{prop:kpz-red-compat}.} for
$\Phi\colon\mathcal V\to\Hom(H,E)$ is
\begin{align}
\Phi_x
  &= -\Phi_{aa}+2C\Phi,
  \label{eq:kpz-red-Ax} \\
\Phi_t
  &= -\tfrac13\Phi_{aaa}+C\Phi_a+(\Lambda+C_a)\Phi.
  \label{eq:kpz-red-At}
\end{align}

\begin{definition}
\label{def:kpz-red-dressing}
We say that $K\colon\mathcal V\to\End(H)$ is
\emph{dressing compatible} with $(\Psi,\Phi)$ if:
\begin{enumerate}[label=\arabic*., leftmargin=*]
\item For all $(t,x,a) \in \mathcal V$,
\begin{align}
K_a
  &= -\Psi\Phi,
  \label{eq:kpz-red-Ka} \\
K_x
  &= -\Psi_a\Phi+\Psi\Phi_a,
  \label{eq:kpz-red-Kx} \\
K_t
  &= \tfrac13(\Psi_{aa}\Phi-\Psi_a\Phi_a+\Psi\Phi_{aa})
     -\Psi C\Phi.
  \label{eq:kpz-red-Kt}
\end{align}
\item There exists $z\in\mathbb{F}$ such that the resolvent
exists for all $(t,x,a) \in \mathcal V$:
\[
R\defeq(I-zK)^{-1}\in\End(H).
\]
\end{enumerate}
\end{definition}

The dressed waves inherit the linear problem with corrected
coefficients.

\begin{proposition}
\label{thm:kpz-red-Darboux}
Let $K$ be dressing compatible with $(\Psi,\Phi)$, with
resolvent $R=(I-zK)^{-1}$.  Define the dressed observable and
dressed waves by
\begin{equation}
\label{eq:kpz-red-A-def}
\mathcal A\defeq z\Phi R\Psi,
\qquad
\widehat\Psi\defeq R\Psi,
\qquad
\widehat\Phi\defeq\Phi R,
\end{equation}
and set
\begin{equation}
\label{eq:kpz-red-hat-CL}
\widehat C\defeq C-\mathcal A_a,
\qquad
\widehat\Lambda\defeq\Lambda+\tfrac12(\mathcal A_x+
\mathcal A_{aa}).
\end{equation}
Then $\widehat\Psi$ and $\widehat\Phi$ satisfy the linear
problem
\eqref{eq:kpz-red-Lx}--\eqref{eq:kpz-red-Lt} and the
adjoint problem
\eqref{eq:kpz-red-Ax}--\eqref{eq:kpz-red-At} respectively,
with coefficients $(\widehat C,\widehat\Lambda)$, and
$\widehat K\defeq KR$ is dressing compatible with
$(\widehat\Psi,\widehat\Phi)$.
\end{proposition}

\begin{proof}
For any coordinate derivative $\pa_\beta$, the resolvent
identity
\begin{equation}
\label{eq:kpz-red-resolvent-alpha}
R_\beta=zRK_\beta R
\end{equation}
follows from differentiating $(I-zK)R=I$.
In particular, \eqref{eq:kpz-red-Ka} gives
$R_a=-zR\Psi\Phi R$, and hence differentiating
$\widehat\Psi=R\Psi$ and $\widehat\Phi=\Phi R$ by the
Leibniz rule gives
\begin{equation}
\label{eq:kpz-red-hatPsi-a}
\widehat\Psi_a=R\Psi_a-\widehat\Psi\mathcal A,
\qquad
\widehat\Phi_a=\Phi_aR-\mathcal A\widehat\Phi.
\end{equation}
The following expressions appear throughout.  Set
\begin{equation}
\label{eq:kpz-red-F-rs}
\mathcal F_{r,s}\defeq
z(\pa_a^r\Phi)R(\pa_a^s\Psi),
\qquad r,s\geq0.
\end{equation}
Then $\mathcal F_{0,0}=\mathcal A$, and differentiating in
$a$ gives
\begin{equation}
\label{eq:kpz-red-Frs-a}
\pa_a\mathcal F_{r,s}
=
\mathcal F_{r+1,s}
-\mathcal F_{r,0}\mathcal F_{0,s}
+\mathcal F_{r,s+1}.
\end{equation}
In particular, evaluating \eqref{eq:kpz-red-Frs-a} at
$r=s=0$ gives
\begin{equation}
\label{eq:kpz-red-Aa-Frs}
\mathcal A_a=
\mathcal F_{1,0}-\mathcal A^2+\mathcal F_{0,1}.
\end{equation}
The following identities express $R\Psi_a$, $R\Psi_{aa}$,
and $R\Psi_{aaa}$ in terms of $\widehat\Psi$ and its
derivatives:
\begin{align}
R\Psi_a
&=\widehat\Psi_a+\widehat\Psi\mathcal A,
\label{eq:kpz-red-RPsi-a} \\
R\Psi_{aa}
&=\widehat\Psi_{aa}
+\widehat\Psi_a\mathcal A
+\widehat\Psi(\mathcal A_a+\mathcal F_{0,1}),
\label{eq:kpz-red-RPsi-aa} \\
R\Psi_{aaa}
&=\widehat\Psi_{aaa}
+\widehat\Psi_{aa}\mathcal A
+\widehat\Psi_a(2\mathcal A_a+
  \mathcal F_{0,1})
+\widehat\Psi(\mathcal A_{aa}+
  \pa_a\mathcal F_{0,1}+\mathcal F_{0,2}).
\label{eq:kpz-red-RPsi-aaa}
\end{align}
The first is a rearrangement of
\eqref{eq:kpz-red-hatPsi-a}.  For the second,
differentiate \eqref{eq:kpz-red-hatPsi-a} in $a$ to give
$\widehat\Psi_{aa} = R_a\Psi_a + R\Psi_{aa}
- \widehat\Psi_a\mathcal A - \widehat\Psi\mathcal A_a$.
Substituting $K_a=-\Psi\Phi$ from \eqref{eq:kpz-red-Ka}
into the resolvent identity
\eqref{eq:kpz-red-resolvent-alpha} and right-multiplying
by $\Psi_a$ gives
$R_a\Psi_a = -z(R\Psi)(\Phi R\Psi_a)
= -\widehat\Psi\mathcal F_{0,1}$, and solving for
$R\Psi_{aa}$ gives \eqref{eq:kpz-red-RPsi-aa}.
Differentiating once more in $a$, evaluating
$R_a\Psi_{aa}=-\widehat\Psi\mathcal F_{0,2}$ by
\eqref{eq:kpz-red-resolvent-alpha} and
\eqref{eq:kpz-red-F-rs} with $(r,s)=(0,2)$, and solving
gives \eqref{eq:kpz-red-RPsi-aaa}.

\noindent\textit{The dressed wave \eqref{eq:kpz-red-Lx}.}
Substituting
$K_x=-\Psi_a\Phi+\Psi\Phi_a$ from
\eqref{eq:kpz-red-Kx} into
\eqref{eq:kpz-red-resolvent-alpha} and regrouping gives
\[
R_x\Psi
= -z(R\Psi_a)(\Phi R\Psi)
  + z(R\Psi)(\Phi_aR\Psi)
= -(R\Psi_a)\mathcal A
  + \widehat\Psi\mathcal F_{1,0}.
\]
By the Leibniz rule,
$\widehat\Psi_x = R_x\Psi + R\Psi_x$.
Using \eqref{eq:kpz-red-Lx} to replace $\Psi_x$ and
\eqref{eq:kpz-red-hatPsi-a} to expand $R\Psi_a$
in $R_x\Psi$ gives
\[
\widehat\Psi_x
= -\widehat\Psi_a\mathcal A
  - \widehat\Psi\mathcal A^2
  + \widehat\Psi\mathcal F_{1,0}
  + R\Psi_{aa}
  - 2\widehat\Psi C.
\]
Substituting \eqref{eq:kpz-red-RPsi-aa} gives
\begin{align*}
\widehat\Psi_x
&= \widehat\Psi_{aa}
   + \widehat\Psi_a\mathcal A
   + \widehat\Psi(\mathcal A_a+\mathcal F_{0,1})
   - \widehat\Psi_a\mathcal A
   - \widehat\Psi\mathcal A^2
   + \widehat\Psi\mathcal F_{1,0}
   - 2\widehat\Psi C.
\end{align*}
The terms $\pm\widehat\Psi_a\mathcal A$ cancel, leaving
\[
\widehat\Psi_x
=\widehat\Psi_{aa}
  +\widehat\Psi(\mathcal A_a+
    \mathcal F_{0,1}-\mathcal A^2+
    \mathcal F_{1,0})-2\widehat\Psi C.
\]
By \eqref{eq:kpz-red-Aa-Frs}, the parenthetical equals
$2\mathcal A_a$, giving
$\widehat\Psi_x = \widehat\Psi_{aa}
-2\widehat\Psi(C-\mathcal A_a)$, which is
\eqref{eq:kpz-red-Lx} with coefficients
$(\widehat C,\widehat\Lambda)$.

\noindent\textit{The dressed wave \eqref{eq:kpz-red-Lt}.}
By the Leibniz rule and the linear problem
\eqref{eq:kpz-red-Lt},
\begin{equation}
\label{eq:kpz-red-hatPsi-t-decomp}
\widehat\Psi_t=R_t\Psi-\tfrac13R\Psi_{aaa}
+(R\Psi_a)C-\widehat\Psi\Lambda.
\end{equation}
The right-hand side involves $R\Psi_a$ and $R\Psi_{aaa}$,
which are expressed in terms of $\widehat\Psi$ by
\eqref{eq:kpz-red-RPsi-a}--\eqref{eq:kpz-red-RPsi-aaa},
and $R_t\Psi$, which requires a separate computation.
The resolvent identity \eqref{eq:kpz-red-resolvent-alpha}
with $\beta=t$ gives $R_t\Psi=zR(K_t\widehat\Psi)$.
Substituting \eqref{eq:kpz-red-Kt} for $K_t$, identifying
$z\Phi\widehat\Psi=\mathcal A$,
$z\Phi_a\widehat\Psi=\mathcal F_{1,0}$, and
$z\Phi_{aa}\widehat\Psi=\mathcal F_{2,0}$ via
\eqref{eq:kpz-red-F-rs}, and applying $R$ from the left
gives
\begin{align}
R_t\Psi
&=\tfrac13\left[(R\Psi_{aa})\mathcal A
-(R\Psi_a)\mathcal F_{1,0}
+\widehat\Psi\mathcal F_{2,0}\right]
-\widehat\Psi C\mathcal A.
\label{eq:kpz-red-RtPsi}
\end{align}
Substituting
\eqref{eq:kpz-red-RPsi-a}--\eqref{eq:kpz-red-RtPsi} into
\eqref{eq:kpz-red-hatPsi-t-decomp} and collecting by
derivatives of $\widehat\Psi$ gives
\[
\widehat\Psi_t
=-\tfrac13\widehat\Psi_{aaa}
+\widehat\Psi_a\theta_1+
\widehat\Psi\theta_0,
\]
where
\[
\theta_1
=C+\tfrac13\mathcal A^2-\tfrac23\mathcal A_a
-\tfrac13\mathcal F_{0,1}-
\tfrac13\mathcal F_{1,0}
=C-\mathcal A_a,
\]
by \eqref{eq:kpz-red-Aa-Frs}.  The zeroth-order coefficient
$\theta_0$ will reduce to
$-\Lambda-\tfrac12(\mathcal A_x+\mathcal A_{aa})$; the
following two identities are needed for the reduction.
\begin{align}
\mathcal A_x
&=\mathcal F_{0,2}-\mathcal F_{2,0}
+\mathcal A\mathcal F_{1,0}
-\mathcal F_{0,1}\mathcal A
+2C\mathcal A-2\mathcal A C,
\label{eq:kpz-red-Ax-Frs} \\
\mathcal A_{aa}
&=\mathcal F_{2,0}+\mathcal F_{0,2}
+2\mathcal F_{1,1}
-\mathcal F_{1,0}\mathcal A
-\mathcal A\mathcal F_{0,1}
-\mathcal A_a\mathcal A
-\mathcal A\mathcal A_a.
\label{eq:kpz-red-Aaa-Frs}
\end{align}
The identities \eqref{eq:kpz-red-Ax-Frs} and
\eqref{eq:kpz-red-Aaa-Frs} follow from
\eqref{eq:kpz-red-Lx}, \eqref{eq:kpz-red-Ax},
\eqref{eq:kpz-red-Frs-a}, and the resolvent identity
\eqref{eq:kpz-red-resolvent-alpha}.
The zeroth-order coefficient extracted from the substitution is
\begin{align*}
\theta_0
=-\Lambda+
\mathcal A C-C\mathcal A
+\tfrac13\mathcal A_a\mathcal A
+\tfrac13\mathcal F_{0,1}\mathcal A
-\tfrac13\mathcal A\mathcal F_{1,0}
+\tfrac13\mathcal F_{2,0}
-\tfrac13\mathcal A_{aa}
-\tfrac13\pa_a\mathcal F_{0,1}
-\tfrac13\mathcal F_{0,2}.
\end{align*}
To verify the reduction, it suffices to check
$\theta_0+\Lambda+\tfrac12(\mathcal A_x+\mathcal A_{aa})=0$.
Substituting \eqref{eq:kpz-red-Ax-Frs} and
\eqref{eq:kpz-red-Aaa-Frs} into
$\tfrac12(\mathcal A_x+\mathcal A_{aa})$ gives
\begin{align*}
\tfrac12(\mathcal A_x+\mathcal A_{aa})
&= \mathcal F_{0,2}
  +\tfrac12\mathcal A\mathcal F_{1,0}
  -\tfrac12\mathcal F_{0,1}\mathcal A
  +C\mathcal A-\mathcal AC \\
&\quad
  +\mathcal F_{1,1}
  -\tfrac12\mathcal F_{1,0}\mathcal A
  -\tfrac12\mathcal A\mathcal F_{0,1}
  -\tfrac12\mathcal A_a\mathcal A
  -\tfrac12\mathcal A\mathcal A_a.
\end{align*}
Adding $\theta_0+\Lambda$, the $[C,\mathcal A]$ terms
$(\mathcal AC-C\mathcal A)+(C\mathcal A-\mathcal AC)$
cancel.  Expanding $\pa_a\mathcal F_{0,1}$ via
\eqref{eq:kpz-red-Frs-a} with $(r,s)=(0,1)$ and
substituting \eqref{eq:kpz-red-Aaa-Frs} for the remaining
$-\tfrac13\mathcal A_{aa}$ term, the $\mathcal F_{0,2}$,
$\mathcal F_{2,0}$, and $\mathcal F_{1,1}$ terms all
cancel in pairs.  The surviving expression is therefore
\[
\theta_0+\Lambda+\tfrac12(\mathcal A_x+\mathcal A_{aa})
=\tfrac16\bigl[
\mathcal A(\mathcal F_{1,0}+\mathcal F_{0,1})
-(\mathcal F_{1,0}+\mathcal F_{0,1})\mathcal A
+\mathcal A_a\mathcal A
-\mathcal A\mathcal A_a\bigr].
\]
By \eqref{eq:kpz-red-Aa-Frs},
$\mathcal F_{1,0}+\mathcal F_{0,1}
=\mathcal A_a+\mathcal A^2$, and the right-hand side
vanishes after expanding.  Therefore
$\theta_0=-\Lambda-\tfrac12(\mathcal A_x+
\mathcal A_{aa})=-\widehat\Lambda$, giving
$\widehat\Psi_t
=-\tfrac13\widehat\Psi_{aaa}
+\widehat\Psi_a\widehat C
-\widehat\Psi\widehat\Lambda$, which is
\eqref{eq:kpz-red-Lt} with coefficients
$(\widehat C,\widehat\Lambda)$.

\noindent\textit{The dressed adjoint \eqref{eq:kpz-red-Ax}.}
The adjoint calculation is the left-handed analogue, with
operators acting from the left on $\widehat\Phi$.  The
left-handed helper identities follow from
\eqref{eq:kpz-red-hatPsi-a} by the same successive
differentiation as
\eqref{eq:kpz-red-RPsi-a}--\eqref{eq:kpz-red-RPsi-aaa},
with $R$ on the right and $\mathcal F_{1,0}$ replacing
$\mathcal F_{0,1}$:
\begin{align}
\Phi_aR
&=\widehat\Phi_a+\mathcal A\widehat\Phi,
  \label{eq:kpz-red-PhiR-a} \\
\Phi_{aa}R
&=\widehat\Phi_{aa}
+\mathcal A\widehat\Phi_a
+(\mathcal F_{1,0}+\mathcal A_a)\widehat\Phi,
  \label{eq:kpz-red-PhiR-aa} \\
\Phi_{aaa}R
&=\widehat\Phi_{aaa}
+\mathcal A\widehat\Phi_{aa}
+(\mathcal F_{1,0}+2\mathcal A_a)\widehat\Phi_a
+(\mathcal A_{aa}
+\pa_a\mathcal F_{1,0}
+\mathcal F_{2,0})\widehat\Phi.
  \label{eq:kpz-red-PhiR-aaa}
\end{align}
By the Leibniz rule,
$\widehat\Phi_x = \Phi_xR + \Phi R_x$.
Substituting \eqref{eq:kpz-red-Ax} gives
$\Phi_xR = -\Phi_{aa}R + 2C\widehat\Phi$.
Substituting \eqref{eq:kpz-red-Kx} into
\eqref{eq:kpz-red-resolvent-alpha} and regrouping gives
\[
\Phi R_x
= -z(\Phi R\Psi_a)(\Phi R)
  + z(\Phi R\Psi)(\Phi_aR)
= -\mathcal F_{0,1}\widehat\Phi
  + \mathcal A(\Phi_aR).
\]
Adding $\Phi_xR$ and $\Phi R_x$ gives
\[
\widehat\Phi_x
= -\Phi_{aa}R + 2C\widehat\Phi
  - \mathcal F_{0,1}\widehat\Phi
  + \mathcal A(\Phi_aR).
\]
Substituting \eqref{eq:kpz-red-PhiR-aa} for $\Phi_{aa}R$
and \eqref{eq:kpz-red-PhiR-a} for $\Phi_aR$ gives
\begin{align*}
\widehat\Phi_x
&= -\widehat\Phi_{aa}
   - \mathcal A\widehat\Phi_a
   - (\mathcal F_{1,0}+\mathcal A_a)\widehat\Phi
   + 2C\widehat\Phi
   - \mathcal F_{0,1}\widehat\Phi
   + \mathcal A\widehat\Phi_a
   + \mathcal A^2\widehat\Phi.
\end{align*}
The terms $\pm\mathcal A\widehat\Phi_a$ cancel.
By \eqref{eq:kpz-red-Aa-Frs}, the remaining
$\widehat\Phi$ coefficient is
$2C-\mathcal F_{1,0}-\mathcal A_a+\mathcal A^2
-\mathcal F_{0,1}=2(C-\mathcal A_a)=2\widehat C$,
giving
$\widehat\Phi_x=-\widehat\Phi_{aa}
+2\widehat C\widehat\Phi$, which is
\eqref{eq:kpz-red-Ax} with coefficients
$(\widehat C,\widehat\Lambda)$.

\noindent\textit{The dressed adjoint \eqref{eq:kpz-red-At}.}
By the Leibniz rule and the adjoint linear problem
\eqref{eq:kpz-red-At},
\[
\widehat\Phi_t
= -\tfrac13\Phi_{aaa}R
  + C(\Phi_aR)
  + (\Lambda+C_a)\widehat\Phi
  + \Phi R_t.
\]
The right-hand side involves $\Phi_aR$, $\Phi_{aaa}R$,
and $\Phi R_t$, which are expressed in terms of
$\widehat\Phi$ by
\eqref{eq:kpz-red-PhiR-a}--\eqref{eq:kpz-red-PhiR-aaa}
and the computation below.
The resolvent identity \eqref{eq:kpz-red-resolvent-alpha}
with $\beta=t$ gives $\Phi R_t=z\widehat\Phi K_tR$.
Substituting \eqref{eq:kpz-red-Kt} for $K_t$ and
identifying
$z\widehat\Phi\Psi=\mathcal A$,
$z\widehat\Phi\Psi_a=\mathcal F_{0,1}$, and
$z\widehat\Phi\Psi_{aa}=\mathcal F_{0,2}$ via
\eqref{eq:kpz-red-F-rs} gives
\[
\Phi R_t
= \tfrac13\left[\mathcal F_{0,2}\widehat\Phi
  - \mathcal F_{0,1}(\Phi_aR)
  + \mathcal A(\Phi_{aa}R)\right]
  - \mathcal A C\widehat\Phi.
\]
Substituting
\eqref{eq:kpz-red-PhiR-a}--\eqref{eq:kpz-red-PhiR-aaa}
and the expression for $\Phi R_t$ and collecting by
derivatives of $\widehat\Phi$ gives
\[
\widehat\Phi_t
=-\tfrac13\widehat\Phi_{aaa}
+\widehat C\widehat\Phi_a+\eta_0\widehat\Phi,
\]
where the first-order coefficient
$\widehat C = C-\mathcal A_a$ matches $\theta_1$ by the
same application of \eqref{eq:kpz-red-Aa-Frs}.  The
zeroth-order coefficient is
\begin{align*}
\eta_0
=
\Lambda+C_a+C\mathcal A-\mathcal A C
-\tfrac13\mathcal A_{aa}
-\tfrac13\pa_a\mathcal F_{1,0}
-\tfrac13\mathcal F_{2,0}
+\tfrac13\mathcal F_{0,2}
-\tfrac13\mathcal F_{0,1}\mathcal A
+\tfrac13\mathcal A\mathcal F_{1,0}
+\tfrac13\mathcal A\mathcal A_a.
\end{align*}
It suffices to verify that
$\eta_0 = \widehat\Lambda+\widehat C_a
= \Lambda+C_a+\tfrac12\mathcal A_x
-\tfrac12\mathcal A_{aa}$,
where the sign of $\mathcal A_{aa}$ differs from
$\theta_0$ because $\widehat C_a = C_a-\mathcal A_{aa}$.
The verification follows the same mechanism as $\theta_0$:
substituting \eqref{eq:kpz-red-Ax-Frs},
\eqref{eq:kpz-red-Aaa-Frs}, and \eqref{eq:kpz-red-Frs-a},
the $[C,\mathcal A]$ terms cancel, then the
$\mathcal F_{0,2}$, $\mathcal F_{2,0}$, and
$\mathcal F_{1,1}$ terms cancel in pairs, and the
remaining six terms vanish after applying
\eqref{eq:kpz-red-Aa-Frs}.  Since
\[
\widehat\Lambda+\widehat C_a
=
\Lambda+C_a+\tfrac12\mathcal A_x-\tfrac12\mathcal A_{aa},
\]
$\widehat\Phi$ satisfies \eqref{eq:kpz-red-At} with
coefficients $(\widehat C,\widehat\Lambda)$.

\noindent\textit{The dressed kernel.}
Since $K$ and $R=(I-zK)^{-1}$ commute, $\widehat K=KR=RK$.
The identity $(I-\eta\widehat K)(I-zK)=I-(z+\eta)K$ shows
that $(I-\eta\widehat K)^{-1}$ exists whenever
$(I-(z+\eta)K)^{-1}$ does; in particular at $\eta=-z$,
where the condition is vacuous.
The resolvent identity
\eqref{eq:kpz-red-resolvent-alpha} gives, for any
coordinate derivative $\pa_\beta$,
$\pa_\beta\widehat K
=K_\beta R+KR_\beta
=K_\beta R+zKRK_\beta R
=RK_\beta R$.
Substituting \eqref{eq:kpz-red-Ka} gives
$\widehat K_a=R(-\Psi\Phi)R
=-\widehat\Psi\widehat\Phi$.
Substituting \eqref{eq:kpz-red-Kx} and regrouping by
\eqref{eq:kpz-red-RPsi-a} and \eqref{eq:kpz-red-PhiR-a}
gives
\begin{align*}
\widehat K_x
=R(-\Psi_a\Phi+\Psi\Phi_a)R
=-(\widehat\Psi_a+\widehat\Psi\mathcal A)\widehat\Phi
  +\widehat\Psi(\widehat\Phi_a+\mathcal A\widehat\Phi)
=-\widehat\Psi_a\widehat\Phi
  +\widehat\Psi\widehat\Phi_a,
\end{align*}
where the $\pm\widehat\Psi\mathcal A\widehat\Phi$ terms
cancel.  Setting $\beta=t$ gives
$\widehat K_t=RK_tR$.  Substituting
\eqref{eq:kpz-red-Kt} and expanding each factor by
\eqref{eq:kpz-red-RPsi-a}--\eqref{eq:kpz-red-RPsi-aaa}
and
\eqref{eq:kpz-red-PhiR-a}--\eqref{eq:kpz-red-PhiR-aaa}
gives
\begin{align*}
RK_tR
=
\tfrac13(
\widehat\Psi_{aa}\widehat\Phi
-\widehat\Psi_a\widehat\Phi_a
+\widehat\Psi\widehat\Phi_{aa})
-\widehat\Psi C\widehat\Phi
+\tfrac13\widehat\Psi
(2\mathcal A_a+\mathcal F_{1,0}
+\mathcal F_{0,1}-\mathcal A^2)
\widehat\Phi.
\end{align*}
On the other hand, evaluating \eqref{eq:kpz-red-Kt} with
$(\widehat\Psi,\widehat\Phi,\widehat C)$ and using
$\widehat C=C-\mathcal A_a$ gives
\[
\tfrac13(
\widehat\Psi_{aa}\widehat\Phi
-\widehat\Psi_a\widehat\Phi_a
+\widehat\Psi\widehat\Phi_{aa})
-\widehat\Psi(C-\mathcal A_a)\widehat\Phi.
\]
Subtracting this from $RK_tR$ leaves
$\widehat\Psi\left[\tfrac13(\mathcal A^2-
\mathcal F_{0,1}-\mathcal F_{1,0}+\mathcal A_a)\right]
\widehat\Phi$,
which vanishes by \eqref{eq:kpz-red-Aa-Frs}.
\end{proof}

Proposition~\ref{thm:kpz-red-Darboux} shows that the dressed
waves satisfy
\eqref{eq:kpz-red-Lx}--\eqref{eq:kpz-red-Lt} and
\eqref{eq:kpz-red-Ax}--\eqref{eq:kpz-red-At} with
coefficients $(\widehat C,\widehat\Lambda)$.
Theorem~\ref{thm:kpz-red-Darboux-compat} next establishes
that $(\widehat C,\widehat\Lambda)$ satisfies the
compatibility equations
\eqref{eq:kpz-red-constraint}--\eqref{eq:kpz-red-evolution},
and that $\mathcal A$ satisfies
\eqref{eq:kpz-shift-background-kp}.

\begin{theorem}
\label{thm:kpz-red-Darboux-compat}
Under the hypotheses of
Proposition~\ref{thm:kpz-red-Darboux}, the dressed
coefficients $(\widehat C,\widehat\Lambda)$ satisfy the
compatibility equations
\eqref{eq:kpz-red-constraint}--\eqref{eq:kpz-red-evolution}.
In particular, the dressed observable $\mathcal A$ satisfies
\begin{align}
0
&=
\mathcal A_{at}
+\tfrac14\mathcal A_{xx}
+\tfrac1{12}\mathcal A_{aaaa}
+\tfrac12(\mathcal A_a^2)_a
+\tfrac12[\mathcal A_a,\mathcal A_x]
\notag \\
&\quad
-(C\mathcal A_a)_a
-\tfrac12[C,\mathcal A_x]
+\tfrac12[C,\mathcal A_{aa}]
+[\mathcal A_a,\Lambda].
\label{eq:kpz-shift-background-kp}
\end{align}
\end{theorem}

The proof uses the following derivative identities for
$\mathcal A$.

\begin{lemma}
\label{lem:kpz-red-derivative-identities}
Under the hypotheses of
Proposition~\ref{thm:kpz-red-Darboux}, the dressed
observable $\mathcal A=z\Phi R\Psi$ satisfies
\begin{align}
\mathcal A_a
&=\mathcal F_{1,0}-\mathcal A^2+\mathcal F_{0,1},
\label{eq:kpz-red-deriv-Aa} \\
\mathcal A_{aa}
&=\mathcal F_{2,0}+\mathcal F_{0,2}
+2\mathcal F_{1,1}
-\mathcal F_{1,0}\mathcal A
-\mathcal A\mathcal F_{0,1}
-\mathcal A_a\mathcal A
-\mathcal A\mathcal A_a,
\label{eq:kpz-red-deriv-Aaa} \\
\mathcal A_x
&=\mathcal F_{0,2}-\mathcal F_{2,0}
+\mathcal A\mathcal F_{1,0}
-\mathcal F_{0,1}\mathcal A
+2C\mathcal A-2\mathcal A C,
\label{eq:kpz-red-deriv-Ax} \\
\mathcal A_t
&=-\tfrac13\mathcal F_{3,0}
  -\tfrac13\mathcal F_{0,3}
  +\tfrac13(\mathcal F_{0,2}\mathcal A
  -\mathcal F_{0,1}\mathcal F_{1,0}
  +\mathcal A\mathcal F_{2,0})
\notag \\
&\quad
  +C\mathcal F_{1,0}+\mathcal F_{0,1}C
  -\mathcal A C\mathcal A
  +(\Lambda+C_a)\mathcal A
  -\mathcal A\Lambda.
\label{eq:kpz-red-deriv-At}
\end{align}
\end{lemma}

\begin{proof}
The identity \eqref{eq:kpz-red-deriv-Aa} is
\eqref{eq:kpz-red-Aa-Frs}, and
\eqref{eq:kpz-red-deriv-Aaa} and
\eqref{eq:kpz-red-deriv-Ax} are
\eqref{eq:kpz-red-Aaa-Frs} and
\eqref{eq:kpz-red-Ax-Frs}, all established in the proof
of Proposition~\ref{thm:kpz-red-Darboux}.
Differentiating $\mathcal A=z\Phi R\Psi$ in $t$ gives
three contributions:
\begin{align*}
z\Phi_tR\Psi
&=-\tfrac13\mathcal F_{3,0}
  +C\mathcal F_{1,0}
  +(\Lambda+C_a)\mathcal A, \\
z\Phi R_t\Psi
&=
\tfrac13(\mathcal F_{0,2}\mathcal A
  -\mathcal F_{0,1}\mathcal F_{1,0}
  +\mathcal A\mathcal F_{2,0})
  -\mathcal A C\mathcal A, \\
z\Phi R\Psi_t
&=-\tfrac13\mathcal F_{0,3}
  +\mathcal F_{0,1}C
  -\mathcal A\Lambda,
\end{align*}
where the first uses \eqref{eq:kpz-red-At}, the second uses
the resolvent identity \eqref{eq:kpz-red-resolvent-alpha}
and \eqref{eq:kpz-red-Kt}, and the third uses
\eqref{eq:kpz-red-Lt}.  Adding them gives
\eqref{eq:kpz-red-deriv-At}.
\end{proof}

The proof of
Theorem~\ref{thm:kpz-red-Darboux-compat} reduces to the
following covariance identity.

\begin{lemma}
\label{lem:kpz-red-residual-identity}
Under the hypotheses of
Proposition~\ref{thm:kpz-red-Darboux}, define the
compatibility residuals
\begin{align}
\mathscr E_1(C,\Lambda)
&\defeq 2\Lambda_a+C_x+C_{aa},
  \label{eq:kpz-E1-residual} \\
\mathscr E_2(C,\Lambda)
&\defeq 2C_t-\Lambda_x+\Lambda_{aa}
+\tfrac23 C_{aaa}-2C_aC+2[C,\Lambda],
  \label{eq:kpz-E2-residual}
\end{align}
so that the compatibility equations
\eqref{eq:kpz-red-constraint}--\eqref{eq:kpz-red-evolution}
are $\mathscr E_1(C,\Lambda)=0$ and
$\mathscr E_2(C,\Lambda)=0$.  The dressed coefficients
$(\widehat C,\widehat\Lambda)$ defined by
\eqref{eq:kpz-red-hat-CL} satisfy
\begin{align}
\mathscr E_1(\widehat C,\widehat\Lambda)
&=\mathscr E_1(C,\Lambda),
  \label{eq:kpz-E1-invariant} \\
\mathscr E_2(\widehat C,\widehat\Lambda)
&=\mathscr E_2(C,\Lambda)
  +[\mathcal A,\mathscr E_1(C,\Lambda)].
  \label{eq:kpz-E2-covariant}
\end{align}
\end{lemma}

\begin{proof}
Substituting \eqref{eq:kpz-red-hat-CL} into
\eqref{eq:kpz-E1-residual} gives
\begin{align*}
\mathscr E_1(\widehat C,\widehat\Lambda)
-\mathscr E_1(C,\Lambda)
&=(\mathcal A_x+\mathcal A_{aa})_a
-(\mathcal A_a)_x-(\mathcal A_a)_{aa}
=0,
\end{align*}
which proves \eqref{eq:kpz-E1-invariant}.

For \eqref{eq:kpz-E2-covariant}, substitute
\eqref{eq:kpz-red-hat-CL} into
\eqref{eq:kpz-E2-residual} and subtract
$\mathscr E_2(C,\Lambda)$.  The difference is
\begin{align}
&\mathscr E_2(\widehat C,\widehat\Lambda)
-\mathscr E_2(C,\Lambda)
\notag \\
&= -2\mathcal A_{at}
-\tfrac12\mathcal A_{xx}
-\tfrac16\mathcal A_{aaaa}
+2\mathcal A_{aa}C
+2C_a\mathcal A_a
-2\mathcal A_{aa}\mathcal A_a
\notag \\
&\quad
+[C,\mathcal A_x+\mathcal A_{aa}]
-2[\mathcal A_a,\Lambda]
-[\mathcal A_a,\mathcal A_x+\mathcal A_{aa}].
\label{eq:kpz-red-E2-difference}
\end{align}
Substituting
\eqref{eq:kpz-red-deriv-Aa}--\eqref{eq:kpz-red-deriv-At}
into
\eqref{eq:kpz-red-E2-difference}, expressing
$\mathcal A_{aaaa}$ and $\mathcal A_{at}$ by
differentiating these identities in $a$, and reducing
$\mathcal A_{xx}$ via \eqref{eq:kpz-red-Lx} and
\eqref{eq:kpz-red-Ax}, the expression reorganizes as
\begin{equation}
\label{eq:kpz-red-E2-covariant-difference}
\mathscr E_2(\widehat C,\widehat\Lambda)
-\mathscr E_2(C,\Lambda)
=[\mathcal A,2\Lambda_a+C_x+C_{aa}]
+\tfrac16(\mathcal N_{aa})_{aa}
+\tfrac12(\mathcal N_x)_x
+2\pa_a\mathcal N_t,
\end{equation}
where $\mathcal N_{aa}$, $\mathcal N_x$, $\mathcal N_t$
denote the differences between the left- and right-hand
sides of
\eqref{eq:kpz-red-deriv-Aaa}--\eqref{eq:kpz-red-deriv-At}
respectively.  Since these identities hold, the defects
$\mathcal N_{aa}$, $\mathcal N_x$, $\mathcal N_t$ all
vanish, and the
surviving terms are
\[
[\mathcal A,2\Lambda_a+C_x+C_{aa}]
=[\mathcal A,\mathscr E_1(C,\Lambda)],
\]
where the last equality uses \eqref{eq:kpz-E1-residual}.
This proves \eqref{eq:kpz-E2-covariant}.
\end{proof}

\begin{proof}[Proof of Theorem~\ref{thm:kpz-red-Darboux-compat}]
Lemma~\ref{lem:kpz-red-residual-identity} gives
$\mathscr E_1(\widehat C,\widehat\Lambda)
=\mathscr E_1(C,\Lambda)=0$ and also
$\mathscr E_2(\widehat C,\widehat\Lambda)
=\mathscr E_2(C,\Lambda)+[\mathcal A,0]=0$,
so $(\widehat C,\widehat\Lambda)$ satisfies
\eqref{eq:kpz-red-constraint}--\eqref{eq:kpz-red-evolution}.
Since both $(C,\Lambda)$ and
$(\widehat C,\widehat\Lambda)$ satisfy
\eqref{eq:kpz-red-constraint}--\eqref{eq:kpz-red-evolution},
the expansion \eqref{eq:kpz-red-E2-difference} vanishes.
Multiplying \eqref{eq:kpz-red-E2-difference} by
$-\tfrac12$ gives
\begin{align*}
0
&=
\mathcal A_{at}
+\tfrac14\mathcal A_{xx}
+\tfrac1{12}\mathcal A_{aaaa}
-\mathcal A_{aa}C
-C_a\mathcal A_a
+\mathcal A_{aa}\mathcal A_a \\
&\quad
-\tfrac12[C,\mathcal A_x+\mathcal A_{aa}]
+[\mathcal A_a,\Lambda]
+\tfrac12[\mathcal A_a,\mathcal A_x+\mathcal A_{aa}].
\end{align*}
The purely $\mathcal A$-dependent nonlinear terms satisfy
\[
\mathcal A_{aa}\mathcal A_a
+\tfrac12[\mathcal A_a,\mathcal A_x+\mathcal A_{aa}]
=
\tfrac12(\mathcal A_a^2)_a
+\tfrac12[\mathcal A_a,\mathcal A_x].
\]
The identity
$-\mathcal A_{aa}C-C_a\mathcal A_a
=-(C\mathcal A_a)_a+[C,\mathcal A_{aa}]$
then gives \eqref{eq:kpz-shift-background-kp}.
\end{proof}

\begin{corollary}
\label{cor:kpz-shift-vacuum-kp}
In the vacuum background $C=0$ and $\Lambda=0$, the dressed
observable satisfies
\begin{equation}
\label{eq:kpz-shift-vacuum-kp}
\mathcal A_{at}
+\tfrac14\mathcal A_{xx}
+\tfrac1{12}\mathcal A_{aaaa}
+\tfrac12(\mathcal A_a^2)_a
+\tfrac12[\mathcal A_a,\mathcal A_x] = 0.
\end{equation}
\end{corollary}

\section{Bilinear Equations}
\label{sec:kpz-red-trace-defect}

The Fredholm determinant $F \defeq \det_H(I - zK)$ encodes
the multipoint distributions of the KPZ fixed point
(Chapter~\ref{ch:the-kpz-fixed-point}).  The proposition
below extracts a bilinear identity for $F$ from the matrix
Darboux theory through the relation
$\pa_a \log F = \operatorname{tr}_E \mathcal A$,
providing a scalar companion to
\eqref{eq:kpz-shift-background-kp}.

\begin{proposition}
\label{prop:kpz-shift-trace-defect}
Under the hypotheses of
Theorem~\ref{thm:kpz-red-Darboux-compat}, assume
additionally that $\dim E<\infty$ and that $H$ is a real
separable Hilbert space with $K\in\mathfrak S_1(H)$ smooth
in trace norm.  Define $F\defeq\det_H(I-zK)$ and assume
that $F\to 1$ and all partial derivatives of $\log F$
vanish as $a\to+\infty$, for each fixed $(x,t)$.
Then
\begin{align}
(D_aD_t+\tfrac14D_x^2+\tfrac1{12}D_a^4)F\cdot F
=
2F^2\left[
\tfrac12\left(
(\operatorname{tr}_E\mathcal A_a)^2
-\operatorname{tr}_E(\mathcal A_a^2)
\right)
+\operatorname{tr}_E(C\mathcal A_a)
\right].
\label{eq:kpz-shift-integrated-Hirota}
\end{align}
\end{proposition}

\begin{proof}
Jacobi's formula and Sylvester's trace identity give
$\pa_a \log F = \operatorname{tr}_E \mathcal A$;
write $q \defeq \operatorname{tr}_E \mathcal A = \pa_a \log F$.
Taking $\operatorname{tr}_E$ of
\eqref{eq:kpz-shift-background-kp} removes all commutator
terms, giving
\begin{equation}
\label{eq:kpz-shift-traced-A-equation}
q_{at}
+\tfrac14q_{xx}
+\tfrac1{12}q_{aaaa}
+\tfrac12\pa_a\operatorname{tr}_E(\mathcal A_a^2)
-\pa_a\operatorname{tr}_E(C\mathcal A_a)
=0.
\end{equation}
Since $q=\pa_a\log F$, every term in
\eqref{eq:kpz-shift-traced-A-equation} is an
$a$-derivative.  Integrating in $a$ gives
\begin{equation}
\label{eq:kpz-shift-differentiated-Hirota}
(\log F)_{at}
+\tfrac14(\log F)_{xx}
+\tfrac1{12}(\log F)_{aaaa}
+\tfrac12\operatorname{tr}_E(\mathcal A_a^2)
-\operatorname{tr}_E(C\mathcal A_a)
=0,
\end{equation}
where the integration constant vanishes by the
boundary assumption $F\to 1$.
The standard Hirota identities
\begin{gather*}
\frac{D_aD_tF\cdot F}{2F^2}=(\log F)_{at},
\qquad
\frac{D_x^2F\cdot F}{2F^2}=(\log F)_{xx},
\qquad
\frac{D_a^4F\cdot F}{2F^2}
=(\log F)_{aaaa}+6(\log F)_{aa}^2
\end{gather*}
convert \eqref{eq:kpz-shift-differentiated-Hirota} into
\eqref{eq:kpz-shift-integrated-Hirota}, because
$(\log F)_{aa}=q_a=\operatorname{tr}_E\mathcal A_a$
and the term
$6(\operatorname{tr}_E\mathcal A_a)^2$ produces
the trace-square difference.
\end{proof}

%% file: chapters/6-discrete-particle-models.tex
\chapter{Discrete Particle Models}
\label{ch:discrete-particle-models}

{
  \setlength{\parskip}{0pt}
}
\section{The Directed-Path Propagator}\label{sec:directed-path}

The admissible propagators required by the product graph
construction of \S\ref{sec:product-graph} arise, for many models in this
chapter, from a common source: a signed sum over
threshold-constrained directed paths through the layer
indices $1, \dotsc, m$, weighted by transition kernels (typically random-walk transition
kernels).  This section
constructs the \emph{directed-path propagator} in abstract
terms and establishes $\mathcal{T}$-covariance and
$\mathcal{T}$-splitting under hypotheses satisfied by all
models treated below. It then derives the
$\mathcal{S}_k$-compatibility condition under additional
hypotheses that hold for the majority of
models.\footnote{The reader may prefer to proceed directly
to the first model verification and refer back to this section as needed.}

Let $\mathcal{G}^m$ be a product graph as in \S\ref{sec:product-graph},
with vertex set $\mathcal{V}^m$, distinguished shift
$\mathcal{T}$, and additional shifts
$\{\mathcal{S}_k\}_{k \in \mathcal I}$ for some index set
$\mathcal I \subseteq \N$.  Throughout this chapter, the
coefficient space is $E = \R$, so all
endomorphism spaces $\End(E^m)$ of the abstract
framework specialize to $\End(\R^m)$.
The construction requires two
pieces of data extracted from each vertex
$u \in \mathcal{V}^m$.

The first is a \emph{threshold vector} recording the
$T$-coordinate of each layer: if $u = (u_1, \dotsc, u_m)$
and the base-graph vertex $u_i$ has $T$-coordinate
$a_i \in \Z$, then define
$\mathbf{a}(u) \defeq (a_1, \dotsc, a_m)$.  Note that since
$\mathcal{T}$ acts diagonally by $T$ in each component, $\mathbf{a}(\mathcal{T}u)
= \mathbf{a}(u) + \mathbf{1}$, where $\mathbf 1 = (1, \dots, 1)$.

The second is a collection of \emph{transition kernels}
$[Q_u]_{\ell, \ell'} \colon \mathbb{Z} \times \mathbb{Z}
\to \mathbb{R}$, one for each pair of layers
$(\ell, \ell')$ with $1 \le \ell < \ell' \le m$.  Two
standing hypotheses are imposed throughout:
\begin{enumerate}[label=(\roman*), leftmargin=*]
\item \emph{Translation invariance:}
$[Q_u]_{\ell,\ell'}(r{+}1, r'{+}1)
= [Q_u]_{\ell,\ell'}(r, r')$ for all
$\ell, \ell', r, r'$.
\item \emph{$\mathcal{T}$-invariance:}
$[Q_{\mathcal Tu}]_{\ell,\ell'}
= [Q_u]_{\ell,\ell'}$ for all $\ell, \ell'$.
\end{enumerate}
For notational convenience, we suppress the subscript $u$
and write $Q_{\ell,\ell'}$ for $[Q_u]_{\ell,\ell'}$.
The propagator is built from directed paths through the
layers $1, \dotsc, m$.

\begin{definition}\label{def:directed-path-propagator}
Fix $u \in \mathcal{V}^m$.  For layer indices $i > j$
and spatial endpoints $r, r' \in \mathbb{Z}$, a
\emph{threshold-constrained directed path} from
$(j, r')$ to $(i, r)$ is a pair of sequences
\[
j = \ell_0 < \ell_1 < \cdots < \ell_k = i,
\qquad
\xi_0 = r', \xi_1, \dotsc, \xi_{k-1},
\xi_k = r,
\]
with $k \ge 1$ steps and internal vertices constrained
by the thresholds:
$\xi_s \le a_{\ell_s}(u)$ for each
$s \in \{1, \dotsc, k{-}1\}$. The
\emph{weight} of such a path $\gamma$ is
\[
w(\gamma) \defeq \prod_{s=0}^{k-1}
\bigl(-Q_{\ell_s, \ell_{s+1}}
(\xi_s, \xi_{s+1})\bigr).
\]
The \emph{directed-path propagator} is the strictly
lower-triangular family
$B_u \colon \Z^m \times \Z^m \to \End(\R^m)$ with
entries for $i > j$
\begin{align}
[B_u]_{ij}(r, r') &\defeq
\sum_{\gamma \colon (j, r') \to (i, r)} w(\gamma) \nonumber \\
&= \sum_{k=1}^{i-j}
\sum_{j = \ell_0 < \ell_1 < \cdots < \ell_k = i}
\sum_{\substack{\xi_1 \le a_{\ell_1}(u) \\ \vdots \\
\xi_{k-1} \le a_{\ell_{k-1}}(u)}}
(-1)^k \prod_{s=0}^{k-1}
Q_{\ell_s, \ell_{s+1}}
(\xi_s, \xi_{s+1}), \label{eq:B-directed-path}
\end{align}
with $\xi_0 = r'$ and $\xi_k = r$, and with
$[B_u]_{ij} = 0$ for $i \le j$.
\end{definition}

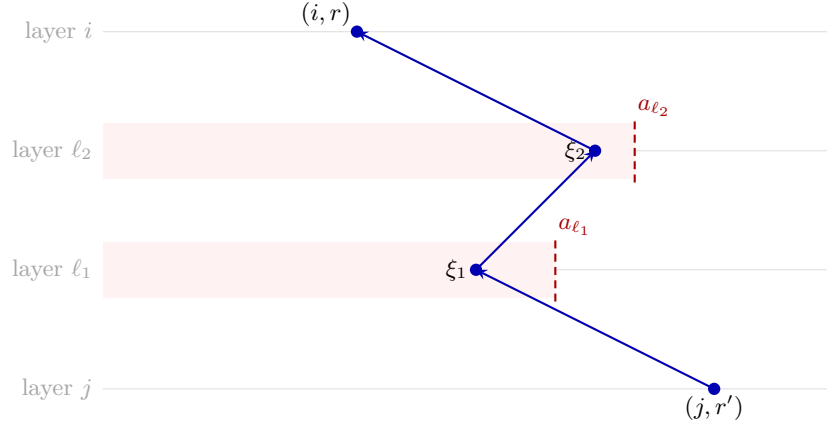
\begin{figure}[ht]
\centering
\begin{tikzpicture}[>=stealth, scale=1.05,
  every node/.style={font=\footnotesize}]

  \def\ls{1.5}
  \def\xmin{-1.2}
  \def\xmax{8}

  \foreach \y in {0,1,2,3} {
    \draw[gray!25] (\xmin, \y*\ls) -- (\xmax, \y*\ls);
  }

  \node[left, gray!70] at (\xmin, 0)
    {layer $j$};
  \node[left, gray!70] at (\xmin, \ls)
    {layer $\ell_1$};
  \node[left, gray!70] at (\xmin, 2*\ls)
    {layer $\ell_2$};
  \node[left, gray!70] at (\xmin, 3*\ls)
    {layer $i$};

  \fill[red!5] (\xmin, \ls-0.35)
    rectangle (4.5, \ls+0.35);
  \fill[red!5] (\xmin, 2*\ls-0.35)
    rectangle (5.5, 2*\ls+0.35);

  \draw[red!65!black, densely dashed,
    line width=0.8pt]
    (4.5, \ls-0.4) -- (4.5, \ls+0.4);
  \node[red!65!black, above right, inner sep=1pt]
    at (4.5, \ls+0.4) {$a_{\ell_1}$};

  \draw[red!65!black, densely dashed,
    line width=0.8pt]
    (5.5, 2*\ls-0.4) -- (5.5, 2*\ls+0.4);
  \node[red!65!black, above right, inner sep=1pt]
    at (5.5, 2*\ls+0.4) {$a_{\ell_2}$};

  \coordinate (P0) at (6.5, 0);
  \coordinate (P1) at (3.5, \ls);
  \coordinate (P2) at (5,   2*\ls);
  \coordinate (P3) at (2,   3*\ls);

  \draw[blue!70!black, thick, -stealth]
    (P0) -- (P1);
  \draw[blue!70!black, thick, -stealth]
    (P1) -- (P2);
  \draw[blue!70!black, thick, -stealth]
    (P2) -- (P3);

  \fill[blue!70!black] (P0) circle (2.2pt);
  \fill[blue!70!black] (P1) circle (2.2pt);
  \fill[blue!70!black] (P2) circle (2.2pt);
  \fill[blue!70!black] (P3) circle (2.2pt);

  \node[below, inner sep=2pt] at (P0)
    {$(j, r')$};
  \node[left, inner sep=3pt] at (P1)
    {$\xi_1$};
  \node[left, inner sep=3pt] at (P2)
    {$\xi_2$};
  \node[above left, inner sep=2pt] at (P3)
    {$(i, r)$};

\end{tikzpicture}
\caption{A threshold-constrained directed path from
$(j, r')$ to $(i, r)$ with $k = 3$ steps through
layers $j < \ell_1 < \ell_2 < i$.  The internal
vertices $\xi_1, \xi_2$ are constrained to lie at or
to the left of the thresholds $a_{\ell_1}, a_{\ell_2}$
(dashed lines); the endpoints $r', r$ are
unconstrained.  Each step carries weight
$(-Q_{\ell_s, \ell_{s+1}})(\xi_s, \xi_{s+1})$.}
\label{fig:directed-path}
\end{figure}

\begin{remark}
    For $x, y \in \mathcal{V}^m$, we write
$[B_u(x, y)]_{ij} \defeq [B_u]_{ij}(a_i(x), a_j(y))$,
identifying $B_u$ with a map
$\mathcal{V}^m \times \mathcal{V}^m \to \End(\R^m)$ as
required by the product graph construction of \S\ref{sec:product-graph}.
In the proofs below, we pass freely between the spatial
entries $[B_u]_{ij}(r, r')$ and the graph-level
evaluations $B_u(x, y)$.
\end{remark}
\begin{remark}
The path length ranges from $k = 1$ (a single step
from layer $j$ directly to layer $i$) to $k = i - j$
(the path visits every intermediate layer), and the sum
over layer sequences is finite.  The sum over internal
spatial positions may be infinite and requires sufficient
decay of the transition kernels; this regularity
condition is verified model by model.  The path
interpretation gives the construction its name: each
term in~\eqref{eq:B-directed-path} corresponds to a
directed path through the layers, weighted by transition
kernels and constrained to remain at or below the
thresholds $\mathbf{a}(u)$ at every internal vertex
(Figure~\ref{fig:directed-path}).
\end{remark}

\begin{remark}[Neumann series]\label{rem:Neumann-series}
Let $L$ be the strictly upper-triangular block operator
on $\ell^2(\{1, \dotsc, m\} \times \mathbb{Z})$ with
blocks $L_{\ell,\ell'} = Q_{\ell,\ell'}$ for
$\ell < \ell'$, and let $\bar{\chi}_{\mathbf{a}}$ be
the block-diagonal projection
$\bar{\chi}_{\mathbf{a}}(\ell, r)
= \mathbf{1}_{r \le a_\ell(u)}$.  The compressed
operator
$\bar{\chi}_{\mathbf{a}} L \bar{\chi}_{\mathbf{a}}$
is nilpotent (strictly upper triangular in the block
indices), so
\[
(I + \bar{\chi}_{\mathbf{a}} L
\bar{\chi}_{\mathbf{a}})^{-1}
= \sum_{k=0}^{m-1}
(-\bar{\chi}_{\mathbf{a}} L
\bar{\chi}_{\mathbf{a}})^k.
\]
Expanding the right-hand side in the block indices,
the $(j,i)$-block with $j < i$ is a sum over chains
$j = \ell_0 < \ell_1 < \cdots < \ell_k = i$ with the
same weights and threshold constraints as
Definition~\ref{def:directed-path-propagator}.  On the
cutoff region $r \le a_i$, $r' \le a_j$, this gives
the identification
\[
\delta_{ij}\delta_{r,r'}
+ [B_u]_{ij}(r, r')
= [(I + \bar{\chi}_{\mathbf{a}} L
\bar{\chi}_{\mathbf{a}})^{-1}]_{ji}(r', r)
\qquad (r \le a_i, r' \le a_j),
\]
where the index reversal ($(j,i)$ on the right vs.\
$(i,j)$ on the left) and the argument reversal
($(r', r)$ vs.\ $(r, r')$) reflect the passage from
the upper-triangular inverse to the lower-triangular
propagator.  The
definition~\eqref{eq:B-directed-path} extends this to
all $r, r' \in \mathbb{Z}$ by removing the endpoint
constraints; this extension is what enters the product
graph construction.
\end{remark}

\subsection{$\mathcal{T}$-covariance and
$\mathcal{T}$-splitting}
\label{sec:T-cov-split}

The two properties established in this subsection hold
for every model in this chapter: they depend only on the
standing hypotheses (i)--(ii) and the combinatorial
structure of the directed-path sum.

\begin{proposition}\label{prop:directed-path-T}
Under the standing hypotheses
\textup{(i)--(ii)}, the directed-path propagator
of Definition~\ref{def:directed-path-propagator}
satisfies $\mathcal{T}$-covariance and
$\mathcal{T}$-splitting on the product graph
$\mathcal{G}^m$.\footnote{The formal rearrangements
in the proof are justified by the regularity condition
of Definition~\ref{def:admissible-prop}, verified
for each model below.}
\end{proposition}

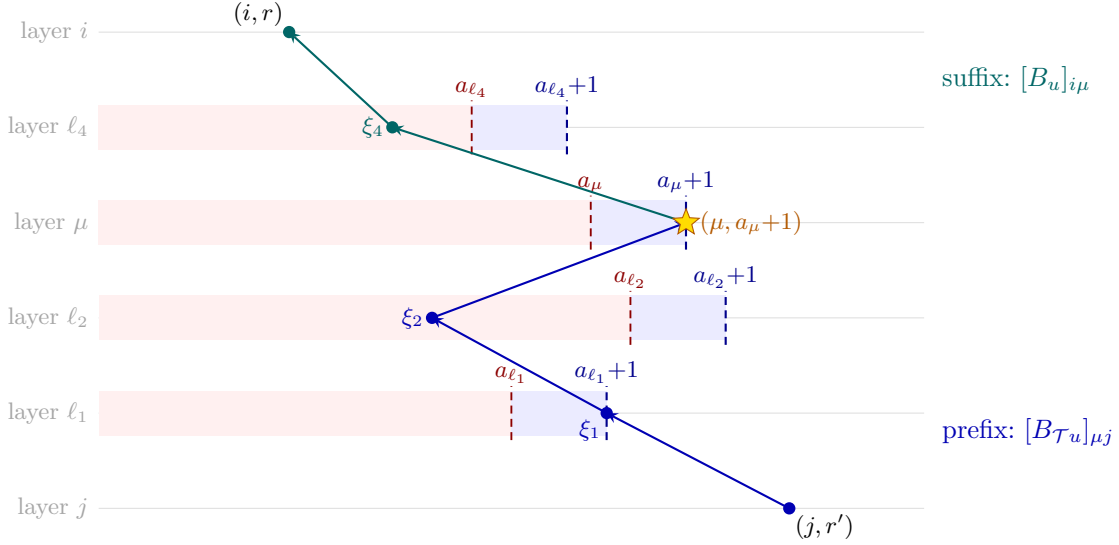
\begin{figure}[ht]
\centering
\begin{tikzpicture}[>=stealth, scale=1.05,
  every node/.style={font=\footnotesize}]
 
  \def\ls{1.2}
  \def\xmin{-1.2}
  \def\xmax{9.2}
 
 
  \foreach \y in {0,1,2,3,4,5} {
    \draw[gray!25] (\xmin, \y*\ls) -- (\xmax, \y*\ls);
  }
 
  \node[left, gray!70] at (\xmin, 0)
    {layer $j$};
  \node[left, gray!70] at (\xmin, \ls)
    {layer $\ell_1$};
  \node[left, gray!70] at (\xmin, 2*\ls)
    {layer $\ell_2$};
  \node[left, gray!70] at (\xmin, 3*\ls)
    {layer $\mu$};
  \node[left, gray!70] at (\xmin, 4*\ls)
    {layer $\ell_4$};
  \node[left, gray!70] at (\xmin, 5*\ls)
    {layer $i$};
 
  \fill[red!6] (\xmin, \ls-0.28)
    rectangle (4.0, \ls+0.28);
  \fill[red!6] (\xmin, 2*\ls-0.28)
    rectangle (5.5, 2*\ls+0.28);
  \fill[red!6] (\xmin, 3*\ls-0.28)
    rectangle (5.0, 3*\ls+0.28);
  \fill[red!6] (\xmin, 4*\ls-0.28)
    rectangle (3.5, 4*\ls+0.28);
 
  \fill[blue!8] (4.0, \ls-0.28)
    rectangle (5.2, \ls+0.28);
  \fill[blue!8] (5.5, 2*\ls-0.28)
    rectangle (6.7, 2*\ls+0.28);
  \fill[blue!8] (5.0, 3*\ls-0.28)
    rectangle (6.2, 3*\ls+0.28);
  \fill[blue!8] (3.5, 4*\ls-0.28)
    rectangle (4.7, 4*\ls+0.28);
 
  \draw[red!55!black, densely dashed,
    line width=0.7pt]
    (4.0, \ls-0.34) -- (4.0, \ls+0.34);
  \draw[red!55!black, densely dashed,
    line width=0.7pt]
    (5.5, 2*\ls-0.34) -- (5.5, 2*\ls+0.34);
  \draw[red!55!black, densely dashed,
    line width=0.7pt]
    (5.0, 3*\ls-0.34) -- (5.0, 3*\ls+0.34);
  \draw[red!55!black, densely dashed,
    line width=0.7pt]
    (3.5, 4*\ls-0.34) -- (3.5, 4*\ls+0.34);
 
  \draw[blue!55!black, line width=0.8pt, densely dashed]
    (5.2, \ls-0.34) -- (5.2, \ls+0.34);
  \draw[blue!55!black, line width=0.8pt, densely dashed]
    (6.7, 2*\ls-0.34) -- (6.7, 2*\ls+0.34);
  \draw[blue!55!black, line width=0.8pt, densely dashed]
    (6.2, 3*\ls-0.34) -- (6.2, 3*\ls+0.34);
  \draw[blue!55!black, line width=0.8pt, densely dashed]
    (4.7, 4*\ls-0.34) -- (4.7, 4*\ls+0.34);
 
  \node[red!55!black, above, inner sep=1pt]
    at (4.0, \ls+0.34) {$a_{\ell_1}$};
  \node[blue!55!black, above, inner sep=1pt]
    at (5.2, \ls+0.34) {$a_{\ell_1}{+}1$};
 
  \node[red!55!black, above, inner sep=1pt]
    at (5.5, 2*\ls+0.34) {$a_{\ell_2}$};
  \node[blue!55!black, above, inner sep=1pt]
    at (6.7, 2*\ls+0.34) {$a_{\ell_2}{+}1$};
 
  \node[red!55!black, above, inner sep=1pt]
    at (5.0, 3*\ls+0.34) {$a_\mu$};
  \node[blue!55!black, above, inner sep=1pt]
    at (6.2, 3*\ls+0.34) {$a_\mu{+}1$};
 
  \node[red!55!black, above, inner sep=1pt]
    at (3.5, 4*\ls+0.34) {$a_{\ell_4}$};
  \node[blue!55!black, above, inner sep=1pt]
    at (4.7, 4*\ls+0.34) {$a_{\ell_4}{+}1$};
 
  \coordinate (P0) at (7.5,  0);         
  \coordinate (P1) at (5.2,  \ls);       
  \coordinate (P2) at (3.0,  2*\ls);     
  \coordinate (P3) at (6.2,  3*\ls);     
  \coordinate (P4) at (2.5,  4*\ls);     
  \coordinate (P5) at (1.2,  5*\ls);     
 
  \draw[blue!70!black, thick, -stealth]
    (P0) -- (P1);
  \draw[blue!70!black, thick, -stealth]
    (P1) -- (P2);
  \draw[blue!70!black, thick, -stealth]
    (P2) -- (P3);
 
  \draw[teal!80!black, thick, -stealth]
    (P3) -- (P4);
  \draw[teal!80!black, thick, -stealth]
    (P4) -- (P5);
 
  \fill[blue!70!black] (P0) circle (2.2pt);
  \fill[blue!70!black] (P1) circle (2.2pt);
  \fill[blue!70!black] (P2) circle (2.2pt);
  \fill[teal!80!black] (P4) circle (2.2pt);
  \fill[teal!80!black] (P5) circle (2.2pt);
 
  \node[star, star points=5, star point ratio=2.3,
    fill=yellow!80!orange, draw=orange!70!black,
    line width=0.4pt,
    inner sep=1.6pt]
    at (P3) {};
 
  \node[below right, inner sep=2pt] at (P0)
    {$(j, r')$};
  \node[above left, inner sep=2pt] at (P5)
    {$(i, r)$};
  \node[below left, inner sep=2pt,
    blue!70!black] at (P1)
    {$\xi_1$};
  \node[left, inner sep=3pt,
    blue!70!black] at (P2)
    {$\xi_2$};
  \node[left, inner sep=3pt, teal!80!black]
    at (P4) {$\xi_4$};
 
  \node[right, inner sep=5pt, orange!70!black]
    at (P3) {$(\mu, a_\mu{+}1)$};
 
  \node[blue!70!black, anchor=west]
    at (\xmax+0.1, 0.8*\ls)
    {\small prefix: $[B_{\mathcal{T}u}]_{\mu j}$};
  \node[teal!80!black, anchor=west]
    at (\xmax+0.1, 4.5*\ls)
    {\small suffix: $[B_u]_{i\mu}$};
 
\end{tikzpicture}
\caption{The $\mathcal{T}$-splitting decomposition.
The difference $B_{\mathcal{T}u} - B_u$ counts paths
that visit the new boundary
$\xi_s = a_{\ell_s} + 1$ (dashed blue lines) at some
internal layer.  Each such path splits at its \emph{last}
visit to the new boundary, at layer~$\mu$ (gold star),
into a \emph{suffix} from $(\mu, a_\mu{+}1)$ to
$(i, r)$ with internal vertices constrained by the
original thresholds $\mathbf a$ (counted
by $[B_u]_{i\mu}$) and a
\emph{prefix} from $(j, r')$ to
$(\mu, a_\mu{+}1)$ with internal vertices
constrained by the relaxed thresholds
$\mathbf a + \mathbf 1$ (counted by
$[B_{\mathcal{T}u}]_{\mu j}$).}
\label{fig:T-splitting}
\end{figure}

\begin{proof}
\noindent\textit{$\mathcal{T}$-covariance.}
For $i \le j$, both sides vanish by strict lower
triangularity.  Fix $i > j$; we must show
\[
[B_{\mathcal{T}u}]_{ij}(r{+}1, r'{+}1)
= [B_u]_{ij}(r, r')
\]
for all $r, r' \in \mathbb{Z}$.  The propagator
$B_{\mathcal{T}u}$ is built from the same transition
kernels as $B_u$ (by $\mathcal{T}$-invariance) but
with thresholds shifted to $\mathbf{a} + \mathbf{1}$.
A path contributing to
$[B_{\mathcal{T}u}]_{ij}(r{+}1, r'{+}1)$ has
endpoints $\xi_0 = r'{+}1$, $\xi_k = r{+}1$ and
internal vertices $\xi_s \le a_{\ell_s} + 1$.  The
uniform translation $\xi_s \mapsto \xi_s - 1$ sends
this to a path with endpoints $\xi_0 = r'$,
$\xi_k = r$ and internal vertices
$\xi_s \le a_{\ell_s}$, contributing to
$[B_u]_{ij}(r, r')$.  Translation invariance of
$Q_{\ell,\ell'}$ ensures that every step weight is
unchanged under the shift, so the map
$\xi_s \mapsto \xi_s - 1$ is a weight-preserving
bijection between the two path sums.

\medskip
\noindent\textit{$\mathcal{T}$-splitting.}
We must show that for all $i > j$ and
$r, r' \in \mathbb{Z}$,
\begin{equation}\label{eq:T-splitting-spatial}
[B_{\mathcal{T}u}]_{ij}(r, r')
- [B_u]_{ij}(r, r')
= \sum_{j < \mu < i}
[B_u]_{i\mu}(r, a_\mu{+}1)
[B_{\mathcal{T}u}]_{\mu j}(a_\mu{+}1, r').
\end{equation}
The propagators $B_{\mathcal{T}u}$ and $B_u$ share
the same transition kernels and differ only in their
thresholds: $\mathbf{a} + \mathbf{1}$ versus
$\mathbf{a}$.  Since $\xi_s \le a_{\ell_s}$ implies
$\xi_s \le a_{\ell_s} + 1$, every
$\mathbf{a}$-admissible path is
$(\mathbf{a}+\mathbf{1})$-admissible.  The
difference on the left therefore sums over exactly
those paths
from $(j, r')$ to $(i, r)$ that are admissible under
the relaxed thresholds but visit the new boundary: at
least one internal vertex satisfies
$\xi_s = a_{\ell_s} + 1$.

Fix such a path $\gamma$, and let $s^*$ be the index
of its last visit to the new boundary:
\[
s^* \defeq \max\bigl\{s \in \{1, \dotsc, k{-}1\} :
\xi_s = a_{\ell_s} + 1\bigr\},
\qquad \mu \defeq \ell_{s^*}.
\]
The path splits at position $s^*$ into two pieces
(Figure~\ref{fig:T-splitting}).

\noindent\textit{Suffix.}
The subpath from $(\mu, a_\mu{+}1)$ to $(i, r)$,
visiting layers $\ell_{s^*}, \dotsc, \ell_k$.  Since
$s^*$ is maximal, every subsequent internal vertex
satisfies $\xi_s \le a_{\ell_s}$, so this subpath is
counted by $[B_u]_{i\mu}(r, a_\mu{+}1)$.

\noindent\textit{Prefix.}
The subpath from $(j, r')$ to $(\mu, a_\mu{+}1)$,
visiting layers $\ell_0, \dotsc, \ell_{s^*}$.  Its
internal vertices satisfy
$\xi_s \le a_{\ell_s} + 1$, so the subpath is counted
by $[B_{\mathcal{T}u}]_{\mu j}(a_\mu{+}1, r')$.

Conversely, given $\mu$ with $j < \mu < i$, any pair
of an $\mathbf{a}$-admissible suffix from
$(\mu, a_\mu{+}1)$ to $(i, r)$ and an
$(\mathbf{a}+\mathbf{1})$-admissible prefix from
$(j, r')$ to $(\mu, a_\mu{+}1)$ concatenates to a
unique path in the difference.  Indeed, the suffix
visits only original thresholds beyond $\mu$, so
$\mu$ is necessarily the last boundary-visiting layer
and the decomposition is injective.  Summing the
product of suffix and prefix weights over $\mu$
gives~\eqref{eq:T-splitting-spatial}.
\end{proof}

\subsection{Admissibility reduction for constant scalar
edge weights}\label{sec:constant-scalar-reduction}

For models whose base-graph edge weights
$(c_k, \lambda_k)$ are constant scalars, the remaining
algebraic admissibility condition, $\mathcal{S}_k$-compatibility, reduces to an easily verifiable
$\mathcal{S}_k$-invariance property;
regularity remains model-specific.  The directed-path propagator
is \emph{$\mathcal{S}_k$-invariant} if for all $x, y, u \in \mathcal{V}^m$,
\begin{equation}\label{eq:Sk-invariance}
B_{\mathcal{S}_k u}
(\mathcal{S}_k x, \mathcal{S}_k y)
= B_u(x, y).
\end{equation}

\begin{proposition}\label{prop:constant-scalar}
Suppose that the base-graph edge weights
$(c_k, \lambda_k)$ are constant scalars, that the
block-diagonal lifts to the product graph are
$C_k = c_k I_m$ and $D_k = \lambda_k I_m$, and that
the directed-path propagator satisfies
$\mathcal{T}$-covariance, $\mathcal{T}$-splitting, and
$\mathcal{S}_k$-invariance~\eqref{eq:Sk-invariance}.
Then:
\begin{enumerate}[label=\textup{(\roman*)},
leftmargin=*]
\item $\mathcal{S}_k$-compatibility holds.

\item $P(\mathcal{S}_k u) = P(u)$, and the corrected
edge weights satisfy
$\Lambda_k(u) = D_k(u) = \lambda_k I_m$.
\end{enumerate}
In particular, the diamond data on $\mathcal{G}^m$
is $(C_k, \Lambda_k) = (c_k I_m, \lambda_k I_m)$.
\end{proposition}

\begin{proof}
\noindent\textit{Part \textup{(i)}.}
Substituting $C_k = c_k I_m$,
$D_k = \lambda_k I_m$, and
$\mathcal{S}_k$-invariance~\eqref{eq:Sk-invariance}
into the $\mathcal{S}_k$-compatibility condition of
Definition~\ref{def:admissible-prop}, the condition
reduces to
\begin{align*}
c_k B_u(x, y)
- \lambda_k B_u(\mathcal{T}x, y)
- c_k B_u(\mathcal{T}x, \mathcal{T}y)
+ \lambda_k B_u(\mathcal{T}x, y)
=
c_k B_u(\mathcal{T}x, \mathcal{T}u)
B_u(u, y).
\end{align*}
The $\lambda_k$-terms cancel, and since $c_k$ is a
common factor on both sides, the identity to verify is
\begin{equation}\label{eq:Sk-compat-reduced}
B_u(x, y) - B_u(\mathcal{T}x, \mathcal{T}y)
= B_u(\mathcal{T}x, \mathcal{T}u) B_u(u, y).
\end{equation}
To recognize~\eqref{eq:Sk-compat-reduced} as an
instance of $\mathcal{T}$-splitting, we rewrite the
two $B_u$ terms using $\mathcal{T}$-covariance:
$B_u(x, y) = B_{\mathcal{T}u}(\mathcal{T}x,
\mathcal{T}y)$ and
$B_u(u, y) = B_{\mathcal{T}u}(\mathcal{T}u,
\mathcal{T}y)$.  Under these substitutions,
\eqref{eq:Sk-compat-reduced} becomes
$\mathcal{T}$-splitting
(Definition~\ref{def:admissible-prop}\textup{(ii)},
Proposition~\ref{prop:directed-path-T}) evaluated at
$(\mathcal{T}x, \mathcal{T}y)$.

\medskip
\noindent\textit{Part \textup{(ii)}.}
Setting $x = y = u$
in~\eqref{eq:Sk-invariance} gives
$B_{\mathcal{S}_k u}(\mathcal{S}_k u,
\mathcal{S}_k u) = B_u(u, u)$, so
$P(\mathcal{S}_k u) = P(u)$.  Since
$D_k = \lambda_k I_m$ commutes with $P(u)$,
\[
\Lambda_k(u)
= P(u) D_k(u) P(\mathcal{S}_k u)^{-1}
= P(u) \lambda_k I_m P(u)^{-1}
= \lambda_k I_m.
\]
\end{proof}

\begin{remark}\label{rem:constant-scalar-scope}
The hypotheses of
Proposition~\ref{prop:constant-scalar} are satisfied
by the majority of models in this chapter, including
all homogeneous Bernoulli and geometric models and
their variants.  Regularity (condition~(iv) of
Definition~\ref{def:admissible-prop}) is verified
separately for each model.  For models whose
propagator is not $\mathcal{S}_k$-invariant, such as
the inhomogeneous Bernoulli model (\S\ref{sec:IBJ}),
$\mathcal{S}_k$-compatibility is verified directly.
\end{remark}
\section{Right Bernoulli Jumps}\label{sec:RBJ}
 
\paragraph{\textbf{System description.}}
We consider the discrete-time totally asymmetric simple
exclusion process with right Bernoulli
jumps, blocking interaction, and sequential update.  The system consists of $N$ particles on
$\mathbb{Z}$.  The configuration at time
$t \in \mathbb{Z}_{\ge 0}$ is denoted by
\[
Y(t) = (Y_1(t) > Y_2(t) > \cdots > Y_N(t)),
\]
where $Y_k(t)$ is the position of the $k$-th particle from
the right.
 
At each discrete time step $t \ge 1$, the particles are
updated sequentially in the order $k = 1, 2, \dots, N$, so the rightmost particle is updated first.  During its update, particle $k$
attempts to jump one step to the right with probability
$p \in (0,1)$, and stays put with probability $q = 1 - p$.
The jump is blocked if the destination site is occupied by
particle $k - 1$, whose position $Y_{k-1}(t)$ has already
been updated to its time-$t$ value by the sequential
ordering.  The evolution is given by
\[
Y_k(t) = \min\{Y_{k-1}(t) - 1, Y_k(t{-}1) + \xi(t,k)\},
\qquad Y_0(t) = \infty,
\]
where $\{\xi(t,k)\}_{t \ge 1, 1 \le k \le N}$ are
independent Bernoulli random variables with
$\mathbb{P}(\xi(t,k) = 1) = p$, and the convention
$Y_0(t) = \infty$ means that the rightmost particle
faces no blocking constraint.

\paragraph{\textbf{Fredholm determinant formula.}}
Fix a time $t \ge 0$ and an initial configuration
$Y(0) = y = (y_1 > y_2 > \cdots > y_N)$.  Fix
$m \in \{1, \dots, N{-}1\}$, particle labels
$\mathbf{n} = (n_1, n_2, \dots, n_m)$ with
$1 \le n_1 < n_2 < \cdots < n_m < N$, and spatial
thresholds $\mathbf{a} = (a_1, \dots, a_m) \in \mathbb{Z}^m$.
The multipoint joint cumulative distribution of the particle
positions is given by a Fredholm determinant on
$\ell^2(\{n_1, \dots, n_m\} \times \mathbb{Z})$~\cite{MatetskiRemenik2023}:
\[
\mathbb{P}_y\Bigl(\bigcap_{i=1}^m
\{Y_{n_i}(t) > a_i\}\Bigr)
= \det(I - \bar{\chi}_{\mathbf{a}} K_t
\bar{\chi}_{\mathbf{a}})_{\ell^2(\{n_1, \dots, n_m\}
\times \mathbb{Z})},
\]
where $\bar{\chi}_{\mathbf{a}}(n_i, x) =
\mathbf{1}_{x \le a_i}$.  The correlation kernel $K_t$ has
block entries
\[
K_t(n_i, x_i; n_j, x_j)
= -Q^{n_j - n_i}(x_i, x_j)\mathbf{1}_{n_i < n_j}
+ \bigl(\mathcal{S}_{-t,-n_i}^*
\overline{\mathcal{S}}_{-t,n_j}^{\operatorname{epi}(y)}
\bigr)(x_i, x_j),
\]
where $A^*$ denotes the transpose kernel,
$A^*(x, y) = A(y, x)$.

Fix an auxiliary parameter $\theta \in (0, 1)$, which
enters the Fredholm determinant representation but does
not appear in the final multipoint equation
\eqref{eq:H5-multipoint}, and let
$\alpha = (1{-}\theta)\theta^{-1}$.  The operators
defining the kernel are as follows.

\medskip
\noindent\textit{The transition matrix $Q$.}
The matrix $Q$ is the transition matrix of a random walk on
$\mathbb{Z}$ taking $\mathrm{Geom}[1{-}\theta]$ steps
strictly to the left:
\[
Q(z_1, z_2) = (1{-}\theta)\theta^{z_1 - z_2 - 1}
\mathbf{1}_{z_1 > z_2},
\]
with its $n$-th power admitting the integral
representation
\[
Q^n(z_1, z_2) = \frac{\alpha^n}{2\pi\mathrm{i}}
\oint_{\gamma_\rho} \diff w\
\frac{\theta^{z_1 - z_2}}{w^{z_1 - z_2 - n + 1}}
\Bigl(\frac{1}{1 - w}\Bigr)^n.
\]

\medskip
\noindent\textit{The scattering operators
$\mathcal{S}_{-t,-n}$ and
$\overline{\mathcal{S}}_{-t,n}$.}\footnote{The kernel
operators $\mathcal{S}_{-t,-n}$ and
$\overline{\mathcal{S}}_{-t,n}$ are unrelated to the
lattice shifts $\mathcal{S}_k$ of the diamond framework
(\S\ref{sec:linear-problem}); the notation follows
Matetski--Remenik~\cite{MatetskiRemenik2023}.}
These operators are built from the probability generating
function $\varphi(w) = q + pw$ of the Bernoulli jump
distribution.  They are defined by the contour integrals
\[
\mathcal{S}_{-t,-n}(z_1, z_2)
= \frac{\alpha^{-n+1}}{2\pi\mathrm{i}}
\oint_{\gamma_\rho} \frac{\diff w}{w}
\frac{\theta^{z_2 - z_1}}{w^{z_2 - z_1 + n}}
(1{-}w)^n (q{+}pw)^t,
\]
\[
\overline{\mathcal{S}}_{-t,n}(z_1, z_2)
= \frac{\alpha^{n-1}}{2\pi\mathrm{i}}
\oint_{\gamma_\delta} \frac{\diff w}{w}
\frac{(1{-}w)^{z_2 - z_1 + n - 1}}
{\theta^{z_2 - z_1} w^{n-1}}
(q{+}p(1{-}w))^{-t},
\]
where $\gamma_\rho$ is a positively oriented circle of
radius $\rho \in (0, 1)$ around the origin, and
\(\gamma_\delta\) is a sufficiently small positively
oriented circle around the origin, chosen so that the only
singularity enclosed is the pole at \(w=0\).

\medskip
\noindent\textit{The epigraph operator
$\overline{\mathcal{S}}_{-t,n}^{\operatorname{epi}(y)}$.}
This is the hitting probability operator defined in terms
of the initial data:
\[
\overline{\mathcal{S}}_{-t,n}^{\operatorname{epi}(y)}
(z_1, z_2)
= \mathbb{E}_{W_0 = z_1}\bigl[
\overline{\mathcal{S}}_{-t,n-\tau}(W_\tau, z_2)
\mathbf{1}_{\tau < n}\bigr],
\]
where $(W_\ell)_{\ell \ge 0}$ is the random walk with
transition matrix $Q$, and
$\tau = \min\{\ell \in \{0, \dots, N{-}1\} :
W_\ell > y_{\ell+1}\}$
is the hitting time of the strict epigraph of the initial
data by the random walk, with the convention
$\tau = \infty$ if the set is empty.

\paragraph{\textbf{References.}}
The sequential-update Bernoulli TASEP above was treated by
Brankov, Priezzhev, and Shelest~\cite{BPS2004}.  The same
model belongs to the family of four basic discrete-time
TASEP variants with blocking or pushing interaction and
Bernoulli or geometric jumps whose determinantal transition
kernels were obtained by
Dieker--Warren~\cite{DiekerWarren2008}.  The Fredholm determinant
formula used below is the Matetski--Remenik formula
(Theorem~1.2, equations~(1.3)--(1.10)
of~\cite{MatetskiRemenik2023}) specialized to the right
Bernoulli case, equivalently to $\kappa = 0$ and
$\varphi(w) = q + pw$ in the notation of Matetski and
Remenik.

\subsection{Kernel reformulation}\label{sec:RBJ-kernel-reform}
 
We rewrite the extended-kernel Fredholm determinant in the
form of the product graph construction of \S\ref{sec:product-graph}.  The
functions $\psi, \phi$ below provide the seed data, while the
lower-triangular operator $B_{\mathbf{a}, \mathbf{n}}$ will be verified
to satisfy the admissibility conditions of Definition~\ref{def:admissible-prop}.
 
\paragraph{\textbf{Seed data.}}
Define
\[
\psi_{t,r,n}(x) \defeq \mathcal{S}_{-t,-n}^*(r,x)
= \mathcal{S}_{-t,-n}(x,r),
\qquad
\phi_{t,r,n}(x) \defeq
\overline{\mathcal{S}}_{-t,n}^{\operatorname{epi}(y)}(x,r),
\]
where we suppress the dependence of $\phi$ on the initial
condition $y$.  The subscript notation $\psi_{t,r,n}$ places
the parameter coordinates $(t, r, n)$ as subscripts,
reserving the function argument for the Hilbert-space
variable $x \in \Z$; the same convention applies throughout
this chapter.

\paragraph{\textbf{Propagator.}}
The propagator is an instance of the directed-path
propagator (Definition~\ref{def:directed-path-propagator})
with transition kernels
$Q_{\ell,\ell'} = Q^{n_{\ell'} - n_\ell}$ for
$\ell < \ell'$.  Explicitly, for $i > j$,
\begin{equation}\label{eq:RBJ-propagator}
[B_{\mathbf{a}, \mathbf{n}}]_{ij}(r, r')
\defeq \sum_{k=1}^{i-j}
  \sum_{j = \ell_0 < \ell_1 < \cdots < \ell_k = i}
  \sum_{\substack{\xi_1 \le a_{\ell_1} \\ \vdots \\
  \xi_{k-1} \le a_{\ell_{k-1}}}}
  \prod_{s=0}^{k-1}
  \Bigl(-Q^{n_{\ell_{s+1}} - n_{\ell_s}}
  (\xi_s, \xi_{s+1})\Bigr),
\end{equation}
with $\xi_0 = r'$, $\xi_k = r$, and $[B_{\mathbf{a}, \mathbf{n}}]_{ij} = 0$ for $i \le j$.
The internal vertices $\xi_1, \dotsc, \xi_{k-1}$ are
constrained by the cutoffs, but the endpoints
$r, r'$ are arbitrary, so
$[B_{\mathbf{a}, \mathbf{n}}]_{ij}(r, r')$ is defined for all
$r, r' \in \mathbb{Z}$.  The propagator depends on $\mathbf{a}$
and $\mathbf{n}$ but not on the time coordinate $t$, which is
accordingly suppressed from the notation.

On the cutoff region $r \le a_i$, $r' \le a_j$, the
entries of $I + B_{\mathbf{a}, \mathbf{n}}$ coincide with those of
$((I + \bar{\chi}_{\mathbf{a}} L
\bar{\chi}_{\mathbf{a}})^{-1})^{\!\top}$, where $L$
is the strictly upper-triangular block operator with
entries $L_{ij}(x, x') = Q^{n_j - n_i}(x, x')
\mathbf{1}_{i < j}$ (Remark~\ref{rem:Neumann-series}).
 
\begin{lemma}\label{lem:RBJ-kernel}
We have
\[
\mathbb{P}_y\Bigl(\bigcap_{i=1}^m
\{Y_{n_i}(t) > a_i\}\Bigr)
= \det(I - K_{t, \mathbf{a}, \mathbf{n}})_{\ell^2(\mathbb{Z})},
\]
where
\[
K_{t, \mathbf{a}, \mathbf{n}}(x, x')
= \sum_{1 \le j \le i \le m}
  \sum_{r \le a_i} \sum_{r' \le a_j}
   \psi_{t, r, n_i}(x)
  \bigl[\delta_{ij}\delta_{r,r'}
  + [B_{\mathbf{a}, \mathbf{n}}]_{ij}(r, r')\bigr]
  \phi_{t, r', n_j}(x').
\]
This identifies $K_{t, \mathbf{a}, \mathbf{n}}$ with the
kernel~\eqref{eq:K-def} on the product graph
$\mathcal{G}^m$, with $\psi, \phi$ as the seed data and
$B_{\mathbf{a}, \mathbf{n}}$ as the propagator.
\end{lemma}

\begin{proof}
We conjugate the extended kernel with block-dependent
exponential weights to obtain a trace-class realization,
then reduce to a single-space determinant via Sylvester's
identity.

Choose \(\beta\in(\theta,1)\), and choose conjugation weights
\[
1<\omega_m<\omega_{m-1}<\cdots<\omega_1<\beta^{-1}.
\]
For a function
\(f\colon\{n_1,\dotsc,n_m\}\times\mathbb Z\to\mathbb R\),
let
\[
(\mathscr Cf)_i(x)\defeq\omega_i^x f_{n_i}(x),
\]
so that for a block kernel \(M\) on
\(\{n_1,\dotsc,n_m\}\times\mathbb Z\),
\[
(\mathscr CM\mathscr C^{-1})_{ij}(x,x')
=
\omega_i^x M_{ij}(x,x')\omega_j^{-x'}.
\]

Next, define
\[
U(n_i,r;x)
\defeq
\mathbf{1}_{r\le a_i}\psi_{t,r,n_i}(x),
\qquad
V(x;n_j,r')
\defeq
\phi_{t,r',n_j}(x)\mathbf{1}_{r'\le a_j}.
\]
Then
\[
\bar{\chi}_{\mathbf a}
\mathcal{S}_{-t,-n_i}^{*}
\overline{\mathcal{S}}_{-t,n_j}^{\operatorname{epi}(y)}
\bar{\chi}_{\mathbf a}
=
UV,
\]
and so the Matetski--Remenik kernel can be written as
\[
I-\bar{\chi}_{\mathbf a}K_t\bar{\chi}_{\mathbf a}
=
I+\bar{\chi}_{\mathbf a}L\bar{\chi}_{\mathbf a}-UV.
\]
Since \(\mathscr C\) multiplies each row of a finite
minor by \(\omega_i^x\) and each column by
\(\omega_j^{-x'}\), the determinant of each minor is
unchanged, and hence the Fredholm series of the
Matetski--Remenik determinant agrees term by term with
that of the conjugated kernel.  The estimates below show
that the conjugated kernel is trace class, so this gives
a trace-class realization of the same extended
determinant.

Set
\[
\widetilde N
\defeq
\mathscr C\bar{\chi}_{\mathbf a}L
\bar{\chi}_{\mathbf a}\mathscr C^{-1},
\qquad
\widetilde U\defeq\mathscr C U,
\qquad
\widetilde V\defeq V\mathscr C^{-1}.
\]
We verify the Schatten-class membership of each piece in turn: \(\widetilde{N}\) is trace class (\(\widetilde{N} \in \mathfrak{S}_1\)), while \(\widetilde{U}\) and \(\widetilde{V}\) are Hilbert--Schmidt (\(\in \mathfrak{S}_2\)).

For \(i<j\), let \(d\defeq n_j-n_i\).  The \((i,j)\)-block
of \(\widetilde N\) has kernel
\[
\omega_i^x Q^d(x,x')\omega_j^{-x'}
\mathbf{1}_{x\le a_i}\mathbf{1}_{x'\le a_j}.
\]
Writing \(s=x-x'\), this becomes
\[
\omega_i^s
\left(\frac{\omega_i}{\omega_j}\right)^{x'}
Q^d(s)
\mathbf{1}_{x'\le a_j}
\mathbf{1}_{x'+s\le a_i}.
\]
Expanding the \(d\)-fold convolution and counting
compositions of \(s\) into \(d\) positive parts gives
\[
Q^d(s)
=
\binom{s-1}{d-1}(1-\theta)^d\theta^{s-d}
\mathbf{1}_{s\ge d}.
\]
Hence, dropping the cutoff
\(\mathbf{1}_{x'+s\le a_i}\) since all summands are
nonnegative,
\[
\sum_{x\le a_i}
\sum_{x'\le a_j}
\left|
\omega_i^x Q^d(x,x')\omega_j^{-x'}
\right|
\le
\sum_{s\ge d}Q^d(s)\omega_i^s
\sum_{x'\le a_j}
\left(\frac{\omega_i}{\omega_j}\right)^{x'}.
\]
The second factor on the right converges because \(i<j\) implies
\(\omega_i/\omega_j>1\).  The first factor converges
because, by the formula above,
\(
\sum_{s\ge d} Q^d(s)\omega_i^s
=
(1-\theta)^d\theta^{-d}
\sum_{s\ge d}
\binom{s-1}{d-1}(\theta\omega_i)^s;
\)
the binomial coefficient grows polynomially in \(s\),
while the geometric factor \((\theta\omega_i)^s\) decays
exponentially since
\(\theta\omega_i<\theta\beta^{-1}<1\).  Thus every
nonzero block of \(\widetilde N\) has absolutely summable
kernel, hence is trace class.  Since there are only
finitely many blocks,
\[
\widetilde N\in\mathfrak S_1.
\]
Moreover \(\widetilde N\) is strictly upper triangular in
the block indices, hence nilpotent.

Next, the contour formula for \(\psi\) gives
\[
\psi_{t,r,n}(x)=0
\qquad\text{unless}\qquad
r-t\le x\le r+n.
\]
Indeed, the integral extracts the coefficient of
\(w^{r-x+n}\) in the polynomial \((1-w)^n(q+pw)^t\),
which has degree \(n+t\); this coefficient vanishes
unless \(0\le r-x+n\le n+t\).  Moreover, the integrand
depends on \(r\) and \(x\) only through \(r-x\), so
\(\psi_{t,r,n}(x)=g(r-x)\) for a fixed function \(g\)
supported on \(\{-n,\dotsc,t\}\).  In particular,
\(\lVert\psi_{t,r,n}\rVert_{\ell^2(\mathbb Z)}
=\lVert g\rVert_{\ell^2}\) is
independent of \(r\).

It follows that
\[
\lVert\widetilde U\rVert_2^2
=
\sum_{i=1}^m
\sum_{r\le a_i}
\omega_i^{2r}\lVert\psi_{t,r,n_i}\rVert_{\ell^2}^2
\le
C
\sum_{i=1}^m
\sum_{r\le a_i}
\omega_i^{2r}
<\infty,
\]
because \(\omega_i>1\).  Thus
\[
\widetilde U\in\mathfrak S_2.
\]

It remains to estimate \(\widetilde V\).  Let \(a_{\max}=\max_{1\le j\le m}a_j\).  We shall only
need the following estimate for \(r\le a_{\max}\) and
\(z>y_N\).  Choose \(\delta>0\) so small that, on
\(\gamma_\delta\),
\[
\eta:=\sup_{w\in\gamma_\delta}
\left|\frac{\theta}{1-w}\right|<\beta .
\]
The contour formula gives
\[
\overline{\mathcal S}_{-t,k}(z,r)
=
\frac{\alpha^{k-1}}{2\pi i}
\oint_{\gamma_\delta}
\frac{\diff w}{w^k}
(1-w)^{k-1}
(1-pw)^{-t}
\left(\frac{\theta}{1-w}\right)^{z-r}.
\]
If \(z\ge r\), then \(z-r\ge0\), and hence
\[
\left|
\left(\frac{\theta}{1-w}\right)^{z-r}
\right|
\le \eta^{z-r}\le \beta^{z-r}.
\]
The remaining factors are uniformly bounded on the fixed
contour, uniformly over \(1\le k\le N\).  If \(z<r\),
the restrictions \(z>y_N\) and \(r\le a_{\max}\)
confine \((z,r)\) to a finite set.  Enlarging the
constant to absorb these finitely many values, we obtain
\[
\left|\overline{\mathcal S}_{-t,k}(z,r)\right|
\le C\beta^{z-r}
\]
for all \(1\le k\le N\), \(z>y_N\), and
\(r\le a_{\max}\).

We now bound \(\phi\).  The hitting-time representation
gives
\[
|\phi_{t,r,n}(x)|
\le
\sum_{\ell=0}^{n-1}
\sum_{z>y_{\ell+1}}
Q^\ell(x,z)
\left|\overline{\mathcal S}_{-t,n-\ell}(z,r)\right|,
\]
where we have replaced
\(\mathbb P(\tau=\ell, W_\ell=z)\) by the larger
quantity \(Q^\ell(x,z)=\mathbb P(W_\ell=z)\).
Applying the contour bound and extending the sum from
\(z>y_{\ell+1}\) to all \(z\) gives
\[
|\phi_{t,r,n}(x)|
\le
C\beta^{-r}
\sum_{\ell=0}^{n-1}
\sum_z Q^\ell(x,z)\beta^z.
\]
Since \(Q\) is the left-jump geometric walk,
\(W_\ell=x-S_\ell\) where \(S_\ell\) is the cumulative
displacement.  Hence
\(\sum_z Q^\ell(x,z)\beta^z
=\beta^x\mathbb E[\beta^{-S_\ell}]\),
and the exponential moment is finite because
\(\beta>\theta\).  Summing over \(\ell\) yields
\(|\phi_{t,r,n}(x)| \le C\beta^{x-r}\).
Since the walk \(W_\ell\) is strictly decreasing,
\(\tau = \infty\) whenever \(x \le y_N\) (the smallest initial position, hence below every \(y_{\ell+1}\)), so
\(\phi_{t,r,n}(x) = 0\) in that range.  Combining:
\begin{equation}\label{eq:RBJ-phi-bound}
|\phi_{t,r,n}(x)|
\le
C\beta^{x-r}\mathbf{1}_{x>y_N}.
\end{equation}
Since \(\beta<1\),
\[
\lVert\phi_{t,r,n}\rVert_{\ell^2}^2
\le
C\beta^{-2r}
\sum_{x>y_N}\beta^{2x}
\le
C\beta^{-2r}.
\]
Consequently
\[
\lVert\widetilde V\rVert_2^2
=
\sum_{j=1}^m
\sum_{r'\le a_j}
\omega_j^{-2r'}
\lVert\phi_{t,r',n_j}\rVert_{\ell^2}^2
\le
C
\sum_{j=1}^m
\sum_{r'\le a_j}
(\omega_j\beta)^{-2r'}
<\infty,
\]
because \(\omega_j\beta<1\).  Hence
\(\widetilde V\in\mathfrak S_2\).

\medskip

With the Schatten-class estimates in hand, we turn to the
determinantal reduction.  The conjugated extended kernel
\(-\widetilde N+\widetilde U\widetilde V\) is
trace class, since \(\widetilde N\in\mathfrak S_1\) and
\(\widetilde U\widetilde V\in\mathfrak S_1\) as a product
of two Hilbert--Schmidt operators.  Factoring out the
nilpotent triangular part:
\[
I+\widetilde N-\widetilde U\widetilde V
=
(I+\widetilde N)
\bigl(
I-(I+\widetilde N)^{-1}\widetilde U\widetilde V
\bigr).
\]
Since \(\widetilde N\) is nilpotent,
\(\det(I+\widetilde N)=1\), and
\((I+\widetilde N)^{-1}
= \mathscr C
(I+\bar{\chi}_{\mathbf a}L\bar{\chi}_{\mathbf a})^{-1}
\mathscr C^{-1}\)
is given by the finite Neumann series.  The product
\((I+\widetilde N)^{-1}\widetilde U\in\mathfrak S_2\)
(since \((I+\widetilde N)^{-1}\) is bounded and
\(\mathfrak S_2\) is an operator ideal), so
Sylvester's identity applies.  Moreover, the conjugation
weights cancel in the product
\(\widetilde V(I+\widetilde N)^{-1}\widetilde U
= V(I+\bar{\chi}_{\mathbf a}L
\bar{\chi}_{\mathbf a})^{-1}U\).
Therefore, setting
\(M\defeq I+\bar{\chi}_{\mathbf a}L
\bar{\chi}_{\mathbf a}\),
\[
\det
\bigl(I-\bar{\chi}_{\mathbf a}K_t
\bar{\chi}_{\mathbf a}\bigr)
=
\det_{\ell^2(\mathbb Z)}
\bigl(
I- VM^{-1}U
\bigr).
\]
Since the Fredholm determinant of a trace-class operator
is invariant under transposition, we may replace the
kernel of \(VM^{-1}U\) by its transpose.  Writing
\([M^{-\top}]_{ij}(r,r') = [M^{-1}]_{ji}(r',r)\),
the transposed kernel is
\[
\sum_{i,j=1}^m
\sum_{r\le a_i}
\sum_{r'\le a_j}
\psi_{t,r,n_i}(x)
[M^{-\top}]_{ij}(r,r')
\phi_{t,r',n_j}(x').
\]

It remains to identify the inverse-transpose entries with
the propagator.  For \(i=j\), the inverse-transpose
contributes the identity \(\delta_{r,r'}\) and
\([B_{\mathbf{a}, \mathbf{n}}]_{ii}=0\), so the combined factor is
\(\delta_{ij}\delta_{r,r'}\).  For \(i<j\), both
\([M^{-\top}]_{ij}\) and \([B_{\mathbf{a}, \mathbf{n}}]_{ij}\) vanish,
the former because the inverse-transpose is block lower
triangular and the latter by the definition
of~$B_{\mathbf{a}, \mathbf{n}}$.  For \(i>j\), expanding the finite Neumann
series in the block indices and restricting to the
cutoff region \(r\le a_i\), \(r'\le a_j\) (which is
the only region over which the kernel formula sums)
gives
\[
[M^{-\top}]_{ij}(r,r')
=
[B_{\mathbf{a}, \mathbf{n}}]_{ij}(r,r').
\]
Indeed, a nonzero term in the Neumann series has the form
\(j=\ell_0<\ell_1<\cdots<\ell_k=i\),
\(\xi_0=r'\), \(\xi_k=r\),
with internal constraints
\(\xi_s\le a_{\ell_s}\) for \(s=1,\dotsc,k-1\), and
weight
\[
(-1)^k
\prod_{s=0}^{k-1}
Q^{n_{\ell_{s+1}} - n_{\ell_s}}
(\xi_s, \xi_{s+1});
\]
the sign \((-1)^k\) is absorbed by the \(k\) factors of
\((-Q)\) in the definition of
\([B_{\mathbf{a}, \mathbf{n}}]_{ij}(r,r')\).  The support condition
\(Q^d(x,x')=0\) unless \(x-x'\ge d\) shows that the
internal sums are finite for fixed endpoints; indeed,
along a nonzero chain,
\(\xi_s-\xi_{s+1}\ge n_{\ell_{s+1}}-n_{\ell_s}\ge1\),
so the spatial chain is strictly decreasing.

\medskip

Thus the transposed kernel is
\[
K_{t,\mathbf a,\mathbf n}(x,x')
=
\sum_{1\le j\le i\le m}
\sum_{r\le a_i}
\sum_{r'\le a_j}
\psi_{t,r,n_i}(x)
\left(
\delta_{ij}\delta_{r,r'}
+
[B_{\mathbf{a}, \mathbf{n}}]_{ij}(r,r')
\right)
\phi_{t,r',n_j}(x').
\]
This is the claimed single-space kernel.
\end{proof}
 For the remainder of this section, we work with the
kernel $K_{t, \mathbf{a}, \mathbf{n}}$ on
$\ell^2(\mathbb{Z})$.
 
\subsection{Base graph and seed data}
\label{sec:RBJ-base-graph}
 
We now identify the lattice structure that governs how
the seed data evolve under shifts of the parameters
$(t, a, n)$, and verify that $\psi, \phi$ satisfy the
linear problems of the diamond framework.
 
Let $\mathcal{V} \subseteq \mathbb{Z}^3$ be the lattice
generated by the shifts
\[
T = e^{\pa_a}, \qquad
S_1 = e^{\pa_n}, \qquad
S_2 = e^{\pa_t + \pa_a},
\]
so that for $u = (t, a, n) \in \mathcal{V}$,
\[
Tu = (t, a{+}1, n), \qquad
S_1 u = (t, a, n{+}1), \qquad
S_2 u = (t{+}1, a{+}1, n).
\]
Take $E = \mathbb{F} = \mathbb{R}$ and define the constant
scalar edge weights
\begin{equation}\label{eq:RBJ-edge-weights}
(c_1, c_2) = \bigl((1{-}\theta)^{-1}, q\bigr),
\qquad
(\lambda_1, \lambda_2)
= \bigl((1{-}\theta)^{-1}\theta, -p\theta\bigr).
\end{equation}
Since $(c_k, \lambda_k)$ are constant scalars, the diamond
equations
\eqref{eq:diamond-C-ij}--\eqref{eq:diamond-mixed-ij} are
satisfied trivially: all three reduce to commutativity in
$\mathbb{R}$.
 
\begin{lemma}\label{lem:RBJ-seed-linear}
Let $H = \ell^2(\mathbb{Z})$.  The seed functions
$\psi_{t,r,n}$ and $\phi_{t,r,n}$ defined in
\S\ref{sec:RBJ-kernel-reform} satisfy the linear problem
\eqref{eq:Linear-Psi-3} and the adjoint linear problem
\eqref{eq:Linear-Phi-Prob} on $\mathcal{V}$ with edge
weights $(c_k, \lambda_k)$.  Explicitly, for each fixed
$r \in \mathbb{Z}$:
\begin{align}
\psi_{t,r,n+1}
&= (1{-}\theta)^{-1}\psi_{t,r+1,n}
   - (1{-}\theta)^{-1}\theta\psi_{t,r,n},
\label{eq:RBJ-psi-S1} \\
\psi_{t+1,r+1,n}
&= q\psi_{t,r+1,n}
   + p\theta\psi_{t,r,n},
\label{eq:RBJ-psi-S2}
\end{align}
and
\begin{align}
\phi_{t,r+1,n}
&= (1{-}\theta)^{-1}\phi_{t,r,n+1}
   - (1{-}\theta)^{-1}\theta\phi_{t,r+1,n+1},
\label{eq:RBJ-phi-S1} \\
\phi_{t,r+1,n}
&= q\phi_{t+1,r+1,n}
   + p\theta\phi_{t+1,r+2,n}.
\label{eq:RBJ-phi-S2}
\end{align}
\end{lemma}
 
\begin{proof}
We verify each identity in turn.

\medskip
\noindent\textit{Verification of \eqref{eq:RBJ-psi-S1}.}
From the contour representation,
\[
\psi_{t,r,n}(x) = \frac{\alpha^{-n+1}}{2\pi\mathrm{i}}
  \oint_{\gamma_\rho} \frac{\diff w}{w}
  \frac{\theta^{r-x}}{w^{r-x+n}}
  (1{-}w)^n (q{+}pw)^t.
\]
The shift $r \to r+1$ contributes a factor $\theta/w$ to
the integrand of $\psi_{t,r+1,n}$ relative to
$\psi_{t,r,n}$.  Factoring out the common integrand, the
right-hand side of \eqref{eq:RBJ-psi-S1} becomes
\begin{align*}
&(1{-}\theta)^{-1}\psi_{t,r+1,n}(x)
- (1{-}\theta)^{-1}\theta\psi_{t,r,n}(x) \\
&\qquad = \frac{\alpha^{-n+1}}{2\pi\mathrm{i}}
  \oint_{\gamma_\rho} \frac{\diff w}{w}
  \frac{\theta^{r-x}}{w^{r-x+n}}
  (1{-}w)^n (q{+}pw)^t
  \cdot\frac{1}{1{-}\theta}
  \Bigl(\frac{\theta}{w} - \theta\Bigr).
\end{align*}
The combined factor simplifies as
\[
\frac{1}{1{-}\theta}
\Bigl(\frac{\theta}{w} - \theta\Bigr)
= \frac{\theta(1{-}w)}{(1{-}\theta) w}
= \frac{1{-}w}{\alpha w},
\]
since $\theta/(1{-}\theta) = \alpha^{-1}$.  Thus
the prefactor becomes $\alpha^{-n}$, while the extra
factor $(1{-}w)/w$ changes the integrand to that of
$\psi_{t,r,n+1}$:
\[
= \frac{\alpha^{-n}}{2\pi\mathrm{i}}
  \oint_{\gamma_\rho} \frac{\diff w}{w}
  \frac{\theta^{r-x}}{w^{r-x+n+1}}
  (1{-}w)^{n+1} (q{+}pw)^t
= \psi_{t,r,n+1}(x).
\]

\medskip
\noindent\textit{Verification of \eqref{eq:RBJ-psi-S2}.}
The right-hand side is
\begin{align*}
q\psi_{t,r+1,n}(x) + p\theta\psi_{t,r,n}(x)
= \frac{\alpha^{-n+1}}{2\pi\mathrm{i}}
  \oint_{\gamma_\rho} \frac{\diff w}{w}
  \frac{\theta^{r-x}}{w^{r-x+n}}
  (1{-}w)^n (q{+}pw)^t
  \Bigl(\frac{q\theta}{w} + p\theta\Bigr).
\end{align*}
The factor $q\theta/w + p\theta = \theta(q + pw)/w$
contributes one extra $(q{+}pw)$ and one extra
$\theta/w$ to the integrand, giving
\[
= \frac{\alpha^{-n+1}}{2\pi\mathrm{i}}
  \oint_{\gamma_\rho} \frac{\diff w}{w}
  \frac{\theta^{(r+1)-x}}{w^{(r+1)-x+n}}
  (1{-}w)^n (q{+}pw)^{t+1}
= \psi_{t+1,r+1,n}(x).
\]

\medskip
\noindent\textit{Free recurrences for
$\overline{\mathcal{S}}$.}
We record three identities for the free operator that will
be applied inside the hitting-time representation of
$\phi$.  For $k \ge 1$,
\begin{equation}\label{eq:Sbar-recurrence-n}
(1{-}\theta)^{-1}\bigl[
\overline{\mathcal{S}}_{-t,k}(z, r)
- \theta\overline{\mathcal{S}}_{-t,k}(z, r{+}1)
\bigr]
= \overline{\mathcal{S}}_{-t,k-1}(z, r{+}1).
\end{equation}
Indeed, replacing $r$ by $r{+}1$ in the contour formula
multiplies the integrand by $(1{-}w)/\theta$.  Factoring
out the common integrand, the left-hand side becomes
\[
\frac{\alpha^{k-1}}{(1{-}\theta)2\pi\mathrm{i}}
\oint_{\gamma_\delta} \frac{\diff w}{w^k}
(1{-}w)^{k-1}(1{-}pw)^{-t}
\Bigl(\frac{\theta}{1{-}w}\Bigr)^{z-r}
\bigl(1 - (1{-}w)\bigr).
\]
The factor $1 - (1{-}w) = w$ reduces
$\diff w/w^k$ to $\diff w/w^{k-1}$, and dividing $\alpha^{k-1}$ by
$(1{-}\theta) = \alpha\theta$ gives the prefactor
$\alpha^{k-2}/\theta$, matching
$\overline{\mathcal{S}}_{-t,k-1}(z, r{+}1)$.

We shall also use
\begin{equation}\label{eq:Sbar-recurrence-t}
q\overline{\mathcal{S}}_{-(t+1),k}(z, r{+}1)
+ p\theta\overline{\mathcal{S}}_{-(t+1),k}(z, r{+}2)
= \overline{\mathcal{S}}_{-t,k}(z, r{+}1).
\end{equation}
In the contour formula for
$\overline{\mathcal{S}}_{-(t+1),k}$, the shift from
$r{+}1$ to $r{+}2$ multiplies the integrand by
$(1{-}w)/\theta$.  The coefficient $p\theta$ on the second
term absorbs the $1/\theta$, contributing $p(1{-}w)$,
while the first term contributes~$q$.  The left-hand side
becomes
\[
\frac{\alpha^{k-1}}{2\pi\mathrm{i}}
\oint_{\gamma_\delta} \frac{\diff w}{w^k}
(1{-}w)^{k-1}(1{-}pw)^{-(t+1)}
\Bigl(\frac{\theta}{1{-}w}\Bigr)^{z-r-1}
\bigl(q + p(1{-}w)\bigr).
\]
The factor $q + p(1{-}w) = 1 - pw$ cancels one power of
$(1{-}pw)^{-(t+1)}$, recovering the integrand of
$\overline{\mathcal{S}}_{-t,k}(z, r{+}1)$.

Finally, we show that $\overline{\mathcal{S}}_{-t,0} \equiv 0$.
When $k = 0$, the measure $\diff w/w^k$ reduces to $\diff w$,
removing the pole at the origin.  The remaining integrand
has singularities only at $w = 1$ (from powers of
$1{-}w$) and $w = 1/p$ (from $(1{-}pw)^{-t}$), both of
which lie outside $\gamma_\delta$.  The integrand is
therefore analytic inside $\gamma_\delta$, and the integral
vanishes by Cauchy's theorem.

\medskip
\noindent\textit{Verification of \eqref{eq:RBJ-phi-S1}.}
The exponential bound on $\overline{\mathcal{S}}$
established in the proof of Lemma~\ref{lem:RBJ-kernel}, together with
the finite exponential moments of the geometric left-jump
walk, implies that the expectations defining $\phi$ are
absolutely convergent.  Hence the free identities above
may be applied inside the hitting-time expectation.

Let $(W_\ell)_{\ell \ge 0}$ denote the random walk with
transition matrix $Q$ and $\tau$ its hitting time of the
strict epigraph of the initial data, as in the definition
of $\phi$.  The right-hand side of \eqref{eq:RBJ-phi-S1}
equals
\begin{align*}
&(1{-}\theta)^{-1}\phi_{t,r,n+1}(x)
- (1{-}\theta)^{-1}\theta\phi_{t,r+1,n+1}(x) \\
&\quad = \mathbb{E}_{W_0 = x}\Bigl[
(1{-}\theta)^{-1}\bigl(
\overline{\mathcal{S}}_{-t,n+1-\tau}(W_\tau, r)
- \theta\overline{\mathcal{S}}_{-t,n+1-\tau}
  (W_\tau, r{+}1)
\bigr)\mathbf{1}_{\tau < n+1}\Bigr].
\end{align*}
We split the expectation according to whether
$\tau \le n{-}1$ or $\tau = n$.  On the event
$\{\tau \le n{-}1\}$, set $k = n{+}1{-}\tau$.  Then
$k \ge 2$, and \eqref{eq:Sbar-recurrence-n} gives
\begin{align*}
(1{-}\theta)^{-1}\bigl(
\overline{\mathcal{S}}_{-t,n+1-\tau}(W_\tau, r)
- \theta\overline{\mathcal{S}}_{-t,n+1-\tau}
  (W_\tau, r{+}1)
\bigr)
= \overline{\mathcal{S}}_{-t,n-\tau}(W_\tau, r{+}1).
\end{align*}
On the boundary event $\{\tau = n\}$, the same recurrence
produces 
\[\overline{\mathcal{S}}_{-t,0}(W_n, r{+}1) = 0.\]
Therefore the boundary contribution vanishes, and the
expression reduces to
\[
\mathbb{E}_{W_0 = x}\bigl[
\overline{\mathcal{S}}_{-t,n-\tau}(W_\tau, r{+}1)
\mathbf{1}_{\tau < n}\bigr]
= \phi_{t,r+1,n}(x).
\]

\medskip
\noindent\textit{Verification of \eqref{eq:RBJ-phi-S2}.}
Here the particle label $n$ is not shifted, so the
indicator $\mathbf{1}_{\tau < n}$ is unchanged.  Applying
\eqref{eq:Sbar-recurrence-t} with $k = n{-}\tau$ on
the event $\{\tau < n\}$, we obtain
\begin{align*}
&q\phi_{t+1,r+1,n}(x)
+ p\theta\phi_{t+1,r+2,n}(x) \\
&\quad = \mathbb{E}_{W_0 = x}\bigl[
\bigl(q\overline{\mathcal{S}}_{-(t+1),n-\tau}
  (W_\tau, r{+}1)
+ p\theta\overline{\mathcal{S}}_{-(t+1),n-\tau}
  (W_\tau, r{+}2)\bigr)
\mathbf{1}_{\tau < n}\bigr] \\
&\quad = \mathbb{E}_{W_0 = x}\bigl[
\overline{\mathcal{S}}_{-t,n-\tau}(W_\tau, r{+}1)
\mathbf{1}_{\tau < n}\bigr]
= \phi_{t,r+1,n}(x).
\end{align*}
\end{proof}

\subsection{Admissible propagators and corrected edge weights}
\label{sec:RBJ-propagators}
 
The product graph $\mathcal{G}^m$ has vertex set
$\mathcal{V}^m$ and is generated by the diagonal shifts
\[
\mathcal{T}u = (t, \mathbf{a} + \mathbf{1}, \mathbf{n}),
\qquad
\mathcal{S}_1 u
= (t, \mathbf{a}, \mathbf{n} + \mathbf{1}),
\qquad
\mathcal{S}_2 u
= (t{+}1, \mathbf{a} + \mathbf{1}, \mathbf{n}),
\]
where $u = (t, \mathbf{a}, \mathbf{n}) \in \mathcal{V}^m$
and $\mathbf{1} = (1, \ldots, 1) \in \mathbb{Z}^m$.  The
block-diagonal edge weights are
$C_k(u) = c_k I_m$ and $D_k(u) = \lambda_k I_m$
(the uncorrected edge weights of
Definition~\ref{def:admissible-prop}).

The propagator \eqref{eq:RBJ-propagator} is an instance of
the directed-path propagator
(Definition~\ref{def:directed-path-propagator}) with
transition kernels $Q_{\ell,\ell'} = Q^{n_{\ell'} - n_\ell}$
for $\ell < \ell'$.  The standing hypotheses of
\S\ref{sec:directed-path} are satisfied: translation
invariance of $Q$ implies translation invariance of
each power $Q^n$, and $\mathcal{T}$-invariance holds
because the transition kernels
$Q_{\ell,\ell'} = Q^{n_{\ell'} - n_\ell}$ depend only
on the particle-label differences, which are unchanged
by~$\mathcal{T}$.

By Proposition~\ref{prop:constant-scalar},
$\mathcal{S}_k$-compatibility follows from the
constant-scalar edge weights once the propagator satisfies
the following invariance property.

\begin{lemma}[$\mathcal{S}_k$-invariance]
\label{lem:RBJ-Sk-invariance}
For $k = 1, 2$ and all $x, y \in \mathcal{V}^m$,
\begin{equation}\label{eq:RBJ-Sk-invariance}
B_{\mathcal{S}_k u}
(\mathcal{S}_k x, \mathcal{S}_k y)
= B_u(x, y).
\end{equation}
\end{lemma}
 
\begin{proof}
The directed-path propagator $[B_v]_{ij}$ depends on $v$
only through the cutoffs $\mathbf{a}(v)$ and the
particle-label differences $n_{\ell'} - n_\ell$.

\noindent\textit{Case $k = 1$.}
Since $\mathcal{S}_1 u =
(t, \mathbf{a}, \mathbf{n}{+}\mathbf{1})$, the cutoffs
$\mathbf{a}$ are unchanged and the particle-label
differences $(n_j{+}1) - (n_i{+}1) = n_j - n_i$ are
unchanged.  Moreover, $\mathcal{S}_1$ does not shift
$\mathbf{a}$, so it acts trivially on spatial
coordinates: $B_{\mathcal{S}_1 u}(\mathcal{S}_1 x,
\mathcal{S}_1 y) = B_u(x, y)$.

\noindent\textit{Case $k = 2$.}
Since $\mathcal{S}_2 u =
(t{+}1, \mathbf{a}{+}\mathbf{1}, \mathbf{n})$, the
particle-label differences are unchanged but the
cutoffs shift to $\mathbf{a}{+}\mathbf{1}$.  Since
$\mathcal{S}_2$ acts on thresholds and spatial
coordinates identically to $\mathcal{T}$, and the
transition kernels do not depend on~$t$,
$\mathcal{S}_2$-invariance reduces to
$\mathcal{T}$-covariance
(Proposition~\ref{prop:directed-path-T}).
\end{proof}

\begin{lemma}\label{lem:RBJ-admissible}
The propagator $B_{\mathbf{a}, \mathbf{n}}$ is an admissible propagator
for the product graph $\mathcal{G}^m$.
\end{lemma}

\begin{proof}
By Proposition~\ref{prop:directed-path-T},
$\mathcal{T}$-covariance and $\mathcal{T}$-splitting hold.
The edge weights $(c_k, \lambda_k)$ are constant scalars
and the propagator satisfies $\mathcal{S}_k$-invariance
(Lemma~\ref{lem:RBJ-Sk-invariance}), so
$\mathcal{S}_k$-compatibility follows from
Proposition~\ref{prop:constant-scalar}.

It remains to verify regularity.  The directed-path sum
\eqref{eq:RBJ-propagator} is finite in the chain length
(at most $m - 1$ steps).  The support condition
$Q^d(x, x') = 0$ unless $x - x' \ge d$ implies that the
spatial chain is strictly decreasing along any nonzero
path, so the internal sums are finite for fixed endpoints.
The propagator entries $[B_{\mathbf{a}, \mathbf{n}}]_{ij}(r, r')$
decay geometrically: each factor
$Q^d(\xi_s, \xi_{s+1})$ is bounded by a polynomial in
$\xi_s - \xi_{s+1}$ times $\theta^{\xi_s - \xi_{s+1}}$,
so for any fixed $\vartheta \in (\theta, 1)$ the
displacements telescope to give
$|[B_{\mathbf{a}, \mathbf{n}}]_{ij}(r, r')| = O(\vartheta^{r'-r})$.
For the seed functions, the contour representation shows
that $\psi_{t,r,n}$ has compact support in~$x$ (specifically,
$\psi_{t,r,n}(x) = 0$ unless $r - t \le x \le r + n$),
while the bound \eqref{eq:RBJ-phi-bound} gives
$|\phi_{t,r,n}(x)| \le C\beta^{x-r}$ for $x > y_N$.
Together, these estimates ensure absolute convergence
of all sums defining $\Psi, \Phi, K$ in
\eqref{eq:Psi-def}--\eqref{eq:K-def}.
\end{proof}
 
\paragraph{\textbf{Corrected edge weights.}}
Since $\mathcal{S}_k$-invariance holds,
Proposition~\ref{prop:constant-scalar}(ii) gives
$P(\mathcal{S}_k u) = P(u)$ and
$\Lambda_k(u) = \lambda_k I_m = D_k(u)$.
The diamond data on $\mathcal{G}^m$ is therefore
$(C_k, \Lambda_k) = (c_k I_m, \lambda_k I_m)$.
Theorem~\ref{thm:product-graph} requires additionally
that the resolvent $(I - zK(u))^{-1}$ exist.  At $z = 1$,
the Fredholm determinant
$F_{t, \mathbf{a}, \mathbf{n}} \defeq
\det(I - K_{t, \mathbf{a}, \mathbf{n}})$ equals the
multipoint gap probability
$\mathbb{P}_y\bigl(\bigcap_i
\{Y_{n_i}(t) > a_i\}\bigr)$.  Since $K$ is trace class,
$F \ne 0$ if and only if the resolvent exists.  At each
time step, each particle attempts to jump one site to the
right, so $Y_{n_i}(t) \le y_{n_i} + t$.  If every jump
attempt for particles $1, \dotsc, n_m$ succeeds (an event
of probability $p^{n_m t}$), the strict integer ordering
prevents blocking and $Y_k(t) = y_k + t$ for
$k \le n_m$.  The resolvent at $z = 1$ therefore exists
precisely on
$\mathcal{R} \defeq \{(t, \mathbf{a}, \mathbf{n}) \in
\mathcal{V}^m : a_i < y_{n_i} + t
\text{ for all } i\}$.
The mixed diamond at a vertex~$u$ requires the resolvent
at the four vertices $u$, $\mathcal{S}_1 u$,
$\mathcal{S}_2 u$, $\mathcal{S}_1\mathcal{S}_2 u$
entering~\eqref{eq:dressed-C-explicit}, together with
the $\mathcal{T}$-shifts $\mathcal{T}u$,
$\mathcal{T}\mathcal{S}_1 u$,
$\mathcal{T}\mathcal{S}_2 u$ from the proof
of~\eqref{eq:M-inverse}.\footnote{The explicit
formula~\eqref{eq:dressed-C-explicit} involves the
resolvent only at the four $\mathcal{S}$-shifted
vertices; a direct verification of the mixed diamond
from this formula would yield the larger domain
$a_i < y_{n_i+1} + t$.}
Since $\mathcal{S}_2$ preserves~$\mathcal{R}$ and
$\mathcal{T}$ tightens each threshold by~$1$, the
binding constraint is
$\mathcal{T}\mathcal{S}_1 u \in \mathcal{R}$, i.e.\
$a_i < y_{n_i+1} + t - 1$ for all~$i$.

\subsection{Multipoint equation}
\label{sec:RBJ-multipoint}
 
The dressed observable at $z = 1$,
\[
\mathcal{M}(u)
= I + \Phi(u)(I - K(u))^{-1}\Psi(u)
\in \End(\mathbb{R}^m),
\]
has dressed edge weights
\begin{equation}\label{eq:RBJ-dressed-Ck}
\mathcal{M}_k(u)
= c_k\mathcal{M}(\mathcal{T}u)^{-1}
  \mathcal{M}(\mathcal{S}_k u).
\end{equation}
For readability, we write
$\mathcal{M}_{t, \mathbf{a}, \mathbf{n}} \defeq
\mathcal{M}(t, \mathbf{a}, \mathbf{n})$.

\begin{theorem}\label{thm:RBJ-H5}
The dressed observable satisfies the mixed diamond
equation
\begin{equation}\label{eq:H5-multipoint}
\bigl(q\mathcal{M}_{t, \mathbf{a}+\mathbf{1},
\mathbf{n}+\mathbf{1}}^{-1}
+ p\mathcal{M}_{t+1, \mathbf{a}+\mathbf{2},
\mathbf{n}}^{-1}\bigr)
\mathcal{M}_{t+1, \mathbf{a}+\mathbf{1},
\mathbf{n}+\mathbf{1}}
- \mathcal{M}_{t, \mathbf{a}+\mathbf{1},
\mathbf{n}}^{-1}
\bigl(q\mathcal{M}_{t+1, \mathbf{a}+\mathbf{1},
\mathbf{n}}
+ p\mathcal{M}_{t, \mathbf{a},
\mathbf{n}+\mathbf{1}}\bigr) = 0,
\end{equation}
at every $(t, \mathbf{a}, \mathbf{n})$ with
$a_i < y_{n_i+1} + t - 1$ for all~$i$.
\end{theorem}

\begin{proof}
By Theorem~\ref{thm:product-graph}, the pair
$(\mathcal{M}_k, \Lambda_k)$ satisfies the diamond
equations on~$\mathcal{R}$.
Since $\Lambda_k = \lambda_k I_m$, the scalar factors
commute with $\mathcal{M}_k$, and the mixed diamond
equation reduces to
\begin{equation}\label{eq:RBJ-mixed-scalar-Lambda}
\lambda_2\bigl[
\mathcal{M}_1(u) - \mathcal{M}_1(\mathcal{S}_2 u)
\bigr]
+ \lambda_1\bigl[
\mathcal{M}_2(\mathcal{S}_1 u)
- \mathcal{M}_2(u)
\bigr] = 0.
\end{equation}
Substituting \eqref{eq:RBJ-dressed-Ck} and the explicit
shifts, every term carries a common factor of
$(1{-}\theta)^{-1}\theta$.  Dividing by this factor
gives
\begin{align*}
-p\bigl[
&\mathcal{M}_{t, \mathbf{a}+\mathbf{1},
\mathbf{n}}^{-1}
\mathcal{M}_{t, \mathbf{a},
\mathbf{n}+\mathbf{1}}
- \mathcal{M}_{t+1, \mathbf{a}+\mathbf{2},
\mathbf{n}}^{-1}
\mathcal{M}_{t+1, \mathbf{a}+\mathbf{1},
\mathbf{n}+\mathbf{1}}
\bigr] \\
+ q\bigl[
&\mathcal{M}_{t, \mathbf{a}+\mathbf{1},
\mathbf{n}+\mathbf{1}}^{-1}
\mathcal{M}_{t+1, \mathbf{a}+\mathbf{1},
\mathbf{n}+\mathbf{1}}
- \mathcal{M}_{t, \mathbf{a}+\mathbf{1},
\mathbf{n}}^{-1}
\mathcal{M}_{t+1, \mathbf{a}+\mathbf{1},
\mathbf{n}}
\bigr] = 0.
\end{align*}
The middle two terms share the right factor
$\mathcal{M}_{t+1, \mathbf{a}+\mathbf{1},
\mathbf{n}+\mathbf{1}}$, and the
first and last share the left factor
$\mathcal{M}_{t, \mathbf{a}+\mathbf{1},
\mathbf{n}}^{-1}$.  Factoring yields
\eqref{eq:H5-multipoint}.
\end{proof}

For $m = 1$, the dressed observable $\mathcal{M}_{t,a,n}$
is scalar.

\begin{corollary}\label{cor:RBJ-scalar}
The one-point Fredholm determinant
$F_{t, a, n} = \det(I - K_{t, a, n})$ satisfies
\begin{equation}\label{eq:RBJ-scalar-HM}
p F_{t, a, n+1} F_{t+1, a+2, n}
+ q F_{t+1, a+1, n} F_{t, a+1, n+1}
- F_{t, a+1, n} F_{t+1, a+1, n+1}
= 0,
\end{equation}
at every $(t, a, n)$ with $a < y_{n+1} + t - 1$.
\end{corollary}

\begin{proof}
Proposition~\ref{prop:scalar-HM} applies at $z = 1$: the
boundary condition $F(T^\ell v) \to 1$ as
$\ell \to -\infty$ holds because the thresholds decrease
to $-\infty$ and the kernel vanishes.  The coefficients are
$\alpha_{12} = c_1 \lambda_2 =
-(1{-}\theta)^{-1}p\theta$
and
$\alpha_{21} = c_2 \lambda_1 =
q(1{-}\theta)^{-1}\theta$.
Substituting the lattice points into
\eqref{eq:HM-variable} and computing
$\alpha_{12} - \alpha_{21}
= -(1{-}\theta)^{-1}\theta$
(using $p + q = 1$), then dividing by
$-(1{-}\theta)^{-1}\theta$, gives
\eqref{eq:RBJ-scalar-HM}.
\end{proof}

\begin{remark}[Structure of the multipoint equation]
\label{rem:RBJ-structure}
Equation \eqref{eq:H5-multipoint} is an
$m \times m$ matrix equation coupling the dressed
observables at six lattice points, with coefficients
$p$ and $q = 1 - p$ given by the jump and stay
probabilities.  For $m \ge 2$, the equation is
intrinsically noncommutative: the matrix inverses and
products do not simplify to a scalar relation.
\end{remark}

\section{Inhomogeneous Bernoulli Jumps}\label{sec:IBJ}

\paragraph{\textbf{System description.}}
The discrete-time totally asymmetric simple
exclusion process with particle- and
time-inhomogeneous Bernoulli jump rates, blocking
interaction, and sequential update
consists of $N$ particles on $\Z$.  The configuration
at time $t \in \Z_{\ge 0}$ is denoted by
\[
Y(t) = (Y_1(t) > Y_2(t) > \cdots > Y_N(t)),
\]
where $Y_k(t)$ is the position of the $k$-th particle
from the right.

At each discrete time step $t \ge 1$, the particles are
updated sequentially in the order
$k = 1, 2, \dotsc, N$, so the rightmost particle is
updated first.  Unlike the homogeneous model
of~\S\ref{sec:RBJ}, the jump probability depends on
both the particle label and the time step.  The
evolution is given by
\[
Y_k(t) = \min\{Y_{k-1}(t) - 1, Y_k(t{-}1)
+ \xi(t,k)\}, \qquad Y_0(t) = \infty,
\]
where $\{\xi(t,k)\}_{t \ge 1, 1 \le k \le N}$ are
independent Bernoulli random variables with
$\PP(\xi(t,k) = 1) = p_t q_k/(1 + p_t q_k)$.  The
parameters $(p_t)_{t \ge 1}$ and $(q_1, \dotsc, q_N)$
are positive with $p_t q_k < 1$ and $q_k > 1$ for
all $t, k$.

\paragraph{\textbf{Fredholm determinant formula.}}
Fix a time $t \ge 0$ and an initial configuration
$Y(0) = y = (y_1 > y_2 > \cdots > y_N)$.  Fix
$m \in \{1, \dotsc, N{-}1\}$, particle labels
$\mathbf{n} = (n_1, \dotsc, n_m)$ with
$1 \le n_1 < \cdots < n_m < N$, and spatial
thresholds $\mathbf{a} = (a_1, \dotsc, a_m) \in \Z^m$.
The multipoint distribution is given by a Fredholm
determinant on
$\ell^2(\{1, \dotsc, m\} \times \Z)$~\cite{BLSZ23}:
\[
\PP_y\Bigl(\bigcap_{i=1}^m
\{Y_{n_i}(t) > a_i\}\Bigr)
= \det\bigl(I - \bar{\chi}_{\mathbf{a}}
K^{\mathrm{Ext}}
\bar{\chi}_{\mathbf{a}}\bigr)
_{\ell^2(\{1, \dotsc, m\} \times \Z)},
\]
where $\bar{\chi}_{\mathbf{a}}(i, x) =
\mathbf{1}_{x \le a_i}$ and the extended kernel has
block entries
\[
K^{\mathrm{Ext}}(i, x; j, x')
= -Q_{(n_i, n_j]}(x, x')\mathbf{1}_{i < j}
+ \bigl(\mathcal{S}_{[1,n_i],(0,t]} \circ
\bar{\mathcal{S}}^{\mathrm{epi}(y)}_{[1,n_j],(0,t]}
\bigr)(x, x').
\]
The operators defining the kernel are as follows.

\medskip
\noindent\textit{The transition matrix $Q_k$.}
For each $k$, the operator $Q_k$ acts on
$\ell^2(\Z)$ by
\[
Q_k(x, y) = q_k^{y-x}\mathbf{1}_{y < x},
\]
and the composed operator is
$Q_{(m,n]} = Q_{m+1} \circ \cdots \circ Q_n$.

\medskip
\noindent\textit{The scattering operators
$\mathcal{S}_{[j,k],(r,t]}$ and
$\bar{\mathcal{S}}_{[j,k],(r,t]}$.}
These are defined by
\[
\mathcal{S}_{[j,k],(r,t]}(x, y)
= \oint_{\Gamma_0} \frac{\diff z}{2\pi\mathrm{i}z}
\frac{\prod_{\ell=j}^k(q_\ell - z)
\prod_{s=r+1}^t(1 + p_s z)}
{z^{x - y + k - j + 1}},
\]
\[
\bar{\mathcal{S}}_{[j,k],(r,t]}(x, y)
= -\prod_{\ell=j}^k(q_\ell - 1)
\oint_{\Gamma_{\mathbf{q}}} \frac{\diff z}{2\pi\mathrm{i}z}
\frac{z^{y - x + k - j + 1}}
{\prod_{\ell=j}^k(q_\ell - z)
\prod_{s=r+1}^t(1 + p_s z)},
\]
where $\Gamma_0$ is a positively oriented contour
enclosing $0$ and excluding
$\{q_j, \dotsc, q_k\}$, and $\Gamma_{\mathbf{q}}$
is a positively oriented contour enclosing
$\{q_j, \dotsc, q_k\}$ and excluding $0$ and
$\{-p_s^{-1} : r < s \le t\}$.

\medskip
\noindent\textit{The epigraph operator
$\bar{\mathcal{S}}^{\mathrm{epi}(y)}_{[1,n],(r,t]}$.}
Let $W$ be the geometric random walk moving strictly
to the left with inhomogeneous transition probabilities
\[
\PP(W_\ell = y \mid W_{\ell-1} = x)
= (q_\ell - 1)q_\ell^{y-x}\mathbf{1}_{y < x},
\qquad 1 \le \ell \le N,
\]
and let
$\tau = \min\{j \in \{0, \dotsc, N\} : W_j > y_{j+1}\}$,
with the convention $y_{N+1} \defeq -\infty$
(so that $\tau \le N$ always holds).  The epigraph
operator is
\[
\bar{\mathcal{S}}^{\mathrm{epi}(y)}_{[1,n],(r,t]}
(x, a)
= \frac{\E_{W_0 = x}\bigl[
\bar{\mathcal{S}}_{[\tau+1,n],(r,t]}(W_\tau, a)
\mathbf{1}_{\tau < n}\bigr]}
{\prod_{\ell=1}^n(q_\ell - 1)}.
\]

\paragraph{\textbf{References.}}
Bisi--Liao--Saenz--Zygouras~\cite{BLSZ23} prove a Fredholm determinant formula for multipoint distributions of the particle- and time-inhomogeneous Bernoulli dTASEP using dual column RSK, intertwinings, and non-intersecting path ensembles.

\subsection{Kernel reformulation}
\label{sec:IBJ-kernel-reform}

The extended-kernel Fredholm determinant admits a
reformulation in the form of the product graph construction
of~\S\ref{sec:product-graph}.
The functions $\psi, \phi$ below provide the seed data,
while the lower-triangular operator $B_{\mathbf{a}, \mathbf{n}}$ will
be verified to satisfy the admissibility conditions
of~Definition~\ref{def:admissible-prop}.

\paragraph{\textbf{Seed data.}}
Define
\begin{align*}
\psi_{t,a,n}(x) &\defeq
\mathcal{S}_{[1,n],(0,t]}(a, x)
= \oint_{\Gamma_0} \frac{\diff z}{2\pi\mathrm{i}z}
\frac{F_n(z)G_t(z)}{z^{a - x + n}}, \\
\phi_{t,a,n}(x) &\defeq
\bar{\mathcal{S}}^{\mathrm{epi}(y)}_{[1,n],(0,t]}
(x, a),
\end{align*}
where $F_n(z) \defeq \prod_{\ell=1}^n(q_\ell - z)$,
$G_t(z) \defeq \prod_{s=1}^t(1 + p_s z)$, and we
suppress the dependence of $\phi$ on the initial
condition $y$.

\paragraph{\textbf{Propagator.}}
The propagator is an instance of the directed-path
propagator
(Definition~\ref{def:directed-path-propagator}) with
transition kernels
$Q_{\ell,\ell'} = Q_{(n_\ell, n_{\ell'}]}$ for
$\ell < \ell'$.  Explicitly, for $i > j$,
\begin{equation}\label{eq:IBJ-propagator}
[B_{\mathbf{a}, \mathbf{n}}]_{ij}(r, r')
\defeq \sum_{k=1}^{i-j}
  \sum_{j = \ell_0 < \ell_1 < \cdots < \ell_k = i}
  \sum_{\substack{\xi_1 \le a_{\ell_1} \\ \vdots \\
  \xi_{k-1} \le a_{\ell_{k-1}}}}
  \prod_{s=0}^{k-1}
  \bigl(-Q_{(n_{\ell_s}, n_{\ell_{s+1}}]}
  (\xi_s, \xi_{s+1})\bigr),
\end{equation}
with $\xi_0 = r'$, $\xi_k = r$, and $[B_{\mathbf{a}, \mathbf{n}}]_{ij} = 0$ for $i \le j$.
The internal vertices $\xi_1, \dotsc, \xi_{k-1}$ are
constrained by the cutoffs, but the endpoints
$r, r'$ are arbitrary, so
$[B_{\mathbf{a}, \mathbf{n}}]_{ij}(r, r')$ is defined for all
$r, r' \in \Z$.  The propagator depends on $\mathbf{a}$
and $\mathbf{n}$ (through the cutoffs $a_{\ell_s}$ and
the operators $Q_{(n_{\ell_s}, n_{\ell_{s+1}}]}$) but
not on $t$.

On the cutoff region $r \le a_i$, $r' \le a_j$, the
entries of $I + B_{\mathbf{a}, \mathbf{n}}$ coincide with those of
$((I + \bar{\chi}_{\mathbf{a}}L
\bar{\chi}_{\mathbf{a}})^{-1})^{\!\top}$, where $L$
is the strictly upper-triangular block operator on
$\ell^2(\{n_1, \dotsc, n_m\} \times \Z)$ with entries
$L_{ij}(x, x') = Q_{(n_i, n_j]}(x, x')
\mathbf{1}_{i < j}$.

\begin{lemma}\label{lem:IBJ-kernel}
We have
\[
\PP_y\Bigl(\bigcap_{i=1}^m
\{Y_{n_i}(t) > a_i\}\Bigr)
= \det(I - K_{t, \mathbf{a}, \mathbf{n}})
_{\ell^2(\Z)},
\]
where
\[
K_{t, \mathbf{a}, \mathbf{n}}(x, x')
= \sum_{1 \le j \le i \le m}
  \sum_{r \le a_i} \sum_{r' \le a_j}
  \psi_{t, r, n_i}(x)
  \bigl[\delta_{ij}\delta_{r,r'}
  + [B_{\mathbf{a}, \mathbf{n}}]_{ij}(r, r')\bigr]
  \phi_{t, r', n_j}(x').
\]
This identifies $K_{t, \mathbf{a}, \mathbf{n}}$ with
the kernel~\eqref{eq:K-def} on the product graph
$\mathcal{G}^m$, with $\psi, \phi$ as the seed data
and $B_{\mathbf{a}, \mathbf{n}}$ as the propagator.
\end{lemma}

\begin{proof}
The argument is analogous to that of
Lemma~\ref{lem:RBJ-kernel}: the block-triangular
decomposition of
$I - \bar{\chi}_{\mathbf{a}}K^{\mathrm{Ext}}
\bar{\chi}_{\mathbf{a}}$,
Sylvester's identity, and the transpose convention
produce the claimed single-space kernel.
The random walk $Q_{(n_i,n_j]}$ and the epigraph
hitting operator are the inhomogeneous counterparts
of~\S\ref{sec:RBJ}; only the scattering operators
$\mathcal{S}_{[1,n],(0,t]}$ and
$\bar{\mathcal{S}}_{[1,n],(0,t]}$ differ, and they
enter only through the seed data $\psi, \phi$.
The Schatten-class estimates are analogous to those of
Lemma~\ref{lem:RBJ-kernel}, with the geometric rate
$\theta$ replaced by $q_{\min}^{-1} \in (0, 1)$.  The
composed transition kernel $Q_{(n_i,n_j]}$ is supported
on $x - x' \ge n_j - n_i$ and satisfies
$|Q_{(n_i,n_j]}(x,x')| \le C\beta^{x-x'}$ for any
$\beta \in (q_{\min}^{-1}, 1)$.  The contour formula for
$\psi$ gives compact support of width $n + t$, and the
contour estimate for $\phi$ gives
$|\phi_{t,a,n}(x)| \le C\beta^{x-a}
\mathbf{1}_{x > y_N}$.  With these bounds, the same
conjugation and Sylvester reduction as in
Lemma~\ref{lem:RBJ-kernel} apply.
\end{proof}

\subsection{Base graph and seed data}
\label{sec:IBJ-base-graph}

The lattice structure governing the seed data under
shifts of the parameters $(t, a, n)$ is identified next,
together with a verification that $\psi$ and $\phi$ satisfy the
linear problems of the diamond framework.

Let $\mathcal{V} \subseteq \Z^3$ be the lattice
generated by the shifts
\[
T = e^{\pa_a}, \qquad
S_1 = e^{\pa_n}, \qquad
S_2 = e^{\pa_t + \pa_a},
\]
so that for $u = (t, a, n) \in \mathcal{V}$,
\[
Tu = (t, a{+}1, n), \qquad
S_1 u = (t, a, n{+}1), \qquad
S_2 u = (t{+}1, a{+}1, n).
\]
Take $E = \mathbb{F} = \R$ and define the
scalar edge weights
\begin{equation}\label{eq:IBJ-edge-weights}
(c_1, c_2) = \bigl(q_{n+1}, 1\bigr),
\qquad
(\lambda_1, \lambda_2)
= \bigl(1, -p_{t+1}\bigr).
\end{equation}
Note that $c_1$ depends on $n$ and $\lambda_2$ depends
on $t$, so the edge weights are vertex-dependent:
$c_1(u) = q_{n+1}$, $\lambda_1(u) = 1$,
$c_2(u) = 1$, $\lambda_2(u) = -p_{t+1}$.

The diamond equations
\eqref{eq:diamond-C-ij}--\eqref{eq:diamond-mixed-ij}
hold for $(i,j) = (1,2)$: \eqref{eq:diamond-C-ij} holds since
$c_1(Tu)c_2(S_1 u) = q_{n+1}
= c_2(Tu)c_1(S_2 u)$, \eqref{eq:diamond-L-ij} since
$\lambda_1(u)\lambda_2(S_1 u) = -p_{t+1}
= \lambda_2(u)\lambda_1(S_2 u)$, and \eqref{eq:diamond-mixed-ij}
since both $c_1(u)\lambda_2(S_1 u)
= -q_{n+1}p_{t+1}
= \lambda_2(Tu)c_1(S_2 u)$ and
$\lambda_1(Tu)c_2(S_1 u) = 1 = c_2(u)\lambda_1(S_2 u)$.

\begin{lemma}\label{lem:IBJ-seed-linear}
Let $H = \ell^2(\Z)$.  The seed functions
$\psi_{t,a,n}$ and $\phi_{t,a,n}$ defined in
\S\ref{sec:IBJ-kernel-reform} satisfy the linear
problem \eqref{eq:Linear-Psi-3} and the adjoint
linear problem \eqref{eq:Linear-Phi-Prob} on
$\mathcal{V}$ with edge weights
$(c_k, \lambda_k)$.  Explicitly:
\begin{align}
\psi_{t,a,n+1}
&= q_{n+1}\psi_{t,a+1,n}
   - \psi_{t,a,n},
\label{eq:IBJ-psi-S1} \\
\psi_{t+1,a+1,n}
&= \psi_{t,a+1,n}
   + p_{t+1}\psi_{t,a,n},
\label{eq:IBJ-psi-S2}
\end{align}
and
\begin{align}
\phi_{t,a+1,n}
&= q_{n+1}\phi_{t,a,n+1}
   - \phi_{t,a+1,n+1},
\label{eq:IBJ-phi-S1} \\
\phi_{t,a+1,n}
&= \phi_{t+1,a+1,n}
   + p_{t+1}\phi_{t+1,a+2,n}.
\label{eq:IBJ-phi-S2}
\end{align}
\end{lemma}

\begin{proof}
Recall $F_n(z) = \prod_{\ell=1}^n(q_\ell - z)$ and
$G_t(z) = \prod_{s=1}^t(1 + p_s z)$.  We verify each
identity in turn.

\medskip
\noindent\textit{Verification of~\eqref{eq:IBJ-psi-S1}.}
Since $F_{n+1}(z) = (q_{n+1} - z)F_n(z)$:
\begin{align*}
\psi_{t,a,n+1}(x)
&= \oint_{\Gamma_0} \frac{\diff z}{2\pi\mathrm{i}z}
  \frac{(q_{n+1} - z)F_n(z)G_t(z)}
  {z^{a - x + n + 1}} \\
&= q_{n+1} \oint_{\Gamma_0}
  \frac{\diff z}{2\pi\mathrm{i}z}
  \frac{F_n(z)G_t(z)}{z^{(a+1) - x + n}}
  - \oint_{\Gamma_0}
  \frac{\diff z}{2\pi\mathrm{i}z}
  \frac{F_n(z)G_t(z)}{z^{a - x + n}} \\
&= q_{n+1}\psi_{t,a+1,n}(x)
  - \psi_{t,a,n}(x).
\end{align*}

\medskip
\noindent\textit{Verification of~\eqref{eq:IBJ-psi-S2}.}
Since $G_{t+1}(z) = (1 + p_{t+1} z)G_t(z)$:
\begin{align*}
\psi_{t+1,a+1,n}(x)
&= \oint_{\Gamma_0} \frac{\diff z}{2\pi\mathrm{i}z}
  \frac{F_n(z)G_t(z)(1 + p_{t+1} z)}
  {z^{(a+1) - x + n}} \\
&= \psi_{t,a+1,n}(x)
  + p_{t+1}\psi_{t,a,n}(x).
\end{align*}

\medskip
\noindent\textit{Verification of~\eqref{eq:IBJ-phi-S1}.}
This requires more care, as the indicator
$\mathbf{1}_{\tau < n}$ in the definition of $\phi$
introduces a boundary term when $n$ is shifted.

We first establish a recurrence for the free operator
$\bar{\mathcal{S}}_{[j,n],(0,t]}$.  Writing
$q_n = (q_n - z) + z$ in the numerator of the contour
representation gives
\begin{equation}\label{eq:IBJ-Sbar-n-rec}
q_n\bar{\mathcal{S}}_{[j,n],(0,t]}(w, a)
- \bar{\mathcal{S}}_{[j,n],(0,t]}(w, a{+}1)
= (q_n - 1)
  \bar{\mathcal{S}}_{[j,n-1],(0,t]}(w, a{+}1).
\end{equation}
Moreover, $\bar{\mathcal{S}}_{[n+1,n],(0,t]}
\equiv 0$: when the lower index exceeds the upper,
the integrand in the contour representation has no
poles inside $\Gamma_{\mathbf{q}}$, and the integral
vanishes.

The right-hand side of~\eqref{eq:IBJ-phi-S1} equals
\begin{align*}
&q_{n+1}\phi_{t,a,n+1}(x) - \phi_{t,a+1,n+1}(x) \\
&\quad = \frac{\E_{W_0 = x}\bigl[
\bigl(q_{n+1}
\bar{\mathcal{S}}_{[\tau+1,n+1],(0,t]}(W_\tau, a)
- \bar{\mathcal{S}}_{[\tau+1,n+1],(0,t]}
(W_\tau, a{+}1)\bigr)
\mathbf{1}_{\tau < n+1}\bigr]}
{\prod_{\ell=1}^{n+1}(q_\ell - 1)}.
\end{align*}
Applying~\eqref{eq:IBJ-Sbar-n-rec} with $n$ replaced
by $n{+}1$ and $j = \tau + 1$ inside the expectation
gives
\[
= \frac{(q_{n+1} - 1)
\E_{W_0 = x}\bigl[
\bar{\mathcal{S}}_{[\tau+1,n],(0,t]}(W_\tau, a{+}1)
\mathbf{1}_{\tau < n+1}\bigr]}
{\prod_{\ell=1}^{n+1}(q_\ell - 1)}.
\]
On the event $\{\tau = n\}$, the factor
$\bar{\mathcal{S}}_{[n+1,n],(0,t]} = 0$ and the
contribution vanishes, so
$\mathbf{1}_{\tau < n+1}$ reduces to
$\mathbf{1}_{\tau < n}$.  The factor
$(q_{n+1} - 1)$ cancels with the last factor of
$\prod_{\ell=1}^{n+1}(q_\ell - 1)$, giving
\[
\frac{\E_{W_0 = x}\bigl[
\bar{\mathcal{S}}_{[\tau+1,n],(0,t]}(W_\tau, a{+}1)
\mathbf{1}_{\tau < n}\bigr]}
{\prod_{\ell=1}^{n}(q_\ell - 1)}
= \phi_{t,a+1,n}(x).
\]

\medskip
\noindent\textit{Verification of~\eqref{eq:IBJ-phi-S2}.}
Since neither side involves a shift in $n$, the
indicator $\mathbf{1}_{\tau < n}$ passes through
unchanged.  Using
$G_t(z)^{-1} - G_{t-1}(z)^{-1}
= -p_t zG_t(z)^{-1}$ in the contour representation
gives
\begin{equation}\label{eq:IBJ-Sbar-t-rec}
\bar{\mathcal{S}}_{[j,n],(0,t-1]}(w, a)
= \bar{\mathcal{S}}_{[j,n],(0,t]}(w, a)
+ p_t\bar{\mathcal{S}}_{[j,n],(0,t]}(w, a{+}1).
\end{equation}
The right-hand side of~\eqref{eq:IBJ-phi-S2} equals
\begin{align*}
&\phi_{t+1,a+1,n}(x)
+ p_{t+1}\phi_{t+1,a+2,n}(x) \\
&\quad = \frac{\E_{W_0 = x}\bigl[
\bigl(\bar{\mathcal{S}}_{[\tau+1,n],(0,t+1]}
(W_\tau, a{+}1)
+ p_{t+1}
\bar{\mathcal{S}}_{[\tau+1,n],(0,t+1]}
(W_\tau, a{+}2)\bigr)
\mathbf{1}_{\tau < n}\bigr]}
{\prod_{\ell=1}^{n}(q_\ell - 1)}.
\end{align*}
Applying~\eqref{eq:IBJ-Sbar-t-rec} with $t$ replaced
by $t{+}1$, $j = \tau + 1$, and $a$ replaced by
$a{+}1$ collapses the two terms inside the
expectation to
$\bar{\mathcal{S}}_{[\tau+1,n],(0,t]}(W_\tau, a{+}1)$,
giving
\[
= \frac{\E_{W_0 = x}\bigl[
\bar{\mathcal{S}}_{[\tau+1,n],(0,t]}(W_\tau, a{+}1)
\mathbf{1}_{\tau < n}\bigr]}
{\prod_{\ell=1}^{n}(q_\ell - 1)}
= \phi_{t,a+1,n}(x).
\]
\end{proof}

\subsection{Admissible propagators and corrected
edge weights}
\label{sec:IBJ-propagators}

The product graph $\mathcal{G}^m$ has vertex set
$\mathcal{V}^m$ and is generated by the diagonal
shifts
\[
\mathcal{T}u
= (t, \mathbf{a} + \mathbf{1}, \mathbf{n}),
\qquad
\mathcal{S}_1 u
= (t, \mathbf{a}, \mathbf{n} + \mathbf{1}),
\qquad
\mathcal{S}_2 u
= (t{+}1, \mathbf{a} + \mathbf{1}, \mathbf{n}),
\]
where $u = (t, \mathbf{a}, \mathbf{n}) \in
\mathcal{V}^m$ and
$\mathbf{1} = (1, \dotsc, 1) \in \Z^m$.
The block-diagonal edge weights on $\mathcal{G}^m$ are
$C_k(u) = \mathrm{Diag}(c_k(u_1), \dotsc, c_k(u_m))$
and
$D_k(u) = \mathrm{Diag}(\lambda_k(u_1), \dotsc,
\lambda_k(u_m))$.  Explicitly,
\begin{gather*}
C_1(u) = \mathrm{Diag}(q_{n_1+1}, \dotsc, q_{n_m+1}),
\qquad
C_2(u) = I_m, \\
D_1(u) = I_m,
\qquad
D_2(u) = -p_{t+1} I_m.
\end{gather*}

The propagator~\eqref{eq:IBJ-propagator} is an instance
of the directed-path propagator
(Definition~\ref{def:directed-path-propagator}) with
transition kernels
$Q_{\ell,\ell'} = Q_{(n_\ell, n_{\ell'}]}$ for
$\ell < \ell'$.  The standing hypotheses of
\S\ref{sec:directed-path} are satisfied: each $Q_k$ is
translation invariant (since
$Q_k(x,y) = q_k^{y-x}\mathbf{1}_{y<x}$ depends only
on $y - x$), and $\mathcal{T}$-invariance holds because
the transition kernels $Q_{(n_\ell, n_{\ell'}]}$ depend
only on the particle labels, which are unchanged
by~$\mathcal{T}$.

To apply Theorem~\ref{thm:product-graph}, we verify
the four admissibility conditions.

\begin{lemma}
\label{lem:IBJ-T-cov-split}
The propagator $B_{\mathbf{a}, \mathbf{n}}$ satisfies
$\mathcal{T}$-covariance and $\mathcal{T}$-splitting
on the product graph $\mathcal{G}^m$.
\end{lemma}

\begin{proof}
The standing hypotheses of \S\ref{sec:directed-path}
are verified above, so
Proposition~\ref{prop:directed-path-T} applies.
\end{proof}

\begin{lemma}
\label{lem:IBJ-S2-inv}
For all $x, y \in \mathcal{V}^m$,
\begin{equation}\label{eq:IBJ-S2-invariance}
B_{\mathcal{S}_2 u}
(\mathcal{S}_2 x, \mathcal{S}_2 y)
= B_u(x, y).
\end{equation}
\end{lemma}

\begin{proof}
Since $\mathcal{S}_2$ acts on thresholds and spatial
coordinates identically to $\mathcal{T}$, and the
transition kernels do not depend on~$t$,
$\mathcal{S}_2$-invariance reduces to
$\mathcal{T}$-covariance
(Proposition~\ref{prop:directed-path-T}).
\end{proof}

Since the propagator is $\mathcal{S}_2$-invariant and
$C_2 = I_m$, the $\mathcal{S}_2$-compatibility
condition reduces
to~\eqref{eq:Sk-compat-reduced} by the same
cancellation as in
Proposition~\ref{prop:constant-scalar}: the factor
$D_2 = -p_{t+1}I_m$ is scalar, and since $\mathcal{T}$
does not shift $t$, both $D_2$-values in the
compatibility condition are equal and cancel.  The
remaining identity follows from
$\mathcal{T}$-covariance and $\mathcal{T}$-splitting.
However, the $\mathcal{S}_1$-direction requires a direct
verification, since $c_1(u) = q_{n+1}$ depends on the
particle label and
Proposition~\ref{prop:constant-scalar} does not apply.

\paragraph{\textbf{Componentwise identities.}}
The following identities are used in the
$\mathcal{S}_1$-compatibility verification.
For each $k$, the operator
$Q_k(x, y) = q_k^{y-x}\mathbf{1}_{y < x}$ has
two-sided inverse
$Q_k^{-1}(x, y) = q_k\delta_{x, y-1}
- \delta_{x, y}$,
and the collapse identity
\begin{equation}\label{eq:IBJ-collapse}
q_kQ_k(x, r) - Q_k(x, r{+}1) = \delta_{x, r+1}
\end{equation}
holds for all $x, r \in \Z$ (being the
evaluation of $Q_k Q_k^{-1} = I$ at $(x, r{+}1)$).
For $m < n$, the edge-shift identity
\begin{equation}\label{eq:IBJ-edge-shift}
Q_{(m+1, n+1]}
= Q_{m+1}^{-1} \circ Q_{(m, n]} \circ Q_{n+1}
\end{equation}
follows by expanding both sides as compositions.

\begin{lemma}
\label{lem:IBJ-S1-compat}
The propagator satisfies the
$\mathcal{S}_1$-compatibility condition
of~Definition~\ref{def:admissible-prop}.
\end{lemma}

\begin{proof}
Since $C_1 = \mathrm{Diag}(q_{n_1+1}, \dotsc, q_{n_m+1})$
is block-diagonal and $D_1 = I_m$, the
$\mathcal{S}_1$-compatibility condition from
Definition~\ref{def:admissible-prop} decomposes into
independent scalar conditions on each $(i,j)$-block.
Write $B = B_{\mathbf{a}, \mathbf{n}}$ and
$B' = B_{\mathbf{a}, \mathbf{n}+\mathbf{1}}$ for the
propagators at labels $\mathbf{n}$ and
$\mathbf{n}+\mathbf{1}$.
For diagonal blocks $i = j$, both sides vanish because
$B$ and $B'$ are strictly lower-triangular.  For
$i > j$, the condition evaluated at spatial arguments
$(r, s)$ reads
\begin{align}
q_{n_i+1}[B']_{ij}(r, s)
&- [B']_{ij}(r{+}1, s)
- [B]_{ij}(r{+}1, s{+}1)q_{n_j+1}
+ [B]_{ij}(r{+}1, s) \nonumber \\
&= \sum_{j < \mu < i}
[B]_{i\mu}(r{+}1, a_\mu{+}1)
q_{n_\mu+1}
[B']_{\mu j}(a_\mu, s),
\label{eq:IBJ-S1-compat-component}
\end{align}
\noindent\textit{Step 1: Endpoint collapse and edge-shift.}
Write $[B']_{ij}(r, s) = \sum_\gamma w'(\gamma)$ for
the directed-path
expansion~\eqref{eq:IBJ-propagator} at shifted labels
$\mathbf{n}{+}\mathbf{1}$, where each path $\gamma$
traverses a chain
$j = \ell_0 < \ell_1 < \cdots < \ell_k = i$
with internal vertices
$\xi_1 \le a_{\ell_1}, \dotsc,
\xi_{k-1} \le a_{\ell_{k-1}}$, initial point
$\xi_0 = s$, endpoint $\xi_k = r$, and weight
\[
w'(\gamma) = \prod_{p=0}^{k-1}
\bigl(-Q_{(n_{\ell_p}+1, n_{\ell_{p+1}}+1]}
(\xi_p, \xi_{p+1})\bigr).
\]
Since the internal vertices are constrained only by
the thresholds $a_{\ell_1}, \dotsc, a_{\ell_{k-1}}$
and not by the endpoint, the sums defining
$[B']_{ij}(r, s)$ and $[B']_{ij}(r{+}1, s)$ range
over the same chains and internal vertices.  For each
chain and choice of internal vertices, there is one
path $\gamma$ ending at $r$ and one path
$\tilde\gamma$ ending at $r{+}1$; their weights
share all factors except the last, which is
$Q_{(n_{\ell_{k-1}}+1, n_i+1]}(\xi_{k-1}, r)$ for
$\gamma$ and
$Q_{(n_{\ell_{k-1}}+1, n_i+1]}(\xi_{k-1}, r{+}1)$
for $\tilde\gamma$.  Applying the collapse
identity~\eqref{eq:IBJ-collapse} to $Q_{n_i+1}$,
composing on the left with
$Q_{(n_{\ell_{k-1}}+1, n_i]}$, and recognizing
$Q_{(n_{\ell_{k-1}}+1, n_i]} Q_{n_i+1}
= Q_{(n_{\ell_{k-1}}+1, n_i+1]}$ gives
\begin{align*}
q_{n_i+1}
&Q_{(n_{\ell_{k-1}}+1, n_i+1]}(\xi_{k-1}, r)
- Q_{(n_{\ell_{k-1}}+1, n_i+1]}(\xi_{k-1}, r{+}1) \\
&= Q_{(n_{\ell_{k-1}}+1, n_i]}(\xi_{k-1}, r{+}1).
\end{align*}
Summing over all chains and internal vertices, the
first two terms
of~\eqref{eq:IBJ-S1-compat-component} therefore
combine into
\begin{align}
&q_{n_i+1}[B']_{ij}(r, s) - [B']_{ij}(r{+}1, s)
= \sum_{k=1}^{i-j}
  \sum_{j = \ell_0 < \ell_1 < \cdots < \ell_k = i}
  \sum_{\substack{\xi_p \le a_{\ell_p} \\ 1 \le p \le k-1}}
(-1)^k \nonumber \\
&\qquad\qquad
\Bigl(\prod_{p=0}^{k-2}
Q_{(n_{\ell_p}+1, n_{\ell_{p+1}}+1]}
(\xi_p, \xi_{p+1})\Bigr)
Q_{(n_{\ell_{k-1}}+1, n_i]}
(\xi_{k-1}, r{+}1),
\label{eq:IBJ-after-collapse}
\end{align}
where the endpoint is pinned at $r{+}1$, the last
edge carries unshifted upper index $n_i$, and the
remaining $k{-}1$ edges retain their shifted labels
$\mathbf{n}{+}\mathbf{1}$.

For a fixed chain
$j = \ell_0 < \ell_1 < \cdots < \ell_k = i$,
the sum over internal vertices
in~\eqref{eq:IBJ-after-collapse} is an operator
composition with cutoffs at internal layers:
\[
(-1)^k \bigl(
Q_{(n_j+1, n_{\ell_1}+1]}
\bar{\chi}_{a_{\ell_1}}
Q_{(n_{\ell_1}+1, n_{\ell_2}+1]}
\bar{\chi}_{a_{\ell_2}}
\cdots
\bar{\chi}_{a_{\ell_{k-1}}}
Q_{(n_{\ell_{k-1}}+1, n_i]}
\bigr)(s, r{+}1),
\]
where the first $k{-}1$ edges carry shifted labels
and the last edge is already collapsed.  Apply the
edge-shift identity~\eqref{eq:IBJ-edge-shift} to
each shifted factor: for $0 \le p \le k{-}2$,
\[
Q_{(n_{\ell_p}+1, n_{\ell_{p+1}}+1]}
= Q_{n_{\ell_p}+1}^{-1}
Q_{(n_{\ell_p}, n_{\ell_{p+1}}]}
Q_{n_{\ell_{p+1}}+1},
\]
and for the collapsed last edge,
$Q_{(n_{\ell_{k-1}}+1, n_i]}
= Q_{n_{\ell_{k-1}}+1}^{-1}
Q_{(n_{\ell_{k-1}}, n_i]}$.
After substitution, the operator product becomes
\begin{gather*}
Q_{n_j+1}^{-1}
Q_{(n_j, n_{\ell_1}]}
Q_{n_{\ell_1}+1}
\bar{\chi}_{a_{\ell_1}}
Q_{n_{\ell_1}+1}^{-1}
Q_{(n_{\ell_1}, n_{\ell_2}]}
Q_{n_{\ell_2}+1}
\bar{\chi}_{a_{\ell_2}}
Q_{n_{\ell_2}+1}^{-1} \\
\cdots
\bar{\chi}_{a_{\ell_{k-1}}}
Q_{n_{\ell_{k-1}}+1}^{-1}
Q_{(n_{\ell_{k-1}}, n_i]}.
\end{gather*}
At each internal vertex $\ell_p$
($1 \le p \le k{-}1$), the factor
$Q_{n_{\ell_p}+1}$ from edge $p{-}1$ and the
factor $Q_{n_{\ell_p}+1}^{-1}$ from edge $p$
appear on either side of the cutoff
$\bar{\chi}_{a_{\ell_p}}$.  Define the vertex operator
$V_{\ell_p} \defeq Q_{n_{\ell_p}+1}
\bar{\chi}_{a_{\ell_p}}
Q_{n_{\ell_p}+1}^{-1}$.
The sum~\eqref{eq:IBJ-after-collapse} therefore
equals
\begin{align}
&\sum_{k=1}^{i-j}
  \sum_{j = \ell_0 < \ell_1 < \cdots < \ell_k = i}
(-1)^k \nonumber \\
&\qquad\bigl(
Q_{n_j+1}^{-1}
Q_{(n_j, n_{\ell_1}]}
V_{\ell_1}
Q_{(n_{\ell_1}, n_{\ell_2}]}
\cdots
V_{\ell_{k-1}}
Q_{(n_{\ell_{k-1}}, n_i]}
\bigr)(s, r{+}1),
\label{eq:IBJ-edgeshift-sum}
\end{align}
where operator products denote kernel compositions
evaluated at $(s, r{+}1)$.

Since $\bar{\chi}_{a_{\ell_p}} Q_{n_{\ell_p}+1}^{-1}
= Q_{n_{\ell_p}+1}^{-1} \bar{\chi}_{a_{\ell_p}}
+ [\bar{\chi}_{a_{\ell_p}}, Q_{n_{\ell_p}+1}^{-1}]$,
the vertex operator decomposes as
$V_{\ell_p}
= \bar{\chi}_{a_{\ell_p}}
+ Q_{n_{\ell_p}+1}
[\bar{\chi}_{a_{\ell_p}}, Q_{n_{\ell_p}+1}^{-1}]$.
Computing the commutator from the definition of
$Q_k^{-1}$ gives
\[
[\bar{\chi}_{a_{\ell_p}}, Q_{n_{\ell_p}+1}^{-1}](x, y)
= q_{n_{\ell_p}+1}\delta_{x, a_{\ell_p}}
\delta_{y, a_{\ell_p}+1},
\]
so
\[
V_{\ell_p}
= \bar{\chi}_{a_{\ell_p}}
+ \mathcal{E}_{\ell_p},
\qquad
\mathcal{E}_{\ell_p}(x, y)
= q_{n_{\ell_p}+1}Q_{n_{\ell_p}+1}(x, a_{\ell_p})
\delta_{y, a_{\ell_p}+1},
\]
where $\mathcal{E}_{\ell_p}$ is rank one.

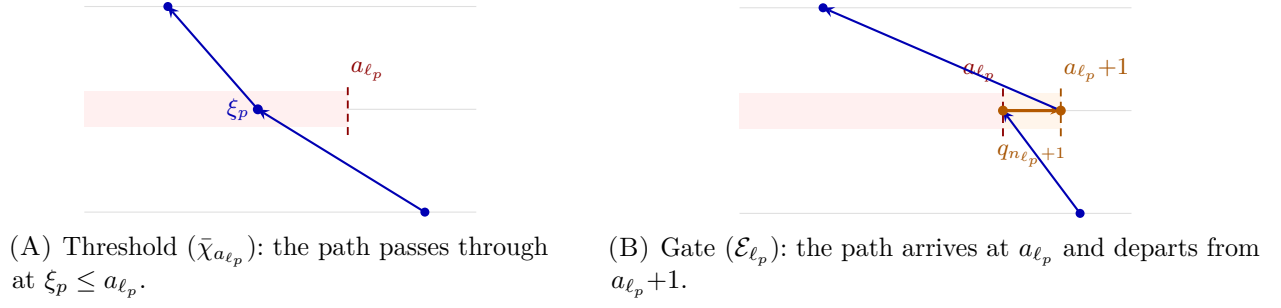
\begin{figure}[!htb]
\centering
\begin{subfigure}[b]{0.43\textwidth}
\centering
\begin{tikzpicture}[>=stealth, scale=0.85,
  every node/.style={font=\footnotesize}]

  \foreach \y in {0, 1.6, 3.2} {
    \draw[gray!25] (-0.3, \y) -- (5.8, \y);
  }

  \fill[red!6] (-0.3, 1.32) rectangle (3.8, 1.88);

  \draw[red!55!black, densely dashed,
    line width=0.7pt]
    (3.8, 1.2) -- (3.8, 2.0);
  \node[red!55!black, above right, inner sep=1pt]
    at (3.8, 2.0) {$a_{\ell_p}$};

  \coordinate (A) at (5.0, 0);
  \coordinate (B) at (2.4, 1.6);
  \coordinate (C) at (1.0, 3.2);

  \draw[blue!70!black, thick, -stealth]
    (A) -- (B);
  \draw[blue!70!black, thick, -stealth]
    (B) -- (C);

  \fill[blue!70!black] (A) circle (2pt);
  \fill[blue!70!black] (B) circle (2.2pt);
  \fill[blue!70!black] (C) circle (2pt);

  \node[left, inner sep=3pt, blue!70!black]
    at (B) {$\xi_p$};

\end{tikzpicture}
\caption{Threshold ($\bar{\chi}_{a_{\ell_p}}$): the
path passes through at
$\xi_p \le a_{\ell_p}$.}
\end{subfigure}
\hfill
\begin{subfigure}[b]{0.52\textwidth}
\centering
\begin{tikzpicture}[>=stealth, scale=0.85,
  every node/.style={font=\footnotesize}]

  \foreach \y in {0, 1.6, 3.2} {
    \draw[gray!25] (-0.3, \y) -- (5.8, \y);
  }

  \fill[red!6] (-0.3, 1.32) rectangle (3.8, 1.88);

  \fill[orange!8] (3.8, 1.32) rectangle (4.7, 1.88);

  \draw[red!55!black, densely dashed,
    line width=0.7pt]
    (3.8, 1.2) -- (3.8, 2.0);
  \node[red!55!black, above left, inner sep=1pt]
    at (3.75, 2.0) {$a_{\ell_p}$};

  \draw[orange!65!black, densely dashed,
    line width=0.7pt]
    (4.7, 1.2) -- (4.7, 2.0);
  \node[orange!65!black, above right, inner sep=1pt]
    at (4.7, 2.0) {$a_{\ell_p}{+}1$};

  \coordinate (A) at (5.0, 0);
  \coordinate (Ba) at (3.8, 1.6);
  \coordinate (Bd) at (4.7, 1.6);
  \coordinate (C) at (1.0, 3.2);

  \draw[blue!70!black, thick, -stealth]
    (A) -- (Ba);
  \draw[orange!70!black, very thick,
    -{Stealth[length=4pt]}]
    (Ba) -- (Bd);
  \draw[blue!70!black, thick, -stealth]
    (Bd) -- (C);

  \fill[blue!70!black] (A) circle (2pt);
  \fill[orange!70!black] (Ba) circle (2.2pt);
  \fill[orange!70!black] (Bd) circle (2.2pt);
  \fill[blue!70!black] (C) circle (2pt);

  \node[orange!65!black, below, inner sep=3pt]
    at (4.25, 1.2) {$q_{n_{\ell_p}+1}$};

\end{tikzpicture}
\caption{Gate ($\mathcal{E}_{\ell_p}$):
the path arrives at $a_{\ell_p}$ and departs from
$a_{\ell_p}{+}1$.}
\end{subfigure}
\caption{The two alternatives at an internal layer
$\ell_p$ after edge-shifting.  The vertex operator
$V_{\ell_p}
= \bar{\chi}_{a_{\ell_p}}
+ \mathcal{E}_{\ell_p}$
either passes the path through below the threshold
(left) or routes it through a gate at the boundary
(right), pinning the path to $a_{\ell_p}$ on entry
and $a_{\ell_p}{+}1$ on exit with coefficient
$q_{n_{\ell_p}+1}$.}
\label{fig:gate-mechanism}
\end{figure}

Expanding each
$V_{\ell_p} = \bar{\chi}_{a_{\ell_p}}
+ \mathcal{E}_{\ell_p}$
in~\eqref{eq:IBJ-edgeshift-sum} produces a sum over
all choices of either $\bar{\chi}_{a_{\ell_p}}$ or
$\mathcal{E}_{\ell_p}$ at each internal layer
(Figure~\ref{fig:gate-mechanism}).

\noindent\textit{Step 2: No-gate sector.}
The no-gate sector of~\eqref{eq:IBJ-edgeshift-sum},
where every internal layer chooses
$\bar{\chi}_{a_{\ell_p}}$, contributes
\begin{gather*}
\sum_{k=1}^{i-j}
  \sum_{j = \ell_0 < \ell_1 < \cdots < \ell_k = i}
(-1)^k \\
\bigl(
Q_{n_j+1}^{-1}
Q_{(n_j, n_{\ell_1}]}
\bar{\chi}_{a_{\ell_1}}
Q_{(n_{\ell_1}, n_{\ell_2}]}
\cdots
\bar{\chi}_{a_{\ell_{k-1}}}
Q_{(n_{\ell_{k-1}}, n_i]}
\bigr)(s, r{+}1).
\end{gather*}
The operator product after $Q_{n_j+1}^{-1}$ is the
unshifted Neumann
series~\eqref{eq:IBJ-propagator} defining
$[B]_{ij}$, so the no-gate sector equals
$(Q_{n_j+1}^{-1} \circ [B]_{ij})(s, r{+}1)$.
Evaluating
$Q_{n_j+1}^{-1}(s, \alpha)
= q_{n_j+1}\delta_{\alpha, s+1}
- \delta_{\alpha, s}$ gives
\[
q_{n_j+1}[B]_{ij}(r{+}1, s{+}1)
- [B]_{ij}(r{+}1, s),
\]
which is the negative of the third and fourth terms
of~\eqref{eq:IBJ-S1-compat-component}.
The left-hand side
of~\eqref{eq:IBJ-S1-compat-component} therefore
reduces to the at-least-one-gate sector
of~\eqref{eq:IBJ-edgeshift-sum}.

\noindent\textit{Step 3: Last-gate factorization.}
For paths with at least one gate crossing,
consider the contribution from paths whose last gate
is at layer $\mu$ with $j < \mu < i$
(Figure~\ref{fig:last-gate}).  The gate operator
$\mathcal{E}_\mu(x, y)
= q_{n_\mu+1}Q_{n_\mu+1}(x, a_\mu)
\delta_{y, a_\mu+1}$
fixes the spatial coordinate at layer $\mu$: the
path arrives at $a_\mu$ and departs from
$a_\mu{+}1$.  This splits the operator product into
two independent pieces.

\begin{figure}[!htb]
\centering
\begin{tikzpicture}[>=stealth, scale=1.05,
  every node/.style={font=\footnotesize}]

  \def\ls{1.2}
  \def\xmin{-1.2}
  \def\xmax{9.2}

  \foreach \y in {0,1,2,3,4,5} {
    \draw[gray!25] (\xmin, \y*\ls) -- (\xmax, \y*\ls);
  }

  \node[left, gray!70] at (\xmin, 0)
    {layer $j$};
  \node[left, gray!70] at (\xmin, \ls)
    {layer $\ell_1$};
  \node[left, gray!70] at (\xmin, 2*\ls)
    {layer $\ell_2$};
  \node[left, gray!70] at (\xmin, 3*\ls)
    {layer $\mu$};
  \node[left, gray!70] at (\xmin, 4*\ls)
    {layer $\ell_4$};
  \node[left, gray!70] at (\xmin, 5*\ls)
    {layer $i$};

  \fill[red!6] (\xmin, \ls-0.28)
    rectangle (4.5, \ls+0.28);
  \fill[red!6] (\xmin, 2*\ls-0.28)
    rectangle (6.0, 2*\ls+0.28);
  \fill[red!6] (\xmin, 3*\ls-0.28)
    rectangle (5.0, 3*\ls+0.28);
  \fill[red!6] (\xmin, 4*\ls-0.28)
    rectangle (3.5, 4*\ls+0.28);

  \fill[orange!8] (4.5, \ls-0.28)
    rectangle (5.3, \ls+0.28);
  \fill[orange!8] (5.0, 3*\ls-0.28)
    rectangle (5.8, 3*\ls+0.28);

  \draw[red!55!black, densely dashed,
    line width=0.7pt]
    (4.5, \ls-0.34) -- (4.5, \ls+0.34);
  \draw[red!55!black, densely dashed,
    line width=0.7pt]
    (6.0, 2*\ls-0.34) -- (6.0, 2*\ls+0.34);
  \draw[red!55!black, densely dashed,
    line width=0.7pt]
    (5.0, 3*\ls-0.34) -- (5.0, 3*\ls+0.34);
  \draw[red!55!black, densely dashed,
    line width=0.7pt]
    (3.5, 4*\ls-0.34) -- (3.5, 4*\ls+0.34);

  \draw[orange!65!black, densely dashed,
    line width=0.7pt]
    (5.3, \ls-0.34) -- (5.3, \ls+0.34);
  \draw[orange!65!black, densely dashed,
    line width=0.7pt]
    (5.8, 3*\ls-0.34) -- (5.8, 3*\ls+0.34);

  \node[red!55!black, above, inner sep=1pt]
    at (4.5, \ls+0.34) {$a_{\ell_1}$};
  \node[orange!65!black, above, inner sep=1pt]
    at (5.3, \ls+0.34) {$a_{\ell_1}{+}1$};

  \node[red!55!black, above, inner sep=1pt]
    at (6.0, 2*\ls+0.34) {$a_{\ell_2}$};

  \node[red!55!black, above, inner sep=1pt]
    at (5.0, 3*\ls+0.34) {$a_\mu$};
  \node[orange!65!black, above, inner sep=1pt]
    at (5.8, 3*\ls+0.34) {$a_\mu{+}1$};

  \node[red!55!black, above, inner sep=1pt]
    at (3.5, 4*\ls+0.34) {$a_{\ell_4}$};

  \coordinate (P0) at (7.5, 0);
  \coordinate (P1a) at (4.5, \ls);
  \coordinate (P1d) at (5.3, \ls);
  \coordinate (P2) at (3.5, 2*\ls);
  \coordinate (P3a) at (5.0, 3*\ls);
  \coordinate (P3d) at (5.8, 3*\ls);
  \coordinate (P4) at (2.5, 4*\ls);
  \coordinate (P5) at (1.2, 5*\ls);

  \draw[blue!70!black, thick, -stealth]
    (P0) -- (P1a);
  \draw[orange!70!black, very thick,
    -{Stealth[length=4pt]}]
    (P1a) -- (P1d);
  \draw[blue!70!black, thick, -stealth]
    (P1d) -- (P2);
  \draw[blue!70!black, thick, -stealth]
    (P2) -- (P3a);

  \draw[teal!80!black, thick, -stealth]
    (P3d) -- (P4);
  \draw[teal!80!black, thick, -stealth]
    (P4) -- (P5);

  \fill[blue!70!black] (P0) circle (2.2pt);
  \fill[teal!80!black] (P5) circle (2.2pt);

  \fill[orange!70!black] (P1a) circle (2.2pt);
  \fill[orange!70!black] (P1d) circle (2.2pt);

  \fill[blue!70!black] (P2) circle (2.2pt);
  \fill[teal!80!black] (P4) circle (2.2pt);

  \node[star, star points=5, star point ratio=2.3,
    fill=yellow!80!orange, draw=orange!70!black,
    line width=0.4pt,
    inner sep=1.6pt] at (P3a) {};
  \draw[orange!70!black, very thick,
    -{Stealth[length=4pt]}]
    (P3a) -- (P3d);
  \fill[teal!80!black] (P3d) circle (2.2pt);

  \node[below right, inner sep=2pt] at (P0)
    {$(j, s)$};
  \node[above left, inner sep=2pt] at (P5)
    {$(i, r{+}1)$};
  \node[left, inner sep=3pt,
    blue!70!black] at (P2)
    {$\xi_2$};
  \node[left, inner sep=3pt,
    teal!80!black] at (P4)
    {$\xi_4$};

  \node[right, inner sep=5pt,
    orange!70!black] at (P3d)
    {$(\mu, a_\mu{+}1)$};

  \node[blue!70!black, anchor=west]
    at (\xmax+0.1, 0.8*\ls)
    {\small prefix: $[B']_{\mu j}$};
  \node[teal!80!black, anchor=west]
    at (\xmax+0.1, 4.5*\ls)
    {\small suffix: $[B]_{i\mu}$};

\end{tikzpicture}
\caption{The last-gate decomposition for
$\mathcal{S}_1$-compatibility.  After endpoint
collapse and edge-shift, each internal layer
carries a vertex operator with two branches: a
threshold (standard cutoff, as at $\ell_2$ and
$\ell_4$) or a gate (orange arrows, as at $\ell_1$
and~$\mu$).  The path splits at its \emph{last}
gate crossing, at layer~$\mu$ (gold star), into a
\emph{prefix} from $(j, s)$ to $(\mu, a_\mu)$
with vertex operators at internal layers
(counted by~$[B']_{\mu j}$) and a \emph{suffix}
from $(\mu, a_\mu{+}1)$ to $(i, r{+}1)$ with
standard cutoffs (counted by~$[B]_{i\mu}$).}
\label{fig:last-gate}
\end{figure}
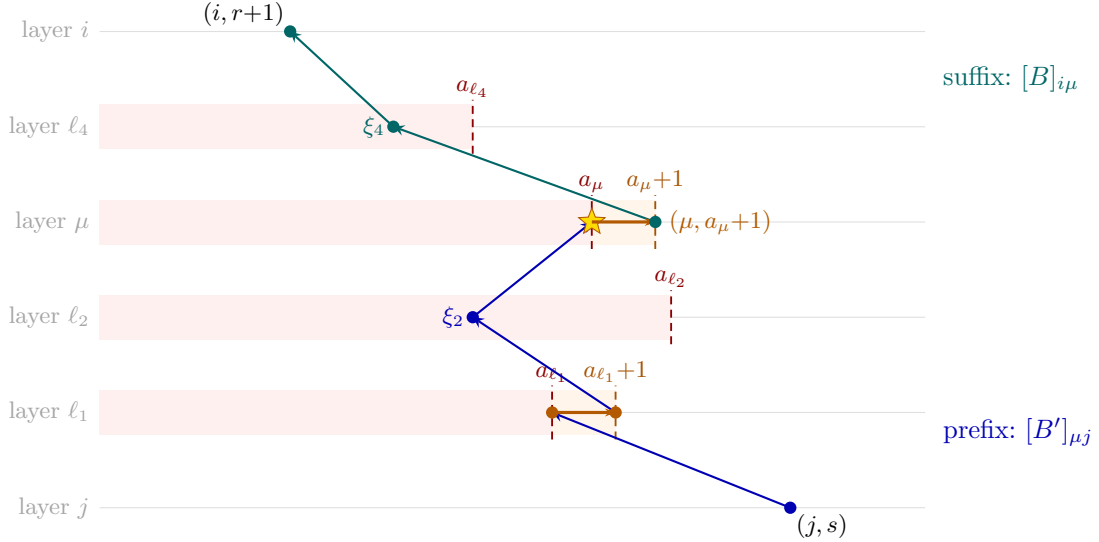

The \emph{suffix}, from layer $\mu$ to layer $i$,
has all internal vertices between $\mu$ and $i$
choosing $\bar{\chi}$ (since $\mu$ is the last gate).
The resulting operator product from $a_\mu{+}1$ to
$r{+}1$ with cutoffs at internal vertices is the
unshifted Neumann
series~\eqref{eq:IBJ-propagator}, giving
$[B]_{i\mu}(r{+}1, a_\mu{+}1)$.

For the \emph{prefix}, from layer $j$ to layer
$\mu$, let $\ell_1, \dotsc, \ell_{g-1}$ denote the
internal layers between $j$ and $\mu$ (so the prefix
has $g$ edges).  These internal vertices retain full
vertex operators $V$, the initial point carries
$Q_{n_j+1}^{-1}$, and the gate contributes
$q_{n_\mu+1}Q_{n_\mu+1}(\cdot, a_\mu)$ at the
endpoint.  Expanding each
$V_{\ell_p}
= Q_{n_{\ell_p}+1}
\bar{\chi}_{a_{\ell_p}}
Q_{n_{\ell_p}+1}^{-1}$
and grouping adjacent factors by the edge-shift
identity~\eqref{eq:IBJ-edge-shift} reconstitutes the
shifted transition kernels:
\begin{align}
Q_{n_j+1}^{-1}
Q_{(n_j, n_{\ell_1}]}
&V_{\ell_1}
Q_{(n_{\ell_1}, n_{\ell_2}]}
\cdots
V_{\ell_{g-1}}
Q_{(n_{\ell_{g-1}}, n_\mu]}
Q_{n_\mu+1}
\nonumber \\
&= Q_{(n_j+1, n_{\ell_1}+1]}
\bar{\chi}_{a_{\ell_1}}
Q_{(n_{\ell_1}+1, n_{\ell_2}+1]}
\cdots
\bar{\chi}_{a_{\ell_{g-1}}}
Q_{(n_{\ell_{g-1}}+1, n_\mu+1]},
\label{eq:IBJ-prefix-edgeshift}
\end{align}
which is the shifted Neumann series defining
$[B']_{\mu j}(a_\mu, s)$.  Summing over all prefix
chains gives
$q_{n_\mu+1}[B']_{\mu j}(a_\mu, s)$.

The prefix and suffix layer sequences range
independently, and the sign factors $(-1)^g$ from
the prefix and $(-1)^{k-g}$ from the suffix combine
to $(-1)^k$.  The contribution for fixed $\mu$
therefore factors as
\[
[B]_{i\mu}(r{+}1, a_\mu{+}1)
q_{n_\mu+1}
[B']_{\mu j}(a_\mu, s).
\]
Summing over $j < \mu < i$ gives the right-hand side
of~\eqref{eq:IBJ-S1-compat-component}.
\end{proof}

Regularity follows as in~\S\ref{sec:RBJ}: the
propagator entries decay geometrically with rate
$q_{\min} = \min(q_1, \dotsc, q_N) > 1$ replacing
the homogeneous parameter, $\psi$ has compact support
in~$x$ from its contour representation, and $\phi$
inherits geometric decay from
$\bar{\mathcal{S}}_{[1,n],(0,t]}$.

\paragraph{\textbf{Corrected edge weights.}}
By Lemma~\ref{lem:Lambda-diamond}, the corrected
edge weights on the product graph are
$\Lambda_k(u) = P(u)D_k(u)
P(\mathcal{S}_k u)^{-1}$.
For $k = 2$, the propagator is
$\mathcal{S}_2$-invariant
(Lemma~\ref{lem:IBJ-S2-inv}), so
$P(\mathcal{S}_2 u) = P(u)$ and
$\Lambda_2(u) = D_2(u) = -p_{t+1}I_m$.
For $k = 1$, the $\mathcal{S}_1$-shift changes
$\mathbf{n} \to \mathbf{n} + \mathbf{1}$ and
$\mathcal{S}_1$-invariance fails, so
$P(\mathcal{S}_1 u) \neq P(u)$ in general; since
$D_1(u) = I_m$, the corrected edge weight is
$\Lambda_1(u) = P(u)P(\mathcal{S}_1 u)^{-1}$.
Theorem~\ref{thm:product-graph} requires additionally
that the resolvent $(I - zK(u))^{-1}$ exist.  At $z = 1$,
the Fredholm determinant
$F_{t, \mathbf{a}, \mathbf{n}} \defeq
\det(I - K_{t, \mathbf{a}, \mathbf{n}})$ equals the
multipoint gap probability
$\PP_y\bigl(\bigcap_i
\{Y_{n_i}(t) > a_i\}\bigr)$.  Since $K$ is trace class,
$F \ne 0$ if and only if the resolvent exists.  At each
time step, each particle attempts to jump one site to the
right, so $Y_{n_i}(t) \le y_{n_i} + t$.  If every jump
attempt for particles $1, \dotsc, n_m$ succeeds (an event
of positive probability), the strict integer ordering
prevents blocking and $Y_k(t) = y_k + t$ for
$k \le n_m$.  The resolvent at $z = 1$ therefore exists
precisely on
$\mathcal{R} \defeq \{(t, \mathbf{a}, \mathbf{n}) \in
\mathcal{V}^m : a_i < y_{n_i} + t
\text{ for all } i\}$.
The mixed diamond at a vertex~$u$ requires the resolvent
at the four vertices $u$, $\mathcal{S}_1 u$,
$\mathcal{S}_2 u$, $\mathcal{S}_1\mathcal{S}_2 u$
entering~\eqref{eq:dressed-C-explicit}, together with
the $\mathcal{T}$-shifts $\mathcal{T}u$,
$\mathcal{T}\mathcal{S}_1 u$,
$\mathcal{T}\mathcal{S}_2 u$ from the proof
of~\eqref{eq:M-inverse}.\footnote{The explicit
formula~\eqref{eq:dressed-C-explicit} involves the
resolvent only at the four $\mathcal{S}$-shifted
vertices; a direct verification of the mixed diamond
from this formula would yield the larger domain
$a_i < y_{n_i+1} + t$.}
Since $\mathcal{S}_2$ preserves~$\mathcal{R}$ and
$\mathcal{T}$ tightens each threshold by~$1$, the
binding constraint is
$\mathcal{T}\mathcal{S}_1 u \in \mathcal{R}$, i.e.\
$a_i < y_{n_i+1} + t - 1$ for all~$i$.

\subsection{Multipoint equation}
\label{sec:IBJ-multipoint}

The dressed observable at $z = 1$,
\[
\mathcal{M}(u)
\defeq I + \Phi(u)(I - K(u))^{-1}\Psi(u)
\in \End(\R^m),
\]
has dressed edge weights
\begin{equation}\label{eq:IBJ-dressed-Ck}
\mathcal{M}_k(u)
= \mathcal{M}(\mathcal{T}u)^{-1}C_k(u)
  \mathcal{M}(\mathcal{S}_k u).
\end{equation}
For readability, we write
$\mathcal{M}_{t, \mathbf{a}, \mathbf{n}} \defeq
\mathcal{M}(t, \mathbf{a}, \mathbf{n})$.

\begin{theorem}
\label{thm:IBJ-multipoint}
The dressed observable satisfies the mixed diamond
equation
\begin{align}\label{eq:IBJ-multipoint}
&\bigl(\Lambda_{\mathbf{a},
\mathbf{n}}
\mathcal{M}_{t, \mathbf{a}+\mathbf{1},
\mathbf{n}+\mathbf{1}}^{-1}
+ p_{t+1}
\mathcal{M}_{t+1, \mathbf{a}+\mathbf{2},
\mathbf{n}}^{-1}
q_{\mathbf{n}+\mathbf{1}}\bigr)
\mathcal{M}_{t+1, \mathbf{a}+\mathbf{1},
\mathbf{n}+\mathbf{1}}
\nonumber \\
&{}-
\mathcal{M}_{t, \mathbf{a}+\mathbf{1},
\mathbf{n}}^{-1}
\bigl(p_{t+1}q_{\mathbf{n}+\mathbf{1}}
\mathcal{M}_{t, \mathbf{a},
\mathbf{n}+\mathbf{1}}
+ \mathcal{M}_{t+1, \mathbf{a}+\mathbf{1},
\mathbf{n}}
\Lambda_{\mathbf{a},
\mathbf{n}}\bigr) = 0,
\end{align}
at every $(t, \mathbf{a}, \mathbf{n})$ with
$a_i < y_{n_i+1} + t - 1$ for all~$i$, where
\[
q_{\mathbf{n}+\mathbf{1}}
\defeq \mathrm{Diag}(q_{n_1+1}, \dotsc, q_{n_m+1}), 
\Lambda_{\mathbf{a}, \mathbf{n}}
\defeq (I + B_{\mathbf{a}, \mathbf{n}}(\mathbf{a},
\mathbf{a}))
(I + B_{\mathbf{a}, \mathbf{n}+\mathbf{1}}(\mathbf{a},
\mathbf{a}))^{-1}.
\]
\end{theorem}

\begin{remark}
By Lemma~\ref{lem:Lambda-diamond},
$(I + B_{\mathbf{a}, \mathbf{n}+\mathbf{1}}
(\mathbf{a}, \mathbf{a}))^{-1}
= I - B_{\mathbf{a}, \mathbf{n}+\mathbf{1}}
(\mathbf{a}+\mathbf{1}, \mathbf{a}+\mathbf{1})$.
\end{remark}

\begin{proof}
By Theorem~\ref{thm:product-graph},
$(\mathcal{M}_k, \Lambda_k)$ satisfies the diamond
equations on~$\mathcal{R}$.  Since
$\Lambda_2 = -p_{t+1}I_m$ is scalar and
$\Lambda_1 = \Lambda_{\mathbf{a}, \mathbf{n}}$
is both $t$-independent and $\mathcal{T}$-covariant,
the mixed diamond
equation~\eqref{eq:diamond-mixed-dressed} reduces to
\begin{equation}\label{eq:IBJ-mixed-simplified}
-p_{t+1}\bigl[\mathcal{M}_1(u)
- \mathcal{M}_1(\mathcal{S}_2 u)\bigr]
+ \Lambda_{\mathbf{a}, \mathbf{n}}
\mathcal{M}_2(\mathcal{S}_1 u)
- \mathcal{M}_2(u)
\Lambda_{\mathbf{a}, \mathbf{n}} = 0.
\end{equation}
Substituting~\eqref{eq:IBJ-dressed-Ck} and the
explicit shifts gives
\begin{align}
&{-}p_{t+1}
\mathcal{M}_{t, \mathbf{a}+\mathbf{1},
\mathbf{n}}^{-1}
q_{\mathbf{n}+\mathbf{1}}
\mathcal{M}_{t, \mathbf{a},
\mathbf{n}+\mathbf{1}}
+ p_{t+1}
\mathcal{M}_{t+1, \mathbf{a}+\mathbf{2},
\mathbf{n}}^{-1}
q_{\mathbf{n}+\mathbf{1}}
\mathcal{M}_{t+1, \mathbf{a}+\mathbf{1},
\mathbf{n}+\mathbf{1}}
\nonumber \\[3pt]
&\quad {}+ \Lambda_{\mathbf{a}, \mathbf{n}}
\mathcal{M}_{t, \mathbf{a}+\mathbf{1},
\mathbf{n}+\mathbf{1}}^{-1}
\mathcal{M}_{t+1, \mathbf{a}+\mathbf{1},
\mathbf{n}+\mathbf{1}}-
\mathcal{M}_{t, \mathbf{a}+\mathbf{1},
\mathbf{n}}^{-1}
\mathcal{M}_{t+1, \mathbf{a}+\mathbf{1},
\mathbf{n}}
\Lambda_{\mathbf{a}, \mathbf{n}} = 0.
\label{eq:IBJ-mixed-expanded}
\end{align}
The second and third terms share the right factor
$\mathcal{M}_{t+1, \mathbf{a}+\mathbf{1},
\mathbf{n}+\mathbf{1}}$, and the first and fourth
share the left factor
$-\mathcal{M}_{t, \mathbf{a}+\mathbf{1},
\mathbf{n}}^{-1}$.  Factoring
yields~\eqref{eq:IBJ-multipoint}.
\end{proof}

\begin{corollary}\label{cor:IBJ-scalar}
The one-point Fredholm determinant
$F_{t, a, n} = \det(I - K_{t, a, n})$ satisfies
\begin{align}
&\frac{q_{n+1} p_{t+1}}{1 + q_{n+1} p_{t+1}}
F_{t, a, n+1}F_{t+1, a+2, n}
+ \frac{1}{1 + q_{n+1} p_{t+1}}
F_{t+1, a+1, n}F_{t, a+1, n+1} \nonumber \\
&\qquad - F_{t, a+1, n}F_{t+1, a+1, n+1} = 0,
\label{eq:IBJ-scalar-HM}
\end{align}
at every $(t, a, n)$ with $a < y_{n+1} + t - 1$.
\end{corollary}

\begin{proof}
Proposition~\ref{prop:scalar-HM} applies at $z = 1$:
the boundary condition $F(T^\ell v) \to 1$ as
$\ell \to -\infty$ holds because the thresholds
decrease to $-\infty$ and the kernel vanishes.
The coefficients are
$\alpha_{12} = c_1 \lambda_2 = -q_{n+1}p_{t+1}$
and
$\alpha_{21} = c_2 \lambda_1 = 1$.
Substituting the lattice points into
\eqref{eq:HM-variable} and computing
$\alpha_{12} - \alpha_{21}
= -(1 + q_{n+1} p_{t+1})$,
then dividing by $-(1 + q_{n+1} p_{t+1})$,
gives~\eqref{eq:IBJ-scalar-HM}.
\end{proof}

\section{Left Bernoulli Jumps}\label{sec:LBJ}

\paragraph{\textbf{System description.}}
The left Bernoulli jump model is a discrete-time totally
asymmetric simple exclusion process with left Bernoulli
jumps, pushing interaction, and sequential update.  The
system consists of $N$ particles on
$\Z$.  The configuration at time
$t \in \Z_{\ge 0}$ is denoted by
\[
Y(t) = (Y_1(t) > Y_2(t) > \cdots > Y_N(t)),
\]
where $Y_k(t)$ is the position of the $k$-th particle
from the right.

At each discrete time step $t \ge 1$, the particles are
updated sequentially in the order $k = 1, 2, \dots, N$,
so the rightmost particle is updated first.  During its
update, particle $k$ attempts to jump one step to the left
with probability $p \in (0,1)$, and stays put with
probability $q = 1 - p$.  If a particle lands on an
occupied site, it pushes that particle and its entire
cluster of nearest neighbours one step to the left,
preserving the strict ordering.  The evolution is given by
\[
Y_k(t) = \min\{Y_k(t{-}1) - \xi(t,k),
Y_{k-1}(t) - 1\},
\qquad Y_0(t) = \infty,
\]
where $\{\xi(t,k)\}_{t \ge 1, 1 \le k \le N}$ are
independent Bernoulli random variables with
$\PP(\xi(t,k) = 1) = p$, and the convention
$Y_0(t) = \infty$ means that the rightmost particle
faces no constraint from above.

\paragraph{\textbf{Fredholm determinant formula.}}
Fix a time $t \ge 0$ and an initial configuration
$Y(0) = y = (y_1 > y_2 > \cdots > y_N)$.  Fix
$m \in \{1, \dots, N{-}1\}$, particle labels
$\mathbf{n} = (n_1, n_2, \dots, n_m)$ with
$1 \le n_1 < n_2 < \cdots < n_m < N$, and spatial
thresholds $\mathbf{a} = (a_1, \dots, a_m) \in \Z^m$.
The multipoint joint cumulative distribution of the particle
positions is given by a Fredholm determinant on
$\ell^2(\{n_1, \dots, n_m\} \times \Z)$~\cite{MatetskiRemenik2023}:
\[
\PP_y\Bigl(\bigcap_{i=1}^m
\{Y_{n_i}(t) > a_i\}\Bigr)
= \det(I - \bar{\chi}_{\mathbf{a}} K_t
\bar{\chi}_{\mathbf{a}})_{\ell^2(\{n_1, \dots, n_m\}
\times \Z)},
\]
where $\bar{\chi}_{\mathbf{a}}(n_i, x) =
\mathbf{1}_{x \le a_i}$.  The correlation kernel $K_t$ has
block entries
\[
K_t(n_i, x_i; n_j, x_j)
= -Q^{n_j - n_i}(x_i, x_j)\mathbf{1}_{n_i < n_j}
+ \bigl(\mathcal{S}_{-t,-n_i}^*
\overline{\mathcal{S}}_{-t,n_j}^{\operatorname{epi}(y)}
\bigr)(x_i, x_j),
\]
where $A^*$ denotes the transpose kernel,
$A^*(x, y) = A(y, x)$.

Fix an auxiliary parameter $\theta \in (0, 1)$, which
enters the Fredholm determinant representation but does
not affect the final multipoint equation
\eqref{eq:H6-multipoint}, and let
$\alpha = (1{-}\theta)\theta^{-1}$.  The operators
defining the kernel are as follows.

\medskip
\noindent\textit{The transition matrix $Q$.}
The matrix $Q$ is the transition matrix of a random walk on
$\Z$ taking $\mathrm{Geom}[1{-}\theta]$ steps
strictly to the left:
\[
Q(z_1, z_2) = (1{-}\theta)\theta^{z_1 - z_2 - 1}
\mathbf{1}_{z_1 > z_2},
\]
with its $n$-th power admitting the integral
representation
\[
Q^n(z_1, z_2) = \frac{\alpha^n}{2\pi\mathrm{i}}
\oint_{\gamma_\rho} \diff w
\frac{\theta^{z_1 - z_2}}{w^{z_1 - z_2 - n + 1}}
\Bigl(\frac{1}{1 - w}\Bigr)^n.
\]

\medskip
\noindent\textit{The scattering operators
$\mathcal{S}_{-t,-n}$ and
$\overline{\mathcal{S}}_{-t,n}$.}\footnote{The kernel
operators $\mathcal{S}_{-t,-n}$ and
$\overline{\mathcal{S}}_{-t,n}$ are unrelated to the
lattice shifts $\mathcal{S}_k$ of the diamond framework
(\S\ref{sec:linear-problem}); the notation follows
Matetski--Remenik~\cite{MatetskiRemenik2023}.}
These operators are built from the probability generating
function $\varphi(w) = q + p/w$ of the left Bernoulli jump
distribution.  The contour integral representations are
\[
\mathcal{S}_{-t,-n}(z_1, z_2)
= \frac{\alpha^{-n+1}}{2\pi\mathrm{i}}
\oint_{\gamma_\rho} \frac{\diff w}{w}
\frac{\theta^{z_2 - z_1}}{w^{z_2 - z_1 + n}}
(1{-}w)^n\Bigl(q + \frac{p}{w}\Bigr)^{\!t},
\]
\[
\overline{\mathcal{S}}_{-t,n}(z_1, z_2)
= \frac{\alpha^{n-1}}{2\pi\mathrm{i}}
\oint_{\gamma_\delta} \frac{\diff w}{w}
\frac{(1{-}w)^{z_2 - z_1 + n - 1}}
{\theta^{z_2 - z_1} w^{n-1}}
\Bigl(q + \frac{p}{1{-}w}\Bigr)^{\!-t},
\]
where $\gamma_\rho$ is a positively oriented circle of
radius $\rho \in (0, 1)$ around the origin, and
\(\gamma_\delta\) is a sufficiently small positively
oriented circle around the origin, chosen so that the only
singularity enclosed is the pole at \(w=0\).

\medskip
\noindent\textit{The epigraph operator
$\overline{\mathcal{S}}_{-t,n}^{\operatorname{epi}(y)}$.}
This is the hitting probability operator defined in terms
of the initial data:
\[
\overline{\mathcal{S}}_{-t,n}^{\operatorname{epi}(y)}
(z_1, z_2)
= \E_{W_0 = z_1}\bigl[
\overline{\mathcal{S}}_{-t,n-\tau}(W_\tau, z_2)
\mathbf{1}_{\tau < n}\bigr],
\]
where $(W_\ell)_{\ell \ge 0}$ is the random walk with
transition matrix $Q$, and
$\tau = \min\{\ell \in \{0, \dots, N{-}1\} :
W_\ell > y_{\ell+1}\}$
is the hitting time of the strict epigraph of the initial
data by the random walk, with the convention
$\tau = \infty$ if the set is empty.

\paragraph{\textbf{References.}}
The left Bernoulli jump model belongs to the family of four
basic discrete-time TASEP variants whose determinantal
transition kernels were obtained by
Dieker and Warren~\cite{DiekerWarren2008}.  The Fredholm determinant
formula used below is the Matetski--Remenik formula
(Theorem~1.2, equations~(1.3)--(1.10)
of~\cite{MatetskiRemenik2023}) specialized to
$\kappa = 0$ and $\varphi(w) = q + p/w$.

\subsection{Kernel reformulation}\label{sec:LBJ-kernel-reform}

We rewrite the extended-kernel Fredholm determinant in the
form of the product graph construction of \S\ref{sec:product-graph}.  The
functions $\psi, \phi$ below provide the seed data, while the
lower-triangular operator $B_{\mathbf{a},\mathbf{n}}$ will be verified
to satisfy the admissibility conditions of Definition~\ref{def:admissible-prop}.

\paragraph{\textbf{Seed data.}}
Define
\[
\psi_{t,r,n}(x) \defeq \mathcal{S}_{-t,-n}^*(r,x)
= \mathcal{S}_{-t,-n}(x,r),
\qquad
\phi_{t,r,n}(x) \defeq
\overline{\mathcal{S}}_{-t,n}^{\operatorname{epi}(y)}(x,r),
\]
where we suppress the dependence of $\phi$ on the initial
condition $y$.

\paragraph{\textbf{Propagator.}}
The propagator is an instance of the directed-path
propagator (Definition~\ref{def:directed-path-propagator})
with transition kernels
$Q_{\ell,\ell'} = Q^{n_{\ell'} - n_\ell}$ for
$\ell < \ell'$.  On the cutoff region $r \le a_i$,
$r' \le a_j$, the entries of
$I + B_{\mathbf{a},\mathbf{n}}$ coincide with those of
$((I + \bar{\chi}_{\mathbf{a}} L
\bar{\chi}_{\mathbf{a}})^{-1})^{\!\top}$, where $L$
is the strictly upper-triangular block operator with
entries $L_{ij}(x, x') = Q^{n_j - n_i}(x, x')
\mathbf{1}_{i < j}$ (Remark~\ref{rem:Neumann-series}).

The propagator depends on $\mathbf{a}$ and $\mathbf{n}$
but not on the time coordinate~$t$, which is accordingly
suppressed from the notation.

\begin{lemma}\label{lem:LBJ-kernel}
We have
\[
\PP_y\Bigl(\bigcap_{i=1}^m
\{Y_{n_i}(t) > a_i\}\Bigr)
= \det(I - K_{t, \mathbf{a}, \mathbf{n}})_{\ell^2(\Z)},
\]
where
\[
K_{t, \mathbf{a}, \mathbf{n}}(x, x')
= \sum_{1 \le j \le i \le m}
  \sum_{r \le a_i} \sum_{r' \le a_j}
   \psi_{t, r, n_i}(x)
  \bigl[\delta_{ij}\delta_{r,r'}
  + [B_{\mathbf{a},\mathbf{n}}]_{ij}(r, r')\bigr]
  \phi_{t, r', n_j}(x').
\]
This identifies $K_{t, \mathbf{a}, \mathbf{n}}$ with the
kernel~\eqref{eq:K-def} on the product graph
$\mathcal{G}^m$, with $\psi, \phi$ as the seed data and
$B_{\mathbf{a},\mathbf{n}}$ as the propagator.
\end{lemma}

\begin{proof}
The argument is completely analogous to that of
Lemma~\ref{lem:RBJ-kernel}: the block-triangular
decomposition of
$I - \bar{\chi}_{\mathbf{a}} K_t
\bar{\chi}_{\mathbf{a}}$,
Sylvester's identity, and the transpose convention
produce the claimed single-space kernel.
The random walk $Q$ and the hitting-time construction
are the same as in \S\ref{sec:RBJ}; the free scattering
operators entering the seed data are different.
The Schatten-class estimates are analogous to those of
Lemma~\ref{lem:RBJ-kernel}, with the same geometric
rate~$\theta$.  The contour formula for $\psi$ gives
compact support of width $n + t$ (the Laurent polynomial
$(1{-}w)^n(q{+}p/w)^t$ has powers ranging from $-t$
to~$n$), and the contour estimate for $\phi$ gives
$|\phi_{t,r,n}(x)| \le C\beta^{x-r}
\mathbf{1}_{x > y_N}$.  With these bounds, the same
conjugation and Sylvester reduction as in
Lemma~\ref{lem:RBJ-kernel} apply.
\end{proof}

\subsection{Base graph and seed data}
\label{sec:LBJ-base-graph}

The lattice structure governing the seed data under
shifts of the parameters $(t, a, n)$ is identified next,
together with a verification that $\psi, \phi$ satisfy the
linear problems of the diamond framework.

Let $\mathcal{V} \subseteq \Z^3$ be the lattice
generated by the shifts
\[
T = e^{\pa_a}, \qquad
S_1 = e^{\pa_n}, \qquad
S_2 = e^{\pa_t},
\]
so that for $u = (t, a, n) \in \mathcal{V}$,
\[
Tu = (t, a{+}1, n), \qquad
S_1 u = (t, a, n{+}1), \qquad
S_2 u = (t{+}1, a, n).
\]
Take $E = \mathbb{F} = \R$ and define the
constant scalar edge weights
\begin{equation}\label{eq:LBJ-edge-weights}
(c_1, c_2) = \bigl((1{-}\theta)^{-1},
p\theta^{-1}\bigr),
\qquad
(\lambda_1, \lambda_2)
= \bigl((1{-}\theta)^{-1}\theta, -q\bigr).
\end{equation}
Since $(c_k, \lambda_k)$ are constant scalars, the diamond
equations
\eqref{eq:diamond-C-ij}--\eqref{eq:diamond-mixed-ij} are
satisfied trivially: all three reduce to commutativity in
$\R$.

\begin{lemma}\label{lem:LBJ-seed-linear}
Let $H = \ell^2(\Z)$.  The seed functions
$\psi_{t,r,n}$ and $\phi_{t,r,n}$ defined in
\S\ref{sec:LBJ-kernel-reform} satisfy the linear problem
\eqref{eq:Linear-Psi-3} and the adjoint linear problem
\eqref{eq:Linear-Phi-Prob} on $\mathcal{V}$ with edge
weights $(c_k, \lambda_k)$.  Explicitly, for each fixed
$r \in \Z$:
\begin{align}
\psi_{t,r,n+1}
&= (1{-}\theta)^{-1}\psi_{t,r+1,n}
   - (1{-}\theta)^{-1}\theta\psi_{t,r,n},
\label{eq:LBJ-psi-S1} \\
\psi_{t+1,r,n}
&= p\theta^{-1}\psi_{t,r+1,n}
   + q\psi_{t,r,n},
\label{eq:LBJ-psi-S2}
\end{align}
and
\begin{align}
\phi_{t,r+1,n}
&= (1{-}\theta)^{-1}\phi_{t,r,n+1}
   - (1{-}\theta)^{-1}\theta\phi_{t,r+1,n+1},
\label{eq:LBJ-phi-S1} \\
\phi_{t,r+1,n}
&= p\theta^{-1}\phi_{t+1,r,n}
   + q\phi_{t+1,r+1,n}.
\label{eq:LBJ-phi-S2}
\end{align}
\end{lemma}

\begin{proof}
Each identity is verified in turn.

\medskip
\noindent\textit{Verification of
\eqref{eq:LBJ-psi-S1}.}
The $S_1$-shift involves only $n$ and the $Q$-power
structure, which is the same as in
\S\ref{sec:RBJ}; the left Bernoulli generating function
$\varphi(w) = q + p/w$ passes through unchanged.  The
computation is identical to the verification of
\eqref{eq:RBJ-psi-S1}: the parenthetical factor
$(1{-}\theta)^{-1}(\theta/w - \theta)
= (1{-}w)/(\alpha w)$ shifts the integrand from $n$
to $n{+}1$.

\medskip
\noindent\textit{Verification of
\eqref{eq:LBJ-psi-S2}.}
From the contour representation,
\[
\psi_{t,r,n}(x) = \frac{\alpha^{-n+1}}{2\pi\mathrm{i}}
  \oint_{\gamma_\rho} \frac{\diff w}{w}
  \frac{\theta^{r-x}}{w^{r-x+n}}
  (1{-}w)^n\Bigl(q + \frac{p}{w}\Bigr)^{\!t}.
\]
The right-hand side of \eqref{eq:LBJ-psi-S2} is
\begin{align*}
p\theta^{-1}\psi_{t,r+1,n}(x)
+ q\psi_{t,r,n}(x)
= \frac{\alpha^{-n+1}}{2\pi\mathrm{i}}
  \oint_{\gamma_\rho} \frac{\diff w}{w}
  \frac{\theta^{r-x}}{w^{r-x+n}}
  (1{-}w)^n\Bigl(q + \frac{p}{w}\Bigr)^{\!t}
  \Bigl(p\theta^{-1} \cdot \frac{\theta}{w}
  + q\Bigr).
\end{align*}
The parenthetical factor equals $q + p/w$, giving
\[
= \frac{\alpha^{-n+1}}{2\pi\mathrm{i}}
  \oint_{\gamma_\rho} \frac{\diff w}{w}
  \frac{\theta^{r-x}}{w^{r-x+n}}
  (1{-}w)^n\Bigl(q + \frac{p}{w}\Bigr)^{\!t+1}
= \psi_{t+1,r,n}(x).
\]

\medskip
\noindent\textit{Verification of
\eqref{eq:LBJ-phi-S1}.}
The $\overline{\mathcal{S}}$ recurrence in~$n$,
\begin{equation}\label{eq:LBJ-Sbar-recurrence-n}
(1{-}\theta)^{-1}\bigl[
\overline{\mathcal{S}}_{-t,n}(z, r)
- \theta\overline{\mathcal{S}}_{-t,n}(z, r{+}1)
\bigr]
= \overline{\mathcal{S}}_{-t,n-1}(z, r{+}1),
\end{equation}
depends only on the $w$-power and $(1{-}w)$-power
structure, not on the probability generating function,
so the verification is the same as for
\eqref{eq:Sbar-recurrence-n}.  The vanishing
$\overline{\mathcal{S}}_{-t,0} \equiv 0$ and the
lifting to the epigraph version follow by the same
argument as in \S\ref{sec:RBJ}.

\medskip
\noindent\textit{Verification of
\eqref{eq:LBJ-phi-S2}.}
Since neither side involves a shift in $n$, the
indicator $\mathbf{1}_{\tau < n}$ passes through
unchanged, and the verification reduces to a
recurrence for the free operator.  From the contour
representation,
\begin{equation}\label{eq:LBJ-Sbar-recurrence-t}
p\theta^{-1}
\overline{\mathcal{S}}_{-(t+1),n}(z, r)
+ q\overline{\mathcal{S}}_{-(t+1),n}(z, r{+}1)
= \overline{\mathcal{S}}_{-t,n}(z, r{+}1).
\end{equation}
Indeed, the left-hand side equals
\begin{align*}
&\frac{\alpha^{n-1}}{2\pi\mathrm{i}}
\oint_{\gamma_\delta} \frac{\diff w}{w}
\frac{(q{+}p/(1{-}w))^{-(t+1)}
      (1{-}w)^{r-z+n-1}}
     {\theta^{r-z} w^{n-1}}
\Bigl(p\theta^{-1}
+ q\frac{1{-}w}{\theta}\Bigr) \\
&\qquad = \frac{\alpha^{n-1}}{2\pi\mathrm{i}}
\oint_{\gamma_\delta} \frac{\diff w}{w}
\frac{(q{+}p/(1{-}w))^{-(t+1)}
      (1{-}w)^{r-z+n-1}}
     {\theta^{r+1-z} w^{n-1}}
\bigl(p + q(1{-}w)\bigr).
\end{align*}
Since
$p + q(1{-}w) = 1 - qw
= (q + p/(1{-}w))(1{-}w)$,
the factor $(q + p/(1{-}w))^{-(t+1)}$ is promoted to
$(q + p/(1{-}w))^{-t}$ and the power of $(1{-}w)$ shifts
by $+1$, giving
$\overline{\mathcal{S}}_{-t,n}(z, r{+}1)$.
Taking the expectation against the hitting-time
representation completes the proof.
\end{proof}

\subsection{Admissible propagators and corrected edge weights}
\label{sec:LBJ-propagators}

The product graph $\mathcal{G}^m$ has vertex set
$\mathcal{V}^m$ and is generated by the diagonal shifts
\[
\mathcal{T}u = (t, \mathbf{a} + \mathbf{1}, \mathbf{n}),
\qquad
\mathcal{S}_1 u
= (t, \mathbf{a}, \mathbf{n} + \mathbf{1}),
\qquad
\mathcal{S}_2 u
= (t{+}1, \mathbf{a}, \mathbf{n}),
\]
where $u = (t, \mathbf{a}, \mathbf{n}) \in \mathcal{V}^m$
and $\mathbf{1} = (1, \dots, 1) \in \Z^m$.  The
block-diagonal edge weights are
$C_k(u) = c_k I_m$ and $D_k(u) = \lambda_k I_m$
(the uncorrected edge weights of
Definition~\ref{def:admissible-prop}).

The propagator is an instance of
the directed-path propagator
(Definition~\ref{def:directed-path-propagator}) with
transition kernels $Q_{\ell,\ell'} = Q^{n_{\ell'} - n_\ell}$
for $\ell < \ell'$.  The standing hypotheses of
\S\ref{sec:directed-path} are satisfied: translation
invariance of $Q^n$ is inherited from the translation
invariance of~$Q$, and $\mathcal{T}$-invariance holds
because the transition kernels depend only on the
particle-label differences, which are unchanged
by~$\mathcal{T}$.

By Proposition~\ref{prop:constant-scalar},
$\mathcal{S}_k$-compatibility follows from the
constant-scalar edge weights once the propagator satisfies
the following invariance property.

\begin{lemma}[$\mathcal{S}_k$-invariance]
\label{lem:LBJ-Sk-invariance}
For $k = 1, 2$ and all $x, y \in \mathcal{V}^m$,
\begin{equation}\label{eq:LBJ-Sk-invariance}
B_{\mathcal{S}_k u}
(\mathcal{S}_k x, \mathcal{S}_k y)
= B_u(x, y).
\end{equation}
\end{lemma}

\begin{proof}
The directed-path propagator $[B_v]_{ij}$ depends on $v$
only through the cutoffs $\mathbf{a}(v)$ and the
particle-label differences $n_{\ell'} - n_\ell$.  Both
quantities are invariant under each $\mathcal{S}_k$.

\noindent\textit{Case $k = 1$.}
Since $\mathcal{S}_1 u =
(t, \mathbf{a}, \mathbf{n}{+}\mathbf{1})$, the cutoffs
$\mathbf{a}$ are unchanged and the particle-label
differences $(n_j{+}1) - (n_i{+}1) = n_j - n_i$ are
unchanged.  Moreover, $\mathcal{S}_1$ does not shift
$\mathbf{a}$, so it acts trivially on spatial
coordinates: $B_{\mathcal{S}_1 u}(\mathcal{S}_1 x,
\mathcal{S}_1 y) = B_u(x, y)$.

\noindent\textit{Case $k = 2$.}
Since $\mathcal{S}_2 u =
(t{+}1, \mathbf{a}, \mathbf{n})$, the cutoffs
$\mathbf{a}$ are unchanged and the particle-label
differences $n_j - n_i$ are unchanged.
The shift $\mathcal{S}_2$ does not translate spatial
coordinates, and the propagator is
independent of~$t$, so
$B_{\mathcal{S}_2 u} = B_u$ trivially.
\end{proof}

\begin{lemma}\label{lem:LBJ-admissible}
The propagator $B_{\mathbf{a},\mathbf{n}}$ is an admissible propagator
for the product graph $\mathcal{G}^m$.
\end{lemma}

\begin{proof}
By Proposition~\ref{prop:directed-path-T},
$\mathcal{T}$-covariance and $\mathcal{T}$-splitting hold.
The edge weights $(c_k, \lambda_k)$ are constant scalars
and the propagator satisfies $\mathcal{S}_k$-invariance
(Lemma~\ref{lem:LBJ-Sk-invariance}), so
$\mathcal{S}_k$-compatibility follows from
Proposition~\ref{prop:constant-scalar}.

The Schatten-class estimates are identical to those of
the right Bernoulli case (\S\ref{sec:RBJ}): the
propagator uses the same transition kernel~$Q$, so the
geometric decay bound carries over unchanged.  The
contour formula for $\psi$ gives compact support
$r \le x \le r + n + t$, and the epigraph estimate gives
$|\phi_{t,r,n}(x)| \le C\beta^{x-r}
\mathbf{1}_{x > y_N}$ as in \eqref{eq:RBJ-phi-bound}.
These bounds ensure absolute convergence of the sums
defining $\Psi$, $\Phi$, and $K$ in
\eqref{eq:Psi-def}--\eqref{eq:K-def}.
\end{proof}

\paragraph{\textbf{Corrected edge weights.}}
Since $\mathcal{S}_k$-invariance holds,
Proposition~\ref{prop:constant-scalar}(ii) gives
$P(\mathcal{S}_k u) = P(u)$ and
$\Lambda_k(u) = \lambda_k I_m = D_k(u)$.
The diamond data on $\mathcal{G}^m$ is therefore
$(C_k, \Lambda_k) = (c_k I_m, \lambda_k I_m)$.
Theorem~\ref{thm:product-graph} requires additionally
that the resolvent $(I - zK(u))^{-1}$ exist.  At $z = 1$,
the Fredholm determinant
$F_{t, \mathbf{a}, \mathbf{n}} \defeq
\det(I - K_{t, \mathbf{a}, \mathbf{n}})$ equals the
multipoint gap probability
$\PP_y\bigl(\bigcap_i
\{Y_{n_i}(t) > a_i\}\bigr)$.  Since $K$ is trace class,
$F \ne 0$ if and only if the resolvent exists.  Particle
positions are non-increasing in time, so
$Y_{n_i}(t) \le y_{n_i}$.  If every jump attempt for
particles $1, \dotsc, n_m$ fails (an event of probability
$q^{n_m t}$), no pushing occurs and $Y_k(t) = y_k$ for
$k \le n_m$.  The resolvent at $z = 1$ therefore exists
precisely on
$\mathcal{R} \defeq \{(t, \mathbf{a}, \mathbf{n}) \in
\mathcal{V}^m : a_i < y_{n_i}
\text{ for all } i\}$.
The mixed diamond at a vertex~$u$ requires the resolvent
at the four vertices $u$, $\mathcal{S}_1 u$,
$\mathcal{S}_2 u$, $\mathcal{S}_1\mathcal{S}_2 u$
entering~\eqref{eq:dressed-C-explicit}, together with
the $\mathcal{T}$-shifts $\mathcal{T}u$,
$\mathcal{T}\mathcal{S}_1 u$,
$\mathcal{T}\mathcal{S}_2 u$ from the proof
of~\eqref{eq:M-inverse}.\footnote{The explicit
formula~\eqref{eq:dressed-C-explicit} involves the
resolvent only at the four $\mathcal{S}$-shifted
vertices; a direct verification of the mixed diamond
from this formula would yield the larger domain
$a_i < y_{n_i+1}$.}
Since $\mathcal{S}_2$ preserves~$\mathcal{R}$ and
$\mathcal{T}$ tightens each threshold by~$1$, the
binding constraint is
$\mathcal{T}\mathcal{S}_1 u \in \mathcal{R}$, i.e.\
$a_i < y_{n_i+1} - 1$ for all~$i$.

\subsection{Multipoint equation}
\label{sec:LBJ-multipoint}

The dressed observable at $z = 1$,
\[
\mathcal{M}(u)
= I + \Phi(u)(I - K(u))^{-1}\Psi(u)
\in \End(\R^m),
\]
has dressed edge weights
\begin{equation}\label{eq:LBJ-dressed-Ck}
\mathcal{M}_k(u)
= c_k\mathcal{M}(\mathcal{T}u)^{-1}
  \mathcal{M}(\mathcal{S}_k u).
\end{equation}
For readability, we write
$\mathcal{M}_{t, \mathbf{a}, \mathbf{n}} \defeq
\mathcal{M}(t, \mathbf{a}, \mathbf{n})$.

\begin{theorem}\label{thm:LBJ-H6}
The dressed observable satisfies the mixed diamond
equation
\begin{equation}\label{eq:H6-multipoint}
\bigl(q\mathcal{M}_{t+1,
\mathbf{a}+\mathbf{1},
\mathbf{n}}^{-1}
+ p\mathcal{M}_{t,
\mathbf{a}+\mathbf{1},
\mathbf{n}+\mathbf{1}}^{-1}\bigr)
\mathcal{M}_{t+1, \mathbf{a},
\mathbf{n}+\mathbf{1}}
- \mathcal{M}_{t,
\mathbf{a}+\mathbf{1},
\mathbf{n}}^{-1}
\bigl(q\mathcal{M}_{t, \mathbf{a},
\mathbf{n}+\mathbf{1}}
+ p\mathcal{M}_{t+1, \mathbf{a},
\mathbf{n}}\bigr) = 0,
\end{equation}
at every $(t, \mathbf{a}, \mathbf{n})$ with
$a_i < y_{n_i+1} - 1$ for all~$i$.
\end{theorem}

\begin{proof}
By Theorem~\ref{thm:product-graph}, the pair
$(\mathcal{M}_k, \Lambda_k)$ satisfies the diamond
equations on~$\mathcal{R}$.
Since $\Lambda_k = \lambda_k I_m$, the scalar factors
commute with $\mathcal{M}_k$, and the mixed diamond
equation reduces to
\begin{equation}\label{eq:LBJ-mixed-scalar-Lambda}
\lambda_2\bigl[
\mathcal{M}_1(u) - \mathcal{M}_1(\mathcal{S}_2 u)
\bigr]
+ \lambda_1\bigl[
\mathcal{M}_2(\mathcal{S}_1 u)
- \mathcal{M}_2(u)
\bigr] = 0.
\end{equation}
Substituting \eqref{eq:LBJ-dressed-Ck} and the explicit
shifts, every term carries a common factor of
$(1{-}\theta)^{-1}$.  Dividing by this factor
gives
\begin{equation}\label{eq:LBJ-mixed-reduced}
\begin{aligned}
-q\bigl[
&\mathcal{M}_{t, \mathbf{a}+\mathbf{1},
\mathbf{n}}^{-1}
\mathcal{M}_{t, \mathbf{a},
\mathbf{n}+\mathbf{1}}
- \mathcal{M}_{t+1,
\mathbf{a}+\mathbf{1},
\mathbf{n}}^{-1}
\mathcal{M}_{t+1, \mathbf{a},
\mathbf{n}+\mathbf{1}}
\bigr] \\
+ p\bigl[
&\mathcal{M}_{t, \mathbf{a}+\mathbf{1},
\mathbf{n}+\mathbf{1}}^{-1}
\mathcal{M}_{t+1, \mathbf{a},
\mathbf{n}+\mathbf{1}}
- \mathcal{M}_{t, \mathbf{a}+\mathbf{1},
\mathbf{n}}^{-1}
\mathcal{M}_{t+1, \mathbf{a},
\mathbf{n}}
\bigr] = 0.
\end{aligned}
\end{equation}
The middle two terms share the right factor
$\mathcal{M}_{t+1, \mathbf{a},
\mathbf{n}+\mathbf{1}}$, and the
first and last share the left factor
$\mathcal{M}_{t, \mathbf{a}+\mathbf{1},
\mathbf{n}}^{-1}$.  Factoring yields
\eqref{eq:H6-multipoint}.
\end{proof}

For $m = 1$, the dressed observable $\mathcal{M}_{t,a,n}$
is scalar.

\begin{corollary}\label{cor:LBJ-scalar}
The one-point Fredholm determinant
$F_{t, a, n} = \det(I - K_{t, a, n})$ satisfies
\begin{equation}\label{eq:LBJ-scalar-HM}
qF_{t, a, n+1}F_{t+1, a+1, n}
+ pF_{t+1, a, n}F_{t, a+1, n+1}
- F_{t, a+1, n}F_{t+1, a, n+1}
= 0,
\end{equation}
at every $(t, a, n)$ with $a < y_{n+1} - 1$.
\end{corollary}

\begin{proof}
Proposition~\ref{prop:scalar-HM} applies at $z = 1$: the
boundary condition $F(T^\ell v) \to 1$ as
$\ell \to -\infty$ holds because the thresholds decrease
to $-\infty$ and the kernel vanishes.  The coefficients are
$\alpha_{12} = c_1 \lambda_2 =
-(1{-}\theta)^{-1}q$
and
$\alpha_{21} = c_2 \lambda_1 =
p(1{-}\theta)^{-1}$.
Substituting the lattice points into
\eqref{eq:HM-variable} and computing
$\alpha_{12} - \alpha_{21}
= -(1{-}\theta)^{-1}$
(using $p + q = 1$), then dividing by
$-(1{-}\theta)^{-1}$, gives
\eqref{eq:LBJ-scalar-HM}.
\end{proof}

\section{Parallel TASEP}\label{sec:PTASEP}

\paragraph{\textbf{System description.}}
The discrete-time totally asymmetric simple
exclusion process with right Bernoulli jumps, blocking
interaction, and parallel
update
consists of $N$ particles on
$\Z$.  The configuration at time
$t \in \Z_{\ge 0}$ is denoted by
\[
Y(t) = (Y_1(t) > Y_2(t) > \cdots > Y_N(t)),
\]
where $Y_k(t)$ is the position of the $k$-th particle
from the right.

At each discrete time step $t \ge 1$, all particles
attempt to jump simultaneously.  During its update from
time $t{-}1$ to $t$, particle $k$ attempts to jump one
step to the right with probability $p \in (0,1)$, and
stays put with probability $q = 1 - p$.  Because the
updates occur in parallel, a particle's jump is blocked
by the position of its right neighbour at time $t{-}1$
rather than at time $t$.  The evolution is given by
\[
Y_1(t) = Y_1(t{-}1) + \xi(t,1),
\]
\[
Y_k(t) = \min\{Y_k(t{-}1) + \xi(t,k),
Y_{k-1}(t{-}1) - 1\},
\qquad k = 2, \dots, N,
\]
where $\{\xi(t,k)\}_{t \ge 1, 1 \le k \le N}$ are
independent Bernoulli random variables with
$\PP(\xi(t,k) = 1) = p$, and the convention
$Y_0(t) = \infty$ means that the rightmost particle
faces no blocking constraint.

\paragraph{\textbf{Fredholm determinant formula.}}
Fix a time $t \ge 0$ and an initial configuration
$Y(0) = y = (y_1 > y_2 > \cdots > y_N)$.  Fix
$m \in \{1, \dots, N{-}1\}$, particle labels
$\mathbf{n} = (n_1, n_2, \dots, n_m)$ with
$1 \le n_1 < n_2 < \cdots < n_m < N$, and spatial
thresholds $\mathbf{a} = (a_1, \dots, a_m) \in
\Z^m$.  The multipoint joint cumulative
distribution of the particle positions is given by a
Fredholm determinant on
$\ell^2(\{n_1, \dots, n_m\} \times \Z)$~\cite{MatetskiRemenik2023}:
\[
\PP_y\Bigl(\bigcap_{i=1}^m
\{Y_{n_i}(t) > a_i\}\Bigr)
= \det(I - \bar{\chi}_{\mathbf{a}} K_t
\bar{\chi}_{\mathbf{a}})_{\ell^2(\{n_1, \dots, n_m\}
\times \Z)},
\]
where $\bar{\chi}_{\mathbf{a}}(n_i, x) =
\mathbf{1}_{x \le a_i}$.  The correlation kernel $K_t$ has
block entries
\[
K_t(n_i, x_i; n_j, x_j)
= -Q^{n_j - n_i}(x_i, x_j)\mathbf{1}_{n_i < n_j}
+ \bigl(\mathcal{S}_{-t,-n_i}^*
\overline{\mathcal{S}}_{-t,n_j}^{\operatorname{epi}(y)}
\bigr)(x_i, x_j),
\]
where $A^*$ denotes the transpose kernel,
$A^*(x, y) = A(y, x)$.

Fix an auxiliary parameter $\theta \in (0, 1)$ satisfying
$\theta p/q < 1$, which enters the Fredholm determinant
representation but does not appear in the final
multipoint equation \eqref{eq:H4-multipoint}, and let
$\alpha = (1{-}\theta)\theta^{-1}(q{+}p\theta)^{-1}$.
The operators defining the kernel are as follows.

\medskip
\noindent\textit{The transition matrix $Q$.}
The matrix $Q$ is the transition matrix of a modified random walk on
$\Z$ strictly to the left:
\[
Q(z_1, z_2) = \frac{1{-}\theta}{q{+}p\theta}
\theta^{z_1 - z_2 - 1}
q_{z_1 - z_2}\mathbf{1}_{z_1 > z_2},
\]
where $q_1 = q$ and $q_i = 1$ for $i > 1$.
Its $n$-th power admits the integral representation
\[
Q^n(z_1, z_2) = \frac{\alpha^n}{2\pi\mathrm{i}}
\oint_{\gamma_\rho} \diff w
\frac{\theta^{z_1 - z_2}}{w^{z_1 - z_2 - n + 1}}
\Bigl(\frac{q{+}pw}{1 - w}\Bigr)^{n}.
\]

\medskip
\noindent\textit{The operators $\mathcal{S}_{-t,-n}$ and
$\overline{\mathcal{S}}_{-t,n}$.}\footnote{The kernel
operators $\mathcal{S}_{-t,-n}$ and
$\overline{\mathcal{S}}_{-t,n}$ are unrelated to the
lattice shifts $\mathcal{S}_k$ of the diamond framework
(\S\ref{sec:linear-problem}); the notation follows
Matetski--Remenik~\cite{MatetskiRemenik2023}.}
These operators encode the Bernoulli
jumps via the probability generating function
$\varphi(w) = q + pw$.  The contour integral representations are
\[
\mathcal{S}_{-t,-n}(z_1, z_2)
= \frac{\alpha^{-n+1}}{2\pi\mathrm{i}}
\oint_{\gamma_\rho} \frac{\diff w}{w}
\frac{\theta^{z_2 - z_1}}{w^{z_2 - z_1 + n}}
(1{-}w)^n(q{+}pw)^{t-n+1},
\]
\[
\overline{\mathcal{S}}_{-t,n}(z_1, z_2)
= \frac{\alpha^{n-1}}{2\pi\mathrm{i}}
\oint_{\gamma_\delta} \frac{\diff w}{w}
\frac{(1{-}w)^{z_2 - z_1 + n - 1}}
{\theta^{z_2 - z_1} w^{n-1}}
(q{+}p(1{-}w))^{n-1-t},
\]
where $\gamma_\rho$ is a positively oriented circle of
radius $\rho \in (0, 1)$ around the origin (with the
additional restriction $\rho < q/p$ when $p \ge 1/2$
and $t < n_m - 1$), and
$\gamma_\delta$ is a small circle around the origin
chosen so that $(q + p(1{-}w))^{n-1-t}$ is analytic
inside.

\medskip
\noindent\textit{The epigraph operator.}
The operator
$\overline{\mathcal{S}}_{-t,n}^{\operatorname{epi}(y)}$
encodes the initial data through a hitting probability:
\[
\overline{\mathcal{S}}_{-t,n}^{\operatorname{epi}(y)}
(z_1, z_2)
= \E_{W_0 = z_1}\bigl[
\overline{\mathcal{S}}_{-t,n-\tau}(W_\tau, z_2)
\mathbf{1}_{\tau < n}\bigr],
\]
where $W_\ell$ is the random walk with transition matrix
$Q$, and
$\tau = \min\{\ell \in \{0, \dots, N{-}1\} :
W_\ell > y_{\ell+1}\}$,
with the convention $\tau = \infty$ if the set is empty.

\paragraph{\textbf{References.}}
The Fredholm determinant formula is due to Matetski--Remenik~\cite{MatetskiRemenik2023} (Prop.~2.3, with $\kappa = 1$ and $\varphi(w) = q + pw$).

\subsection{Kernel reformulation}
\label{sec:PTASEP-kernel-reform}

We rewrite the extended-kernel Fredholm determinant in the
form of the product graph construction of
\S\ref{sec:product-graph}.  The functions $\psi, \phi$
below provide the seed data, while the lower-triangular
operator $B_{\mathbf{a}, \mathbf{n}}$ will be verified to
satisfy the admissibility conditions of
Definition~\ref{def:admissible-prop}.

\paragraph{\textbf{Seed data.}}
The seed functions are
\[
\psi_{t,r,n}(x) \defeq \mathcal{S}_{-t,-n}^*(r,x)
= \mathcal{S}_{-t,-n}(x,r),
\qquad
\phi_{t,r,n}(x) \defeq
\overline{\mathcal{S}}_{-t,n}^{\operatorname{epi}(y)}
(x,r),
\]
with the dependence of $\phi$ on the initial
condition $y$ suppressed.

\paragraph{\textbf{Propagator.}}
The propagator is an instance of the directed-path
propagator (Definition~\ref{def:directed-path-propagator})
with transition kernels
$Q_{\ell,\ell'} = Q^{n_{\ell'} - n_\ell}$ for
$\ell < \ell'$.  On the cutoff region $r \le a_i$,
$r' \le a_j$, the entries of
$I + B_{\mathbf{a}, \mathbf{n}}$ coincide with those of
$((I + \bar{\chi}_{\mathbf{a}} L
\bar{\chi}_{\mathbf{a}})^{-1})^{\!\top}$, where $L$
is the strictly upper-triangular block operator with
entries $L_{ij}(x, x') = Q^{n_j - n_i}(x, x')
\mathbf{1}_{i < j}$ (Remark~\ref{rem:Neumann-series}).
The propagator depends on $\mathbf{a}$ and $\mathbf{n}$
but not on the time coordinate $t$, which is accordingly
suppressed from the notation.

\begin{lemma}\label{lem:PTASEP-kernel}
We have
\[
\PP_y\Bigl(\bigcap_{i=1}^m
\{Y_{n_i}(t) > a_i\}\Bigr)
= \det(I - K_{t, \mathbf{a}, \mathbf{n}})_{\ell^2(\Z)},
\]
where
\[
K_{t, \mathbf{a}, \mathbf{n}}(x, x')
= \sum_{1 \le j \le i \le m}
  \sum_{r \le a_i} \sum_{r' \le a_j}
   \psi_{t, r, n_i}(x)
  \bigl[\delta_{ij}\delta_{r,r'}
  + [B_{\mathbf{a}, \mathbf{n}}]_{ij}(r, r')\bigr]
  \phi_{t, r', n_j}(x').
\]
This identifies $K_{t, \mathbf{a}, \mathbf{n}}$ with the
kernel~\eqref{eq:K-def} on the product graph
$\mathcal{G}^m$, with $\psi, \phi$ as the seed data and
$B_{\mathbf{a}, \mathbf{n}}$ as the propagator.
\end{lemma}

\begin{proof}
The argument is identical to that of
Lemma~\ref{lem:RBJ-kernel}: the block-triangular
decomposition, Sylvester's identity, and the transpose
identification produce the claimed single-space kernel.
The Schatten-class estimates require one modification,
because the factor $(q{+}pw)^{t-n+1}$ in the contour
formula for $\psi$ has a negative exponent when
$t < n - 1$.

The modified one-step kernel satisfies
$Q(x,z) \le C\theta^{x-z}\mathbf{1}_{x > z}$, since
$q_s \le 1$ for every $s \ge 1$.  The $d$-fold
convolution inherits geometric tails: for any
$\vartheta \in (\theta, 1)$,
\begin{equation}\label{eq:PTASEP-Q-tail}
Q^d(x,z) \le C_{d,\vartheta}\vartheta^{x-z}
\mathbf{1}_{x - z \ge d}.
\end{equation}

The contour formula for $\psi$ gives
\[
\psi_{t,r,n}(x) = \alpha^{-n+1}\theta^{r-x}
[w^{r-x+n}](1{-}w)^n(q{+}pw)^{t-n+1}.
\]
Write $k = r - x + n$.  If $t \ge n - 1$, the exponent
$t - n + 1 \ge 0$ makes $(q{+}pw)^{t-n+1}$ a polynomial,
and the coefficient vanishes outside
$0 \le k \le t + 1$, so $\psi$ has compact support
in~$x$.  If $M \defeq n - 1 - t > 0$, the factor
$(q{+}pw)^{-M} = q^{-M}(1 + (p/q)w)^{-M}$ has
coefficients bounded by $C(1{+}k)^{M-1}(p/q)^k$ by the
binomial expansion, and since $(1{-}w)^n$ is a polynomial
of fixed degree, the product has the same exponential
bound.  With
$\gamma \defeq \theta p/q$, the prefactor
$\theta^{r-x} = \theta^{k-n}$ combines with $(p/q)^k$ to
give
\begin{equation}\label{eq:PTASEP-psi-tail}
|\psi_{t,r,n}(x)| \le C(1{+}k)^{M-1}\gamma^k
\mathbf{1}_{x \le r + n}.
\end{equation}
Since $\theta p/q < 1$, we have $\gamma < 1$, and
$\psi_{t,r,n} \in \ell^2(\Z)$ whether $t \ge n - 1$
or $t < n - 1$.  The
integrand depends on $r$ and $x$ only through $r - x$, so
the $\ell^2$ norm of $\psi_{t,r,n}$ is independent
of~$r$.

The contour estimate for $\phi$ follows the same argument
as in Lemma~\ref{lem:RBJ-kernel}.  Choose
$\beta \in (\theta, 1)$ and $\gamma_\delta$ so that
$\sup_{w \in \gamma_\delta}|\tfrac{\theta}{1{-}w}| < \beta$,
giving
$|\overline{\mathcal{S}}_{-t,k}(z,r)|
\le C\beta^{z-r}$ for $z > y_N$ and
$r \le a_{\max}$.  Since $Q$ has geometric tails with
rate $\theta < \beta$, the exponential moment
$\E[\beta^{-S_\ell}]$ is finite, where $S_\ell$ is the
cumulative displacement of the $Q$-walk after
$\ell$ steps.  The hitting
representation gives
\begin{equation}\label{eq:PTASEP-phi-tail}
|\phi_{t,r,n}(x)| \le C\beta^{x-r}
\mathbf{1}_{x > y_N}.
\end{equation}

With these bounds, the same conjugation and Sylvester
reduction as in Lemma~\ref{lem:RBJ-kernel} apply.
\end{proof}

\subsection{Base graph and seed data}
\label{sec:PTASEP-base-graph}

We now identify the lattice structure that governs how
the seed data evolve under shifts of the parameters
$(t, a, n)$, and verify that $\psi, \phi$ satisfy the
linear problems of the diamond framework.

Let $\mathcal{V} \subseteq \Z^3$ be the lattice
generated by the shifts
\[
T = e^{\partial_a}, \qquad
S_1 = e^{\partial_t + \partial_n}, \qquad
S_2 = e^{\partial_t + \partial_a},
\]
so that for $u = (t, a, n) \in \mathcal{V}$,
\[
Tu = (t, a{+}1, n), \qquad
S_1 u = (t{+}1, a, n{+}1), \qquad
S_2 u = (t{+}1, a{+}1, n).
\]
Take $E = \mathbb{F} = \R$ and define the
constant scalar edge weights
\begin{equation}\label{eq:PTASEP-edge-weights}
(c_1, c_2) = \bigl(
(q{+}p\theta)(1{-}\theta)^{-1}, q\bigr),
\qquad
(\lambda_1, \lambda_2)
= \bigl(\theta(q{+}p\theta)(1{-}\theta)^{-1},
-p\theta\bigr).
\end{equation}
Since $(c_k, \lambda_k)$ are constant scalars, the
diamond equations
\eqref{eq:diamond-C-ij}--\eqref{eq:diamond-mixed-ij}
are satisfied trivially.

\begin{lemma}\label{lem:PTASEP-seed-linear}
Let $H = \ell^2(\Z)$.  The seed functions
$\psi_{t,r,n}$ and $\phi_{t,r,n}$ defined in
\S\ref{sec:PTASEP-kernel-reform} satisfy the linear
problem \eqref{eq:Linear-Psi-3} and the adjoint
linear problem \eqref{eq:Linear-Phi-Prob} on
$\mathcal{V}$ with edge weights
$(c_k, \lambda_k)$.  Explicitly, for each fixed
$r \in \Z$:
\begin{align}
\psi_{t+1,r,n+1}
&= (q{+}p\theta)(1{-}\theta)^{-1}
   \psi_{t,r+1,n}
   - \theta(q{+}p\theta)(1{-}\theta)^{-1}
   \psi_{t,r,n},
\label{eq:PTASEP-psi-S1} \\
\psi_{t+1,r+1,n}
&= q\psi_{t,r+1,n}
   + p\theta\psi_{t,r,n},
\label{eq:PTASEP-psi-S2}
\end{align}
and
\begin{align}
\phi_{t,r+1,n}
&= (q{+}p\theta)(1{-}\theta)^{-1}
   \phi_{t+1,r,n+1}
   - \theta(q{+}p\theta)(1{-}\theta)^{-1}
   \phi_{t+1,r+1,n+1},
\label{eq:PTASEP-phi-S1} \\
\phi_{t,r+1,n}
&= q\phi_{t+1,r+1,n}
   + p\theta\phi_{t+1,r+2,n}.
\label{eq:PTASEP-phi-S2}
\end{align}
\end{lemma}

\begin{proof}
Each identity is verified in turn.

\medskip
\noindent\textit{Verification of
\eqref{eq:PTASEP-psi-S1}.}
From the contour representation,
\[
\psi_{t,r,n}(x) = \frac{\alpha^{-n+1}}{2\pi\mathrm{i}}
  \oint_{\gamma_\rho} \frac{\diff w}{w}
  \frac{\theta^{r-x}}{w^{r-x+n}}
  (1{-}w)^n(q{+}pw)^{t-n+1}.
\]
The right-hand side of \eqref{eq:PTASEP-psi-S1} is
\begin{align*}
c_1\psi_{t,r+1,n}(x)
- \lambda_1\psi_{t,r,n}(x)
= \frac{\alpha^{-n+1}}{2\pi\mathrm{i}}
  \oint \frac{\diff w}{w}
  \frac{\theta^{r-x}}{w^{r-x+n}}
  (1{-}w)^n(q{+}pw)^{t-n+1}
  \Bigl(c_1\frac{\theta}{w}
  - \lambda_1\Bigr).
\end{align*}
The parenthetical factor simplifies as
\[
c_1\frac{\theta}{w} - \lambda_1
= \frac{(q{+}p\theta)\theta}{(1{-}\theta)w}
- \frac{\theta(q{+}p\theta)}{1{-}\theta}
= \frac{\theta(q{+}p\theta)}{1{-}\theta}
\frac{1{-}w}{w}
= \frac{1{-}w}{\alpha w},
\]
since
$\alpha^{-1} = \theta(q{+}p\theta)/(1{-}\theta)$.
Multiplying through, the prefactor shifts from
$\alpha^{-n+1}$ to $\alpha^{-n}$, the power of
$(1{-}w)$ increases from $n$ to $n{+}1$, and the
denominator gains one additional factor of $w$.  The power of
$(q{+}pw)$ remains $t - n + 1 = (t{+}1) - (n{+}1) + 1$,
since $S_1$ increments both $t$ and $n$.
This gives $\psi_{t+1,r,n+1}(x)$.

\medskip
\noindent\textit{Verification of
\eqref{eq:PTASEP-psi-S2}.}
The verification is completely analogous to that of
\eqref{eq:RBJ-psi-S2} in \S\ref{sec:RBJ}: the
right-hand side produces the parenthetical factor
$q\theta/w + p\theta = \theta(q + pw)/w$,
which shifts $t$ by $+1$ and $r$ by $+1$.  The power
of $(q + pw)$ changes from $t - n + 1$ to
$t - n + 2 = (t{+}1) - n + 1$, as required.

\medskip
\noindent\textit{Verification of
\eqref{eq:PTASEP-phi-S1}.}
The argument is analogous to the proof of
\eqref{eq:RBJ-phi-S1} in \S\ref{sec:RBJ}.  The
$\overline{\mathcal{S}}$ recurrence in $n$ is
\begin{equation}\label{eq:PTASEP-Sbar-recurrence-n}
c_1\bigl[
\overline{\mathcal{S}}_{-(t+1),n+1}(z, r)
- \theta\overline{\mathcal{S}}_{-(t+1),n+1}
  (z, r{+}1)
\bigr]
= \overline{\mathcal{S}}_{-t,n}(z, r{+}1).
\end{equation}
In the contour formula for
$\overline{\mathcal{S}}_{-(t+1),n+1}$, the shift from
$r$ to $r{+}1$ multiplies the integrand by
$(1{-}w)/\theta$.  On the left-hand side of
\eqref{eq:PTASEP-Sbar-recurrence-n}, the bracket
produces the factor $1 - (1{-}w) = w$, which reduces
$dw/w^{n+1}$ to $dw/w^n$.  Multiplying by
$c_1 = (q{+}p\theta)/(1{-}\theta)$ and using
$c_1 \alpha = 1/\theta$, the prefactor
$c_1 \alpha^n / \theta^{r-z}$ becomes
$\alpha^{n-1}/\theta^{r+1-z}$, matching
$\overline{\mathcal{S}}_{-t,n}(z, r{+}1)$.  The power
of $(q{+}p(1{-}w))$ is unchanged since
$(n{+}1) - 1 - (t{+}1) = (n{-}1) - t$.

We show that
$\overline{\mathcal{S}}_{-t,0} \equiv 0$.  When $n = 0$,
the measure $dw/w^n$ reduces to $dw$, removing the pole
at the origin.  The remaining integrand is analytic
inside $\gamma_\delta$, so the integral vanishes by
Cauchy's theorem.

The lifting to the epigraph version proceeds as in
\S\ref{sec:RBJ}: the expectation splits at
$\tau \le n{-}1$ (where
\eqref{eq:PTASEP-Sbar-recurrence-n} applies) and
$\tau = n$ (where
$\overline{\mathcal{S}}_{-t,0} = 0$ kills the
boundary term).

\medskip
\noindent\textit{Verification of
\eqref{eq:PTASEP-phi-S2}.}
The verification is analogous to that of
\eqref{eq:RBJ-phi-S2}: the recurrence
$q\overline{\mathcal{S}}_{-(t+1),n}(z, r{+}1)
+ p\theta\overline{\mathcal{S}}_{-(t+1),n}
  (z, r{+}2)
= \overline{\mathcal{S}}_{-t,n}(z, r{+}1)$
follows because the shift $r \mapsto r{+}1$ introduces
a factor $\theta^{-1}(1{-}w)$ from the kernel, which
combined with the coefficient $p\theta$ yields
$q + p\theta\theta^{-1}(1{-}w) = q + p(1{-}w)$,
promoting the power of $(q{+}p(1{-}w))$ by $+1$.
Since $S_2$ does not shift the particle label~$n$,
the hitting time $\tau$ is unchanged and the lifting
to the epigraph version is immediate.
\end{proof}

\subsection{Admissible propagators and corrected
edge weights}
\label{sec:PTASEP-propagators}

The product graph $\mathcal{G}^m$ has vertex set
$\mathcal{V}^m$ and is generated by the diagonal shifts
\[
\mathcal{T}u
= (t, \mathbf{a} + \mathbf{1}, \mathbf{n}),
\qquad
\mathcal{S}_1 u
= (t{+}1, \mathbf{a}, \mathbf{n} + \mathbf{1}),
\qquad
\mathcal{S}_2 u
= (t{+}1, \mathbf{a} + \mathbf{1}, \mathbf{n}).
\]
The block-diagonal edge weights are
$C_k(u) = c_k I_m$ and $D_k(u) = \lambda_k I_m$
(the uncorrected edge weights of
Definition~\ref{def:admissible-prop}).

The propagator is an instance of
the directed-path propagator
(Definition~\ref{def:directed-path-propagator}) with
transition kernels $Q_{\ell,\ell'} = Q^{n_{\ell'} - n_\ell}$
for $\ell < \ell'$.  The standing hypotheses of
\S\ref{sec:directed-path} are satisfied: translation
invariance of $Q^n$ is inherited from the translation
invariance of~$Q$, and $\mathcal{T}$-invariance holds
because the transition kernels depend only on the
particle-label differences, which are unchanged
by~$\mathcal{T}$.

By Proposition~\ref{prop:constant-scalar},
$\mathcal{S}_k$-compatibility follows from the
constant-scalar edge weights once the propagator satisfies
the following invariance property.

\begin{lemma}[$\mathcal{S}_k$-invariance]
\label{lem:PTASEP-Sk-invariance}
For $k = 1, 2$ and all $x, y \in \mathcal{V}^m$,
\begin{equation}\label{eq:PTASEP-Sk-invariance}
B_{\mathcal{S}_k u}
(\mathcal{S}_k x, \mathcal{S}_k y)
= B_u(x, y).
\end{equation}
\end{lemma}

\begin{proof}
The directed-path propagator $[B_v]_{ij}$ depends on $v$
only through the cutoffs $\mathbf{a}(v)$ and the
particle-label differences $n_{\ell'} - n_\ell$.

\noindent\textit{Case $k = 1$.}
Since $\mathcal{S}_1 u =
(t{+}1, \mathbf{a}, \mathbf{n}{+}\mathbf{1})$, the
cutoffs $\mathbf{a}$ are unchanged and the particle-label
differences $(n_j{+}1) - (n_i{+}1) = n_j - n_i$ are
unchanged.  Moreover, $\mathcal{S}_1$ does not shift
$\mathbf{a}$, so it acts trivially on spatial
coordinates: $B_{\mathcal{S}_1 u}(\mathcal{S}_1 x,
\mathcal{S}_1 y) = B_u(x, y)$.

\noindent\textit{Case $k = 2$.}
Since $\mathcal{S}_2 u =
(t{+}1, \mathbf{a}{+}\mathbf{1}, \mathbf{n})$, the
particle-label differences are unchanged but the
cutoffs shift to $\mathbf{a}{+}\mathbf{1}$.  Since
$\mathcal{S}_2$ acts on thresholds and spatial
coordinates identically to $\mathcal{T}$, and the
transition kernels do not depend on~$t$,
$\mathcal{S}_2$-invariance reduces to
$\mathcal{T}$-covariance
(Proposition~\ref{prop:directed-path-T}).
\end{proof}

\begin{lemma}\label{lem:PTASEP-admissible}
The propagator $B_{\mathbf{a}, \mathbf{n}}$ is an
admissible propagator for the product graph
$\mathcal{G}^m$.
\end{lemma}

\begin{proof}
By Proposition~\ref{prop:directed-path-T},
$\mathcal{T}$-covariance and $\mathcal{T}$-splitting hold.
The edge weights $(c_k, \lambda_k)$ are constant scalars
and the propagator satisfies $\mathcal{S}_k$-invariance
(Lemma~\ref{lem:PTASEP-Sk-invariance}), so
$\mathcal{S}_k$-compatibility follows from
Proposition~\ref{prop:constant-scalar}.

It remains to verify regularity.  The directed-path sum
is finite in the chain length, and the support condition
$Q^d(x,z) = 0$ unless $x - z \ge d$ makes the internal
spatial sums finite for fixed endpoints.  The tail
estimate \eqref{eq:PTASEP-Q-tail} gives geometric decay
of the propagator entries.  Together with
\eqref{eq:PTASEP-psi-tail} and
\eqref{eq:PTASEP-phi-tail}, these bounds prove absolute
convergence of the sums defining $\Psi$, $\Phi$,
and $K$ in
\eqref{eq:Psi-def}--\eqref{eq:K-def}.
\end{proof}

\paragraph{\textbf{Corrected edge weights.}}
Since $\mathcal{S}_k$-invariance holds,
Proposition~\ref{prop:constant-scalar}(ii) gives
$P(\mathcal{S}_k u) = P(u)$ and
$\Lambda_k(u) = \lambda_k I_m = D_k(u)$.
The diamond data on $\mathcal{G}^m$ are therefore
$(C_k, \Lambda_k) = (c_k I_m, \lambda_k I_m)$.
Theorem~\ref{thm:product-graph} requires additionally
that the resolvent $(I - zK(u))^{-1}$ exist.  At $z = 1$,
the Fredholm determinant
$F_{t, \mathbf{a}, \mathbf{n}} \defeq
\det(I - K_{t, \mathbf{a}, \mathbf{n}})$ equals the
multipoint gap probability
$\PP_y\bigl(\bigcap_i
\{Y_{n_i}(t) > a_i\}\bigr)$.  Since $K$ is trace class,
$F \ne 0$ if and only if the resolvent exists.  Under
parallel update, particle $k$ is blocked by the position
of particle $k - 1$ at time $t - 1$, giving
$Y_k(t) \le h_k(t)$ where
\[
h_1(t) = y_1 + t, \qquad
h_k(t) = \min\{h_k(t{-}1) + 1, h_{k-1}(t{-}1) - 1\},
\quad k \ge 2,
\]
with $h_k(0) = y_k$.  If every jump attempt succeeds (an
event of probability $p^{Nt}$), then $Y_k(t) = h_k(t)$
for all~$k$.  The resolvent at $z = 1$ therefore exists
precisely on
$\mathcal{R} \defeq \{(t, \mathbf{a}, \mathbf{n}) \in
\mathcal{V}^m : a_i < h_{n_i}(t)
\text{ for all } i\}$.
The mixed diamond at a vertex~$u$ requires the resolvent
at the four vertices $u$, $\mathcal{S}_1 u$,
$\mathcal{S}_2 u$, $\mathcal{S}_1\mathcal{S}_2 u$
entering~\eqref{eq:dressed-C-explicit}, together with
the $\mathcal{T}$-shifts $\mathcal{T}u$,
$\mathcal{T}\mathcal{S}_1 u$,
$\mathcal{T}\mathcal{S}_2 u$ from the proof
of~\eqref{eq:M-inverse}.\footnote{The explicit
formula~\eqref{eq:dressed-C-explicit} involves the
resolvent only at the four $\mathcal{S}$-shifted
vertices; a direct verification of the mixed diamond
from this formula would yield the larger domain
$a_i < h_{n_i+1}(t+1)$.}
Since $h_k$ is non-decreasing in $t$ and
$h_{k+1}(t) \le h_k(t) - 1$ for all $k$ and $t$, the
binding constraint is
$\mathcal{T}\mathcal{S}_1 u \in \mathcal{R}$, i.e.\
$a_i < h_{n_i+1}(t+1) - 1$ for all~$i$.

\subsection{Multipoint equation}
\label{sec:PTASEP-multipoint}

The dressed observable at $z = 1$,
\[
\mathcal{M}(u)
= I + \Phi(u)(I - K(u))^{-1}\Psi(u)
\in \End(\R^m),
\]
has dressed edge weights
\begin{equation}\label{eq:PTASEP-dressed-Ck}
\mathcal{M}_k(u)
= c_k\mathcal{M}(\mathcal{T}u)^{-1}
  \mathcal{M}(\mathcal{S}_k u).
\end{equation}
For readability, we write
$\mathcal{M}_{t, \mathbf{a}, \mathbf{n}} \defeq
\mathcal{M}(t, \mathbf{a}, \mathbf{n})$.

\begin{theorem}\label{thm:PTASEP-H4}
The dressed observable satisfies the mixed diamond
equation
\begin{equation}\label{eq:H4-multipoint}
\begin{aligned}
&\bigl(q\mathcal{M}_{t+1, \mathbf{a}+\mathbf{1},
\mathbf{n}+\mathbf{1}}^{-1}
+ p\mathcal{M}_{t+1, \mathbf{a}+\mathbf{2},
\mathbf{n}}^{-1}\bigr)
\mathcal{M}_{t+2, \mathbf{a}+\mathbf{1},
\mathbf{n}+\mathbf{1}} \\
&\qquad - \mathcal{M}_{t, \mathbf{a}+\mathbf{1},
\mathbf{n}}^{-1}
\bigl(q\mathcal{M}_{t+1, \mathbf{a}+\mathbf{1},
\mathbf{n}}
+ p\mathcal{M}_{t+1, \mathbf{a},
\mathbf{n}+\mathbf{1}}\bigr) = 0,
\end{aligned}
\end{equation}
at every $(t, \mathbf{a}, \mathbf{n})$ with
$a_i < h_{n_i+1}(t+1) - 1$ for all~$i$.
\end{theorem}

\begin{proof}
By Theorem~\ref{thm:product-graph}, the pair
$(\mathcal{M}_k, \Lambda_k)$ satisfies the diamond
equations on~$\mathcal{R}$.
Since $\Lambda_k = \lambda_k I_m$, the scalar factors
commute with $\mathcal{M}_k$, and the mixed diamond
equation reduces to
\begin{equation}\label{eq:PTASEP-mixed-scalar-Lambda}
\lambda_2\bigl[
\mathcal{M}_1(u) - \mathcal{M}_1(\mathcal{S}_2 u)
\bigr]
+ \lambda_1\bigl[
\mathcal{M}_2(\mathcal{S}_1 u)
- \mathcal{M}_2(u)
\bigr] = 0.
\end{equation}
Substituting \eqref{eq:PTASEP-dressed-Ck} and the explicit
shifts, every term carries a common factor of
$\theta(q{+}p\theta)(1{-}\theta)^{-1}$
(from $c_1 \lambda_2 =
-p\theta(q{+}p\theta)(1{-}\theta)^{-1}$
and $c_2 \lambda_1 =
q\theta(q{+}p\theta)(1{-}\theta)^{-1}$).
Dividing by this factor gives
\begin{multline}\label{eq:PTASEP-mixed-reduced}
-p\bigl[
\mathcal{M}_{t, \mathbf{a}+\mathbf{1},
\mathbf{n}}^{-1}
\mathcal{M}_{t+1, \mathbf{a},
\mathbf{n}+\mathbf{1}}
- \mathcal{M}_{t+1,
\mathbf{a}+\mathbf{2},
\mathbf{n}}^{-1}
\mathcal{M}_{t+2, \mathbf{a}+\mathbf{1},
\mathbf{n}+\mathbf{1}}
\bigr] \\
+ q\bigl[
\mathcal{M}_{t+1, \mathbf{a}+\mathbf{1},
\mathbf{n}+\mathbf{1}}^{-1}
\mathcal{M}_{t+2, \mathbf{a}+\mathbf{1},
\mathbf{n}+\mathbf{1}}
- \mathcal{M}_{t, \mathbf{a}+\mathbf{1},
\mathbf{n}}^{-1}
\mathcal{M}_{t+1, \mathbf{a}+\mathbf{1},
\mathbf{n}}
\bigr] = 0.
\end{multline}
The middle two terms share the right factor
$\mathcal{M}_{t+2, \mathbf{a}+\mathbf{1},
\mathbf{n}+\mathbf{1}}$, and the
first and last share the left factor
$\mathcal{M}_{t, \mathbf{a}+\mathbf{1},
\mathbf{n}}^{-1}$.  Factoring yields
\eqref{eq:H4-multipoint}.
\end{proof}

For $m = 1$, the dressed observable $\mathcal{M}_{t,a,n}$
is scalar.

\begin{corollary}\label{cor:PTASEP-scalar}
The one-point Fredholm determinant
$F_{t, a, n} = \det(I - K_{t, a, n})$ satisfies
\begin{equation}\label{eq:PTASEP-scalar-HM}
pF_{t+1, a, n+1}F_{t+1, a+2, n}
+ qF_{t+1, a+1, n}F_{t+1, a+1, n+1}
- F_{t, a+1, n}F_{t+2, a+1, n+1}
= 0,
\end{equation}
at every $(t, a, n)$ with $a < h_{n+1}(t+1) - 1$.
\end{corollary}

\begin{proof}
Proposition~\ref{prop:scalar-HM} applies at $z = 1$: the
boundary condition $F(T^\ell v) \to 1$ as
$\ell \to -\infty$ holds because the thresholds decrease
to $-\infty$ and the kernel vanishes.  The coefficients are
$\alpha_{12} = c_1 \lambda_2 =
-p\theta(q{+}p\theta)(1{-}\theta)^{-1}$
and
$\alpha_{21} = c_2 \lambda_1 =
q\theta(q{+}p\theta)(1{-}\theta)^{-1}$.
Substituting the lattice points into
\eqref{eq:HM-variable} and computing
$\alpha_{12} - \alpha_{21}
= -\theta(q{+}p\theta)(1{-}\theta)^{-1}$
(using $p + q = 1$), then dividing by
$-\theta(q{+}p\theta)(1{-}\theta)^{-1}$, gives
\eqref{eq:PTASEP-scalar-HM}.
\end{proof}

\section{Right Geometric Jumps}
\label{sec:RGP}

\paragraph{\textbf{System description.}}
The right geometric jump model is a discrete-time totally
asymmetric simple exclusion process with geometric jumps,
blocking interaction, and parallel update.  The system consists of
$N$ particles on $\Z$.  The configuration at time
$t \in \Z_{\ge 0}$ is denoted by
\[
Y(t) = (Y_1(t) > Y_2(t) > \cdots > Y_N(t)),
\]
where $Y_k(t)$ is the position of the $k$-th particle
from the right.

At each discrete time step $t \ge 1$, each particle $k$
attempts to jump to the right by a distance sampled from a
$\operatorname{Geom}[p]$ distribution, $p \in (0,1)$,
(supported on $\{0, 1, 2, \dots\}$).  If the destination site is
greater than or equal to the position of the
$(k{-}1)$-th particle at time $t{-}1$, the jump is
blocked and particle $k$ arrives at the site immediately
to the left of particle $k{-}1$.  The evolution is given
by
\[
Y_1(t) = Y_1(t{-}1) + \xi(t,1),
\]
\[
Y_k(t) = \min\{Y_k(t{-}1) + \xi(t,k),
Y_{k-1}(t{-}1) - 1\},
\qquad k = 2, \dots, N,
\]
where $\{\xi(t,k)\}_{t \ge 1, 1 \le k \le N}$ are
independent $\operatorname{Geom}[p]$ random variables and
$Y_0(t) = \infty$.  Write $q = 1 - p$.

\paragraph{\textbf{Fredholm determinant formula.}}
Fix a time $t \ge 0$ and an initial configuration
$Y(0) = y = (y_1 > y_2 > \cdots > y_N)$.  Fix
$m \in \{1, \dots, N{-}1\}$, particle labels
$\mathbf{n} = (n_1, n_2, \dots, n_m)$ with
$1 \le n_1 < n_2 < \cdots < n_m < N$, and spatial
thresholds $\mathbf{a} = (a_1, \dots, a_m) \in
\Z^m$.  The multipoint distribution is given by
\[
\PP_y\Bigl(\bigcap_{i=1}^m
\{Y_{n_i}(t) > a_i\}\Bigr)
= \det(I - \bar{\chi}_{\mathbf{a}} K_t
\bar{\chi}_{\mathbf{a}})_{\ell^2(\{n_1, \dots, n_m\}
\times \Z)},
\]
where $\bar{\chi}_{\mathbf{a}}(n_i, x) =
\mathbf{1}_{x \le a_i}$.  The correlation kernel $K_t$
has block entries
\[
K_t(n_i, x_i; n_j, x_j)
= -Q^{n_j - n_i}(x_i, x_j)\mathbf{1}_{n_i < n_j}
+ \bigl(\mathcal{S}_{-t,-n_i}^*
\overline{\mathcal{S}}_{-t,n_j}^{\operatorname{epi}(y)}
\bigr)(x_i, x_j),
\]
where $A^*$ denotes the transpose kernel,
$A^*(x, y) = A(y, x)$.

Fix an auxiliary parameter $\theta \in (0, 1)$, which
enters the Fredholm determinant representation but does
not appear in the final multipoint equation
\eqref{eq:H7-multipoint}, and let
$\alpha = (1{-}\theta)\theta^{-1}$.  The operators
defining the kernel are as follows.

\medskip
\noindent\textit{The transition matrix $Q$.}
The matrix $Q$ is the transition matrix of a random walk on
$\Z$ taking $\operatorname{Geom}[1{-}\theta]$ steps
strictly to the left:
\[
Q(z_1, z_2) = (1{-}\theta)\theta^{z_1 - z_2 - 1}
\mathbf{1}_{z_1 > z_2},
\]
with its $n$-th power admitting the integral
representation
\[
Q^n(z_1, z_2) = \frac{\alpha^n}{2\pi\mathrm{i}}
\oint_{\gamma_\rho} \diff w
\frac{\theta^{z_1 - z_2}}{w^{z_1 - z_2 - n + 1}}
\Bigl(\frac{1}{1 - w}\Bigr)^n.
\]

\medskip
\noindent\textit{The scattering operators
$\mathcal{S}_{-t,-n}$ and
$\overline{\mathcal{S}}_{-t,n}$.}\footnote{The kernel
operators $\mathcal{S}_{-t,-n}$ and
$\overline{\mathcal{S}}_{-t,n}$ are unrelated to the
lattice shifts $\mathcal{S}_k$ of the diamond framework
(\S\ref{sec:linear-problem}); the notation follows
Matetski--Remenik~\cite{MatetskiRemenik2023}.}
These operators are built from the probability generating
function $\varphi(w) = p/(1{-}qw)$ of the geometric jump
distribution.  They are defined by the contour integrals
\[
\mathcal{S}_{-t,-n}(z_1, z_2)
= \frac{\alpha^{-n+1}}{2\pi\mathrm{i}}
\oint_{\gamma_\rho} \frac{\diff w}{w}
\frac{\theta^{z_2 - z_1}}{w^{z_2 - z_1 + n}}
(1{-}w)^n \Bigl(\frac{p}{1{-}qw}\Bigr)^{\!t},
\]
\[
\overline{\mathcal{S}}_{-t,n}(z_1, z_2)
= \frac{\alpha^{n-1}}{2\pi\mathrm{i}}
\oint_{\gamma_\delta} \frac{\diff w}{w}
\frac{(1{-}w)^{z_2 - z_1 + n - 1}}
{\theta^{z_2 - z_1} w^{n-1}}
\Bigl(\frac{p}{1{-}q(1{-}w)}\Bigr)^{\!-t},
\]
where $\gamma_\rho$ is a positively oriented circle of
radius $\rho \in (0, 1)$ around the origin, and
$\gamma_\delta$ is a sufficiently small positively
oriented circle around the origin, chosen so that the only
singularity enclosed is the pole at $w=0$.

\medskip
\noindent\textit{The epigraph operator
$\overline{\mathcal{S}}_{-t,n}^{\operatorname{epi}(y)}$.}
This is the hitting probability operator defined in terms
of the initial data:
\[
\overline{\mathcal{S}}_{-t,n}^{\operatorname{epi}(y)}
(z_1, z_2)
= \E_{W_0 = z_1}\bigl[
\overline{\mathcal{S}}_{-t,n-\tau}(W_\tau, z_2)
\mathbf{1}_{\tau < n}\bigr],
\]
where $(W_\ell)_{\ell \ge 0}$ is the random walk with
transition matrix $Q$, and
$\tau = \min\{\ell \in \{0, \dots, N{-}1\} :
W_\ell > y_{\ell+1}\}$
is the hitting time of the strict epigraph of the initial
data by the random walk, with the convention
$\tau = \infty$ if the set is empty.

\paragraph{\textbf{References.}}
The right geometric jump model belongs to the family of four
basic discrete-time TASEP variants whose determinantal
transition kernels were obtained by
Dieker--Warren~\cite{DiekerWarren2008}.  The Fredholm
determinant formula used below is the Matetski--Remenik
formula~\cite{MatetskiRemenik2023} specialized to
$\kappa = 0$ and $\varphi(w) = p/(1{-}qw)$.

\subsection{Kernel reformulation}
\label{sec:RGP-kernel-reform}

We rewrite the extended-kernel Fredholm determinant in
the form of the product graph construction of
\S\ref{sec:product-graph}.  The functions $\psi, \phi$
below provide the seed data, while the lower-triangular
operator $B_{\mathbf{a}, \mathbf{n}}$ will be verified to satisfy the
admissibility conditions of Definition~\ref{def:admissible-prop}.

\paragraph{\textbf{Seed data.}}
Define
\[
\psi_{t,r,n}(x) \defeq \mathcal{S}_{-t,-n}^*(r,x)
= \mathcal{S}_{-t,-n}(x,r),
\qquad
\phi_{t,r,n}(x) \defeq
\overline{\mathcal{S}}_{-t,n}^{\operatorname{epi}(y)}(x,r),
\]
where the dependence of $\phi$ on the initial
condition $y$ is suppressed.

\paragraph{\textbf{Propagator.}}
The propagator is an instance of the directed-path
propagator (Definition~\ref{def:directed-path-propagator})
with transition kernels
$Q_{\ell,\ell'} = Q^{n_{\ell'} - n_\ell}$ for
$\ell < \ell'$.  On the cutoff region
$r \le a_i$, $r' \le a_j$, the entries of
$I + B_{\mathbf{a}, \mathbf{n}}$ coincide with those of
$((I + \bar{\chi}_{\mathbf{a}} L
\bar{\chi}_{\mathbf{a}})^{-1})^{\!\top}$, where $L$
is the strictly upper-triangular block operator with
entries $L_{ij}(x, x') = Q^{n_j - n_i}(x, x')
\mathbf{1}_{i < j}$ (Remark~\ref{rem:Neumann-series}).

\begin{lemma}\label{lem:RGP-kernel}
We have
\[
\PP_y\Bigl(\bigcap_{i=1}^m
\{Y_{n_i}(t) > a_i\}\Bigr)
= \det(I - K_{t, \mathbf{a}, \mathbf{n}})_{\ell^2(\Z)},
\]
where
\[
K_{t, \mathbf{a}, \mathbf{n}}(x, x')
= \sum_{1 \le j \le i \le m}
  \sum_{r \le a_i} \sum_{r' \le a_j}
   \psi_{t, r, n_i}(x)
  \bigl[\delta_{ij}\delta_{r,r'}
  + [B_{\mathbf{a}, \mathbf{n}}]_{ij}(r, r')\bigr]
  \phi_{t, r', n_j}(x').
\]
This identifies $K_{t, \mathbf{a}, \mathbf{n}}$ with the
kernel~\eqref{eq:K-def} on the product graph
$\mathcal{G}^m$, with $\psi, \phi$ as the seed data and
$B_{\mathbf{a}, \mathbf{n}}$ as the propagator.
\end{lemma}

\begin{proof}
The argument is identical to that of
Lemma~\ref{lem:RBJ-kernel}: the block-triangular
decomposition, Sylvester's identity, and the transpose
identification produce the claimed single-space kernel.
The random walk $Q$ and the triangular-block estimate are
unchanged from \S\ref{sec:RBJ}; the estimate for $\psi$
requires modification because the generating function
$(p/(1{-}qw))^t$ is not a polynomial for $t \ge 1$.

The contour formula for $\psi$ gives
\[
\psi_{t,r,n}(x) = \alpha^{-n+1}\theta^{r-x}
[w^{r-x+n}](1{-}w)^n
\Bigl(\frac{p}{1{-}qw}\Bigr)^{\!t}.
\]
Write $k = r - x + n$.  The coefficient vanishes for
$k < 0$, so $\psi_{t,r,n}(x) = 0$ for $x > r + n$.
Taking the coefficient on a circle $|w| = \rho$ with
$\theta < \rho < 1$ and setting
$\lambda = \theta/\rho$ gives
\begin{equation}\label{eq:RGP-psi-tail}
|\psi_{t,r,n}(x)| \le C\lambda^k
\mathbf{1}_{x \le r + n}.
\end{equation}
Since $\lambda < 1$,
$\psi_{t,r,n} \in \ell^2(\Z)$.  The integrand depends
on $r$ and $x$ only through $r - x$, so the $\ell^2$
norm is independent of~$r$.

The contour estimate for $\phi$ follows the same argument
as in Lemma~\ref{lem:RBJ-kernel}.  Choose
$\beta \in (\theta, 1)$ and $\gamma_\delta$ so that
$\sup_{w \in \gamma_\delta}|\theta/(1{-}w)| < \beta$.
The factor $((p{+}qw)/p)^t$ in the
$\overline{\mathcal{S}}$ contour formula is bounded on
$\gamma_\delta$, so
$|\overline{\mathcal{S}}_{-t,k}(z,r)|
\le C\beta^{z-r}$ for $z > y_N$ and
$r \le a_{\max}$.  Since the $Q$-walk is the same as in
\S\ref{sec:RBJ}, the exponential moment
$\E[\beta^{-S_\ell}]$ is finite, where $S_\ell$ is the
cumulative displacement of the $Q$-walk after
$\ell$ steps.  The hitting representation gives
\begin{equation}\label{eq:RGP-phi-tail}
|\phi_{t,r,n}(x)| \le C\beta^{x-r}
\mathbf{1}_{x > y_N}.
\end{equation}

With these bounds, the same conjugation and Sylvester
reduction as in Lemma~\ref{lem:RBJ-kernel} apply.
\end{proof}

\subsection{Base graph and seed data}
\label{sec:RGP-base-graph}

We now identify the lattice structure that governs how
the seed data evolve under shifts of the parameters
$(t, a, n)$, and verify that $\psi, \phi$ satisfy the
linear problems of the diamond framework.

Let $\mathcal{V} \subseteq \Z^3$ be the lattice
generated by the shifts
\[
T = e^{\partial_a}, \qquad
S_1 = e^{\partial_n}, \qquad
S_2 = e^{-\partial_t + \partial_a},
\]
so that for $u = (t, a, n) \in \mathcal{V}$,
\[
Tu = (t, a{+}1, n), \qquad
S_1 u = (t, a, n{+}1), \qquad
S_2 u = (t{-}1, a{+}1, n).
\]
Take $E = \mathbb{F} = \R$ and define the
constant scalar edge weights
\begin{equation}\label{eq:RGP-edge-weights}
(c_1, c_2) = \bigl((1{-}\theta)^{-1},
p^{-1}\bigr),
\qquad
(\lambda_1, \lambda_2)
= \bigl((1{-}\theta)^{-1}\theta,
q\theta p^{-1}\bigr).
\end{equation}
Since $(c_k, \lambda_k)$ are constant scalars, the
diamond equations are satisfied trivially.

\begin{lemma}\label{lem:RGP-seed-linear}
Let $H = \ell^2(\Z)$.  The seed functions
$\psi_{t,r,n}$ and $\phi_{t,r,n}$ satisfy the linear
problem \eqref{eq:Linear-Psi-3} and the adjoint
linear problem \eqref{eq:Linear-Phi-Prob} on
$\mathcal{V}$ with edge weights
$(c_k, \lambda_k)$.  Explicitly, for each fixed
$r \in \Z$:
\begin{align}
\psi_{t,r,n+1}
&= (1{-}\theta)^{-1}\psi_{t,r+1,n}
   - (1{-}\theta)^{-1}\theta\psi_{t,r,n},
\label{eq:RGP-psi-S1} \\
\psi_{t-1,r+1,n}
&= p^{-1}\psi_{t,r+1,n}
   - q\theta p^{-1}\psi_{t,r,n},
\label{eq:RGP-psi-S2}
\end{align}
and
\begin{align}
\phi_{t,r+1,n}
&= (1{-}\theta)^{-1}\phi_{t,r,n+1}
   - (1{-}\theta)^{-1}\theta\phi_{t,r+1,n+1},
\label{eq:RGP-phi-S1} \\
\phi_{t,r+1,n}
&= p^{-1}\phi_{t-1,r+1,n}
   - q\theta p^{-1}\phi_{t-1,r+2,n}.
\label{eq:RGP-phi-S2}
\end{align}
\end{lemma}

\begin{proof}
We verify each identity in turn.

\medskip
\noindent\textit{Verification of
\eqref{eq:RGP-psi-S1}.}
The $S_1$-shift involves only $n$ and the
$Q$-power structure, which is the same as in
\S\ref{sec:RBJ}; the probability generating function
$(p/(1{-}qw))^t$ passes through unchanged.  The
computation is identical to the verification of
\eqref{eq:RBJ-psi-S1}: the factor
$(1{-}\theta)^{-1}(\theta/w - \theta)
= (1{-}w)/(\alpha w)$ shifts the integrand from $n$
to $n{+}1$.

\medskip
\noindent\textit{Verification of
\eqref{eq:RGP-psi-S2}.}
From the contour representation,
\[
\psi_{t,r,n}(x) = \frac{\alpha^{-n+1}}{2\pi\mathrm{i}}
  \oint_{\gamma_\rho} \frac{\diff w}{w}
  \frac{\theta^{r-x}}{w^{r-x+n}}
  (1{-}w)^n \Bigl(\frac{p}{1{-}qw}\Bigr)^{\!t}.
\]
The right-hand side of \eqref{eq:RGP-psi-S2} is
\begin{align*}
p^{-1}\psi_{t,r+1,n}(x)
- q\theta p^{-1}\psi_{t,r,n}(x)
= \frac{\alpha^{-n+1}}{2\pi\mathrm{i}}
  \oint_{\gamma_\rho} \frac{\diff w}{w}
  \frac{\theta^{r-x}}{w^{r-x+n}}
  (1{-}w)^n \Bigl(\frac{p}{1{-}qw}\Bigr)^{\!t}
  \cdot\frac{1}{p}\Bigl(\frac{\theta}{w}
  - q\theta\Bigr).
\end{align*}
The combined factor simplifies as
\[
\frac{1}{p}\Bigl(\frac{\theta}{w}
- q\theta\Bigr)
= \frac{\theta(1{-}qw)}{pw}.
\]
Since
$(p/(1{-}qw))^t (1{-}qw)/p
= (p/(1{-}qw))^{t-1}$,
and the $\theta/w$ shifts $r$ by $+1$, this gives
\[
= \frac{\alpha^{-n+1}}{2\pi\mathrm{i}}
  \oint_{\gamma_\rho} \frac{\diff w}{w}
  \frac{\theta^{(r+1)-x}}{w^{(r+1)-x+n}}
  (1{-}w)^n \Bigl(\frac{p}{1{-}qw}\Bigr)^{\!t-1}
= \psi_{t-1,r+1,n}(x).
\]

\medskip
\noindent\textit{Verification of
\eqref{eq:RGP-phi-S1}.}
The argument follows the proof of
\eqref{eq:RBJ-phi-S1} in \S\ref{sec:RBJ}.  The
$\overline{\mathcal{S}}$ recurrence in $n$,
\[
(1{-}\theta)^{-1}\bigl[
\overline{\mathcal{S}}_{-t,n}(z, r)
- \theta\overline{\mathcal{S}}_{-t,n}(z, r{+}1)
\bigr]
= \overline{\mathcal{S}}_{-t,n-1}(z, r{+}1),
\]
depends only on the $(1{-}w)$-power and $w$-power
structure, not on the choice of probability generating
function, and its verification is completely analogous.
We show that $\overline{\mathcal{S}}_{-t,0} \equiv 0$.
When $n = 0$, the measure $dw/w^n$ reduces to $dw$,
removing the pole at the origin.  The remaining integrand
is analytic inside $\gamma_\delta$ (since
$(p/(1{-}q(1{-}w)))^{-t} = ((p{+}qw)/p)^t$ is
polynomial in $w$), so the integral vanishes by
Cauchy's theorem.
The lifting to the epigraph version is identical to
\S\ref{sec:RBJ}.

\medskip
\noindent\textit{Verification of
\eqref{eq:RGP-phi-S2}.}
Since neither side involves a shift in $n$, the
indicator $\mathbf{1}_{\tau < n}$ passes through
unchanged, and the verification reduces to a
recurrence for the free operator.  We claim
\begin{equation}\label{eq:RGP-Sbar-recurrence-t}
p^{-1}
\overline{\mathcal{S}}_{-(t-1),n}(z, r)
- q\theta p^{-1}
\overline{\mathcal{S}}_{-(t-1),n}(z, r{+}1)
= \overline{\mathcal{S}}_{-t,n}(z, r).
\end{equation}
Indeed, the left-hand side equals
\begin{align*}
&\frac{\alpha^{n-1}}{2p\pi\mathrm{i}}
\oint_{\gamma_\delta} \frac{\diff w}{w}
\frac{(p/(1{-}q(1{-}w)))^{-(t-1)}
      (1{-}w)^{r-z+n-1}}
     {\theta^{r-z} w^{n-1}}
\bigl(1 - q(1{-}w)\bigr).
\end{align*}
Since
$1 - q(1{-}w) = p + qw$, and
$(p/(1{-}q(1{-}w)))^{-(t-1)}
(1{-}q(1{-}w))/p
= (p/(1{-}q(1{-}w)))^{-t}$,
the left-hand side becomes
$\overline{\mathcal{S}}_{-t,n}(z, r)$.
Evaluating \eqref{eq:RGP-Sbar-recurrence-t} at
$r{+}1$ and taking the expectation against the
hitting-time representation completes the proof.
\end{proof}

\subsection{Admissible propagators and corrected
edge weights}
\label{sec:RGP-propagators}

The product graph $\mathcal{G}^m$ has vertex set
$\mathcal{V}^m$ and is generated by the diagonal shifts
\[
\mathcal{T}u
= (t, \mathbf{a} + \mathbf{1}, \mathbf{n}),
\qquad
\mathcal{S}_1 u
= (t, \mathbf{a}, \mathbf{n} + \mathbf{1}),
\qquad
\mathcal{S}_2 u
= (t{-}1, \mathbf{a} + \mathbf{1}, \mathbf{n}).
\]
The block-diagonal edge weights are
$C_k(u) = c_k I_m$ and $D_k(u) = \lambda_k I_m$.

The propagator is an instance of the directed-path
propagator (Definition~\ref{def:directed-path-propagator})
with transition kernels
$Q_{\ell,\ell'} = Q^{n_{\ell'} - n_\ell}$ for
$\ell < \ell'$.  The standing hypotheses of
\S\ref{sec:directed-path} are satisfied: translation
invariance of $Q^n$ is inherited from the translation
invariance of~$Q$, and $\mathcal{T}$-invariance holds
because the transition kernels depend only on the
particle-label differences, which are unchanged
by~$\mathcal{T}$.

By Proposition~\ref{prop:constant-scalar},
$\mathcal{S}_k$-compatibility follows from the
constant-scalar edge weights once the propagator satisfies
the following invariance property.

\begin{lemma}[$\mathcal{S}_k$-invariance]
\label{lem:RGP-Sk-invariance}
For $k = 1, 2$ and all $x, y \in \mathcal{V}^m$,
\begin{equation}\label{eq:RGP-Sk-invariance}
B_{\mathcal{S}_k u}
(\mathcal{S}_k x, \mathcal{S}_k y)
= B_u(x, y).
\end{equation}
\end{lemma}

\begin{proof}
\noindent\textit{Case $k = 1$.}
The shift gives $\mathcal{S}_1 u =
(t, \mathbf{a}, \mathbf{n}{+}\mathbf{1})$: the
cutoffs $\mathbf{a}$ and the differences
$n_j - n_i$ are unchanged, and $\mathcal{S}_1$
does not shift $\mathbf{a}$, so the propagator entries
$B_u(x, y) = B_{\mathbf{a}, \mathbf{n}}(x, y)$ are unchanged.

\noindent\textit{Case $k = 2$.}
The shift gives $\mathcal{S}_2 u =
(t{-}1, \mathbf{a}{+}\mathbf{1}, \mathbf{n})$:
the cutoffs shift to $\mathbf{a}{+}\mathbf{1}$
and the label differences are unchanged.  Since
$\mathcal{S}_2$ acts on $\mathbf{a}$ identically
to $\mathcal{T}$, the invariance follows from the
same translation-invariance argument as
$\mathcal{T}$-covariance (the propagator is
$t$-independent).
\end{proof}

\begin{lemma}\label{lem:RGP-admissible}
The propagator $B_{\mathbf{a}, \mathbf{n}}$ is admissible for
$\mathcal{G}^m$.
\end{lemma}

\begin{proof}
By Proposition~\ref{prop:directed-path-T},
$\mathcal{T}$-covariance and $\mathcal{T}$-splitting hold.
The edge weights $(c_k, \lambda_k)$ are constant scalars
and the propagator satisfies $\mathcal{S}_k$-invariance
(Lemma~\ref{lem:RGP-Sk-invariance}), so
$\mathcal{S}_k$-compatibility follows from
Proposition~\ref{prop:constant-scalar}.

It remains to verify regularity.  The transition matrix
$Q$ and directed-path propagator are the same as in
\S\ref{sec:RBJ}, so the propagator entries decay
geometrically.  Together with
\eqref{eq:RGP-psi-tail} and
\eqref{eq:RGP-phi-tail}, these bounds prove absolute
convergence of the sums defining $\Psi$, $\Phi$,
and $K$ in \eqref{eq:Psi-def}--\eqref{eq:K-def}.
\end{proof}

\paragraph{\textbf{Corrected edge weights.}}
Since $\mathcal{S}_k$-invariance holds,
Proposition~\ref{prop:constant-scalar}(ii) gives
$P(\mathcal{S}_k u) = P(u)$ and
$\Lambda_k(u) = \lambda_k I_m = D_k(u)$.
The diamond data on $\mathcal{G}^m$ are therefore
$(C_k, \Lambda_k) = (c_k I_m, \lambda_k I_m)$.
Theorem~\ref{thm:product-graph} requires additionally
that the resolvent $(I - zK(u))^{-1}$ exist.  At $z = 1$,
the Fredholm determinant
$F_{t, \mathbf{a}, \mathbf{n}} \defeq
\det(I - K_{t, \mathbf{a}, \mathbf{n}})$ equals the
multipoint gap probability
$\PP_y\bigl(\bigcap_i
\{Y_{n_i}(t) > a_i\}\bigr)$.  Since $K$ is trace class,
$F \ne 0$ if and only if the resolvent exists.  Geometric
jumps are unbounded, but blocking propagates finite upper
bounds down the particle chain: $Y_k(t) \le h_k(t)$ where
\[
h_k(t) \defeq
\begin{cases}
+\infty, & k \le t, \\
y_{k-t} - t, & k > t.
\end{cases}
\]
For $k > t$, iterating the blocking constraint
$Y_k(s) \le Y_{k-1}(s{-}1) - 1$ gives
$Y_k(t) \le y_{k-t} - t$.  For $k \le t$, the
unboundedness of geometric jumps propagates through the
blocking chain in $t$ steps.  If each jump exceeds its
blocking gap (an event of positive probability), then
$Y_k(t) = h_k(t)$ for $k > t$.  The resolvent at $z = 1$
therefore exists precisely on
$\mathcal{R} \defeq \{(t, \mathbf{a}, \mathbf{n}) \in
\mathcal{V}^m : a_i < h_{n_i}(t)
\text{ for all } i\}$.
The mixed diamond at a vertex~$u$ requires the resolvent
at the four vertices $u$, $\mathcal{S}_1 u$,
$\mathcal{S}_2 u$, $\mathcal{S}_1\mathcal{S}_2 u$
entering~\eqref{eq:dressed-C-explicit}, together with
the $\mathcal{T}$-shifts $\mathcal{T}u$,
$\mathcal{T}\mathcal{S}_1 u$,
$\mathcal{T}\mathcal{S}_2 u$ from the proof
of~\eqref{eq:M-inverse}.
The stencil also requires $t \ge 1$, since
$\mathcal{S}_2$ shifts $t \mapsto t{-}1$.
Since $h_k$ is non-decreasing in $t$ and
$h_{k+1}(t) \le h_k(t) - 1$, the binding constraint is
$\mathcal{S}_1\mathcal{S}_2 u \in \mathcal{R}$, i.e.\
$a_i < h_{n_i+1}(t{-}1) - 1$ for all~$i$.

\subsection{Multipoint equation}
\label{sec:RGP-multipoint}

The dressed observable at $z = 1$,
\[
\mathcal{M}(u)
= I + \Phi(u)(I - K(u))^{-1}\Psi(u)
\in \End(\R^m),
\]
has dressed edge weights
\begin{equation}\label{eq:RGP-dressed-Ck}
\mathcal{M}_k(u)
= c_k\mathcal{M}(\mathcal{T}u)^{-1}
  \mathcal{M}(\mathcal{S}_k u).
\end{equation}
For readability, we write
$\mathcal{M}_{t, \mathbf{a}, \mathbf{n}} \defeq
\mathcal{M}(t, \mathbf{a}, \mathbf{n})$.

\begin{theorem}\label{thm:RGP-H7}
The dressed observable satisfies the mixed diamond
equation
\begin{equation}\label{eq:H7-multipoint}
\bigl(\mathcal{M}_{t, \mathbf{a}+\mathbf{1},
\mathbf{n}+\mathbf{1}}^{-1}
- q\mathcal{M}_{t-1, \mathbf{a}+\mathbf{2},
\mathbf{n}}^{-1}\bigr)
\mathcal{M}_{t-1, \mathbf{a}+\mathbf{1},
\mathbf{n}+\mathbf{1}}
- \mathcal{M}_{t, \mathbf{a}+\mathbf{1},
\mathbf{n}}^{-1}
\bigl(\mathcal{M}_{t-1, \mathbf{a}+\mathbf{1},
\mathbf{n}}
- q\mathcal{M}_{t, \mathbf{a},
\mathbf{n}+\mathbf{1}}\bigr) = 0,
\end{equation}
at every $(t, \mathbf{a}, \mathbf{n})$ with $t \ge 1$ and
$a_i < h_{n_i+1}(t{-}1) - 1$ for all~$i$.
\end{theorem}

\begin{proof}
By Theorem~\ref{thm:product-graph}, the pair
$(\mathcal{M}_k, \Lambda_k)$ satisfies the diamond
equations on~$\mathcal{R}$.
Since $\Lambda_k = \lambda_k I_m$, the scalar factors
commute with $\mathcal{M}_k$, and the mixed diamond
equation reduces to
\begin{equation}\label{eq:RGP-mixed-scalar-Lambda}
\lambda_2\bigl[
\mathcal{M}_1(u) - \mathcal{M}_1(\mathcal{S}_2 u)
\bigr]
+ \lambda_1\bigl[
\mathcal{M}_2(\mathcal{S}_1 u)
- \mathcal{M}_2(u)
\bigr] = 0.
\end{equation}
Substituting \eqref{eq:RGP-dressed-Ck} and the explicit
shifts $\mathcal{T}u =
(t, \mathbf{a}{+}\mathbf{1}, \mathbf{n})$,
$\mathcal{S}_1 u =
(t, \mathbf{a}, \mathbf{n}{+}\mathbf{1})$,
$\mathcal{S}_2 u =
(t{-}1, \mathbf{a}{+}\mathbf{1}, \mathbf{n})$,
every term carries a common factor of
$\theta/(p(1{-}\theta))$.  Dividing by this factor
gives
\begin{align*}
q\bigl[
&\mathcal{M}_{t, \mathbf{a}+\mathbf{1},
\mathbf{n}}^{-1}
\mathcal{M}_{t, \mathbf{a},
\mathbf{n}+\mathbf{1}}
- \mathcal{M}_{t-1,
\mathbf{a}+\mathbf{2},
\mathbf{n}}^{-1}
\mathcal{M}_{t-1, \mathbf{a}+\mathbf{1},
\mathbf{n}+\mathbf{1}}
\bigr] \\
+ \bigl[
&\mathcal{M}_{t, \mathbf{a}+\mathbf{1},
\mathbf{n}+\mathbf{1}}^{-1}
\mathcal{M}_{t-1, \mathbf{a}+\mathbf{1},
\mathbf{n}+\mathbf{1}}
- \mathcal{M}_{t, \mathbf{a}+\mathbf{1},
\mathbf{n}}^{-1}
\mathcal{M}_{t-1, \mathbf{a}+\mathbf{1},
\mathbf{n}}
\bigr] = 0.
\end{align*}
The middle two terms share the right factor
$\mathcal{M}_{t-1, \mathbf{a}+\mathbf{1},
\mathbf{n}+\mathbf{1}}$, and the
first and last share the left factor
$\mathcal{M}_{t, \mathbf{a}+\mathbf{1},
\mathbf{n}}^{-1}$.  Factoring yields
\eqref{eq:H7-multipoint}.
\end{proof}

For $m = 1$, the dressed observable $\mathcal{M}_{t,a,n}$
is scalar.

\begin{corollary}\label{cor:RGP-scalar}
The one-point Fredholm determinant
$F_{t, a, n} = \det(I - K_{t, a, n})$ satisfies
\begin{equation}\label{eq:RGP-scalar-HM}
q F_{t, a, n+1} F_{t-1, a+2, n}
- F_{t-1, a+1, n} F_{t, a+1, n+1}
+ p F_{t, a+1, n} F_{t-1, a+1, n+1}
= 0,
\end{equation}
at every $(t, a, n)$ with $t \ge 1$ and
$a < h_{n+1}(t{-}1) - 1$.
\end{corollary}

\begin{proof}
Proposition~\ref{prop:scalar-HM} applies at $z = 1$: the
boundary condition $F(T^\ell v) \to 1$ as
$\ell \to -\infty$ holds because the thresholds decrease
to $-\infty$ and the kernel vanishes.  The coefficients are
$\alpha_{12} = c_1 \lambda_2 =
q\theta/(p(1{-}\theta))$
and
$\alpha_{21} = c_2 \lambda_1 =
\theta/(p(1{-}\theta))$.
Substituting the lattice points into
\eqref{eq:HM-variable} and computing
$\alpha_{12} - \alpha_{21}
= (q{-}1)\theta/(p(1{-}\theta))
= -\theta/(1{-}\theta)$
(using $p + q = 1$), then dividing by
$\theta/(p(1{-}\theta))$, gives
\eqref{eq:RGP-scalar-HM}.
\end{proof}

\begin{remark}[Structure of the multipoint equation]
\label{rem:RGP-structure}
Equation \eqref{eq:H7-multipoint} is an
$m \times m$ matrix equation coupling the dressed
observables at six lattice points, with coefficient $q$
the failure parameter of the geometric jump distribution.  For $m \ge 2$, the equation is intrinsically
noncommutative: the matrix inverses and products do not
simplify to a scalar relation.
\end{remark}

\section{Left Geometric Jumps}\label{sec:LG}

\paragraph{\textbf{System description.}}
The discrete-time totally asymmetric simple
exclusion process with left geometric jumps, pushing
interaction, and sequential update consists of $N$
particles on
$\Z$.  The configuration at time
$t \in \Z_{\ge 0}$ is denoted by
\[
Y(t) = (Y_1(t) > Y_2(t) > \cdots > Y_N(t)),
\]
where $Y_k(t)$ is the position of the $k$-th particle
from the right.

At each discrete time step $t \ge 1$, the particles are
updated sequentially from right to left (i.e., for
$k = 1, 2, \dots, N$).  Each particle makes a jump to
the left with a size sampled from a
$\mathrm{Geom}[p]$ distribution (supported on
$\{0, 1, 2, \dots\}$), pushing all particles in its
way so that the strict ordering is preserved.  The
evolution equations are:
\[
Y_1(t) = Y_1(t{-}1) - \xi(t,1),
\]
\[
Y_k(t) = \min\{Y_k(t{-}1),
Y_{k-1}(t) - 1\} - \xi(t,k),
\qquad k = 2, \dots, N,
\]
where $\{\xi(t,k)\}_{t \ge 1, 1 \le k \le N}$ are
independent $\mathrm{Geom}[p]$ random variables and
$Y_0(t) = \infty$.

\paragraph{\textbf{Fredholm determinant formula.}}
Fix a time $t \ge 0$ and an initial configuration
$Y(0) = y = (y_1 > y_2 > \cdots > y_N)$.  Fix
$m \in \{1, \dots, N{-}1\}$, particle labels
$\mathbf{n} = (n_1, n_2, \dots, n_m)$ with
$1 \le n_1 < n_2 < \cdots < n_m < N$, and spatial
thresholds $\mathbf{a} = (a_1, \dots, a_m) \in
\Z^m$.  The multipoint joint cumulative distribution
of the particle positions is given by
\[
\PP_y\Bigl(\bigcap_{i=1}^m
\{Y_{n_i}(t) > a_i\}\Bigr)
= \det(I - \bar{\chi}_{\mathbf{a}} K_t
\bar{\chi}_{\mathbf{a}})_{\ell^2(\{n_1, \dots, n_m\}
\times \Z)},
\]
where $\bar{\chi}_{\mathbf{a}}(n_i, x) =
\mathbf{1}_{x \le a_i}$.  The correlation kernel $K_t$
has block entries
\[
K_t(n_i, x_i; n_j, x_j)
= -Q^{n_j - n_i}(x_i, x_j) \mathbf{1}_{n_i < n_j}
+ \bigl(\mathcal{S}_{-t,-n_i}^*
\overline{\mathcal{S}}_{-t,n_j}^{\operatorname{epi}(y)}\bigr)
(x_i, x_j),
\]
where $A^*$ denotes the transpose kernel,
$A^*(x, y) = A(y, x)$.

Fix an auxiliary parameter $\theta \in (q, 1)$ (where
$q = 1 - p$), which enters the Fredholm determinant
representation but does not appear in the final
multipoint equation \eqref{eq:H8-multipoint}, and let
$\alpha = (1{-}\theta)\theta^{-1}$.  The operators
defining the kernel are as follows.

\medskip
\noindent\textit{The transition matrix $Q$.}
The matrix $Q$ is the transition matrix of a random walk
on $\Z$ taking $\mathrm{Geom}[1{-}\theta]$ steps
strictly to the left:
\[
Q(z_1, z_2) = (1{-}\theta) \theta^{z_1 - z_2 - 1}
\mathbf{1}_{z_1 > z_2},
\]
with its $n$-th power admitting the integral
representation
\[
Q^n(z_1, z_2) = \frac{\alpha^n}{2\pi\mathrm{i}}
\oint_{\gamma_\rho} \diff w
\frac{\theta^{z_1 - z_2}}{w^{z_1 - z_2 - n + 1}}
\Bigl(\frac{1}{1 - w}\Bigr)^n.
\]

\medskip
\noindent\textit{The scattering operators
$\mathcal{S}_{-t,-n}$ and
$\overline{\mathcal{S}}_{-t,n}$.}\footnote{The kernel
operators $\mathcal{S}_{-t,-n}$ and
$\overline{\mathcal{S}}_{-t,n}$ are unrelated to the
lattice shifts $\mathcal{S}_k$ of the diamond framework
(\S\ref{sec:linear-problem}); the notation follows
Matetski--Remenik~\cite{MatetskiRemenik2023}.}
These operators encode the left
geometric jumps via the probability generating function
$\varphi(w) = p/(1{-}q/w) = pw/(w{-}q)$.  They are
defined by the contour integrals
\[
\mathcal{S}_{-t,-n}(z_1, z_2)
= \frac{\alpha^{-n+1}}{2\pi\mathrm{i}}
\oint_{\gamma_\rho} \frac{\diff w}{w}
\frac{\theta^{z_2 - z_1}}{w^{z_2 - z_1 + n}}
(1{-}w)^n
\Bigl(\frac{pw}{w{-}q}\Bigr)^{\!t},
\]
\[
\overline{\mathcal{S}}_{-t,n}(z_1, z_2)
= \frac{\alpha^{n-1}}{2\pi\mathrm{i}}
\oint_{\gamma_\delta} \frac{\diff w}{w}
\frac{(1{-}w)^{z_2 - z_1 + n - 1}}
{\theta^{z_2 - z_1} w^{n-1}}
\Bigl(\frac{p(1{-}w)}{1{-}w{-}q}\Bigr)^{\!-t},
\]
where $\gamma_\rho$ is a positively oriented circle of
radius $\rho \in (q, 1)$ around the origin, and
$\gamma_\delta$ is a sufficiently small positively
oriented circle around the origin, chosen so that the
only singularity enclosed is the pole at $w = 0$.

\medskip
\noindent\textit{The epigraph operator
$\overline{\mathcal{S}}_{-t,n}^{\operatorname{epi}(y)}$.}
This is the hitting probability operator defined in terms
of the initial data:
\[
\overline{\mathcal{S}}_{-t,n}^{\operatorname{epi}(y)}
(z_1, z_2)
= \E_{W_0 = z_1}\bigl[
\overline{\mathcal{S}}_{-t,n-\tau}(W_\tau, z_2)
\mathbf{1}_{\tau < n}\bigr],
\]
where $(W_\ell)_{\ell \ge 0}$ is the random walk with
transition matrix $Q$, and
$\tau = \min\{\ell \in \{0, \dots, N{-}1\} :
W_\ell > y_{\ell+1}\}$
is the hitting time of the strict epigraph of the
initial data by the random walk, with the convention
$\tau = \infty$ if the set is empty.

\paragraph{\textbf{References.}}
The model belongs to the family of four basic
discrete-time TASEP variants whose determinantal
transition kernels were obtained by
Dieker--Warren~\cite{DiekerWarren2008}.  The Fredholm
determinant formula used below is the Matetski--Remenik
formula~\cite{MatetskiRemenik2023} specialized to
$\kappa = 0$ and $\varphi(w) = pw/(w{-}q)$.

\subsection{Kernel reformulation}
\label{sec:LG-kernel-reform}

We rewrite the extended-kernel Fredholm determinant in
the form of the product graph construction of
\S\ref{sec:product-graph}.  The functions $\psi, \phi$
below provide the seed data, while the
lower-triangular operator $B_{\mathbf{a},\mathbf{n}}$
will be verified to satisfy the admissibility
conditions of Definition~\ref{def:admissible-prop}.

\paragraph{\textbf{Seed data.}}
Define
\[
\psi_{t,r,n}(x) \defeq \mathcal{S}_{-t,-n}^*(r,x)
= \mathcal{S}_{-t,-n}(x,r),
\qquad
\phi_{t,r,n}(x) \defeq
\overline{\mathcal{S}}_{-t,n}^{\operatorname{epi}(y)}
(x,r).
\]

\paragraph{\textbf{Propagator.}}
The propagator is an instance of the directed-path
propagator (Definition~\ref{def:directed-path-propagator})
with transition kernels
$Q_{\ell,\ell'} = Q^{n_{\ell'} - n_\ell}$ for
$\ell < \ell'$.  On the cutoff region $r \le a_i$,
$r' \le a_j$, the entries of
$I + B_{\mathbf{a},\mathbf{n}}$ coincide with those of
$((I + \bar{\chi}_{\mathbf{a}} L
\bar{\chi}_{\mathbf{a}})^{-1})^{\!\top}$, where $L$
is the strictly upper-triangular block operator with
entries $L_{ij}(x, x') = Q^{n_j - n_i}(x, x')
\mathbf{1}_{i < j}$ (Remark~\ref{rem:Neumann-series}).
The propagator depends on $\mathbf{a}$ and $\mathbf{n}$
but not on the time coordinate~$t$, which is accordingly
suppressed from the notation.

\begin{lemma}\label{lem:LG-kernel}
We have
\[
\PP_y\Bigl(\bigcap_{i=1}^m
\{Y_{n_i}(t) > a_i\}\Bigr)
= \det(I - K_{t, \mathbf{a}, \mathbf{n}})_{\ell^2(\Z)},
\]
where
\[
K_{t, \mathbf{a}, \mathbf{n}}(x, x')
= \sum_{1 \le j \le i \le m}
  \sum_{r \le a_i} \sum_{r' \le a_j}
   \psi_{t, r, n_i}(x)
  \bigl[\delta_{ij}\delta_{r,r'}
  + [B_{\mathbf{a}, \mathbf{n}}]_{ij}(r, r')\bigr]
  \phi_{t, r', n_j}(x').
\]
This identifies $K_{t, \mathbf{a}, \mathbf{n}}$ with the
kernel~\eqref{eq:K-def} on the product graph
$\mathcal{G}^m$, with $\psi, \phi$ as the seed data and
$B_{\mathbf{a}, \mathbf{n}}$ as the propagator.
\end{lemma}

\begin{proof}
The argument is identical to that of
Lemma~\ref{lem:RBJ-kernel}: the block-triangular
decomposition, Sylvester's identity, and the transpose
identification produce the claimed single-space kernel.
The random walk $Q$ and the triangular-block estimate are
unchanged from \S\ref{sec:RBJ}; the estimate for $\psi$
requires modification because the generating function
$(pw/(w{-}q))^t$ has a pole at $w = q$ inside the
contour $\gamma_\rho$.

The contour formula for $\psi$ gives
\[
\psi_{t,r,n}(x) = \frac{\alpha^{-n+1}}{2\pi\mathrm{i}}
  \oint_{\gamma_\rho} \frac{\diff w}{w}
  \frac{\theta^{r-x}}{w^{r-x+n}}
  (1{-}w)^n \Bigl(\frac{pw}{w{-}q}\Bigr)^{\!t}.
\]
Expanding the contour to $|w| \to \infty$, the integrand
decays as $|w|^{x-r-1}$, so $\psi_{t,r,n}(x) = 0$ for
$x < r$.  For $x \ge r$, taking the Cauchy estimate on
$|w| = \rho$ with $q < \rho < \theta$ gives
\begin{equation}\label{eq:LG-psi-tail}
|\psi_{t,r,n}(x)| \le C(\rho/\theta)^{x-r}
\mathbf{1}_{x \ge r}.
\end{equation}
Since $\rho/\theta < 1$,
$\psi_{t,r,n} \in \ell^2(\Z)$.

The contour estimate for $\phi$ follows the same argument
as in Lemma~\ref{lem:RBJ-kernel}.  Choose
$\beta \in (\theta, 1)$ and $\gamma_\delta$ so that
$\sup_{w \in \gamma_\delta}|\theta/(1{-}w)| < \beta$.
The factor $((p{-}w)/(p(1{-}w)))^t$ in the
$\overline{\mathcal{S}}$ contour formula is bounded on
$\gamma_\delta$, so
$|\overline{\mathcal{S}}_{-t,k}(z,r)|
\le C\beta^{z-r}$ for $z > y_N$ and
$r \le a_{\max}$.  Since the $Q$-walk is the same as in
\S\ref{sec:RBJ}, the exponential moment
$\E[\beta^{-S_\ell}]$ is finite, where $S_\ell$ is the
cumulative displacement of the $Q$-walk after
$\ell$ steps.  The hitting representation gives
\begin{equation}\label{eq:LG-phi-tail}
|\phi_{t,r,n}(x)| \le C\beta^{x-r}
\mathbf{1}_{x > y_N}.
\end{equation}

With these bounds, the same conjugation and Sylvester
reduction as in Lemma~\ref{lem:RBJ-kernel} apply.
\end{proof}

\subsection{Base graph and seed data}
\label{sec:LG-base-graph}

We now identify the lattice structure that governs how
the seed data evolve under shifts of the parameters
$(t, a, n)$, and verify that $\psi, \phi$ satisfy the
linear problems of the diamond framework.

Let $\mathcal{V} \subseteq \Z^3$ be the lattice
generated by the shifts
\[
T = e^{\partial_a}, \qquad
S_1 = e^{\partial_n}, \qquad
S_2 = e^{-\partial_t},
\]
so that for $u = (t, a, n) \in \mathcal{V}$,
\[
Tu = (t, a{+}1, n), \qquad
S_1 u = (t, a, n{+}1), \qquad
S_2 u = (t{-}1, a, n).
\]
Take $E = \mathbb{F} = \R$ and define the
constant scalar edge weights
\begin{equation}\label{eq:LG-edge-weights}
(c_1, c_2) = \bigl((1{-}\theta)^{-1},
-q(p\theta)^{-1}\bigr),
\qquad
(\lambda_1, \lambda_2)
= \bigl((1{-}\theta)^{-1}\theta, -p^{-1}\bigr).
\end{equation}
Since $(c_k, \lambda_k)$ are constant scalars, the
diamond equations
\eqref{eq:diamond-C-ij}--\eqref{eq:diamond-mixed-ij}
are satisfied trivially: all three reduce to
commutativity in $\R$.

\begin{lemma}\label{lem:LGJ-seed-linear}
Let $H = \ell^2(\Z)$.  The seed functions
$\psi_{t,r,n}$ and $\phi_{t,r,n}$ satisfy the linear
problem \eqref{eq:Linear-Psi-3} and the adjoint
linear problem \eqref{eq:Linear-Phi-Prob} on
$\mathcal{V}$ with edge weights
$(c_k, \lambda_k)$.  Explicitly, for each fixed
$r \in \Z$:
\begin{align}
\psi_{t,r,n+1}
&= (1{-}\theta)^{-1} \psi_{t,r+1,n}
   - (1{-}\theta)^{-1}\theta \psi_{t,r,n},
\label{eq:LG-psi-S1} \\
\psi_{t-1,r,n}
&= -q(p\theta)^{-1} \psi_{t,r+1,n}
   + p^{-1} \psi_{t,r,n},
\label{eq:LG-psi-S2}
\end{align}
and
\begin{align}
\phi_{t,r+1,n}
&= (1{-}\theta)^{-1} \phi_{t,r,n+1}
   - (1{-}\theta)^{-1}\theta \phi_{t,r+1,n+1},
\label{eq:LG-phi-S1} \\
\phi_{t,r+1,n}
&= -q(p\theta)^{-1} \phi_{t-1,r,n}
   + p^{-1} \phi_{t-1,r+1,n}.
\label{eq:LG-phi-S2}
\end{align}
\end{lemma}

\begin{proof}
Each identity is verified in turn.

\medskip
\noindent\textit{Verification of
\eqref{eq:LG-psi-S1}.}
The computation is identical to the verification of
\eqref{eq:RBJ-psi-S1}: the $S_1$-shift involves only
$n$ and the $(1{-}w)^n/w^{r-x+n}$ structure, while the
probability generating function $(pw/(w{-}q))^t$
depends only on $t$ and is unchanged.

\medskip
\noindent\textit{Verification of
\eqref{eq:LG-psi-S2}.}
From the contour representation,
\[
\psi_{t,r,n}(x) = \frac{\alpha^{-n+1}}{2\pi\mathrm{i}}
  \oint_{\gamma_\rho} \frac{\diff w}{w}
  \frac{\theta^{r-x}}{w^{r-x+n}}
  (1{-}w)^n \Bigl(\frac{pw}{w{-}q}\Bigr)^{\!t}.
\]
The right-hand side of \eqref{eq:LG-psi-S2} is
\begin{align*}
-q(p\theta)^{-1} \psi_{t,r+1,n}(x)
+ p^{-1} \psi_{t,r,n}(x)
= \frac{\alpha^{-n+1}}{2\pi\mathrm{i}}
  \oint_{\gamma_\rho} \frac{\diff w}{w}
  \frac{\theta^{r-x}}{w^{r-x+n}}
  (1{-}w)^n \Bigl(\frac{pw}{w{-}q}\Bigr)^{\!t}
  \Bigl(\frac{-q}{p\theta}
  \cdot\frac{\theta}{w}
  + \frac{1}{p}\Bigr).
\end{align*}
The parenthetical factor simplifies as
\[
\frac{-q}{pw} + \frac{1}{p}
= \frac{w - q}{pw}.
\]
Since
$(pw/(w{-}q))^t (w{-}q)/(pw)
= (pw/(w{-}q))^{t-1}$,
this gives $\psi_{t-1,r,n}(x)$.

\medskip
\noindent\textit{Verification of
\eqref{eq:LG-phi-S1}.}
The argument is completely analogous to the proof of
\eqref{eq:RBJ-phi-S1} in \S\ref{sec:RBJ}.  The
$\overline{\mathcal{S}}$ recurrence in $n$ depends
only on the $(1{-}w)$-power and $w$-power structure
and is identical.

We show that $\overline{\mathcal{S}}_{-t,0} \equiv 0$.
When $n = 0$, the measure $dw/w^n$ reduces to $dw$,
removing the pole at the origin.  The remaining
integrand is analytic inside $\gamma_\delta$ (since
$(p(1{-}w)/(1{-}w{-}q))^{-t} = ((p{-}w)/(p(1{-}w)))^t$
is analytic in a neighbourhood of $w = 0$), so the integral
vanishes by Cauchy's theorem.

\medskip
\noindent\textit{Verification of
\eqref{eq:LG-phi-S2}.}
Since neither side involves a shift in $n$, the
indicator $\mathbf{1}_{\tau < n}$ passes through
unchanged.  We claim
\begin{equation}\label{eq:LG-Sbar-recurrence-t}
-q(p\theta)^{-1}
\overline{\mathcal{S}}_{-(t-1),n}(z, r)
+ p^{-1}
\overline{\mathcal{S}}_{-(t-1),n}(z, r{+}1)
= \overline{\mathcal{S}}_{-t,n}(z, r{+}1).
\end{equation}
Writing
$(p(1{-}w)/(1{-}w{-}q))^{-t}
= ((p{-}w)/(p(1{-}w)))^t$,
the left-hand side of
\eqref{eq:LG-Sbar-recurrence-t} combines the
integrands with factor
\[
c_2 - \lambda_2 (1{-}w)/\theta
= -\frac{q}{p\theta}
  + \frac{1{-}w}{p\theta}
= \frac{p - w}{p\theta},
\]
where we used $1 - w - q = p - w$.  Since
\[
\frac{p-w}{p\theta} = \frac{p-w}{p(1{-}w)} \cdot \frac{1{-}w}{\theta},
\]
this raises the exponent from $t{-}1$ to $t$, recovering
$\overline{\mathcal{S}}_{-t,n}(z, r{+}1)$.
Taking the expectation against the hitting-time
representation completes the proof.
\end{proof}

\subsection{Admissible propagators and corrected
edge weights}
\label{sec:LG-propagators}

The product graph $\mathcal{G}^m$ has vertex set
$\mathcal{V}^m$ and is generated by the diagonal shifts
\[
\mathcal{T}u
= (t, \mathbf{a} + \mathbf{1}, \mathbf{n}),
\qquad
\mathcal{S}_1 u
= (t, \mathbf{a}, \mathbf{n} + \mathbf{1}),
\qquad
\mathcal{S}_2 u
= (t{-}1, \mathbf{a}, \mathbf{n}).
\]
The block-diagonal edge weights are
$C_k(u) = c_k I_m$ and $D_k(u) = \lambda_k I_m$
(the uncorrected edge weights of
Definition~\ref{def:admissible-prop}).

The propagator is an instance of the directed-path
propagator (Definition~\ref{def:directed-path-propagator})
with transition kernels
$Q_{\ell,\ell'} = Q^{n_{\ell'} - n_\ell}$ for
$\ell < \ell'$.  The standing hypotheses of
\S\ref{sec:directed-path} are satisfied: translation
invariance of $Q^n$ is inherited from the translation
invariance of~$Q$, and $\mathcal{T}$-invariance holds
because the transition kernels depend only on the
particle-label differences, which are unchanged
by~$\mathcal{T}$.

By Proposition~\ref{prop:constant-scalar},
$\mathcal{S}_k$-compatibility follows from the
constant-scalar edge weights once the propagator satisfies
the following invariance property.

\begin{lemma}[$\mathcal{S}_k$-invariance]
\label{lem:LG-Sk-invariance}
For $k = 1, 2$ and all $x, y \in \mathcal{V}^m$,
$B_{\mathcal{S}_k u}
(\mathcal{S}_k x, \mathcal{S}_k y)
= B_u(x, y)$.
\end{lemma}

\begin{proof}
Neither $\mathcal{S}_1$ nor $\mathcal{S}_2$ shifts
$\mathbf{a}$, and both preserve the label differences
$n_j - n_i$.  Since the propagator depends only on
$\mathbf{a}$ and the label differences and is
independent of~$t$, invariance under both shifts is
immediate.
\end{proof}

\begin{lemma}\label{lem:LG-admissible}
The propagator $B_{\mathbf{a},\mathbf{n}}$ is an
admissible propagator for the product graph
$\mathcal{G}^m$.
\end{lemma}

\begin{proof}
By Proposition~\ref{prop:directed-path-T},
$\mathcal{T}$-covariance and $\mathcal{T}$-splitting hold.
The edge weights $(c_k, \lambda_k)$ are constant scalars
and the propagator satisfies $\mathcal{S}_k$-invariance
(Lemma~\ref{lem:LG-Sk-invariance}), so
$\mathcal{S}_k$-compatibility follows from
Proposition~\ref{prop:constant-scalar}.

It remains to verify regularity.  The transition matrix
$Q$ and directed-path propagator are the same as in
\S\ref{sec:RBJ}, so the propagator entries decay
geometrically.  Together with
\eqref{eq:LG-psi-tail} and
\eqref{eq:LG-phi-tail}, these bounds prove absolute
convergence of the sums defining $\Psi$, $\Phi$,
and $K$ in \eqref{eq:Psi-def}--\eqref{eq:K-def}.
\end{proof}

\paragraph{\textbf{Corrected edge weights.}}
Since $\mathcal{S}_k$-invariance holds,
Proposition~\ref{prop:constant-scalar}(ii) gives
$P(\mathcal{S}_k u) = P(u)$ and
$\Lambda_k(u) = \lambda_k I_m = D_k(u)$.
The diamond data on $\mathcal{G}^m$ is therefore
$(C_k, \Lambda_k) = (c_k I_m, \lambda_k I_m)$.
Theorem~\ref{thm:product-graph} requires additionally
that the resolvent $(I - zK(u))^{-1}$ exist.  At $z = 1$,
the Fredholm determinant
$F_{t, \mathbf{a}, \mathbf{n}} \defeq
\det(I - K_{t, \mathbf{a}, \mathbf{n}})$ equals the
multipoint gap probability
$\PP_y\bigl(\bigcap_i
\{Y_{n_i}(t) > a_i\}\bigr)$.  Since $K$ is trace class,
$F \ne 0$ if and only if the resolvent exists.  Particle
positions are non-increasing in time, so
$Y_{n_i}(t) \le y_{n_i}$.  If every jump size for
particles $1, \dotsc, n_m$ is zero (an event of
probability $p^{n_m t}$), no pushing occurs and
$Y_k(t) = y_k$ for $k \le n_m$.  The resolvent at
$z = 1$ therefore exists precisely on
$\mathcal{R} \defeq \{(t, \mathbf{a}, \mathbf{n}) \in
\mathcal{V}^m : a_i < y_{n_i}
\text{ for all } i\}$.
The mixed diamond at a vertex~$u$ requires the resolvent
at the four vertices $u$, $\mathcal{S}_1 u$,
$\mathcal{S}_2 u$, $\mathcal{S}_1\mathcal{S}_2 u$
entering~\eqref{eq:dressed-C-explicit}, together with
the $\mathcal{T}$-shifts $\mathcal{T}u$,
$\mathcal{T}\mathcal{S}_1 u$,
$\mathcal{T}\mathcal{S}_2 u$ from the proof
of~\eqref{eq:M-inverse}.\footnote{The explicit
formula~\eqref{eq:dressed-C-explicit} involves the
resolvent only at the four $\mathcal{S}$-shifted
vertices; a direct verification of the mixed diamond
from this formula would yield the larger domain
$a_i < y_{n_i+1}$.}
Since $\mathcal{S}_2$ preserves~$\mathcal{R}$ and
$\mathcal{T}$ tightens each threshold by~$1$, the
binding constraint is
$\mathcal{T}\mathcal{S}_1 u \in \mathcal{R}$, i.e.\
$a_i < y_{n_i+1} - 1$ for all~$i$.

\subsection{Multipoint equation}
\label{sec:LG-multipoint}

The dressed observable at $z = 1$,
\[
\mathcal{M}(u)
= I + \Phi(u)(I - K(u))^{-1}\Psi(u)
\in \End(\R^m),
\]
has dressed edge weights
\begin{equation}\label{eq:LG-dressed-Ck}
\mathcal{M}_k(u)
= c_k\mathcal{M}(\mathcal{T}u)^{-1}
  \mathcal{M}(\mathcal{S}_k u).
\end{equation}
For readability, we write
$\mathcal{M}_{t, \mathbf{a}, \mathbf{n}} \defeq
\mathcal{M}(t, \mathbf{a}, \mathbf{n})$.

\begin{theorem}\label{thm:LG-H8}
The dressed observable satisfies the mixed diamond
equation
\begin{equation}\label{eq:H8-multipoint}
\bigl(\mathcal{M}_{t, \mathbf{a}+\mathbf{1},
\mathbf{n}}^{-1}
- q\mathcal{M}_{t+1, \mathbf{a}+\mathbf{1},
\mathbf{n}+\mathbf{1}}^{-1}\bigr)
\mathcal{M}_{t, \mathbf{a},
\mathbf{n}+\mathbf{1}}
- \mathcal{M}_{t+1, \mathbf{a}+\mathbf{1},
\mathbf{n}}^{-1}
\bigl(\mathcal{M}_{t+1, \mathbf{a},
\mathbf{n}+\mathbf{1}}
- q\mathcal{M}_{t, \mathbf{a},
\mathbf{n}}\bigr) = 0,
\end{equation}
at every $(t, \mathbf{a}, \mathbf{n})$ with $n_m < N$ and
$a_i < y_{n_i+1} - 1$ for all~$i$.
\end{theorem}

\begin{proof}
By Theorem~\ref{thm:product-graph}, the pair
$(\mathcal{M}_k, \Lambda_k)$ satisfies the diamond
equations on~$\mathcal{R}$.
Since $\Lambda_k = \lambda_k I_m$, the scalar factors
commute with $\mathcal{M}_k$, and the mixed diamond
equation reduces to
\begin{equation}\label{eq:LG-mixed-scalar-Lambda}
\lambda_2\bigl[
\mathcal{M}_1(u) - \mathcal{M}_1(\mathcal{S}_2 u)
\bigr]
+ \lambda_1\bigl[
\mathcal{M}_2(\mathcal{S}_1 u)
- \mathcal{M}_2(u)
\bigr] = 0.
\end{equation}
Substituting \eqref{eq:LG-dressed-Ck} and the explicit
shifts, every term carries a common factor of
$-(p(1{-}\theta))^{-1}$.  Dividing by this factor
and evaluating at the vertex
$(t{+}1, \mathbf{a}, \mathbf{n})$ gives
\begin{equation}\label{eq:LG-mixed-reduced}
\begin{aligned}
\bigl[
&\mathcal{M}_{t+1, \mathbf{a}+\mathbf{1},
\mathbf{n}}^{-1}
\mathcal{M}_{t+1, \mathbf{a},
\mathbf{n}+\mathbf{1}}
- \mathcal{M}_{t,
\mathbf{a}+\mathbf{1},
\mathbf{n}}^{-1}
\mathcal{M}_{t, \mathbf{a},
\mathbf{n}+\mathbf{1}}
\bigr] \\
+ q\bigl[
&\mathcal{M}_{t+1, \mathbf{a}+\mathbf{1},
\mathbf{n}+\mathbf{1}}^{-1}
\mathcal{M}_{t, \mathbf{a},
\mathbf{n}+\mathbf{1}}
- \mathcal{M}_{t+1, \mathbf{a}+\mathbf{1},
\mathbf{n}}^{-1}
\mathcal{M}_{t, \mathbf{a},
\mathbf{n}}
\bigr] = 0.
\end{aligned}
\end{equation}
The middle two terms share the right factor
$\mathcal{M}_{t, \mathbf{a},
\mathbf{n}+\mathbf{1}}$, and the
first and last share the left factor
$\mathcal{M}_{t+1, \mathbf{a}+\mathbf{1},
\mathbf{n}}^{-1}$.  Factoring yields
\eqref{eq:H8-multipoint}.
\end{proof}

For $m = 1$, the dressed observable $\mathcal{M}_{t,a,n}$
is scalar.

\begin{corollary}\label{cor:LG-scalar}
The one-point Fredholm determinant
$F_{t, a, n} = \det(I - K_{t, a, n})$ satisfies
\begin{equation}\label{eq:LG-scalar-HM}
F_{t+1, a, n+1} F_{t, a+1, n}
- q F_{t, a, n} F_{t+1, a+1, n+1}
- p F_{t+1, a+1, n} F_{t, a, n+1}
= 0,
\end{equation}
at every $(t, a, n)$ with $a < y_{n+1} - 1$.
\end{corollary}

\begin{proof}
Proposition~\ref{prop:scalar-HM} applies at $z = 1$: the
boundary condition $F(T^\ell v) \to 1$ as
$\ell \to -\infty$ holds because the thresholds decrease
to $-\infty$ and the kernel vanishes.  The coefficients are
$\alpha_{12} = c_1 \lambda_2 =
-(p(1{-}\theta))^{-1}$
and
$\alpha_{21} = c_2 \lambda_1 =
-q(p(1{-}\theta))^{-1}$.
Substituting the lattice points into
\eqref{eq:HM-variable} and computing
$\alpha_{12} - \alpha_{21}
= -(1{-}q)(p(1{-}\theta))^{-1}
= -(1{-}\theta)^{-1}$
(using $p + q = 1$), then dividing by
$-(1{-}\theta)^{-1}$ and multiplying by $p$, gives
\eqref{eq:LG-scalar-HM}.
\end{proof}

\section{Geometric Last Passage Percolation with Boundary}
\label{sec:GLPP}

\paragraph{\textbf{System description.}}
We consider geometric last passage percolation with a
deterministic boundary condition.  The corner growth
function $G(m,n)$ for $m, n \ge 0$ satisfies the
recurrence
\[
G(m,n) = \max\{G(m{-}1,n), G(m,n{-}1)\} + \omega_{m,n}
\]
for $m, n \ge 1$, where
$\{\omega_{i,j}\}_{i,j \ge 1}$ are independent geometric
random variables with
$\PP(\omega_{i,j} = k) = (1{-}q)q^{k-1}$ for
$k = 1, 2, 3, \dots$ and $q \in (0,1)$.  The boundary
conditions are $G(0,n) = 0$ and $G(m,0) = x_m$, where
$0 = x_0 < x_1 < x_2 < \cdots$ is a strictly increasing
integer sequence.

\paragraph{\textbf{Fredholm determinant formula.}}
Fix $n \ge 1$, row indices
$\mathbf{m} = (m_1, \dots, m_k)$ with
$1 \le m_1 < \cdots < m_k$,
and integer thresholds
$\mathbf{a} = (a_1, \dots, a_k)$.  The multipoint
distribution is given by a Fredholm determinant on
$\ell^2(\{1, \dots, k\} \times \Z)$:
\[
\PP\Bigl(\bigcap_{i=1}^k
\{G(m_i,n) < a_i\}\Bigr)
= \det(I - \chi_{\mathbf{a}} K \chi_{\mathbf{a}}),
\]
where $\chi_{\mathbf{a}}(i, z) = \mathbf{1}_{z \ge a_i}$.
The correlation kernel $K$ has block entries
\[
K(i, z_1; j, z_2)
= -Q^{m_j - m_i}(z_1, z_2)\mathbf{1}_{m_i < m_j}
+ \bigl(\mathcal{S}_{n,-m_i}
\mathcal{S}_{n,m_j}^{\mathrm{hypo}(x)}\bigr)(z_1, z_2).
\]

Fix an auxiliary parameter $\theta \in (q, 1)$ and let
$\alpha = (1{-}\theta)\theta^{-1}$.  Define the rational factor
$\varphi(w) = (1{-}q)/(w{-}q)$.  The operators
defining the kernel are as follows.

\medskip
\noindent\textit{The transition matrix $Q$.}
The matrix $Q$ is the transition matrix of a random walk
on $\Z$ taking $\mathrm{Geom}[1{-}\theta]$ steps strictly
to the right:
\[
Q(z_1, z_2)
= (1{-}\theta)\theta^{z_2 - z_1 - 1}
\mathbf{1}_{z_2 > z_1},
\]
with $m$-th power admitting the integral representation
\[
Q^m(z_1, z_2)
= \frac{1}{2\pi\mathrm{i}}
\oint_{|w|=r} \diff w
\frac{\theta^{z_2 - z_1}}{w^{z_2 - z_1 - m + 1}}
\Bigl(\frac{\alpha}{1 - w}\Bigr)^m,
\qquad r < 1.
\]

\medskip
\noindent\textit{The scattering operators
$\mathcal{S}_{n,-m}$ and $\overline{\mathcal{S}}_{n,m}$.}\footnote{The kernel operators $\mathcal{S}_{n,-m}$ and $\overline{\mathcal{S}}_{n,m}$ are unrelated to the lattice shifts $\mathcal{S}_k$ of the diamond framework (\S\ref{sec:linear-problem}); the notation follows Rahman~\cite{Rahman2025}.}
These operators are built from the rational factor
$\varphi(w) = (1{-}q)/(w{-}q)$.  They are defined by
the contour
integrals
\[
\mathcal{S}_{n,-m}(z_1, z_2)
= \frac{1}{2\pi\mathrm{i}}
\oint_{|w|=r} \diff w
\frac{\theta^{z_2 - z_1}}
{w^{z_2 - z_1 + m + 1}}
\Bigl(\frac{1{-}w}{\alpha}\Bigr)^{\!m}
\varphi(w)^n,
\qquad r \in (q, 1),
\]
\[
\overline{\mathcal{S}}_{n,m}(z_1, z_2)
= \frac{1}{2\pi\mathrm{i}}
\oint_{|w|=\delta} \diff w
\frac{\theta^{z_2 - z_1}
(1{-}w)^{z_1 - z_2 + m - 1}}
{(w/\alpha)^m}
\varphi(1{-}w)^{-n},
\qquad \delta < 1.
\]

\medskip
\noindent\textit{The hypograph operator
$\mathcal{S}_{n,m}^{\mathrm{hypo}(x)}$.}
This is the hitting probability operator defined in terms
of the boundary data:
\[
\mathcal{S}_{n,m}^{\mathrm{hypo}(x)}(z_1, z_2)
= \E_{B_0 = z_1}\bigl[
\overline{\mathcal{S}}_{n,m-\tau}(B_\tau, z_2)
\mathbf{1}_{\tau < m}\bigr],
\]
where $(B_\ell)_{\ell \ge 0}$ is the random walk with
transition matrix $Q$, and
$\tau = \min\{\ell \ge 0 : B_\ell < x_{\ell+1}\}$
is the hitting time of the boundary.

\paragraph{\textbf{References.}}
Rahman~\cite{Rahman2025} proves the fixed-time Fredholm determinant for geometric LPP with deterministic boundary and gives the associated kernel.  The proof passes through Rahman's duality framework and applies Matetski--Remenik's Theorem~1.2~\cite{MatetskiRemenik2023}.

\subsection{Kernel reformulation}
\label{sec:GLPP-kernel-reform}

We rewrite the extended-kernel Fredholm determinant in
the form of the product graph construction of \S\ref{sec:product-graph}.

\paragraph{\textbf{Seed data.}}
Define
\[
\psi_{n,m,a}(x) \defeq \mathcal{S}_{n,-m}(a, x),
\qquad
\phi_{n,m,a}(x) \defeq
\mathcal{S}_{n,m}^{\mathrm{hypo}(x)}(x, a),
\]
where we suppress the dependence of $\phi$ on the
boundary condition.

\paragraph{\textbf{Propagator.}}
The propagator is a directed-path sum with transition
kernels $Q_{\ell,\ell'} = Q^{m_{\ell'} - m_\ell}$ for
$\ell < \ell'$ and upper-tail internal constraints.
For $i > j$,
\begin{equation}\label{eq:GLPP-propagator}
[B_{\mathbf{a}, \mathbf{m}}]_{ij}(r, r')
\defeq \sum_{k=1}^{i-j}
  \sum_{j = \ell_0 < \ell_1 < \cdots < \ell_k = i}
  \sum_{\substack{\xi_1 \ge a_{\ell_1} \\ \vdots \\
  \xi_{k-1} \ge a_{\ell_{k-1}}}}
  \prod_{s=0}^{k-1}
  \Bigl(-Q^{m_{\ell_{s+1}} - m_{\ell_s}}
  (\xi_s, \xi_{s+1})\Bigr),
\end{equation}
with $\xi_0 = r'$, $\xi_k = r$, and $[B_{\mathbf{a}, \mathbf{m}}]_{ij} = 0$ for $i \le j$.
The internal vertices are constrained by the
upper-tail cutoffs; the endpoints $r, r'$ are
arbitrary.  The propagator depends on $\mathbf{a}$
and $\mathbf{m}$ but not on $n$.

On the cutoff region $r \ge a_i$, $r' \ge a_j$, the
entries of $I + B_{\mathbf{a}, \mathbf{m}}$ coincide with those of
$((I + \chi_{\mathbf{a}} L
\chi_{\mathbf{a}})^{-1})^{\!\top}$, where $L$
is the strictly upper-triangular block operator with
entries $L_{ij}(z, z') = Q^{m_j - m_i}(z, z')
\mathbf{1}_{i < j}$; the identification follows from
the same Neumann expansion as
Remark~\ref{rem:Neumann-series}, with the upper-tail
projection $\chi_{\mathbf{a}}$ in place of
$\bar{\chi}_{\mathbf{a}}$.

\begin{lemma}\label{lem:GLPP-kernel}
We have
\[
\PP\Bigl(\bigcap_{i=1}^k
\{G(m_i, n) < a_i\}\Bigr)
= \det(I - K_{n, \mathbf{m}, \mathbf{a}})
_{\ell^2(\Z)},
\]
where
\[
K_{n, \mathbf{m}, \mathbf{a}}(x, x')
= \sum_{1 \le j \le i \le k}
  \sum_{r \ge a_i} \sum_{r' \ge a_j}
   \psi_{n, m_i, r}(x)
  \bigl[\delta_{ij}\delta_{r,r'}
  + [B_{\mathbf{a}, \mathbf{m}}]_{ij}(r, r')\bigr]
  \phi_{n, m_j, r'}(x').
\]
This identifies $K_{n, \mathbf{m}, \mathbf{a}}$ with the
kernel~\eqref{eq:K-def} on the product graph
$\mathcal{G}^k$, with $\psi, \phi$ as the seed data and
$B_{\mathbf{a}, \mathbf{m}}$ as the propagator.
\end{lemma}

\begin{proof}
The argument is identical to that of
Lemma~\ref{lem:RBJ-kernel}: the block-triangular
decomposition, Sylvester's identity, and the transpose
identification produce the claimed single-space kernel.
The random walk $Q$ moves strictly to the right;
the estimate for $\psi$ requires modification because
$\varphi(w)^n$ has a pole at $w = q$ inside the contour.

The contour formula for $\psi$ gives
\[
\psi_{n,m,a}(x) = \frac{1}{2\pi\mathrm{i}}
  \oint_{|w|=\rho} \diff w
  \frac{\theta^{x-a}}{w^{x-a+m+1}}
  \Bigl(\frac{1{-}w}{\alpha}\Bigr)^{\!m}
  \varphi(w)^n,
\qquad q < \rho < \theta.
\]
On $|w| > q$, the factor $\varphi(w)^n$ contributes
only negative powers of $w$, starting from $w^{-n}$,
so $\psi_{n,m,a}(x) = 0$ for $x > a - n$.
Taking the Cauchy estimate on $|w| = \rho$ gives
\begin{equation}\label{eq:GLPP-psi-tail}
|\psi_{n,m,a}(x)| \le C(\rho/\theta)^{a-x}
\mathbf{1}_{x \le a - n}.
\end{equation}
Since $\rho < \theta$,
$\psi_{n,m,a} \in \ell^2(\Z)$ with norm bounded
uniformly in $a$.

The contour estimate for $\phi$ follows the same
argument as in Lemma~\ref{lem:RBJ-kernel}.  Choose
$\beta \in (\theta, 1)$ and $\gamma_\delta$ so that
$\sup_{w \in \gamma_\delta}|\theta/(1{-}w)| < \beta$.
The factor $\varphi(1{-}w)^{-n}$ in the
$\overline{\mathcal{S}}$ contour formula is bounded on
$\gamma_\delta$, so
$|\overline{\mathcal{S}}_{n,k}(z, r)|
\le C\beta^{r-z}$ for $r \ge z$.  The exponential moment
$\E[\beta^{-S_\ell}]$ is finite for $\beta > \theta$,
where $S_\ell$ is the cumulative displacement of the
$Q$-walk after $\ell$ steps.  The hitting
representation gives
\begin{equation}\label{eq:GLPP-phi-tail}
|\phi_{n,m,a}(x)| \le C\beta^{a-x}
\mathbf{1}_{x < x_m}.
\end{equation}

Choose conjugation weights
$\theta < \omega_1 < \cdots < \omega_k < 1$.
With these bounds, the same conjugation and Sylvester
reduction as in Lemma~\ref{lem:RBJ-kernel} apply.
\end{proof}

\subsection{Base graph and seed data}
\label{sec:GLPP-base-graph}

We now identify the lattice structure that governs how
the seed data evolve under shifts of the parameters
$(n, m, a)$, and verify that $\psi, \phi$ satisfy the
linear problems of the diamond framework.

Let $\mathcal{V} \subseteq \Z^3$ be the lattice
generated by the shifts
\[
T = e^{-\partial_a}, \qquad
S_1 = e^{\partial_m}, \qquad
S_2 = e^{-\partial_n - \partial_a},
\]
so that for $u = (n, m, a) \in \mathcal{V}$,
\[
Tu = (n, m, a{-}1), \qquad
S_1 u = (n, m{+}1, a), \qquad
S_2 u = (n{-}1, m, a{-}1).
\]
Take $E = \mathbb{F} = \R$ and define the
constant scalar edge weights
\begin{equation}\label{eq:GLPP-edge-weights}
(c_1, c_2)
= \bigl((1{-}\theta)^{-1},
-q(1{-}q)^{-1}\bigr),
\qquad
(\lambda_1, \lambda_2)
= \bigl((1{-}\theta)^{-1}\theta,
-(1{-}q)^{-1}\theta\bigr).
\end{equation}
Since $(c_j, \lambda_j)$ are constant scalars, the
diamond equations
\eqref{eq:diamond-C-ij}--\eqref{eq:diamond-mixed-ij}
are satisfied trivially.

\begin{lemma}\label{lem:GLPP-seed-linear}
Let $H = \ell^2(\Z)$.  The seed functions
$\psi_{n,m,a}$ and $\phi_{n,m,a}$ satisfy the
linear problem \eqref{eq:Linear-Psi-3} and the
adjoint linear problem \eqref{eq:Linear-Phi-Prob}
on $\mathcal{V}$ with edge weights
$(c_k, \lambda_k)$.  Explicitly:
\begin{align}
\psi_{n,m+1,a}
&= (1{-}\theta)^{-1}\psi_{n,m,a-1}
   - (1{-}\theta)^{-1}\theta\psi_{n,m,a},
\label{eq:GLPP-psi-S1} \\
\psi_{n-1,m,a-1}
&= -q(1{-}q)^{-1}\psi_{n,m,a-1}
   + (1{-}q)^{-1}\theta\psi_{n,m,a},
\label{eq:GLPP-psi-S2}
\end{align}
and
\begin{align}
\phi_{n,m,a-1}
&= (1{-}\theta)^{-1}\phi_{n,m+1,a}
   - (1{-}\theta)^{-1}\theta\phi_{n,m+1,a-1},
\label{eq:GLPP-phi-S1} \\
\phi_{n,m,a-1}
&= -q(1{-}q)^{-1}\phi_{n-1,m,a-1}
   + (1{-}q)^{-1}\theta\phi_{n-1,m,a-2}.
\label{eq:GLPP-phi-S2}
\end{align}
\end{lemma}

\begin{proof}
We verify each identity in turn.

\medskip
\noindent\textit{Verification of
\eqref{eq:GLPP-psi-S1}.}
The $S_1$-shift involves only $m$ and the
$Q$-power structure, which is the same as in
\S\ref{sec:RBJ}; the factor $\varphi(w)^n$ passes through
unchanged.  The computation is identical to the
verification of \eqref{eq:RBJ-psi-S1}: the factor
$(1{-}\theta)^{-1}(\theta/w - \theta)
= (1{-}w)/(\alpha w)$ shifts the integrand from
$m$ to $m{+}1$.

\medskip
\noindent\textit{Verification of
\eqref{eq:GLPP-psi-S2}.}
From the contour representation,
\[
\psi_{n,m,a}(x) = \frac{1}{2\pi\mathrm{i}}
  \oint \diff w
  \frac{\theta^{x-a}}{w^{x-a+m+1}}
  \Bigl(\frac{1{-}w}{\alpha}\Bigr)^{\!m}
  \varphi(w)^n.
\]
The right-hand side of \eqref{eq:GLPP-psi-S2} is
\begin{align*}
&-q(1{-}q)^{-1}\psi_{n,m,a-1}(x)
+ (1{-}q)^{-1}\theta\psi_{n,m,a}(x) \\
&\qquad = \frac{1}{2\pi\mathrm{i}}
  \oint \diff w
  \frac{\theta^{x-a}}{w^{x-a+m+1}}
  \Bigl(\frac{1{-}w}{\alpha}\Bigr)^{\!m}
  \varphi(w)^n
  \frac{\theta}{1{-}q}
  \Bigl(1 - \frac{q}{w}\Bigr).
\end{align*}
The parenthetical factor simplifies as
\[
\frac{\theta}{1{-}q}\Bigl(1 - \frac{q}{w}\Bigr)
= \frac{\theta(w{-}q)}{(1{-}q)w}.
\]
Since $\varphi(w)^{-1} = (w{-}q)/(1{-}q)$, the
$\varphi$-power shifts from $n$ to $n{-}1$.  The factor
$\theta/w$ shifts $a$ to $a{-}1$:
\[
= \frac{1}{2\pi\mathrm{i}}
  \oint \diff w
  \frac{\theta^{x-(a-1)}}{w^{x-(a-1)+m+1}}
  \Bigl(\frac{1{-}w}{\alpha}\Bigr)^{\!m}
  \varphi(w)^{n-1}
= \psi_{n-1,m,a-1}(x).
\]

\medskip
\noindent\textit{Verification of
\eqref{eq:GLPP-phi-S1}.}
The argument is completely analogous to the proof of
\eqref{eq:RBJ-phi-S1} in \S\ref{sec:RBJ}.  The
$\overline{\mathcal{S}}$ recurrence in $m$ is
\begin{equation}\label{eq:GLPP-Sbar-recurrence-m}
(1{-}\theta)^{-1}\bigl[
\overline{\mathcal{S}}_{n,m+1}(z, r{+}1)
- \theta\overline{\mathcal{S}}_{n,m+1}(z, r)
\bigr]
= \overline{\mathcal{S}}_{n,m}(z, r),
\end{equation}
whose verification depends only on the
$(1{-}w)$-power and $w$-power structure and is
identical to that of \eqref{eq:Sbar-recurrence-n}.

We show that
$\overline{\mathcal{S}}_{n,0} \equiv 0$.  When $m = 0$,
the measure $dw/w^m$ reduces to $dw$, removing the
pole at the origin.  The remaining integrand is
analytic inside $|w| = \delta$, so the integral
vanishes by Cauchy's theorem.

The lifting to the hypograph version proceeds as in
\S\ref{sec:RBJ}: the expectation splits at
$\tau \le m{-}1$ (where
\eqref{eq:GLPP-Sbar-recurrence-m} applies) and
$\tau = m$ (where $\overline{\mathcal{S}}_{n,0} = 0$
kills the boundary term).

\medskip
\noindent\textit{Verification of
\eqref{eq:GLPP-phi-S2}.}
Since neither side involves a shift in $m$, the
indicator $\mathbf{1}_{\tau < m}$ passes through
unchanged, and the verification reduces to a
recurrence for the free operator:
\begin{equation}\label{eq:GLPP-Sbar-recurrence-n}
-q(1{-}q)^{-1}
\overline{\mathcal{S}}_{n-1,m}(z, r)
+ (1{-}q)^{-1}\theta
\overline{\mathcal{S}}_{n-1,m}(z, r{-}1)
= \overline{\mathcal{S}}_{n,m}(z, r).
\end{equation}
In the contour formula for
$\overline{\mathcal{S}}_{n-1,m}$, the shift from $r$
to $r{-}1$ multiplies the integrand by
$\theta^{-1}(1{-}w)$.  On the left-hand side, the
bracket produces the factor
$(-q + (1{-}w))/(1{-}q)
= (1{-}q{-}w)/(1{-}q) = \varphi(1{-}w)^{-1}$,
which promotes $\varphi(1{-}w)^{-(n-1)}$ to
$\varphi(1{-}w)^{-n}$.  Taking the expectation
against the hitting-time representation completes
the proof.
\end{proof}

\subsection{Admissible propagators and corrected
edge weights}
\label{sec:GLPP-propagators}

The product graph $\mathcal{G}^k$ has vertex set
$\mathcal{V}^k$ and is generated by the diagonal shifts
\[
\mathcal{T}u
= (n, \mathbf{m}, \mathbf{a} - \mathbf{1}),
\qquad
\mathcal{S}_1 u
= (n, \mathbf{m} + \mathbf{1}, \mathbf{a}),
\qquad
\mathcal{S}_2 u
= (n{-}1, \mathbf{m}, \mathbf{a} - \mathbf{1}),
\]
where $u = (n, \mathbf{m}, \mathbf{a}) \in
\mathcal{V}^k$ and
$\mathbf{1} = (1, \ldots, 1) \in \Z^k$.
The block-diagonal edge weights are
$C_j(u) = c_j I_k$ and $D_j(u) = \lambda_j I_k$.

Set $\boldsymbol{\beta} = -\mathbf{a}$ and
$\widetilde{Q}_{\ell,\ell'}(r, r')
= Q^{m_{\ell'} - m_\ell}(-r, -r')$.
The propagator~\eqref{eq:GLPP-propagator} is the
pullback
$[B_{\mathbf{a}, \mathbf{m}}]_{ij}(r, r')
= [\widetilde{B}_{-\mathbf{a}, \mathbf{m}}]_{ij}(-r, -r')$
of the lower-tail directed-path propagator
(Definition~\ref{def:directed-path-propagator}) built
from $(\boldsymbol{\beta}, \widetilde{Q})$.  The
standing hypotheses of \S\ref{sec:directed-path} are
satisfied: $\boldsymbol{\beta}(\mathcal{T}u)
= \boldsymbol{\beta}(u) + \mathbf{1}$, and
$\widetilde{Q}$ is translation invariant and
independent of $\boldsymbol{\beta}$.

\begin{lemma}[$\mathcal{S}_j$-invariance]
\label{lem:GLPP-Sk-invariance}
For $j = 1, 2$ and all $x, y \in \mathcal{V}^k$,
$B_{\mathcal{S}_j u}
(\mathcal{S}_j x, \mathcal{S}_j y)
= B_u(x, y)$.
\end{lemma}

\begin{proof}
\noindent\textit{Case $j = 1$.}
$\mathcal{S}_1 u =
(n, \mathbf{m}{+}\mathbf{1}, \mathbf{a})$: the
cutoffs $\mathbf{a}$ and the differences
$m_\ell - m_i$ are unchanged, and $\mathcal{S}_1$
does not shift $\mathbf{a}$, so it acts trivially
on spatial coordinates.

\medskip
\noindent\textit{Case $j = 2$.}
$\mathcal{S}_2 u =
(n{-}1, \mathbf{m}, \mathbf{a}{-}\mathbf{1})$:
the cutoffs shift to $\mathbf{a}{-}\mathbf{1}$
and the label differences are unchanged.  Since
$\mathcal{S}_2$ acts on $\mathbf{a}$ identically
to $\mathcal{T}$, the invariance follows from the
same translation-invariance argument as
$\mathcal{T}$-covariance (the propagator is
$n$-independent).
\end{proof}

\begin{lemma}\label{lem:GLPP-admissible}
The propagator $B_{\mathbf{a}, \mathbf{m}}$ is admissible for
$\mathcal{G}^k$.
\end{lemma}

\begin{proof}
By Proposition~\ref{prop:directed-path-T} applied to
the reflected propagator
$\widetilde{B}_{\boldsymbol{\beta}, \mathbf{m}}$,
$\mathcal{T}$-covariance and $\mathcal{T}$-splitting
hold for $B_{\mathbf{a}, \mathbf{m}}$.
The edge weights $(c_j, \lambda_j)$ are constant scalars
and the propagator satisfies $\mathcal{S}_j$-invariance
(Lemma~\ref{lem:GLPP-Sk-invariance}), so
$\mathcal{S}_j$-compatibility follows from
Proposition~\ref{prop:constant-scalar}.

It remains to verify regularity.  The one-sided support
of $Q$ ensures that the propagator entries decay
geometrically.  Together with
\eqref{eq:GLPP-psi-tail} and \eqref{eq:GLPP-phi-tail},
these bounds prove absolute convergence of the sums
defining $\Psi$, $\Phi$, and $K$ in
\eqref{eq:Psi-def}--\eqref{eq:K-def}.
\end{proof}

\paragraph{\textbf{Corrected edge weights.}}
Since $\mathcal{S}_j$-invariance holds,
Proposition~\ref{prop:constant-scalar}(ii) gives
$P(\mathcal{S}_j u) = P(u)$ and
$\Lambda_j(u) = \lambda_j I_k = D_j(u)$.
The diamond data on $\mathcal{G}^k$ is therefore
$(C_j, \Lambda_j) = (c_j I_k, \lambda_j I_k)$.

Theorem~\ref{thm:product-graph} requires additionally
that the resolvent $(I - zK(u))^{-1}$ exist.  At $z = 1$,
the Fredholm determinant
$F_{n, \mathbf{m}, \mathbf{a}} \defeq
\det(I - K_{n, \mathbf{m}, \mathbf{a}})$ equals the
multipoint gap probability
$\PP\bigl(\bigcap_i
\{G(m_i, n) < a_i\}\bigr)$.  Since $K$ is trace class,
$F \ne 0$ if and only if the resolvent exists.  The
geometric weights satisfy $\omega_{i,j} \ge 1$, so
$G(m, n) \ge x_m + n$ for $m, n \ge 1$.  If
$a_i > x_{m_i} + n$ for every~$i$, all weights in
$\{1, \dots, m_k\} \times \{1, \dots, n\}$ can equal
one, giving $G(m_i, n) = x_{m_i} + n < a_i$ for
all~$i$.  The resolvent at $z = 1$ therefore exists
precisely on
$\mathcal{R} \defeq \{(n, \mathbf{m}, \mathbf{a}) \in
\mathcal{V}^k : n \ge 1,\
a_i > x_{m_i} + n \text{ for all } i\}$.
The mixed diamond at a vertex~$u$ requires the resolvent
at the four vertices $u$, $\mathcal{S}_1 u$,
$\mathcal{S}_2 u$, $\mathcal{S}_1\mathcal{S}_2 u$
entering~\eqref{eq:dressed-C-explicit}, together with
the $\mathcal{T}$-shifts $\mathcal{T}u$,
$\mathcal{T}\mathcal{S}_1 u$,
$\mathcal{T}\mathcal{S}_2 u$ from the proof
of~\eqref{eq:M-inverse}.\footnote{The explicit
formula~\eqref{eq:dressed-C-explicit} involves the
resolvent only at the four $\mathcal{S}$-shifted
vertices; a direct verification of the mixed diamond
from this formula would yield the larger domain
$a_i > x_{m_i+1} + n$.}
The $\mathcal{S}_2$-shift preserves the margin
$a_i - x_{m_i} - n$ but lowers $n$ by~$1$,
requiring $n \ge 2$.  Since $\mathcal{T}$ tightens
each threshold by~$1$, the binding constraint on $a_i$
is $\mathcal{T}\mathcal{S}_1 u \in \mathcal{R}$, i.e.\
$a_i > x_{m_i+1} + n + 1$ for all~$i$.

\begin{remark}\label{rem:GLPP-resolvent-shifts}
The condition $n \ge 2$ ensures that the
$\mathcal{S}_2$-shifted vertices have $n - 1 \ge 1$,
as required by the Fredholm determinant formula.  If the
determinant formula is extended to $n = 0$ using the
deterministic boundary $G(m, 0) = x_m$, the condition
relaxes to $n \ge 1$.
\end{remark}

\subsection{Multipoint equation}
\label{sec:GLPP-multipoint}

The dressed observable at $z = 1$,
\[
\mathcal{M}(u)
= I + \Phi(u)(I - K(u))^{-1}\Psi(u)
\in \End(\R^k),
\]
has dressed edge weights
\begin{equation}\label{eq:GLPP-dressed-Ck}
\mathcal{M}_j(u)
= c_j\mathcal{M}(\mathcal{T}u)^{-1}
  \mathcal{M}(\mathcal{S}_j u).
\end{equation}
For readability, we write
$\mathcal{M}_{n, \mathbf{m}, \mathbf{a}} \defeq
\mathcal{M}(n, \mathbf{m}, \mathbf{a})$.

\begin{theorem}\label{thm:GLPP-H4p}
The dressed observable satisfies the mixed diamond
equation
\begin{equation}\label{eq:H4p-multipoint}
\begin{aligned}
&\bigl(\mathcal{M}_{n-1, \mathbf{m},
\mathbf{a}-2\cdot\mathbf{1}}^{-1}
- q\mathcal{M}_{n, \mathbf{m}+\mathbf{1},
\mathbf{a}-\mathbf{1}}^{-1}\bigr)
\mathcal{M}_{n-1, \mathbf{m}+\mathbf{1},
\mathbf{a}-\mathbf{1}} \\
&\qquad - \mathcal{M}_{n, \mathbf{m},
\mathbf{a}-\mathbf{1}}^{-1}
\bigl(\mathcal{M}_{n, \mathbf{m}+\mathbf{1},
\mathbf{a}}
- q\mathcal{M}_{n-1, \mathbf{m},
\mathbf{a}-\mathbf{1}}\bigr) = 0,
\end{aligned}
\end{equation}
at every $(n, \mathbf{m}, \mathbf{a})$ with $n \ge 2$ and
$a_i > x_{m_i+1} + n + 1$ for all~$i$.
\end{theorem}

\begin{proof}
By Theorem~\ref{thm:product-graph}, the pair
$(\mathcal{M}_j, \Lambda_j)$ satisfies the diamond
equations on~$\mathcal{R}$.
Since $\Lambda_j = \lambda_j I_k$, the scalar factors
commute with $\mathcal{M}_j$, and the mixed diamond
equation reduces to
\begin{equation}\label{eq:GLPP-mixed-scalar-Lambda}
\lambda_2\bigl[
\mathcal{M}_1(u) - \mathcal{M}_1(\mathcal{S}_2 u)
\bigr]
+ \lambda_1\bigl[
\mathcal{M}_2(\mathcal{S}_1 u)
- \mathcal{M}_2(u)
\bigr] = 0.
\end{equation}

Substituting \eqref{eq:GLPP-dressed-Ck} and the explicit
shifts, the coefficients are
$c_1 \lambda_2
= -(1{-}\theta)^{-1}(1{-}q)^{-1}\theta$
and
$c_2 \lambda_1
= -q(1{-}q)^{-1}(1{-}\theta)^{-1}\theta$.
Dividing
\eqref{eq:GLPP-mixed-scalar-Lambda} by
$-(1{-}\theta)^{-1}(1{-}q)^{-1}\theta$ gives
\begin{multline}\label{eq:GLPP-mixed-reduced}
\bigl[
\mathcal{M}_{n, \mathbf{m},
\mathbf{a}-\mathbf{1}}^{-1}
\mathcal{M}_{n, \mathbf{m}+\mathbf{1},
\mathbf{a}}
- \mathcal{M}_{n-1, \mathbf{m},
\mathbf{a}-2\cdot\mathbf{1}}^{-1}
\mathcal{M}_{n-1, \mathbf{m}+\mathbf{1},
\mathbf{a}-\mathbf{1}}
\bigr] \\
+ q\bigl[
\mathcal{M}_{n, \mathbf{m}+\mathbf{1},
\mathbf{a}-\mathbf{1}}^{-1}
\mathcal{M}_{n-1, \mathbf{m}+\mathbf{1},
\mathbf{a}-\mathbf{1}}
- \mathcal{M}_{n, \mathbf{m},
\mathbf{a}-\mathbf{1}}^{-1}
\mathcal{M}_{n-1, \mathbf{m},
\mathbf{a}-\mathbf{1}}
\bigr] = 0.
\end{multline}
The middle two terms share the right factor
$\mathcal{M}_{n-1, \mathbf{m}+\mathbf{1},
\mathbf{a}-\mathbf{1}}$, and the
first and last share the left factor
$\mathcal{M}_{n, \mathbf{m},
\mathbf{a}-\mathbf{1}}^{-1}$.  Factoring and
multiplying by~$-1$ yields
\eqref{eq:H4p-multipoint}.
\end{proof}

For $k = 1$, the dressed observable
$\mathcal{M}_{n,m,a}$ is scalar.

\begin{corollary}\label{cor:GLPP-scalar}
The one-point Fredholm determinant
$F_{n, m, a} = \det(I - K_{n, m, a})$ satisfies
\begin{equation}\label{eq:GLPP-scalar-HM}
F_{n, m+1, a}F_{n-1, m, a-2}
- qF_{n-1, m, a-1}F_{n, m+1, a-1}
- (1{-}q)F_{n, m, a-1}F_{n-1, m+1, a-1}
= 0,
\end{equation}
at every $(n, m, a)$ with $n \ge 2$ and
$a > x_{m+1} + n + 1$.
\end{corollary}

\begin{proof}
Proposition~\ref{prop:scalar-HM} applies at $z = 1$: the
boundary condition $F(T^\ell v) \to 1$ as
$\ell \to -\infty$ holds because the thresholds increase
to $\infty$ and the kernel vanishes.  The coefficients are
$\alpha_{12} = c_1 \lambda_2 =
-(1{-}\theta)^{-1}(1{-}q)^{-1}\theta$
and
$\alpha_{21} = c_2 \lambda_1 =
-q(1{-}q)^{-1}(1{-}\theta)^{-1}\theta$.
Substituting the lattice points into
\eqref{eq:HM-variable} and computing
$\alpha_{12} - \alpha_{21}
= -(1{-}\theta)^{-1}\theta$,
then dividing by
$-(1{-}\theta)^{-1}(1{-}q)^{-1}\theta$, gives
\eqref{eq:GLPP-scalar-HM}.
\end{proof}

\section{Stochastic Six-Vertex Model}\label{sec:S6V}

\paragraph{\textbf{System description.}}
We use the interacting particle system formulation of
the stochastic six-vertex model, following
Borodin--Corwin--Gorin \cite[Sec.~2.2]{BorodinCorwinGorin_S6V}.
Fix stochastic six-vertex parameters $0 < b_2 < b_1 < 1$
and set
\[
\tau \defeq \frac{b_2}{b_1} \in (0,1),
\qquad
\alpha \defeq
\tau^{-1}\frac{1 - b_1}{1 - b_2}.
\]
The system consists of infinitely many particles on
$\Z$.  The configuration at discrete time
$t \in \Z_{\ge 0}$ is
$X(t) = (x_1(t) < x_2(t) < \cdots) \in \Z^\infty$,
with step initial data $X(0) = (1, 2, 3, \dotsc)$.
At each time step, the positions $x_i(t{+}1)$ are
determined sequentially in the order $i = 1, 2, \dotsc$,
with the convention $x_0(t{+}1) = -\infty$.
The transition probabilities depend on whether the
$(i{-}1)$-st particle has been pushed into the position
of particle~$i$.  If $x_{i-1}(t{+}1) < x_i(t)$
(particle~$i$ has not been pushed), then
\begin{align*}
&\PP\bigl(x_i(t{+}1) = x_i(t) + k
  \mid X(t), x_{i-1}(t{+}1)\bigr) \\
&\qquad=
\begin{cases}
b_1,
  & k = 0, \\
(1 - b_1)(1 - b_2) b_2^{k-1},
  & 0 < k < x_{i+1}(t) - x_i(t), \\
(1 - b_1) b_2^{x_{i+1}(t) - x_i(t) - 1},
  & k = x_{i+1}(t) - x_i(t).
\end{cases}
\end{align*}
If $x_{i-1}(t{+}1) = x_i(t)$ (particle~$i$ has been
pushed from its position), then particle~$i$ must jump
at least one step to the right:
\begin{align*}
&\PP\bigl(x_i(t{+}1) = x_i(t) + k
  \mid X(t), x_{i-1}(t{+}1)\bigr) \\
&\qquad=
\begin{cases}
(1 - b_2) b_2^{k-1},
  & 0 < k < x_{i+1}(t) - x_i(t), \\
b_2^{x_{i+1}(t) - x_i(t) - 1},
  & k = x_{i+1}(t) - x_i(t).
\end{cases}
\end{align*}
In both cases, the value
$k = x_{i+1}(t) - x_i(t)$ corresponds to particle~$i$
reaching the current position of particle~$i{+}1$.
When this occurs, the subsequent update of
particle~$i{+}1$ falls into the second case above,
so particle~$i{+}1$ is forced to move at least one
site during its own update.
Informally, each particle flips a $b_1$-biased coin to
decide whether to stay; if it moves, it jumps to the
right by a $\mathrm{Geom}(b_2)$-distributed number of
steps, truncated by the next particle.
For $x \in \Z$, the height function is
\[
N_x(t) \defeq
\#\{i \ge 1 : x_i(t) \le x\}.
\]
For $\zeta \in \C \setminus \R_{\ge 0}$,
define the $\tau$-Laplace transform
\begin{equation}\label{eq:S6V-L-def}
\mathcal{L}_{t,x}(\zeta)
\defeq
\E_{\mathrm{step}}
\left[
\frac{1}{(\zeta\tau^{N_x(t)};\tau)_\infty}
\right],
\qquad
(a;\tau)_\infty
\defeq \prod_{k=0}^\infty (1 - a\tau^k).
\end{equation}

\paragraph{\textbf{Fredholm determinant formula.}}
The following is due to
Borodin--Corwin--Gorin~\cite[Thm.~4.16]{BorodinCorwinGorin_S6V}.
For $t \in \Z_{>0}$, $x \in \Z$, and
$\zeta \in \C \setminus \R_{\ge 0}$,
\begin{equation}\label{eq:S6V-fredholm}
\mathcal{L}_{t,x}(\zeta)
=
\det\bigl(I + K_\zeta^{b_1,b_2}[t,x]\bigr)_{L^2(C_r)},
\end{equation}
where $C_r$ is the positively oriented circle of radius $r$
satisfying $\tau < r < 1/\alpha$.  Since
$1/\alpha = \tau(1 - b_2)/(1 - b_1)$, this is the
contour condition in
Borodin--Corwin--Gorin~\cite{BorodinCorwinGorin_S6V}.
The kernel is given by
\begin{equation}\label{eq:s6v_kernel_main}
K_\zeta^{b_1,b_2}[t,x](w,w')
=
\frac{1}{2\mathrm{i}}
\int_{\frac{1}{2} + \mathrm{i}\R}
\frac{(-\zeta)^s}{\sin(\pi s)}
\frac{g_{t,x}(w)}{g_{t,x}(\tau^s w)}
\frac{\diff s}{\tau^s w - w'},
\end{equation}
with
\begin{equation}\label{eq:S6V-g-def}
g_{t,x}(z)
=
(1 + \alpha z)^t(1 + \tau^{-1}z)^{-x}.
\end{equation}
The $s$-contour is the vertical line
$\operatorname{Re}(s) = 1/2$ oriented bottom-to-top,
and $(-\zeta)^s$ is defined using the principal branch
of $\log(-\zeta)$.

\paragraph{\textbf{References.}}
The one-point Fredholm determinant formula
\eqref{eq:S6V-fredholm} used in this section is
Borodin--Corwin--Gorin~\cite[Thm.~4.16]{BorodinCorwinGorin_S6V}.

\begin{remark}[Fredholm determinant convention]
\label{rem:S6V-fredholm-conv}
The Fredholm determinant
$\det(I + K_\zeta^{b_1,b_2}[t,x])$ on $L^2(C_r)$ is
defined by the contour series of
Borodin--Corwin--Gorin~\cite[Def.~4.13]{BorodinCorwinGorin_S6V},
with the contour measure
$\frac{\diff w}{2\pi\mathrm{i}}$ on $C_r$.  The
combined factor $(-\zeta)^s/\sin(\pi s)$ decays
exponentially on the vertical $s$-contour for each
$\zeta \in \C \setminus \R_{\ge 0}$, giving uniform
convergence of the Mellin--Barnes integral on
$C_r \times C_r$.  The kernel is therefore continuous
and bounded on the contour square, and the Fredholm
series converges absolutely by the Hadamard estimate
of~\cite[Cor.~4.15]{BorodinCorwinGorin_S6V}.  Since
$(-\zeta)^s$ is analytic in $\zeta$ on
$\C \setminus \R_{\ge 0}$ (principal branch), the same
estimates give analyticity of
$\mathcal{L}_{t,x}(\zeta)$ in $\zeta$ on this domain.
When the determinant is nonzero, the resolvent kernel
exists by classical Fredholm theory.  The framework
results invoked below (gauge invariance, the Woodbury
determinant ratio, the Darboux theorem, and the scalar
Hirota--Miwa reduction) hold at the kernel and series
level: their proofs require only pointwise kernel
identities, the resolvent integral equation, and the
rank-one perturbation formula for Fredholm series.
\end{remark}

\subsection{Kernel reformulation}\label{sec:S6V-kernel-reform}

The Fredholm determinant \eqref{eq:S6V-fredholm} depends
on the physical parameters $(t, x)$ and the spectral
variable $\zeta$.  An auxiliary integer parameter
$m \in \Z$, which discretizes the spectral variable
through $\zeta_m = \tau^m \zeta$, promotes $(t, x, m)$
to a three-dimensional lattice.  A gauge transformation
of the kernel on this lattice then produces seed
functions and dressing-compatible data to which the
Darboux theorem (Theorem~\ref{thm:Darboux}) applies.

\paragraph{\textbf{Gauge transformation.}}
Define
\[
G_{t,x,m}(z) \defeq z^{-m} g_{t,x}(z)
= z^{-m}(1 + \alpha z)^t(1 + \tau^{-1}z)^{-x},
\]
and set
\[
\zeta_m \defeq \tau^m \zeta,
\qquad
\mathcal{L}_{t,x,m}(\zeta)
\defeq \mathcal{L}_{t,x}(\zeta_m),
\]
\[
K_{t,x,m}(w,w')
\defeq
G_{t,x,m}^{-1}(w)
K_{\zeta_m}^{b_1,b_2}[t,x](w,w')
G_{t,x,m}(w').
\]
The radius condition $\tau < r < 1/\alpha$ ensures
$G_{t,x,m}$ is nonzero on $C_r$, so gauge invariance
of the contour Fredholm determinant gives
\begin{equation}\label{eq:S6V-gauge-det}
\det(I + K_{t,x,m})_{L^2(C_r)}
= \det(I + K_{\zeta_m}^{b_1,b_2}[t,x])_{L^2(C_r)}
= \mathcal{L}_{t,x,m}(\zeta).
\end{equation}
After the gauge transformation, the kernel takes the form
\begin{equation}\label{eq:S6V-gauged-kernel}
K_{t,x,m}(w,w')
= \frac{w^m g_{t,x}(w')}{(w')^m}
\frac{1}{2\mathrm{i}}
\int_{\frac{1}{2}+\mathrm{i}\R}
\frac{(-\zeta_m)^s}
     {\sin(\pi s) g_{t,x}(\tau^s w)}
\frac{\diff s}{\tau^s w - w'}.
\end{equation}

\paragraph{\textbf{Seed data.}}
Define
\begin{equation}\label{eq:S6V-J-def}
J_{t,x,m}(w)
\defeq
\frac{1}{2\mathrm{i}}
\int_{\frac{1}{2} + \mathrm{i}\R}
\frac{(-\zeta_m)^s}{\sin(\pi s)}
\frac{g_{t,x}(w)}{g_{t,x}(\tau^s w)} \diff s,
\end{equation}
and set
\begin{equation}\label{eq:S6V-seed-def}
\psi_{t,x,m}(w)
\defeq G_{t,x,m}^{-1}(w) J_{t,x,m}(w),
\qquad
\phi_{t,x,m}(w')
\defeq -G_{t,x,m+1}(w').
\end{equation}
Explicitly,
\begin{align*}
\psi_{t,x,m}(w)
&= w^m(1 + \alpha w)^{-t}(1 + \tau^{-1}w)^x
J_{t,x,m}(w), \\
\phi_{t,x,m}(w')
&= -(w')^{-(m+1)}(1 + \alpha w')^t
(1 + \tau^{-1}w')^{-x}.
\end{align*}
The product $\psi_{t,x,m}\phi_{t',x',m'}$ denotes the
rank-one operator on $H$ with integral kernel
$\psi_{t,x,m}(w)\phi_{t',x',m'}(w')$ under the contour
measure of Remark~\ref{rem:S6V-fredholm-conv}.

\subsection{Base graph and seed data}
\label{sec:S6V-base-graph}

The lattice structure governing the seed data under
shifts of $(t, x, m)$ is identified next, and the
seed functions $\psi$, $\phi$ are verified to satisfy
the linear problem of the diamond framework.

Let $\mathcal{V} \subseteq \Z^3$ be the lattice
generated by the shifts
\[
T = e^{-\pa_m},
\qquad
S_1 = e^{\pa_x - \pa_m},
\qquad
S_2 = e^{-\pa_t - \pa_m},
\]
so that for $u = (t, x, m) \in \mathcal{V}$,
\[
Tu = (t, x, m{-}1),
\qquad
S_1 u = (t, x{+}1, m{-}1),
\qquad
S_2 u = (t{-}1, x, m{-}1).
\]
Take $E = \mathbb{F} = \C$ and define the
constant scalar edge weights
\begin{equation}\label{eq:S6V-edge-weights}
(c_1, c_2) = (1, 1),
\qquad
(\lambda_1, \lambda_2)
= (-\tau^{-1}, -\alpha).
\end{equation}
Since $(c_k, \lambda_k)$ are constant scalars, the
diamond equations
\eqref{eq:diamond-C-ij}--\eqref{eq:diamond-mixed-ij}
are satisfied trivially: all three reduce to
commutativity of constant scalars in $\C$.

\begin{lemma}\label{lem:S6V-seed-linear}
Let $H = L^2(C_r)$.  The seed functions
$\psi_{t,x,m}$ and $\phi_{t,x,m}$ defined in
\eqref{eq:S6V-seed-def} satisfy the linear problem
\eqref{eq:Linear-Psi-3} and the adjoint linear problem
\eqref{eq:Linear-Phi-Prob} on $\mathcal{V}$ with edge
weights $(c_k, \lambda_k)$.  Explicitly:
\begin{align}
\psi_{t,x+1,m-1}
&= \psi_{t,x,m-1} + \tau^{-1}\psi_{t,x,m},
\label{eq:S6V-psi-S1} \\
\psi_{t-1,x,m-1}
&= \psi_{t,x,m-1} + \alpha\psi_{t,x,m},
\label{eq:S6V-psi-S2}
\end{align}
and
\begin{align}
\phi_{t,x,m-1}
&= \phi_{t,x+1,m-1}
   + \tau^{-1}\phi_{t,x+1,m-2},
\label{eq:S6V-phi-S1} \\
\phi_{t,x,m-1}
&= \phi_{t-1,x,m-1}
   + \alpha\phi_{t-1,x,m-2}.
\label{eq:S6V-phi-S2}
\end{align}
\end{lemma}

\begin{proof}
Both $\psi$-identities rest on the relation
$(-\zeta_{m-1})^s = \tau^{-s}(-\zeta_m)^s$,
which follows from
$\zeta_{m-1} = \tau^{-1}\zeta_m$.
We establish $J$-level recurrences and then lift them
to $\psi$ via the gauge function $G$.
Recall that $g_{t,x}(z) = (1 + \alpha z)^t(1 + \tau^{-1}z)^{-x}$.
The Mellin--Barnes integrals below converge absolutely
and uniformly for $w \in C_r$, so all manipulations may
be performed under the integral sign.

\medskip
\noindent\textit{Verification of
\eqref{eq:S6V-psi-S1}.}
We first establish the $J$-level identity
\begin{equation}\label{eq:S6V-J-S1}
(1 + \tau^{-1}w) J_{t,x+1,m-1}(w)
- J_{t,x,m-1}(w)
= \tau^{-1}w J_{t,x,m}(w).
\end{equation}
The contour representation \eqref{eq:S6V-J-def} reads
\begin{align*}
J_{t,x,m}(w)
&= \frac{1}{2\mathrm{i}}
\int_{\frac{1}{2} + \mathrm{i}\R}
\frac{(-\zeta_m)^s}{\sin(\pi s)} \cdot
\frac{g_{t,x}(w)}{g_{t,x}(\tau^s w)} \diff s.
\end{align*}
To evaluate the first term on the left-hand side of
\eqref{eq:S6V-J-S1}, substitute this representation
for $J_{t,x+1,m-1}$ and apply the $\zeta$-shift
$(-\zeta_{m-1})^s = \tau^{-s}(-\zeta_m)^s$ together
with the $g$-shift
$g_{t,x+1}(z) = (1 + \tau^{-1}z)^{-1} g_{t,x}(z)$.
Since 
\[\frac{g_{t,x+1}(w)}{g_{t,x+1}(\tau^s w)}
= \frac{(1 + \tau^{s-1}w)}{(1 + \tau^{-1}w)}
\frac{g_{t,x}(w)}{g_{t,x}(\tau^s w)},\] 
the result is
\[
(1 + \tau^{-1}w) J_{t,x+1,m-1}(w)
= \frac{1}{2\mathrm{i}}
\int_{\frac{1}{2}+\mathrm{i}\R}
\frac{\tau^{-s}(-\zeta_m)^s}{\sin(\pi s)}
\frac{g_{t,x}(w)}{g_{t,x}(\tau^s w)}
(1 + \tau^{s-1}w) \diff s.
\]
Subtracting the contour representation of
$J_{t,x,m-1}(w)$, the combined integrand carries the
factor $\tau^{-s}(1 + \tau^{s-1}w) - \tau^{-s}
= \tau^{-1}w$, and the remaining integral gives
$\tau^{-1}w J_{t,x,m}(w)$ on the right-hand side,
confirming \eqref{eq:S6V-J-S1}.
To establish \eqref{eq:S6V-psi-S1}, recall that
$\psi_{t,x,m}(w) = G_{t,x,m}^{-1}(w) J_{t,x,m}(w)$.
Since
$G_{t,x+1,m-1}^{-1}(w)
= G_{t,x,m-1}^{-1}(w)(1 + \tau^{-1}w)$,
\begin{align*}
\psi_{t,x+1,m-1}(w)
&= G_{t,x,m-1}^{-1}(w)
   (1 + \tau^{-1}w) J_{t,x+1,m-1}(w) \\
&= G_{t,x,m-1}^{-1}(w)
   \bigl(J_{t,x,m-1}(w)
   + \tau^{-1}w J_{t,x,m}(w)\bigr) \\
&= \psi_{t,x,m-1}(w)
   + \tau^{-1}\psi_{t,x,m}(w),
\end{align*}
where the second equality applies \eqref{eq:S6V-J-S1}
and the third uses
$w G_{t,x,m-1}^{-1}(w) = G_{t,x,m}^{-1}(w)$.

\medskip
\noindent\textit{Verification of
\eqref{eq:S6V-psi-S2}.}
We establish the $J$-level identity
\begin{equation}\label{eq:S6V-J-S2}
(1 + \alpha w) J_{t-1,x,m-1}(w)
- J_{t,x,m-1}(w)
= \alpha w J_{t,x,m}(w).
\end{equation}
To evaluate the first term on the left-hand side of
\eqref{eq:S6V-J-S2}, substitute the contour
representation for $J_{t-1,x,m-1}$ and apply the
$\zeta$-shift together with the $g$-shift
$g_{t-1,x}(z) = (1 + \alpha z)^{-1} g_{t,x}(z)$.
Since $g_{t-1,x}(w)/g_{t-1,x}(\tau^s w)
= \frac{(1 + \alpha\tau^s w)}{(1 + \alpha w)}
\frac{g_{t,x}(w)}{g_{t,x}(\tau^s w)}$, the result is
\[
(1 + \alpha w) J_{t-1,x,m-1}(w)
= \frac{1}{2\mathrm{i}}
\int_{\frac{1}{2}+\mathrm{i}\R}
\frac{\tau^{-s}(-\zeta_m)^s}{\sin(\pi s)}
\frac{g_{t,x}(w)}{g_{t,x}(\tau^s w)}
(1 + \alpha\tau^s w) \diff s.
\]
Subtracting the contour representation of
$J_{t,x,m-1}(w)$, the combined integrand carries the
factor
$\tau^{-s}(1 + \alpha\tau^s w) - \tau^{-s}
= \alpha w$, and the remaining integral gives
$\alpha w J_{t,x,m}(w)$ on the right-hand side,
confirming \eqref{eq:S6V-J-S2}.
To establish \eqref{eq:S6V-psi-S2}, note that
$G_{t-1,x,m-1}^{-1}(w)
= G_{t,x,m-1}^{-1}(w)(1 + \alpha w)$, so
\begin{align*}
\psi_{t-1,x,m-1}(w)
&= G_{t,x,m-1}^{-1}(w)
   (1 + \alpha w) J_{t-1,x,m-1}(w) \\
&= G_{t,x,m-1}^{-1}(w)
   \bigl(J_{t,x,m-1}(w)
   + \alpha w J_{t,x,m}(w)\bigr) \\
&= \psi_{t,x,m-1}(w)
   + \alpha\psi_{t,x,m}(w),
\end{align*}
where the second equality applies \eqref{eq:S6V-J-S2}
and the third uses
$w G_{t,x,m-1}^{-1}(w) = G_{t,x,m}^{-1}(w)$.

\medskip
\noindent\textit{Verification of
\eqref{eq:S6V-phi-S1} and \eqref{eq:S6V-phi-S2}.}
Since 
\[\phi_{t,x,m}(w')
= -G_{t,x,m+1}(w')
= -(w')^{-(m+1)}(1 + \alpha w')^t
(1 + \tau^{-1}w')^{-x},\]
both follow by substituting the explicit formula
for $\phi$ and factoring the common prefactor.
For \eqref{eq:S6V-phi-S1}:
\begin{align*}
\phi_{t,x+1,m-1}(w')
&+ \tau^{-1}\phi_{t,x+1,m-2}(w') \\
&= -(w')^{-m}(1 + \alpha w')^t
   (1 + \tau^{-1}w')^{-(x+1)}
   \bigl(1 + \tau^{-1}w'\bigr) \\
&= \phi_{t,x,m-1}(w').
\end{align*}
For \eqref{eq:S6V-phi-S2}:
\begin{align*}
\phi_{t-1,x,m-1}(w')
&+ \alpha\phi_{t-1,x,m-2}(w') \\
&= -(w')^{-m}(1 + \alpha w')^{t-1}
   (1 + \tau^{-1}w')^{-x}
   \bigl(1 + \alpha w'\bigr) \\
&= \phi_{t,x,m-1}(w'). \qedhere
\end{align*}
\end{proof}

\subsection{Dressing compatibility}
\label{sec:S6V-dressing}

We verify that the gauged kernel $K_{t,x,m}$ satisfies
the rank-one difference identities of
\S\ref{sec:Darboux-transformations}.

\begin{lemma}\label{lem:S6V-dressing}
The kernel $K_{t,x,m}$ satisfies the rank-one
difference identities
\begin{align}
K_{t,x,m-1} - K_{t,x,m}
&= \psi_{t,x,m-1}\phi_{t,x,m-1},
\label{eq:S6V-K-T} \\
K_{t,x+1,m-1} - K_{t,x,m}
&= \psi_{t,x,m-1}\phi_{t,x+1,m-1},
\label{eq:S6V-K-S1} \\
K_{t-1,x,m-1} - K_{t,x,m}
&= \psi_{t,x,m-1}\phi_{t-1,x,m-1}.
\label{eq:S6V-K-S2}
\end{align}
\end{lemma}

In the notation of the framework, these are
$K(Tu) - K(u) = \Psi(Tu) \Phi(Tu)$,
$K(S_1 u) - K(u) = \Psi(Tu) c_1 \Phi(S_1 u)$,
and
$K(S_2 u) - K(u) = \Psi(Tu) c_2 \Phi(S_2 u)$
respectively, with $c_1 = c_2 = 1$.

\begin{proof}
All three identities follow from the same
Cauchy-cancellation mechanism, demonstrated first for
the $T$-shift.

\medskip
\noindent\textit{Verification of \eqref{eq:S6V-K-T}.}
Substituting the gauged kernel
\eqref{eq:S6V-gauged-kernel} at $(t, x, m{-}1)$ and
$(t, x, m)$ and applying the $\zeta$-shift
$(-\zeta_{m-1})^s = \tau^{-s}(-\zeta_m)^s$, the
difference $K_{t,x,m-1} - K_{t,x,m}$ equals
\[
\frac{w^{m-1} g_{t,x}(w')}{(w')^{m-1}}
\frac{1}{2\mathrm{i}}
\int_{\frac{1}{2}+\mathrm{i}\R}
\frac{(-\zeta_m)^s}
     {\sin(\pi s) g_{t,x}(\tau^s w)}
\left(
\frac{\tau^{-s}}{\tau^s w - w'}
- \frac{w}{w'(\tau^s w - w')}
\right) \diff s.
\]
The bracketed term in the integrand equals
$(w'\tau^{-s} - w)/[w'(\tau^s w - w')]$.
Since $w'\tau^{-s} - w = -\tau^{-s}(\tau^s w - w')$,
the factor $(\tau^s w - w')$ cancels between
numerator and denominator, and the difference
reduces to
\[
K_{t,x,m-1} - K_{t,x,m}
= -\frac{w^{m-1} g_{t,x}(w')}{(w')^m}
\frac{1}{2\mathrm{i}}
\int_{\frac{1}{2}+\mathrm{i}\R}
\frac{\tau^{-s}(-\zeta_m)^s}
     {\sin(\pi s) g_{t,x}(\tau^s w)} \diff s.
\]
The prefactor $-g_{t,x}(w')/(w')^m$ equals
$-G_{t,x,m}(w') = \phi_{t,x,m-1}(w')$, and the
remaining factor $w^{m-1}$ times the integral equals
$G_{t,x,m-1}^{-1}(w) J_{t,x,m-1}(w)
= \psi_{t,x,m-1}(w)$
by comparison with \eqref{eq:S6V-J-def},
confirming \eqref{eq:S6V-K-T}.

\medskip
\noindent\textit{Verification of \eqref{eq:S6V-K-S1}.}
Since $g_{t,x+1}(z) = (1 + \tau^{-1}z)^{-1} g_{t,x}(z)$ and
$(-\zeta_{m-1})^s = \tau^{-s}(-\zeta_m)^s$, we
express $K_{t,x+1,m-1}$ in terms of $g_{t,x}$:
\begin{align*}
K_{t,x+1,m-1}(w,w')
= \frac{w^{m-1} g_{t,x}(w')}
       {(w')^{m-1}(1+\tau^{-1}w')}
\frac{1}{2\mathrm{i}}
\int_{\frac{1}{2}+\mathrm{i}\R}
\frac{\tau^{-s}(-\zeta_m)^s
      (1+\tau^{s-1}w)}
     {\sin(\pi s) g_{t,x}(\tau^s w)}
\frac{\diff s}{\tau^s w - w'}.
\end{align*}
Comparing with the gauged kernel
\eqref{eq:S6V-gauged-kernel} at $(t, x, m)$, the
integrand of $K_{t,x+1,m-1} - K_{t,x,m}$ contains
the factor
\[
\frac{w^{m-1}}{(w')^{m-1}}
\left(
\frac{\tau^{-s}(1+\tau^{s-1}w)}
     {(1+\tau^{-1}w')(\tau^s w - w')}
- \frac{w}{w'(\tau^s w - w')}
\right).
\]
Taking the common denominator
$w'(1+\tau^{-1}w')(\tau^s w - w')$,
the numerator is
\begin{align*}
w'\tau^{-s}(1+\tau^{s-1}w)
&- w(1+\tau^{-1}w') \\
&= w'\tau^{-s} + \tau^{-1}ww'
   - w - \tau^{-1}ww' \\
&= w'\tau^{-s} - w.
\end{align*}
Since
$w'\tau^{-s} - w = -\tau^{-s}(\tau^s w - w')$,
the factor $(\tau^s w - w')$ cancels between
numerator and denominator, giving
\begin{align*}
K_{t,x+1,m-1} - K_{t,x,m}
&= -\frac{w^{m-1} g_{t,x}(w')}
        {(w')^m(1+\tau^{-1}w')}
\frac{1}{2\mathrm{i}}
\int_{\frac{1}{2}+\mathrm{i}\R}
\frac{\tau^{-s}(-\zeta_m)^s}
     {\sin(\pi s) g_{t,x}(\tau^s w)} \diff s \\
&= \psi_{t,x,m-1}(w) \phi_{t,x+1,m-1}(w'),
\end{align*}
confirming \eqref{eq:S6V-K-S1}.

\medskip
\noindent\textit{Verification of \eqref{eq:S6V-K-S2}.}
Since $g_{t-1,x}(z) = (1 + \alpha z)^{-1} g_{t,x}(z)$ and
$(-\zeta_{m-1})^s = \tau^{-s}(-\zeta_m)^s$, we
express $K_{t-1,x,m-1}$ in terms of $g_{t,x}$:
\begin{align*}
K_{t-1,x,m-1}(w,w')
= \frac{w^{m-1} g_{t,x}(w')}
       {(w')^{m-1}(1+\alpha w')}
\frac{1}{2\mathrm{i}}
\int_{\frac{1}{2}+\mathrm{i}\R}
\frac{\tau^{-s}(-\zeta_m)^s
      (1+\alpha\tau^s w)}
     {\sin(\pi s) g_{t,x}(\tau^s w)}
\frac{\diff s}{\tau^s w - w'},
\end{align*}
where we used
$g_{t-1,x}(w') = (1+\alpha w')^{-1} g_{t,x}(w')$
and
$g_{t-1,x}(\tau^s w)^{-1}
= (1+\alpha\tau^s w) g_{t,x}(\tau^s w)^{-1}$.
Comparing with the gauged kernel
\eqref{eq:S6V-gauged-kernel} at $(t, x, m)$, the
integrand of $K_{t-1,x,m-1} - K_{t,x,m}$ contains
the factor
\[
\frac{w^{m-1}}{(w')^{m-1}}
\left(
\frac{\tau^{-s}(1+\alpha\tau^s w)}
     {(1+\alpha w')(\tau^s w - w')}
- \frac{w}{w'(\tau^s w - w')}
\right).
\]
The numerator over the common denominator
$w'(1+\alpha w')(\tau^s w - w')$ is
\begin{align*}
w'\tau^{-s}(1+\alpha\tau^s w)
- w(1+\alpha w')
&= w'\tau^{-s} + \alpha ww'
   - w - \alpha ww'
= w'\tau^{-s} - w.
\end{align*}
Since
$w'\tau^{-s} - w = -\tau^{-s}(\tau^s w - w')$,
the factor $(\tau^s w - w')$ cancels between
numerator and denominator, giving
\begin{align*}
K_{t-1,x,m-1} - K_{t,x,m}
&= -\frac{w^{m-1} g_{t,x}(w')}
        {(w')^m(1+\alpha w')}
\frac{1}{2\mathrm{i}}
\int_{\frac{1}{2}+\mathrm{i}\R}
\frac{\tau^{-s}(-\zeta_m)^s}
     {\sin(\pi s) g_{t,x}(\tau^s w)} \diff s \\
&= \psi_{t,x,m-1}(w) \phi_{t-1,x,m-1}(w'),
\end{align*}
confirming \eqref{eq:S6V-K-S2}.
\end{proof}

\subsection{One-point equations}
\label{sec:S6V-equations}

The dressed observable
\[
\mathcal{M}(u)
= 1 + z \phi(u) (I - zK(u))^{-1} \psi(u)
\in \C
\]
has dressed edge weights
\[
\mathcal{M}_k(u)
= \mathcal{M}(Tu)^{-1} \mathcal{M}(S_k u)
\]
(since $c_1 = c_2 = 1$).  We write
$\mathcal{M}_{t,x,m} \defeq \mathcal{M}(t, x, m)$
throughout.

\begin{theorem}\label{thm:S6V-H9}
The dressed edge weights satisfy the mixed diamond
equation
\begin{align}\label{eq:S6V-H9}
\bigl(\mathcal{M}_{t,x-1,m}^{-1}
- \tfrac{1-b_2}{1-b_1}
\mathcal{M}_{t+1,x,m}^{-1}\bigr)
\mathcal{M}_{t,x,m}
- \mathcal{M}_{t+1,x-1,m+1}^{-1}
\bigl(\mathcal{M}_{t+1,x,m+1}
- \tfrac{1-b_2}{1-b_1}
\mathcal{M}_{t,x-1,m+1}\bigr) = 0
\end{align}
at each $z \in \C$ and
$\zeta \in \C \setminus \R_{\ge 0}$ for which all
displayed $\mathcal{M}$-values are defined.
\end{theorem}

\begin{proof}
The hypothesis ensures the resolvent exists at the
shifted vertices entering the diamond stencil, so by
Theorem~\ref{thm:Darboux}
(Remark~\ref{rem:resolvent-subset}) the pair
$(\mathcal{M}_k, \lambda_k)$ satisfies the diamond
equations.
Since $\Lambda_k = \lambda_k$ are constants, the
mixed diamond equation for $(i, j) = (1, 2)$ reduces to
\[
\lambda_2\bigl[\mathcal{M}_1(u)
- \mathcal{M}_1(S_2 u)\bigr]
+ \lambda_1\bigl[\mathcal{M}_2(S_1 u)
- \mathcal{M}_2(u)\bigr] = 0.
\]
Substituting the explicit shifts and dividing by
$-\alpha$ gives
\begin{align*}
\bigl[\mathcal{M}_{t,x,m-1}^{-1}
&\mathcal{M}_{t,x+1,m-1}
- \mathcal{M}_{t-1,x,m-2}^{-1}
\mathcal{M}_{t-1,x+1,m-2}\bigr] \\
+ \tfrac{1-b_2}{1-b_1}
&\bigl[\mathcal{M}_{t,x+1,m-2}^{-1}
\mathcal{M}_{t-1,x+1,m-2}
- \mathcal{M}_{t,x,m-1}^{-1}
\mathcal{M}_{t-1,x,m-1}\bigr] = 0.
\end{align*}
Evaluating at the vertex $(t{+}1, x{-}1, m{+}2)$
in place of $(t, x, m)$, the middle two terms share
the right factor $\mathcal{M}_{t,x,m}$ and the
first and last share the left factor
$\mathcal{M}_{t+1,x-1,m+1}^{-1}$.  Factoring and
multiplying by $-1$ yields \eqref{eq:S6V-H9}.
\end{proof}

\begin{remark}\label{rem:S6V-H9-continuation}
At $z = -1$, Corollary~\ref{cor:Woodbury} gives
\[
\mathcal{M}_{t,x,m}
= \frac{\mathcal{L}_{t,x,m+1}}
       {\mathcal{L}_{t,x,m}},
\]
so \eqref{eq:S6V-H9} becomes a rational identity in
$\tau$-Laplace transforms.
Corollary~\ref{cor:S6V-scalar} establishes the
corresponding bilinear identity via the scalar
Hirota--Miwa reduction.
\end{remark}

\begin{corollary}\label{cor:S6V-scalar}
For $t \ge 1$ and
$\zeta \in \C \setminus \R_{\ge 0}$, the Fredholm
determinant $\mathcal{L}_{t,x,m}(\zeta)
= \det(I + K_{t,x,m})$ satisfies
\begin{align}\label{eq:S6V-scalar-HM}
&(1 - b_1)
\mathcal{L}_{t,x+1,m+1}
\mathcal{L}_{t-1,x,m}
- (1 - b_2)
\mathcal{L}_{t-1,x,m+1}
\mathcal{L}_{t,x+1,m} \notag\\
&\qquad+ (b_1 - b_2)
\mathcal{L}_{t,x,m+1}
\mathcal{L}_{t-1,x+1,m}
= 0.
\end{align}
\end{corollary}

\begin{proof}
Set $F_z(u) \defeq \det(I - zK(u))$ for $z \in \C$;
the Fredholm series defines an entire function of $z$
for each $u$.  For $|z|$ sufficiently small, the
resolvents at the diamond stencil and the $T$-orbit
exist, and the boundary condition
$F_z(T^\ell v) \to 1$ as $\ell \to -\infty$ holds for
every $z$ because $\tau^{m-\ell}\zeta \to 0$.  In the
notation of Proposition~\ref{prop:scalar-HM}, the
Hirota coefficients are
$\alpha_{12} = c_1\lambda_2 = -\alpha$ and
$\alpha_{21} = c_2\lambda_1 = -\tau^{-1}$; these are
distinct since $b_1 \ne b_2$.
Proposition~\ref{prop:scalar-HM} therefore gives
\eqref{eq:S6V-scalar-HM} with $\mathcal{L}$ replaced
by $F_z$, after substituting the lattice points into
\eqref{eq:HM-variable}, multiplying by
$-\tau(1 - b_2)$, and shifting $m \mapsto m + 2$.
Each term in this identity is a product of two
$F_z$-values, hence entire in $z$.  Since the identity
holds for $|z|$ small, the identity theorem gives it
for all $z \in \C$.  At $z = -1$,
$F_{-1}(u) = \mathcal{L}_{t,x,m}(\zeta)$ by
\eqref{eq:S6V-gauge-det}, recovering
\eqref{eq:S6V-scalar-HM}.
\end{proof}

%% file: chapters/7-euclidean-division-vertex-models-and-positive-polymer-models.tex
\chapter{Euclidean Division, Vertex Models, and Positive Polymer Models}
\label{ch:euclidean-division-vertex-models-and-positive-polymer-models}

In the stochastic six-vertex model (\S\ref{sec:S6V}),
the gauge function has one linear factor per shift
direction, and each dressing identity follows from a
single rank-one Cauchy cancellation.  Other models have
many linear factors per shift and correspondingly many
rank-one Darboux steps that fuse into a single
matrix-valued move.  The natural algebraic framework
for this fusion is the polynomial quotient
$E = \mathbb{F}[w]/(\Pi)$, where the polynomial $\Pi$
collects all the linear factors.

\section{Euclidean Division}\label{sec:PQ}

This section establishes that the Euclidean division
operators of $E = \mathbb{F}[w]/(\Pi)$ satisfy the diamond
equations of the framework, and derives the Bezoutian
identities used in the seed and dressing verifications
of the model sections that
follow.\footnote{The reader may prefer to proceed
directly to \S\ref{sec:HSEP} and refer back to this
section as needed.}

\paragraph{\textbf{Quotient algebra.}}
Fix a field $\mathbb{F}$ and a polynomial
$\Pi \in \mathbb{F}[w]$ of degree $N \geq 1$, and set
\begin{equation}\label{eq:PQ-E-def}
E \defeq \mathbb{F}[w]/(\Pi).
\end{equation}
Recall that $E$ is the ring whose elements are equivalence
classes of polynomials in $\mathbb{F}[w]$: two
polynomials $f, g$ represent the same element of $E$ if
and only if $\Pi \mid (f - g)$.

Polynomial long division writes any
$g \in \mathbb{F}[w]$ uniquely as $g = q\Pi + r$ with
$\deg r < N$.  Since $\Pi \mid (g - r)$, the
polynomials $g$ and $r$ represent the same element of
$E$.  The representative $r$ is unique: if two
polynomials $r_1, r_2$ of degree less than $N$ represent
the same element, then $\Pi \mid (r_1 - r_2)$, but
$\deg(r_1 - r_2) < N = \deg \Pi$, forcing
$r_1 = r_2$.

Each element of $E$ therefore has a unique polynomial
representative of degree less than $N$, and the
monomials $1, w, \ldots, w^{N-1}$ form a basis for $E$
over $\mathbb{F}$.  Multiplication in $E$ is polynomial
multiplication followed by reduction modulo $\Pi$.
When the roots of $\Pi$ are all distinct, $E$ is
semisimple as an $\mathbb{F}$-algebra.

\subsection{Diamond equations from Euclidean division}
\label{sec:PQ-quotient}

Fix $\Pi \in \mathbb{F}[w]$ of degree $N$ and
$E = \mathbb{F}[w]/(\Pi)$ as above.
For any polynomial $p \in \mathbb{F}[w]$ with
$\deg p \leq N$ and any $f \in E$, define
$C_p f, \Lambda_p f \in E$ by
\begin{equation}\label{eq:PQ-division}
p(w) f(w) = \Pi(w)(C_p f)(w) - (\Lambda_p f)(w).
\end{equation}
That is, $C_p f$ is the quotient and $-\Lambda_p f$
the remainder of Euclidean division of $pf$ by $\Pi$.
Since $\deg(pf) < 2N$, the quotient and remainder both
have degree less than $N$, so
$C_p, \Lambda_p \colon E \to E$ are well-defined
endomorphisms.

\medskip

The commutativity of polynomial multiplication forces
the Euclidean division edge weights to satisfy the
diamond equations.

\begin{lemma}
\label{thm:PQ-bare-diamond}
For any $p, q \in \mathbb{F}[w]$ with
$\deg p, \deg q \leq N$, the operators
$C_p, \Lambda_p$ and $C_q, \Lambda_q$ satisfy
\begin{align}
C_p C_q &= C_q C_p,
\label{eq:PQ-CC} \\
\Lambda_p \Lambda_q &= \Lambda_q \Lambda_p,
\label{eq:PQ-LL} \\
\Lambda_p C_q + C_p \Lambda_q
&=
\Lambda_q C_p + C_q \Lambda_p.
\label{eq:PQ-mixed}
\end{align}
These are the diamond equations
\eqref{eq:diamond-C-ij}--\eqref{eq:diamond-mixed-ij}.
\end{lemma}

\begin{proof}
The product $pqf$ is computed in two orders.  Applying
\eqref{eq:PQ-division} first with $q$ and then with $p$
gives
\[
pqf
=
\Pi^2 C_p C_q f
-
\Pi(\Lambda_p C_q + C_p \Lambda_q) f
+
\Lambda_p \Lambda_q f.
\]
Interchanging $p$ and $q$ gives
\[
qpf
=
\Pi^2 C_q C_p f
-
\Pi(\Lambda_q C_p + C_q \Lambda_p) f
+
\Lambda_q \Lambda_p f.
\]
Since $pqf = qpf$, subtracting the two expansions and
collecting powers of $\Pi$ gives
\begin{align*}
\Pi^2 (C_pC_q - C_qC_p) f
- \Pi (\Lambda_pC_q + C_p\Lambda_q
  - \Lambda_qC_p - C_q\Lambda_p) f
+ (\Lambda_p\Lambda_q - \Lambda_q\Lambda_p) f
= 0.
\end{align*}
Each of the three bracketed operators maps $E$ to $E$,
so for each $f \in E$ this is a polynomial identity of
the form $\Pi^2 a - \Pi b + d = 0$ with
$a, b, d \in E$.  Such a representation is unique: if
$\Pi^2 a - \Pi b + d = 0$ with $\deg a, \deg b, \deg d
< N$, then reducing modulo $\Pi$ forces $d = 0$;
dividing by $\Pi$ and reducing again forces $b = 0$;
and then $a = 0$.  The three bracketed operators
therefore vanish identically on $E$, which gives
\eqref{eq:PQ-CC}--\eqref{eq:PQ-mixed}.
\end{proof}

\begin{remark}\label{rem:PQ-commutativity}
The diamond equations hold for any $\Pi$ of degree $N$;
no hypothesis on the roots is required.
\end{remark}

\begin{example}\label{ex:PQ-diamond}
Let $\Pi(w) = w^3 - 6w^2 + 11w - 6 = (w-1)(w-2)(w-3)$,
so $E$ is three-dimensional with monomial basis
$\{1, w, w^2\}$.  Take $p(w) = w - 1$ and
$q(w) = w^2 + 1$; note that $q$ does not divide $\Pi$.
In the monomial basis,
\[
C_p = \begin{pmatrix} 0 & 0 & 1 \\ 0 & 0 & 0 \\
0 & 0 & 0 \end{pmatrix}, \qquad
\Lambda_p = \begin{pmatrix} 1 & 0 & -6 \\ -1 & 1 & 11 \\
0 & -1 & -5 \end{pmatrix},
\]
\[
C_q = \begin{pmatrix} 0 & 1 & 6 \\ 0 & 0 & 1 \\
0 & 0 & 0 \end{pmatrix}, \qquad
\Lambda_q = \begin{pmatrix} -1 & -6 & -36 \\ 0 & 10 & 60 \\
-1 & -6 & -26 \end{pmatrix}.
\]
The individual products $C_p \Lambda_q$ and
$\Lambda_q C_p$ are both nonzero:
\[
C_p \Lambda_q = \begin{pmatrix} -1 & -6 & -26 \\
0 & 0 & 0 \\ 0 & 0 & 0 \end{pmatrix}, \qquad
\Lambda_q C_p = \begin{pmatrix} 0 & 0 & -1 \\
0 & 0 & 0 \\ 0 & 0 & -1 \end{pmatrix},
\]
so $[C_p, \Lambda_q] \neq 0$.  Yet all three diamond
equations hold: $C_p C_q = C_q C_p$,
$\Lambda_p \Lambda_q = \Lambda_q \Lambda_p$, and
$\Lambda_p C_q + C_p \Lambda_q
= \Lambda_q C_p + C_q \Lambda_p$.
\end{example}

\begin{remark}\label{rem:PQ-rank}
For nonzero $p$, the image of $C_p$ consists of
polynomials of degree less than $\deg p$, so
$\operatorname{rank} C_p = \deg p$.
The map $\Lambda_p$ is minus multiplication by $p$ in
$E$, so $\operatorname{rank} \Lambda_p
= N - \deg\gcd(p, \Pi)$; in particular, $\Lambda_p$ is
invertible if and only if $\gcd(p, \Pi) = 1$.

When $\Pi$ splits over $\mathbb{F}$ with distinct roots
$\alpha_1, \dotsc, \alpha_N$, the Chinese remainder
theorem identifies $E \cong \bigoplus_{i=1}^N \mathbb{F}$.
Under this identification, $\Lambda_p$ acts diagonally
with eigenvalues $-p(\alpha_i)$, giving
$\det \Lambda_p = (-1)^N \prod_{i=1}^N p(\alpha_i)$.
\end{remark}

\subsection{Bezoutian identities}
\label{sec:PQ-Bezoutian}

The seed and dressing verifications in the model sections
below require the action of $C_p$ and $\Lambda_p$ on the
Bezoutian polynomial of $\Pi$.  The two identities
established here are used in both the HSEP and IS6V
proofs.

For any polynomial $P$, the Bezoutian is
\begin{equation}\label{eq:PQ-Bezoutian-def}
B_P(x, y) \defeq
\frac{P(x) - P(y)}{x - y},
\end{equation}
a polynomial of degree $\deg P - 1$ in each variable.
In particular, $B_\Pi(\cdot, y) \in E$ for each
$y$.

\begin{lemma}
\label{lem:PQ-Bezoutian}
Let $p \in \mathbb{F}[w]$ with $\deg p \leq N$.  Then
\begin{align}
C_p\bigl[B_\Pi(\cdot, y)\bigr](x)
&= B_p(x, y),
\label{eq:PQ-C-Bezoutian} \\
\Lambda_p\bigl[B_\Pi(\cdot, y)\bigr](x)
&= \Pi(y) B_p(x, y) - p(y)
B_\Pi(x, y).
\label{eq:PQ-Lambda-Bezoutian}
\end{align}
\end{lemma}

\begin{proof}
\noindent\textit{Identity \eqref{eq:PQ-C-Bezoutian}.}
The Euclidean division of $p(x) B_\Pi(x, y)$ by
$\Pi(x)$ reads $p(x) B_\Pi(x, y) = \Pi(x) Q(x) + R(x)$
with $\deg R < N$; we show $Q = B_p(x, y)$ by verifying
that $p(x) B_\Pi(x, y) - \Pi(x) B_p(x, y)$ has degree
less than $N$:
\begin{align*}
p(x) B_\Pi(x, y)
- \Pi(x) B_p(x, y)
&=
\frac{p(x)\bigl(\Pi(x) - \Pi(y)\bigr)
     - \Pi(x)\bigl(p(x) - p(y)\bigr)}
{x - y} \\
&=
\frac{\Pi(x) p(y) - p(x) \Pi(y)}
{x - y}.
\end{align*}
The numerator vanishes at $x = y$, so $(x - y)$ divides
it and the right-hand side is a polynomial in $x$ of
degree less than $N$, giving \eqref{eq:PQ-C-Bezoutian}.

\medskip
\noindent\textit{Identity \eqref{eq:PQ-Lambda-Bezoutian}.}
By \eqref{eq:PQ-division} and \eqref{eq:PQ-C-Bezoutian}, we have 
$\Lambda_p[B_\Pi(\cdot, y)](x)
= \Pi(x) B_p(x, y) - p(x) B_\Pi(x, y)$,
which is the negative of the remainder computed above
and equals $(p(x) \Pi(y) - \Pi(x) p(y))/(x - y)$.
Expanding the right-hand side of
\eqref{eq:PQ-Lambda-Bezoutian} in the same way gives
the same expression.
\end{proof}

\section{Higher-Spin Exclusion Process}\label{sec:HSEP}

\paragraph{\textbf{System description.}}
Fix $0 < q < 1$, $\alpha > 0$, and
$I, J \in \Z_{\geq 1}$.
The $J$-higher-spin exclusion process of
Corwin--Petrov~\cite{CorwinPetrov2016}
is a discrete-time Markov chain on strictly ordered
particle configurations
$\vec{x}(t) = (x_1(t) > x_2(t) > \cdots) \in \Z^\infty$.
At each time step, particles update sequentially
starting from $x_1$, each jumping right by up to $J$
sites subject to exclusion; the transition
probabilities are given by the stochastic higher-spin
vertex weights
of~\cite[Sec.~2--3]{CorwinPetrov2016}.
Step initial data $x_i(0) = -i$ is imposed
throughout.

The parameters
\begin{equation}\label{eq:HSEP-nu-qI}
\nu \defeq q^I, \qquad \beta \defeq \alpha q^J,
\qquad I, J \in \Z_{\geq 1},
\end{equation}
correspond to Case~(1)
of~\cite[Prop.~2.2]{CorwinPetrov2016}.
For $\zeta \in \C \setminus \R_{\geq 0}$, define the
$e_q$-Laplace transform
\begin{equation}\label{eq:HSEP-L-def}
\mathcal{L}_{n,t}(\zeta)
\defeq
\E_{\mathrm{step}}
\left[
\frac{1}{(\zeta q^{x_n(t)+n};q)_\infty}
\right],
\qquad
(a;q)_\infty
\defeq \prod_{k=0}^\infty (1 - aq^k).
\end{equation}

\medskip
\paragraph{\textbf{Fredholm determinant formula.}}
The following is due to
Corwin--Petrov \cite[Thm.~4.2]{CorwinPetrov2016}.
For $n \geq 1$, $t \geq 1$, and
$\zeta \in \C \setminus \R_{\geq 0}$,
\begin{equation}\label{eq:HSEP-fredholm}
\mathcal{L}_{n,t}(\zeta)
=
\det\bigl(I + K_\zeta[n,t]\bigr)_{L^2(C_r)},
\end{equation}
where $C_r$ is a positively oriented circle centred at
$1$ with radius $r$ small enough to exclude $0$,
$q^{-1}$, and $q^{-I}$, and the Mellin--Barnes kernel is
\begin{equation}\label{eq:HSEP-kernel-original}
K_\zeta[n,t](w,w')
=
\frac{g_{n,t}(w)}{2\pi \mathrm{i}}
\int_{\frac{1}{2}+\mathrm{i}\R}
\frac{\pi}{\sin(-\pi s)}
(-\zeta)^s
\frac{1}{g_{n,t}(q^s w)}
\frac{\diff s}{q^s w - w'},
\end{equation}
with
\begin{equation}\label{eq:HSEP-g-def}
g_{n,t}(z)
=
\left(
\frac{(\nu z;q)_\infty}{(z;q)_\infty}
\right)^{\!n}
\left(
\frac{(-\beta z;q)_\infty}{(-\alpha z;q)_\infty}
\right)^{\!t}
\frac{1}{(\nu z;q)_\infty}.
\end{equation}
The $s$-contour is the vertical line
$\operatorname{Re}(s) = 1/2$ oriented bottom-to-top,
and $(-\zeta)^s$ is defined using the principal branch
of $\log(-\zeta)$.

\medskip
\paragraph{\textbf{References.}}
The Fredholm determinant formula
\eqref{eq:HSEP-fredholm} is
Corwin--Petrov \cite[Thm.~4.2]{CorwinPetrov2016},
specialized to $\nu = q^I$ and $\beta = \alpha q^J$.

\begin{remark}[Fredholm determinant convention]
\label{rem:HSEP-fredholm-conv}
The Fredholm determinant
$\det(I + K_\zeta[n,t])$ is taken on $L^2(C_r)$ with
the contour measure
$\diff w / (2\pi\mathrm{i})$.  The combined factor
$\pi(-\zeta)^s/\sin(-\pi s)$ decays exponentially on
the vertical $s$-contour for each
$\zeta \in \C \setminus \R_{\geq 0}$, giving uniform
convergence of the Mellin--Barnes integral on
$C_r \times C_r$.  The kernel is therefore continuous
and bounded on the contour square, and the Fredholm
series converges absolutely by the Hadamard estimate.
Since $(-\zeta)^s$ is analytic in $\zeta$ on
$\C \setminus \R_{\geq 0}$ (principal branch), the
same estimates give analyticity of
$\mathcal{L}_{n,t}(\zeta)$ in $\zeta$ on this domain.
When the determinant is nonzero, the resolvent kernel
exists by classical Fredholm theory.  The framework
results invoked below (gauge invariance, the Woodbury
determinant ratio, and the Darboux theorem) hold at
the kernel and series level: their proofs require only
pointwise kernel identities, the resolvent integral
equation, and the rank-one perturbation formula for
Fredholm series.
\end{remark}

\subsection{Kernel reformulation}\label{sec:HSEP-kernel-reform}

The Fredholm determinant \eqref{eq:HSEP-fredholm}
depends on the physical parameters $(n, t)$ and the
spectral variable $\zeta$.  An auxiliary integer
parameter $m \in \Z$, which discretizes the spectral
variable through $\zeta_m \defeq q^m \zeta$, promotes
$(n, t, m)$ to a lattice on which the Darboux theorem
applies.  For $u = (n, t, m)$, define $K_u$ by
replacing $\zeta$ in \eqref{eq:HSEP-kernel-original}
by $\zeta_m$, so that
\begin{equation}\label{eq:HSEP-F-def}
F(u) \defeq \det(I + K_u)_{L^2(C_r)}
=
\E_{\mathrm{step}}
\left[
\frac{1}{(\zeta_m q^{x_n(t)+n};q)_\infty}
\right].
\end{equation}

\paragraph{\textbf{Gauge transformation.}}
The gauge transformation
$G_u(z) \defeq z^{-m} g_{n,t}(z)$ conjugates the
kernel to
\[
K_u(w,w')
= G_u(w)^{-1} K_{\zeta_m}[n,t](w,w') G_u(w').
\]
Since $G_u$ is bounded and nonzero on $C_r$, gauge
invariance gives $\det(I + K_u) = F(u)$.  Equivalently,
\begin{equation}\label{eq:HSEP-gauged-kernel}
K_u(w,w')
=
\frac{w^m g_{n,t}(w')}{(w')^m}
\frac{1}{2\pi\mathrm{i}}
\int_{\frac{1}{2}+\mathrm{i}\R}
\frac{\pi}{\sin(-\pi s)}
\frac{(-\zeta_m)^s}{g_{n,t}(q^s w)}
\frac{\diff s}{q^s w - w'}.
\end{equation}
The notation $K_u$ denotes the gauged kernel from this
point forward.  Separating the integrand into factors
depending on $(w, s)$ and on $w'$ isolates how the
lattice shifts act: the gauged kernel factors as
\begin{equation}\label{eq:HSEP-XY-factor}
K_u(w,w') = \int_\Gamma X_u(w,s)
\frac{Y_u(w')}{q^s w - w'} \diff s,
\end{equation}
where $\Gamma = \frac{1}{2} + \mathrm{i}\R$ and
\[
X_u(w,s)
\defeq
w^m
\frac{\pi}{2\pi\mathrm{i}\sin(-\pi s)}
\frac{(-\zeta_m)^s}{g_{n,t}(q^s w)},
\qquad
Y_u(w') \defeq \frac{g_{n,t}(w')}{(w')^m}.
\]

\paragraph{\textbf{Shift polynomials.}}
Each physical shift collapses one of the $q$-Pochhammer
ratios in $g_{n,t}$ to a finite polynomial.  Since
$\nu = q^I$, the particle shift gives
\begin{equation}\label{eq:HSEP-p1}
\frac{g_{n+1,t}(z)}{g_{n,t}(z)}
= \frac{(\nu z;q)_\infty}{(z;q)_\infty}
= \frac{1}{p_1(z)},
\qquad
p_1(z) \defeq \prod_{r=0}^{I-1}(1 - q^r z),
\end{equation}
by the telescoping identity
$(az;q)_\infty / (aqz;q)_\infty = 1 - az$.
The same identity with $\beta = \alpha q^J$ gives
\begin{equation}\label{eq:HSEP-p2}
\frac{g_{n,t+1}(z)}{g_{n,t}(z)}
= \frac{(-\beta z;q)_\infty}{(-\alpha z;q)_\infty}
= \frac{1}{p_2(z)},
\qquad
p_2(z) \defeq \prod_{r=0}^{J-1}(1 + \alpha q^r z).
\end{equation}
The spectral shift $m \to m + 1$ contributes the
additional linear factor $z$.  The fusion polynomial
\begin{equation}\label{eq:HSEP-Pi}
\Pi(z) \defeq z p_1(z) p_2(z),
\qquad
\deg \Pi = I + J + 1,
\end{equation}
is the product that defines the quotient algebra
$E = \C[z]/(\Pi)$ of \S\ref{sec:PQ}.

The three lattice shifts are
\begin{gather}\label{eq:HSEP-lattice-shifts}
T(n,t,m) = (n{+}1, t{+}1, m{+}1),
\quad
S_1(n,t,m) = (n{+}1, t, m),
\quad
S_2(n,t,m) = (n, t{+}1, m).
\end{gather}
The $X$--$Y$ factors of \eqref{eq:HSEP-XY-factor}
transform under each shift by its associated polynomial:
\begin{equation}\label{eq:HSEP-shift-rule}
\begin{aligned}
X_{S_k u}(w,s) &= X_u(w,s) p_k(q^s w), &\qquad
Y_{S_k u}(w') &= \frac{Y_u(w')}{p_k(w')}, &\qquad
k &\in \{1,2\}, \\
X_{Tu}(w,s) &= X_u(w,s) \Pi(q^s w), &\qquad
Y_{Tu}(w') &= \frac{Y_u(w')}{\Pi(w')}.
\end{aligned}
\end{equation}

\subsection{Base graph and seed data}
\label{sec:HSEP-base-graph}

Let $\mathcal{V} \subseteq \Z^3$ be the lattice
generated by the shifts
\eqref{eq:HSEP-lattice-shifts}.  The quotient algebra
$E = \C[z]/(\Pi)$ of \S\ref{sec:PQ} has dimension
$I + J + 1$.  Write
$(C_k, \Lambda_k) \defeq (C_{p_k}, \Lambda_{p_k})$
for the Euclidean division edge weights of
\eqref{eq:PQ-division}; these satisfy the bare diamond
equations by Lemma~\ref{thm:PQ-bare-diamond}.

Let $H = L^2(C_r)$.  Define the wave function
$\Psi(u) \colon E \to H$ by
\begin{equation}\label{eq:HSEP-Psi-def}
(\Psi(u) f)(w)
\defeq
\int_\Gamma X_u(w,s) \frac{f(q^s w)}{\Pi(q^s w)} \diff s.
\end{equation}
At the vertex $Tu$, the factor $\Pi(q^s w)$ from the
shift rule \eqref{eq:HSEP-shift-rule} cancels the
denominator in \eqref{eq:HSEP-Psi-def}:
\begin{equation}\label{eq:HSEP-Psi-Tu}
(\Psi(Tu) f)(w)
= \int_\Gamma X_u(w,s) f(q^s w) \diff s.
\end{equation}
Define the adjoint wave function
$\Phi(u) \colon H \to E$ by
\begin{equation}\label{eq:HSEP-Phi-def}
(\Phi(u) h)(x)
\defeq
\int_{C_r} B_\Pi(x,y) Y_u(y) h(y)
\frac{\diff y}{2\pi\mathrm{i}}.
\end{equation}
Since $B_\Pi(\cdot, y)$ has degree $I + J$ in $x$,
$\Phi(u) h$ is a polynomial of degree less than
$\dim E$ and hence an element of $E$.

\begin{lemma}\label{lem:HSEP-seed-linear}
The seed functions $\Psi$ and $\Phi$ satisfy the linear
problem \eqref{eq:Linear-Psi-3} and adjoint linear
problem \eqref{eq:Linear-Phi-Prob} on $\mathcal{V}$
with edge weights $(C_k, \Lambda_k)$.  Explicitly, for
$k \in \{1, 2\}$:
\begin{equation}\label{eq:HSEP-Psi-recurrence}
\Psi(S_k u) = \Psi(Tu) C_k - \Psi(u) \Lambda_k,
\end{equation}
and
\begin{equation}\label{eq:HSEP-Phi-recurrence}
\Phi(Tu) = C_k \Phi(S_k u) - \Lambda_k \Phi(T S_k u).
\end{equation}
\end{lemma}

\begin{proof}
The Mellin--Barnes integrals below converge absolutely
and uniformly for $w \in C_r$, so all manipulations
may be performed under the integral sign.

\noindent\textit{Verification of
\eqref{eq:HSEP-Psi-recurrence}.}
The shift rule \eqref{eq:HSEP-shift-rule} applied to
\eqref{eq:HSEP-Psi-def} gives
\[
(\Psi(S_k u) f)(w)
=
\int_\Gamma X_u(w,s)
\frac{p_k(q^s w) f(q^s w)}{\Pi(q^s w)} \diff s.
\]
Dividing \eqref{eq:PQ-division} by $\Pi(q^s w)$ gives
\[
\frac{p_k(q^s w) f(q^s w)}{\Pi(q^s w)}
= (C_k f)(q^s w)
- \frac{(\Lambda_k f)(q^s w)}{\Pi(q^s w)},
\]
so the integral separates into
$(\Psi(Tu) C_k f)(w) - (\Psi(u) \Lambda_k f)(w)$
by \eqref{eq:HSEP-Psi-Tu} and
\eqref{eq:HSEP-Psi-def}.

\medskip
\noindent\textit{Verification of
\eqref{eq:HSEP-Phi-recurrence}.}
Write $\Phi(u; x, y) \defeq B_\Pi(x, y) Y_u(y)$ for
the integral kernel of $\Phi(u)$, so that
$(\Phi(u)h)(x) = \int_{C_r} \Phi(u; x, y) h(y)
\frac{\diff y}{2\pi\mathrm{i}}$.
The recurrence \eqref{eq:HSEP-Phi-recurrence} reduces to
an identity of these kernels, since the shift rule
\eqref{eq:HSEP-shift-rule} acts only on the $Y_u(y)$
factor.  Evaluating at the vertices $Tu$, $S_k u$, and
$TS_k u$:
\begin{gather*}
\Phi(Tu; x, y) = B_\Pi(x,y)
\frac{Y_u(y)}{\Pi(y)},
\qquad
\Phi(S_k u; x, y) = B_\Pi(x,y)
\frac{Y_u(y)}{p_k(y)}, \\
\Phi(TS_k u; x, y) = B_\Pi(x,y)
\frac{Y_u(y)}{\Pi(y) p_k(y)}.
\end{gather*}
Since $Y_u(y)/p_k(y)$ does not depend on $x$, it
factors from $C_k$ and $\Lambda_k$.  Substituting
$C_k[B_\Pi(\cdot, y)] = B_{p_k}(\cdot, y)$ from
\eqref{eq:PQ-C-Bezoutian} and
$\Lambda_k[B_\Pi(\cdot, y)]
= \Pi(y) B_{p_k}(\cdot, y) - p_k(y) B_\Pi(\cdot, y)$
from \eqref{eq:PQ-Lambda-Bezoutian}:
\begin{align*}
C_k[\Phi(S_k u; \cdot, y)](x)
&- \Lambda_k[\Phi(TS_k u; \cdot, y)](x) \\
&=
\frac{Y_u(y)}{p_k(y)} B_{p_k}(x,y)
-
\frac{Y_u(y)}{p_k(y)} B_{p_k}(x,y)
+
\frac{Y_u(y)}{\Pi(y)} B_\Pi(x,y) \\
&= \Phi(Tu; x, y). \qedhere
\end{align*}
\end{proof}

\subsection{Dressing compatibility}
\label{sec:HSEP-dressing}

\begin{lemma}\label{lem:HSEP-dressing}
The gauged kernel $K_u$ satisfies the dressing
compatibility conditions:
\begin{align}
K(Tu)(w,w') - K(u)(w,w')
&= (\Psi(Tu) \Phi(Tu))(w,w'),
\label{eq:HSEP-K-T} \\
K(S_k u)(w,w') - K(u)(w,w')
&= (\Psi(Tu) C_k \Phi(S_k u))(w,w'),
\qquad k \in \{1, 2\}.
\label{eq:HSEP-K-Sk}
\end{align}
\end{lemma}

\begin{proof}
\noindent\textit{Verification of \eqref{eq:HSEP-K-T}.}
The shift rule \eqref{eq:HSEP-shift-rule} gives
$X_{Tu} = X_u \Pi(q^s w)$ and
$Y_{Tu} = Y_u / \Pi(w')$, so
\begin{align*}
K(Tu)(w,w') - K(u)(w,w')
&=
\frac{Y_u(w')}{\Pi(w')}
\int_\Gamma X_u(w,s)
\frac{\Pi(q^s w) - \Pi(w')}{q^s w - w'} \diff s \\
&=
\frac{Y_u(w')}{\Pi(w')}
\int_\Gamma X_u(w,s) B_\Pi(q^s w, w') \diff s.
\end{align*}
The composition $\Psi(Tu)\Phi(Tu) \in \End(H)$ has kernel
\begin{align*}
(\Psi(Tu)\Phi(Tu))(w, w')
&= (\Psi(Tu)[\Phi(Tu; \cdot, w')])(w) \\
&= \frac{Y_u(w')}{\Pi(w')}
\int_\Gamma X_u(w,s) B_\Pi(q^s w, w') \diff s,
\end{align*}
since for each fixed $w'$ the kernel
$\Phi(Tu; \cdot, w') \in E$, the scalar
$Y_u(w')/\Pi(w')$ factors from $\Psi(Tu)$, and the
remaining integral is \eqref{eq:HSEP-Psi-Tu} applied to
$B_\Pi(\cdot, w')$.  This is the expression computed for
$K(Tu)(w,w') - K(u)(w,w')$, confirming
\eqref{eq:HSEP-K-T}.

\medskip
\noindent\textit{Verification of \eqref{eq:HSEP-K-Sk}.}
The shift rule \eqref{eq:HSEP-shift-rule} gives
$X_{S_k u} = X_u p_k(q^s w)$ and
$Y_{S_k u} = Y_u / p_k(w')$, so
\begin{align*}
K(S_k u)(w,w') - K(u)(w,w')
&=
\frac{Y_u(w')}{p_k(w')}
\int_\Gamma X_u(w,s)
\frac{p_k(q^s w) - p_k(w')}{q^s w - w'} \diff s \\
&=
\frac{Y_u(w')}{p_k(w')}
\int_\Gamma X_u(w,s) B_{p_k}(q^s w, w') \diff s.
\end{align*}
By \eqref{eq:PQ-C-Bezoutian},
$B_{p_k}(x, w') = C_k[B_\Pi(\cdot, w')](x)$.  The
composition $\Psi(Tu) C_k \Phi(S_k u) \in \End(H)$ has
kernel
\begin{align*}
(\Psi(Tu) C_k \Phi(S_k u))(w, w')
&= (\Psi(Tu)[C_k \Phi(S_k u; \cdot, w')])(w) \\
&= \frac{Y_u(w')}{p_k(w')}
(\Psi(Tu)[B_{p_k}(\cdot, w')])(w) \\
&= \frac{Y_u(w')}{p_k(w')}
\int_\Gamma X_u(w,s) B_{p_k}(q^s w, w') \diff s,
\end{align*}
where
\[
C_k[\Phi(S_k u; \cdot, w')](x)
= C_k[B_\Pi(\cdot, w')](x) \frac{Y_u(w')}{p_k(w')}
= B_{p_k}(x, w') \frac{Y_u(w')}{p_k(w')},
\]
the scalar factors from $\Psi(Tu)$, and the integral is
\eqref{eq:HSEP-Psi-Tu}.  This matches the expression
above, confirming \eqref{eq:HSEP-K-Sk}.
\end{proof}

\subsection{One-point equations}
\label{sec:HSEP-equations}

Define the dressed observable
$\mathcal{M}(u) \in \End(E)$ by
\begin{equation}\label{eq:HSEP-M-def}
\mathcal{M}(u)
\defeq I - \Phi(u) (I + K(u))^{-1} \Psi(u),
\end{equation}
corresponding to $z = -1$ in the general dressed
observable $I + z\Phi(u)(I - zK(u))^{-1}\Psi(u)$.
The dressed edge weights are
\begin{equation}\label{eq:HSEP-dressed-weights}
\mathcal{M}_k(u) \defeq \mathcal{M}(Tu)^{-1}
C_{p_k} \mathcal{M}(S_k u),
\qquad k \in \{1, 2\}.
\end{equation}
Write $\mathcal{M}_{n,t,m} \defeq \mathcal{M}(n,t,m)$.

\begin{theorem}\label{thm:HSEP-diamond}
Suppose that $K_u$ is trace class
and that $(I + K_u)^{-1}$ exists at all relevant lattice
vertices.  Then
\begin{align}\label{eq:HSEP-mixed-diamond}
&\mathcal{M}_{n+1,t+1,m+1}^{-1}
\bigl(C_{p_1}\mathcal{M}_{n+1,t,m}\Lambda_{p_2}
- C_{p_2}\mathcal{M}_{n,t+1,m}\Lambda_{p_1}\bigr)
\nonumber \\
&{}+
\bigl(\Lambda_{p_1}\mathcal{M}_{n+2,t+1,m+1}^{-1}
C_{p_2}
- \Lambda_{p_2}\mathcal{M}_{n+1,t+2,m+1}^{-1}
C_{p_1}\bigr)\mathcal{M}_{n+1,t+1,m}
= 0.
\end{align}
\end{theorem}

\begin{proof}
The bare diamond equations are established in
Lemma~\ref{thm:PQ-bare-diamond}, the seed linear
problems in Lemma~\ref{lem:HSEP-seed-linear}, and the
dressing compatibility in
Lemma~\ref{lem:HSEP-dressing}.  Together with the
resolvent hypothesis, all conditions of
Theorem~\ref{thm:Darboux} are satisfied.  The mixed
diamond equation \eqref{eq:diamond-mixed-ij} with
$(i,j) = (1,2)$ is \eqref{eq:HSEP-mixed-diamond}.
\end{proof}

The Woodbury identity (Corollary~\ref{cor:Woodbury})
gives
\begin{equation}\label{eq:HSEP-Woodbury}
F(u) = F(Tu) \det_E\bigl(\mathcal{M}(Tu)\bigr),
\end{equation}
so $\det_E \mathcal{M}(Tu) = F(u)/F(Tu)$ expresses the
scalar content of the matrix equation.

\medskip

In the monomial basis $\{1, z, z^2, \dotsc, z^{I+J}\}$,
the operators $C_{p_k}$ and $\Lambda_{p_k}$ are explicit
$(I{+}J{+}1) \times (I{+}J{+}1)$ matrices determined
by polynomial long division with respect to $\Pi$.

\begin{example}\label{ex:HSEP-IJ2}
Take $I = J = 2$.  The shift polynomials are
\[
p_1(z) = (1-z)(1-qz), \qquad
p_2(z) = (1+\alpha z)(1+\alpha qz),
\]
so $\Pi(z) = z p_1(z) p_2(z)$ has degree $5$ and
$E \cong \C^5$ with monomial basis
$\{1, z, z^2, z^3, z^4\}$.  In this basis:
\[
C_{p_1} =
\begin{pmatrix}
0 & 0 & 0 & \frac{1}{\alpha^2 q}
& -\frac{1+q}{\alpha^3 q^2} \\[3pt]
0 & 0 & 0 & 0 & \frac{1}{\alpha^2 q} \\
0 & 0 & 0 & 0 & 0 \\
0 & 0 & 0 & 0 & 0 \\
0 & 0 & 0 & 0 & 0
\end{pmatrix},
\qquad
C_{p_2} =
\begin{pmatrix}
0 & 0 & 0 & \frac{1}{q}
& \frac{1+q}{q^2} \\[3pt]
0 & 0 & 0 & 0 & \frac{1}{q} \\
0 & 0 & 0 & 0 & 0 \\
0 & 0 & 0 & 0 & 0 \\
0 & 0 & 0 & 0 & 0
\end{pmatrix},
\]
\[
\Lambda_{p_1} =
\begin{pmatrix}
-1 & 0 & 0 & 0 & 0 \\[3pt]
1{+}q & -1 & 0 & \frac{1}{\alpha^2 q}
& -\frac{1+q}{\alpha^3 q^2} \\[3pt]
-q & 1{+}q & -1
& \frac{(\alpha-1)(1+q)}{\alpha^2 q}
& \frac{(1+q)^2-\alpha(1+q+q^2)}{\alpha^3 q^2} \\[3pt]
0 & -q & 1{+}q
& \frac{q-\alpha(1+q)^2}{\alpha^2 q}
& \frac{(1+q)(\alpha(1+q+q^2)-q)}{\alpha^3 q^2} \\[3pt]
0 & 0 & -q & \frac{1+q}{\alpha}
& -\frac{1+q+q^2}{\alpha^2 q}
\end{pmatrix},
\]
\[
\Lambda_{p_2} =
\begin{pmatrix}
-1 & 0 & 0 & 0 & 0 \\[3pt]
-\alpha(1{+}q) & -1 & 0 & \frac{1}{q}
& \frac{1+q}{q^2} \\[3pt]
-\alpha^2 q & -\alpha(1{+}q) & -1
& \frac{(\alpha-1)(1+q)}{q}
& \frac{\alpha(1+q)^2-(1+q+q^2)}{q^2} \\[3pt]
0 & -\alpha^2 q & -\alpha(1{+}q)
& \frac{\alpha(\alpha q-(1+q)^2)}{q}
& \frac{\alpha(1+q)(\alpha q-(1+q+q^2))}{q^2} \\[3pt]
0 & 0 & -\alpha^2 q & -\alpha^2(1{+}q)
& -\frac{\alpha^2(1+q+q^2)}{q}
\end{pmatrix}.
\]
These matrices have
$\operatorname{rank} C_{p_k} = 2$ and
$\operatorname{rank} \Lambda_{p_k} = 3$.
The $C$- and $\Lambda$-commutators vanish, as required
by the bare diamond equations:
$C_{p_1} C_{p_2} = C_{p_2} C_{p_1}$ and
$\Lambda_{p_1} \Lambda_{p_2}
= \Lambda_{p_2} \Lambda_{p_1}$.
The mixed relation is the nontrivial one:
\[
[C_{p_1}, \Lambda_{p_2}] =
\begin{pmatrix}
0 & -1 & 0 & \frac{1+\alpha^2}{\alpha^2 q}
& \frac{(\alpha^3-1)(1+q)}{\alpha^3 q^2} \\[3pt]
0 & 0 & -1 & -\frac{(\alpha-1)(1+q)}{\alpha q}
& -\frac{(1+\alpha^2)(1+q+q^2)}{\alpha^2 q^2} \\[3pt]
0 & 0 & 0 & 1 & 0 \\
0 & 0 & 0 & 0 & 1 \\
0 & 0 & 0 & 0 & 0
\end{pmatrix},
\]
while $[\Lambda_{p_1}, C_{p_2}] = -[C_{p_1}, \Lambda_{p_2}]$.
The bare mixed diamond equation is thus the cancellation
of two rank-four commutator defects; the rank is uniform
in the parameters, since the submatrix of
$[C_{p_1}, \Lambda_{p_2}]$ formed by rows $1$--$4$ and
columns $2$--$5$ has determinant~$1$.
\end{example}

\section{Inhomogeneous Stochastic Six-Vertex Model}\label{sec:IS6V}

\paragraph{\textbf{System description.}}
The stochastic six-vertex model with inhomogeneous
spectral parameters admits a Fredholm determinant formula
for the $q$-Laplace transform of the height function.
Fix $0 < q < 1$ and two sequences of spectral parameters
$u_1, \dotsc, u_y$ (row parameters) and
$v_1, \dotsc, v_x$ (column parameters), with
$0 < u_i v_j < 1$ for all $i, j$.
Step initial data is imposed throughout.

The height function $h_{6V}(x{+}1, y)$ counts the number
of paths that have crossed the horizontal level $x{+}1$
after $y$ rows.  For
$\zeta \in \C \setminus \R_{\geq 0}$, define the
$q$-Laplace transform
\begin{equation}\label{eq:IS6V-L-def}
\mathcal{L}_{x,y}(\zeta)
\defeq
\E_{\mathrm{step}}
\left[
\frac{1}{(\zeta q^{h_{6V}(x+1,y)};q)_\infty}
\right],
\qquad
(a;q)_\infty
\defeq \prod_{k=0}^\infty (1 - aq^k).
\end{equation}

\medskip
\paragraph{\textbf{Fredholm determinant formula.}}
Bufetov--Mucciconi--Petrov \cite[Thm.~9.2]{BufetovMucciconiPetrov2019}
establish that for any
$\zeta \in \C \setminus \R_{\geq 0}$,
\begin{equation}\label{eq:IS6V-fredholm}
\mathcal{L}_{x,y}(\zeta)
=
\det\bigl(I + K_\zeta[x,y]\bigr)_{L^2(C)},
\end{equation}
where $C$ is a positively oriented contour encircling $0$
and $u_1, \dotsc, u_y$, and the Mellin--Barnes kernel is
\begin{equation}\label{eq:IS6V-kernel-original}
K_\zeta[x,y](w,w')
=
\frac{1}{2\mathrm{i}}
\int_{d+\mathrm{i}\R}
\frac{(-\zeta)^r}{\sin(\pi r)}
\frac{f_{x,y}(w)}{f_{x,y}(q^r w)}
\frac{\diff r}{q^r w - w'},
\end{equation}
with $d \in (0, 1)$ and the gauge function
\begin{equation}\label{eq:IS6V-f-def}
f_{x,y}(w)
\defeq
\prod_{i=1}^{y} (w - u_i)^{-1}
\prod_{j=1}^{x} (1 - v_j w).
\end{equation}
The contour $C$ is chosen so that for all
$r \in d + \mathrm{i}\R$, the contour $C$ contains
$q^r C$ but does not contain $v_j^{-1} q^{-r}$ for
any $j$.

\medskip
\paragraph{\textbf{References.}}
The Fredholm determinant formula
\eqref{eq:IS6V-fredholm} is due to 
Bufetov, Mucciconi, Petrov \cite[Thm.\ 9.2]{BufetovMucciconiPetrov2019}.

\begin{remark}[Fredholm determinant convention]
\label{rem:IS6V-fredholm-conv}
The Fredholm determinant
$\det(I + K_\zeta[x,y])$ is taken on $L^2(C)$ with
the contour measure
$\diff w / (2\pi\mathrm{i})$.  The combined factor
$(-\zeta)^r/\sin(\pi r)$ decays exponentially on the
vertical $r$-contour for each
$\zeta \in \C \setminus \R_{\geq 0}$, giving uniform
convergence of the Mellin--Barnes integral on
$C \times C$.  The kernel is therefore continuous and
bounded on the contour square, and the Fredholm series
converges absolutely by the Hadamard estimate.  Since
$(-\zeta)^r$ is analytic in $\zeta$ on
$\C \setminus \R_{\geq 0}$ (principal branch), the
same estimates give analyticity of
$\mathcal{L}_{x,y}(\zeta)$ in $\zeta$ on this domain.
When the determinant is nonzero, the resolvent kernel
exists by classical Fredholm theory.  The framework
results invoked below (gauge invariance, the Woodbury
determinant ratio, and the Darboux theorem) hold at
the kernel and series level: their proofs require only
pointwise kernel identities, the resolvent integral
equation, and the rank-one perturbation formula for
Fredholm series.
\end{remark}

\subsection{Kernel reformulation}\label{sec:IS6V-kernel-reform}

The Fredholm determinant \eqref{eq:IS6V-fredholm}
depends on the physical parameters $(x, y)$ and the
spectral variable $\zeta$.  An auxiliary integer
parameter $m \in \Z$, which discretizes the spectral
variable through $\zeta_m \defeq q^m \zeta$, promotes
$(x, y, m)$ to a lattice on which the Darboux theorem
applies.  For $u = (x, y, m)$,
define $K_u$ by replacing $\zeta$ in
\eqref{eq:IS6V-kernel-original} by $\zeta_m$, so that
\begin{equation}\label{eq:IS6V-F-def}
F(u) \defeq \det(I + K_u)_{L^2(C)}
=
\E_{\mathrm{step}}
\left[
\frac{1}{(\zeta_m q^{h_{6V}(x+1,y)};q)_\infty}
\right].
\end{equation}

\paragraph{\textbf{Gauge transformation.}}
The gauge transformation
$G_u(w) \defeq w^{-m} f_{x,y}(w)$ conjugates the
kernel to
\[
K_u(w,w')
= G_u(w)^{-1} K_{\zeta_m}[x,y](w,w') G_u(w').
\]
Since $G_u$ is bounded and nonzero on $C$, gauge
invariance gives $\det(I + K_u) = F(u)$.  Equivalently,
\begin{equation}\label{eq:IS6V-gauged-kernel}
K_u(w,w')
=
\frac{w^m f_{x,y}(w')}{(w')^m}
\frac{1}{2\mathrm{i}}
\int_{d+\mathrm{i}\R}
\frac{(-\zeta_m)^r}{\sin(\pi r)}
\frac{1}{f_{x,y}(q^r w)}
\frac{\diff r}{q^r w - w'}.
\end{equation}
The symbol $K_u$ refers to the gauged kernel henceforth.
The gauged kernel factors as
\begin{equation}\label{eq:IS6V-XY-factor}
K_u(w,w') = \int_\Gamma X_u(w,r)
\frac{Y_u(w')}{q^r w - w'} \diff r,
\end{equation}
where $\Gamma = d + \mathrm{i}\R$ and
\[
X_u(w,r)
\defeq
w^m
\frac{1}{2\mathrm{i}}
\frac{(-\zeta_m)^r}{\sin(\pi r)}
\frac{1}{f_{x,y}(q^r w)},
\qquad
Y_u(w') \defeq \frac{f_{x,y}(w')}{(w')^m}.
\]

The rational structure of $f_{x,y}$ determines the
shift polynomials.  Each row parameter $u_i$
contributes a pole of $f_{x,y}$ at $w = u_i$, and
each column parameter $v_j$ contributes a zero at
$w = v_j^{-1}$.  The three types of shift polynomial
are
\begin{align}
p_i^{(\mathrm{row})}(w) &\defeq (w - u_i),
\qquad i = 1, \dotsc, y,
\label{eq:IS6V-p-row} \\
p_j^{(\mathrm{col})}(w) &\defeq (1 - v_j w),
\qquad j = 1, \dotsc, x,
\label{eq:IS6V-p-col} \\
p_0(w) &\defeq w.
\label{eq:IS6V-p-spec}
\end{align}
The fusion polynomial collects all $x + y + 1$ linear
factors:
\begin{equation}\label{eq:IS6V-Pi}
\Pi_{x,y}(w)
\defeq
w \prod_{i=1}^{y}(w - u_i)
\prod_{j=1}^{x}(1 - v_j w),
\qquad
\deg \Pi_{x,y} = x + y + 1.
\end{equation}

\paragraph{\textbf{The multiplicity lattice.}}
The $x + y + 1$ linear factors of $\Pi_{x,y}$ each
appear independently in the gauge function
$G_u(w) = w^{-m} f_{x,y}(w)$: the spectral factor
$w^{-m}$, the $y$ row poles $(w - u_i)^{-1}$, and
the $x$ column zeros $(1 - v_j w)$.  Define a lattice
$\mathcal{V} \cong \Z^{x+y+1}$ with coordinates
$\mathbf{a} = (a_0, a_1, \dotsc, a_{x+y})$ recording
how many times each elementary shift has been applied
from the base point.  At a general lattice point the
gauge function is
\begin{equation}\label{eq:IS6V-gauge-lattice}
G_{\mathbf{a}}(w)
= w^{-(m+a_0)}
\prod_{i=1}^{y} (w - u_i)^{-(1+a_i)}
\prod_{j=1}^{x} (1 - v_j w)^{1-a_{y+j}},
\end{equation}
with the base model at the origin
$\mathbf{a} = \mathbf{0}$.  Each lattice point defines
a gauged kernel and hence a Fredholm determinant; the
diamond equations established below are recurrences
connecting these Fredholm determinants.

The coordinate directions correspond to moves in the
inhomogeneous stochastic six-vertex model.  A positive
row step $a_i \mapsto a_i + 1$ adds one row of vertices
with spectral parameter $u_i$.  A positive column step
$a_{y+j} \mapsto a_{y+j} + 1$ lowers the exponent of
$(1 - v_j w)$ in the gauge function by one; at the
base point this removes column $j$, corresponding to
the specialization $v_j = 0$.  The spectral coordinate
$a_0$ shifts the $q$-Laplace variable
$\zeta_m \mapsto q^{a_0} \zeta_m$.

The elementary shifts $S_{p_k}$ each increment a
single coordinate:
\begin{equation}\label{eq:IS6V-lattice-shifts}
S_{p_k}(\mathbf{a})
= \mathbf{a} + \mathbf{e}_k,
\qquad k = 0, 1, \dotsc, x + y,
\end{equation}
where $\mathbf{e}_k$ is the $k$-th standard basis
vector of $\Z^{x+y+1}$.  For a divisor
$p = \prod_{k \in I} p_k$ of $\Pi_{x,y}$, with
$I \subseteq \{0, 1, \dotsc, x+y\}$, the
corresponding shift increments the coordinates
indexed by $I$:\footnote{When row or column parameters
coincide, distinct subsets $I$ may produce the same
polynomial $\prod_{k \in I} p_k$.  The shift
rule~\eqref{eq:IS6V-shift-rule} and the Euclidean
operators $C_p, \Lambda_p$ depend only on the
polynomial $p$, so the identities below hold for any
such factorization.}
\begin{equation}\label{eq:IS6V-general-shift}
S_p(\mathbf{a})
= \mathbf{a} + \sum_{k \in I} \mathbf{e}_k.
\end{equation}
In particular, the distinguished shift
$T = S_{\Pi_{x,y}}$ (with
$I = \{0, 1, \dotsc, x+y\}$) increments all
coordinates:
\begin{equation}\label{eq:IS6V-T-def}
T(\mathbf{a}) = \mathbf{a} + \mathbf{1}.
\end{equation}
The $X$--$Y$ factors of \eqref{eq:IS6V-XY-factor}
transform under $S_p$ as
\begin{equation}\label{eq:IS6V-shift-rule}
X_{S_p u}(w,r) = X_u(w,r) p(q^r w),
\qquad
Y_{S_p u}(w') = \frac{Y_u(w')}{p(w')}.
\end{equation}

\subsection{Base graph and seed data}
\label{sec:IS6V-base-graph}

Let $\mathcal{V} \cong \Z^{x+y+1}$ be the lattice
generated by the elementary shifts
\eqref{eq:IS6V-lattice-shifts}.  The quotient algebra
$E \defeq E_{x,y} = \C[w]/(\Pi_{x,y}(w))$ of
\S\ref{sec:PQ} has dimension $x + y + 1$.  The Euclidean
division edge weights $C_p, \Lambda_p \in \End(E)$ are
defined by \eqref{eq:PQ-division} with
$\Pi = \Pi_{x,y}$; these satisfy the bare diamond
equations by Lemma~\ref{thm:PQ-bare-diamond}.

Let $H = L^2(C)$.  Define the wave function
$\Psi(u) \colon E \to H$ by
\begin{equation}\label{eq:IS6V-Psi-def}
(\Psi(u) f)(w)
\defeq
\int_\Gamma X_u(w,r)
\frac{f(q^r w)}{\Pi_{x,y}(q^r w)} \diff r.
\end{equation}
At the vertex $Tu$, the shift rule
\eqref{eq:IS6V-shift-rule} with $p = \Pi_{x,y}$ cancels
the denominator:
\begin{equation}\label{eq:IS6V-Psi-Tu}
(\Psi(Tu) f)(w)
= \int_\Gamma X_u(w,r) f(q^r w) \diff r.
\end{equation}
Define the adjoint wave function
$\Phi(u) \colon H \to E$ by
\begin{equation}\label{eq:IS6V-Phi-def}
(\Phi(u) h)(w)
\defeq
\int_C B_{\Pi_{x,y}}(w, w') Y_u(w') h(w')
\frac{\diff w'}{2\pi\mathrm{i}}.
\end{equation}
Write 
\[\Phi(u; w, w') \defeq B_{\Pi_{x,y}}(w, w') Y_u(w')\]
for the integral kernel.  Since $B_{\Pi_{x,y}}(\cdot, w')$
has degree $x + y$ in the first variable,
$\Phi(u) h \in E$.  These are the seed functions of
\S\ref{sec:HSEP-base-graph} with $\Pi$ replaced by
$\Pi_{x,y}$, the integration contour $C_r$ by $C$, and
the Mellin--Barnes variable $s$ by $r$.

\begin{lemma}\label{lem:IS6V-seed-linear}
The seed functions $\Psi$ and $\Phi$ satisfy the linear
problem on $\mathcal{V}$ with edge weights
$(C_p, \Lambda_p)$.  For any divisor $p \mid \Pi_{x,y}$
with $1 \leq \deg p < \deg \Pi_{x,y}$:
\begin{equation}\label{eq:IS6V-Psi-recurrence}
\Psi(S_p u) = \Psi(Tu) C_p - \Psi(u) \Lambda_p,
\end{equation}
and
\begin{equation}\label{eq:IS6V-Phi-recurrence}
\Phi(Tu) = C_p \Phi(S_p u) - \Lambda_p \Phi(T S_p u).
\end{equation}
\end{lemma}

\begin{proof}
The argument parallels the proof of
Lemma~\ref{lem:HSEP-seed-linear} with $\Pi$ replaced by
$\Pi_{x,y}$ and $p_k$ by an arbitrary divisor $p$.
The Mellin--Barnes integrals below converge absolutely
and uniformly for $w \in C$, so all manipulations may
be performed under the integral sign.

\medskip
\noindent\textit{Verification of
\eqref{eq:IS6V-Psi-recurrence}.}
The shift rule gives
$X_{S_p u}(w,r) = X_u(w,r) p(q^r w)$, so
\[
(\Psi(S_p u) f)(w)
=
\int_\Gamma X_u(w,r)
\frac{p(q^r w) f(q^r w)}{\Pi_{x,y}(q^r w)} \diff r.
\]
The Euclidean division identity \eqref{eq:PQ-division}
gives $pf/\Pi_{x,y} = C_p f - \Lambda_p f / \Pi_{x,y}$,
separating the integral into
$(\Psi(Tu) C_p f)(w) - (\Psi(u) \Lambda_p f)(w)$
by \eqref{eq:IS6V-Psi-Tu} and \eqref{eq:IS6V-Psi-def}.

\medskip
\noindent\textit{Verification of
\eqref{eq:IS6V-Phi-recurrence}.}
The recurrence reduces to an identity of integral
kernels, since the shift rule acts only on the
$Y_u(w')$ factor.  The kernels at $Tu$, $S_p u$, and
$TS_p u$ are
\begin{gather*}
\Phi(Tu; w, w') = B_{\Pi_{x,y}}(w,w')
\frac{Y_u(w')}{\Pi_{x,y}(w')},
\quad
\Phi(S_p u; w, w') = B_{\Pi_{x,y}}(w,w')
\frac{Y_u(w')}{p(w')}, \\
\Phi(TS_p u; w, w') = B_{\Pi_{x,y}}(w,w')
\frac{Y_u(w')}{\Pi_{x,y}(w') p(w')}.
\end{gather*}
Since $Y_u(w')/p(w')$ does not depend on $w$, it
factors from $C_p$ and $\Lambda_p$.  Substituting
$C_p[B_{\Pi_{x,y}}(\cdot, w')] = B_p(\cdot, w')$ from
\eqref{eq:PQ-C-Bezoutian} and
$\Lambda_p[B_{\Pi_{x,y}}(\cdot, w')]
= \Pi_{x,y}(w') B_p(\cdot, w')
- p(w') B_{\Pi_{x,y}}(\cdot, w')$
from \eqref{eq:PQ-Lambda-Bezoutian}:
\begin{align*}
C_p[\Phi(S_p u; \cdot, w')](w)
&- \Lambda_p[\Phi(TS_p u; \cdot, w')](w) \\
&=
\frac{Y_u(w')}{p(w')} B_p(w,w')
-
\frac{Y_u(w')}{p(w')} B_p(w,w') \\
&\quad +
\frac{Y_u(w')}{\Pi_{x,y}(w')} B_{\Pi_{x,y}}(w,w') \\
&= \Phi(Tu; w, w'). \qedhere
\end{align*}
\end{proof}

\subsection{Dressing compatibility}
\label{sec:IS6V-dressing}

\begin{lemma}\label{lem:IS6V-dressing}
The gauged kernel $K_u$ satisfies the dressing
compatibility conditions: for any divisor
$p \mid \Pi_{x,y}$ with
$1 \leq \deg p < \deg \Pi_{x,y}$,
\begin{align}
K(Tu)(w,w') - K(u)(w,w')
&= (\Psi(Tu) \Phi(Tu))(w,w'),
\label{eq:IS6V-K-T} \\
K(S_p u)(w,w') - K(u)(w,w')
&= (\Psi(Tu) C_p \Phi(S_p u))(w,w').
\label{eq:IS6V-K-Sp}
\end{align}
\end{lemma}

\begin{proof}
The argument parallels the proof of
Lemma~\ref{lem:HSEP-dressing}, with $\Pi$ replaced by
$\Pi_{x,y}$ and $p_k$ by an arbitrary divisor $p$.

\medskip
\noindent\textit{Verification of \eqref{eq:IS6V-K-T}.}
The shift rule gives
$X_{Tu} = X_u \Pi_{x,y}(q^r w)$ and
$Y_{Tu} = Y_u / \Pi_{x,y}(w')$, so
\begin{align*}
K(Tu)(w,w') - K(u)(w,w')
&=
\frac{Y_u(w')}{\Pi_{x,y}(w')}
\int_\Gamma X_u(w,r)
\frac{\Pi_{x,y}(q^r w) - \Pi_{x,y}(w')}{q^r w - w'}
\diff r \\
&=
\frac{Y_u(w')}{\Pi_{x,y}(w')}
\int_\Gamma X_u(w,r) B_{\Pi_{x,y}}(q^r w, w') \diff r.
\end{align*}
The composition $\Psi(Tu)\Phi(Tu) \in \End(H)$ has kernel
\begin{align*}
(\Psi(Tu)\Phi(Tu))(w, w')
&= (\Psi(Tu)[\Phi(Tu; \cdot, w')])(w) \\
&= \frac{Y_u(w')}{\Pi_{x,y}(w')}
\int_\Gamma X_u(w,r) B_{\Pi_{x,y}}(q^r w, w') \diff r,
\end{align*}
since the scalar $Y_u(w')/\Pi_{x,y}(w')$ factors from
$\Psi(Tu)$ and the integral is \eqref{eq:IS6V-Psi-Tu}
applied to $B_{\Pi_{x,y}}(\cdot, w')$.

\medskip
\noindent\textit{Verification of \eqref{eq:IS6V-K-Sp}.}
The shift rule gives
$X_{S_p u} = X_u p(q^r w)$ and
$Y_{S_p u} = Y_u / p(w')$, so
\begin{align*}
K(S_p u)(w,w') - K(u)(w,w')
&=
\frac{Y_u(w')}{p(w')}
\int_\Gamma X_u(w,r) B_p(q^r w, w') \diff r.
\end{align*}
By \eqref{eq:PQ-C-Bezoutian},
$B_p(w, w') = C_p[B_{\Pi_{x,y}}(\cdot, w')](w)$.
The composition
$\Psi(Tu) C_p \Phi(S_p u) \in \End(H)$ has kernel
\begin{align*}
(\Psi(Tu) C_p \Phi(S_p u))(w, w')
&= \frac{Y_u(w')}{p(w')}
(\Psi(Tu)[B_p(\cdot, w')])(w) \\
&= \frac{Y_u(w')}{p(w')}
\int_\Gamma X_u(w,r) B_p(q^r w, w') \diff r,
\end{align*}
confirming \eqref{eq:IS6V-K-Sp}.
\end{proof}

\subsection{One-point equations}
\label{sec:IS6V-equations}

Define the dressed observable
$\mathcal{M}(u) \in \End(E_{x,y})$ by
\begin{equation}\label{eq:IS6V-M-def}
\mathcal{M}(u)
\defeq I - \Phi(u) (I + K(u))^{-1} \Psi(u),
\end{equation}
corresponding to $z = -1$ in the resolvent
$(I - zK)^{-1}$.
The dressed edge weights are
\begin{equation}\label{eq:IS6V-dressed-weights}
\mathcal{M}_p(u) \defeq \mathcal{M}(Tu)^{-1}
C_p \mathcal{M}(S_p u).
\end{equation}

\begin{theorem}\label{thm:IS6V-diamond}
Suppose that $K_u$ is trace class and that
$(I + K_u)^{-1}$ exists at all relevant lattice
vertices.  Then for any two distinct divisors
$p, p' \mid \Pi_{x,y}$ with
$1 \leq \deg p, \deg p' < \deg \Pi_{x,y}$,
the dressed edge weights
satisfy the mixed diamond equation
\begin{equation}\label{eq:IS6V-mixed-diamond}
\mathcal{M}_p(u) \Lambda_{p'}
+
\Lambda_p \mathcal{M}_{p'}(S_p u)
-
\mathcal{M}_{p'}(u) \Lambda_p
-
\Lambda_{p'} \mathcal{M}_p(S_{p'} u)
= 0.
\end{equation}
\end{theorem}

\begin{proof}
The bare diamond equations are established in
Lemma~\ref{thm:PQ-bare-diamond}, the seed linear
problem in Lemma~\ref{lem:IS6V-seed-linear}, and the
dressing compatibility in
Lemma~\ref{lem:IS6V-dressing}.
Together with the resolvent hypothesis, all conditions of
Theorem~\ref{thm:Darboux} are satisfied.  The mixed
diamond equation
\eqref{eq:diamond-mixed-ij} for the pair $(p, p')$
is \eqref{eq:IS6V-mixed-diamond}.
\end{proof}

The Woodbury identity (Corollary~\ref{cor:Woodbury})
gives
\begin{equation}\label{eq:IS6V-Woodbury}
F(u) = F(Tu)
\det_{E_{x,y}}\bigl(\mathcal{M}(Tu)\bigr),
\end{equation}
so the $(x{+}y{+}1) \times (x{+}y{+}1)$ determinant
$\det_{E_{x,y}} \mathcal{M}(Tu) = F(u)/F(Tu)$ expresses
the scalar content of the matrix equation.

\medskip

The $x + y + 1$ linear factors of $\Pi_{x,y}$ produce
$\binom{x+y+1}{2}$ mixed diamond equations, one per pair
of shift directions.  General divisors
$p \mid \Pi_{x,y}$ of degree $2$ or higher package
several elementary shifts into a single compound move
and yield additional diamond equations by
Theorem~\ref{thm:IS6V-diamond}.

The monomial basis $\{1, w, \dotsc, w^{x+y}\}$
identifies $E_{x,y}$ with $\C^{x+y+1}$.  In this
basis, $C_p$ and $\Lambda_p$ become explicit
$(x{+}y{+}1) \times (x{+}y{+}1)$ matrices determined
by polynomial long division with respect to $\Pi_{x,y}$.

\begin{example}\label{ex:IS6V-overlapping}
Take $x = 2$, $y = 1$.  The two degree-two divisors
\[
p(w) \defeq (w-u_1)(1-v_1 w), \qquad
p'(w) \defeq (w-u_1)(1-v_2 w)
\]
both divide $\Pi_{2,1}$ but are not complementary:
$pp' = (w-u_1)^2(1-v_1 w)(1-v_2 w) \neq \Pi_{2,1}$.
The two shifts overlap in the row factor $(w-u_1)$ and
differ by which column factor is included.
In the monomial basis $\{1, w, w^2, w^3\}$:
\[
C_p =
\begin{pmatrix}
0 & 0 & -v_2^{-1} & -v_2^{-2} \\
0 & 0 & 0 & -v_2^{-1} \\
0 & 0 & 0 & 0 \\
0 & 0 & 0 & 0
\end{pmatrix},
\qquad
C_{p'} =
\begin{pmatrix}
0 & 0 & -v_1^{-1} & -v_1^{-2} \\
0 & 0 & 0 & -v_1^{-1} \\
0 & 0 & 0 & 0 \\
0 & 0 & 0 & 0
\end{pmatrix},
\]
\[
\Lambda_p =
\begin{pmatrix}
u_1 & 0 & 0 & 0 \\
-1-v_1 u_1 & u_1 & u_1 v_2^{-1}
& u_1 v_2^{-2} \\[3pt]
v_1 & -1-v_1 u_1 & -(1+v_1 u_1)v_2^{-1}
& -(1+v_1 u_1)v_2^{-2} \\[3pt]
0 & v_1 & v_1 v_2^{-1} & v_1 v_2^{-2}
\end{pmatrix},
\]
\[
\Lambda_{p'} =
\begin{pmatrix}
u_1 & 0 & 0 & 0 \\
-1-v_2 u_1 & u_1 & u_1 v_1^{-1}
& u_1 v_1^{-2} \\[3pt]
v_2 & -1-v_2 u_1 & -(1+v_2 u_1)v_1^{-1}
& -(1+v_2 u_1)v_1^{-2} \\[3pt]
0 & v_2 & v_2 v_1^{-1} & v_2 v_1^{-2}
\end{pmatrix}.
\]
All four matrices have rank two for generic
$v_1 \neq v_2$.  The commutators $[C_p, C_{p'}]$ and
$[\Lambda_p, \Lambda_{p'}]$ vanish.  The mixed
commutator is
\[
[C_p, \Lambda_{p'}] =
\begin{pmatrix}
-1 & u_1 & \frac{(v_1+v_2)u_1}{v_1 v_2}
& (v_1^{-2}+v_2^{-2})u_1 \\[3pt]
0 & -1 & -\frac{v_1+v_2+v_1 v_2 u_1}{v_1 v_2}
& -v_1^{-2}-v_2^{-2} \\[3pt]
0 & 0 & 1 & -u_1 \\
0 & 0 & 0 & 1
\end{pmatrix},
\]
with $[\Lambda_p, C_{p'}] = -[C_p, \Lambda_{p'}]$ and
$\det[C_p, \Lambda_{p'}] = 1$.  The bare mixed diamond
equation is the cancellation of two full-rank commutator
defects, uniform in the parameters.
\end{example}

\section{Log-Gamma Polymer}
\label{sec:LG-onepoint-finiteE}

\paragraph{\textbf{System description.}}
Fix integers $N,M\geq 1$ and parameters
\[
a=(a_1,\ldots,a_N),\qquad
\alpha=(\alpha_1,\ldots,\alpha_M),
\qquad
\min_j\alpha_j-\max_i a_i>0.
\]
Let $\{d_{m,n}\}$ be independent inverse-gamma
random variables of shape $\alpha_m-a_n$.  The log-gamma polymer
partition function is
\[
Z^{M,N}
\defeq
\sum_{\pi\colon(1,1)\to(M,N)}
\prod_{(m,n)\in\pi} d_{m,n},
\]
where the sum is over up/right lattice paths from $(1,1)$ to $(M,N)$.

\paragraph{\textbf{Fredholm determinant formula.}}
For $M\geq N$ and $\mathrm{Re}(s)>0$, the
Borodin--Corwin--Remenik formula gives
\[
\E[e^{-sZ^{M,N}}]
=\det(I+K^{\mathrm{BCR}}_{s,N,M})_{L^2(C_{\delta_1})},
\]
where $C_{\delta_1}$ is the positively oriented circle of
radius~$\delta_1$ centred at the origin, with kernel
\begin{equation}\label{eq:LG-BCR-kernel}
K^{\mathrm{BCR}}_{s,N,M}(v,v')
=
\frac{1}{2\pi\mathrm{i}}\int_{\delta_2+\mathrm{i}\R}
\frac{\pi}{\sin\pi(v-w)}
\prod_{n=1}^{N}\frac{\Gamma(v-a_n)}{\Gamma(w-a_n)}
\prod_{m=1}^{M}
\frac{\Gamma(\alpha_m-w)}{\Gamma(\alpha_m-v)}
\frac{s^{w-v}}{w-v'}\diff w.
\end{equation}
The parameters $\delta_1,\delta_2$ satisfy
$0<\delta_2<1$ and
$0<\delta_1<\min\{\delta_2,1-\delta_2\}$, with
$\max_i a_i<\delta_1$ and $\min_m\alpha_m>\delta_2$.
The notation $s$ for the Laplace variable avoids collision
with the lattice vertex.

\paragraph{\textbf{References.}}
The log-gamma polymer was introduced by
Sepp\"al\"ainen \cite{SeppPolymerBoundary}.
The Fredholm determinant formula \eqref{eq:LG-BCR-kernel}
is due to Borodin, Corwin, and
Remenik~\cite[Corollary~1.8]{BorodinCorwinRemenik2013}.

\subsection{Kernel reformulation}
\label{sec:LG-kernel-reformulation}

The Fredholm determinant \eqref{eq:LG-BCR-kernel} depends on
the parameters $(a, \alpha)$ through the Gamma products in
the integrand.  Integer shifts of these parameters produce
polynomial factors in the gauge, promoting the parameter
space to a lattice on which the Darboux theorem applies.
Write $u = (a, \alpha)$ for the lattice vertex.

\paragraph{\textbf{Gauge transformation.}}
The gauge transformation
\[
G_u(z) \defeq s^{-z}
\prod_{n=1}^{N}\Gamma(z-a_n)
\prod_{m=1}^{M}\Gamma(\alpha_m-z)^{-1}
\]
defines the gauged kernel
\[
K_u(v,y) \defeq G_u(v)^{-1}
K^{\mathrm{BCR}}_{s,N,M}(v,y) G_u(y).
\]
The gauged kernel factors as
\begin{equation}\label{eq:LG-gauged-kernel}
K_u(v,y)
=
\int_\ell X_u(v,w)\frac{Y_u(y)}{w-y}\diff w,
\end{equation}
where
\begin{equation}\label{eq:LG-XY}
X_u(v,w)\defeq
\frac{\pi}{2\pi\mathrm{i}\sin\pi(v-w)}\frac{1}{G_u(w)},
\qquad
Y_u(y)\defeq G_u(y).
\end{equation}
The integral in \eqref{eq:LG-gauged-kernel} and all
$w$-integrals that follow are taken over a contour
$\ell$; the admissibility hypotheses on $\ell$,
principally trace class of the gauged kernel and
nonvanishing of the polynomial denominators, are
specified in \S\ref{sec:LG-quotient-data}, after
$\Psi$ and $\Phi$ are defined.  For the BCR vertical
line $\ell = \delta_2 + \mathrm{i}\R$, the gauge $G_u$
is holomorphic and nonvanishing in a neighbourhood of
$C_{\delta_1}$, and gauge invariance gives
\[
F_u(s) \defeq \det(I+K_u(s))_{L^2(C_{\delta_1})}
= \E[e^{-sZ^{M,N}}];
\]
see Remark~\ref{rem:LG-contour-regime}.

\paragraph{\textbf{Shift polynomials.}}
Integer shifts of the parameters $(a, \alpha)$
transform the gauge by division by affine factors.
By the Gamma recurrence $\Gamma(x+1) = x\Gamma(x)$,
a unit shift $a_i \mapsto a_i + 1$ divides $G_u$ by
$(z - a_i - 1)$, and a unit shift
$\alpha_j \mapsto \alpha_j + 1$ divides $G_u$ by
$(\alpha_j - z)$.  For nonnegative integer vectors
$I = (I_1, \ldots, I_N)$ and
$J = (J_1, \ldots, J_M)$, the shift
$a_i \mapsto a_i + I_i$,
$\alpha_j \mapsto \alpha_j + J_j$ has shift polynomial
\begin{equation}\label{eq:LG-shift-polynomial}
p_{I,J}(z) \defeq
\prod_{i=1}^{N}\prod_{q=1}^{I_i}(z-a_i-q)
\prod_{j=1}^{M}\prod_{q=0}^{J_j-1}(\alpha_j+q-z).
\end{equation}
Write $S_p$ for the parameter shift whose polynomial
is $p$.  A distinguished polynomial $\Pi$ of degree
$D \geq 1$, together with polynomials
$p_2, \ldots, p_r$ of degree at most $D$, all products
of affine factors from
\eqref{eq:LG-shift-polynomial}, determine shift
directions $T \defeq S_\Pi$ and
$S_j \defeq S_{p_j}$, $j \geq 2$, acting on disjoint
sets of parameters.  These shifts generate a lattice
$\mathcal{V}' \cong \Z^r$; each vertex determines a
log-gamma model with shifted parameters.  Since each
direction shifts parameters untouched by the others,
the $X$--$Y$ ratios
\begin{equation}\label{eq:LG-shift-ratios}
\frac{X_{S_p u}}{X_u}=p(w),
\qquad
\frac{Y_{S_p u}}{Y_u}=p(y)^{-1}
\end{equation}
hold at every vertex.

The Laplace variable $s > 0$ provides a continuous
direction: since $G_u(z)$ contains $s^{-z}$,
\[
s\partial_s X_u(v,w) = w X_u(v,w),
\qquad
s\partial_s Y_u(y) = -y Y_u(y).
\]
Setting $\xi \defeq \log s$, the full lattice
$\mathcal{V} \defeq \R \times \mathcal{V}'$ with
$\partial_1 \defeq \partial_\xi = s\partial_s$
instantiates the semi-discrete framework of
Proposition~\ref{prop:sc-diamond}.  The diamond
equations are recurrences connecting the Fredholm
determinants of these shifted models.

\subsection{Base graph and seed data}
\label{sec:LG-quotient-data}

The quotient algebra $E = \C[z]/(\Pi)$ of \S\ref{sec:PQ} has
dimension $D = \deg\Pi$.  The Euclidean division operators
$C_p, \Lambda_p \in \End(E)$ are defined by
\begin{equation}\label{eq:LG-division}
p(z)f(z) = \Pi(z)(C_pf)(z) - (\Lambda_pf)(z),
\qquad f\in E.
\end{equation}
Write $C_1 \defeq C_z$ and $\Lambda_1 \defeq \Lambda_z$,
and $(C_j, \Lambda_j) \defeq (C_{p_j}, \Lambda_{p_j})$ for
$j \geq 2$.  These edge weights satisfy the bare diamond
equations by Lemma~\ref{thm:PQ-bare-diamond}; since they are
independent of $\xi$, they also satisfy the semi-discrete
diamond equations of Proposition~\ref{prop:sc-diamond}.

Let $H = L^2(C_{\delta_1})$.  The wave function
$\Psi(u) \colon E \to H$ is defined by
\begin{equation}\label{eq:LG-Psi-def}
(\Psi(u)f)(v)
\defeq
\int_\ell X_u(v,w)\frac{f(w)}{\Pi(w)}\diff w,
\qquad f\in E.
\end{equation}
At the vertex $Tu$, the shift rule
\eqref{eq:LG-shift-ratios} with $p = \Pi$ cancels the
denominator in \eqref{eq:LG-Psi-def}:
\begin{equation}\label{eq:LG-Psi-Tu}
(\Psi(Tu)f)(v)
= \int_\ell X_u(v,w) f(w) \diff w.
\end{equation}
The adjoint wave function
$\Phi(u) \colon H \to E$ is defined by
\begin{equation}\label{eq:LG-Phi-def}
(\Phi(u)h)(z)
\defeq
\int_{C_{\delta_1}} B_{\Pi}(z,y)Y_u(y)h(y)\diff y.
\end{equation}
Since $B_{\Pi}(\cdot,y)$ has degree $D-1$ in $z$,
the image $\Phi(u)h$ lies in $E$.  The contour $\ell$
and the radius $\delta_1$ are chosen so that $\Pi$ and
each $p_j$ are nonvanishing on $C_{\delta_1}$, that
$\Pi$ is nonvanishing on $\ell$, and that $K(v)$ is
trace class on $H$ at every vertex $v$ of the required
stencil.  Recall that
$\partial_1 = \partial_\xi = s\partial_s$ where
$\xi = \log s$.

\begin{lemma}\label{lem:LG-seed-linear}
The seed functions $\Psi$ \eqref{eq:LG-Psi-def} and $\Phi$
\eqref{eq:LG-Phi-def} satisfy the semi-discrete linear
problem \eqref{eq:sc-linear-1}--\eqref{eq:sc-linear-j} and
adjoint problem
\eqref{eq:sc-adjoint-1}--\eqref{eq:sc-adjoint-j} with edge
weights $(C_k, \Lambda_k)$.  Explicitly:
\[
s\partial_s\Psi(u)
=\Psi(Tu)C_1-\Psi(u)\Lambda_1,
\qquad
\Psi(S_ju)
=\Psi(Tu)C_j-\Psi(u)\Lambda_j,
\]
and
\[
s\partial_s\Phi(Tu)
=\Lambda_1\Phi(Tu)-C_1\Phi(u),
\qquad
\Phi(Tu)
=C_j\Phi(S_ju)-\Lambda_j\Phi(TS_ju).
\]
\end{lemma}

\begin{proof}
The discrete recurrences follow by the same argument as
Lemma~\ref{lem:HSEP-seed-linear}: the shift rule
\eqref{eq:LG-shift-ratios} applied to \eqref{eq:LG-Psi-def}
multiplies the numerator by $p_j(w)$, the Euclidean division
\eqref{eq:LG-division} separates the integral into
$\Psi(Tu)C_jf - \Psi(u)\Lambda_jf$ by
\eqref{eq:LG-Psi-Tu} and \eqref{eq:LG-Psi-def}, and the
discrete adjoint follows from the same kernel
calculation as \eqref{eq:HSEP-Phi-recurrence},
substituting \eqref{eq:PQ-C-Bezoutian} and
\eqref{eq:PQ-Lambda-Bezoutian} into the integral kernel of
$\Phi$.

For the continuous recurrences, $s\partial_sX_u(v,w) = wX_u(v,w)$
gives the $\partial_1$-equation for $\Psi$ by the same
division with $p_1(w) \defeq w$.  For the continuous adjoint,
$B_z(z,y) = 1$ gives $C_1B_{\Pi}(\cdot,y)=1$ and
$\Lambda_1B_{\Pi}(\cdot,y) = \Pi(y)-yB_{\Pi}(\cdot,y)$, so
the kernel of $\Lambda_1\Phi(Tu)-C_1\Phi(u)$ is
\[
\bigl[
\frac{\Pi(y)-yB_{\Pi}(z,y)}{\Pi(y)}
-1
\bigr]Y_u(y)
=
-y\frac{B_{\Pi}(z,y)}{\Pi(y)}Y_u(y),
\]
which is the kernel of $s\partial_s\Phi(Tu)$ since
$s\partial_sY_u(y)=-yY_u(y)$.
\end{proof}

\begin{lemma}\label{lem:LG-dressing}
The gauged kernel $K_u$ is dressing compatible with
$(\Psi, \Phi)$ in the sense of
Definition~\ref{def:sc-dressing}:
\begin{align}
K(Tu)-K(u)&=\Psi(Tu)\Phi(Tu),
\label{eq:LG-KT}\\
s\partial_sK(u)&=\Psi(Tu)C_1\Phi(u),
\label{eq:LG-Kd}\\
K(S_ju)-K(u)&=\Psi(Tu)C_j\Phi(S_ju).
\label{eq:LG-KS}
\end{align}
\end{lemma}

\begin{proof}
\noindent\textit{Verification of \eqref{eq:LG-KT}.}
The shift rule \eqref{eq:LG-shift-ratios} gives
\begin{align*}
K(Tu)(v,y)-K(u)(v,y)
&=
\int X_u(v,w)Y_u(y)
\frac{\Pi(w)\Pi(y)^{-1}-1}{w-y}\diff w \\
&=
\int X_u(v,w)\frac{Y_u(y)}{\Pi(y)}
B_{\Pi}(w,y)\diff w.
\end{align*}
The right-hand side equals $\Psi(Tu)\Phi(Tu)$ by
\eqref{eq:LG-Psi-Tu} and \eqref{eq:LG-Phi-def}.

\medskip
\noindent\textit{Verification of \eqref{eq:LG-KS}.}
The same argument with $p_j$ replacing $\Pi$ gives
\begin{align*}
K(S_ju)(v,y)-K(u)(v,y)
&=
\int X_u(v,w)\frac{Y_u(y)}{p_j(y)}
B_{p_j}(w,y)\diff w.
\end{align*}
By \eqref{eq:PQ-C-Bezoutian}, this is
$\Psi(Tu)C_j\Phi(S_ju)$.

\medskip
\noindent\textit{Verification of \eqref{eq:LG-Kd}.}
Since $M \geq N$, the integrand of
\eqref{eq:LG-gauged-kernel} decays exponentially, so
the integral over $\ell$ converges absolutely and
differentiation in $\xi$ may be performed under the
integral sign.
Differentiating \eqref{eq:LG-gauged-kernel} in $\xi$:
\[
s\partial_s[X_u(v,w)Y_u(y)]
=(w-y)X_u(v,w)Y_u(y).
\]
The factor $(w-y)$ cancels the denominator in
\eqref{eq:LG-gauged-kernel}, so
\[
s\partial_sK(u)(v,y)
=\int X_u(v,w)Y_u(y)\diff w.
\]
Since $C_1B_{\Pi}(\cdot,y)=1$, the last
display equals $\Psi(Tu)C_1\Phi(u)$.
\end{proof}

\subsection{One-point equations}
\label{sec:LG-matrix-equation}

Let $R(u)\defeq(I+K(u))^{-1}$ and define the dressed
observable\footnote{The convention $\det(I+K)$
corresponds to Darboux spectral parameter $-1$ in
Theorem~\ref{thm:sc-Darboux}.}
\begin{equation}\label{eq:LG-M-def}
\mathcal{M}(u)
\defeq
I_E-\Phi(u)R(u)\Psi(u)
\end{equation}
and the dressed edge weights
\[
\mathcal{M}_1(u)
\defeq
\mathcal{M}(Tu)^{-1}C_1\mathcal{M}(u),
\qquad
\mathcal{M}_j(u)
\defeq
\mathcal{M}(Tu)^{-1}C_j\mathcal{M}(S_ju),
\quad j\geq2.
\]

\begin{theorem}\label{thm:LG-onepoint-matrix}
Suppose the resolvents $(I+K_u)^{-1}$ exist at every
vertex in the required stencil.  Then the dressed edge
weights $(\mathcal{M}_k, \Lambda_k)$ satisfy the
semi-discrete diamond equations of
Proposition~\ref{prop:sc-diamond}.  In particular,
\begin{align}
\mathcal{M}_1(u)\Lambda_j
+\Lambda_1\mathcal{M}_j(u)
-\mathcal{M}_j(u)\Lambda_1
-\Lambda_j\mathcal{M}_1(S_ju)
=s\partial_s\mathcal{M}_j(u).
\label{eq:LG-matrix-HM}
\end{align}
\end{theorem}

\begin{proof}
Lemma~\ref{thm:PQ-bare-diamond} gives the bare diamond
equations, Lemma~\ref{lem:LG-seed-linear} gives the seed
data, and Lemma~\ref{lem:LG-dressing} gives dressing
compatibility.  Theorem~\ref{thm:sc-Darboux}
(Remark~\ref{rem:sc-domain-restriction}) gives the
dressed diamond equations on the finite stencil for
$(\mathcal{M}_k, \Lambda_k)$.  Since $\Lambda_k$ is
constant, the mixed equation
\eqref{eq:sc-mixed-dressed} reduces to
\eqref{eq:LG-matrix-HM}.
\end{proof}

At each vertex $v$, the Fredholm determinant
$F_v(s) \defeq \det_H(I + K(v))$ is the log-gamma
Laplace transform for the parameters at $v$.
The Woodbury identity
(Corollary~\ref{cor:sc-Woodbury}) gives
\begin{equation}\label{eq:LG-Woodbury}
F_u(s)=F_{Tu}(s)\det_E\mathcal{M}(Tu).
\end{equation}
The matrix equation \eqref{eq:LG-matrix-HM} and the
scalar relation \eqref{eq:LG-Woodbury} are recurrences
connecting log-gamma Laplace transforms at different
parameter values.

\begin{remark}\label{rem:LG-contour-regime}
For stencils using only $\alpha$-shifts (the
$a$-parameters fixed), the vertical line
$\ell = \delta_2 + \mathrm{i}\R$ satisfies the standing
hypothesis at every vertex: upward $\alpha$-shifts
preserve $\alpha_m > \delta_2$, while the $a$-parameters
remain in $[0, \delta_1)$, so the BCR simple-contour
conditions hold and $F_v = \E[e^{-sZ_v}]$ at each
vertex.  Stencils including upward $a$-shifts require
the augmented contour described
in~\cite[Remark after Theorem~1.4]{BorodinCorwinRemenik2013}:
the vertical line is shifted to the right and
supplemented by residue loops around the crossed integer
points, after which the BCR identity gives
$F_v = \E[e^{-sZ_v}]$ at the shifted vertices as well.
The shift ratios \eqref{eq:LG-shift-ratios} and the
Euclidean division operators are independent of the
contour geometry, so the algebraic identities above hold
on any admissible $\ell$.
\end{remark}

\section{O'Connell--Yor Polymer}
\label{sec:OY-dc-verification}

\paragraph{\textbf{System description.}}
Fix $N\geq 1$, $\tau>0$, and a real drift vector
$a=(a_1,\ldots,a_N)$.  Let $B_1,\ldots,B_N$ be independent
Brownian motions, with $B_i$ having drift~$a_i$.  The
O'Connell--Yor partition function is
\[
Z_N^a(\tau)
\defeq
\int_{0<s_1<\cdots<s_{N-1}<\tau}
\exp\Bigl(\sum_{i=1}^{N}\bigl(B_i(s_i)-B_i(s_{i-1})\bigr)\Bigr)
\diff s_1\cdots\diff s_{N-1},
\]
where $s_0\defeq 0$ and $s_N\defeq\tau$.

\paragraph{\textbf{Fredholm determinant formula.}}
The formula of
Borodin, Corwin, and Ferrari~\cite[Theorem~1.17]{BorodinCorwinFerrari2014}
gives, for $\operatorname{Re}u>0$,
\begin{equation}
\label{eq:OY-source-Fredholm}
\E\left[e^{-uZ_N^a(\tau)}\right]
=\det(I+K_u^a)_{L^2(\mathcal C_{\alpha_{\mathrm{ct}},\varphi})},
\end{equation}
where $\alpha_{\mathrm{ct}}>\max_i a_i$,
$\varphi\in(0,\pi/4)$, and
$\mathcal C_{\alpha_{\mathrm{ct}},\varphi}$ is the contour of
Borodin, Corwin, and Ferrari~\cite[Definition~1.16]{BorodinCorwinFerrari2014}.
The kernel is
\begin{equation}
\label{eq:OY-source-kernel}
K_u^a(v,v')
=
\frac{1}{2\pi i}\int_{\mathcal D_v}
\Gamma(-s)\Gamma(1+s)
\prod_{m=1}^{N}
\frac{\Gamma(v-a_m)}{\Gamma(s+v-a_m)}
\frac{u^s e^{v\tau s+\tau s^2/2}}{v+s-v'}
\diff s,
\end{equation}
where $v+\mathcal D_v$ does not intersect
$\mathcal C_{\alpha_{\mathrm{ct}},\varphi}$.

\paragraph{\textbf{References.}}
The O'Connell--Yor polymer was introduced by
O'Connell and Yor~\cite{OConnellYor2001}.  The Fredholm
determinant formula \eqref{eq:OY-source-Fredholm} is due to
Borodin, Corwin, and Ferrari~\cite[Theorem~1.17]{BorodinCorwinFerrari2014}.

\subsection{Kernel reformulation}
\label{sec:OY-kernel-reform}

Setting $u=e^r$ with $r\in\R$, the source kernel
\eqref{eq:OY-source-kernel} rewrites as
\begin{equation}
\label{eq:OY-Cauchy-source}
K_u^a(v,v')
=
\frac{1}{2\pi i}\int_{v+\mathcal D_v}
\frac{\pi}{\sin\pi(v-w)}
\prod_{m=1}^{N}
\frac{\Gamma(v-a_m)}{\Gamma(w-a_m)}
\frac{e^{r(w-v)+\tau(w^2-v^2)/2}}{w-v'}
\diff w,
\end{equation}
after substituting $w=v+s$, applying the reflection
formula $\Gamma(-s)\Gamma(1+s)=\pi/\sin\pi(v-w)$, and
completing the square in the exponent.
The $\Gamma$-ratio and exponential in
\eqref{eq:OY-Cauchy-source} combine into the ratio
$\bar g_a(v;r,\tau)/\bar g_a(w;r,\tau)$, where
\begin{equation}
\label{eq:OY-gauge}
\bar g_a(z;r,\tau)
\defeq
e^{-rz-\tau z^2/2}\prod_{m=1}^{N}\Gamma(z-a_m).
\end{equation}
The gauge transformation
$K_a(v,v')\defeq\bar g_a(v')\bar g_a(v)^{-1}K_u^a(v,v')$
absorbs this ratio and factors the kernel as
\begin{equation}
\label{eq:OY-gauged-kernel}
K_a(v,v')
=
\int X_a(v,w)\frac{Y_a(v')}{w-v'}\diff w,
\end{equation}
where
\begin{align}
X_a(v,w)
&\defeq
\frac{\pi}{2\pi i\sin\pi(v-w)}
\frac{1}{\bar g_a(w;r,\tau)},
\label{eq:OY-XY} \\
Y_a(v')&\defeq \bar g_a(v';r,\tau). \notag
\end{align}
All $w$-integrals use the contour $v+\mathcal D_v$.
Since the gauge is a conjugation, the row and column
factors cancel in every $n$-point minor and
$\det(I + K_a) = \det(I + K_u^a)$
on $L^2(\mathcal C_{\alpha_{\mathrm{ct}},\varphi})$;
absolute convergence of the Fredholm series follows from the
Borodin--Corwin--Ferrari trace-class
estimates~\cite[Theorem~1.17]{BorodinCorwinFerrari2014}.

\subsection{Recurrences and kernel identities}
\label{sec:OY-recurrences}

Fix a nonzero multi-index
$I=(I_1,\ldots,I_N)\in\mathbb Z_{\geq0}^N$ and set
$D\defeq I_1+\cdots+I_N$.  The shift $S_I$ acts on the
drift vector by
\[
S_Ia=(a_1+I_1,\ldots,a_N+I_N).
\]
The Gamma recurrence gives the gauge ratio
\begin{equation}
\label{eq:OY-string-gauge-ratio}
\frac{\bar g_{S_Ia}(z;r,\tau)}{\bar g_a(z;r,\tau)}
=\frac{1}{p(z)},
\qquad
p(z)\defeq
\prod_{i=1}^{N}\prod_{q=1}^{I_i}(z-a_i-q),
\end{equation}
so that
\begin{equation}\label{eq:OY-XY-shift}
X_{S_Ia}(v,w)=p(w)X_a(v,w),
\qquad
Y_{S_Ia}(v')=p(v')^{-1}Y_a(v').
\end{equation}
The contour parameter is chosen so that
$\alpha_{\mathrm{ct}}>\max_i(a_i+I_i)$; the roots of $p$
then lie to the left of $\alpha_{\mathrm{ct}}$, so
$Y_{S_Ia}$ is well-defined on the contour.

The quotient algebra $E \defeq \C[z]_{<D} \cong \C[z]/(p)$
has dimension $D$.  The Euclidean division identity
\eqref{eq:PQ-division} with $\Pi = p$ and multiplier $z$
gives the companion matrix
$M_p \defeq -\Lambda_z \in \End(E)$:
\begin{equation}
\label{eq:OY-companion-division}
zf(z)=p(z)(C_zf)+(M_pf)(z),
\qquad f \in E.
\end{equation}
The Bezoutian identities of
Lemma~\ref{lem:PQ-Bezoutian} give
\begin{equation}\label{eq:OY-Bez-adjoint}
(v'I-M_p)\operatorname{Bez}_p(\cdot,v')=p(v'),
\qquad
C_z\operatorname{Bez}_p(\cdot,v')=1,
\end{equation}
where $\operatorname{Bez}_p(z,v') \defeq (p(z)-p(v'))/(z-v')
\in E$.

The O'Connell--Yor polymer sits in the parabolic framework of
\S\ref{sec:dc-framework-clean} with
\[
\pa_1 = \pa_\tau, \qquad \pa_2 = \pa_r, \qquad S_3 = S_I,
\]
and edge weights
\[
C_1 = 0, \quad \Lambda_1 = 0, \qquad
C_2 = I_E, \quad \Lambda_2 = M_p.
\]
All weights are constant and $C_1 = \Lambda_1 = 0$, so the
parabolic diamond equations of
Proposition~\ref{prop:dc-clean-diamond} are identically satisfied.

Let $H = L^2(\mathcal{C}_{\alpha_{\mathrm{ct}},\varphi})$.
The wave function $\Psi_a \colon E \to H$ and adjoint
$\Phi_a \colon H \to E$ are defined by
\begin{equation}\label{eq:OY-Psi-def}
(\Psi_a f)(v)\defeq \int X_a(v,w)f(w)\diff w,
\qquad f\in E,
\end{equation}
and
\begin{equation}\label{eq:OY-Phi-def}
(\Phi_a h)\defeq \int Y_a(v')h(v')\diff v' \cdot 1_E,
\end{equation}
where $1_E \in E$ is the constant polynomial.  At the
shifted vertex $S_Ia$, the division identity
\eqref{eq:OY-companion-division} factors out $C_zf$:
\begin{equation}\label{eq:OY-Psi-shifted}
(\Psi_{S_Ia}f)(v)\defeq
(C_zf)\int X_{S_Ia}(v,w)\diff w.
\end{equation}
The shifted adjoint inserts the Bezoutian:
\begin{equation}\label{eq:OY-Phi-shifted}
(\Phi_{S_Ia}h)(z)\defeq
\int \operatorname{Bez}_p(z,v')Y_{S_Ia}(v')h(v')\diff v'.
\end{equation}

\begin{lemma}\label{lem:OY-seed-linear}
The waves $\Psi_a$ \eqref{eq:OY-Psi-def} and $\Phi_a$
\eqref{eq:OY-Phi-def} satisfy the parabolic linear and
adjoint linear problems
\eqref{eq:dc-clean-linear-d1}--\eqref{eq:dc-clean-linear-S3}
and
\eqref{eq:dc-clean-adjoint-d1}--\eqref{eq:dc-clean-adjoint-S3}.
Explicitly:
\begin{align}
\pa_\tau\Psi_a&=\tfrac{1}{2}\pa_r^2\Psi_a,
\label{eq:OY-Psi-heat}\\
\Psi_{S_Ia}&=\pa_r\Psi_a-\Psi_aM_p,
\label{eq:OY-Psi-edge}
\shortintertext{and}
\pa_\tau\Phi_a&=-\tfrac{1}{2}\pa_r^2\Phi_a,
\label{eq:OY-Phi-heat}\\
\Phi_a&=-\pa_r\Phi_{S_Ia}-M_p\Phi_{S_Ia}.
\label{eq:OY-Phi-edge}
\end{align}
\end{lemma}

\begin{proof}
The gauge \eqref{eq:OY-gauge} gives
\[
\pa_rX_a(v,w)=wX_a(v,w),
\qquad
\pa_rY_a(v')=-v'Y_a(v'),
\]
\[
\pa_\tau X_a(v,w)=\tfrac{1}{2}w^2 X_a(v,w),
\qquad
\pa_\tau Y_a(v')=-\tfrac{1}{2}v'^2 Y_a(v').
\]
The continuous equations follow directly.  For $f\in E$,
\[
(\pa_\tau\Psi_a f)(v)
=\tfrac{1}{2}\int w^2X_a(v,w)f(w)\diff w
=\tfrac{1}{2}(\pa_r^2\Psi_a f)(v).
\]
For $\Phi_a$,
\[
\pa_\tau\Phi_a h
=-\tfrac{1}{2}\int v'^2Y_a(v')h(v')\diff v' \cdot 1_E
=-\tfrac{1}{2}\pa_r^2\Phi_a h.
\]

For \eqref{eq:OY-Psi-edge}, the division identity
\eqref{eq:OY-companion-division} gives
\[wf(w) - (M_pf)(w) = (C_zf)p(w),\] so
\begin{align*}
\bigl((\pa_r\Psi_a-\Psi_aM_p)f\bigr)(v)
&=
\int X_a(v,w)\bigl(wf(w)-(M_pf)(w)\bigr)\diff w \\
&=
(C_zf)\int X_a(v,w)p(w)\diff w \\
&=(\Psi_{S_Ia}f)(v).
\end{align*}
For \eqref{eq:OY-Phi-edge}, the Bezoutian identity
\eqref{eq:OY-Bez-adjoint} and
$Y_{S_Ia}=Y_a/p$ give
\begin{align*}
\bigl(-\pa_r\Phi_{S_Ia}-M_p\Phi_{S_Ia}\bigr)h
&=
\int (v'I-M_p)\operatorname{Bez}_p(\cdot,v')
Y_{S_Ia}(v')h(v')\diff v' \\
&=
\int p(v')Y_{S_Ia}(v')h(v')\diff v' \cdot 1_E \\
&=\Phi_a h.
\end{align*}
\end{proof}

\begin{lemma}\label{lem:OY-dressing}
Suppose the gauged kernels $K_a$ are trace class on $H$
and the resolvents $(I+K_a)^{-1}$ exist on the Darboux
stencil.  Then $K_a$ is dressing compatible with
$(\Psi_a, \Phi_a)$ in the sense of
Definition~\ref{def:dc-clean-dressing}:
\begin{align}
\pa_rK_a&=\Psi_a\Phi_a,
\label{eq:OY-K-r}\\
\pa_\tau K_a
&=
\tfrac{1}{2}\bigl((\pa_r\Psi_a)\Phi_a
-\Psi_a(\pa_r\Phi_a)\bigr),
\label{eq:OY-K-tau}\\
K_{S_Ia}-K_a&=\Psi_a\Phi_{S_Ia}.
\label{eq:OY-K-shift}
\end{align}
\end{lemma}

\begin{proof}
In $X_a(v,w)$, the factor $1/\bar g_a(w)$ contains
$e^{\tau w^2/2}$, which decays along the angled rays
of the $w$-contour because
$\operatorname{Re}(w^2) \to -\infty$ on these rays.
In $Y_a(v') = \bar g_a(v')$, the factor
$e^{-\tau v'^2/2}$ gives Gaussian decay along the
$v'$-contour.  Together these ensure that the integrals
converge absolutely and differentiation in $r$ and
$\tau$ may be performed under the integral sign.

\noindent\textit{Verification of \eqref{eq:OY-K-r}.}
Differentiating \eqref{eq:OY-gauged-kernel} in $r$, the
factor $(w-v')$ from $\pa_r(X_aY_a)=(w-v')X_aY_a$ cancels
the denominator:
\[
\pa_rK_a(v,v')=\int X_a(v,w)Y_a(v')\diff w.
\]
Since $\Phi_a$ inserts the constant polynomial, this is
$(\Psi_a\Phi_a)(v,v')$.

\medskip
\noindent\textit{Verification of \eqref{eq:OY-K-tau}.}
Differentiating \eqref{eq:OY-gauged-kernel} in $\tau$, the
factorization $(w^2-v'^2)/(w-v')=w+v'$ gives
\[
\pa_\tau K_a(v,v')
=\tfrac{1}{2}
\int X_a(v,w)Y_a(v')(w+v')\diff w.
\]
Since $\Phi_a$ inserts the constant polynomial,
the operator $(\pa_r\Psi_a)\Phi_a$ has kernel
$Y_a(v')\int wX_a(v,w)\diff w$
and $\Psi_a(\pa_r\Phi_a)$ has kernel
$-v'Y_a(v')\int X_a(v,w)\diff w$.
Their difference is
\[
\int X_a(v,w)Y_a(v')(w+v')\diff w = 2\pa_\tau K_a(v,v').
\]

\medskip
\noindent\textit{Verification of \eqref{eq:OY-K-shift}.}
The shift rule \eqref{eq:OY-XY-shift} gives
\begin{align*}
K_{S_Ia}(v,v')-K_a(v,v')
&=
\int X_a(v,w)Y_a(v')
\frac{p(w)p(v')^{-1}-1}{w-v'}\diff w \\
&=
\int X_a(v,w)\operatorname{Bez}_p(w,v')
\frac{Y_a(v')}{p(v')}\diff w \\
&=(\Psi_a\Phi_{S_Ia})(v,v').
\end{align*}
\end{proof}

\subsection{One-point equation}
\label{sec:OY-onepoint}

Define the dressed observable\footnote{The convention $\det(I+K)$
corresponds to spectral parameter $z=-1$ in
Theorem~\ref{thm:dc-clean-general-darboux}; thus
$R_a=(I+K_a)^{-1}$ and
$\mathcal{A}_a = z\Phi_a R_a\Psi_a = -\Phi_a R_a\Psi_a$.}
\begin{equation}\label{eq:OY-A-def}
\mathcal{A}_a \defeq
-\Phi_a (I + K_a)^{-1} \Psi_a
\in \End(E).
\end{equation}

\begin{theorem}\label{thm:OY-finite-string-equation}
Suppose the resolvents $(I+K_a)^{-1}$ exist on the
Darboux stencil.  The dressed observable $\mathcal{A}_a$
satisfies
\begin{equation}
\label{eq:OY-A-equation}
\begin{aligned}
&\pa_\tau(\mathcal{A}_a - \mathcal{A}_{S_3a})
+\tfrac{1}{2}\pa_r^2(\mathcal{A}_a + \mathcal{A}_{S_3a}) \\
&\quad -(\pa_r\mathcal{A}_a)(M_p + \mathcal{A}_a - \mathcal{A}_{S_3a})
+(M_p + \mathcal{A}_a - \mathcal{A}_{S_3a})\pa_r\mathcal{A}_{S_3a}
= 0.
\end{aligned}
\end{equation}
Here $S_3a = S_Ia$ is the single drift shift on
the stencil $\{a, S_Ia\}$.
\end{theorem}

\begin{proof}
The trace-class property and trace-norm differentiability
in $r$ and $\tau$ follow from the Gaussian decay of the
gauged kernel along the BCF contours; dressing
compatibility is Lemma~\ref{lem:OY-dressing}.
Theorem~\ref{thm:dc-clean-general-darboux} applies: the
dressed $\Lambda$-weights
$\widehat\Lambda_1 = \pa_r\mathcal{A}_a$,
$\widehat\Lambda_2 = M_p + \mathcal{A}_a - \mathcal{A}_{S_3a}$
satisfy the parabolic $\Lambda$-diamond equation
\eqref{eq:dc-clean-L-diamond}.  With $C_1 = 0$ and
$C_2 = I_E$, this equation reads
\[
(\pa_\tau - \tfrac{1}{2}\pa_r^2)\widehat\Lambda_2
+ \pa_r\widehat\Lambda_1
- [\widehat\Lambda_1, \widehat\Lambda_2]_{S_3} = 0.
\]
Since $M_p$ is independent of $r$ and $\tau$,
$(\pa_\tau - \frac{1}{2}\pa_r^2)\widehat\Lambda_2
= (\pa_\tau - \frac{1}{2}\pa_r^2)(\mathcal{A}_a - \mathcal{A}_{S_3a})$
and $\pa_r\widehat\Lambda_1 = \pa_r^2\mathcal{A}_a$.
The second-derivative terms combine:
\[
-\tfrac{1}{2}\pa_r^2(\mathcal{A}_a - \mathcal{A}_{S_3a})
+ \pa_r^2\mathcal{A}_a
= \tfrac{1}{2}\pa_r^2(\mathcal{A}_a + \mathcal{A}_{S_3a}).
\]
Substituting into the commutator form gives
\eqref{eq:OY-A-equation}.
\end{proof}

\begin{corollary}\label{cor:OY-trace-defect}
Let $F_a \defeq \det_H(I + K_a)$ and
$\mathcal{B} \defeq \mathcal{A}_a - \mathcal{A}_{S_3a}$.
Then
\begin{equation}\label{eq:OY-trace-defect}
\frac{
\left(D_\tau - \frac{1}{2}D_r^2\right)
F_{S_3a} \cdot F_a
}{
F_{S_3a} F_a
}
=
\frac{1}{2}\operatorname{tr}_E(\mathcal{B}^2)
- \frac{1}{2}(\operatorname{tr}_E\mathcal{B})^2
+ \operatorname{tr}_E(\mathcal{B} M_p).
\end{equation}
\end{corollary}

\begin{proof}
Proposition~\ref{prop:dc-clean-normalized-trace-defect-direct}
applies with $C_1 = 0$, $C_2 = I_E$, $\Lambda_2 = M_p$.
The trace-defect quantity is
$\mathcal{B} = \Delta_2 C_2^{-1}
= [\mathcal{A}, C_2]_{S_3} = \mathcal{A}_a - \mathcal{A}_{S_3a}$,
and $C_1 + \Lambda_2 C_2^{-1} = M_p$.
\end{proof}

\begin{corollary}\label{cor:OY-degree-one}
If $D = 1$, let $j$ be the unique index with $I_j = 1$.
Then $E = \C$, $p(z) = z - b$ with $b = a_j + 1$, and
$M_p = b$.  The quadratic trace-defect term in
\eqref{eq:OY-trace-defect} vanishes ($\mathcal{B}$ is
scalar), and
Corollary~\ref{cor:dc-clean-scalar-trace-defect} gives
\begin{equation}\label{eq:OY-scalar-degree-one}
\left(D_\tau - bD_r - \frac{1}{2}D_r^2\right)
F_{S_3a} \cdot F_a = 0.
\end{equation}
\end{corollary}

\begin{proof}
Since $C_1 + \Lambda_2 C_2^{-1} = M_p = bI$ is scalar,
Corollary~\ref{cor:dc-clean-scalar-trace-defect} applies.
The quadratic trace-defect vanishes for $\dim E = 1$, while
$\operatorname{tr}_E(\mathcal{B}M_p) = b\operatorname{tr}_E\mathcal{B}$
contributes the term $bD_r$.
\end{proof}

%% file: chapters/8-semi-discrete-model-verifications.tex
\chapter{Semi-Discrete Model Verifications}
\label{ch:semi-discrete-model-verifications}

{
  \setlength{\parskip}{0pt}
}

\section{Continuous-Time TASEP}\label{sec:CT-TASEP}

\paragraph{\textbf{System description.}}
The continuous-time totally asymmetric simple exclusion
process (TASEP) is an interacting particle system on $\Z$
with at most one particle per site.  Given a right-finite,
strictly decreasing initial configuration
$Y(0) = y = (y_1 > y_2 > y_3 > \cdots)$,
the dynamics run in continuous time: each particle carries
an independent rate-one exponential clock, and when
particle~$n$'s clock rings it attempts to jump one step
to the right.  The jump is performed only if the
destination site is empty; otherwise it is suppressed.
After each particle's (attempted) jump, its independent
clock is instantaneously reset.

\paragraph{\textbf{Fredholm determinant formula.}}
Fix a time $t \geq 0$, an initial configuration $y$,
finitely many particle labels
$\mathbf{n} = (n_1, \dotsc, n_m)$ with
$1 \leq n_1 < n_2 < \cdots < n_m$, and spatial
thresholds $\mathbf{a} = (a_1, \dotsc, a_m) \in \Z^m$.
The multipoint joint cumulative distribution of the particle
positions is given by the Fredholm determinant formula of Matetski--Quastel--Remenik~\cite{MQR21} on
$\ell^2(\{n_1, \dotsc, n_m\} \times \Z)$:
\[
\PP_y\Bigl(\bigcap_{i=1}^m
\{Y_{n_i}(t) > a_i\}\Bigr)
= \det(I - \bar{\chi}_{\mathbf{a}} K_t
\bar{\chi}_{\mathbf{a}})_{\ell^2(\{n_1, \dotsc, n_m\}
\times \Z)},
\]
where $\bar{\chi}_{\mathbf{a}}(n_i, x) =
\mathbf{1}_{x \leq a_i}$.  The correlation kernel $K_t$ has
block entries
\[
K_t(n_i, x_i; n_j, x_j)
= -Q^{n_j - n_i}(x_i, x_j)\mathbf{1}_{n_i < n_j}
+ \bigl(\mathcal{S}_{-t,-n_i}^*
\overline{\mathcal{S}}_{-t,n_j}^{\operatorname{epi}(y)}
\bigr)(x_i, x_j),
\]
where $A^*$ denotes the transpose kernel,
$A^*(x, y) = A(y, x)$.  The operators defining the kernel
are as follows.

\medskip
\noindent\textit{The transition matrix $Q$.}
The transition matrix $Q$ governs a random walk on $\Z$
taking $\mathrm{Geom}[1/2]$ steps strictly to the left:
\[
Q(z_1, z_2) = 2^{-(z_1 - z_2)}\mathbf{1}_{z_1 > z_2},
\]
with its $n$-th power given by
\[
Q^n(z_1, z_2) = \binom{z_1 - z_2 - 1}{n - 1}
2^{-(z_1 - z_2)}\mathbf{1}_{z_1 - z_2 \geq n}
\]
and admitting the integral representation, for any positively oriented circle $\gamma_\rho$ of radius $\rho \in (0,1)$ around the origin,
\[
Q^n(z_1, z_2) = \frac{1}{2\pi\mathrm{i}}
\oint_{\gamma_\rho} \diff w
\frac{2^{-(z_1 - z_2)}}{w^{z_1 - z_2 - n + 1}}
\Bigl(\frac{1}{1 - w}\Bigr)^n.
\]

\medskip
\noindent\textit{The scattering operators
$\mathcal{S}_{-t,-n}$ and
$\overline{\mathcal{S}}_{-t,n}$.}\footnote{The kernel
operators $\mathcal{S}_{-t,-n}$ and
$\overline{\mathcal{S}}_{-t,n}$ are unrelated to the
lattice shifts $\mathcal{S}_k$ of the diamond framework
(\S\ref{sec:linear-problem}); the notation follows
Matetski--Quastel--Remenik~\cite{MQR21}.}
These operators are defined by the contour integrals
\[
\mathcal{S}_{-t,-n}(z_1, z_2)
= \frac{1}{2\pi\mathrm{i}}
\oint_{\gamma_\rho} \frac{\diff w}{w}
\frac{(1-w)^n}{2^{z_2 - z_1} w^{n + z_2 - z_1}}
e^{t(w - 1/2)},
\]
\[
\overline{\mathcal{S}}_{-t,n}(z_1, z_2)
= \frac{1}{2\pi\mathrm{i}}
\oint_{\gamma_\delta} \frac{\diff w}{w}
\frac{(1-w)^{z_2 - z_1 + n - 1}}{2^{z_1 - z_2} w^{n-1}}
e^{t(w - 1/2)},
\]
where $\gamma_\rho$ is a positively oriented circle of
radius $\rho \in (0, 1)$ around the origin, and
$\gamma_\delta$ is a sufficiently small positively
oriented circle around the origin, chosen so that the only
singularity enclosed is the pole at $w = 0$.

\medskip
\noindent\textit{The epigraph operator
$\overline{\mathcal{S}}_{-t,n}^{\operatorname{epi}(y)}$.}
This is the hitting probability operator defined in terms
of the initial data:
\[
\overline{\mathcal{S}}_{-t,n}^{\operatorname{epi}(y)}
(z_1, z_2)
= \E_{W_0 = z_1}\bigl[
\overline{\mathcal{S}}_{-t,n-\tau}(W_\tau, z_2)
\mathbf{1}_{\tau < n}\bigr],
\]
where $(W_\ell)_{\ell \geq 0}$ is the random walk with
transition matrix $Q$, and
$\tau = \min\{\ell \geq 0 :
W_\ell > y_{\ell+1}\}$
is the hitting time of the strict epigraph of the initial
data by the random walk, with the convention
$\tau = \infty$ if the set is empty.

\paragraph{\textbf{References.}}
The Fredholm determinant formula above is the
Matetski--Quastel--Remenik formula
(Theorem~2.6, equations~(2.27)--(2.31)
of~\cite{MQR21}) for continuous-time TASEP.
Continuous-time TASEP is the totally asymmetric
specialization of the exclusion processes introduced by
Spitzer~\cite{Spitzer1970}.

\subsection{Kernel reformulation}\label{sec:CT-TASEP-kernel-reform}

The extended-kernel Fredholm determinant admits a reduction to a single-space Fredholm determinant on $\ell^2(\Z)$.  Define
\[
\psi_{t,r,n}(x) \defeq \mathcal{S}_{-t,-n}(x, r),
\qquad
\phi_{t,r,n}(x) \defeq
\overline{\mathcal{S}}_{-t,n}^{\operatorname{epi}(y)}(x, r),
\]
where we suppress the dependence of $\phi$ on the initial
condition $y$.  The propagator is an instance of the
directed-path propagator
(Definition~\ref{def:directed-path-propagator}) with
transition kernels
$Q_{\ell,\ell'} = Q^{n_{\ell'} - n_\ell}$ for
$\ell < \ell'$.  Explicitly, for $i > j$,
\begin{equation}\label{eq:CT-TASEP-propagator}
[B_{\mathbf{a},\mathbf{n}}]_{ij}(r, r')
\defeq \sum_{k=1}^{i-j}
  \sum_{j = \ell_0 < \ell_1 < \cdots < \ell_k = i}
  \sum_{\substack{\xi_1 \leq a_{\ell_1} \\ \vdots \\
  \xi_{k-1} \leq a_{\ell_{k-1}}}}
  \prod_{s=0}^{k-1}
  \Bigl(-Q^{n_{\ell_{s+1}} - n_{\ell_s}}
  (\xi_s, \xi_{s+1})\Bigr),
\end{equation}
with $\xi_0 = r'$ and $\xi_k = r$,
and $[B_{\mathbf{a},\mathbf{n}}]_{ij} = 0$ for $i \leq j$.
The internal vertices $\xi_1, \dotsc, \xi_{k-1}$ are
constrained by the cutoffs, but the endpoints $r, r'$ are
arbitrary.

\begin{lemma}\label{lem:CT-TASEP-kernel}
We have
\[
\PP_y\Bigl(\bigcap_{i=1}^m
\{Y_{n_i}(t) > a_i\}\Bigr)
= \det(I - K_{t, \mathbf{a}, \mathbf{n}})_{\ell^2(\Z)},
\]
where
\[
K_{t, \mathbf{a}, \mathbf{n}}(x, x')
= \sum_{1 \leq j \leq i \leq m}
  \sum_{r \leq a_i} \sum_{r' \leq a_j}
  \psi_{t, r, n_i}(x)
  \bigl[\delta_{ij}\delta_{r,r'}
  + [B_{\mathbf{a},\mathbf{n}}]_{ij}(r, r')\bigr]
  \phi_{t, r', n_j}(x').
\]
\end{lemma}

\begin{proof}
The proof follows the same strategy as
Lemma~\ref{lem:RBJ-kernel}: conjugation by exponential
weights to obtain a trace-class realization, followed by
Sylvester's identity and identification of the
inverse-transpose entries with the propagator.  We provide
the details that differ from the right Bernoulli case and
indicate where the remaining arguments transfer.

Choose $\beta \in (1/2, 1)$, and choose real numbers
\[
1 < \omega_m < \omega_{m-1} < \cdots < \omega_1
< \beta^{-1}.
\]
Since $\beta > 1/2$, we have $\beta^{-1} < 2$, so
$\omega_i \in (1, 2)$ for each $i$.
For a function
$f \colon \{n_1, \dotsc, n_m\} \times \Z \to \R$, let
$(\mathscr{C}f)_i(x) \defeq \omega_i^x f_{n_i}(x)$,
so that for a block kernel $M$ on
$\{n_1, \dotsc, n_m\} \times \Z$,
$(\mathscr{C}M\mathscr{C}^{-1})_{ij}(x, x')
= \omega_i^x M_{ij}(x, x')\omega_j^{-x'}$.

Define $U(n_i, r; x) \defeq
\mathbf{1}_{r \leq a_i}\psi_{t,r,n_i}(x)$ and
$V(x; n_j, r') \defeq
\phi_{t,r',n_j}(x)\mathbf{1}_{r' \leq a_j}$.
The Matetski--Quastel--Remenik kernel decomposes as
$I - \bar{\chi}_{\mathbf{a}} K_t \bar{\chi}_{\mathbf{a}}
= I + \bar{\chi}_{\mathbf{a}} L \bar{\chi}_{\mathbf{a}}
- UV$,
where $L$ is the strictly upper-triangular block operator
with entries $L_{ij}(x, x') = Q^{n_j - n_i}(x, x')
\mathbf{1}_{i < j}$.  Set
$\widetilde{N} \defeq
\mathscr{C}\bar{\chi}_{\mathbf{a}} L
\bar{\chi}_{\mathbf{a}}\mathscr{C}^{-1}$,
$\widetilde{U} \defeq \mathscr{C}U$, and
$\widetilde{V} \defeq V\mathscr{C}^{-1}$.
Since $\mathscr{C}$ multiplies each row of a finite minor
by $\omega_i^x$ and each column by $\omega_j^{-x'}$, the
Fredholm series agrees term by term with that of the
conjugated kernel.

\medskip
\noindent\textit{Trace-class estimate for $\widetilde{N}$.}
For $i < j$, let $d = n_j - n_i$.  The $(i,j)$-block of
$\widetilde{N}$ has kernel
\[
\omega_i^x Q^d(x, x')\omega_j^{-x'}
\mathbf{1}_{x \leq a_i}\mathbf{1}_{x' \leq a_j}.
\]
The explicit formula
$Q^d(s) = \binom{s-1}{d-1}2^{-s}\mathbf{1}_{s \geq d}$
for $s = x - x'$ gives the bound
\[
\sum_{x \leq a_i}\sum_{x' \leq a_j}
\bigl|\omega_i^x Q^d(x,x')\omega_j^{-x'}\bigr|
\leq \sum_{s \geq d} Q^d(s)\omega_i^s
\sum_{x' \leq a_j}
\Bigl(\frac{\omega_i}{\omega_j}\Bigr)^{x'}.
\]
The second sum converges because $\omega_i/\omega_j > 1$
for $i < j$.  The first sum converges because
$\sum_{s \geq d} Q^d(s)\omega_i^s
= \sum_{s \geq d}\binom{s-1}{d-1}
(\omega_i/2)^s$,
and the binomial coefficient grows polynomially in $s$
while the geometric factor $(\omega_i/2)^s$ decays
exponentially since $\omega_i < 2$.  Thus every nonzero
block of $\widetilde{N}$ has absolutely summable kernel,
hence is trace class.  Since there are only finitely many
blocks, $\widetilde{N} \in \mathfrak{S}_1$.  Moreover,
$\widetilde{N}$ is strictly upper triangular in the block
indices, hence nilpotent.

\medskip
\noindent\textit{Hilbert--Schmidt estimate for
$\widetilde{U}$.}
The contour formula for $\psi$ gives
\[
\psi_{t,r,n}(x) = \frac{1}{2\pi\mathrm{i}}
\oint_{\gamma_\rho} \frac{\diff w}{w}
\frac{(1-w)^n}{2^{r-x} w^{n+r-x}} e^{t(w-1/2)}.
\]
The integral extracts the coefficient of $w^{n+r-x}$ in
$(1-w)^n e^{tw}$, up to a constant.  Since $(1-w)^n$ is a
polynomial of degree $n$ and $e^{tw}$ is a power series in
non-negative powers of $w$, this coefficient vanishes when
$n + r - x < 0$, giving $\psi_{t,r,n}(x) = 0$ for
$x > r + n$.  Moreover, the integrand depends on $r$ and
$x$ only through $r - x$, so $\psi_{t,r,n}(x) = g(r - x)$
for a fixed function $g$ supported on $\{-n, -n+1,
-n+2, \dotsc\}$.  For large $s$, $g(s)$ is the coefficient
of $w^{n+s}$ in $(1-w)^n e^{tw}$ (up to a constant), which
decays as $t^s/s!$ since the dominant contribution for
$s \gg n$ comes from $e^{tw}$ alone.  In particular,
$\lVert g\rVert_{\ell^2} < \infty$, and
$\lVert\psi_{t,r,n}\rVert_{\ell^2(\Z)} = \lVert g\rVert_{\ell^2}$ is
independent of $r$.  It follows that
\[
\lVert\widetilde{U}\rVert_2^2
= \sum_{i=1}^m \sum_{r \leq a_i}
\omega_i^{2r}\lVert\psi_{t,r,n_i}\rVert_{\ell^2}^2
\leq C\sum_{i=1}^m \sum_{r \leq a_i} \omega_i^{2r}
< \infty,
\]
since $\omega_i > 1$.  Thus $\widetilde{U} \in
\mathfrak{S}_2$.

\medskip
\noindent\textit{Hilbert--Schmidt estimate for
$\widetilde{V}$.}
Let $a_{\max} = \max_{1 \leq j \leq m} a_j$.
Choose $\delta > 0$ small enough that
$\eta \defeq \sup_{w \in \gamma_\delta}
|1/(2(1-w))| < \beta$.
This is possible since $|1/(2(1-w))| \leq 1/(2(1-\delta))$
on $\gamma_\delta$, which is less than any $\beta > 1/2$
for $\delta$ sufficiently small.  The contour formula for
$\overline{\mathcal{S}}_{-t,k}(z, r)$ then gives
\[
\bigl|\overline{\mathcal{S}}_{-t,k}(z, r)\bigr|
\leq C\beta^{z-r}
\]
for all $1 \leq k \leq n_m$, $z > y_{n_m}$, and
$r \leq a_{\max}$, by the same argument as in the proof of
Lemma~\ref{lem:RBJ-kernel}: for $z \geq r$, the key factor
is $2^{-(z-r)}|1-w|^{-(z-r)} = |1/(2(1-w))|^{z-r}
\leq \eta^{z-r}$; for $z < r$, the restrictions
$z > y_{n_m}$ and $r \leq a_{\max}$ confine $(z,r)$ to a
finite set.

The bound on $\phi$ follows exactly as in
Lemma~\ref{lem:RBJ-kernel}.  The hitting-time
representation gives
$|\phi_{t,r,n}(x)| \leq C\beta^{x-r}$ for
$x > y_{n_m}$, and $\phi_{t,r,n}(x) = 0$ for
$x \leq y_{n_m}$ (since the random walk $(W_\ell)$ is
strictly decreasing, $\tau = \infty$ whenever
$x \leq y_{n_m}$).  Combining:
\begin{equation}\label{eq:CT-TASEP-phi-bound}
|\phi_{t,r,n}(x)| \leq C\beta^{x-r}\mathbf{1}_{x > y_{n_m}}.
\end{equation}
Since $\beta < 1$,
$\lVert\phi_{t,r,n}\rVert_{\ell^2}^2 \leq C\beta^{-2r}$, and
\[
\lVert\widetilde{V}\rVert_2^2
= \sum_{j=1}^m \sum_{r' \leq a_j}
\omega_j^{-2r'}\lVert\phi_{t,r',n_j}\rVert_{\ell^2}^2
\leq C\sum_{j=1}^m \sum_{r' \leq a_j}
(\omega_j\beta)^{-2r'} < \infty,
\]
since $\omega_j\beta < \beta^{-1}\beta = 1$.
Thus $\widetilde{V} \in \mathfrak{S}_2$.

\medskip
\noindent\textit{Sylvester reduction and propagator
identification.}
The conjugated extended kernel
$-\widetilde{N} + \widetilde{U}\widetilde{V}$ is trace
class, since $\widetilde{N} \in \mathfrak{S}_1$ and
$\widetilde{U}\widetilde{V} \in \mathfrak{S}_1$ as a
product of two Hilbert--Schmidt operators.  Factoring out
the nilpotent triangular part:
\[
I + \widetilde{N} - \widetilde{U}\widetilde{V}
= (I + \widetilde{N})\bigl(I - (I + \widetilde{N})^{-1}
\widetilde{U}\widetilde{V}\bigr).
\]
Since $\widetilde{N}$ is nilpotent,
$\det(I + \widetilde{N}) = 1$, and
$(I + \widetilde{N})^{-1}$ is given by the finite Neumann
series.  The product
$(I + \widetilde{N})^{-1}\widetilde{U} \in \mathfrak{S}_2$
(since $(I + \widetilde{N})^{-1}$ is bounded and
$\mathfrak{S}_2$ is an operator ideal), so Sylvester's
identity applies.  Moreover, the conjugation weights cancel
in the product
$\widetilde{V}(I + \widetilde{N})^{-1}\widetilde{U}
= V(I + \bar{\chi}_{\mathbf{a}} L
\bar{\chi}_{\mathbf{a}})^{-1}U$.
Setting $M \defeq I + \bar{\chi}_{\mathbf{a}} L
\bar{\chi}_{\mathbf{a}}$:
\[
\det\bigl(I - \bar{\chi}_{\mathbf{a}} K_t
\bar{\chi}_{\mathbf{a}}\bigr)
= \det_{\ell^2(\Z)}\bigl(I - VM^{-1}U\bigr).
\]
Since the Fredholm determinant of a trace-class operator is
invariant under transposition, we may replace the kernel of
$VM^{-1}U$ by its transpose.  Writing
$[M^{-\top}]_{ij}(r, r') = [M^{-1}]_{ji}(r', r)$, the
transposed kernel is
\[
\sum_{i,j=1}^m \sum_{r \leq a_i} \sum_{r' \leq a_j}
\psi_{t,r,n_i}(x) [M^{-\top}]_{ij}(r, r')
\phi_{t,r',n_j}(x').
\]

It remains to identify the inverse-transpose entries with
the propagator.  For $i = j$, the inverse-transpose
contributes $\delta_{r,r'}$ and
$[B_{\mathbf{a},\mathbf{n}}]_{ii} = 0$, so the combined factor is
$\delta_{ij}\delta_{r,r'}$.  For $i < j$, both
$[M^{-\top}]_{ij}$ and $[B_{\mathbf{a},\mathbf{n}}]_{ij}$ vanish,
the former because the inverse-transpose is block lower
triangular.  For $i > j$, expanding the finite Neumann
series in the block indices and restricting to the cutoff
region $r \leq a_i$, $r' \leq a_j$ gives
\[
[M^{-\top}]_{ij}(r, r') = [B_{\mathbf{a},\mathbf{n}}]_{ij}(r, r').
\]
Indeed, a nonzero term in the Neumann series has the form
$j = \ell_0 < \ell_1 < \cdots < \ell_k = i$,
$\xi_0 = r'$, $\xi_k = r$, with internal constraints
$\xi_s \leq a_{\ell_s}$ for $s = 1, \dotsc, k-1$, and
weight
$(-1)^k\prod_{s=0}^{k-1}
Q^{n_{\ell_{s+1}} - n_{\ell_s}}(\xi_s, \xi_{s+1})$;
the sign $(-1)^k$ is absorbed by the $k$ factors of
$(-Q)$ in the definition of
$[B_{\mathbf{a},\mathbf{n}}]_{ij}(r, r')$.  The support condition
$Q^d(x, x') = 0$ unless $x - x' \geq d$ shows that the
internal sums are finite for fixed endpoints: along any
nonzero chain,
$\xi_s - \xi_{s+1} \geq n_{\ell_{s+1}} - n_{\ell_s}
\geq 1$,
so the spatial chain is strictly decreasing.

The transposed kernel is therefore the claimed single-space
kernel $K_{t, \mathbf{a}, \mathbf{n}}$.
\end{proof}

For the remainder of this section, we work with the kernel
$K_{t, \mathbf{a}, \mathbf{n}}$ on $\ell^2(\Z)$, with
ordered labels $1 \leq n_1 < \cdots < n_m$ as in the
Fredholm formula.  Write
$\mathbf{1} = (1, \ldots, 1) \in \Z^m$.

\subsection{Recurrences and kernel identities}
\label{sec:CT-TASEP-recurrences}

Continuous-time TASEP naturally sits in the dual lattice
presentation of the semi-discrete framework
(\S\ref{sec:sc-dual}): the seed recurrences take the form
$\pa_t\psi(Tu)$ rather than $\pa_t\psi(u)$, with
a discrete particle-label shift in the $S_2$ direction.
The product graph architecture of \S\ref{sec:product-graph}
provides the computational scaffolding: directed-path
propagator $B_{\mathbf{a},\mathbf{n}}$, distinguished shift
$e^{\pa_a}$, and product-graph wave functions
$\Psi, \Phi$, with $\pa_t$ replacing one discrete
shift.  The framework identification proceeds in three
stages: the seed functions satisfy the dual linear problem
(Lemma~\ref{lem:CT-TASEP-seed-linear}), the product-graph
wave functions inherit this structure
(Lemma~\ref{lem:CT-TASEP-WF}), and the multipoint equation
follows from the dual Darboux theorem
(Theorem~\ref{thm:sc-dual-Darboux}).

\begin{lemma}\label{lem:CT-TASEP-seed-linear}
The seed functions $\psi_{t,r,n}$ and $\phi_{t,r,n}$
defined in \S\ref{sec:CT-TASEP-kernel-reform} satisfy the
dual linear problem
\eqref{eq:sc-dual-linear-1}--\eqref{eq:sc-dual-linear-j}
and its adjoint
\eqref{eq:sc-dual-adjoint-1}--\eqref{eq:sc-dual-adjoint-j}
with shifts
$T = e^{\pa_r}$, $\pa_1 = \pa_t$,
$S_2 = e^{\pa_n}$ and constant scalar edge weights
$(c_1, \lambda_1) = (-\tfrac{1}{2}, -\tfrac{1}{2})$,
$(c_2, \lambda_2) = (2, 1)$.
Explicitly:
\begin{align}
\pa_t \psi_{t,r+1,n}
&= -\tfrac{1}{2}\bigl(\psi_{t,r+1,n}
   - \psi_{t,r,n}\bigr),
\label{eq:CT-psi-t} \\
\psi_{t,r,n+1}
&= 2\psi_{t,r+1,n} - \psi_{t,r,n},
\label{eq:CT-psi-n}
\end{align}
and
\begin{align}
\pa_t \phi_{t,r+1,n}
&= -\tfrac{1}{2}\bigl(\phi_{t,r+2,n}
   - \phi_{t,r+1,n}\bigr),
\label{eq:CT-phi-t} \\
\phi_{t,r+1,n}
&= 2\phi_{t,r,n+1} - \phi_{t,r+1,n+1}.
\label{eq:CT-phi-n}
\end{align}
\end{lemma}

\begin{proof}
\noindent\textit{Verification of \eqref{eq:CT-psi-t}.}
From the contour representation,
\[
\psi_{t,r,n}(x) = \frac{1}{2\pi\mathrm{i}}
\oint_{\gamma_\rho} \frac{\diff w}{w}
\frac{(1-w)^n}{2^{r-x} w^{n+r-x}}
e^{t(w - 1/2)}.
\]
Differentiating in $t$ produces the factor
$(w - 1/2)$ in the integrand:
\[
\pa_t \psi_{t,r,n}(x)
= \frac{1}{2\pi\mathrm{i}}
\oint_{\gamma_\rho} \frac{\diff w}{w}
\frac{(1-w)^n}{2^{r-x} w^{n+r-x}}
(w - \tfrac{1}{2})
e^{t(w - 1/2)}.
\]
On the other hand, a shift of $r \to r{-}1$ in
$\psi_{t,r,n}$ multiplies the integrand by $2w$, so
$-\frac{1}{2}(\psi_{t,r,n} - \psi_{t,r-1,n})$ produces
the same factor $-\frac{1}{2}(1 - 2w) = w - 1/2$.
The identity holds at every $r$; evaluating at $r + 1$
gives \eqref{eq:CT-psi-t}.

\medskip
\noindent\textit{Verification of \eqref{eq:CT-psi-n}.}
A shift of $r \to r{+}1$ in $\psi_{t,r,n}$ multiplies
the integrand by $1/(2w)$.  The right-hand side of
\eqref{eq:CT-psi-n} therefore produces the factor
$2 \cdot 1/(2w) - 1 = (1-w)/w$, giving
\[
2\psi_{t,r+1,n}(x) - \psi_{t,r,n}(x)
= \frac{1}{2\pi\mathrm{i}}
\oint_{\gamma_\rho} \frac{\diff w}{w}
\frac{(1-w)^{n+1}}{2^{r-x} w^{n+1+r-x}}
e^{t(w - 1/2)}
= \psi_{t,r,n+1}(x).
\]

\medskip
\noindent\textit{Free recurrences for
$\overline{\mathcal{S}}$.}
We record two identities for the free operator that will
be applied inside the hitting-time representation of
$\phi$.  For $k \geq 1$,
\begin{equation}\label{eq:CT-Sbar-recurrence-n}
2\overline{\mathcal{S}}_{-t,k}(z, r{-}1)
- \overline{\mathcal{S}}_{-t,k}(z, r)
= \overline{\mathcal{S}}_{-t,k-1}(z, r).
\end{equation}
Indeed, the $r$-dependent factors in the contour formula
are $\tfrac{(1{-}w)^{r-z+k-1}}{2^{z-r}}$.  A shift of
$r \to r{-}1$ and the prefactor of $2$ together replace
$(1{-}w)^{r-z+k-1}$ by $(1{-}w)^{r-z+k-2}$.
The left-hand side of \eqref{eq:CT-Sbar-recurrence-n}
therefore produces the factor
$1 - (1{-}w) = w$, reducing $w^{-(k-1)}$ to
$w^{-(k-2)}$ and confirming
\eqref{eq:CT-Sbar-recurrence-n}.

We will also use
\begin{equation}\label{eq:CT-Sbar-recurrence-t}
\pa_t \overline{\mathcal{S}}_{-t,k}(z, r)
= -\tfrac{1}{2}\bigl(
\overline{\mathcal{S}}_{-t,k}(z, r{+}1)
- \overline{\mathcal{S}}_{-t,k}(z, r)\bigr).
\end{equation}
Differentiating $\overline{\mathcal{S}}_{-t,k}(z, r)$
in $t$ produces the factor $(w - 1/2)$ in the integrand.
On the other hand, a shift of $r \to r{+}1$ multiplies
the integrand by $2(1{-}w)$, so the right-hand side of
\eqref{eq:CT-Sbar-recurrence-t} produces the factor
$-\frac{1}{2}[2(1{-}w) - 1] = w - 1/2$,
confirming \eqref{eq:CT-Sbar-recurrence-t}.

Finally, $\overline{\mathcal{S}}_{-t,0} \equiv 0$.
When $k = 0$, the factor $w^{-(k-1)} = w$ cancels the
$w^{-1}$ in the measure, removing the pole at the
origin.  The remaining integrand is analytic inside
$\gamma_\delta$ (the only other possible singularity, at
$w = 1$, lies outside $\gamma_\delta$), and the integral
vanishes by Cauchy's theorem.

\medskip
\noindent\textit{Verification of \eqref{eq:CT-phi-t}.}
The exponential bound on $\overline{\mathcal{S}}$
established in Lemma~\ref{lem:CT-TASEP-kernel}, together
with the finite exponential moments of the geometric
left-jump walk, implies that the expectations defining
$\phi$ are absolutely convergent.  Hence the free
identities above may be applied inside the hitting-time
expectation.

Let $(W_\ell)_{\ell \geq 0}$ denote the random walk with
transition matrix $Q$ and $\tau$ its hitting time of the
strict epigraph of the initial data, as in the definition
of $\phi$.  The particle label $n$ is not shifted, so the
indicator $\mathbf{1}_{\tau < n}$ is unchanged.  Applying
\eqref{eq:CT-Sbar-recurrence-t} with $k = n - \tau$ on
the event $\{\tau < n\}$, we obtain
\begin{align*}
\pa_t \phi_{t,r,n}(x)
&= \E_{W_0 = x}\bigl[
\pa_t \overline{\mathcal{S}}_{-t,n-\tau}
(W_\tau, r)\mathbf{1}_{\tau < n}\bigr] \\
&= -\tfrac{1}{2}\E_{W_0 = x}\bigl[
\bigl(\overline{\mathcal{S}}_{-t,n-\tau}(W_\tau, r{+}1)
- \overline{\mathcal{S}}_{-t,n-\tau}(W_\tau, r)\bigr)
\mathbf{1}_{\tau < n}\bigr] \\
&= -\tfrac{1}{2}\bigl(\phi_{t,r+1,n}(x)
- \phi_{t,r,n}(x)\bigr).
\end{align*}
The interchange of $\pa_t$ and the expectation is
justified by dominated convergence: on $\gamma_\delta$
the factor $(w - 1/2)$ is bounded, so
$|\pa_t \overline{\mathcal{S}}_{-t,k}(z, r)|
\leq C\beta^{z-r}$ by the same contour estimate used
above.  The exponential bound on
$\overline{\mathcal{S}}$ and the finite exponential
moments of the geometric left-jump walk ensure that
this dominating function is integrable.
The identity holds at every $r$; evaluating at $r + 1$
gives \eqref{eq:CT-phi-t}.

\medskip
\noindent\textit{Verification of \eqref{eq:CT-phi-n}.}
The identity
$\phi_{t,r,n-1} = 2\phi_{t,r-1,n} - \phi_{t,r,n}$
holds at every $(r, n)$; evaluating at
$(r{+}1, n{+}1)$ gives \eqref{eq:CT-phi-n}.  Indeed,
\begin{align*}
2\phi_{t,r-1,n}(x) - \phi_{t,r,n}(x)
= \E_{W_0 = x}\Bigl[
\bigl(2\overline{\mathcal{S}}_{-t,n-\tau}(W_\tau, r{-}1)
- \overline{\mathcal{S}}_{-t,n-\tau}(W_\tau, r)\bigr)
\mathbf{1}_{\tau < n}\Bigr].
\end{align*}
The expectation splits according to whether
$\tau \leq n{-}2$ or $\tau = n{-}1$.  On the event
$\{\tau \leq n{-}2\}$, set $k = n - \tau \geq 2$, and
\eqref{eq:CT-Sbar-recurrence-n} gives
\begin{align*}
2\overline{\mathcal{S}}_{-t,n-\tau}(W_\tau, r{-}1)
- \overline{\mathcal{S}}_{-t,n-\tau}(W_\tau, r)
= \overline{\mathcal{S}}_{-t,n-1-\tau}(W_\tau, r).
\end{align*}
On the boundary event $\{\tau = n{-}1\}$, the same
recurrence produces
\[\overline{\mathcal{S}}_{-t,0}(W_{n-1}, r) = 0.\]
Therefore the boundary contribution vanishes, and the
expression reduces to
\[
\E_{W_0 = x}\bigl[
\overline{\mathcal{S}}_{-t,n-1-\tau}(W_\tau, r)
\mathbf{1}_{\tau < n-1}\bigr]
= \phi_{t,r,n-1}(x).
\]
\end{proof}

\medskip

The propagator $B_{\mathbf{a},\mathbf{n}}$ defined by
\eqref{eq:CT-TASEP-propagator} is an instance of the
directed-path propagator
(Definition~\ref{def:directed-path-propagator}) with
transition kernels
$Q_{\ell,\ell'} = Q^{n_{\ell'} - n_\ell}$ for
$\ell < \ell'$.  The following properties govern its
dependence on the parameters $(\mathbf{a}, \mathbf{n})$.

\begin{lemma}\label{lem:CT-TASEP-propagator}
The propagator $B_{\mathbf{a},\mathbf{n}}$ satisfies:
\begin{enumerate}[label=\textup{(\roman*)},
leftmargin=*]
\item For all $i > j$ and $r, r' \in \Z$,
\begin{equation}\label{eq:CT-Ba-covariance}
[B_{\mathbf{a}+\mathbf{1},\mathbf{n}}]_{ij}(r{+}1, r'{+}1)
= [B_{\mathbf{a},\mathbf{n}}]_{ij}(r, r').
\end{equation}

\item For all $i > j$ and $r, r' \in \Z$,
\begin{align}\label{eq:CT-Ba-splitting}
[B_{\mathbf{a}+\mathbf{1},\mathbf{n}}]_{ij}(r, r')
- [B_{\mathbf{a},\mathbf{n}}]_{ij}(r, r')
= \sum_{j < \mu < i}
[B_{\mathbf{a},\mathbf{n}}]_{i\mu}(r, a_\mu{+}1)
[B_{\mathbf{a}+\mathbf{1},\mathbf{n}}]_{\mu j}
(a_\mu{+}1, r').
\end{align}

\item The propagator is $t$-independent and
invariant under uniform shifts of the particle
labels: replacing each $n_i$ by $n_i + 1$ leaves
$[B_{\mathbf{a},\mathbf{n}}]_{ij}(r, r')$ unchanged.
\end{enumerate}
\end{lemma}

\begin{proof}
The standing hypotheses of
\S\ref{sec:directed-path} are satisfied: translation
invariance of $Q^n$ is inherited from that of $Q$,
and $\mathcal{T}$-invariance holds because the
transition kernels $Q^{n_{\ell'} - n_\ell}$ depend
only on the particle-label differences, which are
unaffected by the threshold shift
$\mathbf{a} \to \mathbf{a} + \mathbf{1}$.  Under
the identification
$B_u \leftrightarrow B_{\mathbf{a},\mathbf{n}}$ and
$B_{\mathcal{T}u} \leftrightarrow
B_{\mathbf{a}+\mathbf{1},\mathbf{n}}$, properties~(i) and~(ii)
are $\mathcal{T}$-covariance and
$\mathcal{T}$-splitting respectively, both
established in
Proposition~\ref{prop:directed-path-T}.

For~(iii), the propagator entries depend on
$u = (t, \mathbf{a}, \mathbf{n})$ only through the
thresholds $\mathbf{a}$ and the particle-label
differences $n_{\ell'} - n_\ell$: each transition
kernel $Q^{n_{\ell'} - n_\ell}$ and each threshold
constraint $\xi_s \leq a_{\ell_s}$ depends only on
these quantities.  Replacing each $n_i$ by $n_i + 1$
preserves all differences, and no factor of the
propagator carries $t$-dependence.
\end{proof}

\medskip

For $t \geq 0$, $\mathbf{a} \in \Z^m$, and
$1 \leq n_1 < \cdots < n_m$, define the
product-graph wave functions
$\Psi^{t,\mathbf{a},\mathbf{n}} \in
\Hom(\R^m, \ell^2(\Z))$ and
$\Phi^{t,\mathbf{a},\mathbf{n}} \in
\Hom(\ell^2(\Z), \R^m)$ componentwise by
\begin{align}
\Psi_p^{t,\mathbf{a},\mathbf{n}}
&\defeq \psi_{t,a_p,n_p}
+ \sum_{j > p}
  \sum_{w \leq a_j - 1}
  [B_{\mathbf{a}-\mathbf{1},\mathbf{n}}]_{jp}(w, a_p)
  \psi_{t,w,n_j},
\label{eq:CT-Psi-def} \\
\Phi_p^{t,\mathbf{a},\mathbf{n}}
&\defeq \phi_{t,a_p,n_p}
+ \sum_{j < p}
  \sum_{z \leq a_j}
  [B_{\mathbf{a},\mathbf{n}}]_{pj}(a_p, z)
  \phi_{t,z,n_j}.
\label{eq:CT-Phi-def}
\end{align}
These are the explicit forms of the abstract
constructions
\eqref{eq:Psi-def}--\eqref{eq:Phi-def} on the product
graph with distinguished shift $e^{\pa_a}$.
Under the identification
$u = (t, \mathbf{a}, \mathbf{n})$, the maps
$u \mapsto \Psi^{t,\mathbf{a},\mathbf{n}}$ and
$u \mapsto \Phi^{t,\mathbf{a},\mathbf{n}}$ become the
wave functions of the dual lattice framework of
\S\ref{sec:sc-dual}, with shifts
$T = e^{\pa_{\mathbf{a}}}$,
$\pa_1 = \pa_t$,
$S_2 = e^{\pa_{\mathbf{n}}}$ acting by
\[
Tu = (t, \mathbf{a}{+}\mathbf{1}, \mathbf{n}),
\qquad
S_2 u = (t, \mathbf{a}, \mathbf{n}{+}\mathbf{1}).
\]
The constant scalar edge weights are
$C_1 = -\tfrac{1}{2}I_m$,
$\Lambda_1 = -\tfrac{1}{2}I_m$,
$C_2 = 2I_m$,
$\Lambda_2 = I_m$.
These form a trivial solution of the dual diamond
equations (Corollary~\ref{cor:sc-dual-diamonds}).

\begin{lemma}\label{lem:CT-TASEP-WF}
The product-graph wave functions satisfy the dual
semi-discrete linear problem
\eqref{eq:sc-dual-linear-1}--\eqref{eq:sc-dual-linear-j}
and its adjoint
\eqref{eq:sc-dual-adjoint-1}--\eqref{eq:sc-dual-adjoint-j}
with the shifts and edge weights above.
Explicitly:
\begin{align}
\pa_t \Psi^{t,\mathbf{a}+\mathbf{1},\mathbf{n}}
&= -\tfrac{1}{2}\bigl(\Psi^{t,\mathbf{a}+\mathbf{1},\mathbf{n}}
   - \Psi^{t,\mathbf{a},\mathbf{n}}\bigr),
\label{eq:CT-WF-psi-t} \\
\Psi^{t,\mathbf{a},\mathbf{n}+\mathbf{1}}
&= 2\Psi^{t,\mathbf{a}+\mathbf{1},\mathbf{n}}
   - \Psi^{t,\mathbf{a},\mathbf{n}},
\label{eq:CT-WF-psi-n}
\end{align}
and
\begin{align}
\pa_t \Phi^{t,\mathbf{a}+\mathbf{1},\mathbf{n}}
&= -\tfrac{1}{2}\bigl(\Phi^{t,\mathbf{a}+\mathbf{2},\mathbf{n}}
   - \Phi^{t,\mathbf{a}+\mathbf{1},\mathbf{n}}\bigr),
\label{eq:CT-WF-phi-t} \\
\Phi^{t,\mathbf{a}+\mathbf{1},\mathbf{n}}
&= 2\Phi^{t,\mathbf{a},\mathbf{n}+\mathbf{1}}
   - \Phi^{t,\mathbf{a}+\mathbf{1},\mathbf{n}+\mathbf{1}}.
\label{eq:CT-WF-phi-n}
\end{align}
\end{lemma}

\begin{proof}
\noindent\textit{Discrete recurrences
\eqref{eq:CT-WF-psi-n}, \eqref{eq:CT-WF-phi-n}.}
Under the identification
$\mathcal{T} = e^{\pa_{\mathbf{a}}}$,
$\mathcal{S}_2 = e^{\pa_{\mathbf{n}}}$,
$(C_2, D_2) = (2I_m, I_m)$, the definitions
\eqref{eq:CT-Psi-def}--\eqref{eq:CT-Phi-def} match
the abstract constructions
\eqref{eq:Psi-def}--\eqref{eq:Phi-def} on the product
graph $\mathcal{G}^m$.  Property~(iii) of
Lemma~\ref{lem:CT-TASEP-propagator} is the
$\mathcal{S}_2$-invariance
condition~\eqref{eq:Sk-invariance}; together
with~(i) and~(ii), the constant-scalar reduction of
Proposition~\ref{prop:constant-scalar} gives
admissibility, $P(\mathcal{S}_2 u) = P(u)$, and
$\Lambda_2 = D_2 = I_m$.  The sums defining
$\Psi_p$ and $\Phi_p$ converge in $\ell^2(\Z)$: the
seed function $\psi_{t,r,n}$ has $\ell^2$ norm
independent of $r$,
\eqref{eq:CT-TASEP-phi-bound} gives
$\lVert\phi_{t,r,n}\rVert_{\ell^2} \leq C\beta^{-r}$,
and the propagator entries decay geometrically.
The seed recurrences
\eqref{eq:CT-psi-n}--\eqref{eq:CT-phi-n} give the
product-graph linear problems
\eqref{eq:seed-psi-product}--\eqref{eq:seed-phi-product}
with these edge weights, so
the identities \eqref{eq:CT-WF-psi-n}, \eqref{eq:CT-WF-phi-n} follow from
Lemma~\ref{lem:Psi-Phi-linear}.

\medskip
\noindent\textit{Continuous recurrences
\eqref{eq:CT-WF-psi-t}, \eqref{eq:CT-WF-phi-t}.}
We suppress the dependence on $t$ and $\mathbf{n}$,
writing $\Psi_p^{\mathbf{a}}$ for
$\Psi_p^{t,\mathbf{a},\mathbf{n}}$.
The propagator $B_{\mathbf{a},\mathbf{n}}$ is $t$-independent
(Lemma~\ref{lem:CT-TASEP-propagator}(iii)), so
$\pa_t$ passes through the propagator sums in
\eqref{eq:CT-Psi-def}.  For the $p$-th component,
applying the seed recurrence \eqref{eq:CT-psi-t}
termwise:
\begin{align*}
\pa_t \Psi_p^{\mathbf{a}}
= -\tfrac{1}{2}(\psi_{t,a_p,n_p}
   - \psi_{t,a_p-1,n_p})
- \tfrac{1}{2}\sum_{j > p}
  \sum_{w \leq a_j - 1}
  [B_{\mathbf{a}-\mathbf{1},\mathbf{n}}]_{jp}(w, a_p)
  (\psi_{t,w,n_j} - \psi_{t,w-1,n_j}).
\end{align*}
It remains to show that this equals
$-\frac{1}{2}(\Psi_p^{\mathbf{a}}
- \Psi_p^{\mathbf{a}-\mathbf{1}})$.  By
definition~\eqref{eq:CT-Psi-def} at
$\mathbf{a} - \mathbf{1}$,
\[
\Psi_p^{\mathbf{a}-\mathbf{1}}
= \psi_{t,a_p-1,n_p}
+ \sum_{j > p}
  \sum_{w \leq a_j - 2}
  [B_{\mathbf{a}-\mathbf{2},\mathbf{n}}]_{jp}(w, a_p{-}1)
  \psi_{t,w,n_j}.
\]
By $\mathcal{T}$-covariance
\eqref{eq:CT-Ba-covariance},
$[B_{\mathbf{a}-\mathbf{2},\mathbf{n}}]_{jp}(w, a_p{-}1)
= [B_{\mathbf{a}-\mathbf{1},\mathbf{n}}]_{jp}(w{+}1, a_p)$.
Substituting $w' = w + 1$ (with range
$w' \leq a_j - 1$) gives
\[
\Psi_p^{\mathbf{a}-\mathbf{1}}
= \psi_{t,a_p-1,n_p}
+ \sum_{j > p}
  \sum_{w' \leq a_j - 1}
  [B_{\mathbf{a}-\mathbf{1},\mathbf{n}}]_{jp}(w', a_p)
  \psi_{t,w'-1,n_j}.
\]
Comparing with \eqref{eq:CT-Psi-def}, the difference
$\Psi_p^{\mathbf{a}} - \Psi_p^{\mathbf{a}-\mathbf{1}}$
has leading term
$\psi_{t,a_p,n_p} - \psi_{t,a_p-1,n_p}$ and
propagator terms
$[B_{\mathbf{a}-\mathbf{1},\mathbf{n}}]_{jp}(w, a_p)
(\psi_{t,w,n_j} - \psi_{t,w-1,n_j})$,
which is exactly $-2\pa_t \Psi_p^{\mathbf{a}}$
as computed above.  This gives
\eqref{eq:CT-WF-psi-t}.

For the adjoint, $\pa_t$ again passes through
the propagator sums in \eqref{eq:CT-Phi-def}.
For the $p$-th component, applying
\eqref{eq:CT-phi-t} termwise:
\begin{align*}
\pa_t \Phi_p^{\mathbf{a}}
= -\tfrac{1}{2}(\phi_{t,a_p+1,n_p}
   - \phi_{t,a_p,n_p})
- \tfrac{1}{2}\sum_{j < p}
  \sum_{z \leq a_j}
  [B_{\mathbf{a},\mathbf{n}}]_{pj}(a_p, z)
  (\phi_{t,z+1,n_j} - \phi_{t,z,n_j}).
\end{align*}
It remains to show that this equals
$-\frac{1}{2}(\Phi_p^{\mathbf{a}+\mathbf{1}}
- \Phi_p^{\mathbf{a}})$.  By
definition~\eqref{eq:CT-Phi-def} at
$\mathbf{a} + \mathbf{1}$,
\[
\Phi_p^{\mathbf{a}+\mathbf{1}}
= \phi_{t,a_p+1,n_p}
+ \sum_{j < p}
  \sum_{z \leq a_j + 1}
  [B_{\mathbf{a}+\mathbf{1},\mathbf{n}}]_{pj}(a_p{+}1, z)
  \phi_{t,z,n_j}.
\]
By $\mathcal{T}$-covariance
\eqref{eq:CT-Ba-covariance},
$[B_{\mathbf{a}+\mathbf{1},\mathbf{n}}]_{pj}(a_p{+}1, z)
= [B_{\mathbf{a},\mathbf{n}}]_{pj}(a_p, z{-}1)$.
Substituting $z' = z - 1$ (with range
$z' \leq a_j$) gives
\[
\Phi_p^{\mathbf{a}+\mathbf{1}}
= \phi_{t,a_p+1,n_p}
+ \sum_{j < p}
  \sum_{z' \leq a_j}
  [B_{\mathbf{a},\mathbf{n}}]_{pj}(a_p, z')
  \phi_{t,z'+1,n_j}.
\]
Comparing with \eqref{eq:CT-Phi-def}, the difference
$\Phi_p^{\mathbf{a}+\mathbf{1}}
- \Phi_p^{\mathbf{a}}$ has leading term
$\phi_{t,a_p+1,n_p} - \phi_{t,a_p,n_p}$ and
propagator terms
$[B_{\mathbf{a},\mathbf{n}}]_{pj}(a_p, z)
(\phi_{t,z+1,n_j} - \phi_{t,z,n_j})$,
which is exactly $-2\pa_t \Phi_p^{\mathbf{a}}$
as computed above.  This gives
\eqref{eq:CT-WF-phi-t}.
\end{proof}

\medskip

Dressing compatibility is the remaining hypothesis of
the dual Darboux theorem
(Theorem~\ref{thm:sc-dual-Darboux}).  The following
lemma verifies it for the kernel
$K_{t,\mathbf{a},\mathbf{n}}$.

\begin{lemma}\label{lem:CT-TASEP-dressing}
The kernel $K_{t,\mathbf{a},\mathbf{n}}$ is dual dressing
compatible (Definition~\ref{def:sc-dual-dressing}) with
$(\Psi, \Phi)$ and $z = 1$ for every $t > 0$,
$\mathbf{a} \in \Z^m$, and
$1 \leq n_1 < \cdots < n_m$.
Explicitly:
\begin{align}
K_{t,\mathbf{a}+\mathbf{1},\mathbf{n}}
- K_{t,\mathbf{a},\mathbf{n}}
&= \Psi^{t,\mathbf{a}+\mathbf{1},\mathbf{n}}
   \Phi^{t,\mathbf{a}+\mathbf{1},\mathbf{n}},
\label{eq:CT-DC-a} \\
\pa_t K_{t,\mathbf{a},\mathbf{n}}
&= -\tfrac{1}{2}\Psi^{t,\mathbf{a},\mathbf{n}}
   \Phi^{t,\mathbf{a}+\mathbf{1},\mathbf{n}},
\label{eq:CT-DC-t} \\
K_{t,\mathbf{a},\mathbf{n}+\mathbf{1}}
- K_{t,\mathbf{a},\mathbf{n}}
&= 2\Psi^{t,\mathbf{a}+\mathbf{1},\mathbf{n}}
   \Phi^{t,\mathbf{a},\mathbf{n}+\mathbf{1}}.
\label{eq:CT-DC-n}
\end{align}
\end{lemma}

\begin{proof}
\noindent\textit{Identities \eqref{eq:CT-DC-a} and
\eqref{eq:CT-DC-n}.}
These are the $\mathcal{T}$- and
$\mathcal{S}_2$-factorization identities of
Lemma~\ref{lem:K-factorization}, applied to the
product graph $\mathcal{G}^m$ with distinguished shift
$\mathcal{T} = e^{\pa_{\mathbf{a}}}$ and
$\mathcal{S}_2 = e^{\pa_{\mathbf{n}}}$.  The
hypotheses are verified: the seed recurrences
\eqref{eq:CT-psi-n}--\eqref{eq:CT-phi-n} give the
product-graph linear problems
\eqref{eq:seed-psi-product}--\eqref{eq:seed-phi-product}
with $(C_2, D_2) = (2I_m, I_m)$; the
$\mathcal{T}$-covariance
\eqref{eq:CT-Ba-covariance},
$\mathcal{T}$-splitting \eqref{eq:CT-Ba-splitting},
and $\mathcal{S}_2$-compatibility (established in
Lemma~\ref{lem:CT-TASEP-propagator} and the
paragraph following it) provide the admissibility
conditions; the definitions
\eqref{eq:CT-Psi-def}--\eqref{eq:CT-Phi-def} match
the abstract constructions
\eqref{eq:Psi-def}--\eqref{eq:Phi-def}; and the
kernel $K_{t,\mathbf{a},\mathbf{n}}$ of
Lemma~\ref{lem:CT-TASEP-kernel} matches the abstract
kernel \eqref{eq:K-def}.  The same seed and propagator
estimates give trace-norm convergence of
$K_{t,\mathbf{a},\mathbf{n}}$, completing the regularity
required by Lemma~\ref{lem:K-factorization}.
For \eqref{eq:CT-DC-a},
the proof of Lemma~\ref{lem:K-factorization}
proceeds by splitting the diagonal and off-diagonal
parts of $K_{\mathbf{a}+\mathbf{1}} - K_{\mathbf{a}}$,
applying $\mathcal{T}$-splitting to the off-diagonal
difference, and factoring the result as
$\Psi(\mathcal{T}u)\Phi(\mathcal{T}u)$.  For
\eqref{eq:CT-DC-n}, the same lemma gives
$K(\mathcal{S}_2 u) - K(u)
= \Psi(\mathcal{T}u)C_2(u)\Phi(\mathcal{S}_2 u)$;
since $C_2 = 2I_m$, this is
$2\Psi^{t,\mathbf{a}+\mathbf{1},\mathbf{n}}
\Phi^{t,\mathbf{a},\mathbf{n}+\mathbf{1}}$.

It remains to verify \eqref{eq:CT-DC-t}.  We suppress
$t$ and $\mathbf{n}$ dependence, writing
$\Psi_p^{\mathbf{a}}$ for
$\Psi_p^{t,\mathbf{a},\mathbf{n}}$.
Since $B_{\mathbf{a},\mathbf{n}}$ is $t$-independent
(Lemma~\ref{lem:CT-TASEP-propagator}(iii)),
$\pa_t$ acts only on the seed functions
$\psi, \phi$ in the kernel formula.  The derivatives
$\pa_t\psi$ and $\pa_t\phi$ satisfy the same
bounds as $\psi$ and $\phi$ (differentiation multiplies
each contour integrand by the bounded factor
$(w - 1/2)$), so the differentiated series converges
in trace norm; in particular,
$t \mapsto K_{t,\mathbf{a},\mathbf{n}}$ is $C^1$ in
trace norm.  We split
$K = K^{\mathrm{diag}} + K^{\mathrm{off}}$ and
handle each part separately.

\medskip
\noindent\textit{Step~1: Diagonal part.}
The diagonal contribution is
\[
K^{\mathrm{diag}} = \sum_{p=1}^m \sum_{r \leq a_p}
\psi_{t,r,n_p}\phi_{t,r,n_p}.
\]
Applying the seed recurrences
\eqref{eq:CT-psi-t}--\eqref{eq:CT-phi-t}:
\begin{align*}
\pa_t K^{\mathrm{diag}}
&= -\tfrac{1}{2}\sum_{p=1}^m \sum_{r \leq a_p}
\bigl[(\psi_{t,r,n_p} - \psi_{t,r-1,n_p})
\phi_{t,r,n_p}
+ \psi_{t,r,n_p}
(\phi_{t,r+1,n_p} - \phi_{t,r,n_p})\bigr] \\
&= -\tfrac{1}{2}\sum_{p=1}^m \sum_{r \leq a_p}
\bigl[\psi_{t,r,n_p}\phi_{t,r+1,n_p}
- \psi_{t,r-1,n_p}\phi_{t,r,n_p}\bigr].
\end{align*}
Since $\psi_{t,r,n_p}(x) = 0$ for $r < x - n_p$, the
inner sum over $r$ telescopes to the boundary term at
$r = a_p$:
\begin{equation}\label{eq:CT-DC-t-diag}
\pa_t K^{\mathrm{diag}}
= -\tfrac{1}{2}\sum_{p=1}^m
\psi_{t,a_p,n_p}\phi_{t,a_p+1,n_p}.
\end{equation}

\noindent\textit{Step~2: Off-diagonal part.}
The off-diagonal contribution is
\[
K^{\mathrm{off}} = \sum_{i > j}
\sum_{r \leq a_i}\sum_{r' \leq a_j}
\psi_{t,r,n_i}
[B_{\mathbf{a},\mathbf{n}}]_{ij}(r, r')
\phi_{t,r',n_j}.
\]
Applying the seed recurrences
\eqref{eq:CT-psi-t} and \eqref{eq:CT-phi-t} to
$\psi_{t,r,n_i}$ and $\phi_{t,r',n_j}$ respectively:
\begin{align*}
\pa_t K^{\mathrm{off}}
= -\tfrac{1}{2}\sum_{i > j}
\sum_{r \leq a_i}\sum_{r' \leq a_j}
\bigl[&\psi_{t,r,n_i}
[B_{\mathbf{a},\mathbf{n}}]_{ij}(r, r')
\phi_{t,r'+1,n_j} \\
&- \psi_{t,r-1,n_i}
[B_{\mathbf{a},\mathbf{n}}]_{ij}(r, r')
\phi_{t,r',n_j}\bigr].
\end{align*}
The second term is reindexed using
$\mathcal{T}$-covariance
\eqref{eq:CT-Ba-covariance}:
$[B_{\mathbf{a},\mathbf{n}}]_{ij}(r, r')
= [B_{\mathbf{a}-\mathbf{1},\mathbf{n}}]_{ij}
(r{-}1, r'{-}1)$; substituting
$r \to r{+}1$, $r' \to r'{+}1$ shifts the
summation range to
$r \leq a_i{-}1$, $r' \leq a_j{-}1$:
\begin{align}
\pa_t K^{\mathrm{off}}
&= -\tfrac{1}{2}\sum_{i > j}
\sum_{r \leq a_i}\sum_{r' \leq a_j}
\psi_{t,r,n_i}
[B_{\mathbf{a},\mathbf{n}}]_{ij}(r, r')
\phi_{t,r'+1,n_j}
\notag \\
&\quad +\tfrac{1}{2}\sum_{i > j}
\sum_{r \leq a_i - 1}\sum_{r' \leq a_j - 1}
\psi_{t,r,n_i}
[B_{\mathbf{a}-\mathbf{1},\mathbf{n}}]_{ij}(r, r')
\phi_{t,r'+1,n_j}.
\label{eq:CT-DC-t-split}
\end{align}
The $\mathcal{T}$-splitting identity
\eqref{eq:CT-Ba-splitting} at
$\mathbf{a}{-}\mathbf{1}$ gives
\[
[B_{\mathbf{a}-\mathbf{1},\mathbf{n}}]_{ij}(r, r')
= [B_{\mathbf{a},\mathbf{n}}]_{ij}(r, r')
- \sum_{j < \mu < i}
[B_{\mathbf{a}-\mathbf{1},\mathbf{n}}]_{i\mu}(r, a_\mu)
[B_{\mathbf{a},\mathbf{n}}]_{\mu j}(a_\mu, r').
\]
Substituting this into the second sum of
\eqref{eq:CT-DC-t-split} splits it into a
$[B_{\mathbf{a},\mathbf{n}}]_{ij}$-weighted sum and a
cross-term.  The
$[B_{\mathbf{a},\mathbf{n}}]_{ij}$ sum has range
$r \leq a_i{-}1$, $r' \leq a_j{-}1$, while the
first sum of \eqref{eq:CT-DC-t-split} has range
$r \leq a_i$, $r' \leq a_j$; the two share
the same summand, so their combination reduces to
the boundary regions
$\{r = a_i, r' \leq a_j\}$ and
$\{r \leq a_i{-}1, r' = a_j\}$.
The off-diagonal derivative therefore separates into
boundary contributions~(A), (B) and
cross-term~(C):
\begin{align}
\pa_t K^{\mathrm{off}}
&= -\tfrac{1}{2}\sum_{i > j}
\sum_{r' \leq a_j}
\psi_{t,a_i,n_i}
[B_{\mathbf{a},\mathbf{n}}]_{ij}(a_i, r')
\phi_{t,r'+1,n_j}
\tag{A} \\
&\quad -\tfrac{1}{2}\sum_{i > j}
\sum_{r \leq a_i - 1}
\psi_{t,r,n_i}
[B_{\mathbf{a},\mathbf{n}}]_{ij}(r, a_j)
\phi_{t,a_j+1,n_j}
\tag{B} \\
&\quad -\tfrac{1}{2}\sum_{i > j}
\sum_{j < \mu < i}
\sum_{r \leq a_i - 1}\sum_{r' \leq a_j - 1}
\psi_{t,r,n_i}
[B_{\mathbf{a}-\mathbf{1},\mathbf{n}}]_{i\mu}(r, a_\mu)
[B_{\mathbf{a},\mathbf{n}}]_{\mu j}(a_\mu, r')
\phi_{t,r'+1,n_j}.
\tag{C}
\end{align}

\noindent\textit{Step~3: Term matching.}
The right-hand side of \eqref{eq:CT-DC-t} is
$-\frac{1}{2}\Psi^{t,\mathbf{a},\mathbf{n}}
\Phi^{t,\mathbf{a}+\mathbf{1},\mathbf{n}}
= -\frac{1}{2}\sum_{p=1}^m
\Psi_p^{t,\mathbf{a},\mathbf{n}}
\Phi_p^{t,\mathbf{a}+\mathbf{1},\mathbf{n}}$.
By $\mathcal{T}$-covariance of the propagator,
$\Phi_p^{\mathbf{a}+\mathbf{1}}
= \phi_{t,a_p+1,n_p}
+ \sum_{j < p}\sum_{z \leq a_j}
[B_{\mathbf{a},\mathbf{n}}]_{pj}(a_p, z)
\phi_{t,z+1,n_j}$.
Expanding the product
$\Psi_p^{\mathbf{a}}
\Phi_p^{\mathbf{a}+\mathbf{1}}$
using \eqref{eq:CT-Psi-def}--\eqref{eq:CT-Phi-def},
summing over $p$, and relabelling indices produces
four terms:
\begin{align}
& \sum_{p=1}^m
\psi_{t,a_p,n_p}\phi_{t,a_p+1,n_p},
\label{eq:CT-target-diag} \\
& \sum_{i > j}\sum_{r' \leq a_j}
\psi_{t,a_i,n_i}
[B_{\mathbf{a},\mathbf{n}}]_{ij}(a_i, r')
\phi_{t,r'+1,n_j},
\label{eq:CT-target-ba} \\
& \sum_{i > j}\sum_{r \leq a_i - 1}
\psi_{t,r,n_i}
[B_{\mathbf{a}-\mathbf{1},\mathbf{n}}]_{ij}(r, a_j)
\phi_{t,a_j+1,n_j},
\label{eq:CT-target-ab} \\
& \sum_{i > \mu > j}
\sum_{r \leq a_i - 1}\sum_{r' \leq a_j}
\psi_{t,r,n_i}
[B_{\mathbf{a}-\mathbf{1},\mathbf{n}}]_{i\mu}(r, a_\mu)
[B_{\mathbf{a},\mathbf{n}}]_{\mu j}(a_\mu, r')
\phi_{t,r'+1,n_j}.
\label{eq:CT-target-bb}
\end{align}
Comparing with Step~1 and Step~2, the diagonal
\eqref{eq:CT-DC-t-diag} matches
$-\frac{1}{2}$\eqref{eq:CT-target-diag},
and~(A) matches
$-\frac{1}{2}$\eqref{eq:CT-target-ba}.  It remains
to verify (B) + (C)
$= -\frac{1}{2}$[\eqref{eq:CT-target-ab}
+ \eqref{eq:CT-target-bb}].

The difference (B) $- (-\frac{1}{2})$\eqref{eq:CT-target-ab} is
\[
-\tfrac{1}{2}\sum_{i > j}
\sum_{r \leq a_i - 1}
\psi_{t,r,n_i}
\bigl([B_{\mathbf{a},\mathbf{n}}]_{ij}
- [B_{\mathbf{a}-\mathbf{1},\mathbf{n}}]_{ij}\bigr)
(r, a_j)\phi_{t,a_j+1,n_j}.
\]
By $\mathcal{T}$-splitting, the bracket equals
$\sum_{j < \mu < i}
[B_{\mathbf{a}-\mathbf{1},\mathbf{n}}]_{i\mu}(r, a_\mu)
[B_{\mathbf{a},\mathbf{n}}]_{\mu j}(a_\mu, a_j)$, so
\begin{equation}\label{eq:CT-DC-t-remainder}
\text{(B)} + \tfrac{1}{2}\eqref{eq:CT-target-ab}
= -\tfrac{1}{2}\sum_{i > \mu > j}
\sum_{r \leq a_i - 1}
\psi_{t,r,n_i}
[B_{\mathbf{a}-\mathbf{1},\mathbf{n}}]_{i\mu}(r, a_\mu)
[B_{\mathbf{a},\mathbf{n}}]_{\mu j}(a_\mu, a_j)
\phi_{t,a_j+1,n_j}.
\end{equation}
The difference (C) $- (-\frac{1}{2})$\eqref{eq:CT-target-bb}
reduces to the single term $r' = a_j$, since
\eqref{eq:CT-target-bb} sums over
$r' \leq a_j$ while (C) sums over
$r' \leq a_j{-}1$:
\begin{equation}\label{eq:CT-DC-t-excess}
\text{(C)} + \tfrac{1}{2}\eqref{eq:CT-target-bb}
= +\tfrac{1}{2}\sum_{i > \mu > j}
\sum_{r \leq a_i - 1}
\psi_{t,r,n_i}
[B_{\mathbf{a}-\mathbf{1},\mathbf{n}}]_{i\mu}(r, a_\mu)
[B_{\mathbf{a},\mathbf{n}}]_{\mu j}(a_\mu, a_j)
\phi_{t,a_j+1,n_j}.
\end{equation}
Since \eqref{eq:CT-DC-t-remainder}
$= -$\eqref{eq:CT-DC-t-excess},
(B) + (C) $= -\frac{1}{2}$[\eqref{eq:CT-target-ab}
+ \eqref{eq:CT-target-bb}],
completing the verification of \eqref{eq:CT-DC-t}.

It remains to verify resolvent existence.  By
Lemma~\ref{lem:CT-TASEP-kernel},
$\det(I - K_{t,\mathbf{a},\mathbf{n}})$ equals the gap
probability
$\PP_y\bigl(\bigcap_{i=1}^m \{Y_{n_i}(t) > a_i\}\bigr)$;
since $K_{t,\mathbf{a},\mathbf{n}}$ is trace class, the
resolvent exists if and only if this determinant is
nonzero.  This probability is strictly positive for
$t > 0$: choose $J$ with $y_{n_i} + J > a_i$ for
each~$i$, and consider the event that particles
$1, 2, \ldots, n_m$ each ring exactly $J$ times before
time~$t$, with all of particle~$k{-}1$'s rings preceding
those of particle~$k$.  The strict integer ordering of the
initial data ensures that no jump is blocked, so
$Y_{n_i}(t) = y_{n_i} + J > a_i$ on this event.  This
event has positive probability, so the determinant is
nonzero.
\end{proof}

\subsection{Multipoint equation}
\label{sec:CT-TASEP-multipoint}

The dressed observable
\[
\mathcal{M}_{t,\mathbf{a},\mathbf{n}}
\defeq I + \Phi^{t,\mathbf{a},\mathbf{n}}
(I - K_{t,\mathbf{a},\mathbf{n}})^{-1}
\Psi^{t,\mathbf{a},\mathbf{n}}
\in \End(\R^m)
\]
is invertible for every $t > 0$,
$\mathbf{a} \in \Z^m$, and
$1 \leq n_1 < \cdots < n_m$, with
dressed edge weights
\begin{align}
\mathcal{M}_1(u)
&= -\tfrac{1}{2}I_m
   + \mathcal{M}(Tu)^{-1}
   \pa_t\mathcal{M}(Tu),
\label{eq:CT-dressed-C1} \\
\mathcal{M}_2(u)
&= 2\mathcal{M}(Tu)^{-1}\mathcal{M}(S_2 u).
\label{eq:CT-dressed-C2}
\end{align}

\begin{theorem}\label{thm:CT-TASEP-multipoint}
The dressed edge weights satisfy the mixed diamond
equation
\begin{equation}\label{eq:CT-TASEP-multipoint}
\begin{aligned}
\mathcal{M}_{t,\mathbf{a}+\mathbf{1},
\mathbf{n}+\mathbf{1}}^{-1}
\pa_t \mathcal{M}_{t,\mathbf{a}+\mathbf{1},
\mathbf{n}+\mathbf{1}}
&- \mathcal{M}_{t,\mathbf{a}+\mathbf{1},
\mathbf{n}}^{-1}
\pa_t \mathcal{M}_{t,\mathbf{a}+\mathbf{1},
\mathbf{n}} \\
\qquad + \mathcal{M}_{t,\mathbf{a}+\mathbf{2},
\mathbf{n}}^{-1}
\mathcal{M}_{t,\mathbf{a}+\mathbf{1},
\mathbf{n}+\mathbf{1}}
&- \mathcal{M}_{t,\mathbf{a}+\mathbf{1},
\mathbf{n}}^{-1}
\mathcal{M}_{t,\mathbf{a},\mathbf{n}+\mathbf{1}} = 0
\end{aligned}
\end{equation}
for every $t > 0$, $\mathbf{a} \in \Z^m$, and
$1 \leq n_1 < \cdots < n_m$.
\end{theorem}

\begin{proof}
By Theorem~\ref{thm:sc-dual-Darboux}, the pair
$(\mathcal{M}_k, \Lambda_k)$ satisfies the dual diamond
equations.  The mixed diamond equation
\eqref{eq:sc-dual-diamond-mixed} for the pair $(1, 2)$
reads
\begin{equation}\label{eq:CT-dressed-mixed}
\pa_t\Lambda_2 + \mathcal{M}_1(u)\Lambda_2
+ \Lambda_1\mathcal{M}_2(Tu)
- \mathcal{M}_2(u)\Lambda_1
- \Lambda_2\mathcal{M}_1(S_2 u) = 0.
\end{equation}
Since $\Lambda_1 = -\frac{1}{2}I_m$ and $\Lambda_2 = I_m$
are constant scalar matrices, $\pa_t\Lambda_2 = 0$,
$\Lambda_2$ acts as the identity, and the two
$\Lambda_1$-terms commute with $\mathcal{M}_2$ and combine
into the half-difference below.  The equation reduces to
\begin{equation}\label{eq:CT-reduced-mixed}
\mathcal{M}_1(u) - \mathcal{M}_1(S_2 u)
+ \tfrac{1}{2}\bigl[\mathcal{M}_2(u)
- \mathcal{M}_2(Tu)\bigr] = 0.
\end{equation}
Substituting \eqref{eq:CT-dressed-C1}--\eqref{eq:CT-dressed-C2}
into \eqref{eq:CT-reduced-mixed}: the constant
$-\frac{1}{2}I_m$ in $\mathcal{M}_1(u)$ and
$\mathcal{M}_1(S_2 u)$ cancel, and the factor
$\frac{1}{2} \cdot 2 = 1$ in the $\mathcal{M}_2$ terms
simplifies.  The result is
\begin{align*}
\mathcal{M}(Tu)^{-1}\pa_t\mathcal{M}(Tu)
&- \mathcal{M}(TS_2 u)^{-1}
   \pa_t\mathcal{M}(TS_2 u) \\
&+ \mathcal{M}(Tu)^{-1}\mathcal{M}(S_2 u)
 - \mathcal{M}(T^2 u)^{-1}\mathcal{M}(TS_2 u) = 0.
\end{align*}
Evaluating at $u = (t, \mathbf{a}, \mathbf{n})$ and using
$Tu = (t, \mathbf{a}{+}\mathbf{1}, \mathbf{n})$,
$S_2 u = (t, \mathbf{a}, \mathbf{n}{+}\mathbf{1})$,
$TS_2 u = (t, \mathbf{a}{+}\mathbf{1},
\mathbf{n}{+}\mathbf{1})$,
$T^2 u = (t, \mathbf{a}{+}\mathbf{2}, \mathbf{n})$:
\begin{align*}
\mathcal{M}_{t,\mathbf{a}+\mathbf{1},\mathbf{n}}^{-1}
\pa_t\mathcal{M}_{t,\mathbf{a}+\mathbf{1},\mathbf{n}}
&- \mathcal{M}_{t,\mathbf{a}+\mathbf{1},
\mathbf{n}+\mathbf{1}}^{-1}
\pa_t\mathcal{M}_{t,\mathbf{a}+\mathbf{1},
\mathbf{n}+\mathbf{1}} \\
&+ \mathcal{M}_{t,\mathbf{a}+\mathbf{1},
\mathbf{n}}^{-1}
\mathcal{M}_{t,\mathbf{a},\mathbf{n}+\mathbf{1}}
- \mathcal{M}_{t,\mathbf{a}+\mathbf{2},
\mathbf{n}}^{-1}
\mathcal{M}_{t,\mathbf{a}+\mathbf{1},
\mathbf{n}+\mathbf{1}} = 0.
\end{align*}
Multiplying by $-1$ gives
\eqref{eq:CT-TASEP-multipoint}.
\end{proof}

For $m = 1$, the dressed observable is scalar.

\begin{corollary}\label{cor:CT-TASEP-scalar}
The one-point Fredholm determinant
$F_{t, a, n} \defeq \det(I - K_{t, a, n})$ satisfies
\begin{equation}\label{eq:CT-TASEP-scalar-HM}
F_{t,a,n}\pa_t F_{t,a,n+1}
- F_{t,a,n+1}\pa_t F_{t,a,n}
+ F_{t,a,n}F_{t,a,n+1}
- F_{t,a+1,n}F_{t,a-1,n+1} = 0.
\end{equation}
\end{corollary}

\begin{proof}
Corollary~\ref{cor:sc-dual-scalar} applies with the data
$(c_2, \lambda_1, \lambda_2) = (2, -\tfrac{1}{2}, 1)$
and shifts $T = e^{\pa_a}$,
$S_2 = e^{\pa_n}$.  The Fredholm determinant
$F_{t,a,n}$ is $C^1$ in $t$ by the trace-norm
differentiability established in
Lemma~\ref{lem:CT-TASEP-dressing}.  For $m = 1$ the
propagator vanishes, so
$K_{t,a,n} = \sum_{r \leq a}
\psi_{t,r,n}\phi_{t,r,n}$.  The estimates from
Lemma~\ref{lem:CT-TASEP-kernel} give
$\lVert K_{t,a,n}\rVert_1 \leq C\sum_{r \leq a}
\beta^{-r} \to 0$ as $a \to -\infty$; differentiating
in $t$ multiplies the contour integrands by the bounded
factor $(w - 1/2)$, so
$\lVert\pa_t K_{t,a,n}\rVert_1 \to 0$ likewise.
Continuity of the Fredholm determinant in trace norm
gives $F_{t,a,n} \to 1$.  For all sufficiently negative
$a$ the resolvent $(I - K_{t,a,n})^{-1}$ is uniformly
bounded, and Jacobi's formula gives
$\pa_t F_{t,a,n} \to 0$.
Equation \eqref{eq:sc-dual-HM} gives
\begin{align*}
c_2\lambda_1\bigl[F_{t,a+1,n}F_{t,a-1,n+1}
- F_{t,a,n}F_{t,a,n+1}\bigr]
+ \lambda_2\bigl[F_{t,a,n}\pa_t F_{t,a,n+1}
- F_{t,a,n+1}\pa_t F_{t,a,n}\bigr] = 0.
\end{align*}
Substituting $c_2\lambda_1 = -1$ and
$\lambda_2 = 1$ and rearranging gives
\eqref{eq:CT-TASEP-scalar-HM}.
\end{proof}

\begin{remark}[Structure of the multipoint equation]
\label{rem:CT-TASEP-structure}
Equation \eqref{eq:CT-TASEP-multipoint} is an
$m \times m$ matrix differential-difference equation
coupling the dressed observable $\mathcal{M}$ at four
lattice points, together with the time derivative
$\pa_t$.  For $m \geq 2$, the equation is
intrinsically noncommutative: the matrix inverses and
products do not simplify to a scalar relation.
\end{remark}

{
  \setlength{\parskip}{0pt}
}

\section{Push-TASEP}\label{sec:Push-TASEP}

\paragraph{\textbf{System description.}}
The continuous-time Push-TASEP is an interacting particle
system on $\Z$ with at most one particle per site.  Given
a right-finite, strictly decreasing initial configuration
$Y(0) = y = (y_1 > y_2 > y_3 > \cdots)$,
the dynamics run in continuous time: each particle carries
an independent rate-one exponential clock, and when
particle~$n$'s clock rings it jumps one step to the left,
with any left neighbours also pushed leftward so as to
preserve exclusion.  After each particle's
jump, its independent clock is instantaneously reset.

\paragraph{\textbf{Fredholm determinant formula.}}
Fix a time $t \geq 0$, an initial configuration $y$,
finitely many particle labels
$\mathbf{n} = (n_1, \dotsc, n_m)$ with
$1 \leq n_1 < n_2 < \cdots < n_m$, and spatial
thresholds $\mathbf{a} = (a_1, \dotsc, a_m) \in \Z^m$.
The multipoint joint cumulative distribution of the particle
positions is given by a Fredholm determinant on
$\ell^2(\{n_1, \dotsc, n_m\} \times \Z)$:
\[
\PP_y\Bigl(\bigcap_{i=1}^m
\{Y_{n_i}(t) > a_i\}\Bigr)
= \det(I - \bar{\chi}_{\mathbf{a}} K_t
\bar{\chi}_{\mathbf{a}})_{\ell^2(\{n_1, \dotsc, n_m\}
\times \Z)},
\]
where $\bar{\chi}_{\mathbf{a}}(n_i, x) =
\mathbf{1}_{x \leq a_i}$.  The correlation kernel $K_t$ has
block entries
\[
K_t(n_i, x_i; n_j, x_j)
= -Q^{n_j - n_i}(x_i, x_j)\mathbf{1}_{n_i < n_j}
+ \bigl(\mathcal{S}_{-t,-n_i}^*
\overline{\mathcal{S}}_{-t,n_j}^{\operatorname{epi}(y)}
\bigr)(x_i, x_j),
\]
where $A^*$ denotes the transpose kernel,
$A^*(x, y) = A(y, x)$.  The operators defining the kernel
are as follows.

\medskip
\noindent\textit{The transition matrix $Q$.}
The transition matrix $Q$ is the same geometric left-jump
kernel as in the continuous-time TASEP
(\S\ref{sec:CT-TASEP}):
\[
Q(z_1, z_2) = 2^{-(z_1 - z_2)}\mathbf{1}_{z_1 > z_2},
\]
with its $n$-th power given by
\[
Q^n(z_1, z_2) = \binom{z_1 - z_2 - 1}{n - 1}
2^{-(z_1 - z_2)}\mathbf{1}_{z_1 - z_2 \geq n}.
\]

\medskip
\noindent\textit{The scattering operators
$\mathcal{S}_{-t,-n}$ and
$\overline{\mathcal{S}}_{-t,n}$.}
These operators are defined by the contour integrals
\begin{equation}\label{eq:Push-psi-contour}
\mathcal{S}_{-t,-n}(z_1, z_2)
= \frac{1}{2\pi\mathrm{i}}
\oint_{\gamma_\rho} \frac{\diff w}{w}
\frac{(1-w)^n}{2^{z_2 - z_1} w^{n + z_2 - z_1}}
e^{t(1/w - 2)},
\end{equation}
\begin{equation}\label{eq:Push-Sbar-contour}
\overline{\mathcal{S}}_{-t,n}(z_1, z_2)
= \frac{1}{2\pi\mathrm{i}}
\oint_{\gamma_\delta} \frac{\diff w}{w}
\frac{(1-w)^{z_2 - z_1 + n - 1}}{2^{z_1 - z_2} w^{n-1}}
e^{t(2 - 1/(1-w))},
\end{equation}
where $\gamma_\rho$ is a positively oriented circle of
radius $\rho \in (0, 1)$ around the origin, and
$\gamma_\delta$ is a sufficiently small positively
oriented circle around the origin, chosen so that the only
singularity enclosed is the pole at $w = 0$.

\medskip
\noindent\textit{The epigraph operator
$\overline{\mathcal{S}}_{-t,n}^{\operatorname{epi}(y)}$.}
This is the hitting probability operator defined in terms
of the initial data:
\[
\overline{\mathcal{S}}_{-t,n}^{\operatorname{epi}(y)}
(z_1, z_2)
= \E_{W_0 = z_1}\bigl[
\overline{\mathcal{S}}_{-t,n-\tau}(W_\tau, z_2)\mathbf{1}_{\tau < n}\bigr],
\]
where $(W_\ell)_{\ell \geq 0}$ is the random walk with
transition matrix $Q$, and
$\tau = \min\{\ell \geq 0 :
W_\ell > y_{\ell+1}\}$
is the hitting time of the strict epigraph of the initial
data by the random walk, with the convention
$\tau = \infty$ if the set is empty.

\paragraph{\textbf{References.}}
The Fredholm determinant formula above is the
Nica--Quastel--Remenik PushASEP formula
(Theorem~4.1 of~\cite{NicaQuastelRemenik2020}, specialized to
$r = 0$, $\ell = 1$).

\subsection{Kernel reformulation}\label{sec:Push-TASEP-kernel-reform}

The extended-kernel Fredholm determinant reduces to a
single-space Fredholm determinant on $\ell^2(\Z)$.  Define
\[
\psi_{t,r,n}(x) \defeq \mathcal{S}_{-t,-n}(x, r),
\qquad
\phi_{t,r,n}(x) \defeq
\overline{\mathcal{S}}_{-t,n}^{\operatorname{epi}(y)}(x, r),
\]
where we suppress the dependence of $\phi$ on the initial
condition $y$.  The propagator is an instance of the
directed-path propagator
(Definition~\ref{def:directed-path-propagator}) with
transition kernels
$Q_{\ell,\ell'} = Q^{n_{\ell'} - n_\ell}$ for
$\ell < \ell'$.  Explicitly, for $i > j$,
\begin{equation}\label{eq:Push-TASEP-propagator}
[B_{\mathbf{a},\mathbf{n}}]_{ij}(r, r')
\defeq \sum_{k=1}^{i-j}
  \sum_{j = \ell_0 < \ell_1 < \cdots < \ell_k = i}
  \sum_{\substack{\xi_1 \leq a_{\ell_1} \\ \vdots \\
  \xi_{k-1} \leq a_{\ell_{k-1}}}}
  \prod_{s=0}^{k-1}
  \Bigl(-Q^{n_{\ell_{s+1}} - n_{\ell_s}}
  (\xi_s, \xi_{s+1})\Bigr),
\end{equation}
with $\xi_0 = r'$ and $\xi_k = r$,
and $[B_{\mathbf{a},\mathbf{n}}]_{ij} = 0$ for $i \leq j$.
The internal vertices $\xi_1, \dotsc, \xi_{k-1}$ are
constrained by the cutoffs, but the endpoints $r, r'$ are
arbitrary.  The propagator is identical to that of the
continuous-time TASEP (\S\ref{sec:CT-TASEP}), since it
depends only on the transition matrix $Q$ and not on the
scattering operators.

\begin{lemma}\label{lem:Push-TASEP-kernel}
We have
\[
\PP_y\Bigl(\bigcap_{i=1}^m
\{Y_{n_i}(t) > a_i\}\Bigr)
= \det(I - K_{t, \mathbf{a}, \mathbf{n}})_{\ell^2(\Z)},
\]
where
\[
K_{t, \mathbf{a}, \mathbf{n}}(x, x')
= \sum_{1 \leq j \leq i \leq m}
  \sum_{r \leq a_i} \sum_{r' \leq a_j}
  \psi_{t, r, n_i}(x)
  \bigl[\delta_{ij}\delta_{r,r'}
  + [B_{\mathbf{a},\mathbf{n}}]_{ij}(r, r')\bigr]
  \phi_{t, r', n_j}(x').
\]
\end{lemma}

\begin{proof}
The argument is identical to that of
Lemma~\ref{lem:CT-TASEP-kernel}: the conjugation by
exponential weights, block-triangular decomposition,
Sylvester's identity, and propagator identification
produce the claimed single-space kernel.
The random walk $Q$ and the threshold constraints are
the same.  The temporal factor
$e^{t(2 - 1/(1-w))}$ in the $\overline{\mathcal{S}}$
contour formula \eqref{eq:Push-Sbar-contour} is bounded
on $\gamma_\delta$, and the spatial factor is the same as
in the continuous-time TASEP.  For any
$\beta \in (1/2, 1)$, the estimate
$|\overline{\mathcal{S}}_{-t,k}(z, r)|
\leq C\beta^{z-r}$ for $z \geq r$ and the $\phi$ bound
\begin{equation}\label{eq:Push-TASEP-phi-bound}
|\phi_{t,r,n}(x)| \leq C\beta^{x-r}\mathbf{1}_{x > y_{n_m}}
\end{equation}
therefore carry over.

The seed estimate for $\psi$ requires modification,
because the exponential factor $e^{t(1/w - 2)}$
introduces an essential singularity at $w = 0$, whereas
the continuous-time TASEP factor $e^{t(w - 1/2)}$ is
entire.  The contour formula
\eqref{eq:Push-psi-contour} gives
\[
\psi_{t,r,n}(x) = \frac{1}{2\pi\mathrm{i}}
\oint_{\gamma_\rho} \frac{\diff w}{w}
\frac{(1-w)^n}{2^{r-x} w^{n+r-x}}
e^{t(1/w - 2)}.
\]
The integrand depends on $r$ and $x$ only through
$r - x$, so $\psi_{t,r,n}(x) = g(r - x)$ for a fixed
function $g$.  The residue computation gives
\[
\psi_{t,r,n}(x)
= e^{-2t}2^{x-r}
\sum_{\substack{0 \leq j \leq n \\ j \geq n + r - x}}
(-1)^j\binom{n}{j}
\frac{t^{j-n-r+x}}{(j-n-r+x)!}.
\]
The summation condition $j \geq n + r - x$ with $j \leq n$
requires $r \leq x$, so $\psi_{t,r,n}(x) = 0$ for $r > x$.
For the lower tail, set $q = x - r$; when $q \geq n$, the
sum has $n + 1$ terms, giving
$|\psi_{t,r,n}(x)| \leq C(2t)^{q-n}/(q-n)!$.
This is super-exponential decay; in particular
$g \in \ell^2(\Z)$, and
$\lVert\psi_{t,r,n}\rVert_{\ell^2(\Z)}
= \lVert g\rVert_{\ell^2}$ is independent of $r$.

With these bounds, the same conjugation and Sylvester
reduction as in Lemma~\ref{lem:CT-TASEP-kernel} apply.
\end{proof}

For the remainder of this section, we work with the kernel
$K_{t, \mathbf{a}, \mathbf{n}}$ on $\ell^2(\Z)$, with
ordered labels $1 \leq n_1 < \cdots < n_m$ as in the
Fredholm formula.  Write
$\mathbf{1} = (1, \ldots, 1) \in \Z^m$.

\subsection{Recurrences and kernel identities}
\label{sec:Push-TASEP-recurrences}

Continuous-time Push-TASEP naturally sits in the standard
presentation of the semi-discrete framework
(\S\ref{sec:sc-linear}): the seed recurrences take the form
$\pa_t\psi(u)$ rather than $\pa_t\psi(Tu)$, with
a discrete particle-label shift in the $S_2$ direction.
The product graph architecture of \S\ref{sec:product-graph}
provides the computational scaffolding: directed-path
propagator $B_{\mathbf{a},\mathbf{n}}$, distinguished shift
$e^{\pa_a}$, and product-graph wave functions
$\Psi, \Phi$, with $\pa_t$ replacing one discrete
shift.  The framework identification proceeds in three
stages: the seed functions satisfy the standard linear
problem (Lemma~\ref{lem:Push-TASEP-seed-linear}), the
product-graph wave functions inherit this structure
(Lemma~\ref{lem:Push-TASEP-WF}), and the multipoint
equation follows from the semi-discrete Darboux theorem
(Theorem~\ref{thm:sc-Darboux}).

\begin{lemma}\label{lem:Push-TASEP-seed-linear}
The seed functions $\psi_{t,r,n}$ and $\phi_{t,r,n}$
defined in \S\ref{sec:Push-TASEP-kernel-reform} satisfy the
semi-discrete linear problem
\eqref{eq:sc-linear-1}--\eqref{eq:sc-linear-j}
and its adjoint
\eqref{eq:sc-adjoint-1}--\eqref{eq:sc-adjoint-j}
with shifts
$T = e^{\pa_r}$, $\pa_1 = \pa_t$,
$S_2 = e^{\pa_n}$ and constant scalar edge weights
$(c_1, \lambda_1) = (2, 2)$,
$(c_2, \lambda_2) = (2, 1)$.
Explicitly:
\begin{align}
\pa_t \psi_{t,r,n}
&= 2\bigl(\psi_{t,r+1,n}
   - \psi_{t,r,n}\bigr),
\label{eq:Push-psi-t} \\
\psi_{t,r,n+1}
&= 2\psi_{t,r+1,n} - \psi_{t,r,n},
\label{eq:Push-psi-n}
\end{align}
and
\begin{align}
\pa_t \phi_{t,r+1,n}
&= 2\bigl(\phi_{t,r+1,n}
   - \phi_{t,r,n}\bigr),
\label{eq:Push-phi-t} \\
\phi_{t,r+1,n}
&= 2\phi_{t,r,n+1} - \phi_{t,r+1,n+1}.
\label{eq:Push-phi-n}
\end{align}
\end{lemma}

\begin{proof}
\noindent\textit{Verification of \eqref{eq:Push-psi-t}.}
From the contour representation,
\[
\psi_{t,r,n}(x) = \frac{1}{2\pi\mathrm{i}}
\oint_{\gamma_\rho} \frac{\diff w}{w}
\frac{(1-w)^n}{2^{r-x} w^{n+r-x}}
e^{t(1/w - 2)}.
\]
Differentiating in $t$ produces the factor
$(1/w - 2)$ in the integrand, so the left-hand side of
\eqref{eq:Push-psi-t} is
\[
\pa_t \psi_{t,r,n}(x)
= \frac{1}{2\pi\mathrm{i}}
\oint_{\gamma_\rho} \frac{\diff w}{w}
\frac{(1-w)^n}{2^{r-x} w^{n+r-x}}
(1/w - 2)
e^{t(1/w - 2)}.
\]
On the other hand, a shift of $r \to r{+}1$ in
$\psi_{t,r,n}$ multiplies the integrand by $1/(2w)$, so
the right-hand side of \eqref{eq:Push-psi-t} produces
the factor $2(1/(2w) - 1) = 1/w - 2$,
confirming \eqref{eq:Push-psi-t}.

\medskip
\noindent\textit{Verification of \eqref{eq:Push-psi-n}.}
A shift of $r \to r{+}1$ in $\psi_{t,r,n}$ multiplies
the integrand by $1/(2w)$.  The right-hand side of
\eqref{eq:Push-psi-n} therefore produces the factor
$2 \cdot 1/(2w) - 1 = (1-w)/w$, giving
\[
2\psi_{t,r+1,n}(x) - \psi_{t,r,n}(x)
= \frac{1}{2\pi\mathrm{i}}
\oint_{\gamma_\rho} \frac{\diff w}{w}
\frac{(1-w)^{n+1}}{2^{r-x} w^{n+1+r-x}}
e^{t(1/w - 2)}
= \psi_{t,r,n+1}(x).
\]

\medskip
\noindent\textit{Free recurrences for
$\overline{\mathcal{S}}$.}
We record two identities for the free operator that will
be applied inside the hitting-time representation of
$\phi$.  For $k \geq 1$,
\begin{equation}\label{eq:Push-Sbar-recurrence-n}
2\overline{\mathcal{S}}_{-t,k}(z, r{-}1)
- \overline{\mathcal{S}}_{-t,k}(z, r)
= \overline{\mathcal{S}}_{-t,k-1}(z, r).
\end{equation}
Indeed, the $r$-dependent factors in the contour formula
are $\tfrac{(1{-}w)^{r-z+k-1}}{2^{z-r}}$.  A shift of
$r \to r{-}1$ and the prefactor of $2$ together replace
$(1{-}w)^{r-z+k-1}$ by $(1{-}w)^{r-z+k-2}$.
The left-hand side of \eqref{eq:Push-Sbar-recurrence-n}
therefore produces the factor
$1 - (1{-}w) = w$, reducing $w^{-(k-1)}$ to
$w^{-(k-2)}$ and confirming
\eqref{eq:Push-Sbar-recurrence-n}.

We will also use
\begin{equation}\label{eq:Push-Sbar-recurrence-t}
\pa_t \overline{\mathcal{S}}_{-t,k}(z, r)
= 2\bigl(\overline{\mathcal{S}}_{-t,k}(z, r)
- \overline{\mathcal{S}}_{-t,k}(z, r{-}1)\bigr).
\end{equation}
Differentiating $\overline{\mathcal{S}}_{-t,k}(z, r)$
in $t$ produces the factor
$(2 - \tfrac{1}{1{-}w}) = \tfrac{1{-}2w}{1{-}w}$ in the integrand.
On the other hand, a shift of $r \to r{-}1$ multiplies
the integrand by $1/(2(1{-}w))$, so the right-hand side of
\eqref{eq:Push-Sbar-recurrence-t} produces the factor
$2(1 - 1/(2(1{-}w))) = (1{-}2w)/(1{-}w)$,
confirming \eqref{eq:Push-Sbar-recurrence-t}.

Finally, $\overline{\mathcal{S}}_{-t,0} \equiv 0$.
When $k = 0$, the factor $w^{-(k-1)} = w$ cancels the
$w^{-1}$ in the measure, removing the pole at the
origin.  The remaining integrand is analytic inside
$\gamma_\delta$ (the only other possible singularity, at
$w = 1$, lies outside $\gamma_\delta$), and the integral
vanishes by Cauchy's theorem.

\medskip
\noindent\textit{Verification of \eqref{eq:Push-phi-t}.}
The $t$-recurrence for $\phi$ follows from the
free recurrence \eqref{eq:Push-Sbar-recurrence-t} by
the same argument as in the continuous-time TASEP.
Since $\pa_t$ does not shift the particle label $n$,
the indicator $\mathbf{1}_{\tau < n}$ is unchanged,
and $\pa_t$ passes into the expectation (justified by
dominated convergence: differentiation multiplies the
contour integrand by the bounded factor $(2 - 1/(1{-}w))$,
so the same dominating function applies).  Inside the
expectation, the free recurrence
\eqref{eq:Push-Sbar-recurrence-t} gives
$\pa_t \overline{\mathcal{S}}_{-t,n-\tau}
(W_\tau, r) = 2(\overline{\mathcal{S}}_{-t,n-\tau}
(W_\tau, r) - \overline{\mathcal{S}}_{-t,n-\tau}
(W_\tau, r{-}1))$, which after taking expectations gives 
$\pa_t \phi_{t,r,n} = 2\phi_{t,r,n}
- 2\phi_{t,r-1,n}$.
The identity holds at every $r$; evaluating at $r + 1$
gives \eqref{eq:Push-phi-t}.

\medskip
\noindent\textit{Verification of \eqref{eq:Push-phi-n}.}
The identity
$\phi_{t,r,n-1} = 2\phi_{t,r-1,n} - \phi_{t,r,n}$
holds at every $(r, n)$; evaluating at
$(r{+}1, n{+}1)$ gives \eqref{eq:Push-phi-n}.
The verification follows from the
free recurrence \eqref{eq:Push-Sbar-recurrence-n} by
precisely the same argument as in the continuous-time TASEP
(Lemma~\ref{lem:CT-TASEP-seed-linear}): applying
\eqref{eq:Push-Sbar-recurrence-n} inside the hitting-time
expectation, splitting according to whether $\tau \leq n-2$
or $\tau = n-1$, and using
$\overline{\mathcal{S}}_{-t,0} \equiv 0$ to eliminate the
boundary contribution.
\end{proof}

\medskip

The propagator $B_{\mathbf{a},\mathbf{n}}$ defined by
\eqref{eq:Push-TASEP-propagator} is identical to the
continuous-time TASEP propagator, since it is built from
the same transition matrix $Q$ with the same threshold
constraints.

\begin{lemma}\label{lem:Push-TASEP-propagator}
The propagator $B_{\mathbf{a},\mathbf{n}}$ satisfies:
\begin{enumerate}[label=\textup{(\roman*)},
leftmargin=*]
\item For all $i > j$ and $r, r' \in \Z$,
\begin{equation}\label{eq:Push-Ba-covariance}
[B_{\mathbf{a}+\mathbf{1}}]_{ij}(r{+}1, r'{+}1)
= [B_{\mathbf{a},\mathbf{n}}]_{ij}(r, r').
\end{equation}

\item For all $i > j$ and $r, r' \in \Z$,
\begin{equation}\label{eq:Push-Ba-splitting}
[B_{\mathbf{a}+\mathbf{1}}]_{ij}(r, r')
- [B_{\mathbf{a},\mathbf{n}}]_{ij}(r, r')
= \sum_{j < \mu < i}
[B_{\mathbf{a},\mathbf{n}}]_{i\mu}(r, a_\mu{+}1)[B_{\mathbf{a}+\mathbf{1}}]_{\mu j}
(a_\mu{+}1, r').
\end{equation}

\item The propagator is $t$-independent and
invariant under uniform shifts of the particle
labels: replacing each $n_i$ by $n_i + 1$ leaves
$[B_{\mathbf{a},\mathbf{n}}]_{ij}(r, r')$ unchanged.
\end{enumerate}
\end{lemma}

\begin{proof}
The proof is identical to
Lemma~\ref{lem:CT-TASEP-propagator}: the propagator
depends only on the transition matrix $Q$ and the
threshold constraints, both of which are the same as in
the continuous-time TASEP.
\end{proof}

\medskip

For $t \geq 0$, $\mathbf{a} \in \Z^m$, and
$1 \leq n_1 < \cdots < n_m$, define the
product-graph wave functions
$\Psi^{t,\mathbf{a},\mathbf{n}} \in
\Hom(\R^m, \ell^2(\Z))$ and
$\Phi^{t,\mathbf{a},\mathbf{n}} \in
\Hom(\ell^2(\Z), \R^m)$ componentwise by
\begin{align}
\Psi_p^{t,\mathbf{a},\mathbf{n}}
&\defeq \psi_{t,a_p,n_p}
+ \sum_{j > p}
  \sum_{w \leq a_j - 1}
  [B_{\mathbf{a}-\mathbf{1},\mathbf{n}}]_{jp}(w, a_p)
  \psi_{t,w,n_j},
\label{eq:Push-Psi-def} \\
\Phi_p^{t,\mathbf{a},\mathbf{n}}
&\defeq \phi_{t,a_p,n_p}
+ \sum_{j < p}
  \sum_{z \leq a_j}
  [B_{\mathbf{a},\mathbf{n}}]_{pj}(a_p, z)
  \phi_{t,z,n_j}.
\label{eq:Push-Phi-def}
\end{align}
These are the explicit forms of the abstract
constructions
\eqref{eq:Psi-def}--\eqref{eq:Phi-def} on the product
graph with distinguished shift $e^{\pa_a}$.
Under the identification
$u = (t, \mathbf{a}, \mathbf{n})$, the maps
$u \mapsto \Psi^{t,\mathbf{a},\mathbf{n}}$ and
$u \mapsto \Phi^{t,\mathbf{a},\mathbf{n}}$ become the
wave functions of the semi-discrete framework of
\S\ref{sec:sc-linear}, with shifts
$T = e^{\pa_\mathbf{a}}$, $\pa_1 = \pa_t$,
$S_2 = e^{\pa_\mathbf{n}}$ and constant scalar edge
weights
$C_1 = 2I_m$,
$\Lambda_1 = 2I_m$,
$C_2 = 2I_m$,
$\Lambda_2 = I_m$.
These form a trivial solution of the semi-discrete diamond
equations (Proposition~\ref{prop:sc-diamond}).

\begin{lemma}\label{lem:Push-TASEP-WF}
The product-graph wave functions satisfy the
semi-discrete linear problem
\eqref{eq:sc-linear-1}--\eqref{eq:sc-linear-j}
and its adjoint
\eqref{eq:sc-adjoint-1}--\eqref{eq:sc-adjoint-j}
with the shifts and edge weights above.
Explicitly:
\begin{align}
\pa_t \Psi^{t,\mathbf{a},\mathbf{n}}
&= 2\bigl(\Psi^{t,\mathbf{a}+\mathbf{1},\mathbf{n}}
   - \Psi^{t,\mathbf{a},\mathbf{n}}\bigr),
\label{eq:Push-WF-psi-t} \\
\Psi^{t,\mathbf{a},\mathbf{n}+\mathbf{1}}
&= 2\Psi^{t,\mathbf{a}+\mathbf{1},\mathbf{n}}
   - \Psi^{t,\mathbf{a},\mathbf{n}},
\label{eq:Push-WF-psi-n}
\end{align}
and
\begin{align}
\pa_t \Phi^{t,\mathbf{a}+\mathbf{1},\mathbf{n}}
&= 2\bigl(\Phi^{t,\mathbf{a}+\mathbf{1},\mathbf{n}}
   - \Phi^{t,\mathbf{a},\mathbf{n}}\bigr),
\label{eq:Push-WF-phi-t} \\
\Phi^{t,\mathbf{a}+\mathbf{1},\mathbf{n}}
&= 2\Phi^{t,\mathbf{a},\mathbf{n}+\mathbf{1}}
   - \Phi^{t,\mathbf{a}+\mathbf{1},\mathbf{n}+\mathbf{1}}.
\label{eq:Push-WF-phi-n}
\end{align}
\end{lemma}

\begin{proof}
\noindent\textit{Discrete recurrences
\eqref{eq:Push-WF-psi-n}, \eqref{eq:Push-WF-phi-n}.}
Under the identification
$\mathcal{T} = e^{\pa_{\mathbf{a}}}$,
$\mathcal{S}_2 = e^{\pa_{\mathbf{n}}}$,
$(C_2, D_2) = (2I_m, I_m)$, the definitions
\eqref{eq:Push-Psi-def}--\eqref{eq:Push-Phi-def} match
the abstract constructions
\eqref{eq:Psi-def}--\eqref{eq:Phi-def} on the product
graph $\mathcal{G}^m$.  Property~(iii) of
Lemma~\ref{lem:Push-TASEP-propagator} is the
$\mathcal{S}_2$-invariance
condition~\eqref{eq:Sk-invariance}; together
with~(i) and~(ii), the constant-scalar reduction of
Proposition~\ref{prop:constant-scalar} gives
admissibility, $P(\mathcal{S}_2 u) = P(u)$, and
$\Lambda_2 = D_2 = I_m$.  The seed recurrences
\eqref{eq:Push-psi-n}--\eqref{eq:Push-phi-n} give the
product-graph linear problems
\eqref{eq:seed-psi-product}--\eqref{eq:seed-phi-product}
with these edge weights, so
the identities \eqref{eq:Push-WF-psi-n}, \eqref{eq:Push-WF-phi-n} follow from
Lemma~\ref{lem:Psi-Phi-linear}.

\medskip
\noindent\textit{Continuous recurrences
\eqref{eq:Push-WF-psi-t}, \eqref{eq:Push-WF-phi-t}.}
We suppress the dependence on $t$ and $\mathbf{n}$,
writing $\Psi_p^{\mathbf{a}}$ for
$\Psi_p^{t,\mathbf{a},\mathbf{n}}$.
The propagator $B_{\mathbf{a},\mathbf{n}}$ is $t$-independent
(Lemma \ref{lem:Push-TASEP-propagator}(iii)), so
$\pa_t$ passes through the propagator sums in
\eqref{eq:Push-Psi-def}; the differentiated series
converges in the same operator topology since
$|\pa_t\psi|$ satisfies the same bounds as $|\psi|$
(the contour integrand gains only the bounded factor
$(1/w - 2)$).
For the $p$-th component,
applying the seed recurrence \eqref{eq:Push-psi-t}:
\begin{align*}
\pa_t \Psi_p^{\mathbf{a}}
= 2(\psi_{t,a_p+1,n_p}
   - \psi_{t,a_p,n_p})
+ 2\sum_{j > p}
  \sum_{w \leq a_j - 1}
  [B_{\mathbf{a}-\mathbf{1},\mathbf{n}}]_{jp}(w, a_p)
  (\psi_{t,w+1,n_j} - \psi_{t,w,n_j}).
\end{align*}
It remains to show that this equals
$2(\Psi_p^{\mathbf{a}+\mathbf{1}}
- \Psi_p^{\mathbf{a}})$.  By
definition~\eqref{eq:Push-Psi-def} at
$\mathbf{a} + \mathbf{1}$,
\[
\Psi_p^{\mathbf{a}+\mathbf{1}}
= \psi_{t,a_p+1,n_p}
+ \sum_{j > p}
  \sum_{w \leq a_j}
  [B_{\mathbf{a},\mathbf{n}}]_{jp}(w, a_p{+}1)
  \psi_{t,w,n_j}.
\]
By $\mathcal{T}$-covariance
\eqref{eq:Push-Ba-covariance},
$[B_{\mathbf{a},\mathbf{n}}]_{jp}(w, a_p{+}1)
= [B_{\mathbf{a}-\mathbf{1},\mathbf{n}}]_{jp}(w{-}1, a_p)$.
Substituting $w' = w - 1$ (with range
$w' \leq a_j - 1$) gives
\[
\Psi_p^{\mathbf{a}+\mathbf{1}}
= \psi_{t,a_p+1,n_p}
+ \sum_{j > p}
  \sum_{w' \leq a_j - 1}
  [B_{\mathbf{a}-\mathbf{1},\mathbf{n}}]_{jp}(w', a_p)
  \psi_{t,w'+1,n_j}.
\]
Comparing with \eqref{eq:Push-Psi-def}, the difference
$\Psi_p^{\mathbf{a}+\mathbf{1}} - \Psi_p^{\mathbf{a}}$
has leading term
$\psi_{t,a_p+1,n_p} - \psi_{t,a_p,n_p}$ and
propagator terms
$[B_{\mathbf{a}-\mathbf{1},\mathbf{n}}]_{jp}(w, a_p)
(\psi_{t,w+1,n_j} - \psi_{t,w,n_j})$,
which is exactly $\frac{1}{2}\pa_t \Psi_p^{\mathbf{a}}$
as computed above.  This gives
\eqref{eq:Push-WF-psi-t}.

For the adjoint, $\pa_t$ again passes through
the propagator sums in \eqref{eq:Push-Phi-def}.
For the $p$-th component, applying
\eqref{eq:Push-phi-t} termwise:
\begin{align*}
\pa_t \Phi_p^{\mathbf{a}}
= 2(\phi_{t,a_p,n_p}
   - \phi_{t,a_p-1,n_p})
+ 2\sum_{j < p}
  \sum_{z \leq a_j}
  [B_{\mathbf{a},\mathbf{n}}]_{pj}(a_p, z)
  (\phi_{t,z,n_j} - \phi_{t,z-1,n_j}).
\end{align*}
It remains to show that this equals
$2(\Phi_p^{\mathbf{a}}
- \Phi_p^{\mathbf{a}-\mathbf{1}})$.  By
definition~\eqref{eq:Push-Phi-def} at
$\mathbf{a} - \mathbf{1}$,
\[
\Phi_p^{\mathbf{a}-\mathbf{1}}
= \phi_{t,a_p-1,n_p}
+ \sum_{j < p}
  \sum_{z \leq a_j - 1}
  [B_{\mathbf{a}-\mathbf{1},\mathbf{n}}]_{pj}(a_p{-}1, z)
  \phi_{t,z,n_j}.
\]
By $\mathcal{T}$-covariance
\eqref{eq:Push-Ba-covariance},
$[B_{\mathbf{a}-\mathbf{1},\mathbf{n}}]_{pj}(a_p{-}1, z)
= [B_{\mathbf{a},\mathbf{n}}]_{pj}(a_p, z{+}1)$.
Substituting $z' = z + 1$ (with range
$z' \leq a_j$) gives
\[
\Phi_p^{\mathbf{a}-\mathbf{1}}
= \phi_{t,a_p-1,n_p}
+ \sum_{j < p}
  \sum_{z' \leq a_j}
  [B_{\mathbf{a},\mathbf{n}}]_{pj}(a_p, z')
  \phi_{t,z'-1,n_j}.
\]
Comparing with \eqref{eq:Push-Phi-def}, the difference
$\Phi_p^{\mathbf{a}}
- \Phi_p^{\mathbf{a}-\mathbf{1}}$ has leading term
$\phi_{t,a_p,n_p} - \phi_{t,a_p-1,n_p}$ and
propagator terms
$[B_{\mathbf{a},\mathbf{n}}]_{pj}(a_p, z)
(\phi_{t,z,n_j} - \phi_{t,z-1,n_j})$,
which is exactly $\frac{1}{2}\pa_t \Phi_p^{\mathbf{a}}$
as computed above.  This gives
\eqref{eq:Push-WF-phi-t}.
\end{proof}

\medskip

Dressing compatibility is the remaining hypothesis of
the semi-discrete Darboux theorem
(Theorem~\ref{thm:sc-Darboux}).  The following
lemma verifies it for the kernel
$K_{t,\mathbf{a},\mathbf{n}}$.

\begin{lemma}\label{lem:Push-TASEP-dressing}
The kernel $K_{t,\mathbf{a},\mathbf{n}}$ is dressing
compatible (Definition~\ref{def:sc-dressing}) with
$(\Psi, \Phi)$ and $z = 1$ for all $t > 0$,
$1 \leq n_1 < \cdots < n_m$, and
$a_i < y_{n_i}$ for every~$i$.
Explicitly:
\begin{align}
K_{t,\mathbf{a}+\mathbf{1},\mathbf{n}}
- K_{t,\mathbf{a},\mathbf{n}}
&= \Psi^{t,\mathbf{a}+\mathbf{1},\mathbf{n}}
   \Phi^{t,\mathbf{a}+\mathbf{1},\mathbf{n}},
\label{eq:Push-DC-a} \\
\pa_t K_{t,\mathbf{a},\mathbf{n}}
&= 2\Psi^{t,\mathbf{a}+\mathbf{1},\mathbf{n}}
   \Phi^{t,\mathbf{a},\mathbf{n}},
\label{eq:Push-DC-t} \\
K_{t,\mathbf{a},\mathbf{n}+\mathbf{1}}
- K_{t,\mathbf{a},\mathbf{n}}
&= 2\Psi^{t,\mathbf{a}+\mathbf{1},\mathbf{n}}
   \Phi^{t,\mathbf{a},\mathbf{n}+\mathbf{1}}.
\label{eq:Push-DC-n}
\end{align}
\end{lemma}

\begin{proof}
\noindent\textit{Identities \eqref{eq:Push-DC-a} and
\eqref{eq:Push-DC-n}.}
These are the $\mathcal{T}$- and
$\mathcal{S}_2$-factorization identities of
Lemma~\ref{lem:K-factorization}, applied to the
product graph $\mathcal{G}^m$ with distinguished shift
$\mathcal{T} = e^{\pa_{\mathbf{a}}}$ and
$\mathcal{S}_2 = e^{\pa_{\mathbf{n}}}$.  The
hypotheses are verified: the seed recurrences
\eqref{eq:Push-psi-n}--\eqref{eq:Push-phi-n} give the
product-graph linear problems with
$(C_2, D_2) = (2I_m, I_m)$; the propagator properties
of Lemma~\ref{lem:Push-TASEP-propagator} provide the
admissibility conditions; and the kernel
$K_{t,\mathbf{a},\mathbf{n}}$ of
Lemma~\ref{lem:Push-TASEP-kernel} matches the abstract
kernel \eqref{eq:K-def}.

It remains to verify \eqref{eq:Push-DC-t}.  We suppress $t$ and $\mathbf{n}$ dependence, writing
$\Psi_p^{\mathbf{a}}$ for
$\Psi_p^{t,\mathbf{a},\mathbf{n}}$.
Since $B_{\mathbf{a},\mathbf{n}}$ is $t$-independent,
$\pa_t$ acts only on the seed functions
$\psi, \phi$ in the kernel formula.  The derivatives
$\pa_t\psi$ and $\pa_t\phi$ satisfy the same
bounds as $\psi$ and $\phi$ (differentiation multiplies
each contour integrand by the bounded factor
$(1/w - 2)$ or $(2 - 1/(1{-}w))$ respectively), so the
differentiated series converges in trace norm; in
particular, $t \mapsto K_{t,\mathbf{a},\mathbf{n}}$ is
$C^1$ in trace norm.  We split
$K = K^{\mathrm{diag}} + K^{\mathrm{off}}$ and
handle each part separately.

\medskip
\noindent\textit{Step~1: Diagonal part.}
The diagonal contribution is
\[
K^{\mathrm{diag}} = \sum_{p=1}^m \sum_{r \leq a_p}
\psi_{t,r,n_p}\phi_{t,r,n_p}.
\]
Applying the seed recurrences
\eqref{eq:Push-psi-t}--\eqref{eq:Push-phi-t}:
\begin{align*}
\pa_t K^{\mathrm{diag}}
&= \sum_{p=1}^m \sum_{r \leq a_p}
\bigl[2(\psi_{t,r+1,n_p} - \psi_{t,r,n_p})
\phi_{t,r,n_p}
+ \psi_{t,r,n_p}
2(\phi_{t,r,n_p} - \phi_{t,r-1,n_p})\bigr] \\
&= 2\sum_{p=1}^m \sum_{r \leq a_p}
\bigl[\psi_{t,r+1,n_p}\phi_{t,r,n_p}
- \psi_{t,r,n_p}\phi_{t,r-1,n_p}\bigr].
\end{align*}
For each $L \leq a_p$, the partial sum over
$L \leq r \leq a_p$ telescopes:
\[
\sum_{r=L}^{a_p}
\bigl[\psi_{t,r+1,n_p}\phi_{t,r,n_p}
- \psi_{t,r,n_p}\phi_{t,r-1,n_p}\bigr]
= \psi_{t,a_p+1,n_p}\phi_{t,a_p,n_p}
- \psi_{t,L,n_p}\phi_{t,L-1,n_p}.
\]
The lower boundary term vanishes in trace norm as
$L \to -\infty$:
$\lVert\psi_{t,L,n_p}\rVert_{\ell^2}
= \lVert g\rVert_{\ell^2}$ is bounded while
\eqref{eq:Push-TASEP-phi-bound} gives
$\lVert\phi_{t,L-1,n_p}\rVert_{\ell^2} \to 0$.
Therefore
\begin{equation}\label{eq:Push-DC-t-diag}
\pa_t K^{\mathrm{diag}}
= 2\sum_{p=1}^m
\psi_{t,a_p+1,n_p}\phi_{t,a_p,n_p}.
\end{equation}

\noindent\textit{Step~2: Off-diagonal part.}
The off-diagonal contribution is
\[
K^{\mathrm{off}} = \sum_{i > j}
\sum_{r \leq a_i}\sum_{r' \leq a_j}
\psi_{t,r,n_i}
[B_{\mathbf{a},\mathbf{n}}]_{ij}(r, r')
\phi_{t,r',n_j}.
\]
Applying the seed recurrences
\eqref{eq:Push-psi-t} and \eqref{eq:Push-phi-t} to
$\psi_{t,r,n_i}$ and $\phi_{t,r',n_j}$ respectively:
\begin{align*}
\pa_t K^{\mathrm{off}}
= 2\sum_{i > j}
\sum_{r \leq a_i}\sum_{r' \leq a_j}
\bigl[&\psi_{t,r+1,n_i}
[B_{\mathbf{a},\mathbf{n}}]_{ij}(r, r')
\phi_{t,r',n_j} \\
&- \psi_{t,r,n_i}
[B_{\mathbf{a},\mathbf{n}}]_{ij}(r, r')
\phi_{t,r'-1,n_j}\bigr].
\end{align*}
The first term is split into the boundary $r=a_i$
and the interior $r\leq a_i-1$; reindexing the
interior and applying $\mathcal{T}$-covariance gives
\[
[B_{\mathbf{a},\mathbf{n}}]_{ij}(r-1,r')
= [B_{\mathbf{a}+\mathbf{1}}]_{ij}(r,r'+1).
\]
In the second term, the substitution $r' \to r'{+}1$
gives
\begin{align}
\pa_t K^{\mathrm{off}}
&= 2\sum_{i > j}\sum_{r' \leq a_j}
\psi_{t,a_i+1,n_i}
[B_{\mathbf{a},\mathbf{n}}]_{ij}(a_i,r')
\phi_{t,r',n_j}
\notag \\
&\quad
+2\sum_{i > j}
\sum_{r \leq a_i}\sum_{r' \leq a_j}
\psi_{t,r,n_i}
[B_{\mathbf{a}+\mathbf{1}}]_{ij}(r,r'+1)
\phi_{t,r',n_j}
\label{eq:Push-DC-t-off-expanded} \\
&\quad
-2\sum_{i > j}
\sum_{r \leq a_i}\sum_{r' \leq a_j-1}
\psi_{t,r,n_i}
[B_{\mathbf{a},\mathbf{n}}]_{ij}(r,r'+1)
\phi_{t,r',n_j}.
\notag
\end{align}
Applying $\mathcal{T}$-splitting
\eqref{eq:Push-Ba-splitting} to the middle propagator:
\begin{align*}
[B_{\mathbf{a}+\mathbf{1}}]_{ij}(r,r'+1)
= [B_{\mathbf{a},\mathbf{n}}]_{ij}(r,r'+1)
+ \sum_{j < \mu < i}
[B_{\mathbf{a},\mathbf{n}}]_{i\mu}(r,a_\mu+1)
[B_{\mathbf{a}+\mathbf{1}}]_{\mu j}(a_\mu+1,r'+1).
\end{align*}
The leading $[B_{\mathbf{a},\mathbf{n}}]$ term cancels the last
line of \eqref{eq:Push-DC-t-off-expanded} on the common
range $r\leq a_i$, $r'\leq a_j-1$.  The remaining
boundary is $r'=a_j$.  Using $\mathcal{T}$-covariance once
more,
\[
[B_{\mathbf{a}+\mathbf{1}}]_{\mu j}(a_\mu+1,r'+1)
= [B_{\mathbf{a},\mathbf{n}}]_{\mu j}(a_\mu,r'),
\]
we obtain
\begin{align}
\pa_t K^{\mathrm{off}}
&= 2\sum_{i > j}\sum_{r' \leq a_j}
\psi_{t,a_i+1,n_i}
[B_{\mathbf{a},\mathbf{n}}]_{ij}(a_i,r')
\phi_{t,r',n_j}
\tag{A} \\
&\quad
+2\sum_{i > j}\sum_{r \leq a_i}
\psi_{t,r,n_i}
[B_{\mathbf{a},\mathbf{n}}]_{ij}(r,a_j+1)
\phi_{t,a_j,n_j}
\tag{B} \\
&\quad
+2\sum_{i > \mu > j}
\sum_{r \leq a_i}\sum_{r' \leq a_j}
\psi_{t,r,n_i}
[B_{\mathbf{a},\mathbf{n}}]_{i\mu}(r,a_\mu+1)
[B_{\mathbf{a},\mathbf{n}}]_{\mu j}(a_\mu,r')
\phi_{t,r',n_j}.
\tag{C}
\end{align}

\noindent\textit{Step~3: Term matching.}
It remains to compare this decomposition with the
right-hand side of \eqref{eq:Push-DC-t},
$2\Psi^{t,\mathbf{a}+\mathbf{1},\mathbf{n}}
\Phi^{t,\mathbf{a},\mathbf{n}}$.
Using \eqref{eq:Push-Psi-def} at
$\mathbf{a}+\mathbf{1}$ gives
\[
\Psi_p^{\mathbf{a}+\mathbf{1}}
= \psi_{t,a_p+1,n_p}
+ \sum_{i > p}\sum_{r \leq a_i}
[B_{\mathbf{a},\mathbf{n}}]_{ip}(r,a_p+1)
\psi_{t,r,n_i},
\]
while \eqref{eq:Push-Phi-def} gives
\[
\Phi_p^{\mathbf{a}}
= \phi_{t,a_p,n_p}
+ \sum_{j < p}\sum_{r' \leq a_j}
[B_{\mathbf{a},\mathbf{n}}]_{pj}(a_p,r')
\phi_{t,r',n_j}.
\]
The product
$\Psi_p^{\mathbf{a}+\mathbf{1}}
\Phi_p^{\mathbf{a}}$ expands into four contributions;
summing over $p$, the diagonal
\eqref{eq:Push-DC-t-diag} matches the leading term
$2\sum_{p=1}^m \psi_{t,a_p+1,n_p}\phi_{t,a_p,n_p}$,
and~(A), (B), (C) match the remaining three.  Hence
\[
\pa_t K_{t,\mathbf{a},\mathbf{n}}
= \pa_t K^{\mathrm{diag}}
+ \pa_t K^{\mathrm{off}}
= 2\Psi^{t,\mathbf{a}+\mathbf{1},\mathbf{n}}
\Phi^{t,\mathbf{a},\mathbf{n}},
\]
which proves \eqref{eq:Push-DC-t}.

It remains to identify the resolvent region.
By Lemma~\ref{lem:Push-TASEP-kernel},
$F_{t,\mathbf{a},\mathbf{n}} \defeq
\det(I - K_{t,\mathbf{a},\mathbf{n}})$ equals the
multipoint gap probability
$\PP_y\bigl(\bigcap_i
\{Y_{n_i}(t) > a_i\}\bigr)$.  Since $K$ is trace class,
$F \ne 0$ if and only if the resolvent exists.
Particle positions are non-increasing in time, so
$Y_{n_i}(t) \le y_{n_i}$.  If no clocks among labels
$1, \dotsc, n_m$ ring during $[0, t]$ (an event of
probability $e^{-n_m t}$), then $Y_k(t) = y_k$ for
$k \le n_m$.  The resolvent at $z = 1$ therefore exists
precisely on
$\mathcal{R} \defeq \{(t, \mathbf{a}, \mathbf{n}) :
t > 0,
1 \le n_1 < \cdots < n_m,
a_i < y_{n_i} \text{ for all } i\}$.

The semi-discrete mixed diamond at a vertex~$u$ requires
the resolvent at $u$ and $S_2 u$ for the dressed edge
weight formulas, together with $TS_2 u$ for the
$\mathcal{M}$-inverse appearing in the displayed equation.
$S_2$ replaces $y_{n_i}$ by the smaller $y_{n_i+1}$
and $T$ adds~$1$ to each threshold, so the binding
constraint is $TS_2 u \in \mathcal{R}$, i.e.\
$a_i < y_{n_i+1} - 1$ for all~$i$.
\end{proof}

\subsection{Multipoint equation}
\label{sec:Push-TASEP-multipoint}

The dressed observable
\[
\mathcal{M}_{t,\mathbf{a},\mathbf{n}}
\defeq I + \Phi^{t,\mathbf{a},\mathbf{n}}
(I - K_{t,\mathbf{a},\mathbf{n}})^{-1}
\Psi^{t,\mathbf{a},\mathbf{n}}
\in \End(\R^m)
\]
has dressed edge weights
\begin{align}
\mathcal{M}_1(u)
&= \mathcal{M}(Tu)^{-1}C_1\mathcal{M}(u)
= 2\mathcal{M}(Tu)^{-1}\mathcal{M}(u),
\label{eq:Push-dressed-C1} \\
\mathcal{M}_2(u)
&= \mathcal{M}(Tu)^{-1}C_2\mathcal{M}(S_2 u)
= 2\mathcal{M}(Tu)^{-1}\mathcal{M}(S_2 u).
\label{eq:Push-dressed-C2}
\end{align}

\begin{theorem}\label{thm:Push-TASEP-multipoint}
The dressed edge weights satisfy the mixed diamond
equation
\begin{align}\label{eq:Push-TASEP-multipoint}
& (\pa_t \mathcal{M}_{t,\mathbf{a},
\mathbf{n}+\mathbf{1}})
\mathcal{M}_{t,\mathbf{a},
\mathbf{n}+\mathbf{1}}^{-1}
- (\pa_t \mathcal{M}_{t,\mathbf{a}+\mathbf{1},
\mathbf{n}})
\mathcal{M}_{t,\mathbf{a}+\mathbf{1},
\mathbf{n}}^{-1} \notag \\
& \qquad + \mathcal{M}_{t,\mathbf{a}+\mathbf{1},
\mathbf{n}}
\mathcal{M}_{t,\mathbf{a}+\mathbf{1},
\mathbf{n}+\mathbf{1}}^{-1}
- \mathcal{M}_{t,\mathbf{a},
\mathbf{n}}
\mathcal{M}_{t,\mathbf{a},
\mathbf{n}+\mathbf{1}}^{-1} = 0
\end{align}
at every $(t, \mathbf{a}, \mathbf{n})$ with $t > 0$,
$1 \le n_1 < \cdots < n_m$, and
$a_i < y_{n_i+1} - 1$ for all~$i$.
\end{theorem}

\begin{proof}
By Theorem~\ref{thm:sc-Darboux}, the pair
$(\mathcal{M}_k, \Lambda_k)$ satisfies the semi-discrete
diamond equations on~$\mathcal{R}$.  The mixed diamond equation
\eqref{eq:sc-diamond-mixed} for the pair $(1, 2)$
reads
\[
\mathcal{M}_1(u)\Lambda_2
+ \Lambda_1\mathcal{M}_2(u)
- \mathcal{M}_2(u)\Lambda_1
- \Lambda_2\mathcal{M}_1(S_2 u)
- \pa_t \mathcal{M}_2(u) = 0.
\]
Since $\Lambda_1 = 2I_m$ and $\Lambda_2 = I_m$
are constant scalar matrices, $\pa_t\Lambda_2 = 0$,
$\Lambda_2$ acts as the identity, and the two
$\Lambda_1$-terms cancel:
$\Lambda_1\mathcal{M}_2 - \mathcal{M}_2\Lambda_1
= 2\mathcal{M}_2 - \mathcal{M}_2 \cdot 2 = 0$.
The equation reduces to
\begin{equation}\label{eq:Push-mixed-substituted}
\mathcal{M}_1(u)
- \mathcal{M}_1(S_2 u)
= \pa_t \mathcal{M}_2(u).
\end{equation}

Substituting \eqref{eq:Push-dressed-C1}--\eqref{eq:Push-dressed-C2}
and evaluating at $u = (t, \mathbf{a}, \mathbf{n})$ with
$Tu = (\mathbf{a}{+}\mathbf{1}, \mathbf{n})$,
$S_2 u = (\mathbf{a}, \mathbf{n}{+}\mathbf{1})$:
\[
2\mathcal{M}_{\mathbf{a}+\mathbf{1},
\mathbf{n}}^{-1}\mathcal{M}_{\mathbf{a}, \mathbf{n}}
- 2\mathcal{M}_{\mathbf{a}+\mathbf{1},
\mathbf{n}+\mathbf{1}}^{-1}\mathcal{M}_{\mathbf{a},
\mathbf{n}+\mathbf{1}}
= \pa_t\bigl(
2\mathcal{M}_{\mathbf{a}+\mathbf{1},
\mathbf{n}}^{-1}\mathcal{M}_{\mathbf{a},
\mathbf{n}+\mathbf{1}}\bigr).
\]
Dividing by $2$ and expanding the time derivative on the
right-hand side using
$\pa_t(A^{-1}B)
= -A^{-1}(\pa_t A)A^{-1}B
+ A^{-1}(\pa_t B)$:
\begin{align*}
&\mathcal{M}_{\mathbf{a}+\mathbf{1},
\mathbf{n}}^{-1}\mathcal{M}_{\mathbf{a}, \mathbf{n}}
- \mathcal{M}_{\mathbf{a}+\mathbf{1},
\mathbf{n}+\mathbf{1}}^{-1}\mathcal{M}_{\mathbf{a},
\mathbf{n}+\mathbf{1}} \\
&\qquad = -\mathcal{M}_{\mathbf{a}+\mathbf{1},
\mathbf{n}}^{-1}(\pa_t \mathcal{M}_{\mathbf{a}+\mathbf{1},
\mathbf{n}})\mathcal{M}_{\mathbf{a}+\mathbf{1},
\mathbf{n}}^{-1}\mathcal{M}_{\mathbf{a},
\mathbf{n}+\mathbf{1}}
+ \mathcal{M}_{\mathbf{a}+\mathbf{1},
\mathbf{n}}^{-1}\pa_t \mathcal{M}_{\mathbf{a},
\mathbf{n}+\mathbf{1}}.
\end{align*}
Left-multiplying by
$\mathcal{M}_{\mathbf{a}+\mathbf{1}, \mathbf{n}}$
and right-multiplying by
$\mathcal{M}_{\mathbf{a}, \mathbf{n}+\mathbf{1}}^{-1}$:
\begin{align*}
\mathcal{M}_{\mathbf{a}, \mathbf{n}}\mathcal{M}_{\mathbf{a},
\mathbf{n}+\mathbf{1}}^{-1}
&- \mathcal{M}_{\mathbf{a}+\mathbf{1},
\mathbf{n}}\mathcal{M}_{\mathbf{a}+\mathbf{1},
\mathbf{n}+\mathbf{1}}^{-1} \\
&= -(\pa_t \mathcal{M}_{\mathbf{a}+\mathbf{1},
\mathbf{n}})\mathcal{M}_{\mathbf{a}+\mathbf{1},
\mathbf{n}}^{-1}
+ (\pa_t \mathcal{M}_{\mathbf{a},
\mathbf{n}+\mathbf{1}})\mathcal{M}_{\mathbf{a},
\mathbf{n}+\mathbf{1}}^{-1}.
\end{align*}
Rearranging:
\begin{align*}
(\pa_t \mathcal{M}_{\mathbf{a},
\mathbf{n}+\mathbf{1}})\mathcal{M}_{\mathbf{a},
\mathbf{n}+\mathbf{1}}^{-1}
&- (\pa_t \mathcal{M}_{\mathbf{a}+\mathbf{1},
\mathbf{n}})\mathcal{M}_{\mathbf{a}+\mathbf{1},
\mathbf{n}}^{-1} \\
&+ \mathcal{M}_{\mathbf{a}+\mathbf{1},
\mathbf{n}}\mathcal{M}_{\mathbf{a}+\mathbf{1},
\mathbf{n}+\mathbf{1}}^{-1}
- \mathcal{M}_{\mathbf{a},
\mathbf{n}}\mathcal{M}_{\mathbf{a},
\mathbf{n}+\mathbf{1}}^{-1} = 0.
\end{align*}
This is \eqref{eq:Push-TASEP-multipoint}.
\end{proof}

For $m = 1$, the dressed observable is scalar.

\begin{corollary}\label{cor:Push-TASEP-scalar}
For $t > 0$ and $a < y_{n+1} - 1$, the one-point Fredholm
determinant
$F_{t, a, n} \defeq \det(I - K_{t, a, n})$ satisfies
\begin{equation}\label{eq:Push-TASEP-scalar-HM}
F_{t,a,n}F_{t,a+1,n+1}
- F_{t,a+1,n}F_{t,a,n+1}
+ F_{t,a,n+1}\pa_t F_{t,a+1,n}
- F_{t,a+1,n}\pa_t F_{t,a,n+1} = 0.
\end{equation}
\end{corollary}

\begin{proof}
Proposition~\ref{prop:sc-scalar-HM}(b) applies with the
standard-basis data $(c_1, \lambda_1) = (2, 2)$,
$(c_2, \lambda_2) = (2, 1)$, $T = e^{\pa_a}$,
$S_2 = e^{\pa_n}$, and $\pa_1 = \pa_t$.
The Fredholm determinant $F_{t,a,n}$ is $C^1$ in $t$
(by the trace-norm differentiability established in
Lemma~\ref{lem:Push-TASEP-dressing}), and the boundary
condition $F_{t,a,n} \to 1$ as $a \to -\infty$ is the
standard $T$-boundary, with
$\pa_t F_{t,a,n} \to 0$ in the same limit.
Equation~\eqref{eq:sc-HM-d1} gives
\begin{align*}
c_1(Tu)\lambda_2(Tu)
&\bigl[F(u)F(TS_2 u) - F(Tu)F(S_2 u)\bigr] \\
&+ c_2(Tu)\bigl[F(S_2 u)\pa_t F(Tu)
  - F(Tu)\pa_t F(S_2 u)\bigr] = 0.
\end{align*}
Substituting $c_1\lambda_2 = 2$ and $c_2 = 2$:
\begin{align*}
2\bigl[F_{t,a,n}F_{t,a+1,n+1}
- F_{t,a+1,n}F_{t,a,n+1}\bigr]
+ 2\bigl[F_{t,a,n+1}\pa_t F_{t,a+1,n}
  - F_{t,a+1,n}\pa_t F_{t,a,n+1}\bigr] = 0.
\end{align*}
Dividing by $2$ gives \eqref{eq:Push-TASEP-scalar-HM}.
\end{proof}

\begin{remark}[Structure of the multipoint equation]
\label{rem:Push-TASEP-structure}
Equation \eqref{eq:Push-TASEP-multipoint} is an
$m \times m$ matrix differential-difference equation
coupling the dressed observable $\mathcal{M}$ at four
lattice points, together with the time derivative
$\pa_t$.  For $m \geq 2$, the equation is
intrinsically noncommutative: the matrix inverses and
products do not simplify to a scalar relation.
\end{remark}

\section{Asymmetric Simple Exclusion Process}
\label{sec:asep-onepoint-verification}

\paragraph{\textbf{System description.}}
The asymmetric simple exclusion process (ASEP) on \(\Z\) is a
continuous-time Markov process with state
\(\eta(t)=\{\eta_x(t)\}_{x\in\Z}\in\{0,1\}^{\Z}\),
where \(\eta_x(t)\) indicates whether site \(x\) is occupied at
time \(t\).  Fix \(p,q>0\) with \(p+q=1\) and \(p<q\).  For each
pair of neighbouring sites \((y,y{+}1)\), the occupation variables
exchange at rate \(p\) when a particle is at \(y\) and a hole at
\(y{+}1\), and at rate \(q\) in the reverse configuration; all
exchanges are driven by independent exponential clocks.  In
particle language: each particle independently attempts to jump one
step to the right at rate \(p\) and one step to the left at
rate \(q\), with jumps to occupied sites suppressed.  Set
\[
\tau\defeq \frac{p}{q},
\qquad
\gamma\defeq q-p.
\]
Since \(p<q\), the asymmetry is leftward: \(\tau<1\) and
\(\gamma>0\).  The process is started from step initial data
\(\eta_x(0)=\mathbf{1}_{x\geq1}\), so that particles occupy the
positive integers.  The height function
\[
N_x(t)\defeq \sum_{y\leq x}\eta_y(t)
\]
records the number of particles at or to the left of \(x\).  For
\(\zeta\in\C\setminus\R_{\geq0}\) and \(m\in\Z\), write
\[
\zeta_m\defeq \tau^m\zeta,
\qquad
\mathcal{L}_{t,x,m}(\zeta)
\defeq
\E\left[\frac{1}{(\zeta_m\tau^{N_x(t)};\tau)_\infty}\right].
\]
The power \((-\zeta_m)^s\) is taken with the principal branch of the
logarithm on \(\C\setminus\R_{\leq0}\).  Since \(\tau>0\) lies on the
positive real axis, \(\tau^s\) is single-valued and the branch
satisfies
\[
(-\zeta_{m+1})^s=\tau^s(-\zeta_m)^s.
\]
\paragraph{\textbf{Fredholm determinant formula.}}
The following formula is due to
Borodin, Corwin, and Sasamoto~\cite{BorodinCorwinSasamoto2014}.
For \(\zeta\in\C\setminus\R_{\geq0}\),
\begin{equation}
\label{eq:asep-fredholm}
\mathcal{L}_{t,x,m}(\zeta)
=\det_H(I+K^{\mathrm{ASEP}}_{\zeta_m}),
\end{equation}
where \(H=L^2(\Gamma_w,\diff w/(2\pi\mathrm{i}))\), with
\(\Gamma_w\) a positively oriented simple closed contour
encircling \(0\) and \(-\tau\) but not \(-1\), and the kernel is
\begin{equation}
\label{eq:asep-BCS-kernel}
K^{\mathrm{ASEP}}_{\zeta_m}(w,w')
=
\frac{1}{2\pi\mathrm{i}}\int_{\mathcal{D}}
\Gamma(-s)\Gamma(1+s)(-\zeta_m)^s
\frac{f_{t,x}(w)}{f_{t,x}(\tau^s w)}
\frac{\diff s}{w'-\tau^s w},
\end{equation}
where \(\mathcal{D}\) is the contour \(D_{R,d}\) of
\cite[Def.~3.5]{BorodinCorwinSasamoto2014}, separating the
poles of \(\Gamma(-s)\) from those of \(\Gamma(1+s)\), with
\begin{equation}
\label{eq:asep-f-def}
f_{t,x}(z)\defeq \exp\bigl(\gamma t r(z)\bigr)r(z)^x,
\qquad
r(z)\defeq \frac{\tau}{\tau+z}.
\end{equation}
The reflection formula
\(\Gamma(-s)\Gamma(1+s)=-\pi/\sin(\pi s)\) reduces the kernel to
\[
K^{\mathrm{ASEP}}_{\zeta_m}(w,w')
=
\frac{1}{2\mathrm{i}}\int_{\mathcal{D}}
\frac{(-\zeta_m)^s}{\sin(\pi s)}
\frac{f_{t,x}(w)}{f_{t,x}(\tau^s w)}
\frac{\diff s}{\tau^s w-w'}.
\]

\paragraph{\textbf{References.}}
The ASEP was introduced by
Spitzer~\cite{Spitzer1970} and arose independently in biology in
the work of MacDonald, Gibbs, and
Pipkin~\cite{MacdonaldGibbsPipkin1968}.  The Fredholm determinant
formula \eqref{eq:asep-fredholm} and the contour conditions for the
kernel are Theorems~1.3 and~5.3 of
Borodin--Corwin--Sasamoto~\cite{BorodinCorwinSasamoto2014}.

\begin{remark}[Fredholm determinant convention]
\label{rem:asep-analytic}
The Fredholm determinant
\(\det_H(I+K^{\mathrm{ASEP}}_{\zeta_m})\) is defined by the
contour series with measure \(\diff w/(2\pi\mathrm{i})\) on
\(\Gamma_w\).  The series converges absolutely under the contour
conditions of
\cite[Thm.~5.3]{BorodinCorwinSasamoto2014}.  When the
determinant is nonzero, the resolvent kernel exists by
classical Fredholm theory.  The framework results invoked below
(gauge invariance of the Fredholm determinant, the semi-discrete
Darboux theorem, and the scalar Fredholm reduction) hold at the
kernel and series level: their proofs require only pointwise
kernel identities, the resolvent integral equation, and the
rank-one perturbation formula for Fredholm series.
\end{remark}

\subsection{Kernel reformulation}
\label{sec:asep-kernel-reform}

Define the gauge function
\(G_{t,x,m}(z)\defeq z^{-m}f_{t,x}(z)\) and the conjugated kernel
\[
K_{t,x,m}(w,w')
\defeq
G_{t,x,m}(w)^{-1}
K^{\mathrm{ASEP}}_{\zeta_m}(w,w')
G_{t,x,m}(w').
\]
Since \(G_{t,x,m}\) is nonvanishing on \(\Gamma_w\), gauge
invariance of the contour Fredholm determinant
(Remark~\ref{rem:asep-analytic}) gives
\[
\mathcal{L}_{t,x,m}(\zeta)
=\det_H(I+K_{t,x,m}).
\]
After the gauge transformation, the kernel takes the form
\begin{equation}
\label{eq:asep-gauged-kernel}
K_{t,x,m}(w,w')
=
\frac{w^m f_{t,x}(w')}{(w')^m}
\frac{1}{2\mathrm{i}}\int_{\mathcal{D}}
\frac{(-\zeta_m)^s}{\sin(\pi s)}
\frac{1}{f_{t,x}(\tau^s w)}
\frac{\diff s}{\tau^s w-w'}.
\end{equation}

The seed maps are
\begin{equation}
\label{eq:asep-seed-def}
\begin{aligned}
\psi_{t,x,m}(w)&\defeq
e^{\gamma t}w^m
\frac{1}{2\mathrm{i}}\int_{\mathcal{D}}
\frac{(-\zeta_m)^s}{\sin(\pi s)}
\frac{1}{f_{t,x}(\tau^s w)}\diff s, \\
\phi_{t,x,m}(w')&\defeq -\tau^{-1}e^{-\gamma t}G_{t,x+1,m}(w').
\end{aligned}
\end{equation}

\subsection{Base graph and seed data}
\label{sec:asep-seed-data}

Take \(E=\C\) and let \(\mathcal{V}\) be the lattice with
continuous direction and shifts
\[
\partial_1=\partial_t,
\qquad
T=e^{-\partial_x},
\qquad
S_2=e^{-\partial_x+\partial_m},
\]
so that for \(u=(t,x,m)\),
\begin{equation}
\label{eq:asep-lattice}
Tu=(t,x-1,m),
\qquad
S_2u=(t,x-1,m+1).
\end{equation}
The constant scalar edge weights are
\begin{equation}
\label{eq:asep-edge-weights}
(c_1,\lambda_1)=(-\gamma,-\gamma),
\qquad
(c_2,\lambda_2)=(-\tau,-\tau).
\end{equation}
The edge weights are constant scalars, so the semi-discrete
diamond equations (Proposition~\ref{prop:sc-diamond}) reduce to
commutativity.

\begin{lemma}
\label{lem:asep-seed-linear}
The seed functions \eqref{eq:asep-seed-def} satisfy the
semi-discrete linear problem
\eqref{eq:sc-linear-1}--\eqref{eq:sc-linear-j} and its adjoint
\eqref{eq:sc-adjoint-1}--\eqref{eq:sc-adjoint-j} with the data
\eqref{eq:asep-edge-weights}.  Explicitly:
\begin{align}
\partial_t\psi_{t,x,m}
&=-\gamma\psi_{t,x-1,m}+\gamma\psi_{t,x,m},
\label{eq:asep-psi-t} \\
\psi_{t,x-1,m+1}
&=-\tau\psi_{t,x-1,m}+\tau\psi_{t,x,m},
\label{eq:asep-psi-S}
\end{align}
and
\begin{align}
\partial_t\phi_{t,x-1,m}
&=-\gamma\phi_{t,x-1,m}+\gamma\phi_{t,x,m},
\label{eq:asep-phi-t} \\
\phi_{t,x-1,m}
&=-\tau\phi_{t,x-1,m+1}
+\tau\phi_{t,x-2,m+1}.
\label{eq:asep-phi-S}
\end{align}
\end{lemma}

\begin{proof}
Recall from \eqref{eq:asep-f-def} that
\[
f_{t,x}(z)=\exp\bigl(\gamma t r(z)\bigr)r(z)^x,
\qquad
r(z)=\frac{\tau}{\tau+z}.
\]
The contour conditions of \cite[Thm.~5.3]{BorodinCorwinSasamoto2014}
ensure that \(\tau^s w\) is uniformly separated from \(w'\) and that
the \(f\)-ratios in the Mellin--Barnes formula are bounded on
\(\mathcal{D} \times \Gamma_w \times \Gamma_w\).  Combined with the
exponential decay of \((-\zeta_m)^s/\sin(\pi s)\) along
\(\mathcal{D}\), the integrals below converge absolutely and
uniformly for \(w \in \Gamma_w\).  All manipulations, including
differentiation in \(t\), may therefore be performed under the
integral sign.

\noindent\textit{Verification of \eqref{eq:asep-psi-t}.}
Since
\[
\partial_t f_{t,x}(z)=\gamma r(z)f_{t,x}(z),
\qquad
f_{t,x-1}(z)=r(z)^{-1}f_{t,x}(z),
\]
differentiating \(\psi_{t,x,m}\) gives
\[
\partial_t\psi_{t,x,m}(w)
=\gamma\psi_{t,x,m}(w)
-\gamma w^m
\frac{e^{\gamma t}}{2\mathrm{i}}\int_{\mathcal{D}}
\frac{(-\zeta_m)^s}{\sin(\pi s)}
\frac{r(\tau^s w)}{f_{t,x}(\tau^s w)}\diff s.
\]
Since \(r(\tau^s w)/f_{t,x}(\tau^s w)=1/f_{t,x-1}(\tau^s w)\), the second term
is \(-\gamma\psi_{t,x-1,m}(w)\), proving \eqref{eq:asep-psi-t}.

\medskip
\noindent\textit{Verification of \eqref{eq:asep-psi-S}.}
The shift \(m\to m+1\) replaces \(w^m\) by \(w^{m+1}\) in
the seed definition.  The branch relation
\((-\zeta_{m+1})^s=\tau^s(-\zeta_m)^s\) and the identity
\(f_{t,x-1}^{-1}=r f_{t,x}^{-1}\) then give
\[
\psi_{t,x-1,m+1}(w)
=e^{\gamma t}w^m
\frac{1}{2\mathrm{i}}\int_{\mathcal{D}}
\frac{(-\zeta_m)^s}{\sin(\pi s)}
\frac{w\tau^s r(\tau^s w)}{f_{t,x}(\tau^s w)}\diff s.
\]
Expanding \(r(z)=\tau/(\tau+z)\) gives
\(w\tau^s r(\tau^s w)=\tau\bigl(1-r(\tau^s w)\bigr)\), so
\[
\psi_{t,x-1,m+1}
=\tau\psi_{t,x,m}-\tau\psi_{t,x-1,m},
\]
which is \eqref{eq:asep-psi-S}.

\medskip
\noindent\textit{Verification of \eqref{eq:asep-phi-t}.}
Since \(\phi_{t,x-1,m}=-\tau^{-1}e^{-\gamma t}G_{t,x,m}\),
the Leibniz rule gives
\[
\partial_t\phi_{t,x-1,m}(z)
=-\tau^{-1}\bigl[-\gamma e^{-\gamma t}G_{t,x,m}(z)
+e^{-\gamma t}\partial_t G_{t,x,m}(z)\bigr].
\]
Since \(\partial_t G_{t,x,m}=\gamma r(z)G_{t,x,m}\), substituting gives
\[
\partial_t\phi_{t,x-1,m}
=-\gamma\phi_{t,x-1,m}
-\tau^{-1}e^{-\gamma t}\gamma r(z)G_{t,x,m}.
\]
The second term is \(\gamma\phi_{t,x,m}\), since
\(r(z)G_{t,x,m}(z)=z^{-m}f_{t,x+1}(z)=G_{t,x+1,m}(z)\).
This proves~\eqref{eq:asep-phi-t}.

\medskip
\noindent\textit{Verification of \eqref{eq:asep-phi-S}.}
Expanding the right-hand side of \eqref{eq:asep-phi-S},
\[
-\tau\phi_{t,x-1,m+1}+\tau\phi_{t,x-2,m+1}
=e^{-\gamma t}\bigl(G_{t,x,m+1}(z)-G_{t,x-1,m+1}(z)\bigr).
\]
Since \(G_{t,x-1,m+1}=z^{-m-1}r(z)^{-1}f_{t,x}(z)\), the
difference becomes
\[
e^{-\gamma t}z^{-m-1}f_{t,x}(z)
\left(1-\frac{\tau+z}{\tau}\right)
=-\tau^{-1}e^{-\gamma t}z^{-m}f_{t,x}(z)
=\phi_{t,x-1,m}(z).
\]
This proves~\eqref{eq:asep-phi-S}.
\end{proof}

\subsection{Dressing compatibility}
\label{sec:asep-dressing}

We verify that the gauged kernel \(K_{t,x,m}\) satisfies the
dressing compatibility conditions of \S\ref{sec:sc-Darboux},
so that the semi-discrete Darboux theorem
(Theorem~\ref{thm:sc-Darboux}) applies.

\begin{lemma}
\label{lem:asep-dressing}
The gauged kernel \eqref{eq:asep-gauged-kernel} satisfies
\begin{align}
K_{t,x-1,m}-K_{t,x,m}
  &=\psi_{t,x-1,m}\phi_{t,x-1,m},
\label{eq:asep-K-T} \\
\partial_tK_{t,x,m}
  &=-\gamma\psi_{t,x-1,m}\phi_{t,x,m},
\label{eq:asep-K-t} \\
K_{t,x-1,m+1}-K_{t,x,m}
  &=-\tau\psi_{t,x-1,m}\phi_{t,x-1,m+1}.
\label{eq:asep-K-S}
\end{align}
\end{lemma}

In the notation of \S\ref{sec:sc-Darboux}, these are
\(K(Tu)-K(u)=\Psi(Tu)\Phi(Tu)\),
\(\partial_1 K(u)=\Psi(Tu)C_1(u)\Phi(u)\), and
\(K(S_2u)-K(u)=\Psi(Tu)C_2(u)\Phi(S_2u)\)
respectively.

\begin{proof}
Each identity is a Cauchy-denominator cancellation,
demonstrated first for the \(T\)-shift.

\medskip
\noindent\textit{Verification of \eqref{eq:asep-K-T}.}
Substituting \eqref{eq:asep-gauged-kernel} at \((t,x-1,m)\) and
\((t,x,m)\) and using \(f_{t,x-1}=r^{-1}f_{t,x}\), the
difference \(K_{t,x-1,m}-K_{t,x,m}\) equals
\[
\frac{w^m f_{t,x}(w')}{(w')^m}
\frac{1}{2\mathrm{i}}\int_{\mathcal{D}}
\frac{(-\zeta_m)^s}{\sin(\pi s) f_{t,x}(\tau^s w)}
\left(
\frac{r(\tau^s w)/r(w')}{\tau^s w-w'}
-\frac{1}{\tau^s w-w'}
\right)\diff s.
\]
Since \(r(z)=\tau/(\tau+z)\), the bracketed term equals
\((w'-\tau^s w)/[(\tau+\tau^s w)(\tau^s w-w')]\).
The factor \((\tau^s w-w')\) cancels between numerator and
denominator, and the difference reduces to
\begin{align*}
K_{t,x-1,m}(w,w')-K_{t,x,m}(w,w')
=
\left[
e^{\gamma t}w^m
\frac{1}{2\mathrm{i}}\int_{\mathcal{D}}
\frac{(-\zeta_m)^s}{\sin(\pi s)}
\frac{r(\tau^s w)}{f_{t,x}(\tau^s w)}\diff s
\right]
\left[
-\tau^{-1}e^{-\gamma t}
\frac{f_{t,x}(w')}{(w')^m}
\right].
\end{align*}
The first bracket is \(\psi_{t,x-1,m}(w)\) (since
\(r/f_{t,x}=1/f_{t,x-1}\)) and the second is
\(\phi_{t,x-1,m}(w')\), confirming \eqref{eq:asep-K-T}.

\medskip
\noindent\textit{Verification of \eqref{eq:asep-K-t}.}
Differentiating \eqref{eq:asep-gauged-kernel} in \(t\) using
\(\partial_t\log f_{t,x}=\gamma r\), the kernel derivative
equals
\[
\frac{w^m f_{t,x}(w')}{(w')^m}
\frac{\gamma}{2\mathrm{i}}\int_{\mathcal{D}}
\frac{(-\zeta_m)^s}{\sin(\pi s) f_{t,x}(\tau^s w)}
\frac{r(w')-r(\tau^s w)}{\tau^s w-w'}\diff s.
\]
Substituting \(r(z)=\tau/(\tau+z)\) gives
\([r(w')-r(\tau^s w)]/(\tau^s w-w')
=\tau/[(\tau+w')(\tau+\tau^s w)]\),
and the Cauchy factor cancels:
\begin{align*}
\partial_tK_{t,x,m}(w,w')
=
\left[
e^{\gamma t}w^m
\frac{1}{2\mathrm{i}}\int_{\mathcal{D}}
\frac{(-\zeta_m)^s}{\sin(\pi s)}
\frac{r(\tau^s w)}{f_{t,x}(\tau^s w)}\diff s
\right](-\gamma)
\left[
-\tau^{-1}e^{-\gamma t}
\frac{r(w')f_{t,x}(w')}{(w')^m}
\right].
\end{align*}
The first bracket is \(\psi_{t,x-1,m}(w)\), and
\(r(w')f_{t,x}(w')=f_{t,x+1}(w')\) identifies the second as
\(\phi_{t,x,m}(w')\), confirming \eqref{eq:asep-K-t}.

\medskip
\noindent\textit{Verification of \eqref{eq:asep-K-S}.}
Under the shift \((x,m)\to(x-1,m+1)\) in
\eqref{eq:asep-gauged-kernel}, the prefactor \(w^m/(w')^m\)
acquires a factor \(w/w'\), the branch relation gives
\((-\zeta_{m+1})^s=\tau^s(-\zeta_m)^s\), and
\(f_{t,x-1}=r^{-1}f_{t,x}\) acts on both \(f\)-factors.
The difference \(K_{t,x-1,m+1}-K_{t,x,m}\) equals
\[
\frac{w^m f_{t,x}(w')}{2\mathrm{i}(w')^{m}}
\int_{\mathcal{D}}
\frac{(-\zeta_m)^s}{\sin(\pi s) f_{t,x}(\tau^s w)}
\left(
\frac{w\tau^s(\tau+w')}{w'(\tau+\tau^s w)(\tau^s w-w')}
-\frac{1}{\tau^s w-w'}
\right)\diff s.
\]
The bracketed numerator is
\(w\tau^s(\tau+w')-w'(\tau+\tau^s w)
=\tau(\tau^s w-w')\),
so the Cauchy factor cancels:
\begin{align*}
K_{t,x-1,m+1}(w,w')-K_{t,x,m}(w,w')
&=
\left[
e^{\gamma t}w^m
\frac{1}{2\mathrm{i}}\int_{\mathcal{D}}
\frac{(-\zeta_m)^s}{\sin(\pi s)}
\frac{r(\tau^s w)}{f_{t,x}(\tau^s w)}\diff s
\right](-\tau) \\
&\quad\times
\left[
-\tau^{-1}e^{-\gamma t}
\frac{f_{t,x}(w')}{(w')^{m+1}}
\right].
\end{align*}
The first bracket is \(\psi_{t,x-1,m}(w)\), the scalar is
\(c_2=-\tau\), and the second is \(\phi_{t,x-1,m+1}(w')\),
confirming \eqref{eq:asep-K-S}.
\end{proof}

\subsection{One-point equations}
\label{sec:asep-onepoint-equations}

The dressed observable\footnote{The sign convention
\(\det_H(I+K^{\mathrm{ASEP}}_{\zeta_m})\) of
\cite{BorodinCorwinSasamoto2014} corresponds to \(z=-1\)
in the framework of \S\ref{sec:sc-Darboux}.}
\[
\mathcal{M}(u)
\defeq 1-\Phi(u)(I+K(u))^{-1}\Psi(u)
\in\C
\]
has dressed edge weights
\[
\mathcal{M}_1(u)=\mathcal{M}(Tu)^{-1}c_1\mathcal{M}(u),
\qquad
\mathcal{M}_2(u)=\mathcal{M}(Tu)^{-1}c_2\mathcal{M}(S_2u).
\]
We write \(\mathcal{M}_{t,x,m}\defeq\mathcal{M}(t,x,m)\)
throughout.

\begin{theorem}
\label{thm:asep-M-equation}
The dressed edge weights satisfy the semi-discrete mixed
diamond equation
\begin{equation}
\label{eq:asep-M-equation}
\partial_t\left(
\mathcal{M}_{t,x,m}^{-1}\mathcal{M}_{t,x,m+1}
\right)
+\gamma\left[
\mathcal{M}_{t,x,m}^{-1}\mathcal{M}_{t,x+1,m}
-
\mathcal{M}_{t,x-1,m+1}^{-1}\mathcal{M}_{t,x,m+1}
\right]=0
\end{equation}
at each \(\zeta \in \C \setminus \R_{\geq 0}\) for which
all displayed \(\mathcal{M}\)-values are defined.
\end{theorem}

\begin{proof}
The hypothesis ensures the resolvent exists at the
shifted vertices entering the diamond stencil, so by
Theorem~\ref{thm:sc-Darboux}
(Remark~\ref{rem:sc-domain-restriction}), the pair
\((\mathcal{M}_k,\lambda_k)\) satisfies the semi-discrete
diamond equations.  Since \(\lambda_1=-\gamma\) and
\(\lambda_2=-\tau\) are constant scalars, the mixed diamond
equation \eqref{eq:sc-mixed-dressed} with
\(\partial_1=\partial_t\) reduces to
$-\tau\bigl[\mathcal{M}_1(u)-\mathcal{M}_1(S_2u)\bigr]
= \partial_t\mathcal{M}_2(u)$.
Substituting
\(\mathcal{M}_1(u)=-\gamma\mathcal{M}(Tu)^{-1}\mathcal{M}(u)\)
and
\(\mathcal{M}_2(u)=-\tau\mathcal{M}(Tu)^{-1}\mathcal{M}(S_2u)\),
then dividing by \(\tau\), gives
\begin{align*}
\partial_t\left(
\mathcal{M}_{t,x-1,m}^{-1}\mathcal{M}_{t,x-1,m+1}
\right)
+\gamma\left[
\mathcal{M}_{t,x-1,m}^{-1}\mathcal{M}_{t,x,m}
-
\mathcal{M}_{t,x-2,m+1}^{-1}\mathcal{M}_{t,x-1,m+1}
\right]=0.
\end{align*}
Replacing \(x\) by \(x+1\) gives \eqref{eq:asep-M-equation}.
\end{proof}

\begin{remark}\label{rem:asep-M-continuation}
At \(z = -1\), Corollary~\ref{cor:sc-Woodbury} gives
\[
\mathcal{M}_{t,x,m}
= \frac{\mathcal{L}_{t,x+1,m}}{\mathcal{L}_{t,x,m}},
\]
so \eqref{eq:asep-M-equation} becomes a rational identity in
\(\tau\)-Laplace transforms.
Corollary~\ref{cor:asep-F-equation} establishes the
corresponding bilinear identity via the scalar reduction.
\end{remark}

\begin{corollary}
\label{cor:asep-F-equation}
For \(\zeta \in \C \setminus \R_{\geq 0}\), the
\(\tau\)-Laplace transform
\(\mathcal{L}_{t,x,m}(\zeta)
=\det_H(I+K_{t,x,m})\) satisfies
\begin{align}\label{eq:asep-F-equation}
\mathcal{L}_{t,x,m}\partial_t\mathcal{L}_{t,x,m+1}
-\mathcal{L}_{t,x,m+1}\partial_t\mathcal{L}_{t,x,m}
+\gamma\left[
\mathcal{L}_{t,x+1,m}\mathcal{L}_{t,x-1,m+1}
-\mathcal{L}_{t,x,m}\mathcal{L}_{t,x,m+1}
\right]=0.
\end{align}
\end{corollary}

The proof of Corollary~\ref{cor:asep-F-equation} requires the following boundary estimate.

\begin{lemma}
\label{lem:asep-boundary}
Fix \(t \geq 0\), \(\zeta \in \C \setminus \R_{\geq 0}\),
and \(m \in \Z\).  Then
\[
\mathcal{L}_{t,x,m}(\zeta) \to 1,
\qquad
\partial_t \mathcal{L}_{t,x,m}(\zeta) \to 0
\]
as \(x \to +\infty\).
\end{lemma}

\begin{proof}[Proof of Lemma~\ref{lem:asep-boundary}]
Write \(h_m(n) \defeq 1/(\zeta_m \tau^n; \tau)_\infty\),
so that \[\mathcal{L}_{t,x,m} = \E[h_m(N_x(t))].\]
Since \(\zeta \notin \R_{\geq 0}\), no factor in the
\(\tau\)-Pochhammer product vanishes and
\(h_m(n) \to 1\) as \(n \to \infty\); in particular
\(h_m\) is bounded.  As \(N_x(t) \to \infty\) almost
surely when \(x \to +\infty\), dominated convergence
gives \(\mathcal{L}_{t,x,m} \to 1\).

For the time derivative, only jumps across the bond
\((x, x{+}1)\) change \(N_x\), so the generator gives
\begin{align*}
\partial_t \mathcal{L}_{t,x,m}
&= \E\bigl[
p\eta_x(t)(1 - \eta_{x+1}(t))
\bigl(h_m(N_x(t) {-} 1) - h_m(N_x(t))\bigr) \\
&\quad + q(1 - \eta_x(t))\eta_{x+1}(t)
\bigl(h_m(N_x(t) {+} 1) - h_m(N_x(t))\bigr)
\bigr].
\end{align*}
Since \(h_m\) converges, the differences
\(h_m(N_x \pm 1) - h_m(N_x) \to 0\) as
\(N_x(t) \to \infty\).  The integrand is therefore
bounded and tends to zero almost surely, so
dominated convergence gives
\(\partial_t \mathcal{L}_{t,x,m} \to 0\).
\end{proof}

\begin{proof}[Proof of Corollary~\ref{cor:asep-F-equation}]
Fix \(t, x, m\), and let \(B(\zeta)\) denote the left-hand
side of \eqref{eq:asep-F-equation}.  Each
\(\mathcal{L}\)-value and \(t\)-derivative appearing in
\eqref{eq:asep-F-equation} is analytic in \(\zeta\) on
\(\C \setminus \R_{\geq 0}\), so \(B\) is holomorphic on
this connected domain.

For \(|\zeta|\) sufficiently small, the bound
\(|\zeta_\mu \tau^n| \leq |\zeta_\mu|\) for \(n \geq 0\)
keeps \(1/(\zeta_\mu \tau^n; \tau)_\infty\) close to \(1\)
uniformly in \(n\), so
\(\mathcal{L}_{t,y,\mu}(\zeta) \neq 0\) for every \(y\) and
\(\mu \in \{m, m{+}1\}\).  The resolvents required by
Proposition~\ref{prop:sc-scalar-HM}\textup{(b)} therefore
exist along the entire backward \(T\)-orbit.
Lemma~\ref{lem:asep-boundary} supplies the boundary
normalization, since
\(T^{-n}(t, x, m) = (t, x{+}n, m)\).  The coefficient ratio
\(\alpha_2 = c_1\lambda_2/c_2 = -\gamma\) is constant.

Proposition~\ref{prop:sc-scalar-HM}\textup{(b)} therefore
gives \eqref{eq:sc-HM-d1} with \(\partial_1 = \partial_t\).
Substituting the lattice identifications
\eqref{eq:asep-lattice} gives
\begin{align*}
\mathcal{L}_{t,x-1,m}\partial_t\mathcal{L}_{t,x-1,m+1}
-\mathcal{L}_{t,x-1,m+1}\partial_t\mathcal{L}_{t,x-1,m}
+\gamma\bigl[
\mathcal{L}_{t,x,m}\mathcal{L}_{t,x-2,m+1}
-\mathcal{L}_{t,x-1,m}\mathcal{L}_{t,x-1,m+1}
\bigr]=0.
\end{align*}
Replacing \(x\) by \(x+1\) yields
\(B(\zeta) = 0\) for \(|\zeta|\) small in the slit plane.
Since \(\C \setminus \R_{\geq 0}\) is connected, the
identity theorem gives \(B \equiv 0\).
\end{proof}

%% file: chapters/9-parabolic-model-verifications.tex
\chapter{Parabolic Model Verifications}
\label{ch:parabolic-model-verifications}

{
  \setlength{\parskip}{0pt}
}

\section{Reflected Brownian Motion}\label{sec:RBM}

\paragraph{\textbf{System description.}}
Let $\mathbf{y} = (y_1, y_2, \dotsc)$ with
$y_1 \geq y_2 \geq \cdots$, and let
$B_1(t), B_2(t), \dotsc$ be independent standard Brownian
motions with $B_k(0) = 0$.  One-sided reflected Brownian
motion (RBM) with initial data $\mathbf{y}$ is defined
recursively by
\begin{align}\label{eq:RBM-Skorokhod}
Y_1(t) &= y_1 + B_1(t), \notag \\
Y_{k+1}(t) &= y_{k+1} + B_{k+1}(t)
- \sup_{0 \leq s \leq t}
[y_{k+1} + B_{k+1}(s) - Y_k(s)]^+,
\quad k \geq 1,
\end{align}
where $[\cdot]^+ = \max(\cdot, 0)$.
The first particle is a free Brownian motion; for
$k \geq 1$, particle~$k+1$ is pushed downward whenever it
meets particle~$k$.

\paragraph{\textbf{Fredholm determinant formula.}}
Fix a time $t > 0$ and initial data $\mathbf{y}$.  Fix
finitely many particle labels
$\mathbf{n} = (n_1, \dotsc, n_m)$ with
$1 \leq n_1 < n_2 < \cdots < n_m$, and spatial
thresholds $\mathbf{a} = (a_1, \dotsc, a_m) \in \R^m$.  The
multipoint joint cumulative distribution of the particle
positions is given by a Fredholm determinant on
$L^2(\{n_1, \dotsc, n_m\} \times \R)$:
\[
\PP_{\mathbf{y}}\Bigl(\bigcap_{i=1}^m
\{Y_{n_i}(t) > a_i\}\Bigr)
= \det(I - \bar{\chi}_{\mathbf{a}} K_t
\bar{\chi}_{\mathbf{a}})_{L^2(\{n_1, \dotsc, n_m\}
\times \R)},
\]
where $\bar{\chi}_{\mathbf{a}}(n_i, u) =
\mathbf{1}_{u \leq a_i}$.
The correlation kernel $K_t$ has block entries
\[
K_t(n_i, u; n_j, v)
= -e^{v - u}\pa^{-(n_j - n_i)}(u, v)
\mathbf{1}_{n_i < n_j}\mathbf{1}_{v < u}
+ \bigl(\mathcal{S}_{-t,-n_i}\bigr)^*
\overline{\mathcal{S}}_{t,n_j}^{\operatorname{epi}(\mathbf{y})}
(u, v),
\]
where $A^*(u, v) = A(v, u)$ denotes the transpose kernel.
The operators defining the kernel are as follows.

\medskip
\noindent\textit{The transition kernel.}
The kernel $e^{v-u}\pa^{-1}(u, v)$ governs a random walk
on $\R$ taking $\mathrm{Exp}(1)$ steps to the left:
\[
e^{v-u}\pa^{-1}(u, v) = e^{v-u}\mathbf{1}_{u > v},
\]
with $k$-step transition kernel
\[
e^{v-u}\pa^{-k}(u, v) = e^{v-u}\frac{(u-v)^{k-1}}{(k-1)!}
\mathbf{1}_{u > v}.
\]

\medskip
\noindent\textit{The scattering operators
$\mathcal{S}_{-t,-n}$ and
$\overline{\mathcal{S}}_{-t,n}$.}
These are built from the Hermite functions.
The $n$-th Hermite polynomial is
$H_n(x) = (-1)^n e^{x^2/2}\frac{\diff^n}{\diff x^n} e^{-x^2/2}$,
and the associated functions are
\begin{align}
\varphi_n(t, x) &\defeq t^{-n/2}\frac{1}{\sqrt{2\pi t}}
e^{-x^2/2t}H_n(x/\sqrt{t}),
\label{eq:RBM-varphi-def} \\
\bar{\varphi}_n(t, x) &\defeq \frac{1}{n!}t^{n/2}
H_n(x/\sqrt{t}).
\label{eq:RBM-varphi-bar-def}
\end{align}
The scattering operators are
\begin{align}
\mathcal{S}_{-t, -n}(u, v) &\defeq e^{u-v}\varphi_n(t, u-v),
\label{eq:RBM-scattering-def} \\
\overline{\mathcal{S}}_{-t, n}(u, v) &\defeq e^{v-u}\bar{\varphi}_{n-1}(t, u-v).
\label{eq:RBM-scattering-bar-def}
\end{align}

\medskip
\noindent\textit{The epigraph operator
$\overline{\mathcal{S}}_{t,n}^{\operatorname{epi}(\mathbf{y})}$.}
The initial data enters through a hitting probability
operator.  Let $(B_k)_{k \geq 0}$ be the random walk with
transition kernel $e^{v-u}\pa^{-1}(u, v)$, and define
$\tau = \inf\{k \geq 0 : B_k \geq y_{k+1}\}$.  The
epigraph operator is
\begin{equation}\label{eq:RBM-epigraph-def}
\overline{\mathcal{S}}_{t, n}^{\operatorname{epi}(\mathbf{y})}(u, v)
\defeq \E_{B_0 = u}\bigl[
\overline{\mathcal{S}}_{-t, n-\tau}(B_\tau, v)
\mathbf{1}_{\tau < n}\bigr].
\end{equation}

\paragraph{\textbf{References.}}
The Fredholm determinant formula above is the
conjugated form of the Nica--Quastel--Remenik formula
(Theorem~2.1 of~\cite{NQR}) for one-sided reflected
Brownian motion.

\subsection{Kernel reformulation}\label{sec:RBM-kernel-reform}

The extended-kernel Fredholm determinant on
$L^2(\{n_1, \dotsc, n_m\} \times \R)$ reduces to a
single-space Fredholm determinant on $L^2(\R)$.  The
seed functions are
\begin{equation}\label{eq:RBM-seed-def}
\psi_{t,r,n}(x) \defeq \mathcal{S}_{-t,-n}(x, r),
\qquad
\phi_{t,r,n}(x) \defeq
\overline{\mathcal{S}}_{t,n}^{\operatorname{epi}(\mathbf{y})}(x, r),
\end{equation}
where the dependence of $\phi$ on the initial data
$\mathbf{y}$ is suppressed.  The propagator is the
continuous analogue of the directed-path propagator
(Definition~\ref{def:directed-path-propagator}), with
transition kernels
$N_{ij}(x, x') = -e^{x'-x}\pa^{-(n_j-n_i)}(x, x')
\mathbf{1}_{n_i < n_j}$.
For $1 \leq i \leq j \leq m$,
\begin{align}\label{eq:RBM-propagator-def}
&[B_{\mathbf{a},\mathbf{n}}]_{ij}(z, w) \defeq \notag \\
&\begin{cases}
\displaystyle\sum_{k=1}^{j-i}
\sum_{\substack{i = i_0 < i_1 < \cdots < i_k = j}}
\int_{-\infty}^{a_{i_1}} \cdots
\int_{-\infty}^{a_{i_{k-1}}}
\prod_{r=0}^{k-1} N_{i_r, i_{r+1}}(\xi_r, \xi_{r+1})
\diff\xi_{k-1} \cdots \diff\xi_1,
& i < j, \\[6pt]
\delta(z - w), & i = j, \\
0, & i > j,
\end{cases}
\end{align}
with $\xi_0 = z$ and $\xi_k = w$.  The internal
vertices are constrained by the cutoffs, but the
endpoints $z, w$ are arbitrary.  The dependence on
$t$ is suppressed because the propagator is
$t$-independent (Lemma~\ref{lem:RBM-propagator}(iii)).

\begin{lemma}\label{lem:RBM-kernel-reform}
The multipoint distribution satisfies
\[
\PP_{\mathbf{y}}\Bigl(\bigcap_{i=1}^m
\{Y_{n_i}(t) > a_i\}\Bigr)
= \det(I - K_{t, \mathbf{a}, \mathbf{n}})_{L^2(\R)},
\]
where
\begin{equation}\label{eq:RBM-single-kernel}
K_{t, \mathbf{a}, \mathbf{n}}(x, x')
= \sum_{1 \leq i \leq j \leq m}
\int_{-\infty}^{a_i} \int_{-\infty}^{a_j}
\psi_{t, w, n_j}(x)
[B_{\mathbf{a},\mathbf{n}}]_{ij}(z, w)\phi_{t, z, n_i}(x')\diff w\diff z.
\end{equation}
\end{lemma}

\begin{proof}
The proof follows the same strategy as the kernel
reformulation of \S\ref{sec:CT-TASEP}
(Lemma~\ref{lem:CT-TASEP-kernel}): conjugation by
exponential weights, followed by Sylvester's identity and
identification of the inverse entries with the
propagator.

\medskip
\noindent\textit{Decomposition and conjugation.}
Write $\mathscr{X} = L^2(\{n_1, \dotsc, n_m\}
\times \R)$ for the extended space and
$H = L^2(\R)$.  The cutoff extended kernel
decomposes as
\[
\bar{\chi}_{\mathbf{a}} K_t \bar{\chi}_{\mathbf{a}}
= L + UV,
\]
where $L \defeq \bar{\chi}_{\mathbf{a}} N
\bar{\chi}_{\mathbf{a}}$ is strictly upper triangular
in the layer index.
The operators $U \colon H \to \mathscr{X}$ and
$V \colon \mathscr{X} \to H$ are defined by
\[
U(n_i, r; x)
\defeq \mathbf{1}_{r \leq a_i}\psi_{t,r,n_i}(x),
\qquad
V(x; n_j, r')
\defeq \phi_{t,r',n_j}(x)\mathbf{1}_{r' \leq a_j},
\]
so that $UV$ is the scattering part of the extended
kernel.  The operator $L$ is nilpotent ($L^m = 0$)
because it is strictly upper triangular in finitely
many block indices.

Choose $\gamma \in (0, 1)$ and set
$\beta \defeq e^{-\gamma}$, so that
$\beta \in (e^{-1}, 1)$.  Choose strictly ordered
conjugation parameters
\[
0 < \alpha_m < \alpha_{m-1} < \cdots < \alpha_1 < \gamma.
\]
For a function on $\{n_1, \dotsc, n_m\} \times \R$,
define $(\mathscr{C}f)_i(x) \defeq e^{\alpha_i x}
f_{n_i}(x)$.  For a block kernel $M$,
\[
(\mathscr{C}M\mathscr{C}^{-1})_{ij}(x, x')
= e^{\alpha_i x}M_{ij}(x, x')e^{-\alpha_j x'}.
\]
The conjugation multiplies each row of a finite minor
by $e^{\alpha_i x}$ and each column by
$e^{-\alpha_j x'}$, so the Fredholm series agrees
term by term with that of the conjugated kernel.
Set
$\widetilde{N} \defeq \mathscr{C}L\mathscr{C}^{-1}$,
$\widetilde{U} \defeq \mathscr{C}U$, and
$\widetilde{V} \defeq V\mathscr{C}^{-1}$.

\medskip
\noindent\textit{Nilpotent triangular part.}
The operator $\widetilde{N}$ is strictly upper triangular
in the $m$ block indices, hence nilpotent:
$\widetilde{N}^m = 0$.  The Neumann inverse
$(I - \widetilde{N})^{-1}
= I + \widetilde{N} + \cdots + \widetilde{N}^{m-1}$
is bounded on $\mathscr{X}$.  Each nonzero block is
Hilbert--Schmidt: for $i < j$ with
$d \defeq n_j - n_i$, the $(i,j)$-block of
$\widetilde{N}$ has modulus
\[
|\widetilde{N}_{ij}(x, x')|
= e^{(\alpha_i - \alpha_j)x'}
e^{-(1 - \alpha_i)(x - x')}
\frac{(x - x')^{d-1}}{(d-1)!}
\mathbf{1}_{x > x'}\mathbf{1}_{x \leq a_i}
\mathbf{1}_{x' \leq a_j},
\]
and
\[
\lVert\widetilde{N}_{ij}\rVert_{\mathfrak{S}_2}^2
\leq \int_{-\infty}^{a_j}
e^{2(\alpha_i - \alpha_j)x'}\diff x'
\int_0^{\infty}
e^{-2(1 - \alpha_i)s}
\frac{s^{2(d-1)}}{((d-1)!)^2}\diff s
< \infty,
\]
since $\alpha_i > \alpha_j$ and
$\alpha_i < \gamma < 1$.
Since there are finitely many nonzero blocks,
$\widetilde{N} \in \mathfrak{S}_2(\mathscr{X})$.
(When $d = 1$, the conjugated Volterra block is not
trace class; the determinant reduction below does
not require this.)

\medskip
\noindent\textit{Hilbert--Schmidt estimate for
$\widetilde{U}$.}
The conjugation acts on the extended-space
variable $r$:
$\widetilde{U}(n_i, r; x)
= e^{\alpha_i r}\mathbf{1}_{r \leq a_i}
\psi_{t,r,n_i}(x)$.  The seed function
$\psi_{t,r,n}(x) = e^{x-r}\varphi_n(t, x{-}r)$
has Gaussian decay inherited from $\varphi_n$:
setting $\xi = x - r$ shows that
$\lVert\psi_{t,r,n}\rVert_{L^2(\R)}^2
= \int_{\R}|e^{\xi}\varphi_n(t,\xi)|^2\diff\xi$
is independent of $r$ and finite for each
$t > 0$, $n \geq 0$, because the Gaussian factor
dominates the polynomial growth of $H_n$.
Therefore
\[
\lVert\widetilde{U}\rVert_{\mathfrak{S}_2}^2
= \sum_{i=1}^m
\int_{-\infty}^{a_i}
e^{2\alpha_i r}
\lVert\psi_{t,r,n_i}\rVert_{L^2(\R)}^2
\diff r
\leq C\sum_{i=1}^m
\frac{e^{2\alpha_i a_i}}{2\alpha_i}
< \infty,
\]
since each $\alpha_i > 0$ and the $L^2$-norm
of $\psi_{t,r,n_i}$ is $r$-independent.  Thus
$\widetilde{U} \in \mathfrak{S}_2$.

\medskip
\noindent\textit{Hilbert--Schmidt estimate for
$\widetilde{V}$.}
The conjugation acts on the extended-space
variable $r'$:
$\widetilde{V}(x; n_j, r')
= \phi_{t,r',n_j}(x)e^{-\alpha_j r'}
\mathbf{1}_{r' \leq a_j}$.
The seed $\phi$ is bounded via the hitting-time
representation \eqref{eq:RBM-epigraph-def}.  The free
scattering operator satisfies
$|\overline{\mathcal{S}}_{-t,k}(b, r)|
= e^{r-b}|\bar\varphi_{k-1}(t, b{-}r)|
\leq C\beta^{b-r}$
for all $1 \leq k \leq n_m$, $b \geq y_{n_m}$, and
$r \leq a_{\max}$: for $b - r \geq 0$, the polynomial
$|\bar\varphi_{k-1}(t, \cdot)|$ grows at most
polynomially while $e^{-(1-\gamma)(b-r)}$ decays
exponentially, giving
$e^{-(b-r)}|\bar\varphi_{k-1}|
\leq Ce^{-\gamma(b-r)} = C\beta^{b-r}$;
for $b - r < 0$, the pair $(b, r)$ ranges over
a bounded set.

Using the triangle inequality inside the expectation
and summing over the finitely many values
$k = n - \tau \in \{1, \dotsc, n\}$:
\[
|\phi_{t,r,n}(x)|
\leq C\beta^{-r}
\sum_{\ell=0}^{n-1}
\E_{B_0 = x}\bigl[\beta^{B_\ell}\bigr].
\]
Since $B_\ell = x - S_\ell$ where
$S_\ell = \sum_{k=1}^{\ell} X_k$ is the sum of
$\ell$ independent $\mathrm{Exp}(1)$ random variables,
$\E_{B_0 = x}[\beta^{B_\ell}]
= \beta^x(\E[\beta^{-X_1}])^\ell$.
The moment generating function of $\mathrm{Exp}(1)$
gives
$\E[\beta^{-X_1}]
= \E[e^{X_1\log(1/\beta)}]
= (1 - \gamma)^{-1}$,
which is finite because $\gamma < 1$.  The finite sum
$\sum_{\ell=0}^{n-1}(1 - \gamma)^{-\ell}$ is
bounded by a constant depending on $n$ and $\gamma$,
giving
\begin{equation}\label{eq:RBM-phi-bound}
|\phi_{t,r,n}(x)|
\leq C\beta^{x-r}\mathbf{1}_{x \geq y_{n_m}},
\end{equation}
where the support restriction follows because
$\{\tau < n\}$ requires $B_0 = x \geq y_{n_m}$
(the walk is strictly decreasing, so the hitting
condition forces $x \geq y_{\tau+1} \geq y_{n_m}$).
Using \eqref{eq:RBM-phi-bound}:
\begin{align*}
\lVert\widetilde{V}\rVert_{\mathfrak{S}_2}^2
= \sum_{j=1}^m \int_{-\infty}^{a_j}
e^{-2\alpha_j r'}
\lVert\phi_{t,r',n_j}\rVert_{L^2(\R)}^2
\diff r'
\leq C^2 \sum_{j=1}^m
\int_{-\infty}^{a_j}
e^{-2\alpha_j r'}\beta^{-2r'}
\diff r'
\int_{y_{n_m}}^{\infty}
\beta^{2x}\diff x.
\end{align*}
The $x$-integral is
$\beta^{2y_{n_m}}/(2\gamma) < \infty$.
The $r'$-integral equals
$\int_{-\infty}^{a_j}
e^{2(\gamma - \alpha_j)r'}\diff r'
= e^{2(\gamma - \alpha_j)a_j}/
\bigl(2(\gamma - \alpha_j)\bigr) < \infty$
because $\gamma > \alpha_j$.
Thus $\widetilde{V} \in \mathfrak{S}_2$.

\medskip
\noindent\textit{Determinant reduction and propagator
identification.}
Since $\mathfrak{S}_2$ is a two-sided ideal and
$(I - \widetilde{N})^{-1}$ is bounded,
$(I - \widetilde{N})^{-1}\widetilde{U}
\in \mathfrak{S}_2$, and
$\widetilde{V}(I - \widetilde{N})^{-1}\widetilde{U}
\in \mathfrak{S}_1(H)$ as a product of two
Hilbert--Schmidt operators.
To pass from the extended-space determinant to $H$,
approximate $\widetilde{N}$ in
$\mathfrak{S}_2(\mathscr{X})$ by finite-rank strictly
upper-triangular operators $\widetilde{N}^{(R)}$,
obtained by approximating each nonzero off-diagonal
block of $\widetilde{N}$ in Hilbert--Schmidt norm by
a finite-rank operator.  For each $R$, the operator
$\widetilde{N}^{(R)} + \widetilde{U}\widetilde{V}$
is trace class on $\mathscr{X}$, the nilpotent
factor satisfies
$\det_{\mathscr{X}}(I - \widetilde{N}^{(R)}) = 1$,
and Sylvester's identity gives
\[
\det_{\mathscr{X}}\bigl(I - \widetilde{N}^{(R)}
- \widetilde{U}\widetilde{V}\bigr)
= \det_H\bigl(I -
\widetilde{V}(I - \widetilde{N}^{(R)})^{-1}
\widetilde{U}\bigr).
\]
Since $\widetilde{N}^{(R)} \to \widetilde{N}$ in
operator norm, the inverses converge, and
$\widetilde{V}(I - \widetilde{N}^{(R)})^{-1}
\widetilde{U}
\to \widetilde{V}(I - \widetilde{N})^{-1}
\widetilde{U}$
in $\mathfrak{S}_1(H)$; continuity of the
Fredholm determinant on $\mathfrak{S}_1$ gives
the right-hand limit
$\det_H(I -
\widetilde{V}(I - \widetilde{N})^{-1}
\widetilde{U})$.
The left-hand side converges to the source
Fredholm determinant
$\det(I - \bar{\chi}_{\mathbf{a}} K_t
\bar{\chi}_{\mathbf{a}})$: the source
Fredholm series converges absolutely (it computes
a probability), and the conjugation $\mathscr{C}$
cancels in every finite minor determinant, so
the source series agrees term by term with the
trace-class Fredholm series of
$I - \widetilde{N}^{(R)}
- \widetilde{U}\widetilde{V}$
up to the replacement of $\widetilde{N}$ by
$\widetilde{N}^{(R)}$.
The $\mathfrak{S}_2$ convergence
$\widetilde{N}^{(R)} \to \widetilde{N}$ gives
term-by-term convergence of the minor determinants,
and the Hadamard bound applied to the source
series provides the domination for the interchange
of limit and summation.
Combining the two limits, and noting that the
conjugation weights cancel in the finite Neumann
series
($\widetilde{V}(I - \widetilde{N})^{-1}
\widetilde{U}
= V(I - L)^{-1}U$):
\begin{equation}\label{eq:RBM-det-reduction}
\det(I - \bar{\chi}_{\mathbf{a}} K_t
\bar{\chi}_{\mathbf{a}})
= \det_{L^2(\R)}\bigl(I -
V(I - L)^{-1}U\bigr).
\end{equation}

It remains to identify the entries of
$(I - L)^{-1}$ with the propagator.  The Neumann
series $(I - L)^{-1} = I + L + \cdots + L^{m-1}$
terminates at $L^{m-1}$.  Since
$L = \bar{\chi}_{\mathbf{a}} N
\bar{\chi}_{\mathbf{a}}$, the $k$-th power $L^k$
in the $(i, j)$-block ($i < j$) sums over strictly
increasing chains
$i = i_0 < i_1 < \cdots < i_k = j$ with internal
constraints $\xi_s \leq a_{i_s}$ and integrands
$\prod_{r=0}^{k-1}
N_{i_r,i_{r+1}}(\xi_r, \xi_{r+1})$.
On the region $z \leq a_i$, $w \leq a_j$ sampled
by $V$ and $U$, the endpoint cutoffs from $L$
impose no further restriction, and the chain
matches \eqref{eq:RBM-propagator-def}:
\[
[(I - L)^{-1}]_{ij}(z, w)
= [B_{\mathbf{a},\mathbf{n}}]_{ij}(z, w),
\qquad z \leq a_i,\ w \leq a_j,
\]
for $1 \leq i \leq j \leq m$, with the diagonal
block equal to the identity.
Substituting into the reduced determinant and using
the upper triangularity of $(I - L)^{-1}$:
\begin{align*}
\bigl[V(I - L)^{-1}U\bigr](x, x')
= \sum_{1 \leq j \leq i \leq m}
\int_{-\infty}^{a_j}\int_{-\infty}^{a_i}
\phi_{t,z,n_j}(x)
[B_{\mathbf{a},\mathbf{n}}]_{ji}(z, w)
\psi_{t,w,n_i}(x')\diff w\diff z.
\end{align*}
The Fredholm determinant is invariant under
transposition of the kernel; relabelling
$i \leftrightarrow j$ gives
\eqref{eq:RBM-single-kernel}.
\end{proof}

The single-space kernel
$K_{t, \mathbf{a}, \mathbf{n}}$ decomposes into diagonal
terms (one propagator endpoint coincides with the seed
evaluation point) and off-diagonal terms (both endpoints
are integrated against seeds through the propagator).
Suppressing the time parameter $t$ from the subscripts:
\begin{align}\label{eq:RBM-kernel-decomp}
K_{t, \mathbf{a}, \mathbf{n}}(x, x')
&= \sum_{i=1}^m \int_{-\infty}^{a_i}
\psi_{r, n_i}(x)\phi_{r, n_i}(x')\diff r \notag \\
&\quad + \sum_{i < j} \int_{-\infty}^{a_i}
\int_{-\infty}^{a_j}
\psi_{w, n_j}(x)
[B_{\mathbf{a},\mathbf{n}}]_{ij}(z, w)\phi_{z, n_i}(x')\diff w\diff z.
\end{align}

\subsection{Recurrences and kernel identities}
\label{sec:RBM-recurrences}

Reflected Brownian motion sits in the parabolic
framework of \S\ref{sec:dc-clean-darboux}: the wave
functions $\Psi$, $\Phi$ solve the linear problem
\eqref{eq:dc-clean-linear-d1}--\eqref{eq:dc-clean-linear-S3},
and the kernel $K$ is dressing compatible with
$(\Psi, \Phi)$ in the sense of
Definition~\ref{def:dc-clean-dressing}, so that
Theorem~\ref{thm:dc-clean-general-darboux} applies.

\begin{lemma}\label{lem:RBM-seed-linear}
Let $H = L^2(\R)$.  The seed functions $\psi_{t,a,n}$
and $\phi_{t,a,n}$ defined in
\S\ref{sec:RBM-kernel-reform} satisfy the linear problem
\eqref{eq:dc-clean-linear-d1}--\eqref{eq:dc-clean-linear-S3}
and its adjoint
\eqref{eq:dc-clean-adjoint-d1}--\eqref{eq:dc-clean-adjoint-S3}
with $\pa_1 = \pa_t$, $\pa_2 = \pa_a$,
$S_3 = e^{\pa_n}$ and constant scalar edge weights
$(c_1, \lambda_1) = (1, -\tfrac{1}{2})$,
$(c_2, \lambda_2) = (1, -1)$.
Explicitly:
\begin{align}
\psi_{t,a,n+1}
&= \pa_a \psi_{t,a,n} + \psi_{t,a,n},
\label{eq:RBM-psi-n} \\
\pa_t \psi_{t,a,n}
&= \tfrac{1}{2}\pa_a^2 \psi_{t,a,n}
   + \pa_a \psi_{t,a,n}
   + \tfrac{1}{2}\psi_{t,a,n},
\label{eq:RBM-psi-t}
\end{align}
and
\begin{align}
\phi_{t,a,n}
&= -\pa_a \phi_{t,a,n+1} + \phi_{t,a,n+1},
\label{eq:RBM-phi-n} \\
\pa_t \phi_{t,a,n}
&= -\tfrac{1}{2}\pa_a^2 \phi_{t,a,n}
   + \pa_a \phi_{t,a,n}
   - \tfrac{1}{2}\phi_{t,a,n}.
\label{eq:RBM-phi-t}
\end{align}
\end{lemma}

\begin{proof}
\noindent\textit{Verification of \eqref{eq:RBM-psi-n}.}
The function $\varphi_n$ \eqref{eq:RBM-varphi-def}
satisfies
$\pa_x \varphi_n(t,x) = -\varphi_{n+1}(t,x)$.
Since $\psi_{t,a,n}(x) = e^{x-a}\varphi_n(t, x{-}a)$,
the Leibniz rule gives
\[
\pa_a \psi_{t,a,n}(x)
= -e^{x-a}\varphi_n(t, x{-}a)
+ e^{x-a}\varphi_{n+1}(t, x{-}a)
= -\psi_{t,a,n}(x) + \psi_{t,a,n+1}(x),
\]
confirming \eqref{eq:RBM-psi-n}.

\medskip
\noindent\textit{Verification of \eqref{eq:RBM-psi-t}.}
The identity
$\pa_t \varphi_n(t,x) = \tfrac{1}{2}\varphi_{n+2}(t,x)$
gives
$\pa_t \psi_{t,a,n}
= \tfrac{1}{2}\psi_{t,a,n+2}$.
Two applications of \eqref{eq:RBM-psi-n} give
$\psi_{t,a,n+2}
= \pa_a^2 \psi_{t,a,n}
+ 2\pa_a \psi_{t,a,n}
+ \psi_{t,a,n}$,
confirming \eqref{eq:RBM-psi-t}.

\medskip
\noindent\textit{Free recurrences for
$\overline{\mathcal{S}}_{-t,k}$.}
The following identities will be applied inside the
hitting-time representation of $\phi$.
The function $\bar\varphi_n$
\eqref{eq:RBM-varphi-bar-def} satisfies the lowering
identity
$\pa_x \bar\varphi_n(t,x) = \bar\varphi_{n-1}(t,x)$.
By the chain rule,
$\pa_a \bar\varphi_n(t, x{-}a)
= -\bar\varphi_{n-1}(t, x{-}a)$.
The heat identity is
$\pa_t \bar\varphi_n(t,x)
= -\tfrac{1}{2}\bar\varphi_{n-2}(t,x)$.
Applying these to the free scattering operator
$\overline{\mathcal{S}}_{-t,k}(b, a)
= e^{a-b}\bar\varphi_{k-1}(t, b{-}a)$ gives
\begin{align}
\pa_a \overline{\mathcal{S}}_{-t,k}(b, a)
&= \overline{\mathcal{S}}_{-t,k}(b, a)
- \overline{\mathcal{S}}_{-t,k-1}(b, a),
\label{eq:RBM-Sbar-a} \\
\pa_t \overline{\mathcal{S}}_{-t,k}(b, a)
&= -\tfrac{1}{2}
\overline{\mathcal{S}}_{-t,k-2}(b, a).
\label{eq:RBM-Sbar-t}
\end{align}
Since $\bar\varphi_n = 0$ for $n < 0$,
$\overline{\mathcal{S}}_{-t,k} = 0$ for $k \leq 0$.

\medskip
\noindent\textit{Verification of \eqref{eq:RBM-phi-n}.}
Since $\tau$ and $B_\tau$ are $a$-independent,
$\pa_a$ passes through the expectation in
\eqref{eq:RBM-epigraph-def}.  Applying
\eqref{eq:RBM-Sbar-a} with $k = n - \tau$ gives
\[
\pa_a \phi_{t,a,n}
= \phi_{t,a,n} - \phi_{t,a,n-1},
\]
where $\overline{\mathcal{S}}_{-t,0} = 0$ absorbs the
$\tau = n - 1$ contribution.  Replacing $n$ by
$n + 1$ gives $\pa_a \phi_{t,a,n+1}
= \phi_{t,a,n+1} - \phi_{t,a,n}$, so
$\phi_{t,a,n}
= -\pa_a \phi_{t,a,n+1} + \phi_{t,a,n+1}$,
which is \eqref{eq:RBM-phi-n}.

\medskip
\noindent\textit{Verification of \eqref{eq:RBM-phi-t}.}
The interchange of $\pa_t$ and the expectation is
justified by dominated convergence, locally uniformly
for $t$ in compact subsets of $(0,\infty)$.  On
$\{\tau < n\}$, only finitely many values
$k = n - \tau \in \{1, \dotsc, n\}$ occur, and
\eqref{eq:RBM-Sbar-t} gives
$\pa_t \overline{\mathcal{S}}_{-t,k}(B_\tau, a)
= -\tfrac{1}{2}\overline{\mathcal{S}}_{-t,k-2}(B_\tau, a)$,
which is $e^{a-B_\tau}$ times a polynomial in $B_\tau - a$
with coefficients locally bounded in $t$.
The random walk $(B_k)_{k \geq 0}$ of
\eqref{eq:RBM-epigraph-def} satisfies
$y_n \leq B_\tau \leq B_0$ on $\{\tau < n\}$,
so the differentiated integrand is bounded in absolute
value by $C\mathbf{1}_{\tau < n}$ for a constant $C$
depending on $B_0$, $a$, $t$, and $\mathbf{y}$.
Applying \eqref{eq:RBM-Sbar-t}
with $k = n - \tau$ inside the expectation
(the boundary cases $k = 1, 2$ contribute
$\overline{\mathcal{S}}_{-t,-1}
= \overline{\mathcal{S}}_{-t,0} = 0$) gives
$\pa_t \phi_{t,a,n}
= -\tfrac{1}{2}\phi_{t,a,n-2}$.
The equivalent form $\phi_{t,a,n-1}
= (1 - \pa_a)\phi_{t,a,n}$ of \eqref{eq:RBM-phi-n},
applied twice, gives
$\phi_{t,a,n-2}
= \pa_a^2 \phi_{t,a,n}
- 2\pa_a \phi_{t,a,n}
+ \phi_{t,a,n}$,
confirming \eqref{eq:RBM-phi-t}.
\end{proof}

The propagator $[B_{\mathbf{a},\mathbf{n}}]_{ij}$ satisfies three
structural properties.

\begin{lemma}\label{lem:RBM-propagator}
The propagator satisfies the following identities.
\begin{enumerate}[label=\textup{(\roman*)},
leftmargin=*]
\item \textup{(Splitting.)}
For $i < p < j$,
\begin{equation}\label{eq:RBM-Ba-splitting}
\pa_{a_p} [B_{\mathbf{a},\mathbf{n}}]_{ij}(z, w)
= [B_{\mathbf{a},\mathbf{n}}]_{ip}(z, a_p)[B_{\mathbf{a},\mathbf{n}}]_{pj}(a_p, w).
\end{equation}

\item \textup{(Flux conservation.)}
For $i < j$,
\begin{equation}\label{eq:RBM-Ba-flux}
\Bigl(\pa_z + \pa_w
+ \sum_{i < p < j} \pa_{a_p}\Bigr)
[B_{\mathbf{a},\mathbf{n}}]_{ij}(z, w) = 0.
\end{equation}

\item The propagator depends on the particle labels
$\mathbf{n}$ only through differences
$n_{\ell'} - n_\ell$, and is $t$-independent.
\end{enumerate}
\end{lemma}

\begin{proof}
\noindent\textit{Splitting.}
In \eqref{eq:RBM-propagator-def}, consider a $k$-step
chain $i = i_0 < i_1 < \cdots < i_k = j$.  If the chain
does not pass through vertex~$p$, then $a_p$ does not
appear in any integration limit, and $\pa_{a_p}$
annihilates the term.  If the chain passes through~$p$,
say $i_s = p$, then $a_p$ appears as the upper limit of
the integral over $\xi_s$.  The fundamental theorem of
calculus evaluates the integrand at $\xi_s = a_p$,
factoring the chain product into
\[
\prod_{r=0}^{s-1}
N_{i_r, i_{r+1}}(\xi_r, \xi_{r+1})
\bigg|_{\xi_s = a_p}
\prod_{r=s}^{k-1}
N_{i_r, i_{r+1}}(\xi_r, \xi_{r+1})
\bigg|_{\xi_s = a_p},
\]
where $a_p$ serves as the upper spatial argument of
the $i$-to-$p$ piece and the lower spatial argument of
the $p$-to-$j$ piece.  The remaining integrations
separate: variables $\xi_1, \dotsc, \xi_{s-1}$ belong
to the first factor and $\xi_{s+1}, \dotsc, \xi_{k-1}$
to the second.  Summing over all chains from $i$ to $p$
and from $p$ to $j$ independently recovers
$[B_{\mathbf{a},\mathbf{n}}]_{ip}(z, a_p)
[B_{\mathbf{a},\mathbf{n}}]_{pj}(a_p, w)$.

\noindent\textit{Flux conservation.}
Fix a $k$-step chain
$i = i_0 < i_1 < \cdots < i_k = j$ in
\eqref{eq:RBM-propagator-def} and write
\[
F(\xi_0, \dotsc, \xi_k)
= \prod_{r=0}^{k-1}
N_{i_r, i_{r+1}}(\xi_r, \xi_{r+1}),
\]
with $\xi_0 = z$, $\xi_k = w$, and
\[
I = \int_{-\infty}^{a_{i_1}} \cdots
\int_{-\infty}^{a_{i_{k-1}}}
F(\xi_0, \dotsc, \xi_k)
\diff\xi_{k-1} \cdots \diff\xi_1.
\]
It suffices to show
\[
\Bigl(\pa_z + \pa_w
+ \sum_{s=1}^{k-1} \pa_{a_{i_s}}\Bigr) I = 0
\]
for each such chain.  The full identity
\eqref{eq:RBM-Ba-flux} follows by summing over all
chains and chain lengths, since $\pa_{a_p}$
annihilates any chain whose internal vertices do not
include~$p$.

Each kernel $N_{ij}(x, y)$ depends on
$x - y$ alone, so
$(\pa_x + \pa_y)N_{ij}(x, y) = 0$ in the sense of
distributions.  Applying the
Leibniz rule to each $\pa_{\xi_p} F$ and regrouping
the sum by factor, each factor
$N_{i_r, i_{r+1}}(\xi_r, \xi_{r+1})$ is
differentiated by exactly the pair
$(\pa_{\xi_r} + \pa_{\xi_{r+1}})$, giving
\[
\sum_{p=0}^{k} \pa_{\xi_p} F
= \sum_{r=0}^{k-1}
(\pa_{\xi_r} + \pa_{\xi_{r+1}})
N_{i_r, i_{r+1}}(\xi_r, \xi_{r+1})
\prod_{\substack{q=0 \\ q \neq r}}^{k-1}
N_{i_q, i_{q+1}}(\xi_q, \xi_{q+1})
= 0.
\]

Next, for each $1 \leq s \leq k-1$, the variable
$a_{i_s}$ appears only as the upper limit of the
$\xi_s$-integral in $I$.  The fundamental theorem of
calculus converts $\pa_{a_{i_s}}$ acting on this upper
limit into $\pa_{\xi_s}$ acting on the integrand (the
boundary at $\xi_s = -\infty$ vanishes by exponential
decay of $N$), giving
\[
\pa_{a_{i_s}} I
= \int_{-\infty}^{a_{i_1}} \cdots
\int_{-\infty}^{a_{i_{k-1}}}
\pa_{\xi_s} F
\diff\xi_{k-1} \cdots \diff\xi_1.
\]
The endpoint variables $\xi_0 = z$ and $\xi_k = w$ are
not integrated, so $\pa_z$ and $\pa_w$ pass inside the
integral as $\pa_{\xi_0}$ and $\pa_{\xi_k}$.
Combining,
\[
\Bigl(\pa_z + \pa_w
+ \sum_{s=1}^{k-1} \pa_{a_{i_s}}\Bigr) I
= \int_{-\infty}^{a_{i_1}} \cdots
\int_{-\infty}^{a_{i_{k-1}}}
\sum_{p=0}^{k} \pa_{\xi_p} F
\diff\xi_{k-1} \cdots \diff\xi_1
= 0,
\]
as required.

Property~(iii) is immediate from
\eqref{eq:RBM-propagator-def}: the transition kernels
$N_{ij}$ depend on particle labels only through
differences $n_j - n_i$, and no factor carries
$t$-dependence.
\end{proof}

For each $(t, \mathbf{a}, \mathbf{n})$, define the
product-graph wave functions
$\Psi^{t,\mathbf{a},\mathbf{n}}
\colon \R^m \to L^2(\R)$ and
$\Phi^{t,\mathbf{a},\mathbf{n}}
\colon L^2(\R) \to \R^m$ componentwise by
\begin{equation}\label{eq:RBM-WF-def}
\Psi_i^{t,\mathbf{a},\mathbf{n}} \defeq
\sum_{j=i}^m [\psi_{t,\mathbf{a}}]_{ij},
\qquad
\Phi_i^{t,\mathbf{a},\mathbf{n}} \defeq
\sum_{k=1}^i [\phi_{t,\mathbf{a}}]_{ki},
\end{equation}
where the propagator-dressed components are
\begin{equation}\label{eq:RBM-Psi-def}
[\psi_{t,\mathbf{a}}]_{ij} \defeq
\int_{-\infty}^{a_j}
[B_{\mathbf{a},\mathbf{n}}]_{ij}(a_i, w)\psi_{t, w, n_j}\diff w,
\quad
[\phi_{t,\mathbf{a}}]_{ij} \defeq
\int_{-\infty}^{a_i}
[B_{\mathbf{a},\mathbf{n}}]_{ij}(z, a_j)\phi_{t, z, n_i}\diff z.
\end{equation}

Under the identification
$u = (t, \mathbf{a}, \mathbf{n})$, the maps
$u \mapsto \Psi^{t,\mathbf{a},\mathbf{n}}$ and
$u \mapsto \Phi^{t,\mathbf{a},\mathbf{n}}$ become the
wave functions of the parabolic framework
(\S\ref{sec:dc-framework-clean}), with
$\pa_1 = \pa_t$,
$\pa_2 = \pa_{\mathbf{a}} \defeq \sum_{p=1}^m \pa_{a_p}$,
$S_3(t, \mathbf{a}, \mathbf{n})
= (t, \mathbf{a}, \mathbf{n}+\mathbf{1})$
and constant scalar edge weights
$C_1 = I_m$, $\Lambda_1 = -\tfrac{1}{2}I_m$,
$C_2 = I_m$, $\Lambda_2 = -I_m$.
For notational convenience, the continuous parameters
$t, \mathbf{a}$ are suppressed from the superscript
hereafter, writing $\Psi^{\mathbf{n}}$ for
$\Psi^{t,\mathbf{a},\mathbf{n}}$.  The square
$\pa_{\mathbf{a}}^2 = (\sum_{p=1}^m \pa_{a_p})^2$
is the square of the total threshold derivative.

\begin{lemma}\label{lem:RBM-WF-recurrences}
The product-graph wave functions satisfy the linear
problem
\eqref{eq:dc-clean-linear-d1}--\eqref{eq:dc-clean-linear-S3}
and its adjoint
\eqref{eq:dc-clean-adjoint-d1}--\eqref{eq:dc-clean-adjoint-S3}
with the shifts and edge weights above.
Explicitly:
\begin{align}
\pa_t \Psi^{\mathbf{n}}
&= \tfrac{1}{2}\pa_{\mathbf{a}}^2 \Psi^{\mathbf{n}}
   + \pa_{\mathbf{a}} \Psi^{\mathbf{n}}
   + \tfrac{1}{2}\Psi^{\mathbf{n}},
\label{eq:RBM-WF-t} \\
\Psi^{\mathbf{n}+\mathbf{1}}
&= \pa_{\mathbf{a}} \Psi^{\mathbf{n}} + \Psi^{\mathbf{n}},
\label{eq:RBM-WF-n}
\end{align}
and
\begin{align}
\pa_t \Phi^{\mathbf{n}}
&= -\tfrac{1}{2}\pa_{\mathbf{a}}^2 \Phi^{\mathbf{n}}
   + \pa_{\mathbf{a}} \Phi^{\mathbf{n}}
   - \tfrac{1}{2}\Phi^{\mathbf{n}},
\label{eq:RBM-WF-phi-t} \\
\Phi^{\mathbf{n}}
&= -\pa_{\mathbf{a}} \Phi^{\mathbf{n}+\mathbf{1}}
   + \Phi^{\mathbf{n}+\mathbf{1}}.
\label{eq:RBM-WF-phi-n}
\end{align}
\end{lemma}

\begin{proof}
\noindent\textit{Verification of \eqref{eq:RBM-WF-n}.}
It suffices to show
$[\psi]_{ij}^{\mathbf{n}+\mathbf{1}}
- [\psi]_{ij}^{\mathbf{n}}
= \pa_{\mathbf{a}}[\psi]_{ij}$ for each dressed component;
\eqref{eq:RBM-WF-n} then follows by summing over
$j \geq i$ via \eqref{eq:RBM-WF-def}.
The propagator depends on particle labels only through
differences (Lemma~\ref{lem:RBM-propagator}(iii)), so
the shift $\mathbf{n} \to \mathbf{n}+\mathbf{1}$ leaves
$B$ unchanged and acts only on the seeds in
\eqref{eq:RBM-Psi-def}.  Substituting the seed
recurrence \eqref{eq:RBM-psi-n} into the shifted
component replaces $\psi_{t,w,n_j+1}$ by
$\pa_w\psi_{t,w,n_j} + \psi_{t,w,n_j}$.  Splitting
the resulting integral, the $\psi_{t,w,n_j}$ term
reproduces $[\psi]_{ij}^{\mathbf{n}}$, leaving
\[
[\psi]_{ij}^{\mathbf{n}+\mathbf{1}}
- [\psi]_{ij}^{\mathbf{n}}
= \int_{-\infty}^{a_j}
[B]_{ij}(a_i, w)\pa_w
\psi_{t, w, n_j}\diff w.
\]
The target is to show that this integral equals
$\pa_{\mathbf{a}}[\psi]_{ij}$.

For the diagonal $i = j$, the propagator is
$\delta(a_i - w)$ and the integral evaluates to
$\pa_{a_i}\psi_{t,a_i,n_i}$.  Since
$[\psi]_{ii} = \psi_{t,a_i,n_i}$ depends on
$\mathbf{a}$ only through $a_i$, this equals
$\pa_{\mathbf{a}}[\psi]_{ii}$.

For $i < j$, integration by parts in $w$ gives
\[
[B]_{ij}(a_i, a_j)\psi_{t,a_j,n_j}
- \int_{-\infty}^{a_j}
\pa_w [B]_{ij}(a_i, w)\psi_{t, w, n_j}\diff w,
\]
where the boundary at $-\infty$ vanishes by
exponential decay of $N$.  Flux conservation
\eqref{eq:RBM-Ba-flux} gives
\begin{align*}
-\pa_w [B]_{ij}(a_i, w)
= \pa_z [B]_{ij}(z, w)\big|_{z = a_i}
+ \sum_{i < p < j} \pa_{a_p}[B]_{ij}(a_i, w),
\end{align*}
and the splitting identity
\eqref{eq:RBM-Ba-splitting} factorizes each
$\pa_{a_p}$ term as
$[B]_{ip}(a_i, a_p)[B]_{pj}(a_p, w)$.

On the other hand, differentiating
\eqref{eq:RBM-Psi-def} directly with respect to
$\mathbf{a}$ gives
\begin{align*}
\pa_{\mathbf{a}}[\psi]_{ij}
&= [B]_{ij}(a_i, a_j)\psi_{t,a_j,n_j}
+ \int_{-\infty}^{a_j}
  \pa_z [B]_{ij}(z, w)\big|_{z = a_i}
  \psi_{t, w, n_j}\diff w \\
&\quad + \sum_{i < p < j}
  \int_{-\infty}^{a_j}
  [B]_{ip}(a_i, a_p)[B]_{pj}(a_p, w)
  \psi_{t, w, n_j}\diff w,
\end{align*}
where $\pa_{a_j}$ acts on the upper limit by the
Leibniz integral rule, $\pa_{a_i}$ acts on the first
argument of $B$, and each $\pa_{a_p}$
($i < p < j$) acts through the splitting identity.
All other $\pa_{a_q}$ annihilate $[\psi]_{ij}$.
This expression coincides with the integration-by-parts
expansion above, so the integral equals
$\pa_{\mathbf{a}}[\psi]_{ij}$, as required.

\noindent\textit{Verification of \eqref{eq:RBM-WF-phi-n}.}
The adjoint identity follows by the same argument
applied to
$[\phi]_{ij} = \int_{-\infty}^{a_i}
[B]_{ij}(z, a_j)\phi_{t,z,n_i}\diff z$, with the
seed recurrence \eqref{eq:RBM-phi-n} in the form
$\phi_{r,n+1} - \phi_{r,n} = \pa_r \phi_{r,n+1}$
replacing \eqref{eq:RBM-psi-n} and integration by
parts in the first spatial variable $z$.  Since $a_i$
is now the upper integration limit and $a_j$ the
spatial argument of $B$, the roles of $\pa_{a_i}$ and
$\pa_{a_j}$ in the decomposition of
$\pa_{\mathbf{a}}[\phi]_{ij}^{\mathbf{n}+\mathbf{1}}$
are interchanged.

\medskip
\noindent\textit{Verification of \eqref{eq:RBM-WF-t}.}
The propagator is $t$-independent and depends on labels
only through differences
(Lemma~\ref{lem:RBM-propagator}(iii)), so $\pa_t$
passes through the integrals in \eqref{eq:RBM-Psi-def}
and acts only on the seeds.  The seed heat identity
$\pa_t \psi_{t,w,n} = \tfrac{1}{2}\psi_{t,w,n+2}$
(from the proof of \eqref{eq:RBM-psi-t}) replaces
$\psi_{t,w,n_j}$ by $\tfrac{1}{2}\psi_{t,w,n_j+2}$
in each dressed component.  The label-difference
invariance of $B$ ensures that the diagonal shift
$\mathbf{n} \to \mathbf{n}+\mathbf{2}$ leaves $B$
unchanged, so
$\pa_t [\psi]_{ij}^{\mathbf{n}}
= \tfrac{1}{2}[\psi]_{ij}^{\mathbf{n}+\mathbf{2}}$.
Summing over components,
$\pa_t \Psi^{\mathbf{n}}
= \tfrac{1}{2}\Psi^{\mathbf{n}+\mathbf{2}}$.
Two applications of \eqref{eq:RBM-WF-n} give
\[
\Psi^{\mathbf{n}+\mathbf{2}}
= (\pa_{\mathbf{a}} + 1)\Psi^{\mathbf{n}+\mathbf{1}}
= (\pa_{\mathbf{a}} + 1)^2 \Psi^{\mathbf{n}}
= \pa_{\mathbf{a}}^2 \Psi^{\mathbf{n}}
+ 2\pa_{\mathbf{a}} \Psi^{\mathbf{n}}
+ \Psi^{\mathbf{n}},
\]
confirming \eqref{eq:RBM-WF-t}.

\noindent\textit{Verification of \eqref{eq:RBM-WF-phi-t}.}
The $t$-independence and label-difference invariance
of $B$ (Lemma~\ref{lem:RBM-propagator}(iii)),
together with the seed identity
$\pa_t \phi_{t,w,n}
= -\tfrac{1}{2}\phi_{t,w,n-2}$
(from the proof of \eqref{eq:RBM-phi-t}), give
$\pa_t \Phi^{\mathbf{n}}
= -\tfrac{1}{2}\Phi^{\mathbf{n}-\mathbf{2}}$.
Rewriting \eqref{eq:RBM-WF-phi-n} as
$\Phi^{\mathbf{n}-\mathbf{1}}
= (1 - \pa_{\mathbf{a}})\Phi^{\mathbf{n}}$
and applying twice,
\[
\Phi^{\mathbf{n}-\mathbf{2}}
= (1 - \pa_{\mathbf{a}})^2 \Phi^{\mathbf{n}}
= \pa_{\mathbf{a}}^2 \Phi^{\mathbf{n}}
- 2\pa_{\mathbf{a}} \Phi^{\mathbf{n}}
+ \Phi^{\mathbf{n}},
\]
confirming \eqref{eq:RBM-WF-phi-t}.
\end{proof}

Dressing compatibility
(Definition~\ref{def:dc-clean-dressing}) requires that the
derivatives of $K$ in each direction factorize through
$\Psi$ and $\Phi$.

\begin{lemma}\label{lem:RBM-dressing}
The kernel $K_{t,\mathbf{a},\mathbf{n}}$ is dressing
compatible
(Definition~\ref{def:dc-clean-dressing}) with
$(\Psi, \Phi)$ and the edge weights above.
Explicitly:
\begin{align}
\pa_{\mathbf{a}} K_{\mathbf{n}}
&= \Psi^{\mathbf{n}} \Phi^{\mathbf{n}},
\label{eq:RBM-DC-a} \\
\pa_t K_{\mathbf{n}}
&= \tfrac{1}{2}\bigl(
\pa_{\mathbf{a}}\Psi^{\mathbf{n}} \Phi^{\mathbf{n}}
- \Psi^{\mathbf{n}}
  \pa_{\mathbf{a}}\Phi^{\mathbf{n}}\bigr)
+ \Psi^{\mathbf{n}}\Phi^{\mathbf{n}},
\label{eq:RBM-DC-t} \\
K_{\mathbf{n}+\mathbf{1}}
- K_{\mathbf{n}}
&= \Psi^{\mathbf{n}} \Phi^{\mathbf{n}+\mathbf{1}}.
\label{eq:RBM-DC-n}
\end{align}
\end{lemma}

\begin{proof}
\noindent\textit{Spatial derivative \eqref{eq:RBM-DC-a}.}
From the decomposition \eqref{eq:RBM-kernel-decomp},
differentiating by $\pa_{a_p}$ produces four terms:
the diagonal boundary term
$\psi_{a_p, n_p}\phi_{a_p, n_p}$,
the off-diagonal boundary terms where $\pa_{a_p}$ acts
on an upper integration limit $a_i$ or $a_j$ (producing
propagator evaluations at $a_p$), and the interior term
where
$\pa_{a_p}$ hits $[B_{\mathbf{a},\mathbf{n}}]_{ij}$ through
the splitting identity \eqref{eq:RBM-Ba-splitting}.
Collecting all four contributions:
\begin{align*}
\pa_{a_p} K_{\mathbf{n}}
&= \psi_{a_p, n_p}\phi_{a_p, n_p}
+ \sum_{p < j} [\psi_{\mathbf{a}}]_{pj}
\phi_{a_p, n_p}
+ \sum_{i < p} \psi_{a_p, n_p}
[\phi_{\mathbf{a}}]_{ip}
+ \sum_{i < p < j}
[\psi_{\mathbf{a}}]_{pj}
[\phi_{\mathbf{a}}]_{ip} \\
&= \sum_{j \geq p} [\psi_{\mathbf{a}}]_{pj}
\sum_{i \leq p} [\phi_{\mathbf{a}}]_{ip}
= \Psi_p^{\mathbf{n}}\Phi_p^{\mathbf{n}}.
\end{align*}
Summing over $p$ gives \eqref{eq:RBM-DC-a}.

\medskip
\noindent\textit{Discrete shift \eqref{eq:RBM-DC-n}.}
By Lemma~\ref{lem:RBM-propagator}(iii), the propagator
depends on particle labels only through differences, so
$B_{\mathbf{n}+\mathbf{1}} = B_{\mathbf{n}}$ and the
shift $\mathbf{n} \to \mathbf{n}+\mathbf{1}$ acts only
on the seeds in \eqref{eq:RBM-kernel-decomp}.  The seed
recurrences
\eqref{eq:RBM-psi-n}--\eqref{eq:RBM-phi-n} give
$\psi_{r,n+1} - \psi_{r,n} = \pa_r \psi_{r, n}$
and
$\phi_{r,n+1} - \phi_{r,n} = \pa_r \phi_{r, n+1}$,
so for each diagonal term
\[
\psi_{r,n_i+1}\phi_{r,n_i+1}
- \psi_{r,n_i}\phi_{r,n_i}
= (\pa_r \psi_{r,n_i})\phi_{r,n_i+1}
+ \psi_{r,n_i}(\pa_r \phi_{r,n_i+1})
= \pa_r(\psi_{r,n_i}\phi_{r,n_i+1}).
\]
Integrating and evaluating at the boundary (the
contribution at $-\infty$ vanishes by exponential
decay),
\[
\sum_{i=1}^m \int_{-\infty}^{a_i}
\pa_r(\psi_{r, n_i}\phi_{r, n_i+1}) \diff r
= \sum_{i=1}^m
\psi_{a_i, n_i}\phi_{a_i, n_i+1}.
\]

For the off-diagonal part, the discrete Leibniz rule
and the seed recurrences
$\psi_{w,n_j+1} - \psi_{w,n_j}
= \pa_w \psi_{w,n_j}$,
$\phi_{z,n_i+1} - \phi_{z,n_i}
= \pa_z \phi_{z,n_i+1}$ give
\begin{align*}
\psi_{w,n_j+1}[B]_{ij}\phi_{z,n_i+1}
- \psi_{w,n_j}[B]_{ij}\phi_{z,n_i}
= (\pa_w\psi_{w,n_j})[B]_{ij}\phi_{z,n_i+1}
+ \psi_{w,n_j}[B]_{ij}(\pa_z\phi_{z,n_i+1}).
\end{align*}
Expanding
$(\pa_z + \pa_w
+ \sum_{i<p<j}\pa_{a_p})
[\psi_{w,n_j}[B]_{ij}\phi_{z,n_i+1}]$
by the Leibniz rule, the derivatives of $[B]_{ij}$
sum to zero by flux conservation
\eqref{eq:RBM-Ba-flux}, leaving only the
$\pa_w\psi$ and $\pa_z\phi$ terms.  The right side
above is therefore a divergence
\[
(\pa_z + \pa_w
+ \textstyle\sum_{i<p<j}\pa_{a_p})
\bigl[\psi_{w,n_j}[B]_{ij}(z,w)
\phi_{z,n_i+1}\bigr].
\]
Integrating over the chain domain, the fundamental
theorem of calculus evaluates each component at its
boundary.  The three contributions are
$[\psi_{\mathbf{a},\mathbf{n}}]_{jj}
[\phi_{\mathbf{a},\mathbf{n}+\mathbf{1}}]_{ij}$
(from $\pa_w$ at $w = a_j$),
$[\psi_{\mathbf{a},\mathbf{n}}]_{ij}
[\phi_{\mathbf{a},\mathbf{n}+\mathbf{1}}]_{ii}$
(from $\pa_z$ at $z = a_i$), and
$[\psi_{\mathbf{a},\mathbf{n}}]_{pj}
[\phi_{\mathbf{a},\mathbf{n}+\mathbf{1}}]_{ip}$ for
$i < p < j$
(from splitting \eqref{eq:RBM-Ba-splitting}),
totalling
$\sum_{p=i}^{j}
[\psi_{\mathbf{a},\mathbf{n}}]_{pj}
[\phi_{\mathbf{a},\mathbf{n}+\mathbf{1}}]_{ip}$.
Collecting diagonal and off-diagonal contributions,
\[
K_{\mathbf{n}+\mathbf{1}} - K_{\mathbf{n}}
= \sum_{i \leq p \leq j}
[\psi_{\mathbf{a},\mathbf{n}}]_{pj}
[\phi_{\mathbf{a},\mathbf{n}+\mathbf{1}}]_{ip}
= \Psi^{\mathbf{n}} \Phi^{\mathbf{n}+\mathbf{1}},
\]
confirming \eqref{eq:RBM-DC-n}.

\medskip
\noindent\textit{Time derivative \eqref{eq:RBM-DC-t}.}
The propagator is $t$-independent
(Lemma~\ref{lem:RBM-propagator}(iii)), so $\pa_t$
acts only on the seeds in
\eqref{eq:RBM-kernel-decomp}.  The seed heat
identities
$\pa_t \psi_{r,n} = \tfrac{1}{2}\psi_{r,n+2}$ and
$\pa_t \phi_{r,n} = -\tfrac{1}{2}\phi_{r,n-2}$
(from \eqref{eq:RBM-psi-t}--\eqref{eq:RBM-phi-t})
give
$\pa_t(\psi_{r,n}\phi_{r,n})
= \tfrac{1}{2}\psi_{r,n+2}\phi_{r,n}
- \tfrac{1}{2}\psi_{r,n}\phi_{r,n-2}$.
Expanding $\psi_{n+2} = (\pa_r + 1)^2 \psi_n$ and
$\phi_{n-2} = (\pa_r - 1)^2 \phi_n$ via two
applications of the seed recurrences, the constant
terms $\psi_n\phi_n$ cancel and the remainder is a
total $\pa_r$-derivative
\[
\pa_t(\psi_{r,n}\phi_{r,n})
= \tfrac{1}{2}\pa_r\bigl[
(\pa_r \psi_{r,n})\phi_{r,n}
+ 2\psi_{r,n}\phi_{r,n}
- \psi_{r,n}(\pa_r \phi_{r,n})\bigr].
\]
Integrating over $(-\infty, a_i]$ and substituting
$\pa_r \psi_{r,n} = \psi_{r,n+1} - \psi_{r,n}$ and
$\pa_r \phi_{r,n} = \phi_{r,n} - \phi_{r,n-1}$ at the
boundary $r = a_i$, the bracket evaluates to
$\psi_{a_i,n_i+1}\phi_{a_i,n_i}
+ \psi_{a_i,n_i}\phi_{a_i,n_i-1}$, giving
\begin{equation}\label{eq:RBM-DC-t-diag}
\pa_t K^{\mathrm{diag}}
= \frac{1}{2}\sum_{i=1}^m \bigl(
\psi_{a_i, n_i+1}\phi_{a_i, n_i}
+ \psi_{a_i, n_i}\phi_{a_i, n_i-1}\bigr).
\end{equation}

For the off-diagonal part, $\pa_t$ acts on the seeds
inside each double integral.  The seed recurrences
give
$\pa_t \psi_{w,n_j}
= \tfrac{1}{2}(\pa_w + 1)^2 \psi_{w,n_j}$ and
$\pa_t \phi_{z,n_i}
= -\tfrac{1}{2}(\pa_z - 1)^2 \phi_{z,n_i}$.
Expanding the squares, the constant terms cancel,
leaving
\begin{align*}
\pa_t\bigl[\psi_{w,n_j}[B]_{ij}\phi_{z,n_i}\bigr]
= \tfrac{1}{2}\bigl[
(\pa_w^2 + 2\pa_w)\psi_{w,n_j}
[B]_{ij}\phi_{z,n_i}
- \psi_{w,n_j}[B]_{ij}
(\pa_z^2 - 2\pa_z)\phi_{z,n_i}\bigr].
\end{align*}
The expanded expression equals
\begin{align*}
\tfrac{1}{2}\bigl(\pa_z + \pa_w
+ \textstyle\sum_{i<p<j}\pa_{a_p}\bigr)
\bigl[(\pa_w \psi_{w,n_j})[B]_{ij}\phi_{z,n_i}
+ 2\psi_{w,n_j}[B]_{ij}\phi_{z,n_i}
- \psi_{w,n_j}[B]_{ij}(\pa_z\phi_{z,n_i})\bigr],
\end{align*}
which can be verified by expanding the divergence:
the Leibniz rule distributes it across $\psi$,
$[B]_{ij}$, and $\phi$, and flux conservation
\eqref{eq:RBM-Ba-flux} eliminates all derivatives
of $[B]_{ij}$.  Splitting
$2\psi_{w,n_j}[B]_{ij}\phi_{z,n_i}$ and applying
the seed identities
$\pa_w\psi_{w,n_j} + \psi_{w,n_j} = \psi_{w,n_j+1}$
and
$\phi_{z,n_i} - \pa_z\phi_{z,n_i} = \phi_{z,n_i-1}$
simplifies the bracket to
\[
\tfrac{1}{2}\bigl(\pa_z + \pa_w
+ \textstyle\sum_{i<p<j}\pa_{a_p}\bigr)
\bigl[\psi_{w,n_j+1}[B]_{ij}\phi_{z,n_i}
+ \psi_{w,n_j}[B]_{ij}\phi_{z,n_i-1}\bigr].
\]
Integrating the divergence over the chain domain,
the fundamental theorem of calculus evaluates each
component at its boundary: $\pa_w$ at $w = a_j$,
$\pa_z$ at $z = a_i$, and each $\pa_{a_p}$ via
splitting \eqref{eq:RBM-Ba-splitting}.  The first
term in the bracket produces
$\sum_{p=i}^{j}
[\psi_{\mathbf{a},\mathbf{n}+\mathbf{1}}]_{pj}
[\phi_{\mathbf{a},\mathbf{n}}]_{ip}$ and the second
produces
$\sum_{p=i}^{j}
[\psi_{\mathbf{a},\mathbf{n}}]_{pj}
[\phi_{\mathbf{a},\mathbf{n}-\mathbf{1}}]_{ip}$.
Collecting diagonal and off-diagonal contributions,
\begin{align*}
\pa_t K_{\mathbf{n}}
= \frac{1}{2}\sum_{i \leq p \leq j}
[\psi_{\mathbf{a},\mathbf{n}+\mathbf{1}}]_{pj}
[\phi_{\mathbf{a},\mathbf{n}}]_{ip}
+ [\psi_{\mathbf{a},\mathbf{n}}]_{pj}
[\phi_{\mathbf{a},\mathbf{n}-\mathbf{1}}]_{ip}
= \frac{1}{2}\bigl(\Psi^{\mathbf{n}+\mathbf{1}}
\Phi^{\mathbf{n}}
+ \Psi^{\mathbf{n}}
\Phi^{\mathbf{n}-\mathbf{1}}\bigr).
\end{align*}
Substituting
$\Psi^{\mathbf{n}+\mathbf{1}}
= \pa_{\mathbf{a}}\Psi^{\mathbf{n}}
+ \Psi^{\mathbf{n}}$ and
$\Phi^{\mathbf{n}-\mathbf{1}}
= \Phi^{\mathbf{n}}
- \pa_{\mathbf{a}}\Phi^{\mathbf{n}}$,
this equals
$\tfrac{1}{2}(\pa_{\mathbf{a}}\Psi^{\mathbf{n}}
\Phi^{\mathbf{n}}
- \Psi^{\mathbf{n}}
\pa_{\mathbf{a}}\Phi^{\mathbf{n}})
+ \Psi^{\mathbf{n}}\Phi^{\mathbf{n}}$,
confirming \eqref{eq:RBM-DC-t}.

The identities
\eqref{eq:RBM-DC-a}--\eqref{eq:RBM-DC-n} hold in
$\mathfrak{S}_1(L^2(\R))$: the seed $L^2$ bounds
from Lemma~\ref{lem:RBM-kernel-reform} control the
differentiated integrands in
\eqref{eq:RBM-kernel-decomp}, justifying the
interchange of derivatives and integrals.

It remains to verify resolvent existence.  By
Lemma~\ref{lem:RBM-kernel-reform},
$\det(I - K_{t,\mathbf{a},\mathbf{n}})
= \PP_{\mathbf{y}}\bigl(\bigcap_{i=1}^m
\{Y_{n_i}(t) > a_i\}\bigr)$;
since $K_{t,\mathbf{a},\mathbf{n}}$ is trace class,
the resolvent exists if and only if this determinant
is nonzero.  The probability is strictly positive
for $t > 0$, so the determinant is nonzero.
\end{proof}

\subsection{Multipoint equation}\label{sec:RBM-multipoint}

The dressed observable
$\mathcal{A}(u) = z\Phi(u) R(u) \Psi(u)$ of
Proposition~\ref{prop:dc-clean-dressed-waves} specializes at $z = 1$
to the $m \times m$ matrix
\begin{equation}\label{eq:RBM-A-def}
\mathcal{A}_{\mathbf{n}} \defeq
\Phi^{\mathbf{n}} R_{\mathbf{n}} \Psi^{\mathbf{n}}
\in \End(\R^m),
\qquad
(\mathcal{A}_{\mathbf{n}})_{ij}
= \langle R_{\mathbf{n}} \Psi_j^{\mathbf{n}},
\Phi_i^{\mathbf{n}} \rangle.
\end{equation}

\begin{theorem}\label{thm:RBM-multipoint}
The dressed observable $\mathcal{A}_{\mathbf{n}}$ satisfies
\begin{align}\label{eq:RBM-multipoint}
\pa_t (\mathcal{A}_{\mathbf{n}+\mathbf{1}}
- \mathcal{A}_{\mathbf{n}})
- \tfrac{1}{2}\pa_{\mathbf{a}}^2(\mathcal{A}_{\mathbf{n}+\mathbf{1}}
+ \mathcal{A}_{\mathbf{n}})
- \pa_{\mathbf{a}}(\mathcal{A}_{\mathbf{n}}
  \mathcal{A}_{\mathbf{n}+\mathbf{1}})
+ \mathcal{A}_{\mathbf{n}+\mathbf{1}}
  \pa_{\mathbf{a}} \mathcal{A}_{\mathbf{n}+\mathbf{1}}
+ \pa_{\mathbf{a}} \mathcal{A}_{\mathbf{n}}
  \mathcal{A}_{\mathbf{n}} = 0.
\end{align}
\end{theorem}

\begin{proof}
Theorem~\ref{thm:dc-clean-general-darboux} applies: the
dressed $\Lambda$-weights
$\widehat\Lambda_1 = -\tfrac{1}{2}I + \pa_{\mathbf{a}} \mathcal{A}_{\mathbf{n}}$,
$\widehat\Lambda_2 = -I + \mathcal{A}_{\mathbf{n}}
- \mathcal{A}_{\mathbf{n}+\mathbf{1}}$
satisfy the parabolic $\Lambda$-diamond equation
\eqref{eq:dc-clean-L-diamond}, where
$\mathcal{H} = \pa_t - \tfrac{1}{2}\pa_{\mathbf{a}}^2$ is the heat
operator of \eqref{eq:dc-clean-L-diamond}.  With
$C_1 = C_2 = I$ (constant scalar), the $\Lambda$-diamond
reduces to
\begin{equation}\label{eq:RBM-dressed-L-diamond}
\mathcal{H}\widehat\Lambda_2
+ \pa_{\mathbf{a}}\widehat\Lambda_1
- \pa_{\mathbf{a}}\widehat\Lambda_2
- [\widehat\Lambda_1,
\widehat\Lambda_2]_{S_3} = 0,
\end{equation}
where $[A, B]_{S_3} = AB - B(S_3 A)$ is the
shifted commutator of
\eqref{eq:dc-clean-L-diamond}, since $C_1 = C_2 = I$ acts as the identity under
right-multiplication in the terms
$\pa_{\mathbf{a}}\widehat\Lambda_1 C_2$ and
$-\pa_{\mathbf{a}}\widehat\Lambda_2 C_1(S_3 u)$ of
\eqref{eq:dc-clean-L-diamond}.

Substituting the dressed weights, the constant
parts contribute
$\mathcal{H}(-I) + \pa_{\mathbf{a}}(-\tfrac{1}{2}I)
- \pa_{\mathbf{a}}(-I) - [-\tfrac{1}{2}I, -I]_{S_3} = 0$
(each term vanishes individually).
The $\mathcal{A}$-dependent parts give
\[
\mathcal{H}(\mathcal{A}_{\mathbf{n}} - \mathcal{A}_{\mathbf{n}+\mathbf{1}})
+ \pa_{\mathbf{a}}^2 \mathcal{A}_{\mathbf{n}}
- \pa_{\mathbf{a}} \mathcal{A}_{\mathbf{n}} (\mathcal{A}_{\mathbf{n}}
- \mathcal{A}_{\mathbf{n}+\mathbf{1}})
+ (\mathcal{A}_{\mathbf{n}} - \mathcal{A}_{\mathbf{n}+\mathbf{1}})
\pa_{\mathbf{a}} \mathcal{A}_{\mathbf{n}+\mathbf{1}} = 0.
\]
Expanding $\mathcal{H}(\mathcal{A}_{\mathbf{n}}
- \mathcal{A}_{\mathbf{n}+\mathbf{1}})
= (\pa_t - \tfrac{1}{2}\pa_{\mathbf{a}}^2)(\mathcal{A}_{\mathbf{n}}
- \mathcal{A}_{\mathbf{n}+\mathbf{1}})$ and collecting:
\begin{align*}
-\pa_t(\mathcal{A}_{\mathbf{n}+\mathbf{1}}
- \mathcal{A}_{\mathbf{n}})
+ \tfrac{1}{2}\pa_{\mathbf{a}}^2(\mathcal{A}_{\mathbf{n}}
+ \mathcal{A}_{\mathbf{n}+\mathbf{1}})
+ \pa_{\mathbf{a}} \mathcal{A}_{\mathbf{n}}
(\mathcal{A}_{\mathbf{n}+\mathbf{1}}
- \mathcal{A}_{\mathbf{n}})
- (\mathcal{A}_{\mathbf{n}+\mathbf{1}}
- \mathcal{A}_{\mathbf{n}})
\pa_{\mathbf{a}} \mathcal{A}_{\mathbf{n}+\mathbf{1}} = 0.
\end{align*}
Multiplying by $-1$ and distributing the products
$(\pa_{\mathbf{a}} \mathcal{A}_{\mathbf{n}})
(\mathcal{A}_{\mathbf{n}+\mathbf{1}}
- \mathcal{A}_{\mathbf{n}})$ and
$(\mathcal{A}_{\mathbf{n}+\mathbf{1}}
- \mathcal{A}_{\mathbf{n}})
\pa_{\mathbf{a}} \mathcal{A}_{\mathbf{n}+\mathbf{1}}$,
the Leibniz rule recombines the cross-terms as
$\pa_{\mathbf{a}}(\mathcal{A}_{\mathbf{n}}
\mathcal{A}_{\mathbf{n}+\mathbf{1}})$,
giving \eqref{eq:RBM-multipoint}.
\end{proof}

\begin{corollary}
\label{cor:RBM-trace-defect}
Let $F_{\mathbf{n}} = \det(I - K_{\mathbf{n}})$.
Then
$\pa_{\mathbf{a}} \log F_{\mathbf{n}}
= -\operatorname{tr}_E(\mathcal{A}_{\mathbf{n}})$,
and
\begin{equation}\label{eq:RBM-Hirota}
\bigl[D_t - \tfrac{1}{2}D_{\mathbf{a}}^2\bigr]
F_{\mathbf{n}+\mathbf{1}} \cdot F_{\mathbf{n}}
= \tfrac{1}{2}F_{\mathbf{n}+\mathbf{1}}
F_{\mathbf{n}}
\bigl(\operatorname{tr}_E((\mathcal{A}_{\mathbf{n}+\mathbf{1}}
- \mathcal{A}_{\mathbf{n}})^2)
- (\operatorname{tr}_E(\mathcal{A}_{\mathbf{n}+\mathbf{1}}
- \mathcal{A}_{\mathbf{n}}))^2\bigr).
\end{equation}
\end{corollary}

\begin{proof}
Proposition~\ref{prop:dc-clean-normalized-trace-defect-direct}
with $C_1 = C_2 = I$ and $\Lambda_2 = -I$ gives
$C_1 + \Lambda_2 C_2^{-1} = I - I = 0$, so
Corollary~\ref{cor:dc-clean-scalar-trace-defect} applies
with $a(u) = 0$.  The quantity
$\mathcal{B} = \Delta_2 C_2^{-1}$
(where $\Delta_2 = \widehat\Lambda_2 - \Lambda_2$)
of Proposition~\ref{prop:dc-clean-normalized-trace-defect-direct} reduces to
$\mathcal{B} = \mathcal{A}_{\mathbf{n}}
- \mathcal{A}_{\mathbf{n}+\mathbf{1}}$, and
\eqref{eq:dc-clean-normalized-scalar-linear-trace-defect-direct}
gives \eqref{eq:RBM-Hirota}.
\end{proof}

\begin{corollary}
\label{cor:RBM-one-point}
For $m = 1$, the Fredholm determinant
$F_{t,a,n} = \det(I - K_{t,a,n})$ satisfies
\begin{equation}\label{eq:RBM-scalar}
\bigl[D_t - \tfrac{1}{2}D_a^2\bigr]
F_{t,a,n+1} \cdot F_{t,a,n} = 0.
\end{equation}
\end{corollary}

\begin{proof}
When $m = 1$, $\mathcal{A}_{\mathbf{n}}$ is a $1 \times 1$
matrix (a scalar).  The trace-defect term
$\operatorname{tr}_E((\mathcal{A}_{n+1}
- \mathcal{A}_n)^2)
- (\operatorname{tr}_E(\mathcal{A}_{n+1}
- \mathcal{A}_n))^2$ vanishes, and
\eqref{eq:RBM-Hirota} reduces to
\eqref{eq:RBM-scalar}.
\end{proof}

\begin{remark}
\label{rem:RBM-alternative}
Equation \eqref{eq:RBM-multipoint} may equivalently be
written as
\begin{align}\label{eq:RBM-multipoint-alt}
\pa_t (\mathcal{A}_{\mathbf{n}+\mathbf{1}}
- \mathcal{A}_{\mathbf{n}})
&- \tfrac{1}{2}\pa_{\mathbf{a}}^2
(\mathcal{A}_{\mathbf{n}+\mathbf{1}}
+ \mathcal{A}_{\mathbf{n}})
+ \tfrac{1}{2}\pa_{\mathbf{a}}
(\mathcal{A}_{\mathbf{n}+\mathbf{1}}
- \mathcal{A}_{\mathbf{n}})^2 \notag \\
&+ \tfrac{1}{2}[\mathcal{A}_{\mathbf{n}+\mathbf{1}}
- \mathcal{A}_{\mathbf{n}},
\pa_{\mathbf{a}}(\mathcal{A}_{\mathbf{n}+\mathbf{1}}
+ \mathcal{A}_{\mathbf{n}})] = 0.
\end{align}
\end{remark}

\section{Brownian Last Passage Percolation with Boundary}
\label{sec:BLPP}

\paragraph{\textbf{System description.}}
Let $b \colon [0, \infty) \to \R$ be continuous with $b(0) = 0$,
and let $W_1, W_2, \dotsc$ be independent standard Brownian motions.
Brownian last passage percolation (BLPP) with boundary~$b$ is the
corner-growth process
\[
G(b; t, n)
= \max_{0 \leq t_0 \leq t_1 \leq \cdots \leq t_n = t}
\Bigl\{b(t_0) + \sum_{i=1}^n
\bigl(W_i(t_i) - W_i(t_{i-1})\bigr)\Bigr\}
\]
for $t \geq 0$ and $n \geq 0$.

\paragraph{\textbf{Fredholm determinant formula.}}
Fix $n \geq 1$, times $0 < t_1 < t_2 < \cdots < t_k$,
and thresholds $a_1, \dotsc, a_k \in \R$.  The multipoint
distribution is
\[
\PP\bigl(G(b; t_i, n) \leq a_i;\ 1 \leq i \leq k\bigr)
= \det(I - \chi_{\mathbf{a}} K
\chi_{\mathbf{a}})_{L^2(\{t_1, \dotsc, t_k\} \times \R)},
\]
where $\chi_{\mathbf{a}}(t_i, z) = \mathbf{1}_{z \geq a_i}$.
The correlation kernel has block entries
\[
K(t_i, u; t_j, v)
= -e^{\frac{(t_j - t_i)}{2}\pa^2}(u, v)
\mathbf{1}_{t_i < t_j}
+ \mathcal{S}_{t_i, n}(u, \cdot)\mathcal{S}_{t_j, n}^{\operatorname{epi}(b)}(\cdot, v),
\]
where $e^{\tfrac{(t_j-t_i)}{2}\pa^2}$ is the heat kernel.
The scattering operators take the form
\begin{align}
\mathcal{S}_{t, n}(u, v) &\defeq
\bar{\varphi}_{n-1}(t, u - v),
\label{eq:BLPP-scattering-def} \\
\mathcal{S}_{t, n}^{\operatorname{epi}(b)}(u, v) &\defeq
(-1)^n \E_{B_0 = u}\bigl[
\varphi_n(t - \tau, -B(\tau) - v)
\mathbf{1}_{\tau \leq t}\bigr],
\label{eq:BLPP-epigraph-def}
\end{align}
where $B(t)$ is a standard Brownian motion and
$\tau = \inf\{s \geq 0 : B(s) \geq b(s)\}$.

\paragraph{\textbf{References.}}
The Fredholm determinant formula above is due to
Rahman~\cite{Rahman2025}.

\begin{remark}\label{rem:BLPP-convention}
Rahman~\cite{Rahman2025} proves the BLPP Fredholm
determinant formula using a hypograph operator
$S^{\mathrm{hypo}(b)}$ built from the lower
hitting time
$\tau = \inf\{s \geq 0 : B(s) \leq b(s)\}$.
The epigraph kernel
\eqref{eq:BLPP-epigraph-def}, with upper hitting
time $\tau = \inf\{s \geq 0 : B(s) \geq b(s)\}$
and the sign $-B(\tau)$ in the scattering operator,
is an equivalent representation: the two formulations
compute the same Fredholm determinant for the same
boundary~$b$.  The equivalence follows from
Brownian reflection, Hermite parity
$(-1)^n\varphi_n(t, -x) = \varphi_n(t, x)$, and a
Sylvester rearrangement of the determinant;
see~\cite[Lemma~3.6]{Rodriguez2025} for the
one-point derivation.
\end{remark}

\subsection{Kernel reformulation}\label{sec:BLPP-kernel-reform}

The reduction to a single-space kernel follows the same
strategy as for RBM, with one structural modification: the
right-tail cutoff $\chi_{a_i}(z) = \mathbf{1}_{z \geq a_i}$
reverses the integration direction.

The seed functions are
\begin{equation}\label{eq:BLPP-seed-def}
\psi_{t_i, z, n}(u) \defeq
\mathcal{S}_{t_i, n}^{\operatorname{epi}(b)}(u, z),
\qquad
\phi_{t_j, w, n}(v) \defeq
\mathcal{S}_{t_j, n}(w, v).
\end{equation}
The transition matrix
$N_{ij}(u, v) = -e^{\tfrac{(t_j - t_i)}{2}\pa^2}(u, v)
\mathbf{1}_{t_i < t_j}$ defines the strictly
upper-triangular block operator
$L \defeq \chi_{\mathbf{a}} N \chi_{\mathbf{a}}$.

The BLPP propagator integrates over $[a_p, \infty)$
at intermediate vertices.  For $1 \leq i \leq j \leq k$,
\begin{align}\label{eq:BLPP-propagator-def}
&[B_{\mathbf{t},\mathbf{a}}]_{ij}(z, w) \defeq \notag \\
&\begin{cases}
\displaystyle\sum_{\ell=1}^{j-i}
\sum_{\substack{i = i_0 < i_1 < \cdots < i_\ell = j}}
\int_{a_{i_1}}^{\infty} \cdots
\int_{a_{i_{\ell-1}}}^{\infty}
\prod_{r=0}^{\ell-1}
N_{i_r, i_{r+1}}(\xi_r, \xi_{r+1})
\diff\xi_{\ell-1} \cdots \diff\xi_1,
& i < j, \\[6pt]
\delta(z - w), & i = j, \\
0, & i > j,
\end{cases}
\end{align}
with $\xi_0 = z$ and $\xi_\ell = w$.

\begin{lemma}\label{lem:BLPP-kernel-reform}
The multipoint distribution satisfies
\[
\det(I - \chi_{\mathbf{a}} K
\chi_{\mathbf{a}})_{L^2(\{t_1, \dotsc, t_k\}
\times \R)}
= \det(I - K_{\mathbf{t}, \mathbf{a}, n})_{L^2(\R)},
\]
where
\begin{equation}\label{eq:BLPP-single-kernel}
K_{\mathbf{t}, \mathbf{a}, n}
= \sum_{1 \leq i \leq j \leq k}
\int_{a_i}^{\infty} \int_{a_j}^{\infty}
\psi_{t_i, z, n}
[B_{\mathbf{t},\mathbf{a}}]_{ij}(z, w)\phi_{t_j, w, n}\diff w\diff z.
\end{equation}
\end{lemma}

\begin{proof}
The argument is completely analogous to
Lemma~\ref{lem:RBM-kernel-reform}, with the
integration direction reversed: $(-\infty, a_p]$
becomes $[a_p, \infty)$, and the Neumann series for
$(I - L)^{-1}$ terminates after $k - 1$ terms since
$L$ is strictly upper triangular in the $k$ time
indices.
\end{proof}

The single-space kernel $K_{\mathbf{t}, \mathbf{a}, n}$
decomposes into diagonal and off-diagonal contributions:
\begin{equation}\label{eq:BLPP-kernel-decomp}
K_{\mathbf{t}, \mathbf{a}, n}
= \sum_{i=1}^k \int_{a_i}^{\infty}
\psi_{t_i, r, n}\phi_{t_i, r, n}\diff r
+ \sum_{i < j} \int_{a_i}^{\infty}
\int_{a_j}^{\infty}
\psi_{t_i, z, n}
[B_{\mathbf{t},\mathbf{a}}]_{ij}(z, w)\phi_{t_j, w, n}\diff w\diff z.
\end{equation}

\begin{remark}\label{rem:BLPP-analytic}
The BLPP seeds lack the exponential prefactor of
the RBM seeds, so the Fredholm determinant
$\det(I - K_n)$ and the dressing identities
\eqref{eq:BLPP-DC-a}--\eqref{eq:BLPP-DC-n} are
understood through the absolutely convergent
Fredholm series.  The required kernel decay follows
from the Gaussian factor in the heat semigroup and
the stopping-time bounds in
Rahman~\cite{Rahman2025}.
\end{remark}

\subsection{Recurrences and kernel identities}
\label{sec:BLPP-recurrences}

Brownian last passage percolation sits in the parabolic
framework of \S\ref{sec:dc-clean-darboux}: the wave
functions $\Psi$, $\Phi$ solve the linear problem
\eqref{eq:dc-clean-linear-d1}--\eqref{eq:dc-clean-linear-S3},
and the kernel $K$ is dressing compatible with
$(\Psi, \Phi)$ in the sense of
Definition~\ref{def:dc-clean-dressing}, so that
Theorem~\ref{thm:dc-clean-general-darboux} applies.
Under the identification
$u = (\mathbf{t}, \mathbf{a}, n)$, the shifts are
$\pa_1 = \pa_{\mathbf{t}}
\defeq \sum_{i=1}^k \pa_{t_i}$,
$\pa_2 = -\pa_{\mathbf{a}}$
(where
$\pa_{\mathbf{a}} \defeq \sum_{i=1}^k \pa_{a_i}$),
$S_3 = e^{\pa_n}$, with constant scalar edge weights
$C_1 = 0$, $\Lambda_1 = 0$,
$C_2 = I$, $\Lambda_2 = 0$.
The sign in $\pa_2 = -\pa_{\mathbf{a}}$ absorbs the
right-tail cutoff convention: differentiating the lower
limit $\int_{a_p}^{\infty}$ produces a sign opposite to
the upper limit $\int_{-\infty}^{a_p}$ of RBM.

\begin{lemma}\label{lem:BLPP-seed-linear}
Let $H = L^2(\R)$.  The seed functions $\psi_{t,z,n}$
and $\phi_{t,w,n}$ defined in
\S\ref{sec:BLPP-kernel-reform} satisfy the linear problem
\eqref{eq:dc-clean-linear-d1}--\eqref{eq:dc-clean-linear-S3}
and its adjoint
\eqref{eq:dc-clean-adjoint-d1}--\eqref{eq:dc-clean-adjoint-S3}
with $\pa_1 = \pa_t$,
$\pa_2 = -\pa_a$,
$S_3 = e^{\pa_n}$ and constant scalar edge weights
$C_1 = 0$, $\Lambda_1 = 0$,
$C_2 = I$, $\Lambda_2 = 0$.
Explicitly:
\begin{align}
\pa_z \psi_{t,z,n} &= -\psi_{t,z,n+1},
\label{eq:BLPP-psi-a} \\
\pa_t \psi_{t,z,n}
&= \tfrac{1}{2}\pa_z^2 \psi_{t,z,n},
\label{eq:BLPP-psi-t}
\end{align}
and
\begin{align}
\pa_w \phi_{t,w,n} &= \phi_{t,w,n-1},
\label{eq:BLPP-phi-a} \\
\pa_t \phi_{t,w,n}
&= -\tfrac{1}{2}\pa_w^2 \phi_{t,w,n}.
\label{eq:BLPP-phi-t}
\end{align}
\end{lemma}

\begin{proof}
The RBM seeds carry the exponential prefactor
$e^{x-a}$, which couples the Hermite raising/lowering
identities to a constant shift term
(Lemma~\ref{lem:RBM-seed-linear}).
The BLPP scattering operators omit this prefactor:
$\mathcal{S}_{t,n}(u,v) = \bar\varphi_{n-1}(t, u{-}v)$
and
$\mathcal{S}_{t,n}^{\operatorname{epi}(b)}$
involves $\varphi_n$ without conjugation.
The spatial recurrences \eqref{eq:BLPP-psi-a}--\eqref{eq:BLPP-phi-a} therefore reduce to
the Hermite raising identity
$\pa_x \varphi_n = -\varphi_{n+1}$ and lowering
identity $\pa_x \bar\varphi_n = \bar\varphi_{n-1}$
directly, with no additional $\pm 1$ terms.
Since $\bar\varphi_n = 0$ for $n < 0$,
the adjoint seed satisfies $\phi_{t,w,0} = 0$
and $\phi_{t,w,1}(v) = \bar\varphi_0(t, w{-}v) = 1$.
The time recurrence \eqref{eq:BLPP-phi-t} for $\phi$
follows from
$\pa_t \bar\varphi_n = -\tfrac{1}{2}\bar\varphi_{n-2}$
and two applications of \eqref{eq:BLPP-phi-a}.
For $\psi$, the heat identity
$\pa_t \varphi_n = \tfrac{1}{2}\varphi_{n+2}$ gives
$\pa_t \psi_{t,z,n}
= \tfrac{1}{2}\psi_{t,z,n+2}
= \tfrac{1}{2}\pa_z^2 \psi_{t,z,n}$,
where the second equality uses two applications of
\eqref{eq:BLPP-psi-a}.  The interchange of $\pa_t$ with the hitting-time
expectation is justified by dominated convergence:
on $\{\tau \leq t - \epsilon\}$ for any $\epsilon > 0$,
the differentiated integrand
$\varphi_{n+2}(t - \tau, -b(\tau) - z)$ is bounded
by a constant depending on $\epsilon$, $n$, and
$\sup_{[0,t]}|b|$.  Since $b$ is continuous,
$\PP(\tau = t) = 0$, so letting $\epsilon \to 0$
recovers the identity on the full event
$\{\tau \leq t\}$.
\end{proof}

The BLPP propagator inherits the structural properties
of its RBM counterpart, with the sign in the splitting
identity \eqref{eq:BLPP-Ba-splitting} arising from the
reversed integration direction.

\begin{lemma}\label{lem:BLPP-propagator}
The propagator satisfies the following identities.
\begin{enumerate}[label=\textup{(\roman*)},
leftmargin=*]
\item \textup{(Splitting with sign.)}
For $i < p < j$,
\begin{equation}\label{eq:BLPP-Ba-splitting}
\pa_{a_p} [B_{\mathbf{t},\mathbf{a}}]_{ij}(z, w)
= -[B_{\mathbf{t},\mathbf{a}}]_{ip}(z, a_p)[B_{\mathbf{t},\mathbf{a}}]_{pj}(a_p, w).
\end{equation}

\item \textup{(Flux conservation.)}
For $i < j$,
\begin{equation}\label{eq:BLPP-Ba-flux}
\Bigl(\pa_z + \pa_w
+ \sum_{i < p < j} \pa_{a_p}\Bigr)
[B_{\mathbf{t},\mathbf{a}}]_{ij}(z, w) = 0.
\end{equation}

\item \textup{(Time conservation.)}
\begin{equation}\label{eq:BLPP-Ba-time}
\sum_{p=1}^k \pa_{t_p}
[B_{\mathbf{t},\mathbf{a}}]_{ij} = 0.
\end{equation}

\item The propagator is $n$-independent.
\end{enumerate}
\end{lemma}

\begin{proof}
The sign in \eqref{eq:BLPP-Ba-splitting} arises because
$\pa_{a_p}$ acts on a lower limit of integration
$a_p$ (rather than an upper limit as in RBM), producing
an extra factor of $-1$.  Flux conservation follows from the telescoping argument
of Lemma~\ref{lem:RBM-propagator}, using
$(\pa_u + \pa_v)N_{ij}(u, v) = 0$ (the heat kernel
depends on $u - v$ alone).  For time conservation, the
integration limits in \eqref{eq:BLPP-propagator-def}
carry no time dependence, so $\sum_p \pa_{t_p}$ passes
inside each chain integral.  Each transition factor
$N_{i_r, i_{r+1}}$ depends on $t_{i_{r+1}} - t_{i_r}$
alone, giving
$(\pa_{t_{i_r}} + \pa_{t_{i_{r+1}}})
N_{i_r, i_{r+1}} = 0$; applying $\sum_p \pa_{t_p}$
to a chain product and regrouping by factor, each
factor is annihilated by its pair of time
derivatives.  Property~(iv) holds because
the transition kernels
$N_{ij}(u, v) = -e^{\tfrac{(t_j-t_i)}{2}\pa^2}(u, v)$
carry no $n$-dependence.
\end{proof}

For each $(\mathbf{t}, \mathbf{a}, n)$, define the
product-graph wave functions
$\Psi^n \colon \R^k \to L^2(\R)$ and
$\Phi^n \colon L^2(\R) \to \R^k$ componentwise by
\begin{align}\label{eq:BLPP-psi-prop}
[\psi_{\mathbf{t},\mathbf{a},n}]_{pi} \defeq
\int_{a_p}^{\infty}
[B_{\mathbf{t},\mathbf{a}}]_{pi}(z, a_i)
\psi_{t_p, z, n}\diff z, \qquad
[\phi_{\mathbf{t},\mathbf{a},n}]_{ip} \defeq
\int_{a_p}^{\infty}
[B_{\mathbf{t},\mathbf{a}}]_{ip}(a_i, w)
\phi_{t_p, w, n}\diff w,
\end{align}
and
\begin{equation}\label{eq:BLPP-WF-def}
\Psi_i^n \defeq \sum_{p=1}^i
[\psi_{\mathbf{t},\mathbf{a},n}]_{pi},
\qquad
\Phi_i^n \defeq \sum_{p=i}^k
[\phi_{\mathbf{t},\mathbf{a},n}]_{ip}.
\end{equation}

The absence of the exponential prefactor simplifies the
BLPP wave function recurrences compared with their RBM
counterparts.

\begin{lemma}\label{lem:BLPP-WF-recurrences}
The product-graph wave functions satisfy the linear
problem
\eqref{eq:dc-clean-linear-d1}--\eqref{eq:dc-clean-linear-S3}
and its adjoint
\eqref{eq:dc-clean-adjoint-d1}--\eqref{eq:dc-clean-adjoint-S3}
with the shifts and edge weights above.
Explicitly:
\begin{alignat}{2}
\pa_{\mathbf{a}} \Psi^n
&= -\Psi^{n+1},
&\qquad
\pa_{\mathbf{a}} \Phi^n
&= \Phi^{n-1},
\label{eq:BLPP-WF-a} \\[4pt]
\pa_{\mathbf{t}} \Psi^n
&= \tfrac{1}{2}\pa_{\mathbf{a}}^2 \Psi^n,
&\qquad
\pa_{\mathbf{t}} \Phi^n
&= -\tfrac{1}{2}\pa_{\mathbf{a}}^2 \Phi^n.
\label{eq:BLPP-WF-t}
\end{alignat}
\end{lemma}

\begin{proof}
The argument is analogous to the proof of
Lemma~\ref{lem:RBM-WF-recurrences}.  The propagator is
$n$-independent (Lemma~\ref{lem:BLPP-propagator}(iv)),
so the shift $n \to n+1$ acts only on the seeds in
\eqref{eq:BLPP-psi-prop}.  The BLPP seed recurrence
\eqref{eq:BLPP-psi-a} gives
$\psi_{t,z,n+1} = -\pa_z \psi_{t,z,n}$ (a pure
derivative, with no constant shift as in RBM), so
the shifted dressed component is
\[
[\psi]_{pi}^{n+1}
= -\int_{a_p}^{\infty}
[B]_{pi}(z, a_i)\pa_z\psi_{t_p,z,n}\diff z.
\]
The integration-by-parts, flux-conservation, and
splitting argument of
Lemma~\ref{lem:RBM-WF-recurrences} identifies the
integral as $\pa_{\mathbf{a}}[\psi]_{pi}$
(with sign modifications from the reversed limits
$[a_p, \infty)$ and the signed splitting identity
\eqref{eq:BLPP-Ba-splitting}), giving
$[\psi]_{pi}^{n+1}
= -\pa_{\mathbf{a}}[\psi]_{pi}$.
Summing over $p \leq i$ via \eqref{eq:BLPP-WF-def}
gives \eqref{eq:BLPP-WF-a} for $\Psi$.
The adjoint identity for $\Phi$ follows by the same
argument applied to $[\phi]_{ip}$, using
\eqref{eq:BLPP-phi-a}.

The time recurrences follow from two applications
of the spatial recurrence:
$\pa_{\mathbf{a}}^2 \Psi^n = \Psi^{n+2}$.
The time conservation
$\pa_{\mathbf{t}} [B]_{ij} = 0$
(Lemma~\ref{lem:BLPP-propagator}(iii)) ensures
$\pa_{\mathbf{t}}$
passes through the integrals and acts only on the
seeds via the heat identity
$\pa_t \psi_{t,z,n}
= \tfrac{1}{2}\pa_z^2 \psi_{t,z,n}$
\eqref{eq:BLPP-psi-t}, giving
$\pa_{\mathbf{t}} \Psi^n
= \tfrac{1}{2}\Psi^{n+2}
= \tfrac{1}{2}\pa_{\mathbf{a}}^2 \Psi^n$.
The adjoint time recurrence follows identically
from \eqref{eq:BLPP-phi-t}.
\end{proof}

Dressing compatibility
(Definition~\ref{def:dc-clean-dressing}) requires that the
derivatives of $K$ in each direction factorize through
$\Psi$ and $\Phi$.

\begin{lemma}\label{lem:BLPP-dressing}
The kernel $K_n$ is dressing compatible
(Definition~\ref{def:dc-clean-dressing}) with
$(\Psi, \Phi)$ and the edge weights above.
Explicitly:
\begin{align}
\pa_{\mathbf{a}} K_n
&= -\Psi^n \Phi^n,
\label{eq:BLPP-DC-a} \\
\pa_{\mathbf{t}} K_n
&= \tfrac{1}{2}\bigl(
-\pa_{\mathbf{a}}\Psi^n \Phi^n
+ \Psi^n \pa_{\mathbf{a}}\Phi^n\bigr),
\label{eq:BLPP-DC-t} \\
K_{n+1} - K_n
&= \Psi^n \Phi^{n+1}.
\label{eq:BLPP-DC-n}
\end{align}
\end{lemma}

\begin{proof}
\noindent\textit{Spatial derivative \eqref{eq:BLPP-DC-a}.}
Differentiating the kernel decomposition
\eqref{eq:BLPP-kernel-decomp} by $\pa_{a_p}$ produces
four contributions, as in the RBM proof of
\eqref{eq:RBM-DC-a}: the diagonal boundary evaluation
at $a_p$, two off-diagonal boundary evaluations where
$\pa_{a_p}$ acts on the lower integration limits
$a_i = a_p$ or $a_j = a_p$, and the interior term
where $\pa_{a_p}$ hits $[B]_{ij}$ through the splitting
identity \eqref{eq:BLPP-Ba-splitting}.  Each
contribution carries a minus sign: the fundamental
theorem of calculus on the lower limit
$\int_{a_p}^{\infty}$ produces $-1$ (versus $+1$ from
the upper limit $\int_{-\infty}^{a_p}$ in RBM), and
the splitting identity \eqref{eq:BLPP-Ba-splitting}
carries its own minus.  The four terms reassemble as
$\pa_{a_p} K_n = -\Psi_p^n \Phi_p^n$ by the wave
function definitions \eqref{eq:BLPP-WF-def}.
Summing over $p$ gives \eqref{eq:BLPP-DC-a}.

\medskip
\noindent\textit{Discrete shift \eqref{eq:BLPP-DC-n}.}
The propagator is $n$-independent, so the shift
$n \to n - 1$ acts only on the seeds.  Since
$\psi_{z,n} = -\pa_z\psi_{z,n-1}$ by
\eqref{eq:BLPP-psi-a}, substituting gives
$\psi_{z,n} - \psi_{z,n-1}
= -\pa_z\psi_{z,n-1} - \psi_{z,n-1}
= -(\pa_z + 1)\psi_{z,n-1}$.
Similarly, $\phi_{w,n-1} = \pa_w\phi_{w,n}$ by
\eqref{eq:BLPP-phi-a} gives
$\phi_{w,n} - \phi_{w,n-1}
= \phi_{w,n} - \pa_w\phi_{w,n}
= (1 - \pa_w)\phi_{w,n}$.
For each diagonal term
\begin{align*}
\psi_{r,n}\phi_{r,n}
- \psi_{r,n-1}\phi_{r,n-1}
= (-(\pa_r{+}1)\psi_{r,n-1})\phi_{r,n}
+ \psi_{r,n-1}((1{-}\pa_r)\phi_{r,n})
= -\pa_r(\psi_{r,n-1}\phi_{r,n}),
\end{align*}
where the constant terms
$-\psi_{r,n-1}\phi_{r,n}$ and
$+\psi_{r,n-1}\phi_{r,n}$ cancel.  Integrating
over $[a_i, \infty)$ gives
$\psi_{a_i,n-1}\phi_{a_i,n}$, since
$\int_{a_i}^{\infty}[-\pa_r f]\diff r
= f(a_i)$ when $f$ vanishes at infinity.

For the off-diagonal, the same cancellation
gives 
$-(\pa_z\psi_{z,n-1})[B]_{ij}\phi_{w,n}
- \psi_{z,n-1}[B]_{ij}(\pa_w\phi_{w,n})$.
Expanding
$(\pa_z + \pa_w
+ \sum_{i<p<j}\pa_{a_p})
[\psi_{z,n-1}[B]_{ij}\phi_{w,n}]$
by the Leibniz rule, flux conservation
\eqref{eq:BLPP-Ba-flux} eliminates all derivatives
of $[B]_{ij}$, confirming the integrand is the
negative divergence
\[
-(\pa_z + \pa_w
+ \textstyle\sum_{i<p<j}\pa_{a_p})
\bigl[\psi_{z,n-1}[B]_{ij}(z,w)
\phi_{w,n}\bigr].
\]
As in the diagonal, the negative divergence and
the reversed limits $[a_p, \infty)$ combine to give
positive boundary evaluations.  The three
contributions are
$[\psi_{n-1}]_{ij}[\phi_n]_{jj}$ (from $\pa_w$ at
$w = a_j$),
$[\psi_{n-1}]_{ii}[\phi_n]_{ij}$ (from $\pa_z$ at
$z = a_i$), and
$[\psi_{n-1}]_{ip}[\phi_n]_{pj}$ for $i < p < j$
(from splitting \eqref{eq:BLPP-Ba-splitting}).
Collecting diagonal and off-diagonal contributions
and shifting $n \to n+1$ gives \eqref{eq:BLPP-DC-n}.

\medskip
\noindent\textit{Time derivative \eqref{eq:BLPP-DC-t}.}
The propagator satisfies
$\pa_{\mathbf{t}} [B]_{ij} = 0$ by
\eqref{eq:BLPP-Ba-time}, so $\pa_{\mathbf{t}}$ acts
only on the seeds.  The driftless seed heat
equations
$\pa_t \psi = \tfrac{1}{2}\pa_z^2 \psi$ and
$\pa_t \phi = -\tfrac{1}{2}\pa_w^2 \phi$
(from \eqref{eq:BLPP-psi-t}--\eqref{eq:BLPP-phi-t})
give, for the diagonal,
\begin{align*}
\pa_t(\psi_{r,n}\phi_{r,n})
&= \tfrac{1}{2}\bigl[(\pa_r^2\psi_{r,n})\phi_{r,n}
- \psi_{r,n}(\pa_r^2\phi_{r,n})\bigr] \\
&= \tfrac{1}{2}\pa_r\bigl[
(\pa_r\psi_{r,n})\phi_{r,n}
- \psi_{r,n}(\pa_r\phi_{r,n})\bigr].
\end{align*}
Integrating over $[a_i, \infty)$ and substituting
$\pa_r\psi_{r,n} = -\psi_{r,n+1}$ and
$\pa_r\phi_{r,n} = \phi_{r,n-1}$ at the boundary
gives
$\tfrac{1}{2}(\psi_{a_i,n+1}\phi_{a_i,n}
+ \psi_{a_i,n}\phi_{a_i,n-1})$.

For the off-diagonal, flux conservation identifies
the integrand as a divergence
\begin{align*}
\pa_{\mathbf{t}}[\psi_{z,n}[B]_{ij}\phi_{w,n}]
= -\tfrac{1}{2}\bigl(\pa_z + \pa_w
+ \textstyle\sum_{i<p<j}\pa_{a_p}\bigr)
\bigl[\psi_{z,n+1}[B]_{ij}\phi_{w,n}
+ \psi_{z,n}[B]_{ij}\phi_{w,n-1}\bigr],
\end{align*}
where the seed recurrences $\pa_z\psi_n
= -\psi_{n+1}$ and $\pa_w\phi_n
= \phi_{n-1}$ simplify the bracket.
As in the discrete shift, the negative divergence
and reversed limits combine to give positive
boundary evaluations, producing
$\sum_{p=i}^j [\psi_{n+1}]_{ip}[\phi_n]_{pj}
+ [\psi_n]_{ip}[\phi_{n-1}]_{pj}$.
Collecting diagonal and off-diagonal contributions,
$\pa_{\mathbf{t}} K_n
= \tfrac{1}{2}(\Psi^{n+1}\Phi^n
+ \Psi^n\Phi^{n-1})$.
Substituting
$\Psi^{n+1} = -\pa_{\mathbf{a}}\Psi^n$ and
$\Phi^{n-1} = \pa_{\mathbf{a}}\Phi^n$,
this equals
$\tfrac{1}{2}(-\pa_{\mathbf{a}}\Psi^n\Phi^n
+ \Psi^n\pa_{\mathbf{a}}\Phi^n)$,
confirming \eqref{eq:BLPP-DC-t}.
\end{proof}

\subsection{Multipoint equation}\label{sec:BLPP-multipoint}

The variational formula gives
$G(b; t_i, n) \geq b(t_i)$ for each $i$.  The
Fredholm determinant
$\det(I - \chi_{\mathbf{a}} K \chi_{\mathbf{a}})$ is
therefore strictly positive for
$a_i > b(t_i)$, $1 \leq i \leq k$, and the resolvent
$R_n = (I - K_n)^{-1}$ exists on this set.

The dressed observable
$\mathcal{A}(u) = z\Phi(u) R(u) \Psi(u)$ of
Proposition~\ref{prop:dc-clean-dressed-waves} specializes at
$z = 1$ to the $k \times k$ matrix
\begin{equation}\label{eq:BLPP-A-def}
\mathcal{A}_n \defeq \Phi^n R_n \Psi^n \in \End(\R^k),
\qquad
(\mathcal{A}_n)_{ij} = \langle R_n \Psi_j^n, \Phi_i^n \rangle.
\end{equation}

\begin{theorem}\label{thm:BLPP-multipoint}
The dressed observable $\mathcal{A}_n$ satisfies
\begin{align}\label{eq:BLPP-multipoint}
\pa_{\mathbf{t}} (\mathcal{A}_{n+1} - \mathcal{A}_n)
- \tfrac{1}{2}\pa_{\mathbf{a}}^2(\mathcal{A}_{n+1}
+ \mathcal{A}_n)
+ \pa_{\mathbf{a}}(\mathcal{A}_n \mathcal{A}_{n+1})
- \mathcal{A}_{n+1} \pa_{\mathbf{a}} \mathcal{A}_{n+1}
- \pa_{\mathbf{a}} \mathcal{A}_n
  \mathcal{A}_n = 0.
\end{align}
\end{theorem}

\begin{proof}
Theorem~\ref{thm:dc-clean-general-darboux} applies
at $z = 1$ with $\pa_2 = -\pa_{\mathbf{a}}$,
$C_1 = 0$, $C_2 = I$,
$\Lambda_1 = \Lambda_2 = 0$.
The dressed $\Lambda$-weights
\eqref{eq:dc-clean-L1-dressed-general}--\eqref{eq:dc-clean-L2-dressed-general}
are
\begin{align*}
\widehat\Lambda_1
= \pa_2\mathcal{A}_n
+ [\mathcal{A}_n, 0]
= -\pa_{\mathbf{a}} \mathcal{A}_n, \qquad
\widehat\Lambda_2
= [\mathcal{A}_n, I]_{S_3}
= \mathcal{A}_n - \mathcal{A}_{n+1}.
\end{align*}

The $\Lambda$-diamond equation
\eqref{eq:dc-clean-L-diamond} with $C_1 = 0$ reduces
to
\begin{equation}\label{eq:BLPP-dressed-L-diamond}
\mathcal{H}\widehat\Lambda_2
+ \pa_2\widehat\Lambda_1 C_2
- [\widehat\Lambda_1,\widehat\Lambda_2]_{S_3} = 0,
\end{equation}
since $C_1 = 0$ kills the term
$\pa_2\widehat\Lambda_2 C_1(S_3 u)$.
Since $C_2 = I$, the term
$\pa_2\widehat\Lambda_1 C_2$ reduces to
$\pa_2\widehat\Lambda_1$, and substituting the
dressed weights into
\eqref{eq:BLPP-dressed-L-diamond}:
the heat term is
$\mathcal{H}(\mathcal{A}_n - \mathcal{A}_{n+1})$;
the spatial term is
$(-\pa_{\mathbf{a}})(-\pa_{\mathbf{a}} \mathcal{A}_n)
= \pa_{\mathbf{a}}^2 \mathcal{A}_n$;
and the shifted commutator is
\begin{align*}
[\widehat\Lambda_1,\widehat\Lambda_2]_{S_3}
&= (-\pa_{\mathbf{a}} \mathcal{A}_n)
(\mathcal{A}_n - \mathcal{A}_{n+1})
- (\mathcal{A}_n - \mathcal{A}_{n+1})
(-\pa_{\mathbf{a}} \mathcal{A}_{n+1}) \\
&= -(\pa_{\mathbf{a}} \mathcal{A}_n)
(\mathcal{A}_n - \mathcal{A}_{n+1})
+ (\mathcal{A}_n - \mathcal{A}_{n+1})
(\pa_{\mathbf{a}} \mathcal{A}_{n+1}).
\end{align*}
Collecting all three contributions and expanding
$\mathcal{H} = \pa_{\mathbf{t}} - \tfrac{1}{2}\pa_{\mathbf{a}}^2$:
\begin{align*}
-\pa_{\mathbf{t}}(\mathcal{A}_{n+1} - \mathcal{A}_n)
+ \tfrac{1}{2}\pa_{\mathbf{a}}^2(\mathcal{A}_n
+ \mathcal{A}_{n+1})
- (\pa_{\mathbf{a}} \mathcal{A}_n)
(\mathcal{A}_{n+1} - \mathcal{A}_n)
+ (\mathcal{A}_{n+1} - \mathcal{A}_n)
(\pa_{\mathbf{a}} \mathcal{A}_{n+1}) = 0.
\end{align*}
Multiplying by $-1$ and distributing the products
$(\pa_{\mathbf{a}} \mathcal{A}_n)
(\mathcal{A}_{n+1} - \mathcal{A}_n)$ and
$(\mathcal{A}_{n+1} - \mathcal{A}_n)
\pa_{\mathbf{a}} \mathcal{A}_{n+1}$,
the Leibniz rule recombines the cross-terms as
$\pa_{\mathbf{a}}(\mathcal{A}_n
\mathcal{A}_{n+1})$,
giving \eqref{eq:BLPP-multipoint}.
\end{proof}

\begin{corollary}
\label{cor:BLPP-trace-defect}
Let $F_n = \det(I - K_n)$.
Then
$\pa_{\mathbf{a}} \log F_n
= \operatorname{tr}_E(\mathcal{A}_n)$,
and
\begin{equation}\label{eq:BLPP-Hirota}
\bigl[D_{\mathbf{t}} - \tfrac{1}{2}D_{\mathbf{a}}^2\bigr]
F_{n+1} \cdot F_n
= \tfrac{1}{2}F_{n+1}
F_n
\bigl(\operatorname{tr}_E((\mathcal{A}_{n+1}
- \mathcal{A}_n)^2)
- (\operatorname{tr}_E(\mathcal{A}_{n+1}
- \mathcal{A}_n))^2\bigr).
\end{equation}
\end{corollary}

\begin{proof}
Proposition~\ref{prop:dc-clean-normalized-trace-defect-direct}
with $C_1 = 0$, $C_2 = I$ and $\Lambda_2 = 0$ gives
$C_1 + \Lambda_2 C_2^{-1} = 0$, so
Corollary~\ref{cor:dc-clean-scalar-trace-defect} applies
with $a(u) = 0$.  The quantity
$\mathcal{B} = \Delta_2 C_2^{-1}$
(where $\Delta_2 = \widehat\Lambda_2 - \Lambda_2$)
reduces to
$\mathcal{B} = \mathcal{A}_n - \mathcal{A}_{n+1}$, and
\eqref{eq:dc-clean-normalized-scalar-linear-trace-defect-direct}
gives \eqref{eq:BLPP-Hirota}.
The logarithmic derivative has the opposite sign
from RBM because $\pa_2 = -\pa_{\mathbf{a}}$:
the framework identity
$\pa_2 \log F = -\operatorname{tr}_E(\mathcal{A})$
gives $\pa_{\mathbf{a}} \log F_n
= \operatorname{tr}_E(\mathcal{A}_n)$.
\end{proof}

\begin{corollary}
\label{cor:BLPP-one-point}
For $k = 1$, the Fredholm determinant
$F_{t,a,n} = \det(I - K_{t,a,n})$ satisfies
\begin{equation}\label{eq:BLPP-scalar}
\bigl[D_t - \tfrac{1}{2}D_a^2\bigr]
F_{t,a,n+1} \cdot F_{t,a,n} = 0.
\end{equation}
\end{corollary}

\begin{proof}
When $k = 1$, $\mathcal{A}_n$ is a $1 \times 1$
matrix (a scalar).  The trace-defect term
$\operatorname{tr}_E((\mathcal{A}_{n+1}
- \mathcal{A}_n)^2)
- (\operatorname{tr}_E(\mathcal{A}_{n+1}
- \mathcal{A}_n))^2$ vanishes, and
\eqref{eq:BLPP-Hirota} reduces to
\eqref{eq:BLPP-scalar}.
\end{proof}

%% file: chapters/10-the-kpz-fixed-point.tex
\chapter{The KPZ Fixed Point}
\label{ch:the-kpz-fixed-point}

{
  \setlength{\parskip}{0pt}
}

\label{sec:FP}

This section derives the multipoint equation for the KPZ
fixed point under one-sided initial data by identifying
the Fredholm determinant kernel
with the continuum framework
(\S\ref{sec:kpz-reduced-framework}).  The derivation
proceeds in three stages: reformulation of the
extended-kernel Fredholm determinant as a single-space
operator on $L^2(\R)$, verification of the seed
recurrences and dressing compatibility, and application of
the Darboux theorem
(Proposition~\ref{thm:kpz-red-Darboux} and
Theorem~\ref{thm:kpz-red-Darboux-compat}) to produce the
matrix potential KP equation.

\paragraph{\textbf{System description.}}
The KPZ fixed point is the universal scaling limit of the KPZ
universality class.  It describes the evolution of an
interface height function $\hh(t,x)$, with $t \geq 0$ and
$x \in \R$, starting from an initial condition
$\hh_0 \colon \R \to \R \cup \{-\infty\}$ that is upper
semicontinuous with $\hh_0(x) \leq C(1 + |x|)$ for
some $C < \infty$.

\paragraph{\textbf{Fredholm determinant formula.}}
Fix $t > 0$, $m \geq 1$ spatial sites
$\mathbf{x} = (x_1, \dotsc, x_m) \in \R^m$, and
thresholds
$\mathbf{a} = (a_1, \dotsc, a_m) \in \R^m$.  For
one-sided initial data ($\hh_0(x) = -\infty$ for all
$x > 0$), the joint distribution is
\begin{equation}\label{eq:FP-Fredholm-extended}
\PP\Bigl(\bigcap_{i=1}^m
\{\hh(t, x_i) \leq a_i\}\Bigr)
= \det(I - \chi_{\mathbf{a}} K_t
\chi_{\mathbf{a}})_{L^2(\{x_1, \dotsc, x_m\}
\times \R)},
\end{equation}
where $\chi_{\mathbf{a}}(x_i, u) =
\mathbf{1}_{u > a_i}$.  The extended correlation kernel
$K_t$ has block entries
\begin{equation}\label{eq:FP-extended-kernel}
K_t(x_i, u; x_j, v)
= -e^{(x_j - x_i)\pa^2}\mathbf{1}_{x_i < x_j}(u, v)
+ \bigl(\mathcal{S}_{t,-x_i}^{
\mathrm{hypo}(\hh_0^-)}\bigr)^*
\mathcal{S}_{t,x_j}(u, v),
\end{equation}
where $A^*(u,v) = A(v,u)$ denotes the transposed kernel.

\medskip
\noindent\textit{The heat semigroup
$e^{s\pa^2}$.}
The integral kernel of $e^{s\pa^2}$ is
\[
e^{s\pa^2}(u, v) = \frac{1}{\sqrt{4\pi s}}
\exp\bigl(-\frac{(u-v)^2}{4s}\bigr),
\qquad s > 0.
\]

\medskip
\noindent\textit{The scattering operator
$\mathcal{S}_{t,x}$.}
The scattering operator is defined by
\begin{equation}\label{eq:FP-S-def}
\mathcal{S}_{t,x}(u, v)
\defeq t^{-1/3} e^{\frac{2x^3}{3t^2}
- \frac{(u-v)x}{t}}
\mathrm{Ai}\bigl(-t^{-1/3}(u-v)
+ t^{-4/3}x^2\bigr),
\end{equation}
where $\mathrm{Ai}(\cdot)$ is the Airy function.

\medskip
\noindent\textit{The hypo-operator
$\mathcal{S}_{t,x}^{\mathrm{hypo}(\hh_0^-)}$.}
Given upper semicontinuous initial data
$\hh_0$, set
$\hh_0^-(\cdot) \defeq \hh_0(-\cdot)$ and define
\begin{equation}\label{eq:FP-Shypo-def}
\mathcal{S}_{t,x}^{\mathrm{hypo}(\hh_0^-)}(u, v)
\defeq \E_{B(0) = u}\bigl[
\mathcal{S}_{t,x-\tau}(B(\tau), v)
\mathbf{1}_{\tau < \infty}\bigr],
\end{equation}
where $B$ is a Brownian motion with
diffusion coefficient $2$ and
$\tau \defeq \inf\{y \geq 0 : B(y) \leq \hh_0^-(y)\}$.

\paragraph{\textbf{References.}}
The Fredholm determinant formula
\eqref{eq:FP-Fredholm-extended} is the one-sided
specialization of the construction of Matetski, Quastel,
and Remenik \cite{MQR21}, who built the KPZ fixed point
as the scaling limit of TASEP for general upper
semicontinuous initial data.  The matrix KP equation
for the multipoint distribution was proved by
Quastel and Remenik~\cite{QR2022} in the same generality.

\section{Kernel reformulation}
\label{sec:FP-kernel-reform}

The extended-kernel Fredholm determinant on
$L^2(\{x_1, \dotsc, x_m\} \times \R)$ reduces to a
single-space Fredholm determinant on $L^2(\R)$.
Throughout this subsection the spatial sites are
ordered: $x_1 < \cdots < x_m$.  For distinct sites
this entails no loss, since relabelling the pairs
$(x_i, a_i)$ permutes the block indices and leaves
the Fredholm determinant unchanged.

Define the seed functions
\begin{equation}\label{eq:FP-seed-def}
\psi_{t,x,a}(u) \defeq \mathcal{S}_{t,x}(u, a),
\qquad
\phi_{t,x,a}(u) \defeq
\mathcal{S}_{t,-x}^{\mathrm{hypo}(\hh_0^-)}(u, a),
\end{equation}
where the dependence of $\phi$ on the initial
condition $\hh_0$ is suppressed.  The extended kernel
decomposes as
$K_t = N + (\mathcal{S}^{\mathrm{hypo}})^*
\mathcal{S}$,
where $N$ is the strictly upper-triangular block
operator with entries
$N_{ij}(u, v) \defeq
-e^{(x_j - x_i)\pa^2}\mathbf{1}_{x_i < x_j}(u, v)$.

The propagator is built from $N$.  For
$1 \leq i, j \leq m$,
\begin{align}\label{eq:FP-propagator}
&[B_{\mathbf{x},\mathbf{a}}]_{ij}(z, w) \defeq \notag \\
&\begin{cases}
\displaystyle\sum_{k=1}^{j-i}
\sum_{\substack{i = i_0 < i_1 < \cdots < i_k = j}}
\int_{a_{i_1}}^{\infty} \cdots
\int_{a_{i_{k-1}}}^{\infty}
\prod_{r=0}^{k-1}
N_{i_r, i_{r+1}}(\xi_r, \xi_{r+1})
\diff\xi_{k-1} \cdots \diff\xi_1,
& i < j, \\[6pt]
\delta(z - w), & i = j, \\
0, & i > j,
\end{cases}
\end{align}
with $\xi_0 = z$ and $\xi_k = w$.  The internal
vertices are constrained by the cutoffs, but the
endpoints $z, w$ are arbitrary.

\begin{lemma}\label{lem:FP-kernel}
The joint distribution satisfies
\begin{equation}\label{eq:FP-single-space}
\PP\Bigl(\bigcap_{i=1}^m
\{\hh(t, x_i) \leq a_i\}\Bigr)
= \det(I - K_{t,\mathbf{x},\mathbf{a}})_{L^2(\R)},
\end{equation}
where
\begin{equation}\label{eq:FP-kernel-formula}
K_{t,\mathbf{x},\mathbf{a}}
= \sum_{1 \leq i \leq j \leq m}
\int_{a_i}^{\infty}
\int_{a_j}^{\infty}
\psi_{t,x_i,z}
[B_{\mathbf{x},\mathbf{a}}]_{ij}(z, w)
\phi_{t,x_j,w}
\diff w \diff z.
\end{equation}
\end{lemma}

\begin{proof}
Write $\mathscr{X} = L^2(\{x_1, \dotsc, x_m\}
\times \R)$ for the extended space and
$H = L^2(\R)$.  The cutoff extended kernel
decomposes as
\[
\chi_{\mathbf{a}} K_t \chi_{\mathbf{a}} = L + Q,
\]
where $L \defeq \chi_{\mathbf{a}} N
\chi_{\mathbf{a}}$ and
$Q \defeq \chi_{\mathbf{a}}
(\mathcal{S}^{\mathrm{hypo}})^*\mathcal{S}
\chi_{\mathbf{a}}$.  Since
$x_1 < \cdots < x_m$, the block operator $N$ is
strictly upper triangular, and $L^m = 0$.  The
Neumann series
$(I - L)^{-1} = \sum_{k=0}^{m-1} L^k$ terminates
and gives a bounded operator on $\mathscr{X}$.

\medskip
\noindent\textit{Conjugated trace-class realization.}
The Fredholm determinant in
\eqref{eq:FP-Fredholm-extended} is understood in
the conjugated trace-class realization of
Matetski--Quastel--Remenik~\cite{MQR21}.  Define the
polynomial weight
$\vartheta_i(u) \defeq (1 + u^2)^{2i}$ and the
block-diagonal conjugation
$(\theta f)_i(u) \defeq \vartheta_i(u) f_i(u)$ on
$\mathscr{X}$.  Since each $\vartheta_i$ depends only
on $u$, the operator $\theta$ commutes with the
threshold projections $\chi_{\mathbf{a}}$, and the
Fredholm determinant is invariant under conjugation:
\[
\det_{\mathscr{X}}(I - \chi_{\mathbf{a}} K_t
\chi_{\mathbf{a}})
= \det_{\mathscr{X}}(I -
\theta \chi_{\mathbf{a}} K_t
\chi_{\mathbf{a}} \theta^{-1}).
\]
Matetski--Quastel--Remenik~\cite{MQR21} establish
that
$\theta \chi_{\mathbf{a}} K_t
\chi_{\mathbf{a}} \theta^{-1}
\in \mathfrak{S}_1(\mathscr{X})$.
Set $\widetilde{L} \defeq \theta L \theta^{-1}$ and
$\widetilde{Q} \defeq \theta Q \theta^{-1}$.

\medskip
\noindent\textit{Trace-class estimate for
$\widetilde{L}$.}
Each off-diagonal block of $\widetilde{L}$ has
kernel
\[
\widetilde{L}_{ij}(u, v)
= -(1 + u^2)^{2i}(1 + v^2)^{-2j}
e^{(x_j - x_i)\pa^2}(u, v),
\qquad i < j,
\]
restricted to $u > a_i$, $v > a_j$; the block is
zero for $i \geq j$.  Fix $i < j$ and set
$s \defeq x_j - x_i$.  Since restriction to a
half-line contracts the trace norm, it suffices to
show that $M_fe^{s\pa^2}M_g$ is trace
class on $L^2(\R)$, where $f(u) = (1+u^2)^{2i}$
and $g(v) = (1+v^2)^{-2j}$.

The semigroup property
$e^{s\pa^2}
= e^{(s/2)\pa^2}e^{(s/2)\pa^2}$ provides a
factorization.  The polynomial balancing weight
$h(w) = (1 + w^2)^{i+j}$ inserted between the two
heat semigroup factors gives
\[
M_fe^{s\pa^2}M_g = AB,
\qquad
A \defeq M_fe^{(s/2)\pa^2}
 M_{h^{-1}},
\quad
B \defeq M_he^{(s/2)\pa^2}
 M_g.
\]
The Hilbert--Schmidt norm of $B$ satisfies
\[
\lVert B \rVert_{\mathfrak{S}_2}^2
= \int_{\R^2}
(1+w^2)^{2(i+j)}
\left(e^{(s/2)\pa^2}(w, v)\right)^2
(1+v^2)^{-4j} \diff w \diff v.
\]
At fixed $v$, the inner integral over $w$ is
$(4\pi s)^{-1/2}$ times a polynomial in $v$ of
degree $4(i+j)$: the Gaussian
$\exp(-(w-v)^2/s)$ has finite moments of all
orders, so convolution against
$(1+w^2)^{2(i+j)}$ produces at most polynomial
growth.  Therefore
\[
\lVert B \rVert_{\mathfrak{S}_2}^2
\leq C(s,i,j)
\int_{\R}
(1 + v^2)^{2(i+j) - 4j} \diff v
= C(s,i,j) \int_{\R}
(1 + v^2)^{-2(j-i)} \diff v
< \infty,
\]
since $j - i \geq 1$ gives the exponent
$-2(j-i) \leq -2$.  The same estimate with the
roles of $f$ and $h^{-1}$ exchanged gives
$\lVert A \rVert_{\mathfrak{S}_2} < \infty$:
the exponent is
$4i - 2(i+j) = -2(j-i)$, the same convergent
integral.  Therefore
$M_fe^{s\pa^2}M_g
\in \mathfrak{S}_1(L^2(\R))$
as the product of two Hilbert--Schmidt operators.
Since $\widetilde{L}$ has at most $\binom{m}{2}$
nonzero blocks, each trace class,
$\widetilde{L}
\in \mathfrak{S}_1(\mathscr{X})$.

\medskip
\noindent\textit{Determinant reduction.}
The operator $\widetilde{L}$ is trace class and
nilpotent, so all eigenvalues vanish and
$\det_{\mathscr{X}}(I - \widetilde{L}) = 1$.
Since
$\widetilde{L} + \widetilde{Q}
\in \mathfrak{S}_1(\mathscr{X})$
by Matetski--Quastel--Remenik~\cite{MQR21},
the scattering part satisfies
$\widetilde{Q}
\in \mathfrak{S}_1(\mathscr{X})$.
Multiplicativity of the Fredholm determinant gives
\begin{align*}
\det_{\mathscr{X}}
(I - \widetilde{L} - \widetilde{Q})
= \det_{\mathscr{X}}(I - \widetilde{L})
\cdot \det_{\mathscr{X}}\bigl(
I - (I - \widetilde{L})^{-1}
\widetilde{Q}\bigr)
= \det_{\mathscr{X}}\bigl(
I - (I - \widetilde{L})^{-1}
\widetilde{Q}\bigr),
\end{align*}
where
$(I - \widetilde{L})^{-1}
\widetilde{Q}
\in \mathfrak{S}_1(\mathscr{X})$ because
$(I - \widetilde{L})^{-1}$ is bounded and
$\widetilde{Q}$ is trace class.  Conjugation
invariance gives
$\det_{\mathscr{X}}(I -
(I - \widetilde{L})^{-1}\widetilde{Q})
= \det_{\mathscr{X}}(I -
(I - L)^{-1}Q)$.

\medskip
\noindent\textit{Propagator identification and
Sylvester's identity.}
Each power $L^k$ has the block structure
\[
(L^k)_{ij}
= \sum_{i = i_0 < i_1 < \cdots < i_k = j}
L_{i_0 i_1} L_{i_1 i_2} \cdots L_{i_{k-1} i_k},
\]
since nonzero contributions require strictly
increasing chains of block indices.  The
$(i,j)$-entry of $(I - L)^{-1}$ is therefore
$\delta_{ij} I + \sum_{k=1}^{j-i}(L^k)_{ij}$ for
$i \leq j$, and zero for $i > j$.  Unfolding the
definitions and using
$L_{ij}(u,v)
= \chi_{a_i}(u) N_{ij}(u,v) \chi_{a_j}(v)$,
the cutoff factors constrain the internal variables
to $\xi_s > a_{i_s}$; within the sandwich
$\chi_{\mathbf{a}}(I-L)^{-1}\chi_{\mathbf{a}}$,
the $(i,j)$-entry coincides with the propagator
$[B_{\mathbf{x},\mathbf{a}}]_{ij}$.

The scattering part $Q$ factors through $H$ as a
product of the column operator
$\mathcal{S}_{t,\mathbf{x}}$ and the row operator
$(\mathcal{S}_{t,-\mathbf{x}}^{\mathrm{hypo}})^*$,
with threshold projections in between.
Sylvester's identity
$\det(I - AB) = \det(I - BA)$, valid here because
the trace-class conditions on both sides are
supplied by the estimates of
Matetski--Quastel--Remenik~\cite{MQR21}, gives
\[
\det_{\mathscr{X}}(I -
(I - L)^{-1}Q)
= \det\bigl(I -
\mathcal{S}_{t,\mathbf{x}}
\chi_{\mathbf{a}}
(I - L)^{-1}
\chi_{\mathbf{a}}
(\mathcal{S}_{t,-\mathbf{x}}^{\mathrm{hypo}})^*
\bigr)_{L^2(\R)},
\]
which is $\det(I - K_{t,\mathbf{x},\mathbf{a}})$
after expanding the matrix product blockwise and
identifying the propagator entries with
\eqref{eq:FP-kernel-formula}.
\end{proof}

\begin{remark}\label{rem:FP-weighted-realization}
The kernel reformulation above takes $H = L^2(\R)$.
For observation windows containing non-positive spatial
sites, the Airy scattering seed $\psi_{t,x,a}$ does
not lie in $L^2(\R)$: the exponential factor
$e^{-(u-a)x/t}$ grows for $x < 0$, and for $x = 0$
the Airy oscillation decays too slowly.  A local
gauge resolves this.

Define
$M_\rho(u) \defeq e^{-\rho\max(u,0)}$ for a parameter
$\rho > 0$.  On any compact subset $\mathcal{U}$ of the
strict ordered chamber, choose
$\rho > \sup_{(t,\mathbf{x})\in\mathcal{U}}\max_i|x_i|/t$.
Replace each seed by its gauged counterpart
$M_\rho\psi_{t,x,a}$ and $M_\rho^{-1}\phi_{t,x,a}$.
The condition $\rho > |x_i|/t$ ensures that the
exponential damping $e^{-\rho u}$ overcomes the
non-integrability of $\psi_{t,x_i,a}$ as
$u \to +\infty$, so the gauged seed
$M_\rho\psi_{t,x_i,a}$ lies in $L^2(\R)$.  The
adjoint seed $M_\rho^{-1}\phi_{t,x,a}$ also lies in
$L^2(\R)$ as $u \to +\infty$: the hitting-time
representation and the Airy decay estimates of
Matetski--Quastel--Remenik~\cite{MQR21}, Appendix~A.1,
provide super-exponential decay that dominates the
factor $e^{\rho u}$ from $M_\rho^{-1}$.

The gauged kernel $M_\rho K M_\rho^{-1}$ is trace
class: each integrand in \eqref{eq:FP-kernel-formula}
contributes a rank-one operator
$(M_\rho\psi_{x_i,z})\otimes(M_\rho^{-1}\phi_{x_j,w})$,
which is trace class with trace norm bounded by the
product of the $L^2$ norms of the gauged seeds, and
the Gaussian decay of the heat-kernel propagator
$[B]_{ij}$ ensures that the $z$- and $w$-integrals
converge in trace norm.

Since $M_\rho$ depends only on the Hilbert-space
variable $u$ and not on the parameters
$(t,\mathbf{x},\mathbf{a})$, parameter derivatives
commute with the gauge, and the dressing identities
\eqref{eq:FP-DC-a}--\eqref{eq:FP-DC-t} and all
subsequent algebra hold for the gauged data.
\end{remark}

\section{Recurrences and kernel identities}
\label{sec:FP-recurrences}

The KPZ fixed point sits in the continuum framework of
\S\ref{sec:kpz-reduced-framework}: the wave functions
$\Psi$, $\Phi$ solve the linear problem
\eqref{eq:kpz-red-Lx}--\eqref{eq:kpz-red-Lt} and its
adjoint \eqref{eq:kpz-red-Ax}--\eqref{eq:kpz-red-At}
with $C = 0$, and the kernel $K$ is dressing compatible
with $(\Psi, \Phi)$ in the sense of
Definition~\ref{def:kpz-red-dressing}, so that
Proposition~\ref{thm:kpz-red-Darboux} and
Theorem~\ref{thm:kpz-red-Darboux-compat} apply.

\begin{lemma}\label{lem:FP-seed-linear}
The seed functions $\psi_{t,x,a}$ and $\phi_{t,x,a}$
\eqref{eq:FP-seed-def} satisfy the linear problem
\eqref{eq:kpz-red-Lx}--\eqref{eq:kpz-red-Lt} and its
adjoint \eqref{eq:kpz-red-Ax}--\eqref{eq:kpz-red-At}
with $C = \Lambda = 0$.  Explicitly:
\begin{align}
\pa_x \psi_{t,x,a}
&= \pa_a^2 \psi_{t,x,a},
\label{eq:FP-seed-heat} \\
\pa_t \psi_{t,x,a}
&= -\tfrac{1}{3}\pa_a^3 \psi_{t,x,a},
\label{eq:FP-seed-airy}
\end{align}
and
\begin{align}
\pa_x \phi_{t,x,a}
&= -\pa_a^2 \phi_{t,x,a},
\label{eq:FP-phi-heat} \\
\pa_t \phi_{t,x,a}
&= -\tfrac{1}{3}\pa_a^3 \phi_{t,x,a}.
\label{eq:FP-phi-airy}
\end{align}
\end{lemma}

\begin{proof}
\noindent\textit{Verification of
\eqref{eq:FP-seed-heat}--\eqref{eq:FP-seed-airy}
for $\psi$.}
Write
\[
\psi_{t,x,a}(u)
= t^{-1/3}e^{f(x,u,a)}\mathrm{Ai}(g(x,u,a))
\]
with
\begin{align*}
f(x,u,a) \defeq \tfrac{2x^3}{3t^2}
- \tfrac{(u-a)x}{t},
\qquad
g(x,u,a) \defeq -t^{-1/3}(u-a) + t^{-4/3}x^2.
\end{align*}
Both $f$ and $g$ are affine in $a$, with
$\pa_a f = x/t$ and $\pa_a g = t^{-1/3}$.

For \eqref{eq:FP-seed-heat}, the chain rule and
$\mathrm{Ai}''(g) = g\mathrm{Ai}(g)$
give
\begin{align*}
\pa_x \psi
&= t^{-1/3}e^f\bigl[(\pa_x f)\mathrm{Ai}(g)
+ (\pa_x g)\mathrm{Ai}'(g)\bigr], \\
\pa_a^2 \psi
&= t^{-1/3}e^f\bigl[\bigl((\pa_a f)^2
+ (\pa_a g)^2g\bigr)\mathrm{Ai}(g)
+ 2(\pa_a f)(\pa_a g)
\mathrm{Ai}'(g)\bigr],
\end{align*}
so $\pa_x\psi = \pa_a^2\psi$ reduces to
\[
\pa_x f = (\pa_a f)^2 + (\pa_a g)^2g,
\qquad
\pa_x g = 2(\pa_a f)(\pa_a g).
\]
Substituting
$\pa_x f = 2x^2/t^2 - (u{-}a)/t$ and
$\pa_x g = 2xt^{-4/3}$, the first identity becomes
$2x^2/t^2 - (u{-}a)/t = x^2/t^2 + t^{-2/3}g$,
which holds by expanding $g$; the second becomes
$2xt^{-4/3} = 2(x/t)(t^{-1/3})$.

For \eqref{eq:FP-seed-airy}, using additionally
$\mathrm{Ai}'''(g)
= \mathrm{Ai}(g) + g\mathrm{Ai}'(g)$:
\begin{align*}
\pa_a^3 \psi
&= t^{-1/3}e^f\bigl[\bigl((\pa_a f)^3
+ 3(\pa_a f)(\pa_a g)^2g
+ (\pa_a g)^3\bigr)\mathrm{Ai}(g) \\
&\qquad\qquad
+ \bigl(3(\pa_a f)^2(\pa_a g)
+ (\pa_a g)^3g\bigr)
\mathrm{Ai}'(g)\bigr], \\[4pt]
\pa_t \psi
&= t^{-1/3}e^f\bigl[
\bigl(\pa_t f - \tfrac{1}{3}t^{-1}\bigr)
\mathrm{Ai}(g)
+ (\pa_t g)\mathrm{Ai}'(g)\bigr],
\end{align*}
where $-\tfrac{1}{3}t^{-1}$ arises from
$\pa_t(t^{-1/3})$.  The identity
$\pa_t\psi = -\tfrac{1}{3}\pa_a^3\psi$ reduces to
\begin{align*}
3\pa_t f - t^{-1} + (\pa_a f)^3
&+ 3(\pa_a f)(\pa_a g)^2g
+ (\pa_a g)^3 = 0, \\
3\pa_t g + 3(\pa_a f)^2(\pa_a g)
+ (\pa_a g)^3g &= 0.
\end{align*}
Substituting
$\pa_t f = -\tfrac{4x^3}{3t^3} + \tfrac{(u{-}a)x}{t^2}$,
$\pa_t g = \tfrac{1}{3}t^{-4/3}(u{-}a) - \tfrac{4}{3}t^{-7/3}x^2$,
and expanding $g$ confirms both.

\medskip
\noindent\textit{Verification for $\phi$.}
By \eqref{eq:FP-Shypo-def},
\[
\phi_{t,x,a}(u) = \E_{B(0)=u}\bigl[
\mathcal{S}_{t,-x-\tau}(B(\tau), a)
\mathbf{1}_{\tau < \infty}\bigr].
\]
Since $\tau$ depends only on the path $B$ and
not on $(a, x, t)$, differentiation under the
expectation is permitted by dominated convergence:
the standard Airy decay bound and the upper-growth
hypothesis on $\hh_0$ provide a locally uniform
integrable majorant
(Matetski--Quastel--Remenik~\cite{MQR21},
Appendix~A.1).
To verify \eqref{eq:FP-phi-heat}, note that
$x$ enters $\mathcal{S}_{t,-x-\tau}(B(\tau),a)$
only through the second subscript, so the chain rule
and \eqref{eq:FP-seed-heat} give
\[
\pa_x \mathcal{S}_{t,-x-\tau}(B(\tau), a)
= -\pa_a^2 \mathcal{S}_{t,-x-\tau}(B(\tau), a),
\]
since $\pa_x(-x-\tau) = -1$.  Therefore
$\pa_x\phi = -\pa_a^2\phi$.
To verify \eqref{eq:FP-phi-airy}, notice that
$t$ does not appear in the second subscript of
$\mathcal{S}_{t,-x-\tau}(B(\tau),a)$; the derivative
$\pa_t$ therefore acts at fixed second subscript,
and \eqref{eq:FP-seed-airy} gives the result.
\end{proof}

The propagator $[B_{\mathbf{x},\mathbf{a}}]_{ij}$
defined by \eqref{eq:FP-propagator} satisfies three
structural properties.

\begin{lemma}\label{lem:FP-propagator}
The propagator satisfies:
\begin{enumerate}[label=\textup{(\roman*)},
leftmargin=*]
\item \textup{(Splitting)}
For $i < p < j$,
\begin{equation}\label{eq:FP-splitting}
\pa_{a_p} [B_{\mathbf{x},\mathbf{a}}]_{ij}(z, w)
= -[B_{\mathbf{x},\mathbf{a}}]_{ip}(z, a_p)
[B_{\mathbf{x},\mathbf{a}}]_{pj}(a_p, w).
\end{equation}

\item \textup{(Flux conservation)}
For $i < j$,
\begin{equation}\label{eq:FP-flux}
\Bigl(\pa_z + \pa_w
+ \sum_{i < p < j} \pa_{a_p}\Bigr)
[B_{\mathbf{x},\mathbf{a}}]_{ij}(z, w) = 0.
\end{equation}

\item \textup{($x$-independence)}
\begin{equation}\label{eq:FP-x-indep}
\sum_{p=1}^m \pa_{x_p}
[B_{\mathbf{x},\mathbf{a}}]_{ij}(z, w) = 0.
\end{equation}
\end{enumerate}
\end{lemma}

\begin{proof}
The propagator \eqref{eq:FP-propagator} has the same
chain-of-kernels structure as in
\S\ref{sec:RBM-recurrences}, with the single
difference that integration is from $a_p$ to $\infty$
rather than from $-\infty$ to $a_p$.  The proofs of
(i) and~(ii) in Lemma~\ref{lem:RBM-propagator}
transfer directly.  Splitting~(i) uses the fundamental
theorem of calculus at the integration limit; the sign
is reversed because
$\pa_{a_p}\int_{a_p}^{\infty} f = -f(a_p)$.  Flux
conservation~(ii) uses only the translation invariance
\begin{equation}\label{eq:FP-N-flux}
(\pa_u + \pa_v) N_{ij}(u, v) = 0,
\end{equation}
which holds because
$N_{ij}(u,v) = -e^{(x_j - x_i)\pa^2}(u,v)$
depends on $u - v$ alone.  Property~(iii) follows
because each factor $N_{ij}(u,v)$ depends on the
spatial sites only through the difference
$x_j - x_i$; a uniform shift of all $x_p$ leaves
every factor unchanged.
\end{proof}

For each $(t, \mathbf{x}, \mathbf{a})$ with $t > 0$,
define the product-graph wave functions
$\Psi^{t,\mathbf{x},\mathbf{a}}
\colon \R^m \to L^2(\R)$ and
$\Phi^{t,\mathbf{x},\mathbf{a}}
\colon L^2(\R) \to \R^m$ componentwise by
\begin{align}
\Psi_p^{t,\mathbf{x},\mathbf{a}}
&\defeq \sum_{k=1}^p
\int_{a_k}^{\infty}
[B_{\mathbf{x},\mathbf{a}}]_{kp}(z, a_p)
\psi_{t,x_k,z}\diff z,
\label{eq:FP-Psi-def} \\
\Phi_p^{t,\mathbf{x},\mathbf{a}}
&\defeq \sum_{k=p}^m
\int_{a_k}^{\infty}
[B_{\mathbf{x},\mathbf{a}}]_{pk}(a_p, w)
\phi_{t,x_k,w}\diff w.
\label{eq:FP-Phi-def}
\end{align}
Under the identification
$u = (t, \mathbf{x}, \mathbf{a})$, the maps
$u \mapsto \Psi^{t,\mathbf{x},\mathbf{a}}$ and
$u \mapsto \Phi^{t,\mathbf{x},\mathbf{a}}$ are the
wave functions of the continuum framework
(\S\ref{sec:kpz-reduced-framework}) with
$C = \Lambda = 0$.  The superscript
$(t, \mathbf{x},\mathbf{a})$ is suppressed when the
vertex is clear from context.  Write
$\pa_{\mathbf{x}} \defeq \sum_i \pa_{x_i}$ and
$\pa_{\mathbf{a}} \defeq \sum_i \pa_{a_i}$.

\begin{lemma}\label{lem:FP-WF}
The product-graph wave functions satisfy the linear
problem
\eqref{eq:kpz-red-Lx}--\eqref{eq:kpz-red-Lt} and its
adjoint \eqref{eq:kpz-red-Ax}--\eqref{eq:kpz-red-At}
with $C = \Lambda = 0$.  Explicitly:
\begin{align}
\pa_{\mathbf{x}} \Psi^{t,\mathbf{x},\mathbf{a}}
&= \pa_{\mathbf{a}}^2 \Psi^{t,\mathbf{x},\mathbf{a}},
\label{eq:FP-WF-heat} \\
\pa_t \Psi^{t,\mathbf{x},\mathbf{a}}
&= -\tfrac{1}{3}\pa_{\mathbf{a}}^3
   \Psi^{t,\mathbf{x},\mathbf{a}},
\label{eq:FP-WF-airy}
\end{align}
and
\begin{align}
\pa_{\mathbf{x}} \Phi^{t,\mathbf{x},\mathbf{a}}
&= -\pa_{\mathbf{a}}^2 \Phi^{t,\mathbf{x},\mathbf{a}},
\label{eq:FP-WF-phi-heat} \\
\pa_t \Phi^{t,\mathbf{x},\mathbf{a}}
&= -\tfrac{1}{3}\pa_{\mathbf{a}}^3
   \Phi^{t,\mathbf{x},\mathbf{a}}.
\label{eq:FP-WF-phi-airy}
\end{align}
\end{lemma}

\begin{proof}
By linearity and the definitions
\eqref{eq:FP-Psi-def}--\eqref{eq:FP-Phi-def}, it
suffices to verify the identities componentwise for
$[\psi_{\mathbf{x},\mathbf{a}}]_{kp}
\defeq \int_{a_k}^{\infty}
[B_{\mathbf{x},\mathbf{a}}]_{kp}(z, a_p)
\psi_{t,x_k,z}\diff z$.

\noindent\textit{Verification of \eqref{eq:FP-WF-heat}.}
The propagator satisfies
$\pa_{\mathbf{x}}[B_{\mathbf{x},\mathbf{a}}]_{kp} = 0$
by \eqref{eq:FP-x-indep}, and the seed $\psi_{t,x_k,z}$
depends only on the single site $x_k$.  The
$\pa_{\mathbf{x}}$-derivative of the integral therefore
acts only through $\pa_{x_k}$ on $\psi$, and the seed
equation \eqref{eq:FP-seed-heat} gives
\[
\pa_{\mathbf{x}} [\psi_{\mathbf{x},\mathbf{a}}]_{kp}
= \int_{a_k}^{\infty}
[B_{\mathbf{x},\mathbf{a}}]_{kp}(z, a_p)
\pa_{x_k} \psi_{t,x_k,z}\diff z
= \int_{a_k}^{\infty}
[B_{\mathbf{x},\mathbf{a}}]_{kp}(z, a_p)
\pa_z^2 \psi_{t,x_k,z}\diff z.
\]
On the other hand, differentiating by
$\pa_{\mathbf{a}}$ and evaluating the lower limit at
$a_k$ by the Leibniz rule gives
\[
\pa_{\mathbf{a}} [\psi_{\mathbf{x},\mathbf{a}}]_{kp}
= -[B_{\mathbf{x},\mathbf{a}}]_{kp}(a_k, a_p)
\psi_{t,x_k,a_k}
+ \int_{a_k}^{\infty}
\Bigl(\sum_{k < s \leq p} \pa_{a_s}
[B_{\mathbf{x},\mathbf{a}}]_{kp}(z, a_p)\Bigr)
\psi_{t,x_k,z}\diff z.
\]
To transfer the derivative from $B$ to $\psi$, flux
conservation \eqref{eq:FP-flux} replaces the
$\pa_{\mathbf{a}}$-sum on $B$ by $-\pa_z B$:
\[
= -[B_{\mathbf{x},\mathbf{a}}]_{kp}(a_k, a_p)
\psi_{t,x_k,a_k}
- \int_{a_k}^{\infty}
\pa_z [B_{\mathbf{x},\mathbf{a}}]_{kp}(z, a_p)
\psi_{t,x_k,z}\diff z.
\]
Integration by parts in $z$ (the boundary at $\infty$
vanishes by decay of $B$) gives
\begin{align*}
= -[B_{\mathbf{x},\mathbf{a}}]_{kp}(a_k, a_p)
\psi_{t,x_k,a_k}
+ [B_{\mathbf{x},\mathbf{a}}]_{kp}(a_k, a_p)
\psi_{t,x_k,a_k}
+ \int_{a_k}^{\infty}
[B_{\mathbf{x},\mathbf{a}}]_{kp}(z, a_p)
\pa_z \psi_{t,x_k,z}\diff z.
\end{align*}
The boundary terms cancel, leaving
\begin{equation}\label{eq:FP-passage}
\pa_{\mathbf{a}} [\psi_{\mathbf{x},\mathbf{a}}]_{kp}
= \int_{a_k}^{\infty}
[B_{\mathbf{x},\mathbf{a}}]_{kp}(z, a_p)
\pa_z \psi_{t,x_k,z}\diff z.
\end{equation}
The derivation of \eqref{eq:FP-passage} used only the
flux identity for $B$ and integration by parts;
$\psi_{t,x_k,z}$ entered as an arbitrary $L^2$
function of $z$.
Replacing $\psi_{t,x_k,z}$ by $\pa_z\psi_{t,x_k,z}$
and repeating the argument gives
\[
\pa_{\mathbf{a}}^2 [\psi_{\mathbf{x},\mathbf{a}}]_{kp}
= \int_{a_k}^{\infty}
[B_{\mathbf{x},\mathbf{a}}]_{kp}(z, a_p)
\pa_z^2 \psi_{t,x_k,z}\diff z
= \pa_{\mathbf{x}}
[\psi_{\mathbf{x},\mathbf{a}}]_{kp},
\]
confirming \eqref{eq:FP-WF-heat}.

\medskip
\noindent\textit{Verification of \eqref{eq:FP-WF-airy}.}
The propagator is $t$-independent, so $\pa_t$ passes
through the integral and acts only on $\psi$.
Equation \eqref{eq:FP-seed-airy} and three applications
of \eqref{eq:FP-passage} give
\begin{align*}
\pa_t [\psi_{\mathbf{x},\mathbf{a}}]_{kp}
&= \int_{a_k}^{\infty}
[B_{\mathbf{x},\mathbf{a}}]_{kp}(z, a_p)
\pa_t \psi_{t,x_k,z}\diff z \\
&= -\tfrac{1}{3}\int_{a_k}^{\infty}
[B_{\mathbf{x},\mathbf{a}}]_{kp}(z, a_p)
\pa_z^3 \psi_{t,x_k,z}\diff z
= -\tfrac{1}{3}\pa_{\mathbf{a}}^3
[\psi_{\mathbf{x},\mathbf{a}}]_{kp},
\end{align*}
confirming \eqref{eq:FP-WF-airy}.

\medskip
\noindent\textit{Verification of
\eqref{eq:FP-WF-phi-heat}--\eqref{eq:FP-WF-phi-airy}.}
The same argument applies to
$[\phi_{\mathbf{x},\mathbf{a}}]_{pk}
= \int_{a_k}^{\infty}
[B_{\mathbf{x},\mathbf{a}}]_{pk}(a_p, w)
\phi_{t,x_k,w}\diff w$
from \eqref{eq:FP-Phi-def}, with the integration
variable $w$ in the second argument of $B$.  The
Leibniz-plus-integration-by-parts argument converts
each $\pa_{\mathbf{a}}$ to $\pa_w$ on $\phi$.
Equations
\eqref{eq:FP-phi-heat}--\eqref{eq:FP-phi-airy} then
give $\pa_{\mathbf{x}}[\phi]_{pk}
= -\pa_{\mathbf{a}}^2[\phi]_{pk}$ and
$\pa_t[\phi]_{pk}
= -\tfrac{1}{3}\pa_{\mathbf{a}}^3[\phi]_{pk}$;
the sign in the heat identity arises because
\eqref{eq:FP-phi-heat} reads
$\pa_x\phi = -\pa_a^2\phi$ rather than
$+\pa_a^2\phi$.
\end{proof}

The kernel $K_{t,\mathbf{x},\mathbf{a}}$ satisfies
three factorization identities that match the dressing
conditions of the continuum framework.

\begin{lemma}
\label{lem:FP-dressing}
The kernel $K_{t,\mathbf{x},\mathbf{a}}$ satisfies:
\begin{align}
\pa_{\mathbf{a}} K
&= -\Psi\Phi,
\label{eq:FP-DC-a} \\
\pa_{\mathbf{x}} K
&= -\pa_{\mathbf{a}}\Psi\Phi
   + \Psi\pa_{\mathbf{a}}\Phi,
\label{eq:FP-DC-x} \\
\pa_t K
&= \tfrac{1}{3}\bigl(
   \pa_{\mathbf{a}}^2\Psi\Phi
   - \pa_{\mathbf{a}}\Psi\pa_{\mathbf{a}}\Phi
   + \Psi\pa_{\mathbf{a}}^2\Phi\bigr).
\label{eq:FP-DC-t}
\end{align}
\end{lemma}

\begin{proof}
\noindent\textit{Verification of \eqref{eq:FP-DC-a}.}
Differentiating \eqref{eq:FP-kernel-formula} by
$\pa_{a_p}$ produces boundary contributions from the
lower limits $a_i$ and $a_j$ via the Leibniz rule, and
interior contributions where $\pa_{a_p}$ acts on $B$
through the splitting identity \eqref{eq:FP-splitting}.
All terms carry a minus sign (from
$\pa_{a_p}\int_{a_p}^{\infty} f = -f(a_p)$ and
from the sign in \eqref{eq:FP-splitting}):
\begin{align*}
\pa_{a_p} K
&= -\psi_{x_p,a_p}\phi_{x_p,a_p}
- \sum_{j > p} \psi_{x_p,a_p}[\phi]_{pj}
- \sum_{i < p} [\psi]_{ip}\phi_{x_p,a_p}
- \sum_{i < p < j} [\psi]_{ip}[\phi]_{pj} \\
&= -\bigl(\sum_{i \leq p} [\psi]_{ip}\bigr)
\bigl(\sum_{j \geq p} [\phi]_{pj}\bigr)
= -\Psi_p\Phi_p,
\end{align*}
where $[\psi]_{ip}$ and $[\phi]_{pj}$ denote the
dressed components from
\eqref{eq:FP-Psi-def}--\eqref{eq:FP-Phi-def}, with
$[\psi]_{pp} = \psi_{x_p,a_p}$ and
$[\phi]_{pp} = \phi_{x_p,a_p}$ from the diagonal
$[B]_{pp} = \delta$.  Summing over $p$ gives
\eqref{eq:FP-DC-a}.

\medskip
\noindent\textit{Verification of \eqref{eq:FP-DC-x}.}
The propagator satisfies
$\pa_{\mathbf{x}}[B]_{ij} = 0$ by
\eqref{eq:FP-x-indep}, and each seed depends only on
its own spatial site ($\psi_{x_i,z}$ on $x_i$,
$\phi_{x_j,w}$ on $x_j$).  Differentiating
\eqref{eq:FP-kernel-formula} and applying
\eqref{eq:FP-seed-heat} and \eqref{eq:FP-phi-heat}
gives
\[
\pa_{\mathbf{x}} K
= \sum_{i \leq j}
\int_{a_i}^{\infty}\int_{a_j}^{\infty}
\bigl[\pa_z^2\psi_{x_i,z}[B]_{ij}(z,w)
\phi_{x_j,w}
- \psi_{x_i,z}[B]_{ij}(z,w)
\pa_w^2\phi_{x_j,w}\bigr]
\diff w\diff z.
\]
To reduce this to boundary evaluations, the flux
identity \eqref{eq:FP-flux} rewrites the integrand as
a total divergence:
\begin{align*}
&\pa_z^2\psi_{x_i,z}[B]_{ij}\phi_{x_j,w}
- \psi_{x_i,z}[B]_{ij}\pa_w^2\phi_{x_j,w} \\
&\quad = \bigl(\pa_z + \pa_w
+ \textstyle\sum_{i<p<j}\pa_{a_p}\bigr)
\bigl[\pa_z\psi_{x_i,z}[B]_{ij}\phi_{x_j,w}
- \psi_{x_i,z}[B]_{ij}\pa_w\phi_{x_j,w}\bigr],
\end{align*}
which holds because expanding the right-hand side by
the Leibniz rule and applying
$(\pa_z + \pa_w + \sum\pa_{a_p})[B]_{ij} = 0$ from
\eqref{eq:FP-flux} leaves only the seed derivatives.
Integrating the divergence, each component evaluates
the full bracket at its boundary: $\pa_z$ at the lower
limit $z = a_i$, $\pa_w$ at $w = a_j$, and each
$\pa_{a_p}$ via the splitting identity
\eqref{eq:FP-splitting}.  For example, the $\pa_z$
component evaluates the first bracket entry at
$z = a_i$:
\[
-\int_{a_j}^{\infty}
\pa_{a_i}\psi_{x_i,a_i}[B]_{ij}(a_i,w)
\phi_{x_j,w}\diff w
= -\pa_{a_i}\psi_{x_i,a_i}[\phi]_{ij}
= -\pa_{\mathbf{a}}[\psi]_{ii}[\phi]_{ij},
\]
where the last equality is
\eqref{eq:FP-passage} at the diagonal $[B]_{ii}
= \delta$.  The $\pa_w$ and $\pa_{a_p}$ components
produce the same structure at indices $p = j$ and
$i < p < j$ respectively.  Collecting all boundary
contributions from both bracket entries:
\begin{align*}
\pa_{\mathbf{x}} K
= \sum_{i \leq p \leq j}
\bigl(-\pa_{\mathbf{a}}[\psi]_{ip}[\phi]_{pj}
+ [\psi]_{ip}\pa_{\mathbf{a}}[\phi]_{pj}\bigr)
= -\pa_{\mathbf{a}}\Psi\Phi
+ \Psi\pa_{\mathbf{a}}\Phi,
\end{align*}
where \eqref{eq:FP-passage}
identifies each boundary integral (with $\pa_z\psi$ or
$\pa_w\phi$ as seed) with the $\pa_{\mathbf{a}}$-derivative
of the corresponding dressed component.

\medskip
\noindent\textit{Verification of \eqref{eq:FP-DC-t}.}
The propagator is $t$-independent, so $\pa_t$ acts
only on $\psi$ and $\phi$.  Equations
\eqref{eq:FP-seed-airy} and \eqref{eq:FP-phi-airy}
give
\[
\pa_t K = -\tfrac{1}{3}\sum_{i \leq j}
\int_{a_i}^{\infty}\int_{a_j}^{\infty}
\bigl[\pa_z^3\psi_{x_i,z}[B]_{ij}(z,w)
\phi_{x_j,w}
+ \psi_{x_i,z}[B]_{ij}(z,w)
\pa_w^3\phi_{x_j,w}\bigr]
\diff w\diff z.
\]
The same divergence argument applies, with the bracket
now containing three entries (one degree higher in
each seed derivative):
\begin{align*}
&\pa_z^3\psi_{x_i,z}[B]_{ij}\phi_{x_j,w}
+ \psi_{x_i,z}[B]_{ij}\pa_w^3\phi_{x_j,w} \\
&\quad = \bigl(\pa_z + \pa_w
+ \textstyle\sum_p\pa_{a_p}\bigr)
\bigl[\pa_z^2\psi_{x_i,z}[B]_{ij}\phi_{x_j,w}
- \pa_z\psi_{x_i,z}[B]_{ij}\pa_w\phi_{x_j,w}
+ \psi_{x_i,z}[B]_{ij}\pa_w^2\phi_{x_j,w}\bigr].
\end{align*}
Integrating and evaluating boundary terms as in the
verification of \eqref{eq:FP-DC-x}, each bracket entry
assembles into the corresponding dressed product: the
entry $\pa_z^2\psi_{x_i,z}[B]_{ij}\phi_{x_j,w}$
gives $\pa_{\mathbf{a}}^2\Psi\Phi$, the entry
$-\pa_z\psi_{x_i,z}[B]_{ij}\pa_w\phi_{x_j,w}$
gives $-\pa_{\mathbf{a}}\Psi\pa_{\mathbf{a}}\Phi$, and
$\psi_{x_i,z}[B]_{ij}\pa_w^2\phi_{x_j,w}$
gives $\Psi\pa_{\mathbf{a}}^2\Phi$, so
\[
\pa_t K = \tfrac{1}{3}\bigl(\pa_{\mathbf{a}}^2\Psi\Phi
- \pa_{\mathbf{a}}\Psi\pa_{\mathbf{a}}\Phi
+ \Psi\pa_{\mathbf{a}}^2\Phi\bigr).
\]
\end{proof}

\section{Multipoint equation}
\label{sec:FP-multipoint}

The dressed observable
$\mathcal{A} = z\Phi R\Psi$ of
\eqref{eq:kpz-red-A-def} specializes at $z = 1$ to
the $m \times m$ matrix
\begin{equation}\label{eq:FP-A-def}
\mathcal{A} \defeq \Phi(I-K)^{-1}\Psi
\in \End(\R^m),
\qquad
(\mathcal{A})_{ij}
= \langle (I-K)^{-1}\Psi_j, \Phi_i
\rangle_{L^2(\R)}.
\end{equation}

By Lemma~\ref{lem:FP-kernel},
$\det_{L^2(\R)}(I - K_{t,\mathbf{x},\mathbf{a}})$
coincides with the extended-space determinant
$\det_{\mathscr{X}}(I - \chi_{\mathbf{a}} K_t
\chi_{\mathbf{a}})$, which Quastel and
Remenik~\cite{QR2022} show is nonzero for all
finite $(t, \mathbf{x}, \mathbf{a})$ with $t > 0$.
The resolvent
$(I - K_{t,\mathbf{x},\mathbf{a}})^{-1}$ therefore
exists throughout.

\begin{theorem}\label{thm:FP-pKP}
At every $(t, \mathbf{x}, \mathbf{a})$ with $t > 0$,
the dressed observable $\mathcal{A}$ satisfies
the matrix potential KP equation
\begin{equation}\label{eq:FP-pKP}
\pa_{t,\mathbf{a}} \mathcal{A}
+ \tfrac{1}{4}\pa_{\mathbf{x}}^2 \mathcal{A}
+ \tfrac{1}{12}\pa_{\mathbf{a}}^4 \mathcal{A}
+ \tfrac{1}{2}\pa_{\mathbf{a}}
  (\pa_{\mathbf{a}}\mathcal{A})^2
+ \tfrac{1}{2}[\pa_{\mathbf{a}}\mathcal{A},
\pa_{\mathbf{x}}\mathcal{A}] = 0.
\end{equation}
\end{theorem}

\begin{proof}
By Lemmas~\ref{lem:FP-WF} and~\ref{lem:FP-dressing},
the triple $(\Psi, \Phi, K)$ satisfies the hypotheses
of Theorem~\ref{thm:kpz-red-Darboux-compat} with
$C = \Lambda = 0$.  Setting $C = \Lambda = 0$ in
\eqref{eq:kpz-shift-background-kp} gives
\eqref{eq:FP-pKP}.
\end{proof}

\begin{corollary}\label{cor:FP-trace-defect}
Let $F_{t,\mathbf{x},\mathbf{a}}
\defeq \det_{L^2(\R)}(I - K_{t,\mathbf{x},\mathbf{a}})$.
Then
$\pa_{\mathbf{a}} \log F = \operatorname{tr}_{\R^m} \mathcal{A}$,
and
\begin{equation}\label{eq:FP-Hirota}
\Bigl[D_{t,\mathbf{a}}
+ \tfrac{1}{4}D_{\mathbf{x}}^2
+ \tfrac{1}{12}D_{\mathbf{a}}^4\Bigr] F \cdot F
= F^2\bigl((\operatorname{tr}_{\R^m}
\pa_{\mathbf{a}}\mathcal{A})^2
- \operatorname{tr}_{\R^m}
((\pa_{\mathbf{a}}\mathcal{A})^2)\bigr),
\end{equation}
where $D_{\mathbf{a}}$, $D_{\mathbf{x}}$, $D_t$
denote Hirota derivatives with respect to the
collective variables.
\end{corollary}

\begin{proof}
Jacobi's formula and \eqref{eq:FP-DC-a} give
$\pa_{\mathbf{a}} \log F
= \operatorname{tr}_{\R^m}\mathcal{A}$
by cyclicity of the trace.  The hypotheses of
Proposition~\ref{prop:kpz-shift-trace-defect} hold:
trace-class and invertibility are inherited from
Matetski--Quastel--Remenik~\cite{MQR21}, and
$F \to 1$ as $\mathbf{a} \to +\infty$ because
$K \to 0$ in trace norm.  Setting $C = 0$ in
\eqref{eq:kpz-shift-integrated-Hirota} gives
\eqref{eq:FP-Hirota}.
\end{proof}

\begin{corollary}\label{cor:FP-one-point}
For $m = 1$, the Fredholm determinant
$F_{t,x,a} = \det_{L^2(\R)}(I - K_{t,x,a})$ satisfies
\begin{equation}\label{eq:FP-scalar-pKP}
\Bigl[D_{t,a} + \tfrac{1}{4}D_x^2
+ \tfrac{1}{12}D_a^4\Bigr] F \cdot F = 0,
\end{equation}
which is KP-II in Hirota bilinear form.
\end{corollary}

\begin{proof}
When $m = 1$, $\pa_a\mathcal{A}$ is scalar, so the
right-hand side of \eqref{eq:FP-Hirota} vanishes.
\end{proof}